\documentclass[10pt]{report}
\usepackage[utf8]{inputenc}
\usepackage{stackengine,latexsym,makeidx,amssymb,amsmath,theorem,url}
\usepackage[makeroom]{cancel}
\usepackage{imakeidx}
\usepackage{yfonts}

\newtheorem{lemma}{Lemma}[section]
\newtheorem{thm}[lemma]{Theorem}
\newtheorem{cor}[lemma]{Corollary}
\newtheorem{prop}[lemma]{Proposition}
\newtheorem{defi}[lemma]{Definition}

\newtheorem{prob}[lemma]{Problem}
\newtheorem{axio}[lemma]{Axiom}
\newtheorem{exer}[lemma]{Exercise}

{\theorembodyfont{\rmfamily}}

\newcommand{\N}{\mathbb{N}}
\newcommand{\Z}{\mathbb{Z}}
\newcommand{\Q}{\mathbb{Q}}
\newcommand{\R}{\mathbb{R}}

\newcommand{\C}{\mathbb{C}}
\newcommand{\e}{\mathrm{e}}

\newcommand{\ep}{\varepsilon}
\newcommand{\de}{\delta}
\newcommand{\al}{\alpha}
\newcommand{\be}{\beta}
\newcommand{\ga}{\gamma}
\newcommand{\duk}{\par\noindent{\bf Proof. }}
\newcommand{\kduk}{\hfill{$\Box$}\bigskip}

\newcommand{\sus}{\subset}

\newcommand{\tec}{\hspace{-1.6mm}{\bf. }}

\newcommand{\sgn}{\mathrm{sgn}}
\newcommand{\cc }{\colon}

\def\ds{\dots}

\DeclareFontFamily{U}{wncy}{}
    \DeclareFontShape{U}{wncy}{m}{n}{<->wncyr10}{}
    \DeclareSymbolFont{mcy}{U}{wncy}{m}{n}
    \DeclareMathSymbol{\Sh}{\mathord}{mcy}{"58} 

\makeindex[title={Index\\
\vspace{8mm}
{\small{\rm Page number in {\em italic} refers to the definition of a~notion or to a~theorem with proof.\\
}}}\vspace{-7mm}]
\begin{document}
\pagenumbering{roman}
\thispagestyle{empty}
\begin{center}
{\Huge {\bf Elementary Real Analysis}}
    
\bigskip\bigskip
{\huge Martin Klazar}

\bigskip
{\large (KAM MFF UK Praha)}

\end{center}
    
\bigskip\bigskip\noindent
Preliminary version of October~31, 2026 

\newpage
\thispagestyle{empty}
\noindent
dedicated to the memory of
Ji\v r\'\i\ Matou\v sek (1963--2015)\index{matousek@Matou\v sek, Ji\v r\'\i}

\newpage
\thispagestyle{empty}
\noindent
{\em Was sich \"uberhaupt sagen l\"asst, l\"asst sich klar sagen;} $\ds$\footnote{We omit the well known 
conclusion of Wittgenstein's quote ``Was sich \"uberhaupt sagen l\"asst, l\"asst sich 
klar sagen; und wovon man nicht reden kann, dar\"uber muss man schweigen''.  ``{\em What can be said at all, can be 
said clearly;} and what one cannot talk about, must be left in silence''.}

\bigskip\noindent
L.~Wittgenstein\index{Wittgenstein, Ludwig} \cite[Vorwort]{witt}

\tableofcontents
\newpage

\section*{Introduction}
\addcontentsline{toc}{chapter}{Introduction}

This book is an extension of the lecture notes for the course {\em Mathematical Analysis~1}, which I was teaching in the School of Computer 
Science of MFF UK\footnote{See the index of notation starting on p.~\pageref{znaceni}.} in Prague for several years. When I am done with it, I hope to produce the Czech version. The book 
has nineteen chapters and six appendices. Each chapter and appendix (except 
Appendix~\ref{doda_reseni}) begins with 
a~summary of main points and is divided into sections. Each section 
in a~chapter contains at least one exercise. Solutions or hints to solutions of 
these exercises can be found in Appendix~\ref{doda_reseni}. 

Chapter~\ref{chap_funkArela} is devoted to a~review of notions and notation 
pertaining to functions and relations. The main point is perhaps our definition of functions and their isomorphisms.
Starting in Appendix~\ref{doda_notation} and 
continuing in Chapters~\ref{sec_ZFCtoN0}--\ref{sec_RstarToR}, we build the real 
numbers. In Appendix~\ref{doda_notation}, we introduce enough ZFC set theory so that  in Chapter~\ref{sec_ZFCtoN0} we 
can build the numeric domain of natural numbers. In Chapters~\ref{sec_N0toZ} and \ref{sec_ZtoQ} we introduce the numeric 
domains of integers and rational numbers. In Chapter~\ref{sec_QtoRstar}, which is 
our original contribution, we show that the numeric domain $\R^*=\R\cup
\{-\infty,+\infty\}$ of extended reals 
arises in a~natural way as a~limit completion of the linear order of fractions introduced in 
Chapter~\ref{sec_ZtoQ}. We finish the buildup of $\R$ in 
Chapter~\ref{sec_RstarToR}. In Appendix~\ref{sec_RtoC}, we extend $\R$ to $\C$. Our main results are the construction of $\R^*$ 
and the characterization theorems for the numeric domains: up to isomorphism, $\N_0$ is the unique simple ordered 
semiring, $\Z$~is the unique simple ordered ring, $\Q$ is the unique 
simple ordered field, $\R^*$ is the unique limit completion of the linear 
order $\Q$, and $\C$ is the unique quadratic extension of $\R$. The 
characterization theorem for $\R$ is well known; it is the unique complete ordered 
field. 

In Chapter~\ref{chap_existLim}

I dedicate the book to the memory of 
my late colleague
Ji\v r\'\i\ Matou\v sek\index{matousek@Matou\v sek, Ji\v r\'\i}. He was one of the greatest Czech mathematicians and 
computer scientists. He was preparing lectures for the analysis course in the school year 
2014/15, but illness did not allow him to deliver them.

\bigskip\bigskip\noindent
Praha and Louny, August 2024 to ??? 2026\hfill{Martin Klazar}

\newpage

\pagenumbering{arabic}
\chapter[Prologue on functions and relations\checkmark]{Prologue on functions and relations}\label{chap_funkArela}

The main actors in the arena of elementary real analysis are real 
functions of one real variable. We therefore begin this book
with a~review of notions pertaining to functions and relations. In the following five 
chapters, we construct the complete field of real numbers. We need so many 
chapters for it because en route we construct, as we have to, the numeric domains of 
natural numbers, integers, fractions, and extended reals. All of them, 
especially the last one, are frequently used in elementary real analysis.

In this chapter, 
Definitions~\ref{def_funkce} and \ref{def_rovnostFci} introduce  
functions and their isomorphisms.
Definition~\ref{def_seWoOp} introduces sequences, words, and operations. 
Definition~\ref{def_homom} defines homomorphisms. Composite functions are 
introduced in
Definition~\ref{def_compFun}. The axiom of choice appears as 
Axiom~\ref{axio_AC}. Equivalence relations and the dual set partitions appear in 
Definitions~\ref{def_relEkv} and \ref{def_rozkl}, respectively.
Definition~\ref{def_LO} introduces linear orders. Suprema and infima in linear orders appear in 
Definition~\ref{def_supInf}. Finally,
Definitions~\ref{def_wellOrd} and \ref{def_orders} review well 
orderings and partial orders, respectively. Many other notions are 
introduced and reviewed outside formal definitions.   

\section{Functions and their isomorphisms}

We define functions and their isomorphisms. See 
Definition~\ref{def_kTup} for $k$-tuples
$$
\langle x_1,\,x_2,\,\ds,\,x_k\rangle\ \  
(k\in\N=\{1,2,\ds\})\,. 
$$
We precisely define the set $\N$ in the next chapter. The 
\underline{Cartesian product\index{Cartesian product of sets|emph}} of sets $A$ and $B$ is the set
$$
A\times B:=\{\langle a,\,b\rangle\cc\;a\in A,\,b\in B\}\,.\label{cartProd}
$$
Any $C\sus A\times B$ is 
a~binary \underline{relation between\index{relation@(binary) relation|emph}} 
$A$ and $B$. Instead of 
$\langle a,\,b\rangle\in C$ 
we write $a\,C\,b$, for example $2<5$. If $A=B$, we say that $C$ is 
a~\underline{relation on}\index{relation@(binary) relation!on 
a~set|emph} $A$. A~relation\index{relation@(binary) relation!functional|emph} 
$C$ between $A$ and $B$ is \underline{functional} if for every 
element $a\in A$ there is exactly one element $b\in B$ such that $a\,C\,b$.

\begin{defi}\label{def_funkce} 
Let $A$ and $B$ be sets.
A~\underline{function} or a~\underline{map\index{function|emph}}
$f$ from $A$ to $B$, written $f\cc A\to B$\label{function}, is  
any ordered triple 
$$
f=\langle A,\,B,\,G_f\rangle 
$$
such that $G_f$, the \underline{graph\index{graph of 
$f$|emph}} of $f$, is a~functional relation between $A$ and $B$. We call 
$A$ the 
\underline{domain\index{function!domain of|emph}} of $f$ and denote it by $M(f)$\index{function!mf@$M(f)$|emph}\label{Mf}. $B$ is the \underline{range\index{function!range|emph}} of $f$. If $a\in A$, we write 
$f(a)$ for the unique $b\in B$ such that $\langle a,b\rangle\in G_f$. We say that $f(a)$ is the 
\underline{value\index{function!value of|emph}} of $f$ 
at the \underline{argument\index{function!argument of|emph}} $a$.
\end{defi}
Thus 
$G_f=\{\langle x,f(x)\rangle\cc\;x\in A\}\label{graphFun}$.
Instead of $f\cc A\to B$ we sometimes write 
$$
A\ni a\mapsto \varphi(a)\in B\,,  
$$
where $\varphi(a)$ is a~formula producing the values $f(a)$. 

Any description of mathematical objects should 
include a~description of their isomorphisms. Two structures are 
isomorphic if they are practically identical, even if they differ as sets. 
Isomorphisms are always
equivalences (see Definition~\ref{def_relEkv}).
We introduce isomorphisms of functions. 

\begin{defi}\label{def_rovnostFci}
Two functions are \underline{isomorphic\index{isomorphism of functions|emph}}, if they differ as triples only 
in their ranges.  
\end{defi}

\begin{exer}\label{ex_rovnFci}
Functions $f=\langle A,B,G_f\rangle$ and $g=\langle C,D,G_g\rangle$ are isomorphic if and only if $G_f=G_g$.
\end{exer}

\section{Partial functions and empty functions}

Let $A$ and $B$ be sets. 
A~relation $C$ between $A$ and $B$ is a~\underline{partial 
function\index{partial function|emph}}
from $A$ to $B$ if for every $a\in A$ there is at most one 
$b\in B$ such that $a\,C\,b$. Partial functions may be undefined at 
some elements of their domains.
They are common in recursion 
theory: if an algorithm does not terminate on an input $x$, then the 
associated partial function is not defined at $x$. Other examples of partial functions are the 
subtraction $m-n$ of natural numbers and the division $a/b$ in a~field $F$. The 
former is not defined if $n>m$, and the latter if $b=0_F$. 
We define both operations in more detail later.

\begin{exer}\label{ex_partFun}
What is the relation of functions and partial functions?    
\end{exer}

What is an empty function? For example, the function
$$
f_{\emptyset}:=\langle\emptyset,\,
\emptyset,\,\emptyset\rangle\label{empFun}\,.
$$
A~function is \underline{empty\index{empty function|emph }\index{function!empty|emph}} if it is isomorphic to $f_{\emptyset}$.

\begin{exer}\label{ex_onEmpFun}
Empty functions are functions $\langle\emptyset,X,\emptyset\rangle$.
\end{exer}

\section{Images, preimages, restrictions}

Let $f\cc A\to B$ be 
a~function and $C$ be {\em any} set. We define sets
\begin{eqnarray*}
f[C]&:=&\{f(a)\cc\;a\in C\cap A\}\ (\sus B)\text{\,---\,the \underline{image\index{function!image of a~set, $f[C]$|emph}}\label{image} of $C$ by $f$\,---\,and}\\
f^{-1}[C]&:=&\{a\in A\cc\;f(a)\in C\}\ (\sus A)\text{\,---\,the \underline{preimage\index{function!preimage of a~set, $f^{-1}[C]$|emph}}\label{preimage} of $C$ by $f$}\,.
\end{eqnarray*}
$C$ is arbitrary, it is practical to drop the usual restrictions $C\sus A$ or $C\sus B$. The set $f[A]$ ($\sus B$) is the \underline{image\index{function!image of|emph}} of the function $f$.

\begin{exer}\label{ex_obrazAvzor}
Is it true that $f^{-1}[f[C]]=C$ and that $f[f^{-1}[C]]=C$?
\end{exer}

Let $f\cc A\to B$ be a~function and $C$ be any set. The 
\underline{restriction\index{function!restriction!to a~set, 
$f\,"|\,C$|emph}} of $f$ to $C$ is the 
function $f\,|\,C\cc A\cap C\to B$ such that for every $x\in A\cap C$,
$$
(f\,|\,C)(x):=f(x)\,.\label{restric}
$$
In this situation, we say that $f$  \underline{extends\index{function!extension of|emph}} $f\,|\,C$. It 
is clear that if functions $f_1$ and $f_2$ are isomorphic, then for any set $C$ the 
restrictions $f_1\,|\,C$ 
and $f_2\,|\,C$ are isomorphic as well. We say that a~function $f$ is 
a~\underline{subfunction\index{function!subfunction|emph}} of a~function $g$ if $G_f\sus G_g$. If $f$ and $g$ are functions, $C$ is any 
set and the restrictions $f\,|\,C$ and $g\,|\,C$ are isomorphic, we say that 
\underline{$f=g$ on} $C$\index{equality of $f$ and $g$ on $C$|emph}.

\begin{exer}\label{ex_restEmpFun}
Every empty function is a~subfunction of every function.     
\end{exer}

\section{Sequences, words, operations}

Let $\omega=\{0,1,\ds\}$ be the set of natural numbers introduced in the next 
chapter. Recall that 
$\N=\omega\setminus\{0\}=\{1,2,\ds\}$.\label{natural} For $n\in\N$, let $[n]:=
\{1,2,\ds,n\}$. We set $[0]:=\emptyset$. Let $X$ be any
set. 
\begin{defi}\label{def_seWoOp}
Three important families of functions are sequences, 
words, and operations.
\begin{enumerate}
\item Functions $a\cc\N\to X$ are \underline{sequences\index{function!sequence, $(a_n)\sus X$|emph}\index{sequence in $X$|emph}} in $X$. 
\item Functions $u\cc[n]\to X$, for $n\in\omega$, are \underline{words\index{function!word|emph}\index{words!over $X$|emph}} over 
the \underline{alphabet\index{alphabet|emph}} $X$.
\item Functions $o\cc X\times X\to X$ are binary \underline{operations\index{function!operation|emph}\index{operation on $X$|emph}} on $X$.
\end{enumerate}
\end{defi}
For a~sequence $a$ in $X$ and an \underline{index\index{sequence in $X$!index|emph}} $n\in\N$, we set $a_n:= a(n)$. We invoke a~sequence in $X$  
by writing 
$$
(a_n)\sus X\,. 
$$
A~sequence $(b_n)\sus X$ is a~\underline{subsequence\index{sequence in $X$!subsequence|emph}} of 
$(a_n)$ if there exists a~sequence of numbers $m_1<m_2<\ds$ in $\N$  
such that $b_n=a_{m_n}$ for every $n\in\N$.

We write a~word $u$ over $X$ as 
$$
u=u_1u_2\ds u_n\,,
$$ 
where $u_i:= u(i)$, $i\in[n]$. A~word $v=v_1v_2\ds v_m$ is 
a~\underline{subword\index{words!subwords of|emph}} of $u$ if $v_1=u_i$, 
$v_2=u_{i+1}$, $\ds$, $v_m=u_{i+m-1}$
for some $i\in[n]$ (hence $m\le n$).
For $n=0$ we have the \underline{empty word\index{empty word|emph}} over $X$
$$
u_{\emptyset}:=\langle\emptyset,\,X,\,\emptyset\rangle\,,
$$
which is an empty function; $u_{\emptyset}$ is a~subword of every 
word. The elements in the alphabet $X$, and sometimes also the terms 
$u_i$ of $u$, are called \underline{letters\index{alphabet!letter of|emph}}. The set of words over 
$X$ is denoted by $X^*$ 
and the \underline{length\index{words!length 
of|emph}} $n$ ($\in\omega$) of the word $u$ is denoted by $|u|$. If 
$u=u_1\ds u_m$ and $v=v_1\ds v_n$ are words in $X^*$, their 
\underline{concatenation\index{words!concatenation of|emph}} is the word $uv:=w_1\ds w_{m+n}\in X^*$ given by
$$
\text{$w_i=u_i$ for $i\in[m]$ and 
$w_i=v_{i-m}$ for $i\in[m+n]\setminus[m]$}\,.
$$
We can iterate the concatenation operation and speak of the 
concatenation of several (but always finitely many) words. In view 
of Exercises~\ref{ex_assConc} and \ref{ex_comConc}, only the order of 
these words is relevant.

The value 
$$
o(\langle a,\,b\rangle)=c
$$
of an operation $o\cc X\times X\to X$
is written as $a\,o\,b=c$, for example $1+1=2$. An element $b\in X$ 
such that $b\,o\,a=a$ for every $a\in X$ is called 
\underline{neutral\index{neutral element|emph}} to~$o$. The operation $o$ is
\underline{commutative\index{function!operation!commutative|emph}} if 
always $a\,o\,b=b\,o\,a$. It is 
\underline{associative\index{function!operation!associative|emph}} if 
always $(a\,o\,b)\,o\,c=a\,o\,(b\,o\,c)$. Another operation $p$ on 
$X$ is \underline{distributive\index{function!operation!distributive|emph}} to the 
operation $o$ if always $a\,p\,(b\,o\,c)=
(a\,p\,b)\,o\,(a\,p\,c)$. 
We follow the usual convention that in the absence of 
brackets, the operation denoted by $\cdot$ has precedence over $+$. The distributivity of $\cdot$ to $+$ is then 
written as $a\cdot(b+c)=a\cdot b+a\cdot c$. We abbreviate $a\cdot b$ by $ab$. A~function 
$o\cc X\to X$ is called a~\underline{unary
\index{function!operation!unary|emph} operation} on $X$. 

\begin{exer}\label{ex_uloNeuPr}
Elements neutral to commutative operations are uniquely determined.   
\end{exer}

\begin{exer}\label{ex_onFamFun}
Count words of length five over a~three-element alphabet.    
\end{exer}

\begin{exer}\label{ex_assConc}
Let $X$ be any alphabet. The concatenation operation on $X^*$ is 
associative.    
\end{exer}

\begin{exer}\label{ex_comConc}
However, it is not commutative.
\end{exer}

\section{Injective and other functions}

Let $f\cc X\to Y$ be a~function. 
\begin{enumerate}
\item We say that $f$ is \underline{injective\index{function!injective|emph}} or an injection if 
$x\ne x'$ $\Rightarrow$ $f(x)\ne f(x')$. 
\item We say that $f$ is \underline{surjective\index{function!surjective, onto|emph}} or onto or a~surjection if $f[X]=Y$.
\item The function $f$ is \underline{bijective\index{function!bijective|emph}} or a~bijection if $f$ is surjective and injective. 
\item The function $f$ is \underline{constant\index{function!constant|emph}} if it has exactly one value.
\item The function $f$ is an \underline{identity\index{function!identity, $\mathrm{id}_X$|emph}} function if $f(x)=x$ for every $x\in X$.
\end{enumerate}
More narrowly, the \underline{identity} on a~set $X$ is the map $\mathrm{id}_X\cc X\to X$ given by $\mathrm{id}_X(x)=x$.\label{identityFun}

\begin{exer}\label{ex_uloNaDrFunk}
When is the identity function from $X$ to $Y$ bijective?   
\end{exer}

\begin{exer}\label{ex_drFciArovn}
Suppose that two functions are isomorphic by
Definition~\ref{def_rovnostFci}. Is it true that they both are, or are not, injective? Same question for
surjective, bijective, constant, and identity functions.
\end{exer}
We reveal one answer: bijective functions are not
preserved by isomorphism. Thus, being a~bijection is not an essential property of functions. 

\section{Homomorphisms}

We will use homomorphisms in the following chapters.

\begin{defi}\label{def_homom}
Suppose that $p_i$ for $i\in I$ are operations on a~set $X$ and  that 
$q_i$ for $i\in I$ are operations on a~set $Y$. 
A~\underline{homomorphism\index{homomorphism|emph}} to these operations is a~map $f\cc X\to Y$ such that for every $a,b\in X$ and every $i\in I$,
$$
f(a\,p_i\,b)=f(a)\,q_i\,f(b)\,.
$$
If $f$ is bijective, we call $f$ 
an~\underline{isomorphism\index{isomorphism|emph}} to the mentioned operations. 
\end{defi}

\begin{exer}\label{ex_naHomom}
Let $f\cc X\to Y$ be an isomorphism to operations $p$ on $X$ and
$q$ on $Y$. If $a\in X$ is neutral to $p$, then $f(a)$ is neutral to $q$. What if $f$ is just a~homomorphism?
\end{exer}

\begin{exer}\label{ex_naHomom1}
Let $f$, $p$, and $q$ be as in the previous exercise. Then also 
$f^{-1}\cc Y\to X$ is an isomorphism to the operations $q$ and $p$.
\end{exer}

\begin{exer}\label{ex_naHomom2}
Let $f$, $p$, and $q$ be as before 
and let $g\cc Y\to Z$ be an isomorphism to operations $q$ and 
$r$, where $r$ is an operation on $Z$.
Then $g(f)$ is an isomorphism to $p$ and $r$.
\end{exer}

\section{Inverse functions and compositions}

Let $f\cc X\to Y$ be an injective map. The 
\underline{inverse\index{function!inverse|emph} function} 
$f^{-1}$ of $f$ is the function
$$
\text{$f^{-1}\cc f[X]\to X$, given by $f^{-1}(y)=x$ $\!\iff\!$ $f(x)=y$}\;.\label{inverse}
$$
Non-injective functions do not have inverses.

\begin{exer}\label{ex_jednoZnac}
Is it a~problem that $f^{-1}[A]$ means both the preimage of $A$ by~$f$, and the image of $A$ by $f^{-1}$? 
\end{exer}

\begin{exer}\label{ex_invInv}
Let $f$ be an injection. Are $f$ and $(f^{-1})^{-1}$ isomorphic? Are they equal as triples? Same questions for $f^{-1}$ and $((f^{-1})^{-1})^{-1}$.
\end{exer}

We introduce the kind of composition of functions that is used in 
mathematical analysis. 

\begin{defi}\label{def_compFun}
Let $g\cc X\to Y$ and $f\cc A\to B$ be any two functions. Their \underline{composition\index{function!composition, $f(g)$|emph}}\label{compos} $f(g)\cc X'\to B$ has the domain
$$
X'=\{x\in X\cc\;g(x)\in A\}\ \ (=g^{-1}[A])
$$
and has the values $f(g)(x):=f(g(x))$. We call
$g$ the \underline{inner\index{function!inner|emph} function}, and $f$ the \underline{outer\index{function!outer|emph} function}.    
\end{defi}
For example, if
$$
g(x)=1-x\cc\R\to\R\,\text{ and }\,f(x)=\sqrt{x}\cc[0,\,
+\infty)\to\R\,,
$$ 
then $f(g)(x)=\sqrt{1-x}\cc(-\infty,1]\to\R$ and $g(f)(x)=1-\sqrt{x}\cc[0,+\infty)\to\R$.

\begin{exer}\label{ex_rovnAskla}
Let $f_1,f_2$ and $g_1,g_2$ be pairs of isomorphic functions. Are 
the compositions $f_1(g_1)$ and $f_2(g_2)$ isomorphic? 
\end{exer}

\begin{exer}\label{ex_slozInjSurj}
Composition of two injections is an injection. The composition
$f(g)$ of surjections
$g\cc X\to Y$ and $f\cc Y\to B$ is a~surjection. In general, 
composition of two surjections need not be a~surjection.  
\end{exer}
The next exercise shows that the composition of functions is associative, in a~strong sense.
\begin{exer}\label{ex_asocSkla}
For every three functions 
$f$, $g$, and $h$, 
the two compositions $f(g(h))$ and $f(g)(h)$ are equal as triples. 
\end{exer}

\begin{exer}\label{ex_slozNaProste}
For every map $h\cc X\to Z$ there exists a~set $Y$ and maps
$g\cc X\to Y$ and $f\cc Y\to Z$ such that $h=f(g)$, $g$ is surjective, and $f$ is injective. 
\end{exer}
The next exercise describes the well known method for proving that a~function is bijective.
\begin{exer}\label{ex_oBijekci}
A~map $f\cc X\to Y$ is a~bijection $\iff$ there exists a~map $g\cc Y\to X$ such 
that $f(g)$ is $\mathrm{id}_Y$ and $g(f)$ is $\mathrm{id}_X$.   
\end{exer}

\begin{exer}\label{ex_kdyInvCon}
Which constant functions can be inverted?    
\end{exer}

\section{Set systems and the axiom of choice}\label{sec_AC}

Set systems are just maps. Let $I\ne\emptyset$ be a~set. 
A~\underline{set system\index{set system|emph}} 
$$
S=\{A_i\cc\;i\in I\}
$$ 
indexed by $I$ is just a~map $f\cc I\to B$ to some set $B$. The set corresponding to 
the index $i\in I$ is $A_i:= f(i)$. 
The \underline{union\index{set system!union of|emph}} of $S$ is 
$$
{\textstyle
\bigcup_{i\in I}A_i=
\bigcup\{A_i:\;i\in I\}:=
\bigcup f[I]\,.
}
$$
The \underline{intersection\index{set system!intersection of|emph}} of $S$ is
$$
{\textstyle
\bigcap_{i\in I}A_i=
\bigcap\{A_i:\;i\in I\}:=
\bigcap f[I]\,.
}
$$

\begin{exer}\label{ex_uloNaSS}
Explain the notation $\bigcup_{i\ge1}A_i$, 
$\bigcup_{i=1}^{\infty}A_i$, and 
$\bigcap_{i=0}^{\infty}A_i$.
\end{exer}

The next choice principle plays in mathematics an important role. 

\begin{axio}[AC]\label{axio_AC}
The \underline{axiom of choice\index{axiom!of choice, 
AC|emph}}\label{AC}, abbreviated {\em AC}, is any of the three following mutually 
equivalent statements. 
\begin{enumerate}
\item For every set system $\{A_i\cc i\in I\}$ with $A_i\ne\emptyset$ there exists 
a~\underline{selector\index{axiom!of choice, AC!selector|emph}}, a~map 
$${\textstyle
S\cc I\to\bigcup_{i\in I}A_i
}
$$ 
such that $S(i)\in A_i$ for every $i\in I$. 
\item For every set $X$ of nonempty and mutually disjoint sets there exists 
a~set~$Y$ such that for every $Z\in X$, the set $Z\cap Y$ has one element. 
\item For every surjection $f\cc A\to B$ there exists a~function $g\cc B\to A$ such that the composition $f(g)$ is $\mathrm{id}_B$.
\end{enumerate}
\end{axio}

\begin{exer}\label{ex_equiAC}
These three formulations of {\em AC} are mutually equivalent, in the sense that each can be easily derived from any other.
\end{exer}

\begin{exer}\label{ex_oneMoreAC}
In part~2 of Axiom~\ref{axio_AC}, the assumption of disjointness cannot be 
omitted.    
\end{exer}
What is AC good for? It formalizes the intuitively clear, but 
practically unrealizable, act
of selecting an element from each of {\em infinitely many} nonempty sets 
$\{A_i\cc\;i\in I\}$. If this set system is finite, we do not need AC 
for making the selection. The existence of the selector then follows 
from other set-theoretic axioms. In Appendices~\ref{app_godel} and \ref{app_cohen} we 
present proofs of the equiconsistency of AC and its negation with other set-theoretic axioms; these 
proofs are due to K.~G\"odel\index{goedel@G\"odel, Kurt} 
and P.~Cohen\index{Cohen, Paul}, respectively.

\section{Equivalence relations and partitions}

A~relation $R$ on a~set $A$ is 
\underline{reflexive\index{relation@(binary) 
relation!reflexive|emph}} if $a\,R\,a$ for every $a$  in $A$. It is 
\underline{irreflexive\index{relation@(binary) 
relation!irreflexive|emph}} if $a\,R\,a$ for no $a\in A$. It is 
\underline{symmetric\index{relation@(binary) relation!symmetric|emph}} if $a\,R\,b$ implies $b\,R\,a$ for every $a,b\in A$. It is  
\underline{transitive\index{relation@(binary) 
relation!transitive|emph}} if $a\,R\,b$ and $b\,R\,c$ imply 
$a\,R\,c$ for every $a,b,c\in A$.

\begin{defi}\label{def_relEkv}
An~\underline{equivalence relation}\index{relation@(binary) relation!equivalence|emph} on a~set 
is a~reflexive, symmetric, and transitive relation on the set.
\end{defi}

\begin{defi}\label{def_rozkl}
A~\underline{partition}\index{partition|emph} of a~set $B$ is any 
set $A$ with nonempty and mutually disjoint elements such that $\bigcup A=B$. 
The elements of $A$ are the \underline{blocks\index{partition!block of|emph}} of the partition.
\end{defi}

Let $R$ be an equivalence relation on a~set $A$. For any $a\in A$ we define 
$$
[a]_R:=\{b\in A\cc\;b\,R\,a\}\ \ (\sus A)\,.\label{block}
$$
We call the set $[a]_R$ the
\underline{$a$-block\index{relation@(binary) relation!equivalence!block of|emph}}  
in the relation $R$. It is clear that if $a\,R\,b$ then $[a]_R=[b]_R$. We set
$$
A/R:=\{[a]_R\cc\;a\in A\}\,.
$$

\begin{exer}\label{ex_bloky}
If $R$ is an equivalence relation on a~set $A$, then $A/R$ is a~partition 
of $A$. Elements $b,c\in A$ are in one block of $A/R$ $\iff$ $b\,R\,c$. 
\end{exer}

Let $X$ be a~partition of $Y$.  We define a~relation $Y/X$ on $Y$ by
$$
x\,(Y/X)\,y\iff (\exists\,Z\in X\cc x\in Z\wedge y\in Z)\,.
$$

\begin{exer}\label{ex_rozklad}
If $X$ is a~partition of $Y$, then $Y/X$ is an equivalence relation on $Y$. Elements $x,y\in Y$ are in one block of $X$ $\iff$ $x\,(Y/X)\,y$.
\end{exer}

\begin{exer}\label{ex_dualita}
Let $R$ be an equivalence relation on $A$ and $B$ be a~partition of $A$. Then
$$
A/(A/R)=R\,\text{ and }\,A/(A/B)=B\;.
$$
\end{exer}

\section{Linear orders, suprema, infima}

A~relation $R$ on a~set $A$ is \underline{trichotomic\index{relation@(binary) 
relation!trichotomic|emph}} if for every $a,b\in A$ it is true that $a\,R\,b$ or $b\,R\,a$ or $a=b$. 
$R$ is \underline{dichotomic\index{relation@(binary) 
relation!dichotomic|emph}} if for every $a,b\in A$ it is true that 
$a\,R\,b$ or $b\,R\,a$. $R$ is \underline{asymmetric\index{relation@(binary) 
relation!asymmetric|emph}} if never $a\,R\,b$ and $b\,R\,a$. Finally,
$R$ is \underline{weakly asymmetric\index{relation@(binary) 
relation!weakly asymmetric|emph}} if $a\,R\,b$ and $b\,R\,a$
imply $a=b$.

\begin{defi}\label{def_LO}
A~relation 
on a~nonempty set is a~\underline{linear order\index{linear order|emph}} 
if it is irreflexive, transitive, and trichotomic. 
\end{defi}

\begin{exer}\label{ex_onLO}
Every linear order is asymmetric.    
\end{exer}

\noindent
We denote a~linear order $<$ on a~set $A$ by a~pair $\langle A,<\rangle$.\label{linearOrd} 
Let $a,b\in A$. Notation $a\le b$ means that $a<b$ or $a=b$. Notation $a>b$ is synonymous with
$b<a$, and $a\ge b$  is synonymous with
$b\le a$. The linear orders $<$ and $>$ are also called \underline{strict\index{linear 
order!strict|emph} linear orders}. Relations $\le$ and $\ge$ are \underline{non-strict\index{linear order!non-strict|emph} linear orders}.

\begin{exer}\label{ex_tranzNeos}
Any non-strict linear order is
reflexive, transitive, dichotomic, and weakly asymmetric.
\end{exer}

\begin{exer}\label{ex_jenJedna}
Let $\langle A,<\rangle$ be a~linear order. Then for every pair of elements
$a,b\in A$ exactly one of $a<b$, $b<a$, and $a=b$ holds. 
\end{exer}

Linear orders $\langle A,<\rangle$ and $\langle B,\prec\rangle$ are 
\underline{isomorphic\index{linear order!isomorphism|emph}} if there 
exists a~bijection $f\cc A\to B$ such that for any $a,a'\in A$ we have $a<a'$ $\iff$ $f(a)\prec f(a')$.

\begin{exer}\label{ex_dalsiUl}
In the previous definition, $\iff$ can be weakened to $\Rightarrow$.     
\end{exer}

Let $\langle A,<\rangle$ 
be a~linear order and $B\sus A$. Then $B$ is \underline{bounded from 
above\index{bounded sets in linear orders!from above|emph}} if for some $h\in 
A$ we have $h\ge b$ for every $b\in B$, and $h$ is an 
\underline{upper bound\index{bounded sets in linear orders!upper bound of|emph}} of $B$. We denote the set of upper bounds of $B$ by
$H(B)$\index{bounded sets in linear orders!hb@$H(B)$|emph}\label{upperB} ($\sus A$). 
We similarly define  \underline{boundedness from below\index{bounded sets in linear orders!from 
below|emph}}, \underline{lower bounds\index{bounded sets in linear orders!lower 
bound of|emph}}, and the set $D(B)$\index{bounded sets in linear orders!db@$D(B)$}\label{lowerB} ($\sus A$) 
of lower bounds of $B$. An element $m\in B$ is the \underline{maximum\index{maximum 
of a~set in a~linear order|emph}} of $B$ if $m\in H(B)$. We similarly 
define the \underline{minimum\index{minimum of 
a~set in a~linear order|emph}} of $B$. We denote these elements by $\max(B)$ and
$\min(B)$.\label{minmax} If $B=\emptyset$, they are not defined.

\begin{exer}\label{ex_minimax}
Show that maxima and minima are unique.
\end{exer}

\begin{exer}\label{ex_minimax2}
Any nonempty finite subset in any linear order has both maximum and minimum.
\end{exer}

We define suprema and infima of sets in linear orders.

\begin{defi}\label{def_supInf}
Let $\langle A,\,<\rangle$ be a~linear order and $B\sus A$. The elements 
$$
\sup(B):=\min(H(B))\ (\in A)\,\text{ and }\,\inf(B):=\max(D(B))\ (\in A)
$$
are called the \underline{supremum\index{supremum 
in a~linear order|emph}}\label{sup} of $B$ and 
the \underline{infimum\index{infimum in 
a~linear order|emph}}\label{inf} of $B$, respectively.
\end{defi}
Suprema and infima need not exist. 
In contrast to maxima and minima, they may lie
outside the considered set. 

\begin{exer}\label{ex_infSup}Suprema and infima are unique.
\end{exer}

\begin{exer}\label{ex_dokAprVl}
Prove the following proposition. State and prove the analogous result for infima. 
\end{exer}

\begin{prop}\label{prop_aprVl}
Let $\langle A,\,<\rangle$ be a~linear order, $B\sus A$, and $c\in A$. Then 
$c=\sup(B)$ 
$\iff$ $b\le c$ for every $b\in B$ and for every $a\in A$ with $a<c$ there exists $b\in B$ such that $a<b$.
\end{prop}

\section{Well orderings and partial orders}

Well orderings form an important and useful family of linear orders. 

\begin{defi}\label{def_wellOrd}
A~linear order $\langle A,<\rangle$ is a~\underline{well ordering\index{well 
orderings|emph}} if every nonempty subset of $A$ has a~minimum element.
\end{defi}

\begin{exer}\label{ex_naWO}
The standard linear order $\langle\omega,<\rangle$ of natural numbers, introduced formally in the next chapter, is a~well 
ordering.    
\end{exer}

\begin{exer}\label{ex_naWO2}
The standard linear order $\langle \Z,<\rangle$ of integers, introduced formally in the chapter after the next chapter, is not a~well 
ordering.    
\end{exer}

\begin{exer}\label{ex_naWO3}
A~linear order $\langle A,<\rangle$ is a~well ordering $\iff$ there do not exist elements 
$a_n\in A$, $n\in\N$, such that $a_1>a_2>\cdots$.
\end{exer}

Partial orders, often called just orders, generalize linear orders.

\begin{defi}\label{def_orders}
A~relation 
on a~set is an \underline{order\index{order|emph}} if it is irreflexive and 
transitive.    
\end{defi}
In an order $\langle A,<\rangle$, we may have (distinct) 
\underline{incomparable elements\index{order!incomparable 
elements in|emph}} $a,b\in A$ such that neither $a<b$ nor $b<a$ 
holds. Like for linear orders, we associate with $\langle A,
<\rangle$ the \underline{non-strict order\index{order!non-strict|emph}} $\langle 
A,\le\rangle$ by setting $a\le b$ iff $a<b$ or $a=b$. If $B\sus A$ and $b\in B$, we say that $b$ is a~\underline{maximal element\index{order!maximal element|emph}} in $B$, if $b<c\in B$ for no $c$.
\underline{Minimal elements\index{order!minimal element|emph}} are defined similarly.

\begin{exer}\label{ex_cvicOrder1}
Every order is asymmetric.
Every non-strict order is
reflexive, transitive and weakly asymmetric.    
\end{exer}

\begin{exer}\label{ex_cvicOrder2}
For every set $x$ and subset $y\sus\mathcal{P}(x)$ of the power 
set of $x$ we have the non-strict order $\langle y,\sus\rangle$.     
\end{exer}
\chapter[From ZFC formulas to natural numbers\checkmark]{From ZFC formulas to natural numbers}\label{sec_ZFCtoN0}
\index{natural numbers|(}

As numeric domains are concerned, elementary real 
analysis mainly works with the real numbers. However, it is not 
possible to introduce them rigorously without rigorously dealing 
with the natural numbers, which is not possible without using formalized ZFC set theory. Since 
we built a sufficient fragment of ZFC in Appendix~\ref{sec_ZFC}, we 
can now rigorously construct natural numbers. The main result is 
the characterization Theorem~\ref{thm_N0is}, which asserts 
that natural numbers form a~simple ordered semiring that is unique up to isomorphism.

In this chapter, Definition~\ref{def_omega} introduces the set $\omega$ of natural numbers. 
Definition~\ref{def_plusMin} introduces
the predecessor and successor of a~natural number. In 
Theorem~\ref{thm_wellOrdNat} we prove that natural numbers form a~well 
ordering. Theorem~\ref{thm_indukce} provides the inductive definition of 
functions on natural numbers. We use it to define addition and 
multiplication of natural numbers.  
In Definition~\ref{def_Nzero} we introduce the algebraic structure $\N_0$ of natural numbers. 
Definition~\ref{def_semiR} introduces ordered 
semirings. Theorem~\ref{thm_N0is} is 
the previously mentioned theorem characterizing $\N_0$ as the unique simple ordered semiring. A~unique homomorphism from 
$\N_0$ to any semiring in Definition~\ref{def_evalM} is used in the 
proof of the characterization theorem.
Proposition~\ref{prop_subtrOm} introduces the partial operation of subtraction of
natural numbers.

\section{The linear order of natural numbers}

From Appendix~\ref{sec_ZFC} we know that a~set $x$ is inductive if $\emptyset\in x$ and for every set $y\in x$, also $y\cup\{y\}\in x$.

\begin{defi}\label{def_omega}
We define the set $\omega$ of \underline{natural numbers\index{natural numbers!set|emph}} by means of Axioms~\ref{axi_separ} 
and \ref{axi_infi} as
$$
\omega:=\{x\in x_0\cc\;\varphi(x)\}\,,
$$
where $x_0$ is an inductive set and $\varphi(x)$ is the quasi-formula\index{quasi-formula} 
$$
\forall y(\text{$y$ is inductive}\rightarrow x\in y)\,.
$$
\end{defi}
The existence of $x_0$ is ensured by Axiom~\ref{axi_infi} of infinity, and 
$\omega$ is defined by Axiom~\ref{axi_separ} of separation.  We set $\N:=\omega\setminus\{0\}$. 

\begin{exer}\label{ex_012}
Show that the sets $0:=\emptyset$, 
$1:=\{\emptyset\}$, and
$2:=\{\emptyset,\{\emptyset\}\}$ are natural numbers.    
\end{exer}

\begin{prop}\label{prop_propOm}
The following holds.
\begin{enumerate}
\item The set $\omega$ is unique because it does not depend on the choice of $x_0$.
\item The set $\omega$ is inductive and is a~subset of every inductive set.
\item If $X$ is an inductive set with $X\sus\omega$, then $X=\omega$. 
\end{enumerate}
\end{prop}
\duk
1. Let $x_0$ and $x_1$ be inductive sets, $\varphi(x)$ and $\omega$ be as in Definition~\ref{def_omega}, and let 
$$
\omega':=\{x\in x_1\cc\;\varphi(x)\}\,.
$$
We claim that $\omega'=\omega$. If $x\in\omega$, then $x\in x_0$ and $\varphi(x)$. The 
latter implies that $x\in x_1$ because $x_1$ is one of the inductive sets $y$  
mentioned in $\varphi(x)$. Thus $x\in\omega'$. The same argument shows that if 
$x\in\omega'$, then $x\in\omega$. Hence $\omega'=\omega$ by the axiom of 
extensionality.

2. This follows from Definition~\ref{def_omega}. 

3. This is immediate from part~2.
\kduk

\noindent
We refer to applications of part~3 by the phrase ``by induction''. Next, we derive 
some properties of natural numbers.

\begin{prop}\label{prop_elOfOm}
Let $m,n\in\omega$. The following holds.
\begin{enumerate}
\item $0\in m$ or $m=0$.
\item If $m\ne0$, then $m=l\cup\{l\}$ for a~unique $l\in\omega$.
\item $m\cup\{m\}\in\omega$.
\item $m\sus\omega$.
\item If $m\in n$, then $m\sus n$.
\item If $m\in n$, then $m\cup\{m\}
\in n$ or $m\cup\{m\}=n$.
\end{enumerate}    
\end{prop}
\duk
1. Let $X$ be the set of $m\in\omega$ with the stated property. Obviously, 
$\emptyset\in X$. Let $m\in X$. If $\emptyset\in m$, then also 
$\emptyset\in m\cup\{m\}$. If $m=\emptyset$, then  also 
$\emptyset\in m\cup\{m\}
=\{\emptyset\}$. Thus $m\cup\{m\}\in X$ and $X=\omega$ by induction.

2. We first prove the uniqueness. If $x\cup\{x\}=y\cup\{y\}$ for sets 
$x\ne y$, then the axiom of extensionality implies that 
$x\in y\in x$, which contradicts the axiom of foundation. Let $X$ be the 
set of $m\in\omega$ such that $m=\emptyset$ or $m$ has the stated 
representation. Then $\emptyset\in X$ by definition. Let $m\in X$. Then 
$m\in\omega$ and $m\cup\{m\}$ has  trivially the stated representation. 
Thus $m\cup\{m\}\in X$ and $X=\omega$ by by induction.

3. This follows from part~2 of 
Proposition~\ref{prop_propOm}.

4. Let $X$ be the set of $m\in\omega$ such that $m\sus\omega$. Then 
$\emptyset\in X$ because $\emptyset$ is a~subset of every set. Let $m\in X$. Then $m\in\omega$, $m\sus\omega$, and we deduce that also 
$m\cup\{m\}\sus\omega$. Thus $m\cup
\{m\}\in X$ and $X=\omega$ by by induction.

5. Let $X$ be the set of $n\in\omega$ such that the implication $m\in n$ 
$\Rightarrow$ $m\sus n$ holds. Then 
$\emptyset\in X$ because $\emptyset$ has no elements. Let $n\in X$ and $m\in n\cup\{n\}$. If $m\in n$ then $m\sus 
n\cup\{n\}$ because $m\sus n$ by the assumption that $n\in X$. If $m=n$ 
then again $m=n\sus n\cup\{n\}$. Thus  $n\cup\{n\}\in X$ and $X=\omega$ by 
by induction.

6. Let $X$ be the set of $n\in\omega$ such that the implication $m\in n$ 
$\Rightarrow$ $m\cup\{m\}\in n$ or $m\cup\{m\}=n$ holds. Then 
$\emptyset\in X$ because $\emptyset$ has no elements. Let $n\in X$ and let $m\in n\cup\{n\}$. If $m\in n$ then $m\cup
\{m\}$ is an element of $n$ or equals $n$ by the assumption that $n\in X$. 
In either case $m\cup
\{m\}\in n\cup\{n\}$. If $m=n$ then   $m\cup\{m\}=n\cup\{n\}$. Thus $n\cup
\{n\}\in X$ and $X=\omega$ by 
by induction.
\kduk

\noindent
We use parts~2 and~3 in the following definition.

\begin{defi}\label{def_plusMin}
Let $m\in\omega$. We define the \underline{successor\index{natural numbers!successor|emph}} $m\oplus1$ as $m\cup\{m\}$. If 
$m>0$, we define the \underline{predecessor\index{natural 
numbers!predecessor|emph}} $m-1$ as the $k$ given by $m=k\cup\{k\}$.    
\end{defi}

\begin{exer}\label{ex_minAndplus}
Let $m\in\omega$. Show that $(m\oplus1)-1=n$ and, if $m>0$, also $(m-1)\oplus1=m$.  
\end{exer}

We define a~\underline{linear order\index{natural numbers!linear order|emph}} $<$ on $\omega$ by setting $m<n$ $\iff$ $m\in n$. 

\begin{prop}
$\langle\omega,<\rangle$ is a~linear order.    
\end{prop}
\duk
For the irreflexivity of $<$ see Exercise~\ref{ex_irrefNat}. We show that 
$<$ is transitive. Let $k,l,m\in\omega$ with $k\in l\in m$. By part~5 of 
Proposition~\ref{prop_elOfOm}, we 
have $l\sus m$, so that $k\in m$. We show that $<$ is 
trichotomic. Let $n\in\omega$ and let $X$ be the set of 
$m\in\omega$ such that $n\in m$ or $n=m$ or $m\in n$. Then 
$\emptyset\in X$ by part~1  of Proposition~\ref{prop_elOfOm}. Let 
$m\in X$ and consider the set 
$p=m\cup\{m\}$. If $n\in m$ or $n=m$ then $n\in p$. If $m\in n$ then $p\in n$ or $p=n$ by part~6 of 
Proposition~\ref{prop_elOfOm}.  Thus $p\in X$ and $X=\omega$ by induction.
\kduk    

\begin{exer}\label{ex_irrefNat}
The relation $<$ on $\omega$ is irreflexive.    
\end{exer}

In order to show that $\langle\omega,<\rangle$ is a~well ordering, we use 
the next proposition.

\begin{prop}\label{prop_WOm}
For every $m\in\omega$, the linear order $\langle m,<\rangle$ is 
a~well ordering.    
\end{prop}
\duk
Let $X$ be the set of $m\in\omega$ with the stated 
property. Then $\emptyset\in X$ because $\emptyset$ has no nonempty 
subset. Let $m\in X$ and $p=m\cup\{m\}$. Let $\emptyset\ne A\sus p$. 
If $A\cap m\ne\emptyset$, then 
$\min(A)=\min(A\cap m)$ by the assumption that $m\in X$ and because 
$m$ is the maximum of $p$. If $A\cap m=\emptyset$, then $A=\{m\}$ and $m$ 
is the minimum of $A$. Thus $p\in X$ and $X=\omega$ by induction.
\kduk

The following theorem is basic. 

\begin{thm}\label{thm_wellOrdNat}
The\index{theorem!well ordering $\omega$|emph} 
linear order $\langle\omega,<\rangle$ is a~well ordering.    
\end{thm}
\duk
Let $\emptyset\ne A\sus\omega$. We show that $A$ has a~minimum element. We take any element $m\in A$. If $m\cap 
A=\emptyset$, then $\min(A)=m$. If $m\cap A\ne\emptyset$, then 
$\min(A)=\min(m\cap A)$ by Proposition~\ref{prop_WOm}.
\kduk

So what object or thing is the {\em infinite} set $\omega$ of natural numbers?
As we explain in Appendix~\ref{expl_infiSets}, one can
view $\omega$ in two ways. Either as {\em only} this ZFC formula:
$$
bla
$$
or as, also, the story that the formula describes a~unique infinite set $\omega$ in the 
intuitive universe of sets such that $\omega$ is fully characterized by three properties: (i) 
$\emptyset\in\omega$, (ii) if
$x\in\omega$, then $x\cup\{x\}\in\omega$, and (iii) $\omega$ has no other elements 
than those generated by (i) and (ii).

\section{Addition and multiplication}

In order to define addition and multiplication of natural numbers, we 
justify inductive definitions of functions $f\cc\omega\to\omega$. 
We set $\omega^2:=\omega\times\omega$.

\begin{thm}\label{thm_indukce}
Let $m_0\in\omega$ and let $F\cc\omega^2\to\omega$. 
There exists a~unique function
$f\cc\omega\to\omega$ such that
$f(0)=m_0$ and for every $n\in\omega$ with $n>0$, 
$$
f(n)=F(n,\,f(n-1))\,.
$$
\end{thm}
\duk
We first prove the uniqueness of $f$ and then the existence. Suppose that 
$f\ne g$ are two functions from $\omega$ to $\omega$ with the 
stated property. By Theorem~\ref{thm_wellOrdNat}, there is   the least number $n\in\omega$ such that $f(n)\ne g(n)$. Clearly, $n\ne0$ because $f(0)=m_0=g(0)$, and 
$f(n-1)=g(n-1)$. We get the contradiction
$$
f(n)=F(n,\,f(n-1))=F(n,\,g(n-1))
=g(n)\,.
$$

In order to prove that $f$ exists, we consider the set $X$ of 
$m\in\omega$ such that there exists 
a~function $f\cc m\to\omega$ with the property that $f(0)=m_0$ and $f(n)=F(n,f(n-1))$
for every $n\in m$ with $n>0$. Then $0\in X$ because the empty function 
to $\omega$ 
has, trivially, the stated property. Let $m\in X$. If $m=0$, then $m\oplus1=1\in
X$ because the function 
$f\cc 1\to\omega$ with $f(0)=m_0$ has the stated property. Let $m\in\N$ and 
$f\cc m\to\omega$ be a~function with the stated property. We define $g\cc 
m\oplus1\to\omega$ as an extension of $f$ by the value 
$$
g(m)=F(m,\,f(m-1))\,. 
$$
It 
follows that $g$ has the stated property. Thus, $m\oplus1\in X$. By 
 induction, $X=\omega$. For $m\in\omega$, we denote by $X_m$ 
the set of maps from $m$ to $\omega$ with the stated property. We have 
just proven that $X_m\ne\emptyset$ for every $m\in\omega$. 

Clearly, (i) if $m\in n$ and $f\in X_n$, then 
$f\,|\,m\in X_m$, and (ii) every $X_m$ is a~one-element set; (ii) follows 
by the initial uniqueness argument.
Using AC,\index{axiom!of choice, AC} we choose for every $m\in\omega$ the function $f_m\in X_m=\{f_m\}$. It 
follows from (i) and (ii) that the desired function has graph
$${\textstyle
G_f:=\bigcup_{m\in\omega}G_{f_m}
}\,.
$$
\kduk

\begin{exer}\label{ex_onIndDef}
If $m_0=0$, $F(n,m)=m$ for odd $n$, and 
$F(n,m)=m+1$ for even $n$, what is $f(n)$?
\end{exer}

Using the previous inductive definition, we define \underline{addition\index{natural numbers!addition|emph}} and \underline{multiplication\index{natural numbers!multiplication|emph}} of 
natural numbers. Addition is iterated  successor $\oplus1$.

\begin{cor}\label{cor_addi}
There exists a~unique operation $+\cc\omega^2\to\omega$ such that for every $m\in\omega$ we have $m+0=m$ and 
$m+n=(m+(n-1))\oplus1$ for $n>0$.
\end{cor}
\duk
We fix $m$ and use Theorem~\ref{thm_indukce} with $m_0=m$ and
$F(k,l)=l\oplus1$.
\kduk

\noindent
Multiplication is iterated addition.

\begin{cor}\label{cor_mult}
There exists a~unique operation $\cdot\cc\omega^2\to\omega$ such that for every $m\in\omega$ we have $m\cdot0=0$ and 
$m\cdot n=(m\cdot(n-1))+m$ for $n>0$.
\end{cor}
\duk
We fix $m$ and use Theorem~\ref{thm_indukce} with $m_0=m$ and
$F(k,l)=l+m$.
\kduk

\section{The characterization theorem}

We begin by defining the algebraic structure of natural numbers. 

\begin{defi}\label{def_Nzero}
Recall the set  
$$
\omega=\{0,\,1,\,2,\,\ds\}=
\{\emptyset,\,
\{\emptyset\},\,\{\emptyset,\,
\{\emptyset\}\},\,\ds\}
$$ 
in Definition~\ref{def_omega}. The
algebraic structure of \underline{natural numbers\index{natural numbers!structure|emph}}  
$$
\N_0:=\langle
\omega,\,0,\,1,\,+,\,\cdot,\,<
\rangle\label{natZer}
$$
consists of the set $\omega$, the elements zero 
$0=\emptyset\in\omega$ and one 
$1=\{\emptyset\}\in\omega$, the operations $+$ and $\cdot$ introduced in 
Corollaries~\ref{cor_addi} and \ref{cor_mult}, and the linear order 
$<$ on $\omega$ defined above by
$m<n$ $\iff$ $m\in n$. 
\end{defi}
For natural numbers, we distinguish in notation 
between the base set $\omega$ and the algebraic structure $\N_0$ on $\omega$. 
For integers, fractions, extended reals, real numbers, and complex numbers, we 
use the same symbol $\Z$, $\Q$, $\R^*$, $\R$, and $\C$, respectively, for algebraic structure  and the base set. 
$\N_0$ is an ordered semiring; we will introduce this algebraic structure next.

\begin{defi}\label{def_semiR}
An \underline{ordered semiring\index{semiring!ordered|emph}} is an algebraic structure
$$
S=\langle S,\,0_S,\,1_S,\,+,\,\cdot
,\,<\rangle
$$
such that $S$ is a~set, $0_S,1_S\in S$ and $0_S\ne1_S$, $+$ and $\cdot$ are operations on $S$, $<$ is a~relation on $S$, and the following holds. $0_S$ is neutral to $+$, and $1_S$ to 
$\cdot$, the operations $+$ and $\cdot$ are 
commutative and associative, $\cdot$ is distributive to $+$, the relation $<$ is 
a~linear order, $0_S<1_S$, and for every $a,b,c\in S$,
$$
a<b\Rightarrow a+c<b+c\,.
$$
We call this implication the \underline{first order 
axiom\index{semiring!ordered!first order 
axiom|emph}}. The comparison $0_S<1_S$ is the \underline{second order axiom\index{semiring!ordered!second order 
axiom|emph}}. By omitting the order $<$ we get a~\underline{semiring\index{semiring|emph}}.
\end{defi}
Both neutral elements in $S$ are unique
(Exercise~\ref{ex_uloNeuPr}). One can deduce from the first order axiom
its non-strict version. 

\begin{exer}\label{ex_nonstrOA}
For every $a,b,c\in S$, 
$$
a\le b\Rightarrow a+c\le b+c\,.
$$    
\end{exer}

\begin{exer}\label{ex_semiRinf}
Show that every ordered semiring is infinite.    
\end{exer}

Let $S$ and $S'$ be semirings. A~\underline{semiring 
homomorphism\index{semiring!homomorphism|emph}} from $S$ to $S'$ is any map $f\cc S\to S'$ that is 
a~homomorphism to addition and multiplication in $S$ and $S'$, and that 
has values $f(0_S)=0_{S'}$ and 
$f(1_S)=1_{S'}$. We know from Exercise~\ref{ex_naHomom} that if $f$ is an 
isomorphism, the last condition may be omitted. If $S$
and $S'$ are ordered and $f$ is a~bijective homomorphism that preserves both orders, we say that
$f$ is an \underline{isomorphism of ordered semirings\index{semiring!ordered!isomorphism|emph}}. A~semiring $S$ is \underline{simple\index{semiring!simple|emph}} if for any set $X\sus S$ the following holds. If 
$0_S\in X$ and if for any $x\in X$ also $x+1_S\in X$, then $X=S$. We prove a~property of semirings.

\begin{lemma}\label{lem_prodZeros}
Let $S$ be a~semiring. Then 
$$
0_S\cdot0_S=0_S\,.
$$
\end{lemma}
\duk
We write $0=0_S$ and $1=1_S$ and compute that
$$
0=0\cdot1=0\cdot(0+1)=0\cdot0+0\cdot1=
0\cdot0+0=0\cdot0\,.
$$
\kduk

Related results are left as an exercise.

\begin{exer}\label{ex_nextL}
Let $S$ be a~semiring. Prove that $0_S\cdot a=a\cdot0_S=0_S$ for every 
$a\in S$ if $S$ is ordered or simple.    
\end{exer}
On the other hand, it is easy to define a~large family of semirings $S$ such 
that $0_S\cdot a=a\cdot0_S\ne0_S$ for some $a\in S\setminus\{0_S\}$. In fact, we can extend any semiring 
$$
S=\langle S,\,0_S,\,1_S,\,+,\,\cdot
\rangle
$$
to such an example. We set $S_{\infty}:=S\cup\{\infty\}$, where $\infty$ is a~new element, and define
for every $a\in S_{\infty}$ that
$$
a+\infty=\infty+a=a\cdot\infty=\infty
\cdot a:=\infty\,.
$$
It is easy to check in the next exercise that $S_{\infty}$ is 
a~semiring. In $S_{\infty}$ we have  $0_S\cdot\infty=\infty\cdot0_S=\infty\ne
0_S$.

\begin{exer}\label{ex_example}
Show that $S_{\infty}$ is a~semiring.     
\end{exer}

We can state the characterization theorem.

\begin{thm}\label{thm_N0is}
The\index{theorem!characterization of natural numbers|emph} 
structure $\N_0$ is an ordered simple semiring. Every two ordered 
simple semirings are isomorphic as ordered semirings.    
\end{thm}
The latter claim is not valid without the adjective ``ordered'' 
because there exist both finite and infinite simple semirings (Exercise~\ref{ex_finSiSemi}).
We prove the former claim in the theorem in Proposition~\ref{prop_NatSemir}, and 
the latter claim in Proposition~\ref{prop_isomorpSR}. Both 
proofs require several lemmas. 

\begin{exer}\label{ex_finSiSemi}
Describe some finite simple semirings.  
\end{exer}

\section{$\N_0$ is an ordered simple semiring}

We need the next four lemmas to prove it in Proposition~\ref{prop_NatSemir}.

\begin{lemma}\label{lem_distrM1}
Let $m,n\in\omega$. Then the following holds.
\begin{enumerate}
\item $0+m=m$.
\item If $m,n>0$, then $(m+n)-1=(m-1)+n=m+(n-1)$.
\end{enumerate}
\end{lemma}
\duk
1. For $m=0$ this follows from the definition of $+$. Let $m>0$. Then, 
by the definition of $+$ and induction on $m$, 
$$
0+m=(0+(m-1))+1=(m-1)+1=m\,.
$$

2. We proceed by induction on $n$. Let $n=1$. Then, by the definition of $+$,  
$$
(m-1)+n=(m-1)+1=m=m+0=m+(1-1)=m+(n-1)\,.
$$
By Exercise~\ref{ex_minAndplus}, $(m+n)-1=(m+1)-1=m$.
Let $n>1$. Then, by the definition of $+$ and induction,
\begin{eqnarray*}
(m-1)+n&=&((m-1)+(n-1))+1\\
&=&(m+((n-1)-1))+1=m+(n-1)\,.
\end{eqnarray*}
By the definition of $+$ and Exercise~\ref{ex_minAndplus}, $(m+n)-1=((m+(n-1))+1)-1=m+(n-1)$.
\kduk

\begin{lemma}\label{lem_dalsiLe}
Let $l,m,n\in\omega$. Then the following holds. \begin{enumerate}
\item $(l+m)+n=(l+n)+m$.
\item If $m,n>0$, then $(l+(m-1))+n=(l+(n-1))+m$.
\end{enumerate}   
\end{lemma}
\duk
1. For fixed $l\in\omega$, we proceed by induction on $m$ and $n$. 
If one of them is $0$, the stated identity holds by the definition of 
$+$. Let $m,n>0$. Then
\begin{eqnarray*}
(l+m)+n&=&((l+(m-1))+1)+n\\
&=&((l+1)+(m-1))+n=(((l+1)+(m-1))+(n-1))+1\\
&=&(((l+1)+(n-1))+(m-1))+1\\
&=&((l+1)+(n-1))+m=((l+(n-1))+1)+m\\
&=&(l+n)+m\,.
\end{eqnarray*}
In the first equality, we use the definition of $+$. In the second 
equality, we use induction on $m$. The third equality follows from the 
definition of $+$. The fourth equality follows from induction on 
$m$ (or $n)$. In the fifth equality, we use the definition of $+$. In the sixth 
equality, we use induction on $n$. The last equality follows from the definition of $+$. 

2. By the definition of $+$ and part~1, 
\begin{eqnarray*}
(l+(m-1))+n&=&((l+(m-1))+(n-1))+1\\
&=&((l+(n-1))+(m-1))+1=((l+(n-1))+m\,.
\end{eqnarray*}
\kduk

\begin{lemma}\label{lem_onMul}
Let $m,n\in\omega$. Then the following holds.
\begin{enumerate}
\item $0\cdot m=0$.
\item If $m,n>0$, then $(n-1)\cdot m+m=n\cdot (m-1)+n$.
\end{enumerate}
\end{lemma}
\duk
1. For $m=0$ this follows from the definition of $\cdot$. Let $m>0$. By 
the definitions of $\cdot$ and $+$, and by induction on $m$,
$$
0\cdot m=(0\cdot(m-1))+0=0+0=0\,.
$$

2. We proceed by induction on $m$. Let $m=1$. Then 
\begin{eqnarray*}
(n-1)\cdot m+m&=&((n-1)\cdot(1-1)+(n-1))+1\\
&=&((n-1)\cdot0+(n-1))+1=
(0+(n-1))+1\\
&=&(n-1)+1=n=n\cdot0+n=n\cdot(m-1)+n\,.
\end{eqnarray*}
In the first equality, we use the definition of $\cdot$. The second 
equality uses that $1-1=0$. The third 
equality uses the definition of $\cdot$. The fourth equality uses 
part~1 of Lemma~\ref{lem_distrM1}.
The fifth equality follows from Exercise~\ref{ex_minAndplus}. The sixth equality uses the definition of 
$\cdot$ and part~1 of Lemma~\ref{lem_distrM1}. The last equality uses that $1-1=0$.

Let $m>1$. Then
\begin{eqnarray*}
(n-1)\cdot m+m&=&((n-1)\cdot(m-1)+(n-1))+m\\
&=&((n-1)\cdot(m-1)+(m-1))+n\\
&=&(n\cdot((m-1)-1)+n)+n\\
&=&n\cdot(m-1)+n\,.
\end{eqnarray*}
The first equality follows from the definition of $\cdot$. The second 
equality follows from part~2 of Lemma~\ref{lem_dalsiLe}. In the third 
equality, we use induction on $m$. In the last, fourth equality, we 
use the definition of $\cdot$.  
\kduk

We leave the proof of the last lemma as an exercise.

\begin{exer}\label{ex_lastLem}
Prove the next lemma.    
\end{exer}

\begin{lemma}\label{lem_posNer}
Let $m,n\in\omega$ and $m<n$. Then $m+1<n+1$.
\end{lemma}

Now we can prove that $\N_0$ is an ordered simple semiring.

\begin{prop}\label{prop_NatSemir}
The algebraic structure 
$$
\N_0=\langle
\omega,\,0,\,1,\,+,\,
\cdot,\,<\rangle
$$ 
introduced in Definition~\ref{def_Nzero} is 
an ordered simple  semiring.    
\end{prop}
\duk
The commutativity of $+$ follows at once if we set $l=0$ in part~1 of 
Lemma~\ref{lem_dalsiLe}. We prove that $\cdot$ is commutative. Let 
$m,n\in\omega$. If one of them is $0$, then $m\cdot n=n\cdot m=0$
by the definition of $\cdot$ and part~1 of Lemma~\ref{lem_onMul}. Let $m,n>0$. 
By the definition of $\cdot$, induction on $n$, and by part~2 of 
Lemma~\ref{lem_onMul},
\begin{eqnarray*}
m\cdot n&=&(m\cdot(n-1))+m=
((n-1)\cdot m)+m\\
&=&(n\cdot(m-1))+n=n\cdot m\,.
\end{eqnarray*}

We prove the neutrality of $0$
and $1$. Let $m\in\omega$. We have that $m+0=m$ and 
$$
m\cdot1=m\cdot(1-1)+m=m\cdot0+m=0+m=m+0=m\,,
$$
where these six equalities are easy to justify. Hence, $0$ and $1$ are neutral to $+$ and $\cdot$, respectively.

We prove that $+$ is associative. Let $k,l,m\in\omega$. If one of them is 
zero, then $(k+l)+m=k+(l+m)$ holds due to the neutrality of $0$. Let $k,l,m>0$. Then
\begin{eqnarray*}
(k+l)+m&=&((k+l)+(m-1))+1\\
&=&(k+(l+(m-1)))+1\\
&=&(k+((l+m)-1))+1=k+(l+m)\,. 
\end{eqnarray*}
The first equality follows from the definition of $+$. In the second 
equality, we use induction on $m$. The third equality follows from part~2 of 
Lemma~\ref{lem_distrM1}. The last equality follows from the definition of $+$.

We prove that $\cdot$ is distributive to $+$. Let $k,l,m\in\omega$. If one 
of $l$ and $m$ is $0$, then $k\cdot(l+m)=k\cdot l+k\cdot m$ holds 
by the neutrality of $0$, and by $0n=n0=0$. Let $l,m>0$. Then
\begin{eqnarray*}
k\cdot(l+m)&=&k\cdot((l+m)-1)+k\\
&=&k\cdot(l+(m-1))+k\\
&=&(k\cdot l+k\cdot(m-1))+k\\
&=&k\cdot l+(k\cdot(m-1)+k)=k\cdot l+k\cdot m\,.
\end{eqnarray*}
The first equality follows from the 
definition of $\cdot$. The second equality follows from part~2 of 
Lemma~\ref{lem_distrM1}. In the third equality, we use induction on $m$. 
The fourth equality follows from the associativity of $+$. The last fifth equality follows from the definition of $\cdot$.

We prove that $\cdot$ is associative. Let $k,l,m\in\omega$. If one of them 
is zero, then $(k\cdot l)\cdot m=k\cdot(l\cdot m)=0$ because 
$0n=n0=0$. Let $k,l,m>0$. Then
\begin{eqnarray*}
(k\cdot l)\cdot m&=&(k\cdot l)\cdot(m-1)+k\cdot l\\
&=&k\cdot(l\cdot(m-1))+k\cdot l\\
&=&k\cdot(l\cdot(m-1)+l)=k\cdot(l\cdot m)\,.
\end{eqnarray*}
The first
equality follows from the definition of $\cdot$. In the second 
equality, we use induction on $m$.
The third
equality follows from the distributivity of $\cdot$. The 
last equality follows from the definition of $\cdot$.

We prove that the two order axioms hold in $\N_0$. Let 
$k,l,m\in\omega$ with $k<l$. If $m=0$, then $k+m<l+m$ by the 
neutrality of $0$. Let $m>0$. Then, by the definition of $+$, induction on $m$, and by Lemma~\ref{lem_posNer},
$$
k+m=(k+(m-1))+1<
(l+(m-1))+1=l+m\,,
$$
which proves the first order axiom. The second order axiom $0<1$ is 
obvious.

Finally, let $0\in X\sus\omega$ and let $X$ be closed to adding $1$.
Since $m+1=m\cup\{m\}$, the set $X$ is inductive and $X=\omega$
by induction. The semiring $\N_0$ is ordered and simple.
\kduk

\section{$\N_0$ is unique}

In order to prove the isomorphism claim of Theorem~\ref{thm_N0is} in Proposition~\ref{prop_isomorpSR}, 
we introduce, for any semiring $S$, a~semiring homomorphism from $\N_0$ to 
$S$.

\begin{defi}\label{def_evalM}
Let $S=\langle
S,0_S,1_S,\oplus,\odot
\rangle$ be a~semiring. We define by induction a~map 
$f_S\cc\omega\to S$. First, $f_S(0):=0_S$. For $m\in\omega$ with $m>0$ we define
$$
f_S(m):=f_S(m-1)\oplus1_S\,.
$$
\end{defi}
A~proof of existence and uniqueness of $f_S$ would be similar to the proof of Theorem~\ref{thm_indukce}.
Clearly, $f_S(0)=0_S$ and $f_S(1)=1_S$.

We prove 
Proposition~\ref{prop_isomorpSR} by means of the next four lemmas. We assume that $S=\langle
S,0_S,1_S,\oplus,\odot\rangle$ 
is a~semiring and that $f_S\cc\omega\to S$ is as in Definition~\ref{def_evalM}. The lemmas show that $f_S$ is 
a~semiring homomorphism, and that if $S$ is ordered and simple, then $f_S$ is 
an isomorphism of ordered semirings. 

\begin{lemma}\label{lem_lema1}
For every $m,n\in\omega$, we have $f_S(m+n)=f_S(m)\oplus f_S(n)$.
\end{lemma}
\duk
By induction on $n$. For $n=0$ the equality holds for every $m\in\omega$ due 
to the neutrality of $0_S$. Let $n>0$. Then, writing $f$ for $f_S$, we have
\begin{eqnarray*}
f(m+n)&=&f((m+1)+(n-1))=f(m+1)\oplus f(n-1)\\
&=&(f(m)\oplus1_S)\oplus f(n-1)\\
&=&f(m)\oplus(f(n-1)\oplus1_S)=
f(m)\oplus f(n)\,.
\end{eqnarray*}
The first equality follows from the properties of $+$.
The second equality follows from induction. The third equality follows 
from the definition of~$f$. The fourth equality follows from the 
properties of $\oplus$. The last equality follows from the 
definition of $f$.
\kduk

\begin{lemma}\label{lem_lema2}
For every $m,n\in\omega$, we have $f_S(m\cdot n)=f_S(m)\odot f_S(n)$.
\end{lemma}
\duk
We again write $f$ for $f_S$ and first prove the equality for $n=0$ by induction on $m$. For $m=0$
we have 
$$
f(0\cdot0)=f(0)=0_S=0_S\odot0_S=
f(0)\odot f(0)\,.
$$
The third equality uses Lemma~\ref{lem_prodZeros}. If $m>0$, then
\begin{eqnarray*}
f(m\cdot0)&=&f(0)=0_S=0_S\oplus0_S=
f(m-1)\odot0_S\oplus1_S\odot0_S\\
&=&(f(m-1)\oplus1_S)\odot0_S=
f(m)\odot f(0)\,.
\end{eqnarray*}
The first equality follows from the
properties of $\cdot$. In the second 
equality, we use the definition of $f$. In the third equality, we use
the neutrality of $0_S$.
The fourth equality follows from induction and the neutrality of $1_S$ and $0_S$. The fifth equality follows from the 
distributive law in $S$. In the last equality, we use the 
definition of $f$. The equality therefore holds for every 
$m\in\omega$ and $n=0$.

For $n>0$ we proceed by induction on $n$. Let $m,n\in\omega$ and $n>0$. Then
\begin{eqnarray*}
f(m\cdot n)&=&
f(m\cdot(n-1)+m)=f(m\cdot(n-1))\oplus
f(m)\\
&=&f(m)\odot f(n-1)\oplus f(m)=
f(m)\odot(f(n-1)\oplus1_S)\\
&=&f(m)\odot f(n)\,.
\end{eqnarray*}
In the first equality, we use the properties of $\cdot$. In the second 
equality, we use Lemma~\ref{lem_lema1}.  In the third 
equality, we use induction. The fourth equality follows from the 
properties of $\odot$. The last equality follows from the definition of $f$.
\kduk

\begin{lemma}\label{lem_lema3}
Let $S$ be an ordered semiring, with the linear order $\prec$. 
Then for every $m,n\in\omega$,  
$$
f_S(m)\prec f_S(n)\iff m<n\,.
$$
\end{lemma}
\duk
We again write $f$ for $f_S$.
By the trichotomy of $<$, it suffices to prove the implication 
$$
m<n\Rightarrow f(m)\prec f(n)\,.
$$
Let $m=0$. We proceed by induction on $n$. If $n=1$ then the implication 
holds because 
$$
f(0)=0_S\prec1_S=f(1)
$$
is the second order axiom. Let $n>1$. From the inductive assumption 
$$
f(0)=0_S\prec f(n-1)
$$
we get by adding $1_S$ and using both order axioms that
$$
f(0)=0_S\prec1_S=0_S\oplus1_S\prec f(n-1)\oplus 1_S=f(n)\,,
$$
so that $f(0)\prec f(n)$ by the transitivity of $\prec$.

For $m>0$ we proceed by induction on $m$. We assume that $0<m<n$. Since $m-1<n-1$ (by the trichotomy and first order axiom for $<$), we get by adding $1_S$ to the inductive assumption
$$
f(m-1)\prec f(n-1)
$$
by using the first order axiom for $\prec$, that
$$
f(m)=f(m-1)\oplus1_S\prec f(n-1)\oplus1_S=f(n)\,.
$$
Thus $f(m)\prec f(n)$.
\kduk

We again leave the proof of the last lemma as an exercise.

\begin{exer}\label{ex_jesteCvic}
Prove the next lemma.    
\end{exer} 

\begin{lemma}\label{lem_lema4}
Let $S$ be a~simple ordered semiring, with the linear order denoted by $\prec$. Then the
map 
$$
f_S\cc\N_0=\langle\omega,\,0,\,1,\,+,\,\cdot,\,<\rangle\to
S=\langle S,\,0_S,\,1_S,\,\oplus,\,\odot,\,\prec\rangle
$$
is an isomorphism of ordered semirings.
\end{lemma}

Now we can prove that every two ordered simple semirings are isomorphic.

\begin{prop}\label{prop_isomorpSR}
Suppose that
$$
S=\langle
S,\,0_S,\,1_S,\,+,\,\cdot,\,<
\rangle\,\text{ and }\,
S'=\langle
S',\,0_{S'},\,1_{S'},\,\oplus,\,\odot,\,\prec
\rangle
$$
are two ordered simple semirings. Then $S$ and $S'$ 
are isomorphic, which means that there exists a~bijection $f\cc S\to S'$
 such that for every $a,b\in S$,
$$
f(a+b)=f(a)\oplus f(b),\ 
f(a\cdot b)=f(a)\odot f(b),\,\text{ and  }\,a<b \iff f(a)\prec f(b)\,.
$$
\end{prop}
\duk
By Lemma~\ref{lem_lema4} and Exercises~\ref{ex_naHomom1} and \ref{ex_naHomom2}, the map
$$
f\cc S\to S',\ f:=f_{S'}(f_S^{-1})\,,
$$
where $f_S$ and $f_{S'}$ are as in Definition~\ref{def_evalM}, is the desired 
isomorphism.
\kduk

\noindent
This concludes the proof of Theorem~\ref{thm_N0is}.

\section{Subtraction}

We introduce the partial operation of \underline{subtraction\index{natural numbers!subtraction|emph}} on $\omega$. 
We use it in the next chapter to define arithmetic on $\Z$.

\begin{prop}\label{prop_subtrOm}
There exists a~unique partial operation $-$ of subtraction on $\omega$ such that 
if $m,n\in\omega$ with $m\ge n$, then 
$n+(m-n)=m$. If 
$m<n$, the value $m-n$ is undefined.
\end{prop}
\duk
Let $n,m\in\omega$ with $n\le m$. We prove that there is 
a~unique number $k\in\omega$ such that $n+k=m$; then we set $m-
n:=k$. We first prove the uniqueness and then the existence. If $k,l\in\omega$ 
are distinct, then, for example, $l<k$. Using the first order axiom, we deduce that
$$
n+l<n+k\,,
$$
which shows that the map $\omega\ni k\mapsto n+k\in\omega$ is injective 
and subtraction is unique. Suppose that $n\le m$, but 
$k$ does not exist, and $n$ is $<$-minimum. Then $n>0$ 
because $0+m=m$ and $(n-1)+l=m$ for some $l\in\omega$. Clearly, $l>0$ 
because $n-1<m$. As we know from part~2 of 
Lemma~\ref{lem_distrM1}, 
$$
n+(l-1)=(n-1)+l=m\,.
$$
\kduk

To prove the distributive law in $\Z$, we need 
the distributivity of multiplication in  $\N_0$ with respect to 
subtraction. We leave the proof as an exercise.

\begin{exer}\label{ex_jesteDk}
Prove the next corollary.    
\end{exer}

\begin{cor}\label{cor_distrSubtr}
Let $l,m,n\in\omega$ with $m\le n$. Then
$$
l\cdot(n-m)=l\cdot n-l\cdot m\,.
$$
\end{cor}
\index{natural numbers|)}
\chapter[From natural numbers to integers\checkmark]{From natural numbers to integers}\label{sec_N0toZ}
\index{integers|(}

In this chapter, we extend the semiring $\N_0$ of natural numbers that we 
constructed in the previous chapter to the ring $\Z$ of integers. The main result is 
the characterization Theorem~\ref{thm_Zis}, 
which asserts that integers form an ordered simple ring that is unique up to 
isomorphism. In the last section, we review the alternative construction of 
integers in the form of difference integers. These are equivalence blocks of 
pairs of natural numbers.

Definitions~\ref{def_integers} and  \ref{def_linOrdZ} introduce the set $\Z$ 
of integers and the standard linear order on it. Addition and multiplication of 
integers are introduced in Definitions~\ref{def_addZ} and 
\ref{def_mulZ} in the elementary school way. The algebraic structure of integers 
$\Z$ is introduced in Definition~\ref{def_algStZ}. 
Definition~\ref{def_domain} reviews the notion of an ordered ring. Then comes the characterization Theorem~\ref{thm_Zis}. In
Definition~\ref{def_mapsfR} we introduce a~ring homomorphism from the ring $\Z$ to 
any given ring. Definitions~\ref{def_operZpr} and 
\ref{def_zscarou} introduce operations with difference integers and their 
algebraic structure $\Z'$. By Corollary~\ref{cor_ZandZpr}, which follows 
from the characterization theorem, $\Z$ and $\Z'$ are isomorphic as ordered rings. 

\section{The linear order of integers}

Recall that $\N=\omega\setminus
\{0\}=\omega\setminus\{\emptyset\}$. For $m\in\N$ we define 
$-m:=\langle0,m\rangle$. Let 
$$
-\N:=\{-m\cc\;m\in\N\}\,.
$$

\begin{defi}\label{def_integers}
The set $\Z$ of \underline{integers\index{integers!set|emph}} is given by the disjoint union
$$
\Z:=-\N\cup\omega\,.
$$
The elements of $-\N$ are the \underline{negative 
integers\index{integers!negative|emph}}. The elements of $\omega$ are 
called, in the context of $\Z$, the \underline{nonnegative 
integers\index{integers!nonnegative|emph}}. For $n\in\omega$, the \underline{absolute 
value\index{integers!absolute value|emph}} of $n$ is defined by $|n|:=n$. If $-n\in-\N$, we set $|-n|:=n$. 
\end{defi}
These integers literally contain the natural numbers because $\omega\sus\Z$.

\begin{exer}\label{ex_disjUni}
Check that the displayed union is indeed disjoint.     
\end{exer}

We remind the usual \underline{linear order\index{integers!linear order|emph}} on $\Z$.

\begin{defi}\label{def_linOrdZ}
Let $m,n\in\Z$. We set $m<n$ $\iff$ if $m,n\in\omega$ and $m<n$ in 
$\N_0$, or if $m\in-\N$ and $n\in\omega$, or if $m,n\in-\N$ and $|m|>|n|$ in $\N_0$. 
\end{defi}

\noindent
Clearly, $m<0<n$ for every $m\in-
\N$ and $n\in\N$.

\begin{exer}\label{ex_linOrdZ}
Show that this is a~linear order on $\Z$ and that on $\omega$ it coincides with the 
linear order in $\N_0$. Is it a~well ordering?    
\end{exer}

\section{Addition and multiplication}

We introduce \underline{addition\index{integers!addition|emph}} and \underline{multiplication\index{integers!multiplication|emph}} of integers. We begin with addition, 
which is more complicated than multiplication. We use the partial 
operation of subtraction on $\omega$ introduced in 
Proposition~\ref{prop_subtrOm}. In the next two definitions, operations on 
the right sides of $:=$ are in $\N_0$. 

\begin{defi}\label{def_addZ}
Let $m,n\in\Z$. We define the sum $m+n$ as follows.
\begin{enumerate}
\item If $m,n\ge0$, then $m+n$ is as in $\N_0$.
\item If $m<0\le n$, then $m+n:=n-|m|$ if $|m|\le n$, and
$m+n:=-(|m|-n)$ if $|m|>n$.
\item If $n<0\le m$, then $m+n:=m-|n|$ if $|n|\le m$, and
$m+n:=-(|n|-m)$ if $|n|>m$.
\item If $m,n<0$, then $m+n:=-(|m|+|n|)$.
\end{enumerate}
\end{defi}

\begin{defi}\label{def_mulZ}
Let $m,n\in\Z$. We define the product $m\cdot n$ as follows.
\begin{enumerate}
\item If $m=0$ or $n=0$, then 
$m\cdot n:=0$.
\item If $m,n<0$ or $m,n>0$, then 
$m\cdot n:=|m|\cdot|n|$.
\item If $m<0<n$ or $n<0<m$, then 
$m\cdot n:=-(|m|\cdot|n|)$.
\end{enumerate}
\end{defi}
It is easy to see that the operations $+$ and $\cdot$ are commutative, 
and that $0$ is neutral to $+$ and $1$ to $\cdot$. Also, $0+0=0$, if 
$m\in\N$ then $m+(-m)=0$, and if $-m\in-\N$ then $-m+m=0$\,---\,every integer has an additive inverse (see below the definition of rings).

\begin{exer}\label{ex_oneExZ}
Compute $3+(-5)$, $(-10)+5$, and $2\cdot(-3)$. Compare the pairs of 
integers $-2$ and $-3$, and $1$ and $-10$. \end{exer}

\section{The characterization theorem}

We define the algebraic structure $\Z$. 

\begin{defi}\label{def_algStZ}
The algebraic \underline{structure\index{integers!structure|emph}} of integers
$$
\Z:=\langle
\Z,\,0_{\Z},\,1_{\Z},\,+,\,\cdot,\,<
\rangle
$$
consists of the set $\Z$ introduced in Definition~\ref{def_integers}, the 
elements $0_{\Z}:=0$ and $1_{\Z}:=1$ in $\omega$, the operations of addition $+$ and 
multiplication $\cdot$ on $\Z$ introduced in Definitions~\ref{def_addZ} and \ref{def_mulZ}, and the linear 
order $<$ on $\Z$ introduced in Definition~\ref{def_linOrdZ}.
\end{defi}
Instead of $0_{\Z}$ and $1_{\Z}$, we usually write just $0$ and $1$. These are 
the integers that we learned in elementary school. In our heads, we add and multiply 
integers according to the rules in Definitions~\ref{def_addZ} and 
\ref{def_mulZ}. In Section~\ref{sec_diffInte}, we outline another construction of integers 
that is based on equivalence blocks of pairs of natural numbers.

We review the algebraic structure of a~ring that fits the integers.
Recall Definition~\ref{def_semiR} of ordered semirings.

\begin{defi}\label{def_domain}
An \underline{ordered ring\index{ring!ordered|emph}}
$$
R=\langle R,\,0_R,\,1_R,\,+,\,\cdot
,\,<\rangle
$$
is an ordered semiring satisfying two conditions. {\em (i) }For every 
$a\in R$ there is $b\in R$ with 
$a+b=0_R$. We call $b$ the \underline{additive inverse\index{ring!additive inverses|emph}} of $a$ and denote it by $-a$. 
{\em (ii) }The \underline{second order axiom\index{ring!ordered!second order axiom|emph}} now takes the form that for every $a,b,c\in R$ with $c>0_R$,
$$
a<b\Rightarrow a\cdot c<b\cdot c\,.
$$
By omitting the order $<$ we get a~\underline{ring\index{ring|emph}}.
\end{defi}

\begin{exer}\label{ex_AddInvUni}
Additive inverses in rings are uniquely determined.    
\end{exer}

\begin{exer}\label{ex_2ndOAx}
Show that the ring second order axiom implies the semiring 
version: $0_R<1_R$ in every ordered ring $R$.    
\end{exer}
The non-strict version of the second order axiom is easily deduced: for every $a,b,c\in R$ with $c\ge0_R$,
$$
a\le b\Rightarrow a\cdot c\le b\cdot c\,.
$$
Subtraction of 
natural numbers is a~partial operation. In rings, \underline{subtraction\index{ring!subtraction|emph}} is total: for any pair 
$a,b\in R$ we define 
$$
a-b:=a+(-b)\,. 
$$
Let $R$ and $R'$ be rings. 
A~\underline{ring homomorphism\index{ring!homomorphism|emph}} from $R$ to $R'$ is a~map $f\cc R\to R'$ that is 
a~homomorphism to the operations $+$ and $\cdot$ in $R$ and in $R'$ and that has 
value $f(1_R)=1_{R'}$. 

\begin{exer}\label{ex_onRiHomo}
Show that ring homomorphisms preserve additive neutral elements and 
additive inverses.    
\end{exer}

If $R$ and $R'$ are ordered and $f$ is a~bijective homomorphism preserving both
orders, we call $f$ an \underline{isomorphism of ordered rings\index{ring!ordered!isomorphism|emph}}. We 
know from semirings that bijective homomorphisms preserve both neutral 
elements. A~ring 
$R$ is \underline{simple\index{ring!simple|emph}} if for any set $X\sus R$ the following holds. If 
$0_R\in X$ and if for any $x\in X$ also $x+1_R\in X$ and $x+(-1_R)\in X$, then 
$X=R$. We show that zero is annihilating in multiplication in every ring.

\begin{lemma}\label{lem_zeroPrRin}
Let $R$ be a~ring. Then $a\cdot0_R=0_R\cdot a=0_R$ 
for every $a\in R$.
\end{lemma}
\duk
Let $a\in R$ be arbitrary. We write $0=0_R$ and compute
\begin{eqnarray*}
0&=&a\cdot0+(-(a\cdot0))=
a\cdot(0+0)+(-(a\cdot0))\\
&=&(a\cdot0+a\cdot0)+(-(a\cdot0))=a\cdot0+(a\cdot0+(-(a\cdot0)))\\
&=&a\cdot0+0=a\cdot0\,.
\end{eqnarray*}
\kduk

\noindent
We state the characterization theorem.

\begin{thm}\label{thm_Zis}
The structure $\Z$ is an ordered simple ring. Every two ordered simple 
rings are isomorphic as ordered rings.  
\end{thm}
Again, the latter claim is not valid without the adjective ``ordered'' 
because there exist both finite and infinite simple rings.

\begin{exer}\label{ex_notDoma}
List some examples of finite simple rings. 
\end{exer}
We prove the former claim in the theorem in Proposition~\ref{prop_ZtRing}, and 
the latter claim in Proposition~\ref{prop_isomorpRi}. Both 
proofs require several lemmas. 

\section{$\Z$ is an ordered simple ring}

The next two lemmas are used in the proof of Proposition~\ref{prop_ZtRing}.

\begin{lemma}\label{lem_aux2}
Let $m,n\in\Z$. Then
$$
m+(n+1)=(m+1)+n=(m+n)+1\,.
$$
\end{lemma}
\duk
Due to the commutativity of $+$, it suffices to show that for every $m,n\in\Z$, 
$$
(m+1)+n=(m+n)+1\,. 
$$
We check this identity by dividing it into eleven cases. We use Definition~\ref{def_addZ} and properties of subtraction in $\N_0$.
\begin{enumerate}
\item Let $m\ge0$ and $n\ge0$. Then $(m+1)+n=(m+n)+1$ by the properties of addition in $\N_0$.
\item Let $m\ge0$ and $-m\le n<0$. Then $(m+1)+n=m+1-|n|=m-|n|+1=
(m+n)+1$.
\item Let $m\ge0$ and $n=-(m+1)$. Then $(m+1)+n=(m+1)-|n|=0=-1+1=(m+n)+1$. 
\item  Let $m\ge0$ and $n<-(m+1)$. Then $(m+1)+n=-(|n|-(m+1))=-(|n|-m-1)=-(|m+n|-1)=(m+n)+1$. 
\item Let $m=-1$ and $n\ge1$. Then 
$(m+1)+n=0+n=n=(n-1)+1=(m+n)+1$.
\item Let $m=-1$ and $n=0$. Then $(m+1)+n=0=-1+1=(m+n)+1$.
\item Let $m=-1$ and $n<0$. Then
$(m+1)+n=n=-(|m+n|-1)=(m+n)+1$.
\item Let $m<-1$ and $n\ge|m|$. Then
$(m+1)+n=-(|m|-1)+n=n-(|m|-1)=(n-|m|)+1=(m+n)+1$.
\item Let $m<-1$ and $n=|m|-1$. Then
$(m+1)+n=-(|m|-1)+n=0=-1+1=(m+n)+1$.
\item Let $m<-1$ and $0\le n<|m|-1$. Then $(m+1)+n=-(|m|-1)+n=-(|m|-1-n)=-(|m|-n-1)=-(|m|-n)+1=(m+n)+1$.
\item Finally, let $m<-1$ and $n<0$. Then $(m+1)+n=-(|m|-1)+n=-(|m|-1+|n|)=-(|m+n|-1)=(m+n)+1$.
\end{enumerate}
\kduk

\begin{exer}\label{ex_idenProp}
Identify the properties of subtraction in $\N_0$ that were used in the previous proofs.    
\end{exer}

\begin{lemma}\label{lem_uspZ}
Let $m,n\in\Z$. Then $m<n$ $\iff$ $m+l=n$ for some $l\in\N$.    
\end{lemma}
\duk
Let $m,n\in\Z$ and $m<n$. If $m,n\in\omega$, then $m+(n-m)=n$ for 
some $n-m\in\omega$ by Proposition~\ref{prop_subtrOm}, and 
$n-m>0$. If $m<0\le n$, then 
$m+(|m|+n)=n$ with 
$|m|\in\N$, so that $|m|+n>0$. If $m,n<0$, then $|n|<|m|$, 
so that $|n|+(|m|-|n|)=|m|$ with $|m|-|n|>0$. We have  
$$
m+(|m|-|n|)=-(|m|-(|m|-|n|))=-|n|=n\,.
$$

Let $m,n\in\Z$ and let $m+l=n$ for some $l\in\N$. We 
prove by induction on $l$ that $m<n$. Let $l=1$. If $m\in\omega$, then 
$m\in m\cup\{m\}=m+1=n$ and $m<n$. If $m=-1$, then $n=m+1=0$ and $m<n$. If 
$m<-1$, then $n=m+1=-(|m|-1)$ and $|n|<|m|$, so that $m<n$. Let $l>1$ and 
$n'=m+(l-1)$. Then $m<n'$ by induction, $n'<n'+1=n$ by 
Lemma~\ref{lem_aux2}  and by the case $l=1$, and $m<n$ by the transitivity 
of $<$.
\kduk

We show that $\Z$ is an ordered simple ring.

\begin{prop}\label{prop_ZtRing}
The algebraic structure 
$$
\Z=\langle
\Z,\,0,\,1,\,+,\,\cdot,\,<\rangle
$$ 
introduced in Definition~\ref{def_algStZ} is 
an ordered simple ring.    
\end{prop}
\duk
We already noted that $+$ and $\cdot$ are commutative, that $0$ is 
neutral to $+$ and $1$ to $\cdot$, and that every integer 
has an additive inverse. The associativity of multiplication follows at once from the associativity of 
multiplication in $\N_0$ because the sign of the result depends only on 
the number of minus signs of the three factors. We prove the distributivity of 
$\cdot$ to $+$. Let $l,m,n\in\Z$. If $m,n<0$ or $m,n\ge0$, then the distributive 
identity 
$$
l\cdot(m+n)=l\cdot m+l\cdot n
$$ 
holds because it holds in $\N_0$. 
If $m<0\le n$ or $n<0\le m$, then the distributive 
identity follows from Corollary~\ref{cor_distrSubtr}. 

Least clear is the associativity of $+$. Suppose for the contrary that  
$l,m,n\in\Z$ are such that
$$
(l+m)+n\ne l+(m+n)
$$
and that $|l|+|m|+|n|$ is minimum. We may assume that one of the three numbers
is positive. Let, for example, $m>0$. 
Then 
$$
(l+(m-1))+n=l+((m-1)+n)\,.
$$
We add $1$ to the equation, use the commutativity of $+$, apply four times
Lemma~\ref{lem_aux2}, and obtain the contradiction that 
$$
(l+m)+n=l+(m+n)\,.
$$ 
The other two cases, $l>0$ and $n>0$, are similar.

We show that the two order axioms hold in $\Z$. 
Let $k,m,n\in\Z$ with $m<n$. By Lemma~\ref{lem_uspZ}, $m+l=n$ for 
some $l\in\N$. Thus, since $+$ is commutative and associative, also $(m+k)+l=n+k$. By the same lemma, $m+k<n+k$, which 
proves the first order axiom. Let also $k>0$. Then, by the distributivity of 
$\cdot$,
$$
m\cdot k+l\cdot k=n\cdot k\,.
$$
Since $l\cdot k>0$ by Definition~\ref{def_mulZ}, we have $m\cdot 
k<n\cdot k$ by Lemma~\ref{lem_uspZ}, which proves the second order axiom.

Finally, we show that the ring $\Z$ is simple. Let $0\in X\sus\Z$ and let $X$ be closed to 
adding $1$ and $-1$. Since $\N_0$ is simple, $\omega\sus X$. Since any 
$-n$ in $-\N$ arises by adding $-1$ exactly $n$ times to $0$, also 
$-\N\sus X$. Thus $X=\Z$.
\kduk

\section{$\Z$ is unique}

In order to prove the isomorphism claim in Theorem~\ref{thm_Zis}, 
we introduce for any given ring $R$ a~ring homomorphisms from $\Z$ to 
$R$.

\begin{defi}\label{def_mapsfR}
Let
$$
R=\langle
R,\,0_R,\,1_R,\,\oplus,\,\odot
\rangle
$$
be a~ring. We define a~map $f_R\cc\Z\to R$ as an extension of the map 
$f_R\cc\omega\to R$ introduced in Definition~\ref{def_evalM} by 
setting
$$
f_R(-m):=-f_R(m)
$$
for every $-m\in-\N$.
\end{defi}
In the following four lemmas, 
$$
R=\langle R,\,0_R,\,1_R,\,\oplus,\,\odot\rangle
$$ 
is a~ring and $f_R\cc\Z\to R$ is 
the map in Definition~\ref{def_mapsfR}.
The lemmas show that $f_R$ is a~ring homomorphism, and that if $R$ is simple 
and ordered, then $f_R$ is an isomorphism of ordered rings. 

\begin{lemma}\label{lem_lemaR1}
For every $m,n\in\Z$, we have $f_R(m+n)=f_R(m)\oplus f_R(n)$.
\end{lemma}
\duk
We write $f$ for $f_R$.
Let $m,n\in\Z$. If $m,n\ge0$, then the equality holds by 
Lemma~\ref{lem_lema1}. If $m<0\le n$ and $|m|>n$, then 

\begin{eqnarray*}
f(m+n)&=&f(-(|m|-n))=-f(|m|-n)=-(f(|m|)\ominus f(n))\\
&=&-f(|m|)\oplus f(n)=f(m)\oplus f(n)\,.    
\end{eqnarray*}
The first equality follows from Definition~\ref{def_addZ}. The second 
equality follows from Definition~\ref{def_mapsfR}. The third equality follows 
from Lemma~\ref{lem_lema1} and the fact that $(m-n)+n=m$. The fourth 
equality follows from properties of subtraction in $R$. In the 
last equality, we use Definition~\ref{def_mapsfR}.

Let $m<0\le n$ and $|m|\le n$. Then
$$
f(m+n)=f(n-|m|)=f(n)\ominus f(|m|)=
-f(|m|)\oplus f(n)=f(m)\oplus f(n)\,, 
$$
with similar justifications. The case $n<0\le m$ is 
symmetric. Finally, if $m,n<0$, then
\begin{eqnarray*}
f(m+n)&=&f(-(|m|+|n|))=-f(|m|+|n|)=
-(f(|m|)\oplus f(|n|))\\
&=&-f(|m|)\oplus(-f(|n|))=f(m)\oplus f(n)\,,
\end{eqnarray*}
again with similar justifications using Lemma~\ref{lem_lema1}.
\kduk

\begin{lemma}\label{lem_lemaR2}
For every $m,n\in\Z$, we have $f_R(m\cdot n)=f_R(m)\odot f_R(n)$.
\end{lemma}
\duk
We write $f$ for $f_R$. Let $m,n\in\Z$. If $m=0$ or $n=0$, then 
$$
f(m\cdot n)=f(0)=0_R=f_R(m)\odot f_R(n)\,.
$$
In the third  equality, we use Lemma~\ref{lem_zeroPrRin}.
If $m,n>0$ or $m,n<0$, then
\begin{eqnarray*}
f(m\cdot n)&=&f(|m|\cdot|n|)
=f(|m|)\odot f(|n|)\\
&=&(\pm f(m))\odot
(\pm f(n))=f(m)\odot f(n)\,,    
\end{eqnarray*}
with all upper or all lower signs. If $m<0<n$ or $n<0<m$, then
\begin{eqnarray*}
f(m\cdot n)&=&f(-(|m|\cdot|n|))
=-f(|m|\cdot|n|)=-(f(|m|)\odot f(|n|))\\
&=&-((\pm f(m))\odot
(\mp f(n)))=f(m)\odot f(n)\,,    
\end{eqnarray*}
again with all upper or all lower signs. We used Lemma~\ref{lem_lema2}. 
\kduk

\begin{lemma}\label{lem_lemaR3}
Let $R$ be an ordered ring, with the linear order $\prec$. Then
for every $m,n\in\Z$,  
$$
m<n\iff f_R(m)\prec f_R(n)\,.
$$
\end{lemma}
\duk
We again write $f$ for $f_R$.
By the trichotomy of $<$, it suffices to prove the implication 
$$
m<n\Rightarrow f(m)\prec f(n)\,.
$$
Let $m,n\in\Z$ and $m<n$. If $m,n\ge0$, the implication holds 
by Lemma~\ref{lem_lema3}. Let $m<0\le n<|m|$. Then 
$$
f(m)=-f(|m|)\prec -f(n)\preceq f(n)\,.
$$
The second inequality follows from $f(n)\prec f(|m|)$ (by 
Lemma~\ref{lem_lema3}) via the first order axiom by adding $-f(n)+(-f(|m|))$. To obtain the last non-strict inequality, we add to 
$0_R=f(0)\preceq f(n)$ (by Lemma~\ref{lem_lema3}) the element $-
f(n)$ and get $-f(n)\preceq0_R$. Now 
$0_R\preceq f(n)$ and we use the transitivity of $\preceq$.

Let $m<0<|m|\le n$. Then by similar arguments,
$$
f(m)=-f(|m|)\prec f(0)\prec f(n)\,.
$$
Finally, let $m<n<0$. Then
$$
f(m)=-f(|m|)\prec -f(|n|)=f(n)
$$
by Lemma~\ref{lem_lema3} and the first order axiom, because $|n|<|m|$.
\kduk

We leave the proof of the last lemma as an exercise.

\begin{exer}\label{ex_asAnExe}
Prove the next lemma.    
\end{exer}

\begin{lemma}\label{lem_lemaR4}
Let $R$ be a~simple ordered ring, with the linear order denoted by $\prec$. Then the
map
$$
f_R\cc\Z=\langle\Z,\,0,\,1,\,+,\,\cdot,\,<\rangle\to
R=\langle R,\,0_R,\,1_R,\,\oplus,\,\odot,\,\prec\rangle
$$
 is an isomorphism of ordered rings.
\end{lemma}

We show that every two ordered simple rings are isomorphic.

\begin{prop}\label{prop_isomorpRi}
Suppose that
$$
R=\langle
R,\,0_R,\,1_R,\,+,\,\cdot,\,<
\rangle\,\text{ and }\,
R'=\langle
R',\,0_{R'},\,1_{R'},\,\oplus,\,\odot,\,\prec
\rangle
$$
are two ordered simple rings. Then $R$ and $R'$ are isomorphic, which means 
that there exists a~bijection $f\cc R\to R'$ such that for every $a,b\in R$,
$$
f(a+b)=f(a)\oplus f(b),\ 
f(a\cdot b)=f(a)\odot f(b),\,\text{ and  }\,a<b \iff f(a)\prec f(b)\,.
$$
\end{prop}
\duk
By Lemma~\ref{lem_lemaR4} and Exercises~\ref{ex_naHomom1} and \ref{ex_naHomom2}, the map
$$
f\cc R\to R',\ f:=f_{R'}(f_R^{-1})\,,
$$
where $f_R$ and $f_{R'}$ are as in Definition~\ref{def_mapsfR}, 
is the desired isomorphism.
\kduk

\noindent
This concludes the proof of Theorem~\ref{thm_Zis}.

\section{Difference integers}\label{sec_diffInte}

As we promised, we outline an alternative construction of integers. 
For $m,n\in\omega$, we set $m\ominus n:=\langle m,n\rangle$ and define 
a~relation $\sim$ on $\omega^2=\omega\times\omega$:
$$
m\ominus n\sim m'\ominus n'\iff
m+n'=m'+n\,.
$$

\begin{exer}\label{ex_rozd1}
Show that $\sim$ is an equivalence relation.     
\end{exer}

We set 
$$
\Z':=\omega^2/\!\sim
$$ 
and call the elements of $\Z'$ \underline{difference 
integers\index{difference integers!set|emph}}. The elementary school integers in 
Definition~\ref{def_integers} are finite, but every difference integer is an infinite 
set. We define addition and multiplication on $\omega^2$, and a~relation $<$ on $\omega^2$.

\begin{defi}\label{def_operZpr}
Let $m\ominus n$ and $m'\ominus n'$ be two difference integers. We define the following. 
\begin{enumerate}
\item $m\ominus n+m'\ominus n':=(m+m')\ominus(n+n')$.
\item $m\ominus n\cdot m'\ominus n':=(mm'+nn')\ominus(mn'+nm')$.
\item $m\ominus n<m'\ominus n'$ $\iff$
$m+n'<m'+n$.
\end{enumerate}
\end{defi}

\begin{exer}\label{ex_rozd2}
It is not hard to show that $+$, $\cdot$, and $<$ are congruent to 
$\sim$, so that $+$ and $\cdot$ are operations on $\Z'$, and $<$ is 
a~relation on $\Z'$.    
\end{exer}

\begin{exer}\label{ex_rozd4}
Show that $\langle\Z',<\rangle$ is a~linear order.    
\end{exer}

\begin{defi}\label{def_zscarou}
In the algebraic structure of \underline{difference integers\index{difference integers!structure|emph}}   
$$
\Z':=\langle
\Z',\,0_{\Z'},\,1_{\Z'},\,+,\,\cdot,\,<\rangle\,,
$$
$\Z'=\omega^2/\!\sim$, $0_{\Z'}:=[0\ominus0]_{\sim}$, $1_{\Z'}:=
[1\ominus0]_{\sim}$, and operations $+$ and $\cdot$, and relation $<$ are as 
in Definition~\ref{def_operZpr}.
\end{defi}
One can prove the following.
\begin{prop}\label{prop_ZprIs}
The algebraic structure $\Z'$ is an ordered simple ring.    
\end{prop}
For example, $\Z'$ is simple because $[m\ominus n]_{\sim}$ equals to the 
sum of $m$ copies of $1_{\Z'}=[1\ominus0]_{\sim}$ and $n$ copies of 
$-1_{\Z'}=[0\ominus1]_{\sim}$. 

\begin{exer}\label{ex_rozd3}
Show that the multiplication in $\Z'$ is associative.     
\end{exer}

\begin{exer}\label{ex_rozd5}
Prove both order axioms in $\Z'$.    
\end{exer}
Proposition~\ref{prop_isomorpRi} yields the following corollary.

\begin{cor}\label{cor_ZandZpr}
The algebraic structures $\Z$ and $\Z'$ are isomorphic as ordered rings.    
\end{cor}
\index{integers|)}

\chapter[From integers to fractions\checkmark]{From integers to fractions}\label{sec_ZtoQ}
\index{fractions|(}

In this chapter, we extend the ring $\Z$ of natural numbers that we constructed in the previous chapter to 
the field $\Q$ of fractions. The main result is 
the characterization Theorem~\ref{thm_Qis}, 
which asserts that fractions form an ordered simple field that is unique up to 
isomorphism. 

Definition~\ref{def_Q} introduces the set of fractions and Definition~\ref{def_ordQ} 
introduce the linear order on so called protofractions. 
Proposition~\ref{prop_itsLO} shows that this is also a~linear order on fractions.
Definition~\ref{def_addMulZ} introduces the addition and multiplication of 
fractions. The  algebraic structure of fractions comes on the stage in 
Definition~\ref{def_racCisla}. In Proposition~\ref{thm_fieEmb} we prove that 
unlike for semirings and rings, for fields there are no non-injective homomorphisms.
Then comes the the characterization Theorem~\ref{thm_Qis}. In 
Definition~\ref{def_mapsfF}, we define an homomorphism from $\Q$ to any {\em 
ordered} field. 

\section{The linear order of fractions}

For $m,n\in\Z$ we define
$\frac{m}{n}=m/n$ to be the pair $\langle m,n\rangle$ and set
$${\textstyle
Z:=\big\{\frac{m}{n}\cc\;m\in\Z,\,
n\in\Z\setminus\{0\}\big\}\,.
}
$$
The elements of $Z$ are called \underline{protofractions\index{protofractions|emph}}. In a~protofraction  
$\frac{m}{n}$, the integer $m$ is the \underline{numerator\index{protofractions!numerator|emph}}, 
and $n$ is the 
\underline{denominator\index{protofractions!denominator|emph}}. We consider an equivalence relation $\sim$ on $Z$ (Exercise~\ref{ex_shodnProQ}), defined by
$${\textstyle
k/l\sim m/n\iff k\cdot n=l\cdot m
}
$$
(the multiplication is in $\Z$). For $a/b\in Z$ we set
$$
[a/b]_{\sim}:=\{c/d\in Z\cc\;a/b\sim c/d\}\,.
$$

\begin{exer}\label{ex_shodnProQ}
Show that $\sim$ is an equivalence relation.    
\end{exer}

\begin{defi}\label{def_Q}
The set of equivalence blocks
$$
\Q:=Z/\!\sim
$$
is the set of \underline{fractions\index{fractions!set|emph}} (rational numbers). 
\end{defi}

\begin{exer}\label{ex_posDeno}
Every fraction can be represented by 
a~protofraction with a~positive denominator.    
\end{exer}

We define a~linear order on 
$\Q$. The definition is more complex than in $\N_0$ and $\Z$. All comparisons by $<$ 
in the next definition are in $\Z$.

\begin{defi}\label{def_ordQ}
We define a~relation 
$<_Z$ on protofractions as follows. Let $a/b,c/d\in Z$. For $bd>0$, we set $a/b<_Z c/d$ $\iff$
$ad<bc$. For $bd<0$, we set $a/b<_Z c/d$ $\iff$ $ad>bc$.    
\end{defi}
We cannot have $ad=0$ because $a,d\ne0$.
We derive some properties of $<_Z$.

\begin{prop}\label{prop_itsLO}
The following holds.
\begin{enumerate}
\item The relation $<_Z$ on $Z$ is irreflexive and transitive, and has the 
property that for every two protofractions $\frac{a}{b}$ and 
$\frac{c}{d}$, exactly one of
$${\textstyle
\frac{a}{b}<_Z\frac{c}{d},\ 
\frac{a}{b}>_Z\frac{c}{d},\,\text{ and }\,\frac{a}{b}\sim\frac{c}{d}
}
$$ 
holds.
\item If $\frac{a}{b}\sim\frac{a'}{b'}$ and $\frac{c}{d}\sim\frac{c'}{d'}$, 
then $\frac{a}{b}<_Z\frac{c}{d}$ $\iff$ $\frac{a'}{b'}<_Z\frac{c'}{d'}$.
\end{enumerate}
\end{prop}
\duk
1. Irreflexivity of $<_Z$ is clear. Let $\frac{a}{b}<_Z\frac{c}{d}$ and 
$\frac{c}{d}<_Z\frac{e}{f}$. Suppose that $bd,df>0$. Then $ad<bc$ and 
$cf<de$. It follows that 
$$
ad^2f<dfbc=bdcf<bd^2e\,\text{ and }\,af<be\,.
$$
Since $bf>0$, this means that $\frac{a}{b}<_Z\frac{e}{f}$. The 
other three cases are similar. 
Suppose, for example, that $bd>0$ but $df<0$. Then $ad<bc$ and $cf>de$. It follows that 
$$
ad^2f>dfbc=bdcf>bd^2e\,\text{ and }\,af>be\,.
$$
Since $bf<0$, this means that again $\frac{a}{b}<_Z\frac{e}{f}$. In the 
remaining two cases we proceed similarly. This proves that $<_Z$ is 
transitive.

Let $\frac{a}{b},\frac{c}{d}\in Z$. If $ad=bc$, then $\frac{a}{b}\sim
\frac{c}{d}$. If $ad\ne bc$, then according to if $ad<bc$ or $ad>bc$, and 
according to the sign of $bd$, we have exactly one of $\frac{a}{b}<_Z
\frac{c}{d}$ and $\frac{a}{b}>_Z
\frac{c}{d}$. This proves the last property of $<_Z$.

2. Let $\frac{a}{b}\sim
\frac{a'}{b'}$,  
$\frac{c}{d}\sim\frac{c'}{d'}$, and let
$\frac{a}{b}<_Z\frac{c}{d}$. It suffices to show that $\frac{a'}{b'}<_Z\frac{c'}{d'}$ as well. Thus 
$ab'=ba'$ and $cd'=dc'$. Let $bd>0$. Then $ad<bc$. If $b'd'>0$, we have by 
the second order axiom that
$$
bda'd'=adb'd'<bcb'd'=bdb'c'\,.
$$
hence, again by the second order axiom,
$a'd'<b'c'$ and $\frac{a'}{b'}<_Z\frac{c'}{d'}$. In the remaining 
three cases when $b'd'<0$ or/and $bd<0$, we argue similarly.
\kduk

\noindent
Using this proposition, we can define a~\underline{linear order\index{fractions!linear order|emph}} $<_{\Q}$ on $\Q$.

\begin{defi}\label{def_LOonQ}
Let $\al,\be\in\Q$ be two distinct fractions. We define
$$
\text{
$\al<_{\Q}\be$ $\iff$ $a/b<_Z c/d$ for any
$a/b\in\al$ and $c/d\in\be$}\,.
$$
\end{defi}
By Proposition~\ref{prop_itsLO}, this definition is correct, and $\langle\Q,<_{\Q}\rangle$ is a~linear order.

\section{Addition and multiplication}

In order to define the standard arithmetic on $\Q$, we begin with the arithmetic on 
protofractions. Let $a/b$ and $c/d$ be two protofractions. We define their sum and 
product by $a/b+c/d:=(ad+cb)/bd$ and $a/b\cdot c/d:=ac/bd$. We show that both 
operations are congruent with $\sim$.

\begin{prop}\label{prop_corrAriQ}
Let $a/b\sim a'/b'$ and $c/d\sim c'/d'$ be two pairs of equivalent protofractions. Then  
$a/b+c/d\sim a'/b'+c'/d'$ and 
$a/b\cdot c/d\sim a'/b'\cdot c'/d'$.    
\end{prop}
\duk
We have $ab'=ba'$ and $cd'=dc'$. Then 
$${\textstyle
\frac{a}{b}+\frac{c}{d}=\frac{ad+bc}{bd}\sim\frac{a'd'+b'c'}{b'd'}=\frac{a'}{b'}+\frac{c'}{d'}
}
$$ 
because 
$$
(ad+bc)b'd'=ab'dd'+cd'bb'
=a'd'bd+b'c'bd=(a'd'+b'c')bd\,. 
$$
Similarly, 
$${\textstyle
\frac{a}{b}\cdot\frac{c}{d}
=\frac{ac}{bd}\sim\frac{a'c'}{b'd'}=\frac{a'}{b'}\cdot
\frac{c'}{d'}
}
$$
because 
$$
acb'd'=ab'cd'
=ba'c'd=a'c'bd\,. 
$$
\kduk

Thus we can define \underline{addition\index{fractions!addition|emph}} and \underline{multiplication\index{fractions!multiplication|emph}} of fractions.

\begin{defi}\label{def_addMulZ}
Let $a/b$ and $c/d$ be two protofractions, and let 
$\al=[a/b]_{\sim}$ and $\be=[c/d]_{\sim}$ be the corresponding fractions. We define
$\al+\be:=[(ad+cb)/bd]_{\sim}$ and
$\al\cdot\be:=[ac/bd]_{\sim}$.
\end{defi}

A~protofraction $a/b$ is in \underline{lowest terms\index{protofractions!in lowest 
terms|emph}} if $a$ and $b$ are coprime, and $b>0$. We denote the set of protofractions in lowest terms by $Z'$.   

\begin{exer}\label{ex_lowTer0}
Show that any two different protofractions in lowest terms are nonequivalent in 
$\sim$. 
\end{exer}

\begin{exer}\label{ex_lowTer}
Show that every fraction contains exactly one protofraction in lowest 
terms.    
\end{exer}
These two exercises show that the map $\Q\ni\al\mapsto a/b\in Z'$, $a/b\in\al$, 
that sends any fraction to the unique protofraction in  lowest terms contained 
in it is a~bijection.

\section{The characterization theorem}

We define the algebraic structure $\Q$.

\begin{defi}\label{def_racCisla}
The  algebraic structure of \underline{fractions\index{fractions!structure|emph}}   
$$
\Q=\langle
\Q,\,0_{\Q},\,1_{\Q},\,+,\,\cdot,\,<_{\Q}\rangle
$$
consists of the set $\Q$ introduced in Definition~\ref{def_Q}, the elements 
$0_{\Q}:=[0/1]_{\sim}$ and $1_{\Q}:=[1/1]_{\sim}$ in $\Q$, the operations 
of addition $+$ and multiplication $\cdot$ on $\Q$ introduced in Definition~\ref{def_addMulZ}, 
and the linear order $<_{\Q}$ on $\Q$ introduced in Definition~\ref{def_ordQ}.
\end{defi}
Instead of $0_{\Q}$ and $1_{\Q}$, we write just $0$ and $1$. Instead of 
$<_{\Q}$, we write just $<$.

\begin{exer}\label{ex_hezkeCvic}
Show that $0_{\Q}\ne 1_{\Q}$.    
\end{exer}

We review the algebraic structure of an ordered field that fits the fractions. Recall 
Definition~\ref{def_domain} of ordered rings.

\begin{defi}\label{def_field}
An \underline{ordered field\index{field!ordered|emph}} 
$$
F=\langle
F,\,0_F,\,1_F,\,+,\,\cdot,\,<
\rangle
$$
is an ordered ring such that for every $a\in F\setminus
\{0_F\}$ there is $b\in F$ satisfying $a\cdot b=1_F$. We call $b$ the 
\underline{multiplicative inverse\index{field!multiplicative inverses|emph}} of $a$ and denote it by $a^{-1}$. Without $<$ we 
have the structure of a~\underline{field\index{field|emph}}.  
\end{defi}

\begin{exer}\label{ex_ordFie2}
Show that multiplicative inverses in a~field are uniquely determined.   
\end{exer}
For $a,b\in F$ with $b\ne0_F$ we define
$$
a/b:=a\cdot b^{-1}
$$
and call $a/b$ the \underline{ratio\index{field!ratios|emph}}
of the elements $a$ and $b$. This partial \underline{division\index{field!division|emph}} operation on $F$ has the usual properties. If $F$ and $G$ are 
fields, a~map $f\cc F\to G$ is a~\underline{field 
homomorphism\index{field!homomorphism|emph}} if it is a~ring 
homomorphism. It is easy to see that field homomorphisms preserve 
multiplicative inverses. We show that  
field homomorphisms are injective. 

\begin{prop}\label{thm_fieEmb}
Every field homomorphism is injective.  
\end{prop}
\duk
Let $f\cc F\to K$ be a~field homomorphism. By Exercises~\ref{ex_AddInvUni} and 
\ref{ex_onRiHomo} we have
$$
0_K=f(0_F)=f(1_F+(-1_F))=f(1_F)+f(-1_F)=
1_K+f(-1_F)
$$
and $f(-1_F)=-1_K$. Since in any field $G$ we have that $-a=(-1_G)\cdot a$, we 
deduce that $f(a-b)=f(a)-f(b)$ for every $a,b\in F$. We show that $f^{-1}(0_K)
=\{0_F\}$. Then the injectivity of $F$ follows: $f(a)=f(b)$ is equivalent to 
$f(a-b)=0_K$.

Suppose that $f(a)=0_K$ for some element $a\in F\setminus\{0_F\}$. Then we obtain the contradiction that
$$
1_K=f(1_F)=f(a^{-1}\cdot a)=f(a^{-1})\cdot
f(a)=f(a^{-1})\cdot 0_K=0_K\,,
$$
because $0_K$ is annihilating in multiplication by Lemma~\ref{lem_zeroPrRin}.
\kduk

\noindent
From now on, we therefore speak of \underline{field embeddings\index{field!embedding|emph}} instead of field homomorphisms. If $F$ and $G$ are ordered 
fields and $f\cc F\to G$ is a~field embedding that preserves both
orders, we say that $f$ is an \underline{isomorphism of ordered 
fields\index{field!ordered!isomorphism|emph}}. A~field $F$ is 
\underline{simple\index{field!simple|emph}} if for any set $X\sus F$ the 
following holds. If $0_F\in X$ and if for any two elements $x,y\in X$ with 
$y\ne0_F$, also $x+1_F\in X$, $x-1_F=x+(-1_F)\in X$, and $x/y=x\cdot y^{-1}\in X$, then $X=F$.

We state the characterization theorem. 

\begin{thm}\label{thm_Qis}
The structure $\Q$ is an ordered simple field. Every two ordered simple 
fields are isomorphic as ordered fields.    
\end{thm}
We prove the former claim in the theorem in Proposition~\ref{prop_onFrac1}, and 
the latter claim in 
Proposition~\ref{prop_onFrac2}. 

\section{$\Q$ is an ordered simple field}

As before, we need several lemmas but we can prove 
Proposition~\ref{prop_onFrac1} right away.

\begin{prop}\label{prop_onFrac1}
The algebraic structure $\Q$ introduced in 
Definition~\ref{def_racCisla} is an ordered simple field.
\end{prop}
\duk
We leave the proof of neutrality of $0$ and $1$ to Exercise~\ref{ex_neutrQ}. The 
commutativity of 
$+$ and $\cdot$ follows from the commutativity of
these operations in $\Z$. The same holds for the associativity of 
$\cdot$. We prove that $+$ is associative. By the distributive law in $\Z$ we have
$$
{\textstyle
\big(\frac{a}{b}+\frac{c}{d}\big)+\frac{e}{f}=\frac{(ad+bc)f+bde}{bdf}=
\frac{adf+b(cf+de)}{bdf}=
\frac{a}{b}+\big(\frac{c}{d}+\frac{e}{f}\big)\,.
}
$$
We prove the distributive law:
$${\textstyle
\frac{a}{b}\cdot\big(\frac{c}{d}+\frac{e}{f}\big)=
\frac{a(cf+de)}{bdf}\sim
\frac{acbf+bdae}{b^2df}=
\frac{a}{b}\cdot\frac{c}{d}+
\frac{a}{b}\cdot\frac{e}{f}\,.
}
$$
The additive inverse of
$\frac{a}{b}$ is $\frac{-a}{b}$: 
$${\textstyle
\frac{a}{b}+\frac{-a}{b}=
\frac{ab+b(-a)}{b^2}=\frac{0}{b^2}\sim
\frac{0}{1}\,. 
}
$$
The multiplicative inverse 
of $\frac{a}{b}\not\sim\frac{0}{1}$, which means that $a\ne0$, is $\frac{b}{a}$: 
$${\textstyle
\frac{a}{b}\cdot\frac{b}{a}=\frac{ab}{ba}\sim\frac{1}{1}\,.
}
$$

We show that the two order axioms hold in $\Q$. Let $\frac{a}{b}<\frac{c}{d}$ 
and $\frac{e}{f}$ be three
protofractions with $b,d,f>0$. Thus $ad<bc$. Then also 
$${\textstyle
\frac{a}{b}+\frac{e}{f}<\frac{c}{d}+
\frac{e}{f}
}
$$
because $(af+be)df<(cf+de)bf$ $\iff$ $adf^2<bcf^2$ $\iff$ $ad<bc$ in $\Z$. Suppose that in addition $\frac{e}{f}>\frac{0}{1}$, which means that $e>0$. Then
$${\textstyle
\frac{a}{b}\cdot\frac{e}{f}<
\frac{c}{d}\cdot
\frac{e}{f}
}
$$
because $aedf<bfce$ $\iff$ $ad<bc$ in $\Z$.

Finally, we show that the field $\Q$ is simple. Let $\al=[m/n]_{\sim}$ 
be any fraction. We can easily obtain the fractions $\be=[m/1]_{\sim}$ and 
$\ga=[n/1]_{\sim}$ by
adding and subtracting $1_{\Q}$ to and from $0_{\Q}$. Then we express $\al$ as $\al=\be/\ga$. 
\kduk

\begin{exer}\label{ex_neutrQ}
Show that $0_{\Q}$ is neutral to $+$ and $1_{\Q}$ to $\cdot$.    
\end{exer}

\section{$\Q$ is unique}

In order to prove the isomorphism claim in Theorem~\ref{thm_Qis}, 
we again introduce, for any given ordered field $F$, a~field embedding 
from $\Q$ to $F$.

\begin{defi}\label{def_mapsfF}
Let
$$
F=\langle
F,\,0_F,\,1_F,\,\oplus,\,\odot,\,\prec
\rangle
$$
be an ordered field. We define a~function $f_F\cc\Q\to F$ by means of the ring homomorphism 
$f_F\cc\Z\to F$ introduced in Definition~\ref{def_mapsfR}. We set
$$
f_F(\al)=f_F(a)\oslash f_F(b):=f_F(a)\odot f_F(b)^{-1},\ \al=[a/b]_{\sim}\in\Q\,.
$$
\end{defi}
Differences appear compared to the structures $\N_0$ and $\Z$. We have 
to employ order because the field $\Q$ does not embed in every other field, only in the ordered ones. We also have to show, in the next two lemmas, that the definition of the map 
$f_F$ is correct.

\begin{lemma}\label{lem_deffFcorr0}
Let $F$ be as in Definition~\ref{def_mapsfF}. Then the ring homomorphism 
$f_F\cc\Z\to F$ has the property that $f_F(m)\ne0_F$ for every $m\in\Z
\setminus\{0\}$.
\end{lemma}
\duk
By Lemma~\ref{lem_lemaR3}, the inequality $m<0$, respectively $0<m$, 
is preserved by $f_F$. We have
$$
\text{$f_F(m)\prec f_F(0)=0_F$, respectively $f_F(0)=0_F\prec 
f_F(m)$}\,,
$$
which means that for nonzero $m\in\Z$ also $f_F(m)\ne0_F$.
\kduk

\begin{lemma}\label{lem_deffFcorr}
Let $F$ and the ring homomorphism $f_F\cc\Z\to F$ be as in Definition~\ref{def_mapsfF}, and let $\frac{a}{b}\sim\frac{c}{d}$ be 
two equivalent protofractions. Then 
$$
f_F(a)\oslash f_F(b)=f_F(c)\oslash f_F(d)\ \ (\in F)\,.
$$
\end{lemma}
\duk
$\frac{a}{b}\sim\frac{c}{d}$ means that $a\cdot d=b\cdot c$,  multiplied 
in $\Z$. By Lemma~\ref{lem_lemaR2},
$$
f_F(a)\odot f_F(d)=f_F(a\cdot d)=f_F(b\cdot c)=f_F(b)\odot f_F(c).
$$
The stated equality follows via multiplication by $f_F(d)^{-1}\odot 
f_F(b)^{-1}$. This element is defined due to the previous lemma.
\kduk

In the next four lemmas,
$$
F=\langle
F,\,0_F,\,1_F,\,\oplus,\,\odot,\,\prec
\rangle
$$
is an ordered field and $f_F\cc\Q\to F$ is as in Definition~\ref{def_mapsfF}.
The lemmas show that $f_F$ is a~field embedding, and that if $F$ is simple, then $f_F$ is an isomorphism of ordered fields.

\begin{lemma}\label{lem_prvnifF}
For every $\al,\be\in\Q$, we have $f_F(\al+\be)=f_F(\al)\oplus f_F(\be)$. 
\end{lemma}
\duk
We write $f$ for $f_F$. Let $\al=[a/b]_{\sim}$ and $\be=[c/d]_{\sim}$ 
be two fractions. We compute
\begin{eqnarray*}
f(\al+\be)&=&f([(ad+bc)/bd]_{\sim})=
f(ad+bc)\oslash f(bd)\\
&=&((f(a)\odot f(d)\oplus f(b)\odot f(c))\oslash (f(b)\odot f(d))\\
&=&(f(a)\oslash f(b))\oplus(f(c)\oslash f(d))=
f(\al)\oplus f(\be)\,.
\end{eqnarray*}
The first equality follows from the definition of $+$ in $\Q$. The second
equality is due to the definition of $f$. In the third equality, we use
Lemmas~\ref{lem_lemaR1} and 
\ref{lem_lemaR2}. The fourth equality is a~computation in the field $F$. In 
the last equality, we use the definition of $f$.
\kduk

\begin{lemma}\label{lem_druhefF}
For every $\al,\be\in\Q$, we have $f_F(\al\cdot\be)=f_F(\al)\odot f_F(\be)$.    
\end{lemma}
\duk
We write $f$ for $f_F$. Let $\al=[a/b]_{\sim}$ and $\be=[c/d]_{\sim}$ 
be two fractions. We compute
\begin{eqnarray*}
f(\al\cdot\be)&=&f([ac/bd]_{\sim})=
f(ac)\oslash f(bd)=((f(a)\odot f(c))\oslash(f(b)\odot f(d))\\
&=&(f(a)\oslash f(b))\odot(f(c)\oslash f(d))=f(\al)\odot f(\be)\,.
\end{eqnarray*}
The first equality follows from the definition of $\cdot$ in $\Q$. The second
equality is due to the definition of $f$. In the third equality, we use
Lemma~\ref{lem_lemaR2}. The fourth equality is a~computation in the field $F$. In 
the last equality, we use the definition of $f$.
\kduk

\begin{lemma}\label{lem_tretifF}
For every $\al,\be\in\Q$ we have  
$$
\al<\be\iff f_F(\al)\prec f_F(\be)\,.
$$
\end{lemma}
\duk
We write $f$ for $f_F$.
By the trichotomy of $<$, it suffices to prove the implication 
$$
\al<\be\Rightarrow f(\al)\prec f(\be)\,.
$$
Let $\al=[a/b]_{\sim}$ and 
$\be=[c/d]_{\sim}$ be two fractions such that $\al<\be$ and $b,d>0$. Thus 
$ad<bc$ and, by Lemma~\ref{lem_lemaR3}, $f(b),f(d)\succ 0_F=f(0)$. We have
$$
f(a)\odot f(d)=f(ad)\prec f(bc)=
f(b)\odot f(c)\,.
$$
The first inequality follows from Lemma~\ref{lem_druhefF}. The second 
inequality  follows from Lemma~\ref{lem_lemaR3}; we view $f$ as 
the ring homomorphism $f_F\cc\Z\to F$. The last inequality follows 
from Lemma~\ref{lem_druhefF}. The second order axiom yields, after 
multiplying the last displayed inequality by the element $f(b)^{-1}\odot f(d)^{-1}\succ0_F$, that
$$
f(\al)=f(a)\oslash f(b)\prec f(c)\oslash f(d)=f(\be)\,.
$$
\kduk

We leave the proof of the last lemma again as an exercise.

\begin{exer}\label{ex_opetLast}
Prove the next lemma.     
\end{exer}

\begin{lemma}\label{lem_ctvrtefF}
Let $F$ be in addition simple. 
Then the map
$$
f_F\cc\Q=\langle\Q,\,0,\,1,\,+,\,\cdot,\,<\rangle\to
F=\langle F,\,0_R,\,1_R,\,\oplus,\,\odot,\,\prec\rangle
$$
is an isomorphism of ordered fields.
\end{lemma}

\begin{prop}\label{prop_onFrac2}
Suppose that
$$
F=\langle
F,\,0_F,\,1_F,\,+,\,\cdot,\,<
\rangle\,\text{ and }\,
G=\langle
G,\,0_G,\,1_G,\,\oplus,\,\odot,\,
\prec\rangle
$$
are two ordered simple fields. Then $F$ and $G$ are isomorphic, which means 
that there exists a~bijection $f\cc F\to G$ such that for every $a,b\in F$,
$$
f(a+b)=f(a)\oplus f(b),\ 
f(a\cdot b)=f(a)\odot f(b),\,\text{ and  }\,a<b \iff f(a)\prec f(b)\,.
$$    
\end{prop}
\duk
By Lemma~\ref{lem_ctvrtefF} and Exercises~\ref{ex_naHomom1} and \ref{ex_naHomom2}, the map
$$
f\cc F\to G,\ f:=f_G(f_F^{-1})\,,
$$
where $f_F$ and $f_G$ are as in Definition~\ref{def_mapsfF}, 
is the desired isomorphism.
\kduk

\noindent
This concludes the proof of Theorem~\ref{thm_Qis}.
\index{fractions|)}

\chapter[From fractions to extended reals]{From fractions to extended reals}\label{sec_QtoRstar}

The reader may have expected that now we would proceed in the usual way and 
extend the ordered field $\Q$ of fractions to the complete ordered field $\R$ of real 
numbers. However, they also see from the title of the chapter that instead, 
something else will 
happen. We will make an original 
turn in the development of the basic numeric domains, and we will show how the domain of extended reals $\R^*=\R\cup
\{-\infty,+\infty\}$ naturally and uniquely arises from the linear order
$\langle\Q,<\rangle$ of fractions. We resolve the problem of how to 
obtain in some canonical way the infinities $\pm\infty$ as limits of 
rational sequences. These cannot be metric limits, which are based on the 
arithmetic operations in $\Q$, because the distance of $\pm\infty$ from any fraction 
is $+\infty$. We use topological limits, which are based on the order 
topology defined by $\langle\Q,<\rangle$. There are two main results in the chapter. 1. 
We first prove the characterization Theorem~\ref{thm_RstarIs} that 
$\R^*=\langle\R^*,<\rangle$ 
is, up to isomorphism of linear orders,  a~unique limit completion of $\Q=\langle\Q,
<\rangle$. This means that $\Q$ is a~suborder of $\R^*$, that every sequence 
$(a_n)\sus\Q$ has a~subsequence with a~limit in $\R^*$, and that the set 
$\R^*\setminus\Q$ is inclusion-wise minimal to the previous property. 2. 
In the latter half of the chapter, we extend, by means of limits, the arithmetic 
in $\Q$ to $\R^*$, and obtain in this way the canonical algebraic structure $\R^*$.

In Definition~\ref{def_openInt}

\section{Limits in linear orders}

Recall that $X=\langle X,<\rangle$ is a~{\em linear order} if $<$ is 
an irreflexive, transitive, and trichotomic relation on a~nonempty set $X$. If 
$Y=\langle Y,\prec\rangle$ is another linear order, an {\em isomorphism} of 
$X$ and $Y$ is any bijection $f\cc X\to Y$ such that for every $a,b\in X$, 
$$
a<b\iff f(a)\prec f(b)\,.
$$
Due to the trichotomy of linear orders, $\iff$ may be weakened to $\Rightarrow$ or 
to $\Leftarrow$.

Let $X=\langle X,<\rangle$ be 
a~linear order. A~set $I\sus X$ is an \underline{interval\index{intervals in linear orders|emph}} in $X$ if for every 
$a,b,c\in X$ with $a<b<c$ and $a,c\in I$,  also $b\in I$. 

\begin{defi}\label{def_openInt}
Let $X=\langle X,<\rangle$ be 
a~linear order. The \underline{open intervals\index{intervals in linear orders!open|emph}} in $X$ are the intervals
\begin{eqnarray*}
\{x\}&\ds&\text{if $X=\{x\}$}\,,\\
(a,\,b)&:=&
\{x\in X\cc\;a<x<b\}\,,\\
\,(a,\,\max(X)]&:=&
\{x\in X\cc\;a<x\le\max(X)\}\,,\text{ and }\\ 
\,[\min(X),\,a)&:=&
\{x\in X\cc\;\min(X)\le x<a\}\,,
\end{eqnarray*}
where $a,b\in X$. The latter  two families of intervals are defined only if $X$ has a~maximum, 
respectively, a~minimum.  
\end{defi}

\begin{exer}\label{ex_topBase}
Recall the two conditions a~set system on a~set has to satisfy to form a~base of a~topology 
on the set. Show that open intervals in a~linear order satisfy both conditions 
and therefore determine the \underline{order topology\index{linear order!order topology|emph}} in the linear order.   
\end{exer}

\begin{prop}\label{prop_isHausd}
The order topology is Hausdorff.    
\end{prop}
\duk
Let $X=\langle X,<\rangle$ be a~linear order and $a<b$ be any two distinct 
elements in $X$. We show that there exist disjoint open intervals $I$ and $J$ such 
that $a\in I$ and $b\in J$; these are the {\em separating intervals}. We call the 
situation when $a<b$ is a~jump in $X$ {\em case {\em (i)}}, and the complementary 
situation when $a<c<b$ for some $c\in X$ {\em case {\em (ii)}}. We distinguish four 
cases. If $a=\min(X)$ and $b=\max(X)$, then in case (i) the separating intervals are 
$[a,b)=\{a\}$ and $(a,b]=\{b\}$. In case (ii) they are $[a,c)$ and $(c,b]$. If $a'<a$
for some $a'\in X$ and $b=\max(X)$, then in case (i) the separating intervals are 
$(a',b)$ and $(a,b]=\{b\}$. In case (ii) they are $(a',c)$ and $(c,b]$. The case when
$a=\min(X)$ and $b<b'$ for some $b'\in X$ is treated symmetrically. Finally, if 
$a'<a$ and $b<b'$ for some $a',b'\in X$, then in case (i) the separating intervals 
are $(a',b)$ and $(a,b')$. In case (ii) they are $(a',c)$ and $(c,b')$.
\kduk

We introduce limits in linear orders.

\begin{defi}\label{def_limLO}
Let $X=\langle X,<\rangle$ be 
a~linear order, $(a_n)\sus X$ be a~sequence, and let $L\in X$. If for every open interval $I$ in $X$ with $L\in I$ there exists an index $n_0$ in $\N$ such that 
$a_n\in I$ for every $n\ge n_0$, we write 
$$
\lim a_n=L\,\text{ or }\,\lim_{n\to\infty}a_n=L
$$ 
and say that  the sequence $(a_n)$ has the limit $L$.
\end{defi}
It follows that limits in linear orders are unique.

\begin{cor}\label{cor_HausdTop}
Let $\langle X,<\rangle$ be a~linear order, $(a_n)\sus X$, and let $a,b\in 
X$ be such that $\lim a_n=a$ and $\lim a_n=b$. Then $a=b$.    
\end{cor}
\duk
This is immediate from Proposition~\ref{prop_isHausd} and Definition~\ref{def_limLO}.
\kduk

We leave the proof of the next lemma as an exercise for the interested reader.

\begin{lemma}\label{lem_onSubseq}
Let $\langle X,<\rangle$ be a~linear order, let $(a_n)\sus X$ have a~limit, 
and let $(b_n)$ be a~subsequence of $(a_n)$. Then $\lim b_n=\lim a_n$.    
\end{lemma}

We prove the next lemma.

\begin{lemma}\label{lem_limOrder}
Let $X=\langle X,<\rangle$ be a~linear order, $(x_n)\sus X$ be such that $\lim x_n=a$, and let $b\in X$. Then the following holds.
\begin{enumerate}
\item If $x_n\le b$ for every $n\ge n_0$, then $a\le b$.
\item If $a<b$, then $x_n<b$ for every $n\ge n_0$.  
\end{enumerate}
\end{lemma}
\duk
1. Suppose, for the contrary, that $b<a$. There exists an open interval $I$
in $X$ of the form $I=(b,\ds$ such that $a\in I$. Then $x_n\in I$ for every 
large $n$ and $x_n>b$ for the same $n$, which is a~contradiction. 

2. Similarly, there exists an open interval $I$ in $X$ of the form 
$I=\ds,b)$ such that $a\in I$. By the definition of limit, $x_n\in I$ for 
every large $n$. Thus $x_n<b$ for the same $n$.
\kduk

\noindent
The same result holds with the two $\le$ in part~1 replaced by $\ge$, and with the two $<$ in part~2 replaced by $>$.

\section{Limit completions of linear orders}

We introduce limit completions of linear orders.

\begin{defi}\label{def_limClos}
A~linear order $Y=\langle Y,<_Y\rangle$ is 
a~limit completion of a~linear order $X=\langle X,<_X\rangle$, if $X\sus Y$, $<_X$ is a~restriction of $<_Y$, that is, 
$$
<_X\,=\,<_Y\cap\, (X\times X)\,, 
$$
and if two conditions hold.
\begin{enumerate}
\item Every sequence $(a_n)\sus X$ has a~subsequence with a~limit in $Y$.
\item The set $Y\setminus X$ is inclusion-wise minimal with respect to condition~1. 
\end{enumerate}
\end{defi}
The next corollary is immediate.

\begin{cor}\label{cor_XdenseY}
If $\langle Y,<\rangle$ is a~limit completion of $\langle X,<\rangle$, 
then for any element $y\in Y$ there is a~sequence $(a_n)\sus X$ such that
$$
\lim_{n\to\infty}a_n=y\,.
$$
\end{cor}

Three natural questions are (i) whether any linear order has a limit completion,
(ii) whether this completion is unique, and (iii) whether this completion is 
itself complete to limits of subsequences. We hope to treat them in 
generality in 
here we 
resolve them only for the linear order $\langle\Q,<\rangle$.

In the Introduction, we have already mentioned cuts on $\langle\Q,<\rangle$. We 
remind them and define the linear order $\langle\R^*,<\rangle$. A set $X\sus\Q$
is a~{\em cut} if (i) for any $a,b\in\Q$ with $a<b\in X$, we also have 
$a\in X$ and (ii) for every $a\in X$,  there exists a~$b\in X$ with $a<b$.

\begin{defi}\label{def_Rstar}
Let $\R^*$ be the set of cuts. We call the elements of $\R^*$ also extended 
reals. We define the linear order $\langle\R^*,<\rangle$ by
$$
X<Y\iff X\sus Y\,.
$$
\end{defi}
It is clear that this is a~linear order. We define the cuts
$$
-\infty:=\emptyset,\ +\infty:=\Q,\,
\text{ and }\,\R:=\R^*
\setminus\{-\infty,\,+\infty\}\,.
$$
We call the last cuts in $\R$ the {\em real numbers}. We embed the linear 
order $\langle\Q,<\rangle$ in $\langle\R^*,<\rangle$ by viewing any 
fraction $a\in\Q$ as the {\em rational cut}  
$$
\{b\in\Q\cc\;b<a\}\,.
$$
Then $\Q\sus\R^*$ and we can state the characterization theorem for extended 
reals.

\section{The characterization theorem}

\begin{thm}\label{thm_RstarIs}
The linear order $\langle\R^*,<\rangle$ is a~limit completion of the linear order $\langle\Q,<\rangle$. Every 
two limit completions of $\Q$ are isomorphic as linear orders.    
\end{thm}
We prove the former claim in Proposition~\ref{prop_Rstar1}, and the 
latter in Proposition~\ref{prop_Rstar2}. We 
need several lemmas. 

\section{$\R^*$ is a~limit completion of fractions}

The next three are used in the proof of 
Proposition~\ref{prop_Rstar1}.  

Recall that a~sequence $(a_n)\sus A$ in
a~linear order $\langle A,<\rangle$ is {\em monotone} if
$$
a_1\le a_2\le\ds\,\text{ or }\,a_1\ge a_2\ge\ds\;.
$$
The next lemma is an infinite version of the famous Erd\H{o}s--Szekeres lemma.

\begin{lemma}\label{lem_ESlemma}
Any sequence in any linear order $\langle A,<\rangle$ has a~monotone 
subsequence.    
\end{lemma}
\duk
Let $(a_n)\sus A$. We define
$$
X=\{n\in\N\cc\;\forall\,m\cc\,m>n
\Rightarrow a_n>a_m\}\ \ (\sus\N)\,.
$$
If the set $X$ is infinite, 
$X=\{n_1<n_2<\ds\}$, then
$$
a_{n_1}>a_{n_2}>\ds
$$
is a~monotone subsequence of $(a_n)$. Suppose that $X$ is finite.  If $X=\emptyset$, we take any $n_1\in\N$.
If $X\ne\emptyset$, we take any $n_1\in\N$ with $n_1>\max(X)$. Since 
$n_1\not\in X$, there exists an index $n_2>n_1$ such that $a_{n_1}\le a_{n_2}$. Since 
$n_2\not\in X$, there exists an index $n_3>n_2$ such that $a_{n_2}\le a_{n_3}$. Continuing in this way, we obtain the monotone subsequence of $(a_n)$
$$
a_{n_1}\le a_{n_2}\le a_{n_3}\le\ds\,.
$$
\kduk

The next lemma shows that $\langle\R^*,<\rangle$ behaves nicely.

\begin{lemma}\label{lem_SupInfinExtRe}
In the linear order $\langle\R^*,<\rangle$, every set has an infimum 
and a~supremum.
\end{lemma}
\duk
Let $U\sus\R^*$ be any set of cuts. We prove the existence of suprema. We denote by $H(U)$ the set of upper bounds of the set $U$ in $\langle\R^*,<\rangle$. If $U=\emptyset$, then 
$$
\sup(U)=\min(H(U))=\min(\R^*)=-\infty\,.
$$
Let $U\ne\emptyset$. We set 
$X:=\bigcup U$ ($\ne\emptyset$). It is easy to see that $X\in\R^*$ and 
that $X\in H(U)$. We show that $\sup(U)=X$. Let $Y$ be any cut with $Y<X$.
We take any fraction $a\in X\setminus Y$. Then $a>b$ for every $b\in Y$. We have 
$a\in Z$ for some cut $Z\in U$. Then $Y<Z$ and we see that $Y\not\in H(U)$. Thus, 
$$
Y=\min(H(U))=\sup(U)\,. 
$$
We leave the existence of infima for 
the interested reader as an exercise.
\kduk

We show that $\Q$ is dense in $\R^*$. 

\begin{lemma}\label{lem_ratCutDens}
Rational cuts are dense in $\R^*$. In more details, for every cut $X$ there 
exists a~sequence $(a_n)\sus\Q$ of rational cuts such that 
$$
\lim_{n\to\infty}a_n=X.
$$
\end{lemma}
\duk
Let $X\in\R^*$ be any cut. If $X=\emptyset=-\infty$, then $\lim(-n)=X$. Similarly, if $X=\Q=+\infty$, then 
$\lim n=X$. Let $X\in\R$. Then there is $b\in\Q$ such that $b\ge 
X$. We define two sequences $(b_n)\sus\Q$ and $(k_n)\sus\N$
and initialize them by setting $b_1:=b$ and $k_1:=1$. If $b_1,\ds,b_n$ and 
$k_1,\ds,k_n$ are already defined, we set $c:=b_n-\frac{1}{k_n}$. If still $c\ge X$, 
we define $b_{n+1}:=c$ and $k_{n+1}:=k_n$. If $c<X$, we set $b_{n+1}:=b_n$ and 
$k_{n+1}:=k_n+1$. 
It is clear that $(b_n)\sus\Q$. We claim that
$$
\lim\,b_n=X\,.
$$
Let $I$ be any open interval in $\langle\R^*,<\rangle$ containing $X$.
Since $X\in\R$, we may assume that $I=(Y,Z)$ with $Y<Z$. We take any two 
fractions $c<d$ such that 
$$
Y<X<c<d<Z\,.
$$ 
It follows from the definition of 
sequences $(b_n)$ and $(k_n)$ and from $X\ne-\infty$, that $(k_n)$ is unbounded. Thus there exists an index 
$m\in\N$ such that $\frac{1}{k_m}<d-c$ and
$b_m-\frac{1}{k_m}<X$. Then 
$$
X\le b_n\le b_m<d\,\text{ and }\,b_n\in I
$$
for every $n\ge m$. Thus, $\lim b_n=X$. 
\kduk

\begin{prop}\label{prop_Rstar1}
The linear order $\langle\R^*,<\rangle$ is a~limit completion of the linear order $\langle\Q,<\rangle$.    
\end{prop}
\duk
We first show that every sequence $(a_n)\sus\Q$, with terms represented by 
rational cuts, has a~subsequence with a~limit in $\R^*$. Let 
$(a_n)\sus\Q$. By Lemma~\ref{lem_ESlemma}, $(a_n)$ has 
a~monotone subsequence $(b_n)$.  
We assume that $b_1\le b_2\le\ds$, use Lemma~\ref{lem_SupInfinExtRe}, and set
$$
X:=\sup(\{b_n\cc\;n\in\N\})\,.
$$
The other case, when $b_1\ge b_2\ge\ds$ and $X$ is the infimum of the stated set, 
is treated similarly. We show that $\lim b_n=X$. Let $I$ be any open interval in 
$\langle\R^*,<\rangle$ containing $X$, so that $I=(Y,Z)\ni X$ or $I=(Y,X]=(Y,+\infty]$. By the definition of 
supremum, there is an $m\in\N$ such that 
$$
Y<b_m\le X\,.
$$
Since $(b_n)$ weakly increases, the displayed inequalities hold also for any 
$n\ge m$ and $b_n\in I$ for these $n$. We see that $\lim b_n= X$.

The minimality condition~2 in Definition~\ref{def_limClos} follows from 
Lemma~\ref{lem_ratCutDens}. If we omit from $\R^*$ any cut $X$, or several 
cuts, then by Lemma~\ref{lem_ratCutDens} there is a~sequence $(a_n)\sus\Q$ (represented by rational 
cuts) such that $\lim a_n=X$. But then $(a_n)$ has no subsequence with a~limit 
in $\R^*\setminus\{X\}$ because every subsequence $(b_n)$ of $(a_n)$ has $\lim 
b_n=X$. We see that $\langle\R^*,<\rangle$ is a~limit completion of $\langle\Q,<\rangle$.
\kduk

\section{$\R^*$ is unique}

We proceed to the proof of the isomorphism claim in 
Proposition~\ref{prop_Rstar2}. We need five lemmas. In them, $\R'=\langle\R',<\rangle$ is a~limit completion of $\langle\Q,<\rangle$. Although 
$$
\Q=\langle\Q,\,0,\,1,\,+,\,\cdot,\,<
\rangle
$$ 
is an ordered field, no arithmetic operations are yet defined on the 
elements in $\R'\setminus\Q$, and we can approach them only by means of the 
linear order $<$. We have Cauchy sequences $(a_n)\sus\Q$, but it makes 
no sense to speak about Cauchy 
sequences in $\R'$. Let $\langle X,<\rangle$ be a~linear order and $a,b\in 
X$. We define the interval $[a,b]$ in $X$ by 
$$
[a,\,b]:=\{c\in X\cc\;a\le c\le b\}\,.
$$

\begin{lemma}\label{lem_pomoc2}
Let $\langle X,<\rangle$ be a~linear order, let $a,b,x\in X$, $(x_n)\sus[a,b]$, 
and let $\lim x_n=x$. Then $x\in[a,b]$. 
\end{lemma}
\duk
This follows from part~1 of Lemma~\ref{lem_limOrder}.
\kduk

\begin{lemma}\label{lem_pomoc0}
For every $\al<\be$ in $\R'$, there exists a~fraction $a\in\Q$ such that $\al<a<\be$. 
\end{lemma}
\duk
If $\al<\be$ are in $\Q$, the claim holds because the linear order $\Q$ is dense. We therefore assume that $\al\in\R'\setminus\Q$. If 
$\be\in\R'\setminus\Q$, we argue similarly. 
We consider three elements $\al<\ga<\be$ in $\R'$ and the open interval 
$(\al,\be)\ni\ga$. By Corollary~\ref{cor_XdenseY} and 
Definition~\ref{def_limLO}, there exists a~fraction $a\in(\al,\be)$. Thus if the 
stated claim does not hold, then there exist elements $\al<\be$ in $\R'$ ($\al\not\in\Q$) such that 
$\al<\ga<\be$ for no $\ga\in\R'$. We set $\R'':=\R'\setminus\{\al\}$ and show that 
every sequence $(a_n)\sus\Q$ has a~subsequence $(b_n)$ with a~limit in 
$\R''$; this contradicts the minimality condition~2 in 
Definition~\ref{def_limClos}. Let $(a_n)\sus\Q$. Since $\R'$ is a~limit completion of $\Q$, there is an element 
$\ga\in\R'$ and a~subsequence $(b_n)$ of $(a_n)$ such that $\lim b_n=\ga$ in $\R'$.
If $\ga\ne\al$, then also $\lim b_n=\ga$ in $\R''$ because every open interval in 
$\R''$ containing $\ga$ remains the same in $\R'$, or is enlarged only by $\al$. 
Let $\ga=\al$. We claim that then $\lim b_n=\be$ in $\R''$. Let $I''$ be an open 
interval in $\R''$ containing $\be$. Since $\al<\be$ is a~jump in the linear order 
$\R'$, the interval $I''$ becomes the open interval $I'=I''\cup\{\al\}$ in $\R'$. Then
$b_n\in I'$ for every large $n$ by 
the assumption on $(b_n)$. It follows that $b_n\in I''$ for the same $n$, and we see 
that indeed $\lim b_n=\be$ in $\R''$.  
\kduk

\begin{lemma}\label{lem_pomoc1}
There exist extremal elements 
$$
-\infty':=\min(\R')\,\text{ and }\,+\infty':=\max(\R')\,.
$$
\end{lemma}
\duk
We only prove the existence of $\max(\R')$; the element $\min(\R')$ is treated similarly.
The sequence $(n)\sus\Q$ has a~subsequence $(a_n)$ with a~limit 
$$
\lim_{n\to\infty}\,a_n=\al\in\R'\setminus\Q\,. 
$$
If $\al$ is not $\max(\R')$, then $\al<\be$ for some $\be\in\R'$. By 
the previous lemma, there exists a~fraction $b\in\Q$ such that $\al<b<\be$. 
Then we can take an open interval $I$ in $\R'$ such that $I=\ds,b)$ and $I\ni\al$.
Since $a_n>b$ for every large $n$, we have $a_n\not\in I$ for the same 
$n$, which is a~contradiction. 
\kduk

\noindent
From now on we denote by $\R''$ the finite elements $\R'
\setminus\{-\infty',+\infty'\}$.

\begin{lemma}\label{lem_1pomoc}
Let $\al\in\R'$, let $a,b\in\Q$, and let $(a_n)\sus\Q$ with $(a_n)\sus[a,b]$ 
be such that $\lim a_n=\al$. 
Then $(a_n)$ is a~Cauchy sequence.     
\end{lemma}
\duk
For $\al\in\Q$, this is a~standard result on convergent sequences in the  ordered field $\Q$ (see Proposition~\ref{prop_equiDefi}). Let $\al\in\R'\setminus\Q$. 
Let $e\in\Q$ with $e>0$ be given. We take $m\in\N$ such that $(b-a)/2^m\le e$.
Lemma~\ref{lem_pomoc2} implies $\al\in[a,b]$. We halve the interval $[a,b]$ into open intervals $I=(a,\frac{1}{2}(a+b))$ and 
$J=(\frac{1}{2}(a+b),b)$.  Then $\al\in I$
or $\al\in J$. By the definition of limits,
$(a_n)\sus I$ or $(a_n)\sus J$ for every large $n$. After $m$ such halvings,
we obtain an interval $K$ with rational endpoints and length at most $e$ such that 
$(a_n)\sus K$ for every large $n$. This means that $(a_n)$ is Cauchy. 
\kduk

\begin{lemma}\label{lem_2pomoc}
Let $(a_n)\sus\Q$ be a~Cauchy sequence, and let $(b_n)$ and  $(b_n')$ be two subsequences of 
$(a_n)$ that have limits in $\R'$. Then $\lim b_n=\lim b_n'$.  
\end{lemma}
\duk
Suppose for the contrary that $\lim b_n=\al<\be=\lim b_n'$, where $\al,\be\in\R'$. For $\al>\be$ we 
argue similarly. Using Lemma~\ref{lem_pomoc0}, we take 
fractions $c,d\in\Q$ such that 
$$
\al<c<d<\be\,.
$$
We take open intervals $I=\ds,c)$ and $J=(d,\ds$ such that $I\ni\al$ and $J\ni\be$.
By the assumption on $(b_n)$ and $(b_n')$, we have $b_n\in I$ for every 
large $n$ and, similarly, $b_n'\in J$ for every large $n$. It follows that 
for every $n_0$ there 
exist indices $m,n\in\N$ such that $m,n\ge n_0$ and $a_m-a_n>d-c$ ($>0$). This
contradicts the assumption that $(a_n)$ is Cauchy. 
\kduk

Now we can prove the uniqueness of the limit completion $\R^*$ of $\langle\Q,
<\rangle$. We prove a~slightly stronger result with $\Q$-isomorphisms. These 
are isomorphisms of linear orders that pointwise fix the common subset $\Q$. 

\begin{prop}\label{prop_Rstar2}
Every limit completion of the linear order of fractions is 
$\Q$-isomorphic as a~linear order to $\R^*$. Hence every two limit completions of 
the linear order of fractions are $\Q$-isomorphic as linear orders.    
\end{prop}
\duk
Let $\R':=\langle\R',\,\prec\rangle$ 
be a~limit 
completion of the linear order $\Q$. We define a~map $F\cc\R^*\to\R'$ such 
that $F$ is onto, $F(a)=a$ for every $a\in\Q$, and for any two cuts $X,Y\in\R^*$, 
$$
X<Y\Rightarrow F(X)\prec F(Y)\,.
$$
Then $F$ is a~$\Q$-isomorphism of the linear orders $\R^*$ and $\R'$. 

By Lemma~\ref{lem_pomoc1}, we can set 
$F(-\infty):=-\infty'$ and $F(+\infty):=+\infty'$. We define $F(X)$ for any real number $X\in\R$. 
By Lemma~\ref{lem_pomoc0}, $X\in(c,d)$ for some fractions 
$c,d\in\Q$. Since $\R^*$ is a~limit completion of $\Q$, by 
Corollary~\ref{cor_XdenseY} there is a~sequence 
$(a_n)\sus\Q$ such that $X=\lim a_n$ (in $\R^*$). By Definition~\ref{def_limLO}, we may assume that $(a_n)\sus[c,d]$. Since $\R'$ is a~limit 
completion of $\Q$, the 
sequence $(a_n)$ has a~subsequence $(b_n)$ with a~limit
$$
\lim_{n\to\infty}\,b_n=:F(X)\ \ (\in\R')\,.
$$
We show that the value $F(X)$ depends neither on the choice of the 
subsequence $(b_n)$ nor on the choice of the original sequence $(a_n)$. The former is 
immediate from Lemmas~\ref{lem_1pomoc} and \ref{lem_2pomoc}. To prove the latter, we 
consider another sequence 
$(a_n')\sus\Q$ with $\lim a_n'=X$ (in $\R^*$) and a~subsequence $(b_n')$ of it 
that has a~limit in $\R'$. But then both $(b_n)$ and $(b_n')$ are 
subsequences of the combined sequence
$$
(a_n''):=(a_1,\,a_1',\,a_2,\,a_2',\,a_3,\,\ds)
$$
that has also $\lim a_n''=X$. Thus $\lim b_n=\lim b_n'$ (in $\R'$) again by 
Lemmas~\ref{lem_1pomoc} and \ref{lem_2pomoc}. We see that the map 
$F\cc\R^*\to\R'$ is correctly defined. 
It is clear that $F(X)=X$ for every rational cut $X=a\in\Q$ because then we 
may take the constant sequence $(a_n)=(a,a,\ds)$.

We prove that $F$ is onto. Let $\al\in\R'$ be
an arbitrary element. By Corollary~\ref{cor_XdenseY}, $\al=\lim a_n$
(in $\R'$) for some sequence $(a_n)\sus\Q$. Since $\R^*$ is a~limit 
completion of $\Q$, we have a~subsequence $(b_n)$ of $(a_n)$  such that $\lim 
b_n=X\in\R^*$. Since $\R'$ is a~limit 
completion of $\Q$, we have a~subsequence $(c_n)$ of $(b_n)$  such that 
$\lim c_n=\be\in\R'$. But $(c_n)$ 
is a~subsequence of $(a_n)$, and therefore $\be=\al$. We see that $F(X)=\be=\al$.

Finally, we show that $F$ increases. Let $X<Y$ be two distinct cuts in $\R^*$. We 
prove that $F(X)\prec F(Y)$. By Lemma~\ref{lem_pomoc0}, $X<a<b<Y$ for 
two fractions $a,b\in\Q$. Using Corollary~\ref{cor_XdenseY}, we take 
sequences $(a_n),(b_n)\sus\Q$ such that $\lim a_n=X$ and $\lim b_n=Y$ (in 
$\R^*$). Let $(a_n')$ be a~subsequence of $(a_n)$ such that $F(X)=\lim a_n'$ 
(in $\R'$). Similarly, let $(b_n')$ be a~subsequence of $(b_n)$ such that 
$F(Y)=\lim b_n'$ (in $\R'$). Since $a_n'\le a$ and $b\le b_n'$ for every
large $n$ (by Definition~\ref{def_limLO}), by part~1 of Lemma~\ref{lem_limOrder} we have
$$
F(X)\preceq a\prec b\preceq F(Y)\,,
$$
and $F(X)\prec F(Y)$. 
\kduk

\noindent
This completes the proof of Theorem~\ref{thm_RstarIs}.

We have answered the earlier mentioned questions (i) and (ii) for $\langle\Q,<\rangle$: the linear 
order $\Q$ has a~limit completion and this completion is unique. Now we prove 
that it is complete to limits of subsequences.  

\begin{prop}\label{prop_RstarLCo}
Every sequence $(\al_n)\sus\R^*$ has a~subsequence with a~limit in the linear order $\R^*$.    
\end{prop}
\duk
This was already proven in the proof of Proposition~\ref{prop_Rstar1}: nothing 
changes in the argument if we take any sequence $(a_n)\sus\R^*$.
\kduk

\section{Addition and multiplication}

In the second part of this section, we extend, using limits, the addition and 
multiplication in the ordered field $\Q$ to the linear order of extended 
reals $\R^*$. We begin by proving the equivalence of the two definitions of 
limits that apply in ordered fields. Let
$$
\langle F,\,0_F,\,1_F,\,+,\,\cdot,\,<\rangle
$$ 
be an ordered field. For a~sequence $(a_n)\sus F$ and $L\in F$, we can 
define that
$$
\lim_{n\to\infty}a_n=L
$$
via Definition~\ref{def_limLO} or we can use the {\em classical definition} $-$ for every $e\in F$ with $e>0_F$ there is an index $n_0$ such for every $n\ge n_0$,
that
$$
|a_n-L|\le e\,.
$$
Recall that the {\em absolute value} $|\cdots|$ is 
defined for $a\in R$ in any ordered ring $R$ as $|a|:=a$ if $a\ge0_R$, 
and as $|a|:=-a$ if $a<0_R$. It has the familiar 
properties: (i) $|a|\ge 0_R$ and $|a|=0_R$ iff $a=0_R$, (ii) $|a\cdot 
b|=|a|\cdot|b|$, and (iii) the triangle inequality $|a+b|\le|a|+|b|$ holds.

\begin{prop}\label{prop_equiDefi}
Let 
$$
F=\langle F,\,0_F,\,1_F,\,+,\,\cdot,\,<\rangle
$$ 
be an ordered field. 
Then the above classical definition of limits in $F$ and 
Definition~\ref{def_limLO} (applied to the linear order $\langle F,<\rangle$) 
are equivalent.     
\end{prop}
\duk
Let $L\in F$ and let $(a_n)\sus F$. Suppose that $\lim a_n=L$ by the 
classical definition, and that $I$ is an open interval in the linear order 
$\langle F,<\rangle$ such that $L\in I$. Thus $I=(a,b)$ for some $a<b$ in 
$F$ because $|F|>1$, and  
$\min(F)$ and $\max(F)$ do not exist. We set 
$$
d:=\min(b-L,\,L-a)\ \ (>0_F)\,\text{ and }\,e:=d/2_F\,.
$$
Note that $2_F\ne 0_F$ because $F$ is ordered, and hence has characteristic 
$0$. There is $n_0\in\N$ such that $|a_n-L|\le e$ for every $n\ge n_0$. 
If $n\ge n_0$ and $a_n\not\in I$, then $a_n\le a$ or $a_n\ge b$. In both 
cases, $|a_n-L|\ge d>e$, which is a~contradiction. It 
follows that for every $n\ge n_0$ we have $a_n\in I$. So $\lim a_n=L$ also 
by Definition~\ref{def_limLO}.

Suppose that $\lim a_n=L$ by Definition~\ref{def_limLO}, and that 
$e\in F$ with $e>0_F$ is given. We consider the open interval 
$I:=(L-e,L+e)$. Since for any 
$c\in F$,
$$
c\in I\iff |c-L|<e\,,
$$
we see that $\lim a_n=L$ also by the classical definition.
\kduk

Let $\langle\R',\prec\rangle$ be any limit completion of $\langle\Q,<\rangle$, such as $\langle\R^*,
<\rangle$, and let $\al$ and $\be$ in $\R'$ be two arbitrary elements. Using Corollary~\ref{cor_XdenseY}, we 
take sequences $(a_n),(b_n)\sus\Q$ such that $\lim a_n=\al$ and $\lim b_n=\be$ 
(in $\R'$). Since $\R'$ is a~limit completion of $\Q$, there exists 
a~subsequence $(c_n)$ of $(a_n+b_n)$ such that we can set
$$
\al+\be:=\lim_{n\to\infty}c_n\ \ (\text{in $\R'$})\,.
$$
We similarly define
$$
\al\cdot\be:=\lim_{n\to\infty}d_n\,,
$$
where $(d_n)$ is a~subsequence of $(a_nb_n)$ that has a~limit in $\R'$.

\begin{defi}\label{def_spraDefi}
We define expressions $\al+\be$ and $\al\cdot\be$ by the above limits, provided that they do not depend 
on the choice of sequences $(a_n)$, $(b_n)$, $(c_n)$, and $(d_n)$. 
Else, the expression has no value. 
\end{defi}
In the next theorem, we determine for every $\al,\be\in\R'$ which case occurs.
Recall that an ordered field is {\em complete} if every nonempty and 
upper-bounded set in it has a~supremum. We use two lemmas.

\begin{lemma}\label{lem_narrInt}
Let $\langle\R',\prec\rangle$ be any limit completion of $\langle\Q,
<\rangle$ and let $\al\in\R''$. Then for every fraction $e>0$ there exist
fractions $a$ and $b$ such that $a\prec\al\prec b$ and $b-a<e$.
\end{lemma}
\duk
The existence of fractions $c$ and $d$ such that $\al\in(c,d)$ follows from 
Lemma~\ref{lem_pomoc0}. The fractions $a$ and $b$ are obtained by repeated 
halving of the interval $(c,d)$.
\kduk

\begin{lemma}\label{lem_euqiOF}
Let 
$$
F=\langle
F,\,0_F,\,1_F,\,+,\,\cdot,\,<
\rangle
$$
be an ordered field. The second order axiom has the equivalent form 
$$
a,\,b>0_F\Rightarrow a\cdot b>0_F\,.
$$
\end{lemma}
\duk
We deduce the new form from the old one. Let $a,b>0_F$. Then, by the old
form,
$$
a>0_F\leadsto a\cdot b>0_F\cdot b=0_F\,.
$$
We deduce the old form from the new one. Let $a<b$ and $c>0_F$. By the 
first order axiom, $b-a>0_F$. By the new form, distributivity, and again the 
first order axiom,
$$
b\cdot c-a\cdot c>0_F\leadsto
a\cdot c<b\cdot c\,.
$$
\kduk

\begin{thm}\label{thm_strExtRea}
Let $\langle\R',\prec\rangle$ be any limit completion of $\langle\Q,
<\rangle$ and let $\al+\be$ and $\al\cdot\be$, with $\al,\be\in\R'$, be 
predefined by the 
above limits. Then the following holds.
\begin{enumerate}
\item $\al+\be$ is always correctly defined, with the exception of two expressions $(-\infty')+(+\infty')$ and 
$(+\infty')+(-\infty')$, which have no value.
\item $\al\cdot\be$ is always correctly defined, with the exception of four 
expressions $0\cdot(+\infty')$, $0\cdot(-\infty')$, $(+\infty')\cdot0$, 
and $(-\infty')\cdot0$, which have no value.
\item The partial operations $+$ and $\cdot$ on $\R'$ extend the addition and multiplication in $\Q$.
\item Recall that $\R''=\R'\setminus\{-\infty',+\infty'\}$.   We have values
\begin{enumerate}
\item $\al+(\pm\infty')=(\pm\infty')+\al=\pm\infty'$ (same signs) for every $\al\in\R''$. 
\item $(\pm\infty')+(\pm\infty')=\pm\infty'$ (same signs).
\item $\al\cdot(\pm\infty')=(\pm\infty')\cdot\al=\pm\infty'$ (same signs) for every $\al\in\R'$ with $\al\succ0$.
\item $\al\cdot(\pm\infty')=(\pm\infty')\cdot\al=\mp\infty'$ (same signs) for every $\al\in\R'$ with $\al\prec0$.
\end{enumerate}
\item The algebraic structure
$$
\R'':=\langle
\R'',\,0,\,1,\,+,\,\cdot,\,\prec
\rangle\ \ (\R''=\R'\setminus\{-\infty',\,+\infty'\})
$$
is a~complete ordered field. 
\end{enumerate}
\end{thm}
\duk
Let $\al$, $\be$, $(a_n)$, and $(b_n)$ be as stated. 1. First, let  
$\al,\be\in\R''$. By Lemma~\ref{lem_pomoc0} and 
Definition~\ref{def_limLO}, there exist fractions $a,b$ such that 
$\al,\be\in[a,b]$ and $(a_n),(b_n)\sus[a,b]$ for every $n$. Thus for every $n$, 
$$
(a_n+b_n)\sus[2a,\,2b]\,.
$$
By Lemmas~\ref{lem_1pomoc} and \ref{lem_2pomoc}, $\al+\be:=\lim c_n$ 
does not depend on the choice of the subsequence $(c_n)$. We can prove that 
$\al+\be$ depends neither on the choice of $(a_n)$ and $(b_n)$ by using the 
combined sequences as in the proof of Proposition~\ref{prop_Rstar2}. Let 
$\al=\pm\infty'$ and $\be\in\R''$ or $\be=\al$. It is not hard to see that 
then $\lim(a_n+b_n)=\lim a_n=\al$. Thus 
$\al+\be$ does not depend on the choice of the subsequence $(c_n)$ by 
Lemma~\ref{lem_onSubseq}, and $\al+\be=\al$. The case when 
$\be=\pm\infty'$ and $\al\in\R''$ or $\al=\be$ is symmetric. Finally, we 
show that the two expressions
$$
(\pm\infty')+(\mp\infty')
$$
have no value. Let $\al=\pm\infty'$ and $(a_n)\sus\Q$ be such that $\lim 
a_n=\al$. Then for any $a\in\Q$ and $(b_n):=(-a_n+a)$,
$$
\lim\,(b_n)=-\al\,\text{ and }\,
\lim(a_n+b_n)=a\,,
$$
so that $\lim c_n=a$ for every subsequence $(c_n)$ of $(a_n+b_n)$. 
Thus, $\al+\be$ now depends on the choice of sequences $(a_n)$ and $(b_n)$,
and cannot be assigned any value.

2. We again begin with the case $\al,\be\in\R''$. The argument is the 
same as in 1, because if $(a_n),(b_n)\sus[-a,a]$, with a~fraction $a>0$, 
then $(a_nb_n)\sus[-a^2,a^2]$. Let $\al=\pm\infty'$ and $\be\succ0$. Again, it 
is not hard to see that 
then $\lim(a_nb_n)=\lim a_n=\al$. We argue as in 1 and deduce that 
$\al\cdot\be=\al$. The remaining three cases, when one of $\al$ and $\be$ is 
$\pm\infty'$ and the other is nonzero, are similar. We show that the four expressions
$$
(\pm\infty')\cdot0\,\text{ and }\,
0\cdot(\pm\infty')
$$
have no value. Let $\al=\pm\infty'$ and let $(a_n)\sus\Q$ be such that $\lim 
a_n=\al$ and $a_n\ne0$ for every $n$. Then for every $a\in\Q$ and $(b_n):=(a/a_n)$, 
$$
\lim\,b_n=0=\be\,\text{ and }\,\lim\,(a_nb_n)=a\,,
$$
so that $\lim d_n=a$ for every subsequence $(d_n)$ of $(a_nb_n)$. 
Similarly if $\al$ and $\be$ are exchanged. Thus, $\al\cdot\be$ now 
depends on the choice of sequences $(a_n)$ and $(b_n)$, and cannot be 
assigned any value.

3. Let $\al,\be\in\Q$. We select constant sequences $(a_n)=
(\al,\al,\ds)$ and $(b_n)=(\be,\be,\ds)$; $\lim a_n=\al$ and
$\lim b_n=\be$. Then $\lim(a_n+b_n)=\al+\be$ and by 
Lemma~\ref{lem_onSubseq}, $\lim c_n=\al+\be$ for every subsequence 
$(c_n)$ of $(a_n+b_n)$. We see that the addition on $\R'$ extends the addition 
in $\Q$. The argument for multiplication is similar.

4. We proved (a) and (b) already in part~1, and (c), (d), and (e) in part~2.

5. We check all properties of an ordered field described in 
Definition~\ref{def_field} (which refers to earlier definitions). In view
of the isomorphism of $\langle\R',\prec\rangle$ and $\langle\R^*,\prec\rangle$,
completeness of the linear order $\langle\R'',\prec\rangle$ was proven 
in Lemma~\ref{lem_SupInfinExtRe}. In order to show that for every 
$\be\in\R''$ we have $0+\be=1\cdot\be=\be$, we just set 
$(a_n):=(0,0,\ds)$, respectively $(a_n):=(1,1,\ds)$, in the definition 
of $\al+\be$, respectively $\al\cdot\be$, and use 
Lemma~\ref{lem_onSubseq}. Commutativity and associativity of $+$ and $\cdot$ on 
$\R''$ follow from the definition and from the fact that they hold in $\Q$. 
The same applies to the distributive law.   

We prove the existence of additive and multiplicative inverses in $\R''$. Let $\al\in\R''$ and $(a_n)\sus\Q$ be such 
that $\lim a_n=\al$. We set $(b_n):=
(-a_n)$. Let $(c_n)=(b_{m_n})$ be a~subsequence of $(b_n)$ with a~limit 
$\be:=\lim c_n$ in $\R'$. It is not hard to show that $\be\in\R''$. By the 
definition of addition,    
$$
\al+\be=\lim\,c_n'=0\,, 
$$
because any subsequence $(c_n')$ of the sequence $(a_{m_n}+c_n)=(0,0,\ds)$ is 
again $(0,0,\ds)$. Let $\al\in\R''$ with $\al\ne0$ and $(a_n)\sus\Q$ be 
such that $\lim a_n=\al$. We may assume that $a_n\ne0$ for every $n$. We set 
$(b_n):=(1/a_n)$ and argue as before, with $+$ replaced by $\cdot$, and 
$(0,0,\ds)$ by $(1,1,\ds)$. We obtain an element $\be\in\R''$ such that 
$\al\cdot\be=1$. We have proven that $\R''$ is a~field.

It remains to prove the two order axioms in $\R''$. Let $\al\prec\be$ and 
$\ga$ be in $\R''$. Using Lemma~\ref{lem_pomoc0} and 
Definition~\ref{def_limLO}, we take sequences $(a_n),(b_n),(c_n)\sus\Q$ and 
fractions $a$ and $b$ such that $\lim a_n=\al$, $\lim b_n=\be$, $\lim 
c_n=\ga$, $\al\prec a\prec b\prec\be$, and that $a_n\preceq a$ and $b\preceq 
b_n$ for every $n$. By Lemma~\ref{lem_narrInt}, there exist 
fractions $c$ and $d$ such that $\ga\in(c,d)$ and $d-c\le
\frac{1}{3}(b-a)$. It follows from the definition of addition in $\R''$ that 
$${\textstyle
\al+\ga\preceq a+\frac{1}{3}(b-a)<b-\frac{1}{3}(b-a)\preceq\be+\ga\,,
}
$$
so that $\al+\ga\prec\be+\ga$. We prove the second order axiom in the 
equivalent form given in Lemma~\ref{lem_euqiOF}. Let 
$\al,\be\in\R''$ be positive elements. We take a~fraction $c>0$ and sequences 
$(a_n),(b_n)\sus\Q$ such that $\lim a_n=\al$, $\lim b_n=\be$, $\al,\be
\succ c$, and $a_n,b_n\ge c$ for every $n$. Then $a_nb_n\ge c^2>0$ for every 
$n$, so that, by the definition of multiplication in $\R''$, we have 
$\al\cdot\be\succeq c^2\succ0$. 
\kduk

\section{Extended reals and ordered fields}

In the last result in this section, we determine which properties of 
a~complete ordered field are preserved or acquired when the 
two infinities are allowed and operations are understood as 
partial. 

\begin{prop}\label{prop_whichPres}
Let $\langle\R',\prec\rangle$ be any limit completion of $\langle\Q,
<\rangle$. The algebraic structure
$$
\R'=\langle\R',\,0,\,1,\,+,\,\cdot,\,
\prec\rangle_{\mathrm{partial}}
$$
retains the properties of a~complete ordered field, with the following 
changes.
\begin{enumerate}
\item The infinities $-\infty'$ and $+\infty'$ have neither an additive nor a~multiplicative inverse. 
\item  In the linear order $\langle\R',\prec\rangle$, every subset has an infimum and a~supremum.
\item Let $\ga=\pm\infty'$. The inequalities $\al\prec\be$ with $\al,\be\ne 
-\ga$ are not preserved under the addition of $\ga$. In all other cases, the first 
order axiom holds if the involved expressions are defined.
\item Let $\ga=+\infty'$. The inequalities $\al\prec\be$ such that $\al,\be\ne0$ and have 
equal signs are not preserved under multiplication by $\ga$. In all other 
cases, the second 
order axiom holds if involved expressions are defined.
\item Every other identity of an ordered field holds in $\R'$, provided that 
involved expressions are defined. This concerns neutral elements, 
commutativity, associativity, and distributivity.
\end{enumerate}
\end{prop}
\duk
1. This is obvious. 2. This was proven in Lemma~\ref{lem_SupInfinExtRe}. 

3. This follows at once from parts~4a and~4b in 
Theorem~\ref{thm_strExtRea}. We show that the first order axiom is not 
violated in any other case. Let $\al\prec\be$ and $\ga$ be in $\R'$, and 
let one of them be $\pm\infty'$. If $\ga\in\R''$, then $\al$ or $\be$ is 
$\pm\infty'$. Since infinity remains unchanged after adding $\ga$, we have 
that $\al+\ga\prec\be+\ga$. If $\ga=\pm\infty'$, then in every case 
different from the stated violating one, the addition of $\ga$ leads to an 
undefined expression.  

4. This follows at once from part~4c in 
Theorem~\ref{thm_strExtRea}. We show that the second order axiom is not 
violated in any other case. Let $\al\prec\be$ and $\ga\succ0$ be in $\R'$, and 
let one of them be $\pm\infty'$. If $\ga\in\R''$, then $\al$ or $\be$ is 
$\pm\infty'$. Since infinity remains unchanged after multiplying it by 
$\ga$, we have that $\al\cdot\ga\prec\be\cdot\ga$. If $\ga=+\infty'$ and $\al\prec0\prec\be$, then 
$\al\cdot\ga=-\infty'\prec+\infty'=\be\cdot\ga$. The 
remaining cases with $\ga=+\infty'$ are either the stated violating case or 
lead to an undefined expression. 

5. It follows from parts~4 and~5 in Theorem~\ref{thm_strExtRea} that 
$0+\al=\al$ and $1\cdot\al=\al$ for every $\al\in\R'$. It also follows from 
these parts that addition and multiplication in $\R'$ are commutative.
We show that they are associative and that the distributive law holds. Let 
$\al$, $\be$, and $\ga$ be in $\R'$, and let one of them be $\pm\infty'$. It 
is clear that the identity
$$
(\al\cdot\be)\cdot\ga=\al\cdot(\be\cdot
\ga)
$$
either holds or neither side is defined. By part~4 in 
Theorem~\ref{thm_strExtRea}, if $\al,\be,\ga\ne0$, then the former case 
occurs, with an infinity being equated to itself. If one of the three elements 
is $0$, then the latter case occurs. We consider the identity
$$
(\al+\be)+\ga=\al+(\be+\ga)\,.
$$
Again, it either holds, or neither side is defined. By part~4 in 
Theorem~\ref{thm_strExtRea}, if the infinities among $\al$, $\be$, and 
$\ga$ have the same sign, then the former case occurs, with an infinity 
being equated to itself. If two infinities with different signs occur, 
then the latter case occurs. Finally, we
consider the identity
$$
\ga\cdot(\al+\be)=\ga\cdot\al+\ga\cdot
\be\,.
$$
Let $\ga=\pm\infty'$. If $\al=0$ or $\be=0$ or $\al+\be=0$ or $\al$ and 
$\be$ have different signs, then one side of the identity is not defined. In 
all other cases, the identity holds, with an infinity being equated to 
itself. Let $\ga\in\R''$. If $\ga=0$ or $\al$ and $\be$ are infinities with 
different signs, then neither side of the identity is defined. In 
all other cases, the identity holds, with an infinity being equated to 
itself.  
\kduk
\chapter[From extended reals to real numbers]{From extended reals to real numbers}\label{sec_RstarToR}

Thanks to the results on extended reals $\R^*$ in Section~\ref{sec_QtoRstar}, this section on real 
numbers $\R$ can be relatively short. We first review the algebraic structure $\R$ and 
then prove the characterization Theorem~\ref{thm_Ris}. 

\section{The real numbers}

\begin{defi}\label{def_realNum}
Recall that the set of real numbers is
$$
\R=\{\al\in\R^*\cc\;-\infty<\al<+\infty\}=\R^*\setminus
\{-\infty,\,+\infty\}\,.
$$
So $\al\in\R$ if and only if $\al$ is a~cut on $\langle\Q,
<\rangle$ different from $\emptyset$ and $\Q$.
\end{defi}

\section{Addition and multiplication}

We recall the definitions of addition and multiplication of real numbers 
given in the previous section.

\begin{defi}\label{def_ariOnR}
Let $\al,\be\in\R$. Then 
$$
\al+\be=\lim_{n\to\infty}c_n\,\text{ and }\,\al\cdot\be=\lim_{n\to\infty}d_n\,,
$$
where $(c_n)$, resp. $(d_n)$, is any subsequence of the sequence of 
fractions $(a_n+b_n)$, resp. $(a_n\cdot b_n)$, that has a~limit in $\langle\R^*,
<\rangle$. The sequence $(a_n)$, resp. $(b_n)$, is any sequence of fractions 
such that $\lim a_n=\al$, resp. $\lim b_n=\be$.
\end{defi}
The last $+$ and $\cdot$ take place in the ordered field $\Q$. The limits are 
taken in the sense of Definition~\ref{def_limLO}. The 
sequences $(a_n)$ and $(b_n)$ exist by Corollary~\ref{cor_XdenseY} or by 
Lemma~\ref{lem_ratCutDens}. For the same reason, that $\langle\R^*,<\rangle$
is a~limit completion of $\langle\Q,<\rangle$, there exist the subsequences 
$(c_n)$ and $(d_n)$. Soon we simplify this definition. $\R$ also denotes the 
following algebraic structure 

\begin{defi}\label{def_RasStru}
The algebraic structure
$$
\R=\langle
\R,\,0_{\R},\,1_{\R},\,+,\,\cdot,\,<
\rangle
$$
is formed by the set $\R$ in Definition~\ref{def_realNum}, by the 
distinct elements
$$
0_{\R}=\{a\in\Q\cc\;a<0\}\,\text{ and }\,1_{\R}=\{a\in\Q\cc\;a<1\}\ \ (\in\R)\,,
$$
by the operations $+$ and $\cdot$ in Definition~\ref{def_ariOnR}, and by the 
linear order $<$, which is just the set inclusion of cuts. 
\end{defi}
As usual, instead of $0_{\R}$ and $1_{\R}$ we write just $0$ and $1$. The 
algebraic structure that fits the real numbers is the structure of an ordered 
field, see Definition~\ref{def_field}. Recall that an ordered field is complete
if every nonempty and upper-bounded set in it has a~supremum. 

\section{The characterization theorem}

Now we can state the well known characterization theorem. 

\begin{thm}\label{thm_Ris}
The structure $\R$ is a~complete ordered field. Every two complete 
ordered fields are isomorphic as ordered fields.    
\end{thm}
The former claim was proven as part~5 of Theorem~\ref{thm_strExtRea}. We prove the latter isomorphism claim in Proposition~\ref{prop_isomR}. 

\section{$\R$ is unique}

We need five lemmas.

\begin{lemma}\label{lem_clOFisAr}
Let $F$ be a~complete ordered field. Then for every element $a\in F$ there 
exists a~sum of ones $s=1_F+1_F+\ds+1_F$ such that
$$
-s\le a\le s\,.
$$
\end{lemma}
\duk
It suffices to prove just the latter inequality. Suppose for the contrary 
that the set $S$ ($\sus F$) of finite sums of $1_F$'s is bounded from above. Ten we can 
set
$$
b:=\sup(S)\ \ (\in F)\,.
$$
By the definition of supremum there exists an element $s\in S$ such that $b-
1_F<s$. We get the contradiction that $b<s+1_F\in S$.
\kduk

The following fact is well known.

\begin{lemma}\label{lem_wellKnF}
Let $F$ be a~complete ordered field and $(a_n)\sus F$ be a~Cauchy sequence. 
Then $(a_n)$ has a~limit in $F$.
\end{lemma}

Let $F$ be an ordered field. Due to the field embedding $f_F\cc\Q\to F$, we can 
view $\Q$ as a~subfield of $F$.

\begin{lemma}\label{lem_ariInR}
Let $F$ be a~complete ordered field and $(a_n),(b_n)\sus\Q$ be sequences that 
have limits in $F$. Then the sequences
$$
(a_n+b_n)\,\text{ and }\,(a_nb_n)
$$
have limits in $F$ as well.
\end{lemma}
\duk
It follows that $(a_n)$ and $(b_n)$ are Cauchy. It is not hard to show that $(a_n+b_n)$ and $(a_nb_n)$ are also 
Cauchy. The result follows from Lemma~\ref{lem_wellKnF}.
\kduk

\noindent
Thus we simplify Definition~\ref{def_ariOnR} as follows.
If $\al,\be\in F$ are two elements in a~complete ordered field and $(a_n),
(b_n)\sus\Q$ are such that $\lim a_n=\al$ and $\lim b_n=\be$, then $\al+\be=\lim(a_n+b_n)$ and $\al\cdot\be=\lim(a_n\cdot b_n)$.

The next lemma is a~generalization of the standard lemma on the arithmetic of 
limits of real sequences, and its proof can be left to the reader.

\begin{lemma}\label{lem_ariLim}
Let $F$ be an ordered field, $a,b\in F$, and let $(a_n),(b_n)\sus F$ be sequences such that $\lim a_n=a$ and $\lim b_n=b$. Then
$$
\lim\,(a_n+b_n)=a+b\,\text{ and }\,
\lim\,(a_n\cdot b_n)=a\cdot b\,.
$$
\end{lemma}

The next lemma establishes that Cauchy sequence are preserved by the 
field embeddings $f_F$.

\begin{lemma}\label{lem_isCauchy}
Let $F=\langle F,0_F,1_F,\oplus,\odot,\prec\rangle$    
be an ordered field, $g=f_F\cc\Q\to F$ be the field embedding from 
Definition~\ref{def_mapsfF}, and $(a_n)\sus\Q$ be a~Cauchy sequence. 
Then $(g(a_n))\sus F$ is a~Cauchy sequence in $F$.
\end{lemma}
\duk
It is important that by Lemma~\ref{lem_tretifF} the embedding 
$g$ respects also the linear orders $<$ and $\prec$; it goes back to Lemma~\ref{lem_lemaR3}. Let an $e\in F$ with 
$e\succ 0_F$ be given. We set $e':=g^{-1}(e)$ ($\in\Q$ and $>0$) and 
take an $n_0\in\N$ such that $m,n\ge n_0$ $\Rightarrow$ $|a_m-a_n|\le e'$, 
i.e
$$
-e'\le a_m-a_n\le e'\,.
$$
We apply $g$ to these two inequalities and get that
$$
-e\preceq g(a_m)-g(a_n)\preceq e\,,
$$
i.e. $|g(a_m)-g(a_n)|\preceq e$ for the same $m$ and $n$.
\kduk

Now we can prove the uniqueness of $\R$.

\begin{prop}\label{prop_isomR}
Every two complete ordered fields are isomorphic as ordered fields.     
\end{prop}
\duk
It suffices to show that any complete ordered field
$$
F=\langle F,\,0_F,\,1_F,\,\oplus,\,
\odot,\,\prec\rangle
$$
is isomorphic to $\R$ as an ordered field. We may assume that $\Q$ is 
a~subfield of $\R$, and we take the field embedding $g=f_F\cc\Q\to F$ in 
Definition~\ref{def_mapsfF}. We 
define a~map $\Phi\cc\R\to F$ and show that it is an isomorphism of ordered 
fields. Let $\al\in\R$. By Lemma~\ref{lem_ratCutDens} or 
Corollary~\ref{cor_XdenseY}, we can take a~sequence $(a_n)\sus\Q$ such that 
$\lim a_n=\al$. (By Proposition~\ref{prop_equiDefi}, we may 
use two equivalent definitions of limits.) It follows that $(a_n)$ is 
Cauchy, and by Lemma~\ref{lem_isCauchy} the sequence $(g(a_n))\sus F$ is also Cauchy. By Lemma~\ref{lem_wellKnF} we
can set $\Phi(\al):=\lim_F g(a_n)$. We first show that the value $\Phi(\al)$ 
does not depend on the choice of $(a_n)$. If $(b_n)\sus\Q$ is another 
sequence with $\lim b_n=\al$ and $\lim_F g(b_n)\ne \lim_F g(a_n)$, we 
consider the combined rational sequence
$$
(c_n):=(a_1,\,b_1,\,a_2,\,b_2,\,a_3,\,\ds)\,.
$$
Again, $\lim c_n=\al$. It follows that the sequence $(g(c_n))$ has two 
different limits, which is a~contradiction (with Corollary~\ref{cor_HausdTop}).

We show that $\Phi$ is an isomorphism of ordered fields. Let $\al,\be\in\R$ 
and let $(a_n),(b_n)\sus\Q$ be such that $\lim a_n=\al$ and $\lim b_n=\be$.
We compute
\begin{eqnarray*}
\Phi(\al+\be)&=&\lim_F g(a_n+b_n)=
\lim_F (g(a_n)+g(b_n))\\
&=&\lim_F g(a_n)\oplus\lim_F g(b_n)=\Phi(\al)\oplus\Phi(\be)\,.
\end{eqnarray*}
The first equality follows from Lemma~\ref{lem_ariInR} and the 
definition of $\Phi$. In the second equality, we use that $g$ is a~field 
embedding. In the third equality, we use Lemma~\ref{lem_ariLim}. The last 
fourth equality follows from  the 
definition of $\Phi$. We similarly compute
\begin{eqnarray*}
\Phi(\al\cdot\be)&=&\lim_F g(a_n\cdot b_n)=\lim_F (g(a_n)\cdot g(b_n))\\
&=&\lim_F g(a_n)\odot\lim_F g(b_n)=\Phi(\al)\odot\Phi(\be)\,,
\end{eqnarray*}
with similar justifications. We assume that $\al<\be$ and deduce that $\Phi(\al)\prec\Phi(\be)$. By Lemma~\ref{lem_odmvR}, 
$\al+\ga^2=\be$ for some $\ga>0$. We have proven that $\Phi$ is a~field 
embedding, and therefore
$$
\Phi(\al)\oplus\Phi(\ga)^2=\Phi(\be)\,.
$$
Since $\Phi(\ga)^2\succ0_F$, we deduce that $\Phi(\al)\prec\Phi(\be)$. In 
particular, $\Phi$ is injective. It remains to show that $\Phi$ is onto. 
Let $y\in F$ be an arbitrary element. Since $F$ is Archimedean by 
Lemma~\ref{lem_clOFisAr}, by halving $g$-rational intervals containing $y$ 
we can obtain a~sequence $(y_n)\sus g[\Q]$ such that $\lim_F y_n=y$. Then 
$(y_n)$ is Cauchy, and so is $(x_n):=(g^{-1}(y_n))\sus\Q$. Using 
Lemma~\ref{lem_wellKnF}, we set $\ga:=\lim x_n$ ($\in\R$). It follows 
that $\Phi(\ga)=y$.
\kduk

\noindent
This completes the proof of Theorem~\ref{thm_Ris}.

We leave the proof of the next useful lemma to the interested reader as 
an exercise.

\begin{lemma}\label{lem_odmvR}
For every real number $a>0$, there exists a~unique real number $b>0$ 
such that $b^2=b\cdot b=a$.     
\end{lemma}
We write $\sqrt{a}:=b$.

\section{$\R$ is uncountable}
\chapter[Limits of real sequences]{Limits of real sequences}\label{chap_existLim}

After the lengthy but inevitable, and hopefully not uninteresting 
preparation in the six previous chapters, we get to the first 
topic of elementary real analysis proper. It is the limit behavior of 
sequences of real numbers.

Section~\ref{podkap_nekoOkolLimi} introduces limits of real sequences 
in the domain $\R^*$ of extended reals. Theorem~\ref{thm_idenDiff} 
establishes two identities for subtraction in $\R^*$. In 
Theorem~\ref{thm_deleVRhvezd} we prove another two identities in 
$\R^*$, now for addition, multiplication, and division.
Definition~\ref{def_vlLimi} introduces limits of real sequences
in $\R^*$ classically in terms of neighborhoods (in 
Chapter~\ref{sec_QtoRstar} 
we have already considered 
limits in the linear order $\langle\R^*,<\rangle$). Definition~\ref{def_robust} 
introduces the notion of robustness of properties of real sequences. In 
Proposition~\ref{prop_Lim_nta_odm_zn} we prove that $n^{1/n}\to1$.

Section~\ref{sec_subseq} is devoted to subsequences.

\section[Infinities, neighborhoods, limits]{Infinities, neighborhoods, limits}\label{podkap_nekoOkolLimi}

We recall the algebraic structure of extended reals that was introduced in 
Then we give already the third definition of limits of real sequences. 

\medskip\noindent
{\em $\bullet$ Notation. }For logical and set-theoretic notation, see 
Appendix~\ref{sec_logST}. Letters $i$, $j$, $k$, $l$, $m$, $m_0$, $m_1$, $\ds$, and $n$, $n_0$, 
$n_1$, $\ds$, possibly with primes, are variables ranging in $\N=\{1,2,\ds\}$. 
By $a$, $b$, $c$, $d$, $e$, $\de$, $\ep$ and $\theta$, possibly with 
indices and primes, we denote 
real numbers. \underline{Always $\de, \ep,\theta>0$}. Extended reals, 
elements in $\R^*=\R\cup\{-\infty,+\infty\}$, are denoted by the letters $A$, $B$, $C$, $D$, $K$, and $L$. Recall that a~sequence 
$(a_n)\sus\R$ is a~function $a\cc\N\to\R$ with $a_n=a(n)$. Sets of real numbers are denoted by $M$ and $N$. Recall that the 
absolute value of $a\in\R$ is $|a|=a$ if $a\ge0$, and $|a|=-a$ if $a<0$. 

\begin{exer}\label{ex_trojNero}
Recall the proof (in ordered fields) of the triangle inequality, abbreviated {\em TI}, that for every $a,b\in\R$,
$$
|a+b|\le|a|+|b|\,.
$$
\end{exer}
Theorem~\ref{thm_AKjeKonv} is
an infinite TI. In 
give two proofs of the TI in $\C$.

\medskip\noindent
{\em $\bullet$ Review of the operations on the extended reals introduced in 
Recall that 
$$
\R^*=\R\cup\{+\infty,\,-\infty\}
$$
is the set of extended reals.
We denote elements in $\R^*$ by $A$, $B$, $C$, and $D$. Recall that the 
linear order $\langle\R^*,<\rangle$ extends $\langle\R,<\rangle$ by the 
comparisons $-\infty<+\infty$, 
$-\infty<a$, and $a<+\infty$ for every $a\in\R$. The addition $+$ in $\R^*$ 
extends the addition in $\R$ by the values $a+A=A+a=A+A=A$ for $a\in\R$ and 
$A=\pm\infty$. The multiplication is extended by the values $a\cdot A=A\cdot 
a=A\cdot A=A$ for real $a>0$ and $A=\pm\infty$, and by $a\cdot A=A\cdot 
a=(-A)\cdot A=-A$ for real $a<0$. Here 
$-(\pm\infty):=\mp\infty$ (all top or all bottom signs). We have the new 
differences $a-A=-A$, $A-a=A$, and $A-(-A)=A$ for every $a\in\R$ and $A=\pm\infty$. We have the new ratios
$A/a=A$ for real $a>0$ and $A=\pm\infty$, $A/a=-A$ for real $a<0$, 
and $a/A=0$ for any $a\in\R$. The remaining undefined expressions, called indefinite expressions, are
$$
{\textstyle
\pm\infty+(\mp\infty),\ 
0\cdot(\pm\infty),\ 
(\pm\infty)\cdot 0,\ 
(\pm\infty)/(\pm\infty),\,\text{ and }\,A/0\,\text{ for }\, A\in\R^*\,,
}
$$
with all top or all bottom signs in the first expression, and all four sign combinations in the fourth one. 

\begin{exer}\label{ex_pocSnek}
Compute $(-\infty)/(-2)$, $(-\infty)-(+\infty)$, $-\infty+10$, and $(+\infty)/0$.
\end{exer}

\noindent
{\em $\bullet$ Review of the algebraic structure $\R^*$ introduced in 
this section for more details. Recall that the original structure
$$
\R=\langle
\R,\,0,\,1,\,+,\,\cdot,\,<
\rangle
$$
is a~complete ordered field 
In the extension 
$$
\R^*=\langle
\R^*,\,0,\,1,\,+,\,\cdot,\,<
\rangle\,, 
$$
the elements $0$ and $1$ remain neutral to $+$ and $\cdot$, respectively, $+$ 
and $\cdot$ remain commutative and associative, and $\cdot$ remains 
distributive to $+$, if the involved expressions are not indefinite. The two 
infinities $\pm\infty$ have neither additive nor multiplicative inverses, 
but the other inverses are as in $\R$. With the exception of cases described 
in parts~10 and~11 in 
the two 
order axioms remain valid in $\R^*$.

\vspace{-3mm}
\begin{exer}\label{ex_supMinusNek}
Find all sets $X\sus\R^*$ such that $\sup(X)=-\infty$. 
\end{exer}

In the next two theorems, we show that some familiar identities for 
differences and ratios remain valid in $\R^*$.

\begin{thm}\label{thm_idenDiff}
For every $A,B,C\in\R^*$ we have
$$
A-(B-C)=(A-B)+C\,\text{ and }\,
(A-B)-C=(A-C)-B\,,
$$    
if the involved expressions are not indefinite.
\end{thm}
\duk
Let $A$ or $B$ or $C$ be infinity, and let all involved expressions be 
defined. We prove the former identity. Let $A=\pm\infty$. If $B-C\in\R$, then 
$B,C\in\R$ and the identity holds as the equality $A=A$. Let $B-C=-A$. If 
$B\in\R$, then $C=A$ and we get $A=A$. If $B=-A$, then $C\in\R$ or $C=A$, and 
we get $A=A$. Let $A\in\R$. If $B=\pm\infty$, then $C\in\R$ or $C=-B$, 
and we get $-B=-B$. If $B\in\R$, then $C=\pm\infty$, and we get $C=C$. This 
exhausts all cases when all involved expressions are defined and infinity is 
present, and the identity holds in each of them. 

We prove the latter identity. Let $A-B\in\R$. Then $A,B\in\R$, 
$C=\pm\infty$, and we get $-C=-C$. Let 
$A-B=\pm\infty$. If $A\in\R$, then $B=\pm\infty$, $C\in\R$ or $C=B$, and 
we get $-B=-B$. If $A=\pm\infty$, then $B\in\R$ or $B=-A$, $C\in\R$ or $C=-A$, 
and we get $A=A$. We are done in the same way as for the former identity.
\kduk

\begin{thm}\label{thm_deleVRhvezd}
For every $A,B,C,D\in\R^*$ we have
$$
A/C+B/C=(A+B)/C\,\text{ and }\,
(A/C)\cdot(B/D)=(A\cdot B)/(C\cdot D)\,,
$$
if the involved expressions are not indefinite.
\end{thm}
\duk
Let $A$ or $B$ or $C$ or $D$ be infinity, and let all involved 
expressions be defined. We prove the former identity. Let $A=\pm\infty$. 
Then $C\in\R\setminus\{0\}$, $B\in\R$ or $B=A$, and we get $A=A$ or $-A=-A$. 
The case $B=\pm\infty$ is symmetric. We are left with the case $A,B\in\R$ and 
$C=\pm\infty$. Then we get $0+0=0$ and are done. 

We prove the latter identity. Let $A=\pm\infty$. Then $C\in\R
\setminus\{0\}$. If $B=\pm\infty$, then
$D\in\R\setminus\{0\}$, and we get $\pm\infty=\pm\infty$, with the sign 
equal to the product of the signs of the four elements $A$, $B$, $C$, and $D$.
If $B\in\R$, then $B,D\ne0$, $D\in\R$, and we again get  
$\pm\infty=\pm\infty$, with the sign 
equal to the product of the four signs. The case $B=\pm\infty$ is symmetric. We 
are left with the case $A,B\in\R$. Since $C,D\ne0$ and one of $C$ and $D$ 
is infinity, we get the valid equality
$a\cdot b=0$, where $a,b\in\R$ and one of them is zero. We are done.
\kduk

\begin{exer}\label{ex_oneMore}
Show that even the more general identity
$$
A/C+B/D=(A\cdot D+C\cdot B)/(C\cdot D)
$$
is valid in $\R^*$.
\end{exer}

\noindent
{\em $\bullet$ Neighborhoods of elements in $\R^*$. }Let $b\in\R$ and 
let $\ep>0$ be given. We define certain intervals in the linear order $\R^*$: 
$$
U(b,\,\ep):=(b-\ep,\,b+\ep)\label{neighb} 
$$
is the
\underline{$\ep$-neighborhood\index{neighborhood!of 
a~point|emph}} of $b$, and 
$${\textstyle
U(-\infty,\,\ep):=[-\infty,\,-\frac{1}{\ep})\;\text{ and }\;U(+\infty,\,\ep):=(\frac{1}{\ep},\,+\infty]\label{neighbInf}}
$$
are \underline{$\ep$-neighborhoods} of 
infinities\index{neighborhood!of an infinity|emph}. Note that 
$$
-\infty\in U(-\infty,\ep)\,\text{ and }\,+\infty\in 
U(+\infty,\ep)\,. 
$$
Basic properties of neighborhoods are reviewed in five exercises.
\index{neighborhood!properties of}

\begin{exer}\label{ex_jeToOpIn}
Every neighborhood $U(L,\ep)$ is an open interval in the linear order $\langle\R^*,<\rangle$.     
\end{exer}

\begin{exer}\label{ex_nekdVlOkoli}
Let $A\in\R^*$, $c\in U(A,\ep)$ and $c<b<A+\ep$ or $A-\ep< b<c$. Then $b\in U(A,\ep)$.
\end{exer}
If $M,N\sus\R^*$, then $M<N$ means that $A<B$ for every $A\in M$ and $B\in N$.
\begin{exer}\label{ex_ulohaNaokoli1}
Let $A<B$ be in $\R^*$. Then there is an $\ep$ such that $U(A,\ep)
<U(B,\ep)$. In particular, $U(A,\ep)\cap U(B,\ep)=\emptyset$ for 
some $\ep$. 
\end{exer}

\begin{exer}\label{ex_ulohaNaokol2}
Let $A\in\R^*$. 
If $\ep\le\de$ then $U(A,\ep)\sus U(A,\de)$.
\end{exer}

\begin{exer}\label{ex_ulohaNaokoli3}
Let $A\in\R^*$. Then
$\bigcap_{k=1}^{\infty}U(A,\,
\frac{1}{k})=\{A\}$. 
\end{exer}

\noindent
{\em $\bullet$ Limits of real sequences. }We denote sequences of real 
numbers by $(a_n)$, $(b_n)$, 
and $(c_n)$. We call them real 
sequences or just sequences. We restate
the specialization of Definition~\ref{def_limLO} of limits in linear orders to the 
linear order $\langle\R^*,<\rangle$.

\begin{defi}\label{def_vlLimi}
Let $(a_n)\sus\R$ and $L\in\R^*$. If for every $\ep>0$ there is an $n_0$ 
such that for every $n\ge n_0$ we have 
$$
a_n\in U(L,\,\ep)\,,
$$
we write that $\lim\,a_n=L$\label{lim} or $\lim_{n\to\infty} a_n=L$ or 
$a_n\to L$, and say that the sequence $(a_n)$ has the \underline{limit\index{limit of a~sequence|emph}} $L$.
\end{defi}

\noindent
If $L\in\R$, we say that $(a_n)$ has a~\underline{finite\index{limit of a~sequence!finite|emph} limit} or that $(a_n)$ \underline{converges\index{sequence!convergent|emph}}. If 
$L=\pm\infty$,  we say that $(a_n)$ has an \underline{infinite\index{limit of 
a~sequence!infinite|emph} limit}. A~sequence \underline{diverges\index{sequence!divergent|emph}} if it has no limit or an infinite limit. An \underline{eventually constant\index{sequence!eventually 
constant|emph}} sequence $(a_n)$ with $a_n=a$ 
for $n\ge n_0$ converges and $\lim 
a_n=a$. We show that the present definition of limits is equivalent to 
the earlier one. It is a~corollary of a~universality property of the 
neighborhoods $U(L,\ep)$.

\begin{prop}\label{prop_eq2defsL}
Let $L\in\R^*$ and $(a_n)\sus\R$. If
$$
\lim_{n\to\infty}a_n=L
$$
by the specialization of Definition~\ref{def_limLO} to the linear order 
$\langle\R^*,<\rangle$,  \underline{then} this limit also holds 
by Definition~\ref{def_vlLimi}, and vice
versa. 
\end{prop}
\duk
Suppose that $\lim a_n=L$ in the linear order $\R^*=\langle\R^*,<\rangle$ by 
Definition~\ref{def_limLO}, and that $\ep$ is given. By Exercise~\ref{ex_jeToOpIn} and the assumption, there is an $n_0$ such that 
$a_n\in U(L,\ep)$ for every $n\ge n_0$.

Suppose that $\lim a_n=L$ by Definition~\ref{def_vlLimi}, and that 
an open interval $I\ni L$ in the linear order $\R^*$ is given. By 
Definition~\ref{def_openInt} of open intervals in linear orders, there exists $\ep$ such that $U(L,\ep)\sus 
I$. Then there is an $n_0$ such that 
$a_n\in U(L,\ep)\sus I$ for every $n\ge n_0$.
\kduk

\begin{exer}\label{ex_blowUp}
Suppose that $(a_n)\sus\R$, $(m_n)\sus\N$ and that the
sequence $(b_n)$ arises from $(a_n)$ by replacing each term $a_n$ 
by $m_n$ copies $a_n,a_n,\ds,a_n$. Then
$\lim a_n=\lim b_n$
whenever one limit exists. 
\end{exer}

\begin{prop}\label{prop_limMinNek}
$\lim a_n=-\infty$ $\iff$ for every $c<0$ there exists $n_0$ such that if $n\ge n_0$, then $a_n\le c$. 
\end{prop}
\duk
For the implication $\Rightarrow$ it suffices to take $\ep>0$ such that 
$-1/\ep\le c$. For the implication $\Leftarrow$ it suffices to take $c<0$ such that $c<-1/\ep$. 
\kduk

\begin{exer}\label{ex_limPlInf}
State and prove the analogous equivalence for the limit $+\infty$.  
\end{exer}

We know from 
that limits of sequences are unique for any 
linear order. For reader's convenience, we prove it again in the instance of 
$\R^*$. See also 

\begin{prop}\label{prop_jednoznLim}
Let $K,L\in\R^*$,\underline{\index{limit of a~sequence!uniqueness of|emph}} 
$\lim a_n=K$ and $\lim a_n=L$. \underline{Then} $K=L$.
\end{prop}
\duk
Let $K$ and $L$ be limits of $(a_n)$ and $\ep$ be arbitrary. By Definition~\ref{def_vlLimi},  
if $n\ge n_0$ then $a_n\in U(K,\ep)$ and $a_n\in 
U(L,\ep)$. Thus $U(K,\ep)\cap 
U(L,\ep)\ne\emptyset$ for every $\ep$. By Exercise~\ref{ex_ulohaNaokoli1},  
$K=L$.
\kduk

\noindent
{\em $\bullet$ Robust properties of sequences. }Let $\R^{\N}$\label{RtoN} be the set of real 
sequences. Any set $V\sus \R^{\N}$ 
is called a~\underline{property\index{property of sequences|emph}} of real 
sequences.

\begin{defi}\label{def_robust}
A~property $V\sus \R^{\N}$ is 
\underline{robust\index{property of sequences!robust|emph}} if for every two sequences $(a_n)$ and $(b_n)$
such that $a_n\ne b_n$ for only finitely many indices $n$, 
we have $(a_n)\in V\iff(b_n)\in V$.
\end{defi}

\begin{exer}\label{ex_jenomCvic}
For every $L\in\R^*$, the property $\{(a_n)\cc\;\lim a_n=L\}$ is robust.
\end{exer}

The next exercise shows that the notion of robustness is itself robust. 

\begin{exer}\label{ex_robRob}
Prove the following.
\begin{enumerate}
\item If $V\sus\R^{\N}$ is robust, then so is $\R^{\N}\setminus V$.
\item If $X\sus\mathcal{P}(\R^{\N})$ is such that every 
$Y\in X$ is robust, 
then $\bigcup X$ is robust.
\item If $\emptyset\ne X\sus\mathcal{P}(\R^{\N})$ is such that every 
$Y\in X$ is robust, then 
$\bigcap X$ is robust.
\end{enumerate}
\end{exer}

\begin{exer}\label{ex_zmenaKonmnoha}
Which of the following properties of sequences $(a_n)$ is robust? 
\begin{enumerate}
\item $(a_n)$ converges.
\item $a_1\le a_2\le\cdots$.
\item There exists an index $m$ such that $a_m\ge a_{m+1}\ge\cdots$.
\item There exist indices $m_1<m_2<\ds$ such that $a_{m_1}>a_{m_2}>\cdots$.
\item $\inf(\{a_n\cc\;n\in\N\})=-1$.
\item $\inf(\{a_n\cc\;n\in\N\})=-\infty$.
\end{enumerate} 
\end{exer}

\noindent
{\em $\bullet$ Three interesting limits. }For $a\in\R$, we define the \underline{upper integer part\index{integer part of 
a~number!upper|emph}} $\lceil a\rceil$ ($\in\Z$) of $a$ as the 
smallest integer $v$ such that $v\ge a$. Similarly, the lower 
\underline{integer part\index{integer part of a~number!lower|emph}} 
$\lfloor a\rfloor$  ($\in\Z$) of $a$ is the largest $v\in\Z$ such that 
$v\le a$. Many\underline{\index{limit of a~sequence!one over en@$\frac{1}{n}\to0$|emph}} limits can be
reduced to the next limit. 

\begin{prop}\label{prop_1overN} 
$\lim\frac{1}{n}=0$.     
\end{prop}
\duk
Let an $\ep>0$ be given and
$n_0=\lceil
\frac{1}{\ep}\rceil+1$. Then for every $n\ge n_0$, 
$${\textstyle
0<\frac{1}{n}\le\frac{1}{n_0}=\frac{1}{\lceil 1/\ep\rceil+1}<
\frac{1}{1/\ep}=\ep\,. 
}
$$
Thus if $n\ge n_0$ then
$\frac{1}{n}\in U(0,\ep)$ and
$\lim\frac{1}{n}=0$.
\kduk
\vspace{-3mm}
\begin{exer}\label{ex_optVal}
Show that $\lceil
\frac{1}{\ep}\rceil+1$ is, for every $\ep>0$, the minimum value of $n_0\in\N$ such that $\frac{1}{n}\in U(0,\ep)$ for every $n\ge n_0$.     
\end{exer}

We compute the next\underline{\index{limit of 
a~sequence!third root@$\sqrt[3]{n}-\sqrt{n}\to-\infty$|emph}} limit via an algebraic transformation.

\begin{prop}\label{prop_2ndExa}
$\lim (\sqrt[3]{n}-\sqrt{n})=-\infty$.  
\end{prop}
\duk
We use Proposition~\ref{prop_limMinNek}. 
Let a~$c<0$ be given. We take any $n_0\ge \max(4c^2,2^6)$. Then for every
$n\ge n_0$ we have
$$
\overbrace{\sqrt[3]{n}-\sqrt{n}\,}^{\text{nontrivial}}=\overbrace{n^{1/2}\cdot\underbrace{(n^{-1/6}-1)}_{\text{ $\ds\le-1/2$}}}^{\text{trivial}}\le\underbrace{-n^{1/2}}_{\text{$\ds\le-2|c|$}}/\,2\le-2|c|/2=c\;,
$$
and the limit is $-\infty$. 
The first upper bracket says that in this form the limit is 
non-trivial, in the indefinite form 
$+\infty-(+\infty)$.
We transform it algebraically in the trivial form  
$(+\infty)\cdot(0-1)=(+\infty)(-1)=-\infty$. Lower brackets show upper bounds for enclosed expressions for $n\ge n_0$.
\kduk

\begin{exer}\label{ex_priklNalimitu}
Find the limit
$\lim_{n\to\infty}\frac{\sqrt[3]{n}-\sqrt{n}}{\sqrt[4]{n}}$.
\end{exer}

Limit of a~sequence is 
\underline{nontrivial\index{limit of a~sequence!nontrivial}} if in the original form it leads to an indefinite expression. Else it is \underline{trivial\index{limit of a~sequence!trivial}}. 
For example, the limits $\lim (2^n+3^n)$ and $\lim\,\frac{4}{5n-3}$
are trivial, but $\lim (2^n-3^n)=+\infty-(+\infty)$ and 
$\lim\frac{4n+7}{5n-3}=\frac{+\infty}{+\infty}$ 
are nontrivial. Non-trivial limits can often, but not always, be 
computed by algebraic transformations into trivial limits, as in Proposition~\ref{prop_2ndExa}. 
The next limit of $n^{1/n}$ is nontrivial because $n\to+\infty$,  
$\frac{1}{n}\to0$ and $(+\infty)^0$ is an indefinite power 
expression (see Exercise~\ref{ex_laterEx}). No 
algebraic transformation works in this case, and we have to compute the 
limit from the definition. We show that the exponent prevails and 
$n^{1/n}\to1$. 
We use the well known 
\underline{binomial theorem\index{binomial theorem|emph}} 
that is reviewed in the next exercise.

\begin{exer}\label{ex_binomVeta}
Let $a,b\in\R$ and $n\in\N_0$. Then 
$${\textstyle
(a+b)^n=\sum_{j=0}^n\binom{n}{j}a^j b^{n-j}\,.
}
$$
Here $\binom{n}{j}:=\frac{1}{j!}n(n-1)\ds(n-j+1)$ for $j\in\N$ 
and $\binom{n}{0}:=1$.\label{binom}   
\end{exer}

We assume that the reader is familiar with the real powers $a^b$ with base 
$a>0$ and exponent $b\in\Q$. We introduce them in Section~\ref{sec_elemenFce}.

\begin{prop}\label{prop_Lim_nta_odm_zn}
\underline{\index{limit of a~sequence!nth root@$\sqrt[n]{n}\to1$|emph}} 
$\lim_{n\to\infty}\,n^{1/n}=1$.
\end{prop}
\duk
Always $n^{1/n}\ge1$. If $n^{1/n}\not\to 1$, there would be 
a~number $c>0$ and a~sequence of integers $2\le n_1<n_2<\ds$  such 
that $n_i^{1/n_i}\ge1+c$  for every $i$ (Exercise~\ref{ex_vysvPodp}). 
Raising this inequality to the power $n_i$ (Exercise~\ref{ex_explIndet}) and using 
Exercise~\ref{ex_binomVeta} we get
\begin{eqnarray*}
n_i&\ge&{\textstyle (1+c)^{n_i}=\sum_{j=0}^{n_i}\binom{n_i}{j}c^j=1+\binom{n_i}{1}c+\binom{n_i}{2}c^2+\ds+\binom{n_i}{n_i}c^{n_i}}\\
&\ge&{\textstyle\binom{n_i}{2}c^2=\frac{1}{2}n_i(n_i-1)\cdot c^2\,.
}
\end{eqnarray*}
So for every $i$, 
$${\textstyle
n_i\ge\frac{1}{2}n_i(n_i-1)\cdot c^2\,\text{ and }\,
1+2/c^2\ge n_i\,.
}
$$
This is impossible, the sequence $n_1<n_2<\ds$ is not bounded from above.
\kduk
\vspace{-3mm}
\begin{exer}\label{ex_vysvPodp}
Explain why there is the sequence $(n_i)$.
\end{exer}

\begin{exer}\label{ex_explIndet}
Explain why $n_i^{1/n_i}\ge1+c$  $\Rightarrow$ $n_i\ge(1+c)^{n_i}$.   
\end{exer}

\noindent
{\em $\bullet$ Limits of monotone sequences. }In the next 
passage, we determine the limits of sequences $(q^n)$ with $q\in\R$. We 
need the existence of 
limits of monotone sequences for this.  
proved it for finite limits. 
We extend it to infinite limits. 

\begin{prop}\label{prop_O_mon1}
Let 
$a_1\le a_2\le\ds$ be real numbers. \underline{Then}    
$$
\lim_{n\to\infty} a_n=\sup(\{a_n\cc\;n\in\N\})\ \ (\in\R^*)\,.
$$
The same holds when $\le$ is replaced with $\ge$, and $\sup$ with 
$\inf$. Infima and suprema are taken in the linear order $\langle\R^*,<\rangle$.
\end{prop}
\duk
We may assume that the supremum is $+\infty$. This means that for any 
given $\ep>0$ there is an $n_0$ such that $a_{n_0}>1/\ep$. By the 
monotonicity of $(a_n)$ we have $a_n\in U(+\infty,\ep)$ for every $n\ge n_0$.
Thus $\lim a_n=+\infty$.
\kduk

\begin{exer}\label{ex_O_mono}
Is the assumption on $(a_n)$ a~robust property of sequences?    
\end{exer}

\begin{exer}\label{ex_O_mono1}
Reduce the $\ge$ case to the $\le$ case.   
\end{exer}

\noindent
{\em $\bullet$ Limits of geometric sequences. }A~\underline{geometric sequence\index{sequence!geometric|emph}} is any sequence of the form
$$
(q^n)=(q,q^2,q^3,\ds),\ q\in\R\,.
$$

\begin{exer}\label{ex_jesteDodskAL}
Let $(a_n)\sus\R$. Then 
$\lim a_n=0$ $\iff$ $\lim |a_n|=0$. 
\end{exer}

\begin{prop}\label{prop_geoPosl}
Let $q\in\R$. The limit $L=\lim q^n$ has the following values.
\begin{enumerate}
\item If $|q|<1$ then $L=0$.
\item If $q=1$ then $L=1$.
\item If $q\le -1$ then $L$ does not exist.
\item If $q>1$ then $L=+\infty$. 
\end{enumerate}
\end{prop}
\duk
1. By Exercise~\ref{ex_jesteDodskAL}, we can assume that 
$q\in(0,1)$. We show that $\lim q^n=0$. Let $L=
\inf(\{q^n\cc\;n\in\N\})$ ($\in[0,1]$). Since 
$q>q^2>q^3>\ds>0$, $\lim q^n=L$ by Proposition~\ref{prop_O_mon1}. It remains to show that 
$L=0$. Let $L>0$. Then 
$L/q>q^n\ge L$ for some $n\in\N$. We get the contradiction $L>q^{n+1}$.

2. The constant sequence $(1,1,\ds)$ has limit $1$. 

3. For 
$q\le-1$, the sequence $(q^n)$ does not have a~limit: 
if $m,n\in\N$ have different parity, then $|q^m-q^n|\ge2$.

4. Let $q>1$. We proceed as in part~1. Let 
$L=\sup(\{q^n\cc\;n\in\N\})$ ($\in(1,+\infty)\cup\{+\infty\}$), taken in the linear order 
$\langle\R^*,<\rangle$. Since 
$q<q^2<q^3<\ds$, 
$\lim q^n=L$ by Proposition~\ref{prop_O_mon1}. It remains to show that $L=+\infty$. Let $L<+\infty$. So 
$L/q<q^n\le L$ for some $n\in\N$. We get the contradiction that $L<q^{n+1}$.
\kduk

\section[Subsequences]{Subsequences}\label{sec_subseq}

We obtain some results on subsequences.  

\medskip\noindent
{\em $\bullet$ Subsequences. }Recall the definition of subsequences after Definition~\ref{def_seWoOp}: $(b_n)$ is a~subsequence of $(a_n)$ if $b_n=a_{m_n}$ for every $n$ and some sequence of indices $m_1<m_2<\ds\;$. We write
$(a_n)\preceq(b_n)$.\label{preceq}
Definition~\ref{def_slabPodpo} introduces weak subsequences. The subsequence 
$(a_k,a_{k+1},\ds)$ of $(a_n)$ is a~\underline{tail\index{sequence!tail
of|emph}} of $(a_n)$.

\begin{exer}\label{ex_preqRefTranz}
The relation $\preceq$ on $\R^{\N}$ is reflexive and transitive. 
\end{exer}

\begin{exer}\label{ex_podpNotWeakAntirefl}
Find distinct sequences $(a_n)$ and $(b_n)$ such that $(a_n)\preceq(b_n)$ 
and $(b_n)\preceq(a_n)$.
\end{exer}

\begin{prop}\label{prop_LimPodpo}
Let $(b_n)\preceq(a_n)$ and $\lim a_n=L$. \underline{Then} $\lim b_n=L$.
\end{prop}
\duk
This follows from Definition~\ref{def_vlLimi} and the
fact that in the above sequence of indices, $m_n\ge n$ for every $n$.
\kduk

\noindent 
{\em $\bullet$ Orderings of sets of natural numbers. }Subsequences 
arise by deleting terms from the given sequences. We can delete only so 
many terms that infinitely 
many of them remain. We formalize it in the next proposition. A~map 
$f\cc A\to\N$, where $A\sus\N$, \underline{increases\index{increasing map to $\N$|emph}} if always $m<n$ $\Rightarrow$ $f(m)<f(n)$.

\begin{prop}\label{prop_ordSet}
Let $B\sus\N$. If $B$ is finite, \underline{then} there exist a~unique number $m\in\omega$ and
a~unique increasing bijection $f\cc[m]\to B$. If $B$ is infinite, \underline{then} there exists a~unique
increasing bijection $f\cc\N\to B$. In both cases, we call $f$ the \underline{ordering} of $B$\index{ordering of $B\sus\N$|emph}. 
\end{prop}
\duk
For finite $B$, the number $m\in\omega$ is the number of elements of $B$. 
The map $f$ is defined by the inductive formula 
$$
f(n)=\min(B\setminus f[\,[n-1]\,])\,,
$$
where $n\in[m]$, respectively $n\in\N$. 
\kduk

We can restate subsequences in terms of orderings of sets, 

\begin{prop}\label{prop_ekviSubs}
$(b_n)\preceq(a_n)$ $\iff$ there is a~unique infinite set $B\sus\N$ such 
that $b_n=a_{f(n)}$ for every $n\in\N$, where $f$ is the ordering of $B$.
\end{prop}
\duk
Suppose that $(b_n)$ is a~subsequence of $(a_n)$, so that $b_n=a_{m_n}$ for a~sequence of natural numbers $1\le 
m_1<m_2<\cdots$. It follows that $B=
\{m_n\cc\;n\in\N\}$. Clearly,  $f(n)=m_n$. 

Suppose that the right-hand side of the equivalence holds. For $n\in\N$ 
we set $m_n=f(n)$. Then $m_1<m_2<\ds$ because $f$ is 
increasing and $b_n=a_{m_n}$. So $(b_n)$
is a~subsequence of $(a_n)$.
\kduk

\noindent
We call the set $B$ ($\sus\N$) the \underline{support\index{support of a~subsequence|emph}} of the subsequence $(b_n)$ of $(a_n)$. 

\begin{defi}\label{def_slabPodpo}
We say that $(b_n)$ is a~\underline{weak 
subsequence\index{weak subsequences|emph}} of $(a_n)$ if there 
is a~sequence $(m_n)\sus\N$ such that $\lim m_n=+\infty$ and $b_n=a_{m_n}$ 
for every $n$. We write $(b_n)\preceq^*(a_n)$. 
\end{defi}

\begin{exer}\label{ex_slabaPodposl}
Generalize Proposition~\ref{prop_LimPodpo}: if 
$(b_n)\preceq^*(a_n)$ and $\lim a_n=L$, then $\lim\,b_n=L$.
\end{exer}

\begin{exer}\label{ex_oSlPodp}
If $(b_n)\preceq^*(a_n)$, then there is a~sequence $(c_n)$ such that 
$(c_n)\preceq(b_n)$ and $(c_n)\preceq(a_n)$.
\end{exer}

\noindent
{\em $\bullet$ Limit dualities. }In Theorem~\ref{thm_oPodposl} and Exercise~\ref{ex_3rdDual}
we present three results to the effect that the non-existence of limits of
certain sequences is equivalent to the existence of limits of certain 
subsequences. We recall a~result from 

\begin{prop}\label{prop_exSUbLim}
Every sequence $(a_n)\sus\R$ has a~subsequence $(b_n)$ such that $(b_n)$
has a~limit in $\R^*$.
\end{prop}
\duk
\kduk

The first two limit dualities are as follows. 

\begin{thm}\label{thm_oPodposl}
Let\index{theorem!two limit dualities|emph} 
$(a_n)\sus\R$ and $A\in\R^*$. The following holds.
\begin{enumerate}
\item $\lim a_n$ does not exist $\iff$ two
subsequences of $(a_n)$ have different limits. 
\item It is not true that $\lim a_n=A$ $\iff$ a~subsequence $(b_n)\preceq(a_n)$ has a~limit 
$\lim b_n\ne A$.
\end{enumerate}
\end{thm}
\duk
1. The implication $\neg\Rightarrow\neg$ follows from Proposition~\ref{prop_LimPodpo}. We prove the implication $\Rightarrow$. Suppose that the sequence 
$(a_n)$ does not have a~limit. By Proposition~\ref{prop_exSUbLim}, there is a~sequence $(b_n)\preceq(a_n)$   
with $\lim b_n=B$. Since $B$ is not a~limit of $(a_n)$, there exists an 
$\ep$ and a~sequence 
$(c_n)\preceq(a_n)$ such that $c_n\not\in U(B,\ep)$ for every $n$. By 
Proposition~\ref{prop_exSUbLim} there is a~sequence $(d_n)\preceq(c_n)$ such that $\lim d_n=K$. 
Then $(d_n)\preceq(a_n)$ and $K\ne B$. The desired subsequences are $(b_n)$ and $(d_n)$.

2. The implication $\neg\Rightarrow\neg$ again follows from Proposition~\ref{prop_LimPodpo}. We prove the implication 
$\Rightarrow$. Let $\neg (\lim a_n=A)$. Hence there is an $\ep$ 
and a~sequence $(b_n)\preceq(a_n)$ such that $b_n\not\in U(A,\ep)$ for 
every $n$. By Proposition~\ref{prop_exSUbLim},  
there is a~sequence $(c_n)\preceq(b_n)$ such that $\lim c_n=B$. Then 
$(c_n)\preceq(a_n)$ 
and $B\ne A$. So $(c_n)$ is the desired subsequence.
\kduk

\noindent
We see that if a~sequence 
does not have a~limit, it is always possible to prove it by presenting two 
subsequences with different limits. For example,  
$$
(a_n)=((-1)^n)=(-1,\,1,\,-1,\,1,\,-1,\,\ds)
$$ 
does not have a~limit because 
$$
(1,\,1,\,\ds)\preceq(a_n)\,\text{ and }\,(-1,\,-1,\,\ds)\preceq(a_n)\,, 
$$
and these constant subsequences have
different limits $1$ and $-1$.

\begin{exer}\label{ex_disjointSub}
In part~1 of the theorem, the two subsequences may have disjoint supports.    
\end{exer}

Here is another limit duality.

\begin{exer}\label{ex_3rdDual}
$(a_n)$ diverges $\iff$ $(a_n)$ has two subsequences with two different limits, or a~subsequence with limit $\pm\infty$.
\end{exer}

\noindent
{\em $\bullet$ Partitions into subsequences. }Let $(a_n)\sus\R$. A~\underline{partition\index{sequence!partitions in subsequences|emph}} 
of $(a_n)$ into subsequences is given by any partition $P$ of $\N$ with 
infinite blocks. To any block $Z\in P$ corresponds the subsequence $(b_{Z,n})$
of $(a_n)$ such that $b_{Z,n}=a_{f(n)}$ for every $n$, where $f\cc\N\to Z$
is the ordering of $Z$.

\begin{thm}\label{thm_finManyBl}
Let $A\in\R^*$. 
If\index{theorem!finite partitions of sequences|emph} 
a~sequence $(a_n)\sus\R$ has a~partition into $k\in\N$
subsequences such that each has limit $A$, \underline{then} $\lim 
a_n=A$ as well.   
\end{thm}
\duk
Let $\neg\lim a_n=A$ and let the subsequences in the partition be 
$(a_{i,n})$ with $i=1,2,\ds,k$. By part~2 in Theorem~\ref{thm_oPodposl},  
there is a~sequence $(b_n)\preceq(a_n)$ such that $\lim b_n=B\ne A$. By the 
pigeonhole principle, there exists a~sequence $(c_n)$ and an index 
$i\in[k]$ such that $(c_n)\preceq(b_n)$
and $(c_n)\preceq(a_{i,n})$. Since $\lim c_n=B\ne 
A=\lim_{n\to\infty}a_{i,n}$, it is a~contradiction.
\kduk

\noindent
We use this theorem in the proofs of Theorems~\ref{thm_altSer} and \ref{thm_DerSlozFce}. For infinite partitions, an opposite result holds.

\begin{thm}\label{thm_infinManyBl}
Every\index{theorem!infinite partitions of sequences|emph} 
sequence $(a_n)\sus\R$ can be partitioned into infinitely many
subsequences with the same limit.   
\end{thm}
\duk
Let $(a_n)\sus\R$ be arbitrary. By Proposition~\ref{prop_exSUbLim}, there 
is an $L\in\R^*$ and an infinite set $B_0\sus\N$ such that $\lim b_n=L$, 
where $(b_n)=(a_{f(n)})$ for the ordering $f$ of $B_0$. By Exercise~\ref{ex_whyAssu}, we may 
assume that the
complement $\N\setminus B_0=
\{c_1<c_2<\ds\}$ is infinite. We take any infinite partition
$\{C_j\cc\;j\in\N\}$ of $B_0$ with infinite blocks $C_j$ and set $B_j=C_j\cup\{c_j\}$. Then
$$
\{B_j\cc\;j\in\N\}
$$
is an infinite partition of $\N$ with infinite blocks. The corresponding subsequences $(b_{j,n})=(a_{f_j(n)})$ of 
$(a_n)$, where $f_j$ is the ordering 
of $B_j$, form the desired partition of $(a_n)$; it is clear (why?\,---\,Exercise~\ref{ex_whyClear}) that $\lim_{n\to\infty}b_{j,n}=L$
for every $j\in\N$.
\kduk
\vspace{-3mm}
\begin{exer}\label{ex_whyAssu}
Why can we assume that the complement $\N\setminus B_0$
is infinite?    
\end{exer}

\begin{exer}\label{ex_whyClear}
Why is it clear?    
\end{exer}

\begin{exer}\label{ex_partBound}
Show that every bounded sequence 
can be partitioned in infinitely many convergent subsequences with the same limit.   
\end{exer}

\begin{exer}\label{ex_naParti}
Show that every sequence that does not have a~limit can be 
partitioned in infinitely many subsequences such that one of them 
has limit $A$, and any other has limit $B\ne A$. 
\end{exer}

\section[Limes inferior and limes superior]{Limes inferior and limes superior}\label{podkap_liminLimsup}

In Latin, this means ``the lowest limit'' and ``the 
highest limit'', respectively. Limes inferior, briefly liminf, and limes superior, 
briefly limsup, of a~real sequence always exist, which is an advantage 
over limits.

\medskip\noindent
{\em $\bullet$ Limit points of a~sequence. }They are limits of subsequences of the sequence. 

\begin{defi}\label{def_hromBody}
We say that $A\in\R^*$ is a~\underline{limit 
point}\index{limit point!of a sequence|emph} of a~sequence 
$(a_n)$ if $A=\lim b_n$ for a~sequence $(b_n)\preceq(a_n)$. 
The set of limit points of $(a_n)$ is denoted by $L(a_n)$\index{limit point!of a 
sequence!lan@$L(a_n)$|emph} ($\sus\R^*$).\label{Lan} 
\end{defi}
For example, the sequence 
$${\textstyle
(a_n)=(n-1+(-1)^n n+\frac{1}{n})
}
$$ 
has two limit points, $L(a_n)=\{-1, +\infty\}$.

\begin{exer}\label{ex_existHB}
Every real sequence has at least one limit point.    
\end{exer}

\noindent
{\em $\bullet$ Limes inferior and limes superior of a~sequence. }We 
have already revealed that these are the smallest and 
largest limit points of the sequence, respectively. 

\begin{defi}\label{def_liminfLimsup}
Let $(a_n)\sus\R$. We define 
$$
\liminf_{n\to\infty} a_n=
\liminf a_n\index{sequence!liminf of|emph}:=\min(L(a_n)),\ 
\limsup_{n\to\infty} a_n=
\limsup a_n\index{sequence!limsup of|emph}:=\max(L(a_n))\,. 
$$
The 
minimum and maximum are taken in the linear order $\langle\R^*,<\rangle$. \label{liminf}\label{limsup}
\end{defi}
We show that liminf and limsup always exist.

\begin{thm}\label{thm_limsup}
The following holds.
\begin{enumerate}
\item 
$L(a_n)\ne\emptyset$\index{theorem!liminf and limsup exist|emph} 
for every sequence $(a_n)\sus\R$.
\item In the linear order $\langle\R^*,<\rangle$, the set 
$L(a_n)$ has minimum and maximum.
\end{enumerate} 
\end{thm}
\duk
1. Exercise~\ref{ex_existHB}. 

2. We show the existence of $\max(L(a_n))$; the 
minimum is treated similarly. Let $A:=\sup(L(a_n))$, in the 
linear order $\langle\R^*,<\rangle$ 
We show that $A\in L(a_n)$. This holds if $A=-\infty$, then
$L(a_n)=\{-\infty\}$. Let $A>-\infty$. There is a~sequence $(b_n)\sus 
L(a_n)\cap\R$ such that $\lim b_n=A$: if $A<+\infty$ or if  
$A=+\infty\not\in L(a_n)$, it follows from the definition of supremum, and if
$A=+\infty\in L(a_n)$, we are done. Since every number $b_n$
is the limit of a~subsequence of $(a_n)$, it is easy to find 
a~subsequence $(a_{m_n})$ of $(a_n)$ such that $a_{m_n}\in 
U(b_n,1/n)$ for every $n$. Then $\lim 
a_{m_n}=\lim b_n=A$ and $A\in L(a_n)$.
\kduk

Clearly, if $\lim a_n$ exists, then $L(a_n)=\{\lim a_n\}$. We obtain
some more properties of liminfs and limsups.

\begin{prop}\label{prop_limiJeLims}
The following holds.
\begin{enumerate}
\item Always $\liminf a_n\le\limsup a_n$.
\item The inequality in part~1 holds as an equality $\iff$ $\lim a_n$ exists. 
Then $\liminf a_n=\limsup\,a_n=\lim a_n$.
\end{enumerate}  
\end{prop}    
\duk
Let $(a_n)\sus\R$. 1. 
This is obvious as $\liminf a_n=\min(L(a_n))$ and $\limsup 
a_n=\max(L(a_n))$. 

2. If equality holds, then $L(a_n)$ has just one 
element and $(a_n)$ does not have two subsequences with different limits. 
By part~2 of Theorem~\ref{thm_oPodposl},  $\lim a_n$ exists and equals to 
$\liminf a_n=\limsup a_n$. If $\liminf 
a_n\ne\limsup a_n$, then the sequence $(a_n)$ has two subsequences with 
different limits and $\lim a_n$ does not exist.
\kduk
\vspace{-3mm}
\begin{thm}\label{thm_VlLimsup}
Let\index{theorem!properties of liminfs and limsups|emph} 
$(a_n)\sus\R$, $A:=\liminf a_n$, and $B:=\limsup a_n$. 
The following holds.
\begin{enumerate}
\item If $A=-\infty$ \underline{then} for every $c$ we have $a_n\le c$ for infinitely 
many $n$. If $A=+\infty$ \underline{then} $\lim a_n=+\infty$.
\item If $A\in\R$ \underline{then} for every $\ep$ we have $a_n\le A+\ep$ for 
infinitely many $n$, and  $a_n\ge A-\ep$ for every $n\ge n_0$.
\item If $B=+\infty$ \underline{then} for every $c$ we have $a_n\ge c$ for infinitely 
many $n$. If $B=-\infty$ \underline{then} $\lim a_n=-\infty$.
\item If $B\in\R$ \underline{then} for every $\ep$ we have $a_n\ge B-\ep$ for 
infinitely many $n$, and $a_n\le B+\ep$ for every $n\ge n_0$.
\end{enumerate}
\end{thm}
\duk
We prove parts 1 and 2. Parts~3 and 4 are left to Exercise~\ref{ex_casti3a4}.
1. Let $A=-\infty$. Then for some sequence $(b_n)\preceq(a_n)$ we have 
$\lim b_n=-\infty$ and the claim follows. Let $A=+\infty$. Then $L(a_n)=\{+\infty\}$ and by 
Proposition~\ref{prop_limiJeLims} we have $\lim a_n=+\infty$.

2. Let $A\in\R$ and $\ep$ be given. Since there is a~sequence 
$(b_n)\preceq(a_n)$ with 
$\lim b_n=A$, the inequality $a_n\le A+\ep$ holds for infinitely many $n$. If we had $a_n<A-\ep$ 
for infinitely many $n$, a~sequence $(c_n)\preceq(a_n)$ would exist with $\lim c_n\le A-\ep$. 
This is impossible because $A=\min(L(a_n))$. Hence $a_n\ge A-\ep$ 
for every $n\ge n_0$.
\kduk

We relate the liminf and limsup to infima and suprema in $\langle\R^*,<\rangle$.

\begin{prop}\label{prop_onLiminf}
Let $(a_n)\sus\R$. \underline{Then}
\begin{eqnarray*}
\liminf_{n\to\infty}a_n&=&\lim_{n\to\infty}\inf(\{a_m\cc\;m\ge n\})\,\text{ and}\\
\limsup_{n\to\infty}a_n&=&\lim_{n\to\infty}\sup(\{a_m\cc\;m\ge n\})\,,
\end{eqnarray*}
where we define
$$
\lim_{n\to\infty} (\pm\infty,\,
\pm\infty,\,\ds):=\pm\infty\,,
$$
with all top or all bottom signs.
\end{prop}
\duk
We prove the former formula, the proof of the latter is similar. Let 
$$
A_n:=\inf(\{a_m\cc\;m\ge n\})\ \ 
(\in\R\cup\{-\infty\})\,.
$$ 
Clearly, 
$A_1\le A_2\le\cdots$. It is easy to see that if $A_1=-\infty$, then $A_n=-
\infty$  for every $n$, and that then $(a_n)$ has a~subsequence $(a_{m_n})$
with $\lim a_{m_n}=-\infty$. So  
$\liminf a_n=-\infty$, and the equality holds.

Let $A_1\in\R$. Then $A_n\in\R$ for every $n$. We set 
$$
A:=\sup(\{A_n\cc\;n\in\N\})\ \  (\in\R\cup\{+\infty\})\,. 
$$
Clearly, if $A=+\infty$, then $\lim a_n=
\lim A_n=+\infty$ (Proposition~\ref{prop_O_mon1}) and the equality again holds.

Finally, let $A\in\R$. Since $A_1\le A_2\le\ds$, we have 
that $\lim A_n=A$ by Proposition~\ref{prop_O_mon1}. We show 
that $\liminf a_n=A$. Let an $\ep>0$ be given. Then $a_n\le A+\ep$ for 
infinitely many $n$, for else we would have $A_n\ge A+\ep$ for $n\ge 
n_0$. Thus, $\liminf a_n\le A+\ep$. On the other hand, $A_{n_0}\ge A-\ep$ 
for some $n_0$, so that $\liminf a_n\ge A-\ep$. We see that $\liminf 
a_n=A=\lim A_n$.
\kduk

We illustrate the use of liminfs and limsups by a~number-theoretic result. Recall that the \underline{function 
$\tau(n)$\index{number of divisors, $\tau(n)$|emph}}\label{tau} counts 
the divisors of $n\in\N$. For example, 
$$
\tau(6)=|\{1,\,2,\,3,\,6\}|=4\,. 
$$
It is known that 
$$
\limsup_{n\to\infty}\frac{\log(\tau(n))}{(\log 2)(\log n)/(\log\log n)}=1\,\text{ and }\,\liminf_{n\to\infty}\tau(n)=2\,.
$$
\begin{exer}\label{ex_lininftau}
Prove the latter equality.  
\end{exer}

\begin{exer}\label{ex_casti3a4}
Prove parts 3 and 4 of the last theorem.    
\end{exer}

\begin{exer}\label{ex_naliminf1}
Find a~sequence $(a_n)\sus\R$ such that $L(a_n)=\R^*$.
\end{exer}

\begin{exer}\label{ex_naliminf2}
Show that $L(a_n)=[-1,1]\setminus\{0\}$ for no sequence $(a_n)\sus\R$.
\end{exer}

\begin{exer}\label{ex_naliminf3}
Find $\liminf a_n$ and $\limsup a_n$ for $a_n=n(1+(-1)^n)$ .
\end{exer}

\section[Existence theorems for limits]{Existence theorems for limits}\label{podkap_ctyriRobustni}

We obtain or review several theorems on the existence of limits of real sequences, namely  
Theorems~\ref{thm_O_mon2}, \ref{thm_kvazimon}, 
\ref{thm_BolzWeier}, and \ref{thm_CauchyPodm}.

\medskip\noindent {\em $\bullet$ Monotonicity and boundedness. }We say that a~sequence $(a_n)$ 
\underline{weakly increases\index{sequence!weakly increases|emph}} (respectively 
\underline{weakly decreases\index{sequence!weakly decreases|emph}}) if for every $n$ 
we have $a_n\le a_{n+1}$ (respectively $a_n\ge a_{n+1}$). It 
\underline{increases\index{sequence!increases|emph}} (respectively \underline{decreases\index{sequence!decreases|emph}}) if for every $n$ we have $a_n<a_{n+1}$ (respectively, $a_n>a_{n+1}$). It is 
\underline{monotone\index{sequence!monotone|emph}} if it weakly decreases or weakly increases. It 
is \underline{strictly monotone\index{sequence!strictly 
monotone|emph}} if it decreases or increases.

We say that $(a_n)$ is \underline{bounded from 
above\index{sequence!bounded from above|emph}} if there is
a~$c$ such that for every $n$ we have $a_n\le c$. Reversing the inequality, we get \underline{ boundedness from
below\index{sequence!bounded from below|emph}}. A~sequence $(a_n)$ 
is \underline{bounded\index{sequence!bounded|emph}} if it is bounded
both from above and from below.

\begin{exer}\label{ex_onMonot}
$(a_n)$ weakly increases iff $m\le n$ $\Rightarrow$ $a_m\le a_n$. 
State and prove analogous equivalences for  weakly 
decreasing, increasing and decreasing sequences.   
\end{exer}

\begin{exer}\label{ex_ekvivOmez}
$(a_n)$ is bounded iff there is a~$c$ such that $|a_n|\le c$ for every $n$.    
\end{exer}

\begin{exer}\label{ex_petRobustnich}
Which of the above nine underlined properties of sequences are robust?
\end{exer}

\noindent
{\em $\bullet$ Limits of monotone sequences. }We proved the existence of limits of monotone sequences in 
Proposition~\ref{prop_O_mon1}. We give a~robust version of this result.

\begin{thm}\label{thm_O_mon2}
Let\index{theorem!monotone sequences|emph} 
$(a_n)$ be a~real sequence. If a~tail $(a_m,a_{m+1},\ds)$ of $(a_n)$ is monotone, \underline{then}  
$$
\lim_{n\to\infty}a_n=\sup(\{a_m,a_{m+1},\ds\}),\,\text{ resp. }\,\lim_{n\to\infty}a_n=
\inf(\{a_m,a_{m+1},\ds\})\,,
$$
depending on whether the tail weakly increases or weakly decreases. Infima and suprema are taken in $\langle\R^*,<\rangle$.
\end{thm}
\duk
Clearly, the limit of the tail $(a_m,a_{m+1},\ds)$ is the limit of $(a_n)$.
\kduk
\vspace{-3mm}
\begin{exer}\label{ex_uloNaMonPos}
The assumption on $(a_n)$ is a~robust property of sequences.
\end{exer}

\noindent 
{\em $\bullet$ Limits of quasi-monotone sequences. }There are generalizations of monotone sequences. A~sequence $(a_n)\sus\R$ \underline{goes up\index{sequence!goes up|emph}} (respectively \underline{goes down\index{sequence!goes down|emph}})
if for every index $n$, the inequality $a_m<a_n$ (respectively $a_m>a_n$) holds for only finitely many indices $m$. We say that $(a_n)$ is
\underline{quasi-monotone\index{sequence!quasi-monotone|emph}} if it goes up or goes down.

\begin{exer}\label{ex_monoJekvazim}
Every monotone sequence is quasi-monotone.    
\end{exer} 

\begin{exer}\label{ex_oKvazimon1}
Find a~quasi-monotone sequence $(a_n)$ such that for no $m$ the tail $(a_m,a_{m+1},\ds)$ is monotone. 
\end{exer}

\begin{exer}\label{ex_oKvazimon2}
Write the property of quasi-monotonicity in terns of quantifiers, logical connectives, brackets, 
variables, and inequalities, which are applied to natural and real 
numbers.
\end{exer}

\begin{thm}\tec\label{thm_kvazimon}
Every\index{theorem!limits of quasi-monotone
sequences|emph} quasi-monotone sequence $(a_n)$ 
has a~limit.
\end{thm}
\duk
We assume that $(a_n)$ goes up, the argument for $(a_n)$ going down is similar. Let 
$A:=\limsup a_n$ and $\ep$ be given. Then 
$(a_n)$ has a~subsequence $(a_{m_n})$ with $\lim a_{m_n}=A$ and we have $a_n<A+\ep$ for every 
$n\ge n_0$. We take an $n'$ such that 
$a_{m_{n'}}\in U(A,\ep)$. Since $(a_n)$ goes up, we can take an $n_1\ge 
n_0$ such that $a_{m_{n'}}\le a_n< A+\ep$ for every $n\ge n_1$. By Exercise~\ref{ex_nekdVlOkoli}, $a_n\in U(A,\ep)$ for the same $n$. Hence 
$\lim a_n=A$.
\kduk

We state a~robust version of this result; we omit the easy proof.

\begin{cor}\label{cor_robuKMP}
Every sequence $(a_n)$ that has a~quasi-monotone tail
has a~limit.
\end{cor}

\begin{exer}\label{ex_robuKVMo}
The assumption in the corollary defines a~robust property of sequences.
\end{exer}

\noindent
Quasi-monotone 
sequences were introduced by {\em Godfrey H. Hardy\index{Hardy, 
Godfrey H.} (1877--1947)}. 

\medskip\noindent
{\em $\bullet$ The Bolzano--Weierstrass theorem. }Part~1 of the next theorem is 
known as the Bolzano--Weierstrass theorem.

\vspace{-3mm}
\begin{thm}\label{thm_BolzWeier}
Let\index{theorem!Bolzano--Weierstrass|emph} 
$(a_n)\sus\R$. \underline{Then} one of three cases occurs.
\begin{enumerate}
\item The sequence $(a_n)$ is bounded and has a~convergent subsequence.
\item The sequence $(a_n)$ is not bounded from above and has a~subsequence $(b_n)$ such that $\lim b_n=+\infty$.
\item The sequence is not bounded from below and has a~subsequence $(b_n)$ such that $\lim b_n=-\infty$.
\end{enumerate}
\end{thm}
\duk
1. Let $(a_n)$ be bounded and $(b_n)\preceq(a_n)$ be a~monotone 
subsequence guaranteed by Lemma~\ref{lem_ESlemma}. 
The sequence $(b_n)$ is bounded. By Proposition~\ref{prop_O_mon1}, in fact already by 
has a~finite limit.

2. Suppose that $(a_n)$ is not bounded from above. We inductively define a~subsequence $(a_{m_n})$ of $(a_n)$
such that $a_{m_n}\ge n$. Then $\lim a_{m_n}=+\infty$.

We take any $m_1$ such that $a_{m_1}\ge1$. Suppose that 
indices $m_1<m_2<\ds<m_n$ are defined such that 
$a_{m_i}\ge i$ for $i=1,2,\ds,n$. Since $(a_n)$ is not bounded from above, there exists $j\in\N$ such that
$$
a_j\ge1+
\max(\{a_1,\,a_2,\,\ds,\,a_{m_n}\}\cup
\{n+1\})\,.
$$
Then $j>m_n$ and $a_j\ge n+1$. We set $m_{n+1}:=j$.

3. Exercise~\ref{ex_part3}. 
\kduk
\vspace{-3mm}
\begin{exer}\label{ex_part3}
Prove part~3 of the theorem.    
\end{exer}
\begin{exer}\label{ex_verzeBWthh}
Let $a\le b$ be real numbers. Then every sequence $(a_n)\sus[a,b]$ 
has a~subsequence $(a_{m_n})$ such that $\lim a_{m_n}\in[a,b]$.
\end{exer}

\noindent
{\em Bernard Bolzano\index{Bolzano, Bernard} (1781--1848)} was a~Czech--German--Italian priest, philosopher, 
and  mathematician. {\em Karl 
Weierstrass\index{Weierstrass, Karl} (1815--1897)} was a~German 
mathematician.

\medskip\noindent {\em $\bullet$ Real Cauchy sequences. }We used rational 
Cauchy sequences to define real numbers 
Cauchy sequences are defined in the 
same way. 

\begin{defi}\label{def_cauchyPosl}\index{sequence!Cauchy!real|emph}
A~sequence $(a_n)\subset\R$ is \underline{Cauchy} if for every $\ep$ there is an $n_0$ such that for every $m,n\ge n_0$ we have 
$$
|a_m-a_n|\le\ep\,.
$$
\end{defi}

\begin{exer}\label{ex_cauchyRobust}
 Cauchy sequences form a~robust property of sequences. 
\end{exer}

\begin{exer}\label{ex_CauchyJeOmez}
Every Cauchy sequence is bounded. 
\end{exer}
\vspace{-3mm}
\begin{thm}\label{thm_CauchyPodm}
A~real\index{theorem!metric completeness of 
$\R$|emph} sequence converges if and only if it is Cauchy.
\end{thm}
\duk
Implication $\Rightarrow$. Let $\lim a_n=a$ and $\ep$ be given. For 
every large $n$ we have $|a_n-a|\le\frac{\ep}{2}$. 
Using TI (Exercise~\ref{ex_trojNero}) we have for every large $m$ and $n$ that
$${\textstyle
|a_m-a_n|\le|a_m-a|+|a-a_n|\le\frac{\ep}{2}+\frac{\ep}{2}=\ep\,.
}
$$
We see that $(a_n)$ is Cauchy.

Implication $\Leftarrow$. Let $(a_n)$ be Cauchy. By Exercise~\ref{ex_CauchyJeOmez}, $(a_n)$ is bounded. By the Bolzano--Weierstrass theorem, $(a_n)$ has a~convergent subsequence
$(a_{m_n})$ with a~limit $a$. Let an $\ep$ be given. There is an $n_0$ such that for every $m,n\ge n_0$, we have $|a_{m_n}-a|\le\frac{\ep}{2}$ and $|a_m-a_n|\le\frac{\ep}{2}$.
By TI, we have for every $n\ge n_0$, since $m_n\ge n$, that
$${\textstyle
|a_n-a|\le|a_n-a_{m_n}|+|a_{m_n}-a|\le\frac{\ep}{2}+\frac{\ep}{2}=\ep\,.
}
$$
We see that $\lim a_n=a$.
\kduk

\noindent
Interestingly, A.-L. Cauchy\index{Cauchy, Augustin-Louis} lived in Prague in the Austrian Empire for several years in political exile. 

\begin{exer}\label{ex_CauchyVeQ}
There is a~Cauchy sequence $(a_n)\sus\Q$ such that $\lim a_n\not\in\Q$. 
\end{exer}
The previous theorem therefore does not hold in the ordered field 
$\Q$. This is not surprising\,---\,we know that 
$\Q$ is not complete. 

\begin{exer}\label{ex_kdeuplnost}
Where is the completeness of $\R$ used in the previous proof?
\end{exer}

\section[Fekete's lemma in combinatorics]{Fekete's lemma in combinatorics}\label{sec_fekete}

Fekete's lemma is another existence theorem for limits. We state its 
additive and multiplicative forms as 
a~theorem and its corollary, respectively. In this section, we are 
interested in applications of this lemma in combinatorics.

\medskip\noindent
{\em $\bullet$ Additive Fekete's lemma. }The lemma is credited to the 
Hungarian--Israeli mathematician {\em Michael Fekete\index{Fekete, Michael} (1886--1957)}.

\begin{exer}\label{ex_fekete}
``fekete'' means $\ds\ $.
\end{exer}

\begin{exer}\label{ex_fekete2}
Solve the previous exercise by methods available in 1984.
\end{exer}

A~sequence $(a_n)\sus\R$ is \underline{superadditive\index{sequence!superadditive|emph}}, respectively 
\underline{subadditive\index{sequence!subadditive|emph}}, if for every two
indices $m$ and $n$ we have $a_{m+n}\ge a_m+a_n$, respectively
$a_{m+n}\le a_m+a_n$.

\begin{thm}[Fekete's lemma]\label{thm_fekete}
Let\index{theorem!Fekete's lemma|emph} 
$(a_n)\sus\R$ and
$M=\{a_n/n\cc\;n\in\N\}$. If $(a_n)$ is superadditive, resp. subadditive, 
\underline{then}
$${\textstyle
\lim\frac{a_n}{n}=\sup(M),\, \text{ resp. }\,\lim\frac{a_n}{n}=\inf(M)\,.
}
$$
The supremum and infimum are taken in $\langle\R^*,<\rangle$.
\end{thm}
\duk
Let $(a_n)$ be superadditive; the other case is similar. Let an $\ep$ be given. We take a~number $c\in 
U(\sup(M),\ep)$ with $c<\sup(M)$. Then $\frac{a_m}{m}>c$ for some $m$. We 
write any $n\in\N$ with $n\ge m$ as $n=km+l$, where $k\in\N$, 
$l\in\N_0$ and $0\le l<m$. We have 
$a_n\ge ka_m+a_l$ (Exercise~\ref{ex_naAFeLe}), so that
$${\textstyle
\frac{a_n}{n}\ge\frac{ka_m}{km+l}+\frac{a_l}{n}=\frac{a_m/m}{1+l/km}+\frac{a_l}{n}\,.
}
$$
For $n\to\infty$ also $k\to\infty$. So $1+\frac{l}{km}\to1$ and 
$\frac{a_l}{n}\to0$. Hence for every $\de$ we have for large $n$ that $\frac{a_n}{n}\ge\frac{a_m}{m}-\de$. Thus there is an $n_0\ge m$ such that for every $n\ge n_0$,
$$
c<a_n/n\le\sup(M)\,. 
$$
By Exercise~\ref{ex_nekdVlOkoli}, $\frac{a_n}{n}\in U(\sup(M),\ep)$ for the same $n$. Hence 
$\frac{a_n}{n}\to\sup(M)$. 
\kduk

\begin{exer}\label{ex_naAFeLe}
If $(a_n)$ is superadditive and $n=km+l$, then $a_n\ge ka_m+a_l$.    
\end{exer}

\noindent
{\em $\bullet$ Multiplicative Fekete's lemma. }We say that 
a~sequence $(a_n)\sus(0,+\infty)$ is \underline{supermultiplicative\index{sequence!supermultiplicative|emph}}, resp. 
\underline{submultiplicative\index{sequence!submultiplicative|emph}}, if for every two
indices $m$ and $n$ we have $a_{m+n}\ge a_ma_n$, resp.
$a_{m+n}\le a_ma_n$. The real 
power $a^b$ and functions $\log x$ and $\exp x$ used in the the next corollary are introduced in Section~\ref{sec_elemenFce}. 

\begin{exer}\label{ex_whySuperad}
Prove the next two lemmas.
\end{exer}

\begin{lemma}\label{lem_forFeke}
Let $(c_n)\sus\R$ and $\lim c_n=K$, where $K\in\R\cup\{+\infty\}$. \underline{Then} $\lim\exp(c_n)=
\exp(K)$, where we define $\exp(+\infty)=+\infty$.
\end{lemma}

\begin{lemma}\label{lem_forFeke1}
Let $\emptyset\ne X\sus\R$ with $A=\sup(X)\in\R\cup\{+\infty\}$ and 
let $Y=\{\exp x\cc\;x\in X\}$. \underline{Then} $\sup(Y)=\exp(A)$.   
\end{lemma}

\begin{cor}[multiplicative Fekete's lemma]\label{cor_mulFeke}
Let $(a_n)\sus(0,+\infty)$ and let
$M=\{a_n^{1/n}\cc\;n\in\N\}$.
If $(a_n)$ is supermultiplicative, resp. submultiplicative, 
\underline{then}
$${\textstyle
\lim a_n^{1/n}=\sup(M),\,\text{ resp. }\,\lim a_n^{1/n}=\inf(M)\,.
}
$$
Suprema and infima are taken in $\langle\R^*,<\rangle$.    
\end{cor}
\duk
We assume that $(a_n)$ is supermultiplicative (the 
submultiplicative case is similar)
and set $b_n=\log a_n$. The sequence $(b_n)$ is superadditive: $b_{m+n}=\log(a_{m+n})
\ge\log(a_ma_n)=\log(a_m)+\log(a_n)=b_m+b_n$. By Theorem~\ref{thm_fekete},
$$
\lim_{n\to\infty}b_n/n=\sup(\{b_n/n\cc\;n\in\N\})=B\ \ (\in\R\cup\{+\infty\})\,.
$$
Since $b_n/n=\log\big(a_n^{1/n}\big)$ and $a_n^{1/n}=\exp(b_n/n)$, using 
Lemma~\ref{lem_forFeke} and \ref{lem_forFeke1} in the second and 
third equality, respectively, we obtain
$$
\lim_{n\to\infty}a_n^{1/n}=
\lim_{n\to\infty}\exp(b_n/n)=\exp(B)=\sup(M)\,.
$$
\kduk

\noindent
{\em $\bullet$ Fekete's lemma in combinatorics. }We present five applications of Fekete's lemma: two 
in extremal combinatorics\index{extremal combinatorics} and three in
enumerative combinatorics\index{enumerative 
combinatorics}.

\medskip\noindent
{\em $\bullet$ Extremal functions of words. }Let $u=a_1a_2\ds a_j$
and $v=b_1b_2\ds b_k$ with
$j,k\in\N$ be two words. We write $u\preceq v$ and
say that $u$ is \underline{contained\index{words!containment of|emph}} in $v$
if there exist $j$ indices $1\le 
i_1<i_2<\ds<i_j\le k$ such that for
every $l,l'\in[j]$ we have
$$
b_{i_l}=b_{i_{l'}}\iff
a_l=a_{l'}\,.
$$
In other words, a~subsequence in $v$ has the same equality pattern 
of letters as the word $u$. Let $r\in\N$. The word $v$ is 
$r$-\underline{sparse\index{words!sparse@$r$-sparse|emph}} if for any two indices
$1\le l<l'\le k$ with $b_l=b_{l'}$ we have $l'-l\ge r$. 

\begin{defi}[$\mathrm{ex}(u,n)$]
Suppose that $u=a_1a_2\ds a_j$ is a~word which uses 
$r=|\{a_1,a_2,\ds,a_j\}|$ 
distinct letters. Its 
\underline{extremal function\index{extremal functions of words|emph}} 
$\mathrm{ex}(u,n)\cc\N\to\N_0$ is defined by
$$
\mathrm{ex}(u,\,n):=
\max(\{|v|\cc\;\text{$v\in[n]^*$ is an $r$-sparse word
such that $u\not\preceq v$}\})\,.\label{exuen}
$$
\end{defi}

\begin{exer}\label{ex_exIsDef}
$\mathrm{ex}(u,n)$ is correctly defined for every $u\ne\emptyset$ and every $n\in\N$.     
\end{exer}

A~word $u$ is \underline{irreducible\index{words!irreducible|emph}} if it cannot be written as 
a~concatenation $u=vw$ of nonempty words $v$ and $w$ over
two disjoint alphabets. For example, $abab$ is an irreducible word, but $aabb$ is not irreducible. Here is the first application of Fekete's lemma.

\begin{prop}[1st application]\label{prop_irrWor}
Let $u\ne\emptyset$ be an irreducible word.
\underline{Then} the limit
$$
L(u):=\lim_{n\to\infty}\frac{\mathrm{ex}(u,\,n)}{n}\ \ (\in[0,\,+\infty)\cup\{+\infty\})
$$
exists. 
\end{prop}
\duk
In view of Theorem~\ref{thm_fekete}, it suffices to show that $\mathrm{ex}(u,n)$ is superadditive. Let $r$ ($\in\N$) be the number of distinct letters in $u$, let $m,n\in\N$, and let $v$, respectively $w$, be 
a~word witnessing the value $\mathrm{ex}(u,m)$, 
respectively $\mathrm{ex}(u,n)$. Let the word $w'$ be obtained from $w$ by
renaming the elements of $[n]$ using the map $i\mapsto 
i+m$. So $w'$ is over the alphabet 
$\{m+1,m+2,\ds,m+n\}$, $u\not\preceq w'$ and $|w'|=|w|$. The concatenated word $vw'$ over $[m+n]$ is $r$-sparse and $u\not\preceq vw'$, due to the irreducibility of $u$. Hence
$$
\mathrm{ex}(u,\,m)+\mathrm{ex}(u,\,n)=|v|+|w'|=|vw'|\le
\mathrm{ex}(u,\,m+n)\,.
$$
\kduk 

For some short words $u$, the extremal function can be determined exactly.

\begin{exer}\label{ex_abab}
Prove that $\mathrm{ex}(abab,n)=
2n-1$.     
\end{exer}
Thus $L(abab)=2$. In contrast, $L(ababa)=+\infty$, 
but this is quite hard to prove. See the survey article \cite{klaz_DSsurv} for more information on extremal 
functions of words.

\medskip\noindent
{\em $\bullet$ Szemer\'edi's theorem. }An \underline{arithmetic 
progression\index{arithmetic progression|emph}} with 
length $k\in\N$, abbreviated $k$-AP, is any set of integers of the form
$$
X=\{a+jd\cc\;j\in[k]\}\ \  (a\in\Z\,,d\in\N)\,.
$$

\begin{defi}[$r_k(n)$]\label{def_funrkn}
Let $k\in\N$. For $n\in\N$ we define 
$$
r_k(n)=\max(
\{|A|\cc\;\text{$A\sus[n]$ and contains no {\em $k$-AP}}\})\,.\label{erken}
$$
\end{defi}

\begin{exer}\label{ex_trivSzem}
Prove that $r_1(n)=0$ and $r_2(n)=1$ for every $n$.    
\end{exer}\label{ex_rTriv}
For $k\ge3$, the problem of determining or estimating
$r_k(n)$ becomes quite hard. The second application of Fekete's lemma is, however, easy to 
obtain.

\begin{prop}[2nd application]\label{prop_2ndAppl}
Let $k\in\N$. \underline{Then} the limit
$$
r_k:=\lim_{n\to\infty}\frac{r_k(n)}{n}\ \ 
(\in[0,\,1])
$$
exists.   
\end{prop}
\duk
In view of Theorem~\ref{thm_fekete}, it suffices to show that $r_k(n)$ is subadditive.
Let $m,n\in\N$ and let $A\sus[m+n]$ be a~set witnessing the value $r_k(m+n)$. We consider the sets $A'=A\cap[m]$ and 
$A''=\{a\in[n]\cc\;a+m\in A\}$. Clearly, $|A|=|A'|+|A''|$,  
$A'\sus[m]$, $A''\sus[n]$ and both $A'$ and $A''$ avoid $k$-AP. Thus,
$$
r_k(m+n)=|A|=|A'|+|A''|\le r_k(m)+r_k(n)\,.
$$
\kduk

The following simply formulated theorem is one of the most important results in combinatorics in the 20th century. 

\begin{thm}[E.~Szemerédi, 1975]\label{thm_Szemeredi}
For every $k\ge3$\index{theorem!Szemerédi's}  
we have $r_k=0$.
\end{thm}
The proof is quite complicated;
see \cite{szem,tao_szem}. {\em Endre Szemerédi\index{Szemeredi@Szemerédi, Endre} (1940)
} is a~Hungarian mathematician.

\medskip\noindent
{\em $\bullet$ Paths aka 
self-avoiding walks. }Let
$G=\langle V,E\rangle$ with $E\sus\binom{V}{2}$\label{graph}
be a~\underline{graph\index{graphs|emph}} with the set of 
\underline{vertices\index{graphs!vertices of|emph}} $V$
and the set of \underline{edges\index{graphs!edges of|emph}} $E$. Here $\binom{V}{2}:=\{e\cc\;e\sus V\wedge|e|
=2\}$. A~\underline{path\index{graphs!paths in|emph}} in $G$ of \underline{length\index{graphs!paths in!length|emph}} $n$ ($\in\N_0$) is any $n+1$-tuple 
$$
w=\langle v_0,\,v_1,\,\ds,\,v_n\rangle
$$
of \underline{mutually distinct} vertices $v_i\in V$ such that 
$\{v_{i-1},v_i\}\in E$ for every $i\in[n]$.  
The 
graph $G=\langle V,E\rangle$ is \underline{locally finite\index{graphs!locally 
finite|emph}} if for every vertex $v\in V$ there exist only finitely many edges $e\in E$ such that 
$v\in e$. If there is $r\in\N_0$ such that for every vertex 
$v\in V$ the number of such edges equals $r$, we call $G$ \underline{$r$-regular\index{graphs!regular@$r$-regular|emph}}.

Let $G=\langle V,E\rangle$ be a~(typically infinite) graph. 
An \underline{automorphism\index{graphs!
automorphisms of|emph}} of $G$ is any bijection 
$$
f\cc V\to V
$$ 
such that for every edge $e\in E$ we have $f[e],
f^{-1}[e]\in E$. We say that $G$ is \underline{transitive\index{graphs!transitive|emph}} if for every two vertices 
$u,v\in V$ the graph $G$ has an automorphism $f$ such that
$f(u)=v$. A~\underline{LFT graph\index{graphs!LFT|emph}}\label{LFT} is a~locally finite and transitive graph.

\begin{exer}\label{ex_saw1}
Every {\em LFT} graph is $r$-regular for some $r\in\N_0$.    
\end{exer}

Let $G=\langle V,E\rangle$ be an LFT graph, $v\in V$ and 
$n\in\N$. We denote by $p_G(v,n)$ ($\in\N_0$)
the number of paths in $G$ of length $n$ starting at $v$.

\begin{exer}\label{ex_saw2}
 Show that the number $p_G(v,n)$ is finite and independent of $v$.
\end{exer}
For a~LFT graph $G$, 
we define $p_G(n)$ ($\in\N_0$) to be the common value of $p_G(v,n)$ for any 
starting vertex $v$.

\begin{exer}\label{ex_bouNumPath}
Let  $r\in\N$ and $G$ be an 
$r$-regular {\em LFT} graph. Show that $p_G(n)\le r(r-1)^{n-1}$ for every $n\in\N$.
\end{exer}

In the third application, Fekete's lemma proves the existence of 
growth constants $\kappa(G)$ of LFT graphs $G$.

\begin{prop}[3rd application]\label{prop_numPaths}
Let $r\in\N$ and let $G=\langle V,E\rangle$ be an $r$-regular {\em LFT} graph. \underline{Then} the limit
$$
\kappa(G):=\lim_{n\to\infty}p_G(n)^{1/n}\ \ (\in[0,\,r-1])
$$
exists.    
\end{prop}
\duk
In view of Corollary~\ref{cor_mulFeke}, it 
suffices to prove that $p_G(n)$ is submultiplicative. The bound $\kappa(G)\in[0,r-1]$ follows from Exercise~\ref{ex_bouNumPath}. We fix an 
initial vertex $v\in V$ of paths and take, for every $u\in V$, an 
automorphism $f_u\cc V\to V$ of $G$ such that $f_u(u)=v$. We denote
by $P_n$, for 
$n\in\N$, the finite set of paths 
in $G$ of length $n$ starting at 
$v$. Let $m,n\in\N$. We consider the map $F\cc P_{m+n}\to P_m\times 
P_n$ given by
$F(w)=\langle w',f_u[w'']\rangle$,
where 
$$
w=\langle u_0=v,\,u_1,\,\ds,\,u_{m+1},\,\ds,\,u_{m+n+1}\rangle\,,
$$ 
$w'=\langle u_0,\ds,u_{m+1}\rangle$, $w''=\langle u_{m+1},\ds,u_{m+n+1}\rangle$
and $u=u_{m+1}$. It is easy to see that $F$ is injective. Hence
$$
p_G(m+n)=|P_{m+n}|\le|P_m\times P_n|=|P_m|\cdot|P_n|=p_G(m)\cdot p_G(n).
$$
\kduk

A~very interesting result in enumerative combinatorics is the determination 
of $\kappa(H)$ for the graph $H$ of the hexagonal lattice. $H$ is obtained by tiling the plane $\R^2$ 
by congruent regular hexagons. In 2012, 
the Russian mathematician {\em Stanislav Smirnov\index{Smirnov, Stanislav} (1970)} together with the 
French mathematician {\em Hugo Duminil-Copin\index{Duminil-Copin, Hugo} (1985)} proved 
\cite{dumi_smir} that
$${\textstyle
\kappa(H)=\sqrt{2+
\sqrt{2}}=2\cos(\pi/8)\,.
}
$$
\begin{exer}\label{ex_hexag}
Define the graph $H$ in set-theoretic terms and 
show that $H$ is a~$3$-regular {\em LFT} graph.    
\end{exer}
\cite{grim_li} is a~survey article on growth constants of graphs.

\medskip\noindent
{\em $\bullet$ Meanders. }A~\underline{matching\index{matching|emph}} is any graph $M=\langle 
V,E\rangle$ such that $V\sus\Z$ 
is finite with even cardinality and $E$ is a~partition of $V$. 
A~matching $M=\langle V,E\rangle$ is \underline{non-crossing\index{matching!non-crossing|emph}} 
if for no two edges $e,e'\in E$ we have
$$
\min(e)<\min(e')<
\max(e)<\max(e')\,.
$$

\begin{defi}[meanders]\label{def_meand}
A~\underline{meander\index{meanders|emph}} is any triple 
$$
M=\langle V,\,E,\,E'\rangle
$$ 
such that $\langle V,E\rangle$ and $\langle V,E'\rangle$ are two
non-crossing matchings on a~common vertex set $V$ such that for every partition 
$\{A,B\}$ of $V$, some edge in $E\cup E'$ joins $A$ and $B$.
\end{defi}
We can draw any meander $\langle V,E,E'\rangle$ in the plane:  the vertices
in $V$ are located on the $x$-axis and the edges in $E$, resp. $E'$, 
are represented by half-circles drawn  above, resp. below, 
the $x$-axis. The last condition in Definition~\ref{def_meand} is (equivalent to) connectivity: the 
resulting union of semicircles is a~connected and closed non-self-intersecting plane curve 
(circuit). It is more common to introduce meanders in these 
geometric terms\,---\,see, for example, the beginning of 
\cite{albe_pate,dele_al}. For $n\in\N$, we define $m(n)$\label{meaNum} to be 
the number of meanders $\langle[2n],E,E'\rangle$. 

\begin{exer}\label{ex_bouMeaCat}
Let 
$C_n:=\frac{1}{n+1}\binom{2n}{n}\le 4^n$\label{Catalan} be the $n$-th \underline{Catalan number\index{Catalan numbers|emph}}. Show that $m(n)\le C_n^2\le 16^n$.
\end{exer}

The fourth application of Fekete's lemma establishes the existence of 
the growth constant $\mu$ of meanders.  

\begin{prop}[4th application]\label{prop_meand}
The limit
$$
\mu:=\lim_{n\to\infty}m(n)^{1/n}\ \ (\in[0,\,16])
$$
exists.
\end{prop}
\duk
In view of Corollary~\ref{cor_mulFeke}, it suffices to show that $m(n)$ is
supermultiplicative. The bound $\mu\in[0,16]$ follows from Exercise~\ref{ex_bouMeaCat}.

We denote by $M_n$ the set of meanders with the vertex set $[2n]$ and define an injection 
$f\cc M_m\times M_n\to M_{m+n}$\,---\,by this we will be done. Let $M\in M_m$ and $M'\in
M_n$,  with $M=\langle[2m],E,E'\rangle$.  
We shift $M'$ by the vector $\langle 2m,0\rangle$ to the meander 
$$
M''=\langle\{2m+1,\,
2m+2,\,\ds,\,2m+2n\},\,E_0,\,
E_0'\rangle\,.
$$
We take the edges $e_1$ and $e_2$
such that $2m\in e_1\in E$ and $2m+1\in e_2\in E_0$, and define new edges $e_3=(e_1\cup 
e_2)\setminus\{2m,2m+1\}$ and $e_4=\{2m,2m+1\}$. We set 
$E_1=((E\cup E_0)\setminus\{e_1,e_2\})\cup\{e_3,e_4\}$. It 
follows that
$$
N=\langle[2m+2n],\,E_1,\,
E'\cup E_0'\rangle\in M_{m+n}\,.
$$
We define $f(M,M')=N$. Since $M$ and $M'$ can be recovered from $N$, the map $f$ is an injection.
\kduk

\begin{exer}\label{ex_howReco}
How does one recover $M$ and 
$M'$ from $N$?
\end{exer}

In the article
\cite{albe_pate} the best currently known bounds
$$
11.380\le\mu\le 12.901
$$
are obtained. Another interesting article on meanders is \cite{dele_al}. It is not known if 
there is an algorithm that computes
the function $n\mapsto m(n)$ in polynomially many (in $n$) bit 
operations. 

\medskip\noindent
{\em $\bullet$ Pattern-free permutations. }Let 
$m\in\N$. An \underline{$m$-permutation\index{permutations!emperm@$m$-permutation|emph}} $p=a_1a_2\ds a_m$ is any word 
over the alphabet $[m]$ of length $m$ that uses every letter 
$i\in[m]$ exactly once. 
A~\underline{permutation\index{permutations|emph}} is an $m$-permutation for some $m$.
We denote the set of $m$-permutations by $S_m$.

\begin{exer}\label{ex_numsPerms}
Show that $|S_m|=m!=\prod_{i=1}^m i$.     
\end{exer}

We introduce a~containment
for permutations. It is similar to the containment of words. If $q=b_1b_2\ds 
b_n$ is an $n$-permutation, we write $p\preceq q$ and say that $p$
is \underline{contained\index{permutations!containment of|emph}} in $q$ if 
there exist $m$ indices $1\le i_1<i_2<\ds<i_m\le n$ such that for 
every $l,l'\in [m]$ we have
$$
b_{i_l}<b_{i_{l'}}\iff
a_l<a_{l'}\,.
$$
In other words, the word $q$ contains a~subsequence 
with the same comparison pattern as $p$. For $n\in\N$ and an 
$m$-permutation $p$, we define $\pi(p,n)$ to be 
the number of $n$-permutations not
containing $p$, that is,
$$
\pi(p,\,n)=|\{q\in S_n\cc\;p\not\preceq q\}|\,.\label{piePn}
$$

In the fifth application of Fekete's lemma, we prove the 
existence of growth constants for pattern-free permutations.
Finiteness of these constants follows from the Marcus--Tardos Theorem~\ref{thm_MarcTard} below.

\begin{prop}[5th application]\label{prop_groConPer}
For every permutation $p$, the finite limit
$$
\pi(p):=\lim_{n\to\infty}\pi(p,\,n)^{1/n}\ \ (\in[0,\,+\infty))
$$
exists.    
\end{prop}
\duk
In view of Corollary~\ref{cor_mulFeke}, it suffices to show that $\pi(p,n)$
is supermultiplicative. The bound $\pi(p)<+\infty$ follows from Theorem~\ref{thm_MarcTard}. 

We denote by $\mathrm{per}(p,n)$ the set of $n$-permutations avoiding the permutation $p$ and define an injection 
$$
f\cc \mathrm{per}(p,m)\times \mathrm{per}(p,n)\to \mathrm{per}(p,m+n) 
$$
---\,by this we will be done. A~permutation $q=b_1b_2\ds b_l$ is 
$\oplus$-irreducible if there is no $i$ with $1\le i<l$ such that  
$b_j<b_{j'}$ whenever $1\le j \le i<j'\le l$. For example, $312$ is 
$\oplus$-irreducible, but $2143$ is not. $\ominus$-irreducible permutations are defined similarly, by replacing the inequality
$b_j<b_{j'}$ with $b_j>b_{j'}$. It is easy to see 
that every permutation is $\oplus$-irreducible or $\ominus$-irreducible. We assume that the 
forbidden permutation $p$ is 
$\ominus$-irreducible; the other case is very similar. 

Let $m,n\in\N$, $q=a_1a_2\ds a_m\in\mathrm{per}(p,m)$ and 
$r=b_1b_2\ds b_n\in\mathrm{per}(p,n)$. Let
$$
s=c_1c_2\ds c_{m+n}\,,
$$
where $c_i=a_i+n$ for $1\le i\le m$ and $c_{i+m}=b_i$ for $1\le 
i\le n$. It follows, due to $\ominus$-irreducibility of $p$, 
that $s\in\mathrm{per}(p,m+n)$. We set $f(q,r)=s$. Since we easily 
recover $q$ and $r$ from $s$, the map $f$ is an injection.
\kduk

\noindent
It is instructive to compare this proof and that of Proposition~\ref{prop_irrWor}.

For some time, it was an open problem to show that the
value $+\infty$ cannot occur as a~permutation growth constant. In 2004 this 
was confirmed in \cite{marc_tard} by the American
mathematician {\em Adam Marcus\index{Marcus, Adam} (1979)} 
and the Hungarian mathematician and computer scientist {\em 
G\'abor Tardos\index{Tardos, G\'abor} (1964)}. 

\begin{thm}[Marcus--Tardos]\label{thm_MarcTard}
For\index{theorem!Marcus and Tardos} 
every permutation $p$ there is a~constant $c\in\N$ such 
that $\pi(p,n)\le c^n$ for every 
$n\in\N$.    
\end{thm}

We close this section with the remark that the existence of the 
(necessarily finite) permutation growth constants
for the generalization of the counting function 
$\pi(p,n)$ with more than one forbidden permutation $p$ is still an open problem. 

\begin{exer}\label{ex_whereFails}
Where does the argument in the proof of Proposition~\ref{prop_groConPer} fail in the case of several forbidden permutations?    
\end{exer}

\section[Arithmetic of limits]{Arithmetic of limits}\label{podkap_limAaritOpe}

We investigate the interplay of limits with arithmetic operations in $\R^*$.

\medskip\noindent
{\em $\bullet$ Arithmetic of limits.\underline{\index{limit of 
a~sequence!arithmetic of (AL)|(emph}} }Recall that $(a_n)$, $(b_n)$ and 
$(c_n)$ denote real sequences, that always $\ep,\de,\theta>0$ and that 
$\R^*=\R\cup\{-\infty,+\infty\}$, with elements denoted by $A$, $B$, $K$ and 
$L$. Recall the arithmetic on $\R^*$ introduced in 
Section~\ref{podkap_nekoOkolLimi}.

\begin{exer}[variants of the triangle inequality]\label{ex_obmeTroj}
For\underline{\index{triangle inequality, $\Delta$-inequality!variants|emph}} 
every $a,b\in\R$, 
$$
|a+b|\ge|a|-|b|\,\text{ and }\,|a-b|\ge|a|-|b|\,.
$$    
\end{exer}
The next theorem is the main tool for computing limits.

\begin{thm}[arithmetic of limits]\label{thm_ari_lim}
Let $(a_n),(b_n)\sus\R$\index{theorem!arithmetic of limits of sequences|emph} be sequences with limits $\lim a_n=K$ and $\lim b_n=L$. \underline{Then} 
$${\textstyle
\text{ $\lim 
(a_n+b_n)=K+L$, $\lim a_nb_n=KL$ and $\lim 
\frac{a_n}{b_n}=\frac{K}{L}$}\,, }
$$
provided that the expression on the right-hand side is not indeterminate.
\end{thm}
\duk
\underline{Sum}. Let $K,L\in\R$ and an $\ep$ be given.
For every large $n$ we have that $|a_n-K|\le\frac{\ep}{2}$ and $|b_n-L|$ $\le\frac{\ep}{2}$. By the TI, we have
$${\textstyle
|(a_n+b_n)-(K+L)|\le|a_n-K|+|b_n-L|\le\frac{\ep}{2}+\frac{\ep}{2}=\ep
}
$$
for the same $n$. Hence $a_n+b_n\to K+L$.

Let $K=L=\pm\infty$ and an $\ep$ be given. For every large $n$ the 
numbers $a_n$ and $b_n$ have the same sign as $K$ and 
$|a_n|,|b_n|\ge\frac{1}{2\ep}$. Thus for these $n$ the sum $a_n+b_n$ has 
the same sign as $K$ and $|a_n+b_n|=|a_n|+|b_n|\ge
\frac{1}{2\ep}+\frac{1}{2\ep}=
\frac{1}{\ep}$. Hence $a_n+b_n\to K+L=K=L$. 

Let $K=\pm\infty$, $L\in\R$ and an $\ep$ be given. For every large $n$ the number $a_n$ 
has the same sign as $K$, $|a_n|\ge\frac{1}{\ep}+|L|+1$ and $|b_n-L|\le1$, thus 
$|b_n|\le|L|+1$. For these $n$ the sum $a_n+b_n$ has the same sign as 
$K$ and, by Exercise~\ref{ex_obmeTroj}, $|a_n+b_n|\ge|a_n|-
|b_n|\ge\frac{1}{\ep}+|L|+1-|L|-1=\frac{1}{\ep}$. Hence $a_n+b_n\to K+L=K$. The cases 
$K\in\R$ and $L=\pm\infty$ follow from the commutativity of 
addition. 

\underline{Product}. Let $K,L\in\R$  and an $\ep\le1$ be given. For every large $n$ one has that 
$|a_n-K|\le\frac{\ep}{2|L|+1}$, thus $|a_n|\le|K|+1$, and $|b_n-L|\le\frac{\ep}{2|K|+2}$. By TI, we have
\begin{eqnarray*}
|a_nb_n-KL|&\le&|a_n|\cdot|b_n-L|+|L|\cdot|a_n-K|\\
&\le&{\textstyle
(|K|+1)\cdot\frac{\ep}{2|K|+2}+|L|\cdot\frac{\ep}{2|L|+1}
\le\frac{\ep}{2}+\frac{\ep}{2}=\ep
}
\end{eqnarray*}
for the same $n$. 
Hence $a_nb_n\to KL$.

Let $K=\pm\infty$, $L=\pm\infty$ and an $\ep$ be given. For every large $n$ the number $a_n$ has the same sign as $K$, $b_n$ as $L$ 
and $|a_n|,|b_n|\ge\frac{1}{\sqrt{\ep}}$. Thus for these $n$ the product $a_nb_n$ has the same sign as $KL$ and
$|a_nb_n|=|a_n|\cdot|b_n|\ge\frac{1}{\sqrt{\ep}}\cdot\frac{1}{\sqrt{\ep}}=\frac{1}{\ep}$. Hence $a_nb_n\to KL$.

Let $K=\pm\infty$, $L\in\R\setminus\{0\}$ ($L=0$ yields an indefinite 
expression) and let an $\ep$ be given. For every large $n$ the 
number $a_n$ has the same sign as $K$, $|a_n|\ge\frac{2}{\ep|L|}$ and $|b_n-L|\le\frac{|L|}{2}$, thus 
$|b_n|\ge\frac{|L|}{2}$. So for these $n$ the product $a_nb_n$ has the same sign as $KL$ and  $|a_nb_n|=|a_n|\cdot|b_n|\ge
\frac{2}{\ep|L|}\cdot\frac{|L|}{2}=\frac{1}{\ep}$. Hence $a_nb_n\to 
KL$. The cases $K\in\R\setminus\{0\}$ and $L=\pm\infty$ follow from the 
commutativity of multiplication.

\underline{Ratio}. Let $K\in\R$, $L\in\R\setminus\{0\}$ ($L=0$ yields an indefinite expression) and an $\ep$ be given. For every large 
$n$ it holds that 
$|a_n-K|\le\min(\{1,\frac{\ep L^2}{4(|L|+1)}\})$ and $|b_n-L|\le\min(\{1, 
\frac{\ep L^2}{4(|K|+1)},\frac{|L|}{2}\})$, thus $|a_n|\le |K|+1$, $|b_n|\le |L|+1$ and $|b_n|\ge\frac{|L|}{2}$. 
By TI, we have 
$${\textstyle
\big|\frac{a_n}{b_n}-\frac{K}{L}\big|=\big|\frac{a_nL-b_nK}{b_nL}\big|
\le\frac{|a_n|\cdot|L-b_n|+|b_n|\cdot|a_n-K|}{|b_n|\cdot|L|}\le\frac{\ep L^2/4+\ep L^2/4}{L^2/2}=\ep
}
$$
for the same $n$. 
Hence $\frac{a_n}{b_n}\to\frac{K}{L}$.

Let $K=\pm\infty$, $L\in\R\setminus\{0\}$ and an $\ep$ be given. For every large $n$ the number $a_n$ 
has the same sign as $K$, $|a_n|\ge\frac{|L|+1}{\ep}$ and $|b_n-L|\le1$, thus 
$|b_n|\le|L|+1$. So for these $n$ the ratio $\frac{a_n}{b_n}$ has the same sign as 
$\frac{K}{L}$ and $\big|\frac{a_n}{b_n}\big|=\frac{|a_n|}{|b_n|}\ge
\frac{|L|+1}{\ep(|L|+1)}=\frac{1}{\ep}$. Hence $\frac{a_n}{b_n}\to \frac{K}{L}$.

Let $K\in\R$, $L=\pm\infty$ and an $\ep$ be given. For every large $n$ one has that $|a_n-K|\le1$, thus $|a_n|\le |K|+1$, and $|b_n|\ge \frac{|K|+1}{\ep}$. For these $n$ it holds that $\big|\frac{a_n}{b_n}-0\big|=\big|\frac{a_n}{b_n}\big|=\frac{|a_n|}{|b_n|}\le\frac{|K|+1}{(|K|+1)/\ep}=\ep$. Hence $\frac{a_n}{b_n}\to\frac{K}{L}=0$.
\kduk

\noindent
If $\lim a_n=K$, $\lim b_n=L$ and $\frac{K}{L}$ is not an indefinite 
expression,  
then $L\ne0$. There are only finitely many $n$ with $b_n=0$. The 
corresponding undefined ratios $\frac{a_n}{b_n}$
may be ignored or defined arbitrarily. 

\medskip\noindent
{\em $\bullet$ Two more results on the arithmetic of limits. }The 
previous theorem does not describe the arithmetic of limits 
completely. Even if one of the two sequences does not have a limit, 
the sequence of sums, products, or ratios may still have a unique
limit. We present six cases when this happens, and leave proofs for them as an exercise.

\begin{exer}\label{ex_dodatek1}
Prove the following proposition.
\end{exer}

\begin{prop}[more on AL]\label{prop_dod1}
Let $(a_n),(b_n)\sus\R$.
The following implications hold.
\begin{enumerate}
\item If $(a_n)$ is bounded and $L=\lim b_n=\pm\infty$ \underline{then} $\lim (a_n+b_n)=L$. 
\item If $(a_n)$ is bounded and $\lim b_n=0$ \underline{then} $\lim a_n b_n=0$.
\item If $a_n\ge c>0$ for every $n\ge n_0$ and $L=\lim b_n=\pm\infty$ \underline{then} $\lim a_n b_n=L$. 
\item If $(a_n)$ is bounded and  $\lim b_n=\pm\infty$ \underline{then} $\lim\frac{a_n}{b_n}=0$.
\item If $a_n\ge c>0$ and $b_n>0$ for every $n\ge n_0$, and if $\lim\,b_n=0$ \underline{then} $\lim \frac{a_n}{b_n}=+\infty$. 
\item If $0<a_n\le c$ for every $n\ge n_0$ and $L=\lim b_n=\pm\infty$ \underline{then} $\lim\frac{b_n}{a_n}=L$.
\end{enumerate}
\end{prop}
In parts 3 and 5 we may also have $a_n\le c<0$, in part 6 we may have $c\le a_n<0$ and in part 
5 we may have $b_n<0$. It is not hard to state precisely and prove these modifications.

If $a_n\to K$, $b_n\to L$ and 
the expression $K+L$ is indefinite, then $\lim (a_n+b_n)$ is far from
being uniquely determined. Similarly for $KL$ and $\frac{K}{L}$. We make
it more precise in the next proposition. 

\begin{prop}[on indefiniteness]\label{prop_jesteAritLimi2}
The following holds.
\begin{enumerate}
\item Indefiniteness of $(+\infty)+(-\infty)$. If $a_n\to+\infty$ and $A\in\R^*$, \underline{then} $a_n+b_n\to A$ for some $(b_n)$ with $b_n\to-\infty$.
\item Indefiniteness of $(+\infty)\cdot0$. If $a_n\to+\infty$ and $A\in\R^*$, \underline{then} $a_nb_n\to A$ for some $(b_n)$ with $b_n\to0$. 
\item Indefiniteness of $\frac{0}{0}$. If $a_n\to0$, $a_n\ne 0$ and $A\in\R^*$, \underline{then} $\frac{a_n}{b_n}\to A$ for some $(b_n)$ with $b_n\to0$. 
\item Indefiniteness of $\frac{A}{0}$
with $A\ne0$. If $a_n\to A$ and $B\in\{-\infty,+\infty\}$, \underline{then} $\frac{a_n}{b_n}\to B$ for some $(b_n)$ with $b_n\to0$.
\item Indefiniteness of $\frac{\pm\infty}{\pm\infty}$. If $a_n\to+\infty$ and $A\in\R^*$, \underline{then} $\frac{a_n}{b_n}\to A$ for some $(b_n)$ with $b_n\to\pm\infty$.  \underline{\index{limit of  a~sequence!arithmetic of (AL)|)emph}}
\end{enumerate}
\end{prop}
\duk
1. If $A=+\infty$, we set $b_n=-\frac{a_n}{2}$. If $A\in\R$, we set $b_n=-a_n+A$. If $A=-\infty$, we set $b_n=-2a_n$.

2. If
$A=+\infty$, we set $b_n=(1+|a_n|)^{-1/2}$. If $A\in\R$, we set 
$b_n=\frac{A}{1+|a_n|}$. 
If $A=-\infty$, we set $b_n=-(1+|a_n|)^{-1/2}$.

3. If 
$A=+\infty$ we set $b_n=\sgn(a_n)a_n^2$. If
$A\in\R\setminus\{0\}$, we set $b_n=\frac{a_n}{A}$. 
If $A=0$, we set $b_n=\sqrt{|a_n|}$.
If $A=-\infty$ we set $b_n=
-\sgn(a_n)a_n^2$. 

4. If
$B=+\infty$, we set $b_n=1$ if $a_n=0$ and 
$b_n=\sgn(a_n)n^{-1}$ otherwise. If $B=-\infty$, we set $b_n=1$ if
$a_n=0$ and $b_n
=-\sgn(a_n)n^{-1}$ otherwise.

5. If
$A=+\infty$, we set $b_n=\sqrt{1+|a_n|}$. If 
$A\in\R\setminus\{0\}$, we set $b_n=\frac{\sqrt{1+a_n^2}}{A}$. If $A=0$, 
we set $b_n= 1+a_n^2$.
If $A=-\infty$, we set $b_n=-\sqrt{1+|a_n|}$.
\kduk

\begin{exer}\label{ex_whatPo}
What can be said in part~3 when we assume instead of 
$a_n\ne0$ for every $n$ that, contrary-wise, $a_n=0$ for infinitely many $n$?     
\end{exer}

\section[Limits of recurrent sequences]{Limits of recurrent sequences}\label{podkap_rekuPosl}

This section begins with an example of the determination of the limit of a~recurrent sequence. We 
define $f$-recurrent sequences and
make explicit the {\em method of equations} that finds their limits. We illustrate this method with several 
examples, including the sequence $(F_{n+1}/F_n)$ of ratios of consecutive Fibonacci numbers. 

\medskip\noindent
{\em $\bullet$ A~recurrent sequence. }We begin with the well known inequality between the arithmetic
and geometric mean.

\begin{exer}[AG inequality]\label{ex_AGnerov}
For every $a,b\ge0$ we have $\frac{a+b}{2}\ge\sqrt{ab}$.    
\end{exer}

\begin{prop}[a~recurrent limit]\label{prop_limReku}
Let $(a_n)\sus\Q$ be given by $a_1=1$ and
$a_n=\frac{a_{n-1}}{2}+\frac{1}{a_{n-1}}$ for $n\ge2$.
\underline{Then} $\lim a_n=\sqrt{2}$.
\end{prop}
\duk
We show that $a_2\ge a_3\ge\ds\ge0$. From this, it follows by 
Theorem~\ref{thm_O_mon2} that $(a_n)$ has a finite 
nonnegative limit. Clearly, $a_n$ is defined for every $n$. For 
$n\ge2$, we get by the AG inequality that 
$${\textstyle
a_n=\frac{a_{n-
1}}{2}+\frac{1}{a_{n-1}}\ge 2\sqrt{\frac{a_{n-1}}{2}\cdot\frac{1}{a_{n-1}}}=\sqrt{2}\,. 
}
$$
Then for $n\ge3$ we have $a_{n-1}\ge a_n$ $\iff$ $\frac{a_{n-1}}{2}\ge\frac{1}{a_{n-1}}$ $\iff$ $a_{n-1}\ge\sqrt{2}$ which is true. Thus $(a_n)$ has a~weakly decreasing tail. Let $a=\lim 
a_n$ ($\ge\sqrt{2}$). By the arithmetic of limits and limits of subsequences, 
$${\textstyle
a=\lim a_n=\frac{\lim a_{n-1}}{2}+\frac{1}{\lim a_{n-1}}=\frac{a}{2}+\frac{1}{a}\,.
}
$$
Hence $\frac{a}{2}=\frac{1}{a}$ $\iff$ $a^2=2$, and $a=\sqrt{2}$.
\kduk

\noindent
Such computations are rigorous
only when the limit is proven to exist. For 
example, the sequence given by $a_1=1$ and $a_n=-a_{n-1}$ for $n\ge2$, does \underline{not} converge to
$0$, although the equation
$L=-L$ has in $\R^*$ the only solution $L=0$. The sequence alternates, $(a_n)=(1,-1,1,-1,\ds)$, and has no limit.

\medskip\noindent
{\em $\bullet$ $f$-recurrent sequences. }What is a~recurrent 
sequence? To explain it, we have to employ functions of several 
variables.

\begin{defi}[$f$-recurrent sequences]\label{def_fRecSeq}
Let $k\in\N$, let
$M\sus\R^k$ and let $f\cc M\to\R$. A~sequence $(a_n)$ is 
\underline{$f$-recurrent\index{sequence!efrec@$f$-
recurrent|emph}} if for every $n\in\N$, 
$$
\langle a_n,\,a_{n+1},\,\ds,\,a_{n+k-1}\rangle\in M\wedge a_{n+k}=
f(a_n,\,a_{n+1},\,\ds,\,a_{n+k-1})\,.
$$    
\end{defi}
We write $f(a_n,a_{n+1},\ds,a_{n+k-1})$ instead of $f(\langle a_n,a_{n+1},\ds,a_{n+k-1}\rangle)$.

\begin{exer}\label{ex_OnFrec}
The {\em P}-recurrent sequences in 
Section~\ref{sec_P_rec} are not $f$-recurrent for any function $f$. Extend 
Definition~\ref{def_fRecSeq} so that it includes {\em P}-recurrent 
sequences. 
\end{exer}

The well known \underline{Fibonacci sequence\index{Fibonacci sequence|emph}} 
$$
(F_n):=(1,\,1,\,2,\,3,\,5,\,8,\,13,\,21,\,34\ds)
$$ 
is an $f$-recurrent sequence; $f=f(x_1,x_2)=x_1+x_2$. 

\begin{exer}\label{ex_fibon}
Prove that $\lim F_n=+\infty$    
\end{exer}

\noindent
{\em $\bullet$ The method of equations. }We are interested in the limits of $f$-recurrent sequences. They can be determined by the method we explain. 

\begin{defi}[limit fix points]\label{def_limFixP}
Let $k$, $M$, and $f$ be as in Definition~\ref{def_fRecSeq} and let $L\in\R^*$. We call $L$ 
a~\underline{limit fix point\index{limit fix point|emph}} of $f$ if for some sequence $(a_n)$ three conditions hold.     
\begin{enumerate}
\item $\lim a_n=L$. 
\item For every $n$ we have 
$\langle a_n,a_{n+1},\ds,a_{n+k-1}\rangle\in M$.
\item
$\lim_{n\to\infty} f(a_n,a_{n+1},\ds,a_{n+k-1})=L$.
\end{enumerate}
We denote the set of limit fix points of $f$ by $\mathrm{LFP}(f)$.\label{LFP} 
\end{defi}
Limit fix points of $f$ are, in a~sense, exactly the solutions $L\in\R^*$ 
of the equation $f(L,L,\ds,L)=L$.

\begin{exer}\label{ex_triviThm}
Prove the next proposition.    
\end{exer}

\begin{prop}[on LFP]\label{prop_simpleO}
Let $L\in\R^*$ and let $(a_n)$ be an $f$-recurrent 
sequence with $\lim a_n=L$. \underline{Then} $L$ is a~limit fix point of $f$.    
\end{prop}

A~function $f\cc M\to\R$, where $M\sus\R^k$, is continuous if for every point $\overline{b}=\langle 
b_1,b_2,\ds,b_k\rangle\in M$ and every $k$-tuple of real sequences 
$$
\langle(a_{n,\,1}),\,(a_{n,\,2}),\,\ds,\,(a_{n,\,k})\rangle
$$
with $\lim a_{n,i}=b_i$, $i=1,2,\ds,k$, we have 
$\lim f(a_{n,1},a_{n,2},\ds,a_{n,k})=f(\overline{b})$.

\begin{exer}\label{ex_epDeDef}
State the equivalent $\ep$-$\de$ definition of continuity of $f$ and prove the equivalence.    
\end{exer}

\begin{thm}[LFP of continuous functions]\label{thm_metEqu}
Let\index{theorem!LFP of continuous functions|emph} 
$k$ be in $\N$, let $M\sus\R^k$ and let $f\cc M\to\R$ be
continuous. \underline{Then}
$$
\{L\in\mathrm{LFP}(f)\cc\;\langle L,\,L,\,\ds,\,L\rangle\in M\}=
\{L\in\R\cc\;f(L,\,L,\,\ds,\,L)=L\}
\,.
$$
\end{thm}
\duk
Suppose that $L$ is in the set on the right-hand side. Then the 
constant sequence $(a_n)=(L,L,\ds)$ shows that $L$ is an 
element of the set on the left-hand side.

Suppose that $L$ is in the set on the left-hand side. Then there is 
a~sequence $(a_n)$ with $\lim a_n=L$ such that
$$
\langle a_n,\,a_{n+1},\,\ds,\,a_{n+k-1}\rangle\in M
$$
for every $n$ and $\lim f(a_n,a_{n+1},\ds,a_{n+k-1})=L$. Since 
$\langle L,L,\ds,L\rangle\in M$ and $f$ is continuous, the last 
limit also equals $f(L,L,\ds,L)$. Hence $L$ is an 
element of the set on the right-hand side.
\kduk

The previous proposition and theorem yield a~method for 
determining finite limits of $f$-recurrent sequences. The proof of 
the corollary is immediate.

\begin{cor}[method of equations]\label{cor_metEqu}
Let $(a_n)\sus\R$ be an $f$-recurrent sequence with $\lim 
a_n=L\in\R$, where $f\cc M\to\R$ with $M\sus\R^k$ is continuous, and let 
$\langle L,L,\ds,L\rangle\in M$. \underline{Then} $L$ is a~solution of the equation
$$
f(L,\,L,\,\ds,\,L)=L\,.
$$
\end{cor}
Infinite limits $L=\pm\infty$, as well as finite limits 
$L\in\R$ lying outside $M$ (that is, $\langle L,L,\ds,L\rangle\not\in M$), have to be handled differently. 

\medskip\noindent
{\em $\bullet$ The initial example. }Equipped with the last proposition and corollary, we 
revisit the initial example. In it we have $k=1$ and 
the continuous function
$${\textstyle
f=f(x)=x/2+1/x\cc
\R\setminus\{0\}\to\R\,.
}
$$ 
It is easy to see that $0\not\in\mathrm{LFP}(f)$ and
$\pm\infty\in\mathrm{LFP}(f)$. So by Theorem~\ref{thm_metEqu} we have 
$$
\mathrm{LFP}(f)=\{-\infty,+\infty, 
-\sqrt{2},\sqrt{2}\}\,,
$$
because the solutions of the equation
$f(L)=L$ $\iff$ $L/2+1/L=L$ are exactly $L=\pm\sqrt{2}$.
If $(a_n)$ is an $f$-recurrent sequence with $\lim a_n=L$,
then by Proposition~\ref{prop_simpleO}, $L=\pm\infty$ or $L=\pm\sqrt{2}$. 

\begin{exer}\label{ex_proGen}
Make it more precise in the next generalization of Proposition~\ref{prop_limReku}.    
\end{exer}

\begin{prop}[a~generalization]\label{prop_limRekuG}
Let $f(x)=\frac{x}{2}+\frac{1}{x}$. For every number 
$a\in\R\setminus\{0\}$ the unique $f$-recurrent sequence $(a_n)$ 
starting from $a_1=a$ has 
$$
\lim a_n=\sgn(a)\cdot\sqrt{2}\,.
$$
\end{prop}

\noindent
{\em $\bullet$ On the Fibonacci sequence. }The Fibonacci sequence $(F_n)$ was recalled above. We determine
the limit $\lim\frac{F_{n+1}}{F_n}$ by the method of equations. 

\begin{exer}\label{ex_fibo1}
Prove by induction on $n$ that $F_{n+1}F_{n-1}-F_n^2=(-1)^n$.    
\end{exer}
Let 
$\phi:=\frac{1}{2}(1+\sqrt{5})$ be the \underline{golden ratio\index{golden ratio|emph}}\label{phi} and 
$\psi:=\frac{1}{2}(1-\sqrt{5})$ ($\in[-0.6,-0.7]$) be its conjugate.

\begin{prop}[golden ratio and Fibonacci]\label{prop_fibRat}
It is true that
$$
\lim_{n\to\infty}\frac{F_{n+1}}{F_n}=\phi\,.
$$
\end{prop}
\duk
For $n\in\N$, let $g_n=\frac{F_{n+1}}{F_n}$. Since 
$g_{n+1}=\frac{F_{n+2}}{F_{n+1}}=1+\frac{F_n}{F_{n+1}}
=1+\frac{1}{g_n}$, the sequence $(g_n)$ is $f$-recurrent for the continuous function $f=f(x)=1+\frac{1}{x}$. By Exercise~\ref{ex_fibo1} ($n\ge2$),
$$
{\textstyle
g_n-g_{n-1}=\frac{F_{n+1}}{F_n}-\frac{F_n}{F_{n-1}}=
\frac{(-1)^n}{F_nF_{n-1}}\,.
}
$$
It follows that  ($n\ge3$)
$$
{\textstyle
g_n-g_{n-2}=g_n-g_{n-1}+g_{n-1}-g_{n-2}=\frac{(-1)^nF_{n-2}+(-1)^{n-1}F_n}{F_nF_{n-1}F_{n-2}}
=\frac{(-1)^{n-1}}{F_nF_{n-2}}\,.
}
$$
Thus
$$
1=g_1<g_3<g_5<\ds<\ds<g_6<g_4<g_2=2
$$
and $g_n-g_{n-1}\to0$. We see that $(g_n)$ converges because it is Cauchy. Since
$f(L)=L$ $\iff$ $L^2-L-1=0$, Corollary~\ref{cor_metEqu} implies that $\lim g_n=\phi$.
\kduk

\noindent
We can obtain the limit $F_{n+1}/F_n\to\phi$ more easily from the next theorem. The point of the proposition is to 
illustrate the use of the method of equations for sequences that are not quasi-monotone; one must then resort to the Cauchy condition. 

\begin{thm}[Binet's formula]\label{thm_binet}
The\index{theorem!Binet's formula|emph} 
$n$-th Fibonacci number
$$
F_n=\frac{\phi^n-\psi^n}{\sqrt{5}}\,.
$$
\end{thm}
\duk
It is easy to check this formula for $n=1$ and $2$. Let $G_n$ be the ratio on the right-hand side. For every $n$ we 
have
$${\textstyle
G_{n+2}-G_{n+1}-G_n=
\frac{(\phi^2-\phi-1)\phi^n-(\psi^2-\psi-1)\psi^n}{\sqrt{5}}
=\frac{0-0}{\sqrt{5}}=0
}
$$
because $(x-\phi)(x-\psi)=x^2-x-1$. Hence $F_n=G_n$ for every $n$.
\kduk

\noindent
{\em Jacques P.~M. Binet (1786–1856)\index{Binet, Jacques P.~M.}} was a French mathematician, physicist, and astronomer.

\medskip\noindent
{\em $\bullet$ Two more examples of applications of the method of equations. }In the first, we determine for $a\ge0$ the value $v(a)$ of the nested radicals
$$
v(a)=\sqrt{a+\sqrt{a+\sqrt{a+\ds}}}\ .
$$
Let $f_a=f_a(x)=\sqrt{a+x}$. We interpret $v(a)$ as $v(a)=\lim a_n$ where $(a_n)$ is the $f_a$-recurrent 
sequence starting from $a_1=0$.

\begin{prop}[nested radicals]\label{prop_nesRad}
Then $v(a)=\frac{1}{2}(1+\sqrt{1+4a})$ for $a>0$ and  $v(0)=0$. 
\end{prop}
\duk
The value $v(0)=0$ is clear. Let $a>0$. We have $k=1$, 
$M=[-a,+\infty)$ (the domain of $f_a$), and $f_a$ is 
continuous. It follows that
$${\textstyle
\mathrm{LFP}(f_a)=\{+\infty,\,r_a\}\,\text{ where }\,
r_a=\frac{1}{2}(1+\sqrt{1+4a})
}
$$
because $L=r_a$ is the only solution of the equation
$$
f_a(L)=L\iff\sqrt{a+L}=L\iff L^2-L-a=0\wedge L\ge0\,.
$$
Let $(a_n)$ be the $f_a$-recurrent sequence with $a_1=0$. Then $a_1<a_2=\sqrt{a}$. If $a_n<a_{n+1}$, then also
$$
a_{n+1}=\sqrt{a+a_n}<\sqrt{a+a_{n+1}}=a_{n+2}
$$
and we see that $(a_n)$ increases. Clearly, $a_1<r_a$.
If $a_n<r_a$, then also $a_{n+1}<r_a$ because 
$$
a_{n+1}=\sqrt{a+a_n}<r_a\iff a+a_n<r_a^2=r_a+a
\iff a_n<r_a\,.
$$
So $(a_n)$ is bounded from above by $r_a$. By Proposition~\ref{prop_O_mon1} and Corollary~\ref{cor_metEqu},
$v(a)=\lim a_n=r_a$.
\kduk

\noindent
Note that $v(1)=\phi$.

In the second example, we find the limits of generalized Fibonacci sequences. 

\begin{prop}[gen. Fibonacci sequences]\label{prop_genFibo}
Let $f=f(x_1,x_2)=x_1+x_2$, $a,b\in\R$ and let $(a_n)$ be the $f$-recurrent sequence starting from $a_1=a$ and $a_2=b$; for $a=b=1$ we get $(a_n)=(F_n)$. The following holds.
\begin{enumerate}
\item If $b<-a/\phi$ \underline{then} $\lim a_n=-\infty$.
\item If $b=-a/\phi$ \underline{then} $\lim a_n=0$.
\item If $b>-a/\phi$ \underline{then} $\lim a_n=+\infty$.
\end{enumerate}
\end{prop}
\duk
We have $k=2$, $M=\R^2$ (the domain of $f$), and $f$
is continuous. It follows that 
$$
\mathrm{LFP}(f)=\{-\infty,\,0,\,+\infty\}
$$
because $f(L,L)=L$ $\iff$ $2L=L$ has the only solution 
$L=0$. We determine for which $a$ and $b$ each of these three possible limits occurs. 

We compute
$$
(a_n)=(a,\,b,\,a+b,\,a+2b,\,2a+3b,\,3a+5b,\,5a+8b,\,\ds)\,.
$$
So for $n\ge3$ induction gives that $a_n=aF_{n-2}+bF_{n-1}$, where $(F_n)$ is the Fibonacci sequence. We use the expression 
$${\textstyle
a_n=F_{n-2}\cdot\big(a+\frac{F_{n-1}}{F_{n-2}}\cdot b\big)\,.
}
$$
Since $\frac{F_{n-1}}{F_{n-2}}\to\phi$ by Proposition~\ref{prop_fibRat} and
$\lim F_{n-2}=+\infty$ by Exercise~\ref{ex_fibon}, we get the limits in cases~1 and~3. Suppose that case~2 occurs and $a+\phi b=0$.
Using Theorem~\ref{thm_binet} and item~1 of Proposition~\ref{prop_geoPosl}, we see that ($n\ge3$) 
$$
{\textstyle
a_n=\frac{1}{\sqrt{5}}\big(
a\phi^{n-2}-a\psi^{n-2}+b\phi^{n-1}-b\psi^{n-1}\big)=
-\frac{1}{\sqrt{5}}\big(a\psi^{n-2}+b\psi^{n-1}\big)\to0
}
$$
because $|\psi|<1$. 
\kduk

\begin{exer}\label{ex_uloNaRek}
$$
\cfrac{1}{1+\cfrac{1}{1+\cfrac{1}{1+\cfrac{1}{\ddots}}}}=\,?
$$    
\end{exer}

\section[Limits and order]{Limits and order}\label{podkap_limUsp}
 
We\underline{\index{limit of 
a~sequence!versus order|(emph}} 
investigate interactions between limits and the linear order $\langle\R^*,<\rangle$. 

\medskip\noindent
{\em $\bullet$ Strengthening a~standard theorem. }If we can compare terms in two sequences,
we can compare the limits, and vice versa. However, which terms in the sequences
are being compared? 

\begin{thm}[limits versus order]\label{thm_limAuspo}
Let\index{theorem!limits versus order 1|emph} $\lim a_n=K$ and $\lim b_n=L$. The following holds.
 \begin{enumerate}
\item If $K<L$ \underline{then} $a_m<b_n$ for every two indices $m,n\ge n_0$.
\item If for every index $n_0$ there exist indices $m,n\ge n_0$ such that $a_m\ge b_n$, \underline{then} $K\ge L$.
 \end{enumerate}
\end{thm}
\duk
1. Let $K<L$. By Exercise~\ref{ex_ulohaNaokoli1} there is an $\ep$ such that $U(K,\ep)<U(L,\ep)$. By the 
definition of limits, $a_m\in U(K,\ep)$ and $b_n\in 
U(L,\ep)$ for every $m,n\ge n_0$. So $m,n\ge n_0$ $\Rightarrow$ 
$a_m<b_n$.

2. The implication $\varphi\Rightarrow\psi$ 
is equivalent to the contrapositive 
$\neg\psi\Rightarrow\neg\varphi$. The contrapositive of the implication in part~1 
is the implication in part~2.
\kduk

\noindent
Strangely, Theorem~\ref{thm_limAuspo} is 
always presented in unnecessarily weak forms. Part~1 as: 
if $K<L$ then there is an $n_0$ such that $a_n<b_n$ for every $n\ge n_0$. Part~2 as: if $a_n\le b_n$ for every $n\ge n_0$, then
$K\le L$. I was teaching these weakish variants of 
Theorem~\ref{thm_limAuspo} for many years.

\begin{exer}\label{ex_procJesiln}
Why is the first part of Theorem~\ref{thm_limAuspo} with possibly 
$m\ne n$ stronger than the standard version with $m=n$?    
\end{exer}

Strict inequalities may not be preserved in limits; they may 
turn into equalities. This is another reason why non-strict 
equalities are
safer than strict ones.\index{less or equal@$\le$ is safer than $<$}

\begin{exer}\label{ex_muzePrejit} 
Find convergent sequences $(a_n)$ and $(b_n)$ such that $a_m<b_n$ for all $m$ and $n$, but $\lim a_n=\lim b_n$.   
\end{exer}

\begin{exer}\label{ex_daleZesil}
Prove the next strengthening of Theorem~\ref{thm_limAuspo}. State the corresponding part~2.    
\end{exer}

\begin{prop}[a~strengthening]\label{prop_limAuspo}
\index{theorem!limits versus order 1!strengthening} Let 
$\lim a_n=K$, $\lim b_n=L$, and let $K<L$. \underline{Then} $a_m\le 
a<b\le b_n$ for every $m,n\ge n_0$ for some index $n_0$ and some real numbers $a<b$. 
\end{prop}

\noindent
{\em $\bullet$ Intervals. }We remind the notion of real 
intervals.

\begin{defi}[intervals]\label{def_interv}
An \underline{interval\index{interval|emph}}, or 
a~real \underline{convex\index{convex set in $\R$|emph}} set, is any set 
$I\sus\R$ such that if $a<b<c$ with $a,c\in I$, then always $b\in I$. An interval is 
\underline{nontrivial\index{interval!nontrivial|emph}} if it has at least two
elements.   
\end{defi}
In the next proposition we work in the linear order $\langle\R,<\rangle$.
\begin{prop}[on intervals]\label{prop_intervaly}
Intervals are exactly the sets
$\emptyset$, $\{a\}$, $\R$, $(a,b)$, 
$(-\infty,a)$, $(a,+\infty)$, $(a,b]$, $[a,b)$, $[a,b]$, $(-\infty,a]$ and  
$[a,+\infty)$, where $a<b$ are real numbers.
\end{prop}
\duk
The transitivity of $<$ shows that all stated sets are convex. 
We show that no other
real convex sets exist. Let $X\sus\R$
be a~convex set different from
$\emptyset$, $\R$ and $\{a\}$, and let $a\in\R\setminus X$. 
The convexity of $X$ implies that $a\in H(X)$ (upper bounds of $X$) or $a\in 
D(X)$ (lower 
bounds of $X$). We discuss only the former case because the latter 
case reduces to the former by the reversal of inequalities.

So let $H(X)\ne\emptyset$. We set $b=\sup(X)$. Let $D(X)=\emptyset$. If $b\in X$ then $X=(-\infty,b]$. If $b\not\in X$ then $X=(-\infty,b)$. Let $D(X)\ne\emptyset$. Then we set $c=\inf(X)$, clearly $c<b$. If $b\not\in X$ and $c\not\in X$ then $X=(c,b)$. If $b\not\in X$ and $c\in X$ then $X=[c,b)$. If $b\in X$ and $c\not\in X$ then $X=(c,b]$. Finally, if $b\in X$ and $c\in X$ then $X=[c,b]$. 
\kduk
\vspace{-3mm}
\begin{exer}\label{ex_konecInterv}
Are there nonempty finite intervals?  
\end{exer}

\noindent
{\em $\bullet$ The squeeze theorem. }In Czech textbooks, it 
is the {\em two-cops-theorem}. We strengthen it by 
employing just one cop.

\begin{thm}[one cop]\label{thm_dvaStraz}
Let\index{theorem!one cop|emph}
$\lim a_n=b\in\R$ and let $(b_n)\sus\R^n$ be such that
$\lim|a_n-b_n|=0$. \underline{Then} also
$\lim b_n=b$. 
\end{thm}
\duk
Let $\ep$ be given. Then $a_n\in U(b,\frac{\ep}{2})$ and
$|a_n-b_n|\le\frac{\ep}{2}$  for every large $n$. By TI we have $b_n\in U(b,\ep)$ for the same 
$n$. 
Hence $\lim b_n=b$.
\kduk 

\noindent
The classical two-cops-theorem is the next corollary. We 
denote by $I(a,b)$ the interval $[a,b]$ if $a\le b$, and the 
interval $[b,a]$ if $b\le a$.

\begin{cor}[two cops]\label{cor_dvaStraz}
Let $\lim a_n=\lim b_n=b$ 
and let $(c_n)$ be such that $c_n\in I(a_n,b_n)$ for every large $n$. 
\underline{Then} also $\lim c_n=b$.    
\end{cor}
\duk
Clearly, $\lim|a_n-b_n|=0$. Since for large $n$ we have 
$$
|a_n-c_n|,\,|b_n-c_n|\le|a_n-b_n|\,,
$$
either of the sequences $(a_n)$ and $(b_n)$ can serve as a~cop
for $(c_n)$ in the previous theorem.
\kduk

\begin{exer} \label{ex_Jedenstr}
State and prove the variants of Theorem~\ref{thm_dvaStraz} with the cop going to $\pm\infty$.
\end{exer}

\chapter[Infinite series]{Infinite series}\label{chap_lim_fce}

Infinite series are an important application of limits of real 
sequences. We treat their theory in three sections. In  
Section~\ref{sec_AKrady} we 
develop a~theory of {\em absolutely convergent set series}. These are maps $h\cc 
X\to\R$ defined on at most countable sets $X$ such that for countable $X$ 
and every bijection $f\cc\N\to X$, the finite (and necessarily unique) 
sum 
$$
\lim_{n\to\infty}{\textstyle
\sum_{j=1}^n h(f(j))}\ \ (\in\R)
$$ 
exists. In Section~\ref{sec_polya} we apply this theory to generalize P\'olya's theorem. This theorem determines the limits of probabilities that 
walks in the grid graph $\Z^d$, starting at the origin, visit the 
given vertex. In 
Section~\ref{sec_classAbsc}, we specialize set series to the classical case 
$X=\N$ and relax the finiteness requirement by allowing unique sums 
$\pm\infty$. We call such series {\em 
abscon series}. The 
last Section~\ref{sec_classCondC} is devoted to {\em classical 
conditionally convergent series} $\sum_{n\in\N}a_n$ that may have non-unique finite sums. 

\section[Absolutely convergent set series]{Absolutely convergent set series}\label{sec_AKrady}

In this section, we consider series of a~general form such that the index set 
may be any at most countable set.

\medskip\noindent
{\em $\bullet$ Absolutely convergent set series. }We begin by defining set series.

\begin{defi}\label{def_setSer}
A~\underline{set series\index{set series|emph}} is any map $h\cc 
X\to\R$ such that the domain $X$ is at most countable. We use notation $\sum_{x\in X}h(x)$.      
\end{defi}
Next we define absolutely convergent set series.

\begin{defi}\label{def_AKrady}
A~set series $\sum_{x\in X}h(x)$  
\underline{absolutely converges\index{set series!absolutely convergent|emph}}
if two equivalent conditions hold.
\begin{enumerate}
\item There exists a~constant $c>0$ such that for every finite set $Y\sus 
X$ we have $\sum_{x\in Y}|h(x)|\le c$.
\item In the case that $X$ is countable, for every bijection 
$f\cc\N\to X$ there exists the finite limit
$$
\lim_{n\to\infty}{\textstyle
\sum_{j=1}^n h(f(j))\ \ (\in\R)\,.
}
$$
\end{enumerate}
\end{defi}

\begin{exer}\label{ex_setSer1}
Every finite set series absolutely converges.     
\end{exer}

\begin{exer}\label{ex_setSer2}
If $\sum_{x\in X}h(x)$ absolutely converges, then so 
does the set series $\sum_{x\in X}|h(x)|$.     
\end{exer}

Let $S=\sum_{x\in X}h(x)$ be a~set series and $Y\sus X$. 
We call $\sum_{x\in Y}h(x)$ the \underline{subseries\index{set series!subseries|emph}} of $S$.

\begin{exer}\label{ex_setSer3}
Every subseries of an absolutely convergent series absolutely converges.     
\end{exer}

\begin{thm}\label{thm_AbCoSeSe}
The two\index{Theorem!absolutely convergent set series|emph} 
conditions in Definition~\ref{def_AKrady} are equivalent.    
\end{thm}
\duk
Let $\sum_{x\in X}h(x)$ be a~set series. The implication 
$1\Rightarrow2$. Let 
condition~1 hold, $X$ be countable, and $f\cc\N\to X$ be 
a~bijection. For $n\in\N$, let 
$${\textstyle
t_n=\sum_{j=1}^n|h(f(j))|\,\text{ and }\,s_n=\sum_{j=1}^n h(f(j))\,. 
}
$$
The sequence $(t_n)$ is convergent 
because it weakly increases and, by the assumption, is bounded 
from above by $c$ 
$(t_n)$ is Cauchy (Theorem~\ref{thm_CauchyPodm}), and 
for any given $\ep$, there is $n_0$ such that for every $m\ge n\ge n_0$ 
we have
$$
{\textstyle
\sum_{j=n+1}^m|h(f(j))|\le\ep\,.
}
$$
Using TI, we obtain for the same $m$ and $n$ the bound
$$
{\textstyle
\big|\sum_{j=n+1}^m h(f(j))\big|\le
\sum_{j=n+1}^m|h(f(j))|\le\ep\,.
}
$$
The sequence $(s_n)$ is Cauchy and therefore convergent (Theorem~\ref{thm_CauchyPodm}).

The implication $\neg1\Rightarrow\neg2$. We assume
that $X$ is countable and that condition~1 does not hold. We find 
a~bijection $f\cc\N\to X$ such that 
$$
\lim_{n\to\infty}{\textstyle
\sum_{j=1}^n h(f(j))=\pm\infty\,.
}
$$
Let $f_0\cc\N\to X$ be any bijection. It follows from the assumption that 
$${\textstyle
\lim\sum_{j=1}^n|h(f_0(n))|=
+\infty\,. 
}
$$
Thus there is an injection $f_1\cc\N\to X$ such that 
$\lim\sum_{j=1}^nh(f_1(j))=
\pm\infty$; the image $f_1[\N]$ is the subset of $X$ on 
which $h$ has the same sign $+$ or $-$. Let $Y=X\setminus 
f_1[\N]$. It is clear that by inserting the elements of $Y$, one by 
one, sufficiently sparsely into the sequence $(f_1(n))$ we 
can transform $f_1$ into a~bijection $f\cc\N\to X$ such that
$$
\lim_{n\to\infty}{\textstyle
\sum_{j=1}^n h(f(j))}=
\lim_{n\to\infty}{\textstyle
\sum_{j=1}^n h(f_1(j))}=\pm\infty\,.
$$
\kduk

\noindent
{\em $\bullet$ Sums. }We define sums of absolutely convergent set series.

\begin{prop}\label{prop_sumSetSer}
Let $\sum_{x\in X}h(x)$ be absolutely convergent and $X$ be countable. 
\underline{Then} there is a~number $s\in\R$ such that for every bijection $f\cc\N\to X$, 
$$
\lim_{n\to\infty}{\textstyle
\sum_{j=1}^n h(f(j))=s}\,.
$$
\end{prop}
\duk
For the contrary, 
let $f,g\cc\N\to X$ be bijections 
and $s,t\in\R$ be numbers such that 
$${\textstyle
s=\lim_{n\to\infty}
\sum_{j=1}^n h(f(j))\ne\lim_{n\to\infty}
\sum_{j=1}^n h(g(j))=t\,.
}
$$
We obtain a~bijection $e\cc\N\to X$ such that the limit
$$
\lim_{n\to\infty}{\textstyle
\sum_{j=1}^n h(e(j))}
$$
does not exist, which contradicts the assumption of absolute convergence. Let $\ep=\frac{1}{3}|s-t|$. We 
define a~bijection $e\cc\N\to X$ and a~sequence of integers 
$0<n_1<n_2<\ds$ such that for every $i\in\N$,
$$
{\textstyle
\sum_{j=1}^{n_{2i-1}} h(e(j))\in U(s,\,\ep)\,\text{ and }\,\sum_{j=1}^{n_{2i}} h(e(j))\in 
U(t,\,\ep)\,.
}
$$
Then $e$ is as desired.  

We take $n_1\in\N$ such that 
$${\textstyle
\sum_{j=1}^{n_1}h(f(j))\in U(s,\,\ep)
}
$$ 
and set $e(m)=f(m)$ for $m\le n_1$. We take $n_2>n_1$ such that 
$${\textstyle
f[\,[n_1]\,]\sus g[\,[n_2]\,]\wedge
\sum_{j=1}^{n_2}h(g(j))\in U(t,\,\ep)\,. 
}
$$
For $m\in[n_2]\setminus[n_1]$, the values $e(m)$ run, in some order, 
through 
$g[\,[n_2]\,]\setminus f[\,[n_1]\,]$. We take $n_3>n_2$ such that 
$${\textstyle
g[\,[n_2]\,]\sus f[\,[n_3]\,]\wedge
\sum_{j=1}^{n_3}h(f(j))\in U(s,\ep)\,. 
}
$$
For $m\in[n_3]\setminus[n_2]$, the values $u(m)$ run, in some 
order, through $f[\,[n_3]\,]\setminus g[\,[n_2]\,]$. Continuing in this way indefinitely, 
we obtain the bijection $e$. 
\kduk

\begin{exer}\label{ex_cvicSumSeSer}
Let the maps $f$ and $g$ be as in the previous proof. We define a~surjection 
$e\cc\N\to X$ by interleaving 
them as 
$$
(e(n))=(f(1),\,g(1),\,f(2),\,g(2),\,\ds)
$$ 
($e$ is not injective). Does $\lim\sum_{j=1}^n h(e(j))$ exist?  
\end{exer}

\begin{defi}\label{def_SeSerSum}
Let $S=\sum_{x\in X}h(x)$ be 
absolutely convergent and $X$ be countable. The 
\underline{sum\index{set series!sum|emph}} of $S$ is the 
unique limit
$$
\lim_{n\to\infty}{\textstyle
\sum_{j=1}^n h(f(j))\,,
}
$$
for any bijection $f\cc\N\to X$. We denote the sum of $S$ again by 
$\sum_{x\in X}h(x)$, and by $h(X)$.
\end{defi}

\noindent
The definition is correct due to Pro;position~\ref{prop_sumSetSer}. 
If $X=\{x_1,x_2,\ds,x_n\}$ ($x_i\ne x_j$ for $i\ne j$) is finite and 
nonempty, we define the sum as the usual finite 
sum 
$${\textstyle
h(X)=\sum_{x\in X}h(x):=
h(x_1)+h(x_2)+\ds+h(x_n)\,.
}
$$
We set $h(\emptyset)=\sum_{x\in \emptyset}h(x):=0$. The adjective ``absolute'' in 
Definition~\ref{def_AKrady} refers not so much to absolute values but to the independence of 
sums on the order of summation.

\begin{exer}\label{ex_setSer4}
Find an absolutely convergent set series 
$S=\sum_{x\in X}h(x)$ and a~subseries $\sum_{x\in Y}h(x)$ of $S$ such that $h(X)=0$ and $h(Y)=1000$.     
\end{exer}

\noindent
{\em $\bullet$ Approximating sums. }We show that sums of 
absolutely convergent set series can be arbitrarily tightly approximated 
by finite sums. 

\begin{prop}\label{prop_aproxSou}
Let $S=\sum_{x\in X}h(x)$ be absolutely convergent. \underline{Then} for every $\ep>0$ 
there exists a~finite set $Y\sus X$, denoted by $Y(\ep,S)$, such that for every finite set $Z$ with $Y\sus Z\sus X$, we 
have
$${\textstyle
|\sum_{x\in Z}h(x)-h(X)|\le\ep\,.
}
$$
\end{prop}
\duk
Let an $\ep>0$ be given. 
For finite $X$ we take the set $Y$ according to 
Exercise~\ref{ex_triviSe}. Let $X$ be countable and $f\cc\N\to X$ be any 
bijection. Let $n_0\in\N$ be so large that 
$$
{\textstyle
\sum_{n>n_0}
|h(f(n))|\le\frac{\ep}{2}\,\text{ and }\,
\big|\sum_{n=1}^{n_0} 
h(f(n))-h(X)\big|\le\frac{\ep}{2}\,\,. 
}
$$
Let $Y=f[\,[n_0]\,]$. Then, for any finite set $Z$ with $Y\sus Z\sus X$, we have 
\begin{eqnarray*}
{\textstyle
\big|\sum_{x\in Z}h(x)-h(X)\big|}&\le&
{\textstyle
\big|\sum_{x\in Y}h(x)-h(X)\big|+
\big|\sum_{x\in Z\setminus Y}h(x)\big|}\\
&\le&{\textstyle
\big|\sum_{x\in Y}h(x)-h(X)\big|+\sum_{x\in Z\setminus Y}|h(x)|}\\
&\le&{\textstyle
\frac{\ep}{2}+\frac{\ep}{2}=\ep\,.
}
\end{eqnarray*}
The first two inequalities follow from TI. The third inequality follows 
from the definition of $Y$. The last fourth equality is trivial.
\kduk
\vspace{-3mm}
\begin{exer}\label{ex_triviSe}
How do we define $Y=Y(\ep,S)$ for finite $X$?    
\end{exer}

\noindent
{\em $\bullet$ Linear combinations of set series. }

\noindent
{\em $\bullet$ The grouping construction. }Recall the 
definition of partition in 
$S=\sum_{x\in X}h(x)$ be a~set series and $Y$ be a~partition of 
$X$. If every subseries $\sum_{x\in Z}h(x)$ for $Z\in 
Y$ absolutely converges, we call the series $S_Y=\sum_{Z\in Y}h(Z)$
the \underline{grouping\index{set series!grouping|emph}} of $S$.

\begin{thm}[sums of groupings]\label{thm_asocAKr}
Let\index{theorem!sums of groupings|emph} $S=\sum_{x\in 
X}h(x)$ be an absolutely convergent set series and $Y$ be 
a~partition of $X$. \underline{Then} the grouping 
series $S_Y$ absolutely converges and has sum $h(X)$.
\end{thm}
\duk
Note that $S_Y$ is correctly defined because $\sum_{x\in Z}h(x)$ for $Z\in Y$ is a~subseries of an absolutely convergent series $S$.
We first show that 
$S_Y$ absolutely converges. 
Let 
$${\textstyle
c\equiv\sup(\{\sum_{x\in Z}|r(x)|\cc\;\text{$Z\sus X$ and is finite}\})\ \ (\in[0,\,+\infty))
}
$$
and let $Y'=\{Z_1,\ds, Z_n\}\sus Y$ be a~finite set. For every $i\in[n]$ we use Proposition~\ref{prop_aproxSou} and take a~finite set $Z_i'\equiv Y(2^{-i},R_{Z_i})$ ($\sus Z_i$). 
We set $Z_0\equiv Z_1'\cup\ds\cup Z_n'$ ($\sus X$); this is a~disjoint union. Then
$\sum_{Z\in Y'}|S(R_Z)|$ is
$$
{\textstyle
\sum_{i=1}^n|S(R_{Z_i})-S(R_{Z_i'})+S(R_{Z_i'})|\le
\sum_{i=1}^n 2^{-i}+\sum_{x\in Z_0}|r(x)|\le 1+c\,.
}
$$
Hence $\sum_{Z\in Y}S(R_Z)\in\mathfrak{S}$.

Let an $\ep$ be given. We show that $|S(R)-S(R')|\le\ep$. We use Proposition~\ref{prop_aproxSou} 
and take finite sets $X'\equiv Y(\frac{\ep}{3},R)$ ($\sus X$)
and $Y'\equiv Y(\frac{\ep}{3},R')$ ($\sus Y$). Then we take a~finite 
superset of blocks $\{Z_1,\ds,Z_n\}$ of $Y'$, that is $Y'\sus\{Z_1,\ds,Z_n\}\sus Y$
and $n\in\N$, such that $X'\sus\bigcup_{i=1}^n Z_i$. For
every $i\in[n]$ we use Proposition~\ref{prop_aproxSou} and 
take a~finite set $Z_i'\equiv 
Y(2^{-i}\frac{\ep}{3},R_{Z_i})$ ($\sus Z_i$). Finally, for every 
$i\in[n]$ we set $Z_i''\equiv Z_i'\cup(X'\cap Z_i)$ ($\sus Z_i$) 
and take the disjoint union $X_0\equiv\bigcup_{i=1}^n Z_i''$ 
($\sus X$). Then $X_0$ is finite,  $X'\sus X_0$, $Z_i'\sus Z_i''$ and 
\begin{eqnarray*}
|S(R)-S(R')|&\le&{\textstyle|S(R)-\sum_{x\in X_0}r(x)|+
\sum_{i=1}^n|\sum_{x\in Z_i''}r(x)-S(R_{Z_i})|\,+}\\
&+&{\textstyle |\sum_{i=1}^n S(R_{Z_i})-S(R')|\le\frac{\ep}{3}+\frac{\ep}{3}+\frac{\ep}{3}=\ep\,.
}
\end{eqnarray*}
Hence $S(R)=S(R')$.
\kduk
\vspace{-3mm}
\begin{exer}\label{ex_vysvEpdel3}
Explain the last three bounds $\ds\le\frac{\ep}{3}$.    
\end{exer}

\noindent
{\em $\bullet$ Congruence of {\em AK} series. }$R=\sum_{x\in X}r(x)$ and $R'=\sum_{x\in 
Y}s(x)$ in $\mathfrak{S}$ are \underline{congruent\index{AK 
series!congruent|emph}}, in symbols $R\sim R'$,\index{congruence!on S@on $\mathfrak{S}$|emph} 
if there is a~bijection $f\cc X\to Y$ such that 
for every $x\in X$ we have $r(x)=s(f(x))$.

\begin{exer}\label{ex_congIsER}
Show that $\sim$ is an equivalence relation on $\mathfrak{S}$. (Equivalence relation on a~class is defined in the same way as on a~set.) 
\end{exer}

\begin{exer}\label{ex_congAKser}
If $R,R'\in\mathfrak{S}$ are congruent then $S(R)=S(R')$ ($\in\R$).    
\end{exer}

\noindent
{\em $\bullet$ Binary sums and products of {\em AK} series. }We 
introduce two binary 
operations on $\mathfrak{S}$.

\begin{thm}[binary sums on $\mathfrak{S}$]\label{thm}
Suppose that\index{theorem!binary sum of AK series|emph} 
$R=\sum_{x\in X}r(x)$ and $R'=\sum_{y\in Y}s(y)$ are in 
$\mathfrak{S}$. We set $Z\equiv X
\times\{0\}\cup Y\times\{1\}$ and define 
$${\textstyle
R+R'=\sum_{z\in Z}t(z) 
}
$$
by setting $t(z)\equiv r(x)$ if $z=(x,0)$, and $t(z)\equiv s(y)$ if $z=(y,1)$. \underline{Then}
$R+R'\in\mathfrak{S}$ and $S(R+R')=S(R)+S(R')$. We call $R+R'$
the \underline{binary sum\index{AK series!binary sum of, $+$|emph}} of $R$ and $R'$.
\end{thm}
\duk
First we show that $R+R'\in\mathfrak{S}$. 
Let $c$ be a~constant witnessing that 
both $R\in\mathfrak{S}$ and $R'\in\mathfrak{S}$, and let $W\sus 
Z$ be a~finite set. Then $W=(X'\times\{0\})\cup(Y'\times\{1\})$, 
where $X'\sus X$ and $Y'\sus Y$ are finite sets, and 
$$
{\textstyle
\sum_{z\in W}|t(z)|=\sum_{x\in X'}|r(x)|+\sum_{y\in Y'}|s(y)|\le c+c=2c\,.
}
$$
Hence $R+R'\in\mathfrak{S}$. 

We prove that $S(R+R')=S(R)+S(R')$.
Let $r\equiv S(R)$, $s\equiv 
S(R')$ and $t\equiv S(R+R')$, 
and let an $\ep$ be given. We show that $|t-(r+s)|\le\ep$. We use  
Proposition~\ref{prop_aproxSou} and take finite sets $X'\equiv 
Y(\frac{\ep}{3},R)$ ($\sus X$), $Y'\equiv Y(\frac{\ep}{3},S)$ ($\sus 
Y$) and $Z'\equiv Y(\frac{\ep}{3},R+S)$ ($\sus X\times\{0\}\cup 
Y\times\{1\}$). We take finite sets $X''$ and 
$Y''$ such that $X'\sus X''\sus X$, $Y'\sus Y''\sus Y$ and $Z'\sus X''\times\{0\}\cup Y''\times
\{1\}\equiv W$. Then $|t-(r+s)|$ is at most
$$
{\textstyle
|t-\sum_{z\in W}t(z)|+|\sum_{x\in X''}r(x)-r|+|\sum_{y\in Y''}s(y)-s|\le
\frac{\ep}{3}+\frac{\ep}{3}+\frac{\ep}{3}=\ep\,.
}
$$
Hence $t=r+s$.
\kduk
\vspace{-3mm}
\begin{exer}\label{ex_jesteTriEp}
Explain the last tree bounds $\ds\le\frac{\ep}{3}$.   
\end{exer}
\begin{exer}\label{ex_congAndBiSu}
If $Q\sim Q'$ and $R\sim R'$ then $Q+R\sim Q'+R'$.    
\end{exer}

\begin{thm}[products on $\mathfrak{S}$]\label{thm_SoucAKrad}
Suppose that\index{theorem!product of AK series|emph} 
$R=\sum_{x\in X}r(x)$ and $R'=\sum_{y\in Y}s(y)$ are in 
$\mathfrak{S}$. We define
$${\textstyle
R\cdot R'\equiv\sum_{(x,\,y)\in X\times Y}r(x)s(y)\,.
}
$$
\underline{Then}
$R\cdot R'\in\mathfrak{S}$ and $S(R\cdot R')=S(R)S(R')$. We call $R\cdot R'$ the \underline{product\index{AK series!product of, $\cdot$|emph}} of $R$ and $R'$.
\end{thm}
\duk
First we show that $R\cdot R'\in\mathfrak{S}$. We take a~constant $c$ 
witnessing that $R$ and $R'$ are AK series. Let $Z\sus X\times Y$ be 
a~finite set. We take finite sets $X'\sus X$ and $Y'\sus Y$ such that $Z\sus X'\times Y'$. Then
$$
{\textstyle
\sum_{(x,\,y)\in Z}|r(x)s(y)|\le
\sum_{x\in X'}|r(x)|\cdot
\sum_{y\in Y'}|s(y)|\le c\cdot c=c^2\,.
}
$$
Hence $R\cdot R'\in\mathfrak{S}$. 

We prove that $S(R\cdot R')=S(R)S(R')$. Let $r\equiv S(R)$, $s\equiv S(R')$ and $t\equiv S(R\cdot R')$, 
and let an $\ep\le1$ be given. We show that $|t-rs|\le\ep$. We use
Proposition~\ref{prop_aproxSou} and take finite sets $X'\equiv 
Y(\frac{\ep}{3(|s|+1)},R)$ ($\sus X$), $Y'\equiv Y(\frac{\ep}{3(|r|+1)},S)$ ($\sus Y$) and 
$Z\equiv Y(\frac{\ep}{3},R\cdot R')$ ($\sus X\times Y$). We take finite sets $X''$ and 
$Y''$ such that $X'\sus X''\sus X$, $Y'\sus Y''\sus Y$ and $Z\sus X''\times Y''$. 
Then $|t-rs|$ is at most
\begin{eqnarray*}
&&{\textstyle
|t-\sum_{(x,\,y)\in X''\times Y''}r(x)s(y)|+|\sum_{x\in X''}r(x)\cdot\sum_{y\in Y''}s(y)-rs|\le
}\\
&&{\textstyle\le\frac{\ep}{3}+|(r+\de)(s+\theta)-rs|\,\text{ where }\,|\de|\le\frac{\ep}{3(|s|+1)}\,\text{ and }\,|\theta|
\le\frac{\ep}{3(|r|+1)}\,.
}
\end{eqnarray*}
Hence $|t-rs|\le\frac{\ep}{3}+\frac{\ep}{3}+\frac{\ep}{3}=\ep$. This holds for every $\ep$ and $t=rs$.
\kduk

\begin{exer}\label{ex_jesteeTriEp}
Explain the last tree bounds $\ds\le\frac{\ep}{3}$.   
\end{exer}

\begin{exer}\label{ex_congAndProd}
If $Q\sim Q'$ and $R\sim R'$ then $Q\cdot R\sim Q'\cdot R'$.    
\end{exer}

\noindent
{\em $\bullet$ A~semiring of factorized {\em AK} series. }We show 
that AK series, when factorized by 
the congruence $\sim$, 
form a~semiring with respect to binary sums and to products. We set 
$\mathfrak{T}
\equiv\mathfrak{S}/\!\sim$
and call this class \underline{factorized AK 
series\index{AK series!factorized|emph}\index{class!t@
$\mathfrak{T}$|emph}}. A~\underline{semiring\index{semiring|emph}} is the ring structure on 
a~set or a~class, with the existence of additive inverses dropped. Let $0_{\mathfrak{T}}\equiv\emptyset$ be 
the empty AK series and $1_{\mathfrak{T}}\equiv\big[\sum_{x\in
\{1\}}1\big]_{\sim}$ be the AK series with just a~single summand $1$. 

\begin{thm}[semiring $\mathfrak{T}_{\mathrm{SR}}$]\label{thm_vlOperaci}
The structure
$$
\mathfrak{T}_{\mathrm{SR}}\equiv\langle\mathfrak{T},\,0_{\mathfrak{T}},\,1_{\mathfrak{T}},\,+,\,\cdot\rangle\label{tsr}
$$
is a~semiring. In more detail, $+$ and $\cdot$ are commutative and associative operations on $\mathfrak{T}$, the elements
$0_{\mathfrak{T}}$ and $1_{\mathfrak{T}}$ of $\mathfrak{T}$ are neutral to $+$ 
and $\cdot$, respectively, and $\cdot$
is distributive to $+$.
\end{thm}
\duk
Exercises~\ref{ex_congAndBiSu} and
\ref{ex_congAndProd} show that $+$ and $\cdot$ operate on the class
$\mathfrak{T}$. First we show that $+$ is commutative. Let 
$R=\sum_{x\in X}r(x)$ and $R'=\sum_{x\in Y}s(x)$ be in $\mathfrak{S}$, and let $Z\equiv X\times
\{0\}\cup Y\times\{1\}$ and $W\equiv X\times\{1\}\cup Y\times\{0\}$. We 
take the bijection $f\cc Z\to W$ that sends $(x,0)\in Z$ to $(x,1)\in W$, and $(y,1)\in Z$ 
to $(y,0)\in W$. Then we see that for $R+R'=\sum_{z\in Z}t(z)$ and 
$R'+R=\sum_{z\in W}t(z')$ it holds for every $z\in Z$ that
$t(z)=t'(f(z))$ because $t(z)=r(x)=t'(f(z))$ if $z=(x,0)$, 
and $t(z)=s(y)=t'(f(z))$ if $z=(y,1)$. Hence $R+R'\sim R'+R$. The
similar proof of the associativity of $+$ is relegated to Exercise~\ref{ex_jakJeBij}. We leave the proofs of commutativity and associativity of $\cdot$ to respective
Exercises~\ref{ex_jakJeBij2} and \ref{ex_jakJeBij3}.

We prove the neutrality of $0_{\mathfrak{T}}$ to $+$, and $1_{\mathfrak{T}}$ to $\cdot$.
Let $R=\sum_{x\in X}r(x)$ be in $\mathfrak{S}$. Since $X\times\{0\}\cup\emptyset\times\{1\}
=X\times\{0\}$, the bijection sending $x\in X$
to $(x,0)\in X\times\{0\}$ proves that $R\sim R+0_{\mathfrak{T}}$. Similarly, the
bijection sending $x\in X$ to $(x,1)\in X\times\{1\}$
proves that $R\sim R\cdot 1_{\mathfrak{T}}$. 

Finally, we show that $\cdot$ is distributive to $+$. Let $R=\sum_{x\in X}r(x)$, 
$R'=\sum_{x\in Y}s(x)$ and $R''=\sum_{x\in Z}t(x)$ be in $\mathfrak{S}$, and let 
$$
W\equiv X\times(Y\times\{0\}\cup 
Z\times\{1\})\,\text{ and }\,W'\equiv(X\times Y)\times\{0\}\cup (X\times Z)\times
\{1\}\,. 
$$
We take the bijection $f\cc W\to W'$ sending $(x,(y,0))\in W$ 
to $((x,y),0)\in W'$, and $(x,(z,1))\in W$ to $((x,z),1)\in W'$. 
Let 
$${\textstyle
\sum_{w\in W}u(w)\equiv 
R\cdot(R'+R'')\,\text{ and }\,\sum_{w\in W'}u'(w)
\equiv (R\cdot R')+(R\cdot R'')\,. 
}
$$
Then for every $w\in W$ we have $u(w)=u'(f(w))$ because 
$u(w)=r(x)s(y)=u'(f(w))$ if $w=(x,(y,0))$, and $u(w)=r(x)t(z)=u'(f(w))$ if $w=(x,(z,1))$. 
Hence $R\cdot(R'+R'')\sim(R\cdot R')+
(R\cdot R'')$.
\kduk
\vspace{-3mm}
\begin{exer}\label{ex_jakJeBij}
Give the bijection proving that $R+(R'+R'')\sim(R+R')+R''$.
\end{exer}

\begin{exer}\label{ex_jakJeBij2}
Give the bijection proving that $R\cdot R'\sim R'\cdot R$.    
\end{exer}

\begin{exer}\label{ex_jakJeBij3}
Give the bijection proving that $R\cdot(R'\cdot R'')\sim(R\cdot R')\cdot R''$.    
\end{exer}

\section[An application of set series: generalizing P\'olya's theorem]{Generalizing P\'olya's theorem}\label{sec_polya}

\section[Classical abscon series]{Classical abscon series}\label{sec_classAbsc}

\noindent
{\em $\bullet$ Some definitions. }A~classical \underline{series\index{series|emph
}} is a~sequence $(a_n)\sus\R$, denoted by 
$$
{\textstyle
\sum_{n=1}^{\infty}a_n,\ 
\sum_{n\ge1}a_n,\ \sum a_n\,\text{ or }\,a_1+a_2+\ds\,.
}
$$
Terms $a_n$ of the sequence are called \underline{summands\index{series!summands|emph}} of the series. For $n\in\N$, the 
$n$-th \underline{partial sum\index{series!partial sum|emph}}
of the series is $s_n:=\sum_{j=1}^n a_j$. If $\lim s_n\in\R$, we say that the series $\sum a_n$ 
\underline{converges\index{series!convergent|emph}}. If $\sum a_n$ is a~series and $f\cc\N\to\N$
is a~bijection, we call the series $\sum a_{f(n)}$ the 
\underline{reordering\index{series!reordering|emph}} of $\sum a_n$. 
A~\underline{subseries\index{series!subseries of|emph}} of $\sum a_n$ is any series $\sum b_n$ such that $(b_n)$ is a~subsequence of $(a_n)$.

\begin{defi}[abscon series]\label{def_abscon}
A~series $\sum a_n$ is an \underline{abscon series\index{abscon series|emph}} if for 
every bijection $f\cc\N\to\N$ the limit    
$$
L=\lim_{n\to\infty}
{\textstyle
\sum_{j=1}^n a_{f(j)}\ \ (\in\R^*)
}
$$
exists and does not depend on $f$.
\end{defi}
In other words, a~series 
is abscon iff in every reordering of it the sequence of partial sums 
has the same limit.
The limit $L$ is then the \underline{sum\index{abscon 
series!sum|emph}} of the abscon series $\sum a_n$. We denote sums 
by the same symbols as series. This is standard ambiguous 
classical notation; we will make an effort to always 
indicate by the words ``series'' and ``sum'' which meaning is 
intended. Abscon series are more general than the classical 
absolutely convergent series because infinite sums are allowed. 
For the same reason, abscon series are a~strictly larger family of series 
than the absolutely convergent set series $\sum_{x\in\N}h(x)$. 

\begin{exer}\label{ex_clAbSer0}
If a~set series $\sum_{n\in\N}a_n$ is absolutely convergent, then $\sum 
a_n$ is an abscon series.    
\end{exer}

\begin{exer}\label{ex_clAbSer1}
Let $\sum a_n$ be a~series with $a_n\in\{0,1\}$. Is it an abscon series?    
\end{exer}

\noindent
{\em $\bullet$ An equivalent definition of abscon series. }We first obtain an analog of condition~1 of Definition~\ref{def_AKrady}. For a~series $\sum a_n$ we set $A=
\{n\in\N\cc\;a_n>0\}$, $B=
\{n\in\N\cc\;a_n<0\}$, and define the set series 
$$
{\textstyle
\sum^+a_n:=\sum_{n\in A}a_n\,\text{ and }\,\sum^-a_n:=\sum_{n\in B}a_n\,.
}
$$

\begin{prop}[equivalent definition]\label{prop_eqDeAbscon}
A~series $\sum a_n$ is abscon $\iff$ at least one of the set 
series $\sum^+a_n$ and $\sum^-a_n$ absolutely converges. 
\end{prop}
\duk
Let $\sum a_n$ be a~series. Implication $\Leftarrow$. Suppose 
that both $\sum^+a_n$ and $\sum^-a_n$ absolutely converge, 
and that $c\ge0$ is the constant in condition 1 of 
Definition \ref{def_AKrady} that works for both $\sum^+a_n$ and 
$\sum^-a_n$. Then for any finite set $X\sus\N$,
$$
{\textstyle
\sum_{n\in X}|a_n|=
\sum_{n\in X\cap A}|a_n|+\sum_{n\in X\cap B}|a_n|\le c+c=2c\,.
}
$$
We see that $\sum_{n\in\N}a_n$ absolutely converges and therefore 
$\sum a_n$ is abscon series by Exercise~\ref{ex_clAbSer0}. 

Suppose that $\sum^+a_n$  absolutely converges, with the corresponding constant $c\ge0$, but that $\sum^-
a_n$ does not; the other case is similar. Let $f\cc\N\to\N$ be any 
bijection and let a~$d<0$ be given. We take sufficiently large 
$n_0$ such that 
$$
{\textstyle
\sum_{n\in f[\,[n_0]\,]\cap B}a_n\le d-c\,.
}
$$
Then for every $n\ge n_0$ we have
$$
{\textstyle
\sum_{j=1}^n a_{f(j)}=
\sum_{n\in f[\,[n]\,]\cap A}a_n+
\sum_{n\in f[\,[n]\,]\cap B}a_n
\le c+(d-c)=d\,.
}
$$
We see that $\lim\sum_{j=1}^n a_{f(j)}=-\infty$.

Implication $\neg\Leftarrow\neg$. We suppose that neither $\sum^+ a_n$ nor $\sum^- a_n$ absolutely 
converges. Then we are done by Exercise~\ref{ex_clAbSer2}, 
\kduk

\begin{exer}\label{ex_clAbSer2}
If neither $\sum^+ a_n$ nor $\sum^- a_n$ absolutely 
converges, then there exist bijections $f,g\cc\N\to\N$ such that 
$${\textstyle
\lim\sum_{j=1}^n a_{f(j)}=+\infty\,\text{ and }\,\lim\sum_{j=1}^n a_{g(j)}=-\infty\,.
}
$$
\end{exer}

\noindent
{\em $\bullet$ Operations with abscon series. }

\noindent
{\em $\bullet$ Abscon series with nonnegative summands. }

\noindent
{\em $\bullet$ Abscon series with positive and negative summands. }

\section[Classical conditionally convergent series]{Classical conditionally convergent series}\label{sec_classCondC}

\noindent
{\em $\bullet$ The Riemann theorem. }We show that if a~series 
has a~non-unique finite sum, then it can be reordered to have any sum in $\R^*$.

\begin{defi}[CC series]\label{def_CCser}
Let $\sum a_n$ be a~classical series. We say that it is \underline{conditionally 
convergent\index{series!conditionally convergent|emph}}, abbreviated {\em CC}, if it is convergent but is not an abscon series.   
\end{defi}
One of the simplest examples of a~CC series is in the next 
exercise.

\begin{exer}\label{ex_CC0}
The series
$$
{\textstyle
1-1+\frac{1}{2}-\frac{1}{2}
+\frac{1}{3}-\frac{1}{3}+\ds
}
$$
 is {\em CC}. 
\end{exer}

\begin{thm}[Riemann]\label{thm_Riemannova}
Let $\sum a_n$ be a~series. The next three claims on $\sum a_n$ are 
equivalent.
\begin{enumerate}
\item $\sum a_n$ is {\em CC} series.
\item The limit $\lim a_n=0$, the sum $\sum^+ a_n=+\infty$ and the sum $\sum^- a_n=-\infty$.
\item For every element $L\in\R^*$ there is a~bijection $f\cc\N\to\N$ such that
$$
\lim_{n\to\infty}{\textstyle
\sum_{j=1}^n a_{f(j)}=L\,.
}
$$
\end{enumerate}
\end{thm}
\duk

\kduk

\medskip\noindent
Theorem~\ref{thm_RiemmanRady}, which we attributed to
Riemann\index{Riemann, Bernhard}, completely characterizes the sequences
$(a_n)\sus\R$ with the property that for every $A\in\R^*$ there is a~bijection 
$\pi\cc\N\to\N$ such that the sum
$$
{\textstyle
\sum_{n=1}^{\infty}a_{\pi(n)}=A\,.
}
$$
In this extending section we obtain an analogous characterization for infinite products, 
that is, we replace addition with multiplication.

\medskip\noindent
{\em $\bullet$ Riemannian infinite products. } Recall 
Definition~\ref{def_infPro} of infinite products. For any sequence $(a_n)\sus\R$
we denote by $(a_{z_n})$ the subsequence 
of terms $a_n$ with $a_n\in(0,1)$, and by $(a_{k_n})$ the subsequence of terms $a_n$ with $a_n\ge1$.

\begin{defi}[Riemannian infinite products]\label{def_riemInfSer}
An infinite product 
$$
\prod_{n=1}^{\infty}a_n
$$ 
is \underline{Riemannian\index{infinite product!Riemannian|emph}} if the 
following conditions hold.
\begin{enumerate}
\item We have $\lim a_n=1$, $a_n\ne0$ for every $n\in\N$ and $a_n<0$ for only finitely many $n\in\N$.
\item $\sum\log(a_{z_n})=-\infty$.
\item $\sum\log(a_{k_n})=+\infty$.
\end{enumerate}
\end{defi}

\begin{exer}\label{ex_naRiemPro}
Logarithms of positive terms in a~Riemannian infinite product form
a~Riemannian series.
\end{exer}

\noindent
{\em $\bullet$ Riemann's theorem for infinite products. }We have devised the following analog of Theorem~\ref{thm_RiemmanRady}.

\begin{thm}[getting any product]\label{thm_RiemmanProd}
Let\index{theorem!getting any product|emph}
$\prod_{n=1}^{\infty}a_n$ be an infinite product. The
following properties are equivalent.
\begin{enumerate}
\item $\prod_{n=1}^{\infty}a_n$ is Riemannian.
\item For every nonnegative element $A\in\R^*$ there exists 
a~bijection $\pi\cc\N\to\N$ such that 
$$
\prod_{n=1}^{\infty}a_{\pi(n)}=A\,,
$$
or the same holds for nonpositive elements of $\R^*$.
\end{enumerate}
There do not exist two bijections $\pi,\rho\cc\N\to\N$ such that
$$
\prod_{n=1}^{\infty}a_{\pi(n)}=\lim_{n\to\infty}\prod_{j=1}^n a_{\pi(j)}<0<
\lim_{n\to\infty}\prod_{j=1}^n a_{\rho(j)}=
\prod_{n=1}^{\infty}a_{\rho(n)}\,.
$$
\end{thm}
\duk
Implication $1\Rightarrow 2$. Let $\prod_{n=1}^{\infty}a_n$ be a~Riemannian 
infinite product. We may assume that the negative terms in it are $a_1$, $a_2$, $\ds$, $a_k$, with $k\in\N_0$. For $n\in\N$ we set
$$
b_n\equiv\log(a_{n+k})\ \ (a_{n+k}>0)\,.
$$
Let 
$$
{\textstyle
c\equiv\prod_{j=1}^k a_j\ \ (\in\R\setminus\{0\})\,,
}
$$ 
where $c\equiv1$ if $k=0$. Let an $A\in\R^*$ be given, with $A\le 0$ if $c<0$ and $A\ge 0$ if $c>0$. We set
$${\textstyle
B\equiv\log(|A|)-\log|c|\,,
}
$$
where $\log(|\pm\infty|)\equiv+\infty$
and $\log 0\equiv-\infty$. By Exercise~\ref{ex_naRiemPro} and Theorem~\ref{thm_RiemmanRady}, there is a~bijection $\rho\cc\N\to\N$ such that we have the 
sum
$${\textstyle
\sum_{n=1}^{\infty}b_{\rho(n)}=B\,.
}
$$
We define the bijection $\pi\cc\N\to\N$
by 
$$
\pi(1)=1,\,\pi(2)=2,\,\ds,\,
\pi(k)=k,\,\pi(k+1)=k+\rho(1),\,
\pi(k+2)=k+\rho(2),\,\ds\,.
$$
Using the continuity of the exponential function and its limits
$$
\lim_{x\to-\infty}\exp x=0\,\text{ and }\,\lim_{x\to+\infty}\exp x=+\infty
$$
we get that $\prod_{n=1}^{\infty}a_{\pi(n)}$ equals
$$
{\textstyle
c\lim_{n\to\infty}\exp\big(\sum_{j=1}^n b_{\rho(j)}\big)=c\exp(\log|A|-\log|c|)=
\frac{c}{|c|}\cdot|A|=A\,,
}
$$
as required. But see Exercise~\ref{ex_ButSee}.

Implication $\neg 1\Rightarrow \neg 2$. Let $\prod_{n=1}^{\infty}a_n$ be an 
infinite product that is not Riemannian.
If $a_n=0$ for some $n\in\N$ then $\prod_{n=1}^{\infty}a_{\pi(n)}=0$ for
every permutation $\pi$ of $\N$. If $a_n\ne0$ for every $n$ but $a_n<0$ for
infinitely many $n$, then for
every permutation $\pi$ of $\N$ the permuted partial products
$$
{\textstyle
\prod_{j=1}^n a_{\pi(j)}
}
$$
change sign infinitely often and cannot have any nonzero limit. If $a_n\ne0$ for every $n$ and $\lim a_n$ does not exist or is not $1$, then this persists  for 
any reordering of $(a_n)$ and 
the proof of Proposition~\ref{prop_NCCprod} shows 
that the permuted partial products cannot have any nonzero limit. Suppose 
that $a_n<0$ for only finitely many $n$. If $\sum\log(a_{z_n})
>-\infty$ then no permutation $\pi$ of $\N$ gives as infinite product $0$, and 
if $\sum\log(a_{k_n})<+\infty$ then  no permutation $\pi$ of $\N$ gives as 
infinite product $\pm\infty$. 

We show that permutations $\pi$ and $\rho$ as stated do not exist. We may
assume that $a_n\ne0$ for every $n$.
If $a_n<0$ for infinitely many $n$ then, as we know, by reordering we can get 
as an infinite product only $0$. If $a_n<0$ for finitely many $n$, let $c$
be the product of the negative terms. Then for every reordering the permuted partial
products have, eventually, the same sign as $c$. Thus we cannot get two infinite 
products with different signs.
\kduk
\vspace{-3mm}
\begin{exer}\label{ex_ButSee}
Check that in the cases $A=0$ and $A=\pm\infty$ the computation in the proof is really correct.     
\end{exer}

\begin{exer}\label{ex_naRieThmPro}
Obtain an analog of Proposition~\ref{prop_ZabSouc} for infinite products.
\end{exer}

\section[Classical infinite series]{Classical infinite series}\label{sec_stanRady}

In Section~\ref{sec_AKrady} we introduced AK series. Now we give a~standard introduction to infinite series. 

\medskip\noindent
{\em $\bullet$ Series in general. }We begin with basic notions of the theory of infinite series.

\begin{defi}[infinite series]\label{def_infSer}
An (infinite) 
\underline{series\index{series|emph}} is a~real sequence $(a_n)$. We denote
it by 
$${\textstyle
\text{$\sum a_n$ or $\sum_{n=1}^{\infty}a_n$ or  $a_1+a_2+\cdots$}\,.
}
$$
The numbers $a_n$ are the 
\underline{summands\index{series!summand of|emph}} of the series. The  
\underline{sum\index{series!sum of|emph}} of $\sum a_n$ is the limit 
$$
\lim_{n\to\infty} (a_1+a_2+\ds+a_n)\ \ (\in\R^*)\,,
$$
if it exists. We denote the sum again by $\sum a_n$ or $\sum_{n=1}^{\infty}a_n$ or
$a_1+a_2+\cdots$.\label{seriesSum} 
Terms of the sequence 
$$
(s_n)\equiv(a_1+\ds+a_n)\label{partSum}
$$ 
are called 
\underline{partial sums\index{series!partial sum 
of|emph}} of the series, so that the sum is $\lim s_n$.
\end{defi} 
A~series with finite sum 
\underline{converges\index{series!converges|emph}}. If it has an infinite or no sum, it
\underline{diverges\index{series!diverges|emph}}. The notation for series generalizes the usual notation
$$
{\textstyle
\sum_{j=1}^n a_j=a_1+a_2+\ds+a_n
}
$$
for finite sums.

Let $S\equiv\sum a_n$ be a~series. For $m\in\N$, the series 
$${\textstyle
\sum_{n=m}^{\infty}a_n\,\text{ or }\,
a_m+a_{m+1}+\ds 
}
$$
is the standard series $\sum b_n$ with $b_n\equiv a_{m+n-1}$, 
$n\in\N$, and is called a~\underline{tail\index{series!tail 
of|emph}} of $S$. A~\underline{subseries\index{series!subseries of|emph}} of $S$ is any series $\sum b_n$ of the form $b_n\equiv a_{f(n)}$, $n\in\N$,
where $f\cc\N\to B$ is the ordering of an infinite set $B\sus\N$. In other words, $\sum b_n$ is a~subsequence of the sequence (series) $S$. 

Every subsequence of a~sequence with a~limit has the same limit. For
series the situation is different. 

\begin{prop}[on subseries]
For any pair of elements $A,B\in\R^*$ there exist a~series 
$\sum a_n$ and a~subseries $\sum b_n$ of it such that the sums are
$${\textstyle
\text{$\sum a_n=A$ and $\sum b_n=B$}\,. 
}
$$
\end{prop}
\duk
For $A,B\in\R$ we set $a_1\equiv 2^{-1}(A-B)$, $a_2\equiv 2^{-1}B$, 
$a_3\equiv 2^{-2}(A-B)$, $a_4\equiv 2^{-2}B$, $\dots$ and $\sum 
b_n\equiv\sum a_{2n}$. For $A=\pm\infty$ and $B\in\R$ we set 
$a_1\equiv \pm1$, $a_2\equiv 2^{-1}B$, $a_3\equiv\pm1$, 
$a_4\equiv 2^{-2}B$, $\dots$ (equal signs) and $\sum b_n\equiv\sum 
a_{2n}$. For positive $A\in\R$ and $B=-\infty$ we set 
$a_1\equiv\frac{A}{1}$, $a_2\equiv-\frac{A}{2}$, $a_3\equiv
\frac{A}{2}$, $a_4\equiv-\frac{A}{3}$, $\ds$ and
$\sum b_n\equiv\sum a_{2n}$. The remaining cases when $A\in\R$ and 
$B=\pm\infty$ are left in the next exercise. For $A=-\infty$ and $B=+\infty$ we set $a_1\equiv-2$, $a_2\equiv1$, $a_3\equiv-2$, $a_4\equiv1$, $\ds$ and
$\sum b_n\equiv\sum a_{2n}$. The remaining cases with
$A,B\in\{-\infty,+\infty\}$ are similar.
\kduk
\vspace{-3mm}
\begin{exer}\label{ex_remaiCase}
Resolve the remaining three cases with $A\in\R$ and $B=\pm\infty$, namely $A>0$, $B=+\infty$
and $A\le0$, $B=\pm\infty$.
\end{exer}

\begin{exer}\label{ex_robKonvDiv}
Convergence of a~series is a~robust property of sequences.    
\end{exer}
However, the sum is sensitive to changes of summands.
\begin{exer}\label{ex_zmeSou}
In every convergent series any change of any single summand changes the sum.    
\end{exer}

\begin{exer}\label{ex_onSubser}
Find a~convergent series that has a~divergent subseries.    
\end{exer}

\begin{exer}\label{ex_onTails}
A~series converges iff its every tail converges. A~series has 
a~sum $\pm\infty$ iff its every tail has the same infinite sum.
\end{exer}

\begin{exer}\label{ex_nezapScitance}
Suppose that the series $S\equiv a_1+a_2+\ds$ is such that $a_n\ge0$
for every $n\ge n_0$. Then $S$ has a~sum that is not
$-\infty$. Similarly,  
every series with almost all summands non-positive has a~sum that is not $+\infty$.    
\end{exer}

\begin{exer}\label{ex_souPlusneko} $\sum_{n=1}^{\infty}1=+\infty$.  
\end{exer}

Let $\sum a_n$ be a~series. We assume that $A$ is a~set such that
$$
A\sus\{n\in\N\cc\;a_n=0\}
$$
and that $f$ is the ordering of $\N\setminus A$. If $\N\setminus A$ is finite and $m=|\N\setminus A|$, then the
$A$-\underline{deletion of 
zeros\index{series!deletion of 
zeros|emph}} from $\sum a_n$ is the finite sum $\sum_{i=1}^m 
a_{f(i)}$. If $\N\setminus A$ is
infinite then it is the series $\sum b_n$ with $b_n=a_{f(n)}$. We show that this operation affect neither existence nor value of the sum. We begin with the finite case.

\begin{prop}[deletion of zeros~1]\label{prop_delZeros1}
Let $\sum a_n$ be a~series and let the finite sum $\sum_{i=1}^m a_{f(i)}$ arise from $\sum a_n$ by $A$-deletion of zeros. \underline{Then} we have the equality of sums
$$
{\textstyle
\sum a_n=\sum_{i=1}^m a_{f(i)}\,.
}
$$
\end{prop}
\duk
The sequence $(s_n)$ of partial sums of the series $\sum a_n$ is evenually constantly equal to $\sum_{i=1}^m a_{f(i)}$. 
\kduk
 
\begin{prop}[deletion of zeros~2]\label{prop_delZeros2}
Let $\sum a_n$ be a~series and let the series $\sum b_n$ arise from $\sum a_n$ by $A$-deletion of zeros. \underline{Then} we have the equality of sums
$$
{\textstyle
\sum a_n=\sum b_n\,,
}
$$
whenever one sum exists.
\end{prop}
\duk
Let $\sum a_n$ and $\sum b_n$ be as stated, $f$ be the ordering of $\N\setminus A$ and $(s_n)$ be partial sums of $\sum a_n$. Then for every $n$ we have 
$${\textstyle
\sum_{j=1}^n b_j=
\sum_{j=1}^n a_{f(j)}=s_{f(n)}\,, 
}
$$
and $s_m=s_{f(n)}$ for every $m$ such
that $f(n)\le m<f(n+1)$. Thus, by Exercise~\ref{ex_blowUp}, 
$$
{\textstyle
\sum b_n=\lim_{n\to\infty}\sum_{j=1}^n b_j=
\lim_{m\to\infty}s_m=\sum a_n\,,
}
$$
whenever either limit exists.
\kduk

Let $\sum a_n$ be a~series. Its \underline{reordering\index{series!reordering of|emph}} is any series $\sum b_n$ such that $b_n=a_{f(n)}$ for a~bijection $f\cc\N\to\N$.

\begin{prop}[commutativity]\label{prop_konZnam}
Let $\sum a_n$ be a~series. If 
$\{n\cc\;a_n<0\}$ or
$\{n\cc\;a_n>0\}$ is finite, \underline{then} all reorderings of 
$\sum a_n$ have the same sum.
\end{prop}
\duk
Suppose that $I=\{n\cc\;a_n>0\}$ is finite. Let $f,g\cc\N\to\N$ be any 
bijections and let 
$${\textstyle
s_n\equiv\sum_{i=1}^n a_{f(i)}\,\text{ and }\,
t_n\equiv\sum_{i=1}^n a_{g(i)}\,.
}
$$
We take $m\in\N$ such that $f[\,[m]\,]\supset I$ and $g[\,[m]\,]
\supset I$. Then both $(s_n)$ and $(t_n)$ weakly decrease starting
from the index $m$, and by  Theorem~\ref{thm_O_mon2} both 
have a~limit. It is not hard to see that for every $n_1\ge m$
there exist two indices $n_2,n_3\ge m$ such that $t_{n_2}\le s_{n_1}$ and $s_{n_3}\le 
t_{n_1}$. It follows that $\lim s_n=\lim t_n$. The other case 
when $\sum a_n$ has only finitely many negative summands is treated similarly.
\kduk

The equality $\lim a_n=0$ is the \underline{necessary convergence condition\index{series!NCC|emph}\index{NCC|emph}}\label{NCC} (NCC) of the series $\sum a_n$.

\begin{prop}[NCC]\label{prop_nutnaPoKon} If a~series $\sum a_n$ converges \underline{then} $\lim a_n=0$.
\end{prop}
\duk
Suppose that $s\equiv\lim s_n=\lim (a_1+\ds+a_n)$ is in $\R$. 
Then $\lim a_n=\lim(s_n-s_{n-1})=\lim s_n-\lim s_{n-1}=s-s=0$.
\kduk

\noindent
Thus if $\lim a_n$ does not exist or is not $0$ then $\sum a_n$ diverges.
Exercise~\ref{ex_souPlusneko} shows
that NCC does not generalize to infinite sums.

\begin{exer}\label{ex_vzsvRovn}
Explain the last four equalities in the previous proof.    
\end{exer}

Let $\sum a_n$ and $\sum b_n$ be series and $c,d\in\R$. Their \underline{linear 
combination\index{series!linear combination of|emph}} 
$${\textstyle
c\sum a_n+d\sum b_n
}
$$ 
is the series $\sum(ca_n+db_n)$.

\begin{prop}[lin. combinations of series]\label{prop_linComSer}
Suppose that the series $\sum a_n$ and $\sum b_n$ have sums $A$ and $B$, 
respectively. \underline{Then} the linear combination 
$${\textstyle
c\sum a_n+d\sum b_n
}
$$ 
has sum $cA+dB$, if this expression is not 
indefinite.
\end{prop}
\duk
This follows at once from Definition~\ref{def_infSer} and Theorem~\ref{thm_ari_lim}.
\kduk
\vspace{-3mm}
\begin{exer}\label{ex_onLinComSer}
L.~Euler\index{Euler, Leonhard} proved that 
$\sum n^{-2}=\frac{\pi^2}{6}$. Using a~linear combination of series find the sum $\sum (-1)^{n+1}n^{-2}$. 
\end{exer}

\noindent
{\em $\bullet$ Infinite products. }We briefly mention the analog of sums of
series for the operation of multiplication of real numbers.

\begin{defi}[infinite products]\label{def_infPro}
For any sequence $(a_n)\sus\R$ we define 
$$
\prod_{n=1}^{\infty}a_n\equiv\lim_{n\to\infty}a_1a_2\ds a_n\ \ (\in\R^*)\,,
\label{infiProd}
$$
if the limit of \underline{partial products\index{partial product|emph}} exists, and call it the \underline{(infinite) product\index{infinite product|emph}} of the numbers $a_n$.
\end{defi}
As for series, $\prod_{n=1}^{\infty}a_n$ also means the
sequence $(a_n)$. If $\prod_{n=1}^{\infty}a_n\in\R$, we say 
that the infinite product \underline{converges\index{infinite 
product!converges|emph}}. We can and will consider more general infinite 
products of sequences $(z_n)\sus\C$, but in the complex domain we do not allow infinite products equal to $\pm\infty$.

\begin{exer}\label{ex_prNaNekSou}
What is $\prod_{n=1}^{\infty}\big(1+\frac{1}{n}\big)$?    
\end{exer}

We have the following NCC\index{NCC|emph} for infinite products.

\begin{prop}[NCC for infinite products]\label{prop_NCCprod}
Let $\prod_{n=1}^{\infty}a_n$ be a~convergent infinite product. \underline{Then} 
$$
\text{$a_n=0$ for some $n$ or $\liminf_{n\to\infty}|a_n|\le1$}\,.
$$
\end{prop}
\duk
If $a_n=0$ for some $n$ then
$\prod_{n=1}^{\infty}a_n=0$. Suppose that 
$$
\lim_{n\to\infty}a_1a_2\ds 
a_n=a\in\R
$$ 
and that $a_n\ne0$ for every $n$. The first case is that $a\ne0$. Then
$${\textstyle
\lim a_n=\lim\frac{a_1a_2\ds a_n}{a_1a_2\ds a_{n-1}}=
\frac{\lim a_1a_2\ds a_n}{\lim a_1a_2\ds a_{n-1}}=\frac{a}{a}=1\,,
}
$$
so that 
$$
\liminf|a_n|=\lim|a_n|=\lim a_n=1\le1\,.
$$
Let $a=0$. Then it is clear that we 
cannot have $|a_n|\le1$ for only finitely many $n$. Thus 
$\liminf|a_n|\le1$.
\kduk

\noindent
We will not hide from the reader that the standard terminology of infinite 
products regards the zero products
$${\textstyle
\prod_{n=1}^{\infty}a_n=0
}
$$
as divergent.

\begin{exer}\label{ex_statStan}
State for infinite 
products the standard {\em 
NCC\index{NCC|emph}}.    
\end{exer}

It could be an interesting project to work out analogues of results about 
series for infinite products. Here we realize it only for Riemann's 
Theorem~\ref{thm_RiemmanRady} in 
Section. 

\medskip\noindent
{\em $\bullet$ The \underline{harmonic series\index{series!harmonic|emph}} is the series $\sum
\frac{1}{n}$. }We show that it has the sum $+\infty$.

\begin{exer}\label{ex_prvPomoc}
If a~sequence $(a_n)$ weakly increases and has a~subsequence with the limit $+\infty$, then $\lim a_n=+\infty$.    
\end{exer}

\begin{exer}\label{ex_druPomoc}
If series $\sum a_n$ and $\sum 
b_n$ satisfy $a_n\ge b_n$ for every $n\ge n_0$ and $\sum 
b_n=+\infty$, then $\sum a_n=+\infty$.    
\end{exer}

\begin{prop}[$\sum\frac{1}{n}=+\infty$]\label{prop_divHarm}
The harmonic series sums to $+\infty$.
\end{prop}
\duk
We consider the series 
$${\textstyle
\sum b_n\equiv\frac{1}{2}+\frac{1}{4}+
\frac{1}{4}+\frac{1}{8}+\frac{1}{8}+\frac{1}{8}+\frac{1}{8}+\frac{1}{16}+\cdots\,, 
}
$$
where in general $b_{2^k}=b_{2^k+1}=\ds=b_{2^{k+1}-1}=
\frac{1}{2^{k+1}}$. Then 
$\frac{1}{n}\ge b_n$ for every $n$. The partial sums $(s_n)$ of $\sum 
b_n$ increase. For every $k\in\N_0$, 
$${\textstyle
s_{2^{k+1}-1}=\frac{1}{2}+2\cdot\frac{1}{4}+
4\cdot\frac{1}{8}+\ds+2^k\cdot\frac{1}{2^{k+1}}=\frac{k+1}{2}\,. 
}
$$
By Exercise~\ref{ex_prvPomoc}, 
$\sum b_n=\lim s_n=+\infty$. By Exercise~\ref{ex_druPomoc}, 
$\sum\frac{1}{n}=+\infty$.
\kduk

\noindent 
By Proposition~\ref{prop_konZnam} every reordering of the
series $\sum\frac{1}{n}$ has sum $+\infty$. The partial sums 
$${\textstyle
(h_n)\equiv\big(\sum_{i=1}^n
\frac{1}{i}\big)=\big(1,\,
\frac{3}{2},\,\frac{11}{6},\,
\frac{25}{12},\,
\frac{137}{60},\,\ds\big)\ 
\ (\sus\Q) 
}
$$
are called \underline{harmonic 
numbers\index{harmonic number,   $h_n$|emph}}.\label{hen} Already in 1350
the French medieval philosopher {\em Nicolas Oresme\index{Oresme, 
Nicolas} (1320 to 1325\,--\,1382)} proved Proposition~\ref{prop_divHarm},  that is, that $h_n\to+\infty$.

\begin{thm}[growth of $h_n$]\label{thm_harm_cial} There is a~$c\ge0$ such that for every $n$\index{theorem!asymptotics of $h_n$}, 
$$
h_n=\log n+\gamma+\Delta(n) 
$$
where $|\Delta(n)|\le\frac{c}{n}$
and $\ga=0.57721\ds$ is so called \underline{Euler's 
constant\index{Euler's constant, $\ga$|emph}}.\label{gamma}
\end{thm}
We prove this asymptotics in lecture~14 with the help of 
integrals. Theorem~\ref{thm_HarmExp}, which 
will be proven in {\em MA~1${}^+$}, provides much more precise
asymptotic estimates of $h_n$. 

\begin{exer}\label{ex_HnNenicele}
Prove that $h_n\in\N$ only for $n=1$. Hint: $m=(2l-1)2^k$. 
\end{exer}

\begin{exer}[an open problem]\label{ex_otProGamma}
Prove that
Euler's\underline{\index{Euler's constant, 
$\ga$!irrational?}} 
constant~$\ga$ is an irrational number.    
\end{exer}

\noindent
{\em $\bullet$ Riemannian series. }In the first 
lecture we encountered the series
$${\textstyle
1-1+\frac{1}{2}-\frac{1}{2}+\frac{1}{3}-\frac{1}{3}+\ds+\frac{1}{n}-\frac{1}{n}+\ds
}
$$
with the sum $0$. We described a~reordering of it with a~positive sum. We show in Theorem~\ref{thm_RiemmanRady}  that this 
and similar series can be reordered to have any sum. We begin with less extensive changes of sums.
For a~series $\sum a_n$ we denote by $k_1<k_2<\ds$ the
indices $n$ such that $a_n\ge0$, and by $z_1<z_2<\ds$ the
indices $n$ such that $a_n<0$. If there are infinitely many 
$k_n$, by Proposition~\ref{prop_konZnam}
all reorderings of the series $\sum a_{k_n}$ have the same sum. The same
holds for the infinite series $\sum a_{z_n}$.

Before we state and prove two propositions and a~theorem on reorderings of series, we
explain two mutually inverse operations related to real sequences, decomposition and composition. 
A~\underline{segment\index{segment|emph}} is 
a~$k$-tuple $U=\langle b_1,\ds,b_k\rangle$ of real numbers $b_i$; we denote its \underline{length\index{segment!length of|emph}} $k$ by $|U|$ ($\in\N$). Let $S=(a_n)$ be a~real sequence. An 
\underline{initial segment\index{segment!initial|emph}} of $S$ is any nonempty segment
$U=\langle a_1,\ds,a_k\rangle$. 
We then write $S\setminus U$ for the tail $(a_{k+1},a_{k+2},\ds)$.
A~\underline{decomposition\index{segment!
decomposition into segments|emph}} of $S$ in segments 
$$
S=U_1U_2\ds
$$ is a~sequence $(U_1,U_2,
\ds)$ of
segments such that $U_1=\langle a_1,\ds,a_{k_1}\rangle$, $U_2=\langle a_{k_1+1},
\ds,a_{k_1+k_2}\rangle$, $\ds$ ($k_i\in\N$).
Reversely, if $U_i=\langle u_{i,1},\ds,u_{i,k_i}\rangle$, $i,k_i\in\N$, are 
segments then 
the real sequence 
$$
(b_n)=U_1U_2\ds
$$ 
obtained by the
\underline{composition\index{segment!composition operation|emph}} of the segments
$U_i$ is defined by setting $b_n\equiv a_{i,\,j}$, where $k_0\equiv0$, the index $i$ is 
given by
$$
k_0+k_1+\ds+k_{i-1}<n\le k_0+k_1+\ds+k_i\,\text{ and }\, j\equiv n-k_0-k_1-\ds-k_{i-1}\,.
$$

\begin{prop}[getting sums 
$\pm\infty$]\label{prop_naSoucNek}
$\sum a_n$ can be  reordered to the sum $-\infty$ $\iff$ the sum
$\sum a_{z_n}=-\infty$. Similarly, $\sum a_n$ has a~reordering with the 
sum $+\infty$ $\iff$ the sum $\sum a_{k_n}=+\infty$.  
\end{prop}
\duk
We prove the former equivalence and leave the latter for Exercise~\ref{ex_druheTvrz}. 
Let $\sum a_{z_n}=-\infty$. In particular the number of
indices $z_n$ is infinite. We may assume that also the number of 
indices $k_n$ is infinite; if their number is finite, 
by Proposition~\ref{prop_konZnam} every reordering of $\sum a_n$ has sum $-\infty$. We define a~bijection 
$f\cc\N\to\N$ for which the sum $\sum a_{f(n)}=-\infty$. It is 
the ``limit'' $\lim_{k\to\infty}P_k$ of certain injective sequences $P_k=(m_{n,k})\sus\N$, 
$k\in\N_0$. Their terms are indices $k_n$ and $z_n$. We set $P_0\equiv(z_n)$. We take an initial segment $U_1$ of 
$P_0$ such that $\sum_{n\in U_1}a_n\le-1-a_{k_1}$. 
We insert in $P_0$ after $U_1$ the index $k_1$ and get the sequence $P_1$. We take an initial segment $U_2$ of $P_1$ such that $|U_2|>|U_1|$ and $\sum_{n\in U_2}a_n
\le-2-a_{k_2}$. We insert in $P_1$ after $U_2$ the index $k_2$ and get the sequence $P_2$. And so on. Since $|U_1|<|U_2|<\ds$, the sequences
$P_0$, $P_1$, $P_2$, $\ds$ converge, in the obvious sense, to the sought 
for bijection $f$.

Suppose that it is not true that $\sum a_{z_n}=-\infty$. Then either 
the number of indices $z_n$ is finite or the sum 
$\sum a_{z_n}\in\R$. In the former case by 
Proposition~\ref{prop_konZnam} all reorderings of $\sum a_n$ have the 
same sum different from $-\infty$. In the latter case there exists a~$c$ 
such that for every $n$ it holds that $\sum_{i=1}^n a_{z_i}\ge c$. It 
follows that no reordering of $\sum a_n$ has the sum $-\infty$. 
\kduk
\vspace{-3mm}
\begin{exer}\label{ex_druheTvrz}
Reduce the latter equivalence to the former. 
\end{exer}

\begin{prop}[getting no sum]\label{prop_ZabSouc}
A~series $\sum a_n$ has a~reordering
 with no sum $\iff$ the sum $\sum a_{z_n}=-\infty$ and the sum $\sum a_{k_n}=+\infty$.  
\end{prop}
\duk
We suppose that $\sum a_{z_n}=-\infty$ and $\sum a_{k_n}=+\infty$. 
We take an initial segment $U_1$ of $(z_n)$ such that $\sum_{n\in 
U_1}a_n\le-1$, and an initial segment $V_1$ of $(k_n)$ such that $\sum_{n\in 
U_1\cup V_1}a_n\ge1$. We take an initial segment $U_2$ of $(z_n)\setminus U_1$ such that 
$${\textstyle
\sum_{n\in U_1\cup V_1\cup  U_2}a_n\le-1\,, 
}
$$
and an initial segment $V_2$ of $(k_n)\setminus V_1$ such that 
$${\textstyle
\sum_{n\in U_1\cup V_1\cup  
U_2\cup V_2}a_n\ge1\,. 
}
$$
We continue in this way indefinitely. The composed sequence
$$
(p_n)\equiv U_1V_1U_2V_2\ds\ \ (\sus\N)
$$
is a~bijection from $\N$ to $\N$ and the series $\sum a_{p_n}$
does not have a~sum because both $\sum_{i=1}^n a_{p_i}\le-1$ and
$\sum_{i=1}^n a_{p_i}\ge1$ hold for infinitely many $n$.

Suppose, for example, that the sum $\sum a_{z_n}\in\R$. For $\sum 
a_{k_n}\in\R$ we argue similarly. If $\sum a_{k_n}=+\infty$, then it
follows that every reordering of $\sum a_n$ has the sum $+\infty$. If
also $\sum a_{k_n}\in\R$, then the series $\sum a_n$ is, in the 
terminology introduced below, an abscon series, and by 
Proposition~\ref{prop_oAbskonR} all reorderings of it have the same finite sum.
\kduk

A~series $\sum a_n$ is \underline{Riemannian\index{series!Riemannian|emph}} if $\lim a_n=0$, the sum $\sum a_{k_n}=+\infty$ and the 
sum $\sum a_{z_n}=-\infty$.

\begin{thm}[getting any sum]\label{thm_RiemmanRady}
Let\index{theorem!getting any sum|emph}
$\sum a_n$ be a~series. The
following claims are equivalent.
\begin{enumerate}
\item $\sum a_n$ is Riemannian.
\item For every 
$A\in\R^*$ there exists a~reordering of $\sum a_n$ with the sum $A$.
\end{enumerate}
\end{thm}
\duk
Implication $1\Rightarrow2$. Let $\sum a_n$ be Riemannian. 
The two cases $A=\pm\infty$ are dealt with in Proposition~\ref{prop_naSoucNek}.
The case of no sum is dealt with in Proposition~\ref{prop_ZabSouc}.
Let $A\in\R$. Both sequences $(k_n)$ and $(z_n)$ are infinite and every 
index $n$ lies in exactly one of them. We define two decompositions
$$
(k_n)=U_1U_2\ds\,\text{ and }\,
(z_n)=V_1V_2\ds
$$
in segments $U_i$ 
and $V_i$. $U_1$ is the shortest initial segment of $(k_n)$ such that $\sum_{n\in U_1}
a_n\ge A$. $V_1$ is the
shortest initial segment of $(z_n)$
such that 
$${\textstyle
\sum_{n\in U_1\cup V_1}a_n\le A\,.
}
$$
$U_2$ is the 
shortest initial segment of $(k_n)\setminus U_1$
such that 
$${\textstyle
\sum_{n\in U_1\cup V_1\cup U_2}a_n\ge A\,. 
}
$$
$V_2$ is the 
shortest initial segment of $(z_n)\setminus V_1$
such that 
$${\textstyle
\sum_{n\in U_1\cup V_1\cup U_2\cup V_2}a_n\le A\,. 
}
$$
We continue in this way indefinitely. We show that the composed sequence
$$
(p_n)\equiv U_1V_1U_2V_2U_3\ds\ \ (\sus\N)
$$
is a~bijection from $\N$ to $\N$ and that the sum $\sum a_{p_n}=A$.

First we remark that the definition of the segments $U_i$ and $V_i$ is correct because
every tail of $\sum a_{k_n}$ (resp. $\sum a_{z_n}$) has sum $+\infty$ (resp. $-\infty$), see Exercise~\ref{ex_onTails}.
It is also clear that $p\cc\N\to\N$ is a~bijection because 
$$
\{(k_1,k_2,\ds),\,
(z_1,z_2,\ds)\}
$$ 
is a~partition of $\N$ in two
infinite sets. We show that $\sum a_{p_n}=A$. Let $k_i\equiv|U_i|$ and
$l_i\equiv|V_i|$, $i\in\N$, and let $u_i$ (resp. $v_i$) be the last term of $U_i$ (resp. 
$V_i$). For $i\ge2$ we set 
$S_i\equiv\sum_{j=1}^{i-1}(k_j+l_j)$. Let $s_n\equiv\sum_{j=1}^n a_{p_j}$.
It follows from the definition of $U_i$ and $V_i$ that for every $i\ge2$ we have
$$
S_i<n<S_i+k_i\Rightarrow A+a_{v_{i-1}}\le 
s_n<A,\ n=S_i+k_i\Rightarrow A\le s_n<A+a_{u_i}
$$
and similarly
$$
S_i+k_i<n<S_{i+1}\Rightarrow A\le s_n<A+a_{u_i},\ n=S_{i+1} \Rightarrow A+a_{v_i}<s_n\le A\,.
$$
But $\lim_{i\to\infty}a_{u_i}=\lim_{i\to\infty}a_{v_i}=0$, hence $\lim s_n=A$ and
$\sum a_{p_n}=A$.

Implication $\neg1\Rightarrow\neg2$. Let $\sum a_n$ be not Riemannian. If the sum
$${\textstyle
\sum_{a_{k_n}}<+\infty
}
$$ then no reordering of the series yields the sum $+\infty$. If the sum
$${\textstyle
\sum_{a_{z_n}}>-\infty
}
$$ then no reordering of the series yields the sum 
$-\infty$. If $\lim a_n$ does not exist or is not $0$, then this persists in any 
reordering and the series cannot be reordered to have a~finite sum
\kduk

\noindent
The previous theorem is due to the German mathematician {\em Bernhard Riemann\index{Riemann, 
Bernhard} (1826--1866)}. 

\medskip\noindent
{\em $\bullet$ Leibnizian series. }The series
$${\textstyle
1-1+\frac{1}{2}-\frac{1}{2}+\frac{1}{3}-\frac{1}{3}+\ds+\frac{1}{n}-\frac{1}{n}+\ds
}
$$
in the first lecture is not only Riemannian, but it is also a~\underline{Leibnizian} series\index{series!Leibnizian|emph}.
These are series of the form $\sum (-1)^{n-1}a_n$, where $a_1\ge a_2\ge\ds\ge0$ and $\lim a_n=0$.

\begin{thm}[on Leibnizian series]\label{thm_altSer}
Every\index{theorem!bounding sums of Leibnizian series|emph} 
Leibnizian series 
$${\textstyle
\sum (-1)^{n-1}a_n
}
$$ 
has a~finite sum
$s\in\R$. Moreover, for every $n$ we have bounds
$$
a_1-a_2+a_3-\ds-a_{2n}\le s\le 
a_1-a_2+a_3-\ds+a_{2n-1}\,.
$$
\end{thm}
\duk
Let $\sum (-1)^{n-1}a_n$ be Leibnizian and 
$s_n\equiv\sum_{j=1}^n(-1)^{j-1}a_j$. We show by 
induction on $k$ and $l$ that always
$$
s_1\ge s_3\ge\ds\ge s_{2k-1}\ge
s_{2l}\ge s_{2l-2}\ge\ds\ge s_2\,.
$$
we have the common limit $s\equiv\lim s_{2n-1}=\lim s_{2n}$.
By Theorem~\ref{thm_finManyBl}, $\lim s_n=s$. Since $s_{2n-1}\ge 
s_m\ge s_{2n}$ for every $m\ge 2n$, using part~2 of 
Theorem~\ref{thm_limAuspo} we get with $m\to\infty$ the stated 
inequalities for $s$. 
\kduk

\noindent
This theorem is credited to the German mathematician, philosopher and diplomat {\em 
Gottfried W. Leibniz (1646--1716)\index{Leibniz, Gottfried 
W.}}. It is a~remarkable result,  partial sums are both lower and upper bounds on the sum. 

\begin{exer}\label{ex_geneLeiSer}
Extend the previous theorem to series with Leibnizian tails.      
\end{exer}

We prove the next theorem with two
examples of sums of Leibnizian series in Chapter~\ref{chap_pr9}, see Corollaries~\ref{cor_alterHarm} and.

\begin{thm}[two sums]\label{thm_twoSums}
We\index{theorem!two sums of Leibnizian series} have sums
$$
{\textstyle
\sum\frac{(-1)^{n-1}}{n}=1-\frac{1}{2}+\frac{1}{3}-\ds=\log 2\,\text{ and }\,\sum\frac{(-1)^{n-1}}{2n-1}=1-\frac{1}{3}+\frac{1}{5}-\ds=\frac{\pi}{4}\,.
}
$$
\end{thm}
See preprints \cite{maxw_schn,schn} for summations like
$$
{\textstyle
\sum\frac{1}{n(4n-3)}=\frac{\pi}{6}+\log 2,\ 
\sum\frac{1}{n(4n-1)}=3\log 2-\frac{\pi}{2}\,\text{ and }\, 
\sum\frac{1}{n(2n-1)(4n-3)}=\frac{\pi}{3}\,.
}
$$

\noindent
{\em $\bullet$ Grouping of summands. }This is a~simple and useful operation on series. If 
$\sum a_n$ is a~series and $S\equiv(m_n)\sus\N$ is a~sequence, we set $m_0\equiv0$, define the series $\sum b_n$ by
$${\textstyle
b_n\equiv\sum_{j=m_0+m_1+\ds+m_{n-1}+1}
^{m_0+m_1+\ds+m_{n-1}+m_n}a_j
}
$$
and call it the $S$-\underline{grouping\index{series!grouping of|emph}} of $\sum a_n$.

\begin{thm}[grouping of series]\label{thm_groupSum}
Suppose\index{theorem!grouping of series|emph} 
that $\sum a_n$ is a~series, 
$$
S\equiv(m_n)\sus\N
$$ 
is a~sequence and that $\sum b_n$ is the $S$-grouping of $\sum a_n$. The
following holds.
\begin{enumerate}
\item If the sum $\sum a_n$ exists, \underline{then} the equality of sums 
$\sum b_n=\sum a_n$ holds.
\item If $\lim a_n=0$, the sequence $S$ is bounded and the sum $\sum b_n$ exists, \underline{then} the equality of sums 
$\sum a_n=\sum b_n$ holds.
\item If $\sum a_n$ and $S$ are bounded sequences and the infinite sum $\sum b_n$ 
exists, 
\underline{then} the equality of sums 
$\sum a_n=\sum b_n$ holds.
\end{enumerate}
\end{thm}
\duk
1. The sequence of partial sums of
$\sum b_n$ is always a~subsequence of the sequence of partial sums of 
$\sum a_n$.

2. We suppose that $\lim a_n=0$, $m_n\le m$ ($\in\N$) for every $n$ 
and that the sum $\sum b_n=s$ $\in\R$. The case that 
$\sum b_n=\pm\infty$ is treated in part~3; the assumptions in 
part~2 imply the assumptions for part~3. Let an $\ep$ be given. We
take an $n_0$ such that for every $n\ge n_0$ we have 
$${\textstyle
m|a_n|\le\frac{\ep}{2}\,\text{ and }\,|\sum_{j=1}^n b_j-s|\le
\frac{\ep}{2}\,.
}
$$
We set $n_1\equiv\sum_{n=1}^{n_0}m_n$. Then 
certainly $n_1\ge n_0$ and for any given $n\ge n_1$ we take the unique
$k\in\N$ with $k\ge n_0$ such that $\sum_{n=1}^k m_n\le n<
\sum_{n=1}^{k+1} m_n$. Then, with
$l\equiv m_1+m_2+\ds+m_k$, we have for every $n\ge n_1$ that
$$
{\textstyle
\big|s-\sum_{j=1}^n a_j\big|\le
\big|s-\sum_{j=1}^k b_j\big|+
\sum_{j=l+1}^n|a_j|\le
\frac{\ep}{2}+m_{k+1}\max_{l+1\le j\le n}|a_j|
}
$$
which is at most $\frac{\ep}{2}+\frac{\ep}{2}=\ep$. Hence the sum $\sum a_n=s$.

3. We suppose that $|a_n|, m_n\le m$ ($\in\N$) for every $n$ and that the
sum $\sum b_n=-\infty$; the case  $\sum b_n=+\infty$ is
treated similarly. Let a~$c$ be given. We take an $n_0$ such that 
for every $n\ge n_0$ we have $\sum_{j=1}^n b_j\le c-m^2$. We 
define $n_1$ as before, and for any given $n\ge n_1$ we define $k=k(n)$ and $l=l(n)$ as before. Then for every $n\ge n_1$,
$$
{\textstyle
\sum_{j=1}^n a_j=\sum_{j=1}^k b_j+
\sum_{j=l+1}^n a_j\le c-m^2+
m_{k+1}\max_{l+1\le j\le n}|a_j|
}
$$
which is at most $c-m^2+m^2=c$. Hence the sum $\sum a_n=-\infty$.
\kduk

\begin{cor}
We have the sums
$${\textstyle
\sum\frac{1}{(2n-1)2n}=\log 2\,\text{ and }\,
\sum\frac{1}{(4n-3)(4n-1))}=\frac{\pi}{8}\,.
}
$$
\end{cor}
\duk
We obtain these sums by applying part~2 of Theorem~\ref{thm_groupSum} with $S=(2,2,\ds)$
to the two series in Theorem~\ref{thm_twoSums}. We divide 
the latter obtained series by $2$ and use Proposition~\ref{prop_linComSer}.
\kduk

\noindent
Cancellation of summands is at work: $\frac{1}{(2n-1)2n}$ is much
smaller than any of the two summands in $\frac{1}{2n-1}-\frac{1}{2n}$. Thus we transform by grouping the Riemannian series
$\sum\frac{(-1)^{n-1}}{n}$ in the better behaved abscon series (defined next) 
$\sum\frac{1}{(2n-1)2n}$. Same for the second series $\sum\frac{(-1)^{n-1}}{2n-1}$.

\begin{exer}\label{ex_naGroup}
Show by an example that part~2 of Theorem~\ref{thm_groupSum} in
general does not hold if we omit the assumption $\lim a_n=0$. 
\end{exer}

\begin{exer}\label{ex_naGroup2}
Show by an example that part~2 of Theorem~\ref{thm_groupSum} in
general does not hold if we omit the assumption of boundedness of $S$. 
\end{exer}

\noindent 
{\em $\bullet$ Abscon series. }In Section~\ref{sec_AKrady} we considered a~more 
general version of this kind of series. Here
we treat them standardly. 

\begin{defi}[abscon series]\label{def_abskonRady}
$\sum a_n$ is an \underline{absolutely convergent} series, or 
an \underline{abscon\index{series!abscon|emph}} series, if the sum 
$${\textstyle
\sum|a_n|<+\infty\,.
}
$$
\end{defi}
Every abscon series is an AK series. In the next theorem we present an
infinite triangle inequality for 
abscon series.  

\begin{thm}[infinite triangle inequality]\label{thm_AKjeKonv} If 
\index{theorem!infinite triangle inequality|emph}\underline{\index{triangle 
inequality, $\Delta$-inequality!infinite|emph}}
$\sum a_n$ is an abscon series, \underline{then} $\sum a_n$ 
converges and
$${\textstyle
|\sum a_n|\le\sum|a_n|\,.
}
$$
\end{thm}
\duk
Let $\sum a_n$ be an abscon series with partial sums $(s_n)$ 
and let $\sum|a_n|$ have partial sums $(t_n)$. By 
Theorem~\ref{thm_CauchyPodm} the sequence $(t_n)$ is Cauchy. Hence 
for any given $\ep$ for every two large indices $m\le n$ it holds that
$$
|t_n-t_m|=\big||a_{m+1}|+|a_{m+2}|+\ds+|a_n|\big|=|a_{m+1}|+|a_{m+2}|+\ds+|a_n|\le\ep
$$
(for $m=n$ these sums are zero). By the $\Delta$-inequality we have for the same indices $m\le n$ that
$$
|s_n-s_m|=|a_{m+1}+a_{m+2}+\ds+a_n|\le|a_{m+1}|+|a_{m+2}|+\ds+|a_n|\le\ep\,.
$$
Hence $(s_n)$ is Cauchy. By Theorem~\ref{thm_CauchyPodm} the 
sequence $(s_n)$ converges. Hence the series $\sum a_n$ converges. 

By the $\Delta$-inequality it holds for every $n$ that $|s_n|\le t_n$, equivalently 
$$
-t_n\le s_n\le t_n\,. 
$$
Sending $n\to\infty$ gives by Theorem~\ref{thm_limAuspo} for sums the inequalities 
$${\textstyle
-\sum|a_n|\le\sum a_n\le\sum|a_n|\,.
}
$$
Hence $\big|\sum a_n\big|\le\sum|a_n|$.
\kduk
\vspace{-3mm}
\begin{exer}\label{ex_prerovAKrad}
Every reordering of an abscon series is an abscon series.
\end{exer}

\begin{exer}\label{ex_subserAKser}
Every subseries of an abscon series is an abscon series.
\end{exer}

\begin{prop}[on abscon series]\label{prop_oAbskonR}
$\sum a_n$ is an abscon series $\iff$ all reorderings of it have the same finite sum.     
\end{prop}
\duk
Suppose that $\sum a_n$ is an abscon series. Let $f,g\cc\N\to\N$ be 
bijections and let an $\ep$ be 
given. We take an $n_0$ such that $\sum_{n\ge n_0}|a_n|\le\ep$. Let $n_1$ be such that $f[\,[n_1]\,]\supset[n_0]$ and
$g[\,[n_1]\,]\supset[n_0]$.
Then for every $m,n\ge n_1$,
$$
{\textstyle
\big|\sum_{j=1}^m a_{f(j)}-
\sum_{j=1}^n a_{g(j)}\big|\le
\sum_{k\in A}|a_k|+\sum_{k\in B}|a_k|
\le\ep+\ep=2\ep\,,
}
$$
where $A,B\sus\N\setminus[n_0]$ are  finite sets. Thus every reordering of $\sum a_n$ has a~finite sum, 
because the sequence of partial sums is Cauchy ($f=g$), and all these 
sums are equal ($f\ne g$). 

If $\sum a_n$ is not an abscon series then the sum $\sum a_{z_n}=-\infty$ 
or the sum $\sum a_{k_n}=+\infty$,  and by 
Proposition~\ref{prop_naSoucNek} some reordering does not have finite sum.
\kduk

\noindent
{\em $\bullet$ The comparison criterion. }Absolute convergence of 
a~series is usually shown by comparing it with another abscon
series.

\begin{prop}[comparison criterion]\label{prop_srovKrit}
Suppose that 
$${\textstyle
\text{$\sum a_n$ and  
$\sum b_n$} 
}
$$
are series such that $a_n\ge0$ for every $n$, $|b_n|\le 
a_n$ for every $n\ge n_0$ and the sum $\sum a_n<+\infty$. \underline{Then} $\sum b_n$ is an abscon series.      
\end{prop}
\duk
Let $s_n\equiv\sum_{i=1}^n a_i$ and $t_n\equiv\sum_{i=1}^n |b_i|$. Then for every $n\ge n_0$,
$$
{\textstyle
0\le t_n\le
\sum_{i=1}^{n_0}|b_i|+s_n\,.
}
$$
Hence there is a~constant $c\ge0$ such that $0\le t_n\le c$ for every 
$n$ and $\sum b_n$ is an abscon series. 
\kduk
\vspace{-3mm}
\begin{exer}\label{ex_exOnCompCr}
Suppose that $(c_n)$ is a~bounded sequence. Prove that the series
$\sum\frac{c_n}{n(n+1)}$ is abscon. 
\end{exer}

\noindent
{\em $\bullet$ The Cauchy 
product of series. }We need this operation to prove the exponential 
identity in Theorem~\ref{thm_expFce}. In this context it is convenient to work  instead of $\N$ with the index 
set $\N_0$ and with series of the form 
$${\textstyle
\sum_{n=0}^{\infty}a_n=a_0+a_1+
\cdots\,. 
}
$$
We understand 
$a_0+a_1+\ds$ as the standard 
series $b_1+b_2+\ds$ in which 
$b_n\equiv a_{n-1}$. The \underline{Cauchy product\index{series!Cauchy 
product|emph}} of the series $\sum_{n=0}^{\infty}a_n$ and 
$\sum_{n=0}^{\infty}b_n$ is the series
$$
{\textstyle
\sum_{n=0}^{\infty}a_n\odot
\sum_{n=0}^{\infty}b_n\equiv
\sum_{n=0}^{\infty}\sum_{j=0}^n a_jb_{n-j}\,.
}\label{CauchyPr}
$$

\begin{thm}[Cauchy product of series]\label{thm_cauSouRad}
Any two\index{theorem!Cauchy product of series|emph} 
abscon series with sums $s$ and $t$ have abscon Cauchy product with the sum $st$.    
\end{thm}
\duk
Let $\sum_{n=0}^{\infty}a_n$ and $\sum_{n=0}^{\infty}b_n$ be abscon 
series with the respective sums $s$ 
and $t$. For $n\in\N_0$ let $s_n\equiv\sum_{j=0}^n a_j$, 
$t_n\equiv\sum_{j=0}^n b_j$, $z_n\equiv\sum_{j=0}^n\sum_{i=0}^j 
a_ib_{j-i}$, $s'\equiv\sum_{n=0}^{\infty}|a_n|$
(sum) and $t'\equiv\sum_{n=0}^{\infty}|b_n|$
(sum). Then
$$
|st-z_n|\le|st-s_nt_n|+|s_nt_n-z_n|\equiv A_n+B_n\,.    
$$
By Definition~\ref{def_infSer} and Theorem~\ref{thm_ari_lim}, 
$\lim A_n=0$. We show that also $\lim B_n=0$. We have 
$$
{\textstyle
B_n=\big|
\sum_{\substack{0\le i,j\le n\\
i+j>n}}a_ib_j\big|\le s'\sum_{\frac{n}{2}<j\le n}|b_j|+
t'\sum_{\frac{n}{2}<i\le n}|a_i|
\equiv s'C_n+t'D_n\,.
}
$$
Since the sequences $(\sum_{i=0}^n|a_i|)$ and 
$(\sum_{j=0}^n|b_j|)$ are Cauchy, we have $\lim C_n=\lim D_n=0$ and $\lim B_n=0$. Thus the sum of $\sum_{n=0}^{\infty}a_n\odot
\sum_{n=0}^{\infty}b_n$ is $st$. It remains to show that 
$z'\equiv\sum_{n=0}^{\infty}|z_n|<+\infty$. Since 
$${\textstyle
|z_n|\le\sum_{j=0}^n\sum_{i=0}^j 
|a_i|\cdot|b_{j-i}|\,, 
}
$$
our previous result implies that $z'\le s't'$ and we are done. 
\kduk
\vspace{-3mm}
\begin{exer}\label{ex_onCauPro}
Let $F(x)=\sum_{n=0}^{\infty}a_nx^n$
and $G(x)=\sum_{n=0}^{\infty}b_nx^n$
be formal power series with real coefficients and let $\sum_{n=0}^{\infty}c_n\equiv\sum_{n=0}^{\infty}a_n\odot\sum_{n=0}^{\infty}b_n$. Then $F(x)\cdot G(x)=\sum_{n=0}^{\infty}c_nx^n$.
\end{exer}

\noindent
{\em $\bullet$ A~geometric series\index{series!geometric|emph} }is any series of the form
$${\textstyle
\sum_{n=0}^{\infty}q^n=1+q+q^2+\ds+q^n+\cdots,\ q\in\R\,. 
}
$$
The number $q$ is the \underline{quotient\index{series!geometric!quotient of|emph}} of the
geometric series.

\begin{thm}[sums of geometric series]\label{thm_geomRada}
The\index{theorem!sums of geometric series|emph} 
geometric series has sum $\frac{1}{1-q}$ if $-1<q<1$ and $+\infty$ if $q\ge 1$. For $q\le-1$ the sum does not exist.
\end{thm}
\duk
For every $q\in\R\setminus\{1\}$ and  $n\in\N$ we have the identity
$${\textstyle
s_n\equiv 1+q+q^2+\ds+q^{n-1}=\frac{1-q^n}{1-q}=\frac{1}{1-q}+\frac{q^n}{q-1}\,.
}
$$
For $q<-1$ Theorem~\ref{thm_ari_lim} shows that $\lim s_{2n-
1}=+\infty$ and $\lim s_{2n}=-\infty$. Hence $\lim s_n$ does not exist. For $q=-1$ similarly $s_{2n-
1}=1$ and $s_{2n}=0$, the sum again does not exist. For $-1<q<1$ 
we have $\lim q^n=0$. Thus by Theorem~\ref{thm_ari_lim}, $\lim 
s_n=\frac{1}{1-q}$. For $q=1$ we have $s_n=n$ and the sum is 
$+\infty$. For $q>1$ we have $\lim q^n=+\infty$ and the sum is again 
$+\infty$ by Theorem~\ref{thm_ari_lim}.
\kduk

\noindent
As an application we get that, for example,
$$
{\textstyle
27.272727\ds=27(1+10^{-2}+10^{-4}+\ds)=27\cdot\frac{1}{1-10^{-2}}
=\frac{27\cdot100}{99}=\frac{300}{11}\;.
}
$$

\begin{exer}\label{ex_zobecGeom}  
Let $q\in(-1,1)$ and $m\in\Z$, with $q\ne0$ if $m<0$. Then the sum $\sum_{n\ge m}q^n=\frac{q^m}{1-q}$.
\end{exer}
\begin{exer}\label{ex_GeomRposlU}
Which geometric series are abskon?
\end{exer}

\noindent
{\em $\bullet$ Root and ratio tests. }These two classical convergence
criteria for series with nonnegative summands use bounds by geometric 
series $1+x+x^2+\cdots$ for $x\in[0,1)$.

\begin{prop}[root test]\label{prop_rootTest}
Let $a_n\ge0$ for every $n$. \underline{Then} the sum
$$
{\textstyle\sum a_n}\left\{
\begin{tabular}{lll}
$=+\infty$&$\ds$&$\limsup a_n^{1/n}>1$\,\text{ and}\\
$<+\infty$&$\ds$&$\limsup a_n^{1/n}<1$\,. 
\end{tabular}
\right.
$$     
\end{prop}
\duk
Let $\limsup a_n^{1/n}>1$. Then for some $c>1$ we have $a_n\ge 
c^n>1$ for infinitely many $n$. By Exercise~\ref{ex_druPomoc}, the sum $\sum 
a_n=+\infty$.

Let $\limsup a_n^{1/n}<1$. Then for some $c\in[0,1)$ and $n_0$ we have 
$0\le a_n\le c^n$ for every $n\ge n_0$. With $s_n\equiv\sum_{j=1}^n 
a_j$ we have
$${\textstyle
s_1\le s_2\le\ds\le s_n\le\ds\le a_1+\ds+a_{n_0-1}+
\frac{1}{1-c}
}
$$
\kduk
\vspace{-3mm}
\begin{exer}\label{ex_onRooTes}
Show that with nonstrict inequalities $\ds\ge1$ or $\ds\le1$ the root test is not valid.   
\end{exer}

\begin{prop}[ratio test]\label{prop_ratioTest}
Let always $a_n>0$. \underline{Then} the sum
$$
{\textstyle\sum a_n}\left\{
\begin{tabular}{lll}
$=+\infty$&$\ds$&${\textstyle\liminf\frac{a_{n+1}}{a_n}>1}$\,\text{ and}\\
$<+\infty$&$\ds$&${\textstyle\limsup \frac{a_{n+1}}{a_n}
<1}$\,. 
\end{tabular}
\right.
$$    
\end{prop}
\duk
Let 
$\liminf\frac{a_{n+1}}{a_n}>1$. Then for some $n_0$ and
$c>1$ we have $\frac{a_{n+1}}{a_n}\ge c$ for every $n\ge n_0$. 
It follows that $a_n\ge a_{n_0}>0$ 
for every $n\ge n_0$. Hence
$\sum a_n=+\infty$ by Exercise~\ref{ex_druPomoc}.

Let $\limsup\frac{a_{n+1}}{a_n}<1$. Then for some $n_0$ and $c\in(0,1)$ 
we have $0<\frac{a_{n+1}}{a_n}\le c$ for every $n\ge n_0$. Thus $0<a_n\le 
c^{n-n_0}a_{n_0}$ for every $n\ge n_0$. With $s_n\equiv\sum_{j=1}^n a_j$ we have
$${\textstyle
s_1\le s_2\le\ds\le s_n\le\ds\le a_1+\ds+a_{n_0-1}+
\frac{a_{n_0}}{1-c}\,.
}
$$
\kduk

\noindent Proposition~\ref{prop_delZeros2}
shows that the ratio test extends to series with nonnegative summands. 

\begin{exer}\label{ex_onRatTes}
Show that with the nonstrict inequalities $\ds\ge1$ or $\ds\le1$ the ratio test is not valid.   
\end{exer}

\begin{exer}\label{ex_limsupRatio}
Find a~convergent series $\sum 
a_n$ with positive summands and $\limsup\frac{a_{n+1}}{a_n}=+\infty$.
\end{exer}

\noindent
{\em $\bullet$ The zeta function (series) $\zeta(s)$. }We use real exponentiation
$a^b$ which we introduce soon in 
Section~\ref{sec_elemenFce}. 

\begin{defi}[series $\zeta(s)$]\label{def_radaZeta}
For $s\in\R$ the \underline{zeta series\index{series!zeta, $\zeta(s)$|emph}} 
is $\zeta(s)\equiv\sum\frac{1}{n^s}$.\label{zeta} 
\end{defi}

We determine the convergence of $\zeta(s)$ by the \underline{Cauchy condensation criterion}.

\begin{thm}[CCC]\label{thm_CKK}
Let\index{theorem!CCC|emph} $a_1\ge a_2\ge\ds\ge0$ be real numbers. \underline{Then} 
$${\textstyle 
\text{$\sum a_n$ converges
iff $R\equiv\sum2^n\cdot a_{2^n}$
converges}\,.
}
$$ 
\end{thm}
\duk
Suppose that $R$ has sum $+\infty$. Hence also the series 
$\frac{1}{2}R=\sum2^{n-1}\cdot a_{2^n}$ has sum $+\infty$. We have the inequalities
$${\textstyle
a_2\ge a_2,\,a_3+a_4\ge2a_4,\,\ds,\,  \sum_{j=2^{k-1}+1}^{2^k}a_j\ge 2^{k-1}a_{2^k},\,\cdots
}
$$
and summing them we get $\sum a_n=+\infty$. 

Suppose that $R$ converges. We have the inequalities
$${\textstyle
a_2+a_3\le2a_2,\, a_4+a_5+a_6+a_7\le 4a_4,\,\ds,\, \sum_{j=2^k}^{2^{k+1}-1}a_j\le 
2^k a_{2^k},\,\cdots
}
$$
and summing them we get that $\sum a_n$ converges. 
\kduk

\noindent
The proof of convergence of $\zeta(s)$ for $s>1$ is a~nice application of CCC.

\begin{thm}[convergence of $\zeta(s)$]\label{thm_zetaFce} 
For\index{theorem!convergence of $\zeta(s)$|emph} 
$s\le1$ the series $\zeta(s)$ has sum $+\infty$. For $s>1$ the zeta series converges.
\end{thm}
\duk
The former claim is Exercise~\ref{zetaSmen1}. Let $s>1$. The series $R$ in CCC for $\zeta(s)$ is 
$$
{\textstyle
\sum\frac{2^n}{(2^n)^s}=\sum\frac{1}{(2^{s-1})^n}\,.
}
$$
Since $0<\frac{1}{2^{s-1}}<1$, by Theorem~\ref{thm_geomRada} this geometric series converges. So by 
Theorem~\ref{thm_CKK} the series $\zeta(s)$ converges.
\kduk

\noindent
In {\em MA~1${}^+$} we show that $\zeta(2)=\frac{\pi^2}{6}$. 
This sum is due to the Swiss mathematician {\em Leonhard
Euler\index{Euler, Leonhard} (1707--1783)}. In {\em MA~1${}^+$} we
also show that the sum $\zeta(3)$ is an irrational number; this was proved by the French mathematician {\em Roger Ap\'ery\index{apery@Ap\'ery, Roger} (1916--1994)} in 1979.

\begin{exer}\label{zetaSmen1}
Prove that for $s\le1$ the series $\zeta(s)$ has sum $+\infty$.
\end{exer}

\begin{exer}\label{ex_aplCKT}
For which real $s$ does the series $\sum_{n=2}^{\infty}n^{-1}(\log n)^s$ converge? 
\end{exer}

\begin{exer}\label{ex_zeta2Easy}
Using bound $n^{-2}\le\frac{1}{n(n-1)}$ for $n\ge 2$, give a~simple proof of convergence of the series $\sum n^{-2}$. 
\end{exer}

\noindent
{\em $\bullet$ Using $\zeta(s)$. }We conclude this section on infinite series with three applications of the zeta function in mathematics and physics. The most 
important connection is to prime numbers. Let $(p_n)$ ($\sus\N$) be 
the ordering (Proposition~\ref{prop_ordSet}) of 
the set of primes, so that
$$
(p_n)=(2,\,3,\,5,\,7,\,11,\,13,\,17,\,19,\,23,\,\ds)\,.\label{pen}
$$

\begin{thm}[Euler's infinite product]\label{thm_EukerPro}
For
\index{theorem!Euler's product}
every real $s>1$ we have
$$
\prod_p\bigg(1-\frac{1}{p^s}\bigg)^{-1}=
\prod_{n=1}^{\infty}\bigg(1-\frac{1}{(p_n)^s}\bigg)^{-1}=
\lim_{n\to\infty}
\prod_{j=1}^n\bigg(1-\frac{1}{(p_j)^s}\bigg)^{-1}=\zeta(s)\,.
$$
\end{thm}

\begin{exer}\label{ex_eulPro}
Prove it.    
\end{exer}

\begin{exer}\label{ex_eulPro2}
Could not the value of Euler's infinite product be changed by reordering
the terms in it?   
\end{exer}

In combinatorics, more precisely in graph theory, we have the following 
interesting result of A.~M.~Frieze\index{Frieze, Alan M.} 
\cite{frie}. Let $n\in\N$ and $w\cc\binom{[n]}{2}\to[0,1]$. 
A~spanning tree\index{spanning tree|emph} $T$ on $[n]$ is a~tree (connected graph without cycles) $T=(V,E)$ such that
$V\sus[n]$ and $|E|=n-1$ ($E\sus\binom{[n]}{2}$).  

\begin{exer}\label{ex_spanTree}
Prove that then $V=[n]$ and $T$ is really ``spanning''.    
\end{exer}
We define the $w$-weight, or simply the weight, of the spanning tree $T$
by $w(T)=\sum_{e\in E}w(e)$. We denote by $M(w)$\label{Mw} the minimum
weight $w(T)$ for $T$ running through
all spanning trees on $[n]$. Recall 
that in probability theory we denote by $\mathbb{E}\,X$\label{expec} the expectation 
of the random variable $X$.

\begin{thm}[Frieze's]\label{thm_Frieze}
For\index{theorem!Frieze's} 
$n$ in $\N$ let $w=w_n$ be the random weight
$w\cc\binom{[n]}{2}\to[0,1]$, that is, $w_n$ is the product of $\binom{n}{2}$ 
independent copies of the uniform distribution on the interval 
$[0,1]$. \underline{Then}
$$
\lim_{n\to\infty}\mathbb{E}\,M(w_n)=\zeta(3)\,.
$$
\end{thm}

Finally, when browsing through the lecture notes {\em Low Temperature 
Physics} \cite{skrb_al} we spotted on page 16 the formula
$${\textstyle
E=\frac{\zeta(5/2)\cdot\Gamma(5/2)}{\zeta(3/2)\cdot\Gamma(3/2)}\cdot N\cdot k_{\mathrm{B}}\cdot T\cdot\big(\frac{T}{T_{\mathrm{B}}}\big)^{3/2}
}
$$
for the energy of the 
Bose\index{Bose, Satyendra Nath} gas\index{Bose gas} at temperature $T$.

\chapter[Limits of real functions]{Limits of real functions}\label{chap_pr5}

We begin with
Section~\ref{podkap_jedostrLimSpobod}
on one-sided limits of functions. Propositions~\ref{prop_ojednLimi} and 
\ref{prop_ojednLimi2} establish various relations between two-sided 
and one-sided limits. Section~\ref{sec_spojFcevBodu} 
introduces continuity function at a~point. In 
Proposition~\ref{prop_spojVbLimitou} 
we characterize it by limits, and in  Exercise~\ref{ex_HeiDefSpo} by 
Heine's definition. In Proposition~\ref{prop_spojVizol} we point out that every function is continuous at every isolated point of its domain.
Section~\ref{podkap_vlLimFci} 
contains Theorem~\ref{thm_limMonFce} on limits of monotone functions,
Theorem~\ref{thm_AritLimFce} on arithmetic of limits of functions, 
Theorem~\ref{thm_limFceUsp} on relations of limits of functions and 
the linear order $(\R^*,<)$, and a~generalization of the  
squeeze theorem for sequences in Theorem~\ref{thm_unexSque}. 
The most interesting result 
of Section~\ref{podkap_vlLimFci}
is the corrected Cauchy's theorem from \cite{brad_sand} in 
Theorem~\ref{thm_CauchyFixed}. 
In Section~\ref{sec_limSloFce} 
we present Theorem~\ref{thm_LimSlozFunkce}
on limits of composite functions. Our version is an equivalence.

Later I returned to this chapter and added the extending 
Section~\ref{sec_limInvFci} on limits of inverse 
functions. The relation of functional limits and 
functional inverses was  so far not considered in the 
literature. Thus the section contains original  results. We obtain 
conditions under which $\lim_{x\to A}f(x)=B$ implies $\lim_{y\to B}f^{\langle-1\rangle}(y)=A$. In 
Theorem~\ref{thm_LimInvFun1} we constrain the definition of limit. 
In Theorems~\ref{thm_LimInvFun2} and \ref{thm_LimInvFun2b} we 
constrain the considered function.
In Theorem~\ref{thm_LimInvFun3a} we give yet another sufficient 
condition for this inversion of limits.

The last two sections are devoted to asymptotic notation. In 
Section~\ref{sec_asympZnac} we explain the meaning of symbols $O$, 
$\ll$, $\gg$, $\Omega$, $\Theta$, $\asymp$, $o$, $\omega$ and $\sim$.
Definition~\ref{def_asympRel} introduces the notion of 
asymptoticity of relations on the set of real functions $\mathcal{R}$.
Relations $f=O(g)$ (on $N$), $f=o(g)$ ($x\to A$) and $f\sim g$ 
($x\to A$), as we define them, are asymptotic. In the extending Section~\ref{sec_asyExp} we explain
asymptotic expansions of
functions and give, without proofs, three examples of them: for $\log(n!)$, for the harmonic numbers $h_n$ and for the
probability of connectedness of
labeled graphs on $n$ vertices.

\section[Limits of functions]{Limits of functions}\label{podkap_limFunkci}

We extend limits of real sequences $(a_n)$, which are functions of the type 
$$
a\cc\N\to\R\,, 
$$
to functions $f\cc M\to\R$ defined on arbitrary sets $M\sus\R$.

\medskip\noindent 
{\em $\bullet$ Deleted neighborhoods and limit points. }Recall that for
$A,\ep\in\R$ with $\ep>0$
we have $\ep$-neighborhoods of points
$$
U(A,\,\ep)=(a-\ep,\,a+\ep)\,,
$$
and that for $A=\pm\infty$ we have $\ep$-neighborhoods of infinities
$${\textstyle
U(-\infty,\,\ep)=(-\infty,\,-\frac{1}{\ep})\,\text{ and }\,
U(+\infty,\,\ep)=(\frac{1}{\ep},\,+\infty)\,.
}
$$
We define the
\underline{deleted\index{neighborhood!deleted|emph} $\ep$-neighborhood} of $A\in\R^*$ by 
$$
P(A,\,\ep)\equiv U(A,\,\ep)\setminus\{A\}\,.\label{PAeps} 
$$
Let $M\sus\R$. An element $L\in\R^*$ is a~\underline{limit point\index{limit point!of 
a set|emph}} of $M$ if for every $\ep$ we have $P(L,\ep)\cap M\ne\emptyset$.
The set of limit points of $M$ is denoted by $L(M)$\index{limit point!of a 
set!lm@$L(M)$|emph} ($\sus\R^*$).\label{LM}

\begin{exer}\label{ex_limBodyLimitou}
Prove the following proposition.
\end{exer}

\begin{prop}[on limit points]\label{prop_OlimBodech}
Let $M\sus\R$ and $A\in\R^*$. The next four claims are mutually equivalent.
\begin{enumerate}
    \item $A\in L(M)$.
    \item There is a~sequence $(a_n)\sus M\setminus\{A\}$ such that $\lim a_n=A$.
    \item There is an injective sequence $(a_n)\sus M$ such that $\lim a_n=A$.
    \item For every $n\in\N$ we have $P(A,\frac{1}{n})\cap M\ne\emptyset$.
\end{enumerate}
\end{prop}

\begin{exer}\label{ex_AjesteJedenPr}
A~set $M\sus\R$ is finite $\iff$ $L(M)=\emptyset$.
\end{exer}

\begin{exer}\label{ex_vyjmutiBodu}
If $M\sus\R$ and $b\in L(M)$ then also
$b\in L(M\setminus\{b\})$.
\end{exer}

\noindent
{\em $\bullet$ Limits of real functions. }We introduce notation we will often use.

\begin{defi}[notation for real functions]\label{def_znacFci}
For $M\sus\R$ we define 
$$
\mathcal{F}(M)\index{function!fm@$\mathcal{F}(M)$|emph}\equiv\{f=\langle M,\,\R,\,G_f\rangle\cc\;f\cc M\to\R\,\text{ and }\,M\sus\R\}\label{FM}
$$
(recall Definition~\ref{def_funkce})
and set
$${\textstyle
\mathcal{R}
\underline{\index{function!r@$\mathcal{R}$|emph}}\equiv\bigcup_{M\sus\R}\mathcal{F}(M)\,.\label{funcR} 
}
$$
For $f\in\mathcal{F}(M)$ we define 
$$
Z(f)\equiv\{b\in M\cc\;f(b)=0\}\ \ 
(\sus M)\,.
\index{function!zf@$Z(f)$|emph}\label{zeroFun}
$$
Recall that the domain of any function $f\cc X\to Y$ is $M(f)=X$.
\end{defi}
Thus $\mathcal{F}(M)$, for $M\sus\R$, is the set of functions with the definition domain 
$M$ and range $\R$, and $\mathcal{R}$ is the set of all such 
functions for all sets of real numbers~$M$. By $Z(f)$ we denote the 
set of zeros of $f$. The \underline{real empty 
function\index{real empty function|emph}\index{function!real 
empty|emph}} is
$$
\emptyset_f\equiv\langle\emptyset,\,\R,\,\emptyset\rangle\ \ (\in\mathcal{R})\,.\label{empReaFun}
$$
It is an empty function.

\begin{defi}[limits of functions]\label{def_limFce}
Let $f$ be in $\mathcal{R}$, $A$ in $L(M(f))$ and $L$ in $\R^*$. If for every $\ep$ there is a~$\de$ such that
\begin{equation}
 f[P(A,\,\de)]\sus U(L,\,\ep)\tag{$*$}\,,   
\end{equation}
we write $\lim_{x\to A}f(x)=L$\label{limFun} and say that $f$ has at $A$ the \underline{limit\index{limit of a function|emph}} $L$.
\end{defi}
(In Definition~\ref{def_strongLim} we consider a~strengthening of this definition.)
Recall that by our definition of image, 
$$
f[P(A,\,\de)]=f[P(A,\,\de)\cap 
M(f)]=\{f(x)\cc\;x\in P(A,\,\de)\cap M(f)\}\,. 
$$
The limit does not depend on the value $f(A)$ and the function $f$ need not be defined at $A$. If $A=\pm\infty$ then it in
fact cannot 
be defined at $A$. For a~sequence $(a_n)\sus\R$, which is a~function 
$a\cc\N\to\R$, we have
$\lim_{x\to+\infty}a(x)=\lim a_n$.

\begin{exer}\label{ex_kteryDalsiLimB}
Find all limit points of the set $\N$ ($\sus\R$).  
\end{exer}
If $A=a\in\R$ and $L=b\in\R$, then 
$\lim_{x\to a}f(x)=b$ means that
$$
\forall\,\ep\,\exists\,\de\,\big(x\in M(f)\wedge 0<|x-a|\le\de\Rightarrow |f(x)-b|\le\ep \big)
$$
(recall that we like $\le$ more than $<$). We stress that 
$$
\text{if $f\in\mathcal{F}(M)$ and $A\not\in L(M)$ then 
$\lim_{x\to A}f(x)$ \underline{is not defined}}\,. 
$$
For then for some $\de$ 
we have $P(A,\de)\cap M=\emptyset$ and
$f[P(A,\de)]=\emptyset$, which means that the above inclusion
$(*)$ holds for every $L$ and
every $\ep$. If $\lim_{x\to A}f(x)$ exists then always $A\in L(M(f))$.

\begin{prop}[locality of limits]\label{prop_locLimFce}
If\underline{\index{limit of a function!locality of|emph}} 
$f,g\in\mathcal{R}$, $A\in\R^*$ and for some $\theta$ we have
$f\,|\,P(A,\theta)=g\,|\,P(A,
\theta)$, \underline{then} 
$$
\lim_{x\to A}f(x)=\lim_{x\to 
A}g(x)\,, 
$$
if one side is defined.   
\end{prop}
\duk
We can take $\de$ in Definition~\ref{def_limFce} such
that $\de\le\theta$. Then
$P(A,\de)\sus P(A,\theta)$ and $f[P(A,\de)]=g[P(A,\de)]$. 
\kduk

\noindent Later we show that continuity and derivatives are local 
too.

\begin{prop}[unique limits]\label{prop_jednLimi}
If $\lim_{x\to K}f(x)=L$ and $\lim_{x\to K}f(x)=L'$ \underline{then} $L=L'$. Thus limits\underline{\index{limit of a function!uniqueness of|emph}} 
of functions are unique.
\end{prop}
\duk
For every $\ep$ there is a~$\de$ such that the nonempty (!) set 
$f[P(K,\de)]$ is contained both in  $U(L,\ep)$ and $U(L',\ep)$. Thus 
$\forall\ep\big(U(L,\ep)\cap U(L',\ep)\ne\emptyset\big)$. 
By Exercise~\ref{ex_ulohaNaokoli1}, $L=L'$.
\kduk

\noindent
We show that limits of restrictions are equal to limits of original 
functions. 

\begin{prop}[limits of restrictions]\label{prop_restrAlimFce}
Let $f\in\mathcal{F}(M)$, $X$ be any set, $A\in L(X\cap M)$ and let $\lim_{x\to A}f(x)=L$. \underline{Then} $\lim_{x\to A}(f\,|\,X)(x)=L$. 
\end{prop}
\duk
Let an $\ep$ be given. Then there is a~$\de$ such that $f[P(A,\de)]\sus 
U(L,\ep)$. From $P(A,\de)\cap(X\cap M)\sus 
P(A,\de)\cap M$ we get
$$
(f\,|\,X)[P(A,\,\de)]\sus
f[P(A,\,\de)]\sus U(L,\,\ep)\,.
$$
Hence $\lim_{x\to A}(f\,|\,X)(x)=L$.
\kduk

\noindent
This is the first of several results on interactions between limits of
functions and operations on $\mathcal{R}$. Operations we
consider later are composition and the arithmetic operations of addition, multiplication,
and division. Functional inverses are relegated to {\em MA~1${}^+$}.

\begin{exer}\label{ex_uvedtePriklad}
Find a~function 
$f\in\mathcal{F}(M)$ and a~set $X$  such that $\lim_{x\to 
A}f(x)$ does not exist but $\lim_{x\to A}(f\,|\,X)(x)$ exists.
\end{exer}

\noindent
{\em $\bullet$ Heine's definition of the limit of a~function. }It is 
clear that limits of sequences are particular cases of limits of 
functions. The German mathematician E. Heine realized that one can go in the opposite direction and reduce the limits of functions to the limits of sequences. 

\begin{thm}[Heine's limit]\label{thm_HeinehoDef}
Let\index{theorem!Heine's definition of limits of functions|emph}\underline{\index{limit of a function!Heine's definition of|emph}} 
$f\in\mathcal{F}(M)$ and $K\in L(M)$. \underline{Then} 
$$
\lim_{x\to K}f(x)=L\iff \forall\,(a_n)\sus M\setminus
\{K\}\,\big(\lim a_n=K\Rightarrow
\lim f(a_n)=L\big)\,.
$$
In words, a~function $f$ defined on $M$ has limit $L$ at $K$ iff for 
every sequence $(a_n)$ in $M$ with $a_n\ne K$ but with $\lim a_n=K$ the 
sequence $(f(a_n))$ has limit $L$.
\end{thm}
\duk
Implication $\Rightarrow$. Suppose that 
$$
\lim_{x\to K}f(x)=L\,, 
$$
that $(a_n)\sus M\setminus\{K\}$ is a~sequence with $\lim a_n=K$ and that an $\ep$ is given. 
Then there is a~$\de$ such that for every $x\in P(K,\de)\cap M$ we have 
$f(x)\in U(L,\ep)$. For this $\de$ there is an $n_0$ such 
that for every $n\ge n_0$ we have $a_n\in P(K,\de)\cap M$. 
Hence if $n\ge n_0$ then $f(a_n)\in U(L,\ep)$ and $f(a_n)\to L$.

Reverse implication $\neg\Rightarrow\neg$. Suppose that
$$
\text{it is not true that 
$\lim_{x\to K}f(x)=L$}\,. 
$$
Then there exists a~number $\ep$ such that for every $\de$ 
there exists a~point $b=b(\de)\in P(K,\de)\cap M$ with 
$$
f(b)\not\in U(L,\ep)\,. 
$$
We set $\de=\frac{1}{n}$ and select points 
$${\textstyle
b_n\equiv b(\frac{1}{n})\in P(K,\frac{1}{n})\cap M
}
$$ 
such that $f(b_n)\not\in U(L,\ep)$ for every $n$. Then 
$(b_n)$ is a~sequence in $M\setminus\{K\}$ with $\lim 
b_n=K$, but the sequence $(f(b_n))$ does not have limit $L$. 
The right side of the equivalence does not hold. 
\kduk

\noindent
In the proof of the reverse implication we used the axiom of choice\index{axiom!of choice, AC}. 

\begin{exer}\label{ex_AC}
How exactly did we use it? 
\end{exer}

\noindent
{\em $\bullet$ A~functional limit. }We get with the help of the identities
$x-y=\frac{x^2-y^2}{x+y}$ a~$\frac{x}{y}=\frac{1}{y/x}$ 
that
\begin{eqnarray*}
{\textstyle
\lim_{x\to+\infty}\big(\sqrt{x+\sqrt{x}}-\sqrt{x}\big)}&=&{\textstyle
\lim_{x\to+\infty}\frac{\sqrt{x}}{\sqrt{x+\sqrt{x}}\,+\,\sqrt{x}}
}\\
&=&{\textstyle
\lim_{x\to+\infty}\frac{1}{\sqrt{1+1/\sqrt{x}}\,+\,1}
}\\
&=&{\textstyle
\frac{1}{\sqrt{1+1/(+\infty)}\,+\,1}=\frac{1}{1+1}=\frac{1}{2}\,.
}
\end{eqnarray*}

\begin{exer}\label{ex_triLimity}
Compute the following limits of functions.
\begin{enumerate}
    \item $\lim_{x\to-\infty}\frac{x}{\sqrt{1+x^2}-1}$.
    \item $\lim_{x\to+\infty}\frac{1}{\sqrt{1+x}-\sqrt{x}}$.
    \item $\lim_{x\to0}\frac{1}{x}$.
    \item $\lim_{x\to-\infty}\frac{1}{x}$.
\end{enumerate}
\end{exer}

\noindent
{\em $\bullet$ A~theorem of A.-L. Cauchy. }The book {\em Cours 
d'analyse}\index{Cours d'analyse}, see \cite{brad_sand} for English 
translation, ``is a~seminal textbook in infinitesimal calculus published 
by Augustin-Louis Cauchy\index{Cauchy, Augustin-Louis} 
in 1821'' \cite{coursWiki}. In the next proposition we present counterexamples to Theorem~I 
in Section~2.3 of it. We illustrate by this limits of functions. 
Cauchy's theorem, quoted from \cite{brad_sand}, says the following.
\begin{quote}
If the difference $f(x+1)-f(x)$
converges towards a~certain limit $k$, for increasing values of $x$, then the fraction
$\frac{f(x)}{x}$ converges at the same time towards the same limit.    
\end{quote}
However, we found the following counterexamples.

\begin{prop}[counterexamples]\label{prop_couExa}
Let $k\in\R$ and 
$$
(a_n)\sus(0,\,+\infty)
$$ 
be any sequence such that $a_1<a_2<\ds$,
$\lim a_n=+\infty$ and $a_n-a_m\not\in\N$ for every pair of indices $m<n$ (see Exercise~\ref{ex_giveExam}). 
\underline{Then} there exists 
a~function $f\cc(0,+\infty)\to\R$ such that
$$
f(x+1)-f(x)=k
$$ 
for every $x>0$, but $f(a_{2n})=0$ and $f(a_{2n-1})=a_{2n-1}$ for
every $n$. We have in particular  
$$
\lim_{x\to+\infty}(f(x+1)-f(x))=k
$$ 
but 
$${\textstyle
\frac{f(a_{2n})}{a_{2n}}=0\,\text{ and }\, 
\frac{f(a_{2n-1})}{a_{2n-1}}=1
}
$$ 
for every $n$, hence
$\lim_{x\to+\infty}\frac{f(x)}{x}$ does not exist.
\end{prop}
\duk
Let $k$ and $(a_n)$ be as stated, and $y\in(0,1]$ be arbitrary. If $y+m\ne a_n$ for 
every $m\in\N_0$ and every $n\in\N$, we define the value $f(y)\in\R$ arbitrarily and then for 
every $l\in\N$ we set 
$$
f(y+l)\equiv f(y)+lk\,.
$$
If $y+m=a_{2n}$ for a~unique $m\in\N_0$ and $n\in\N$, we define $f(y)\equiv-mk$ and then for
every $l\in\N$ we set 
$$
f(y+l)\equiv f(y)+lk\,.
$$
Finally, if $y+m=a_{2n-1}$ for
a~unique $m\in\N_0$ and $n\in\N$, we define $f(y)\equiv a_{2n-1}-mk$ and 
then for every $l\in\N$ we set 
$$
f(y+l)\equiv f(y)+lk\,.
$$
By the assumption on $(a_n)$ this definition of $f$ is correct. It is 
clear that $M(f)=(0,+\infty)$ and that the function $f$ has the stated 
properties.
\kduk
\vspace{-3mm}
\begin{exer}\label{ex_giveExam}
Give examples of real sequences $(a_n)$
with the properties stated in the proposition.
\end{exer}
The next exercise requires access to \cite{brad_sand}.
\begin{exer}\label{ex_CauchyError}
Find an error in Cauchy's proof of Theorem I in Section~2.3 of the {\em Cours d'analyse} {\em (\cite{brad_sand})} by which if 
$f\cc(0,+\infty)\to\R$ has 
$$
\lim_{x\to+\infty}(f(x+1)-f(x))=k\ \ 
(\in\R)\,, 
$$
then $\lim_{x\to+\infty}\frac{f(x)}{x}=k$.     
\end{exer}
We fix Cauchy's theorem in Theorem~\ref{thm_CauchyFixed} by adding monotonicity assumption.

\section[One-sided limits]{One-sided limits}\label{podkap_jedostrLimSpobod}

The\underline{\index{limit of a function!one-sided|(}} 
complement of any point $a$ to the real axis, 
$$
\R\setminus\{a\}=(-\infty,\,a)
\cup(a,\,+\infty)\,,
$$ consists of two separated intervals. In the 
plane $\R^2$ we can go around a~point but in $\R$ this is impossible. So we consider 
left-sided and right-sided limits of functions. First we define one-sided 
neighborhoods.

\medskip\noindent
{\em $\bullet$ One-sided neighborhoods and one-sided limit points. }A~\underline{left\index{neighborhood!left|emph}}, resp. a~\underline{right\index{neighborhood!right|emph}}, \underline{$\ep$-neighborhood} of a~point $b\in\R$ is
$$
\text{$U^-(b,\,\ep)\equiv(b-\ep,\,b]$, resp. $U^+(b,\,\ep)\equiv[b,\,b+\ep)$}\,.\label{neighbPm} 
$$
A~\underline{left\index{neighborhood!left deleted|emph}}, resp. a~\underline{right\index{neighborhood!right deleted|emph}, deleted $\ep$-neighborhood} of a~point $b\in\R$ is 
$$
\text{$P^-(b,\,\ep)\equiv(b-\ep,\,b)$, resp. 
$P^+(b,\,\ep)\equiv(b,\,b+\ep)$}\,.
\label{delNeiPm} 
$$

A~point $b\in\R$ is a~\underline{left\index{limit 
point!left|emph}}, resp. a~\underline{right\index{limit 
point!right|emph}, limit point} of $M\sus\R$ if for every $\ep$ we have
$$
\text{$P^-(b,\,\ep)\cap M\ne\emptyset$, resp.
$P^+(b,\,\ep)\cap M\ne\emptyset$}\,. 
$$
The set of these points is denoted by \underline{$L^-(M)$\index{limit 
point!lmminus@$L^-(M)$|emph}}, resp. \underline{$L^+(M)$\index{limit 
point!lmplus@$L^+(M)$|emph}} 
($\sus\R$).\label{LMpm} 
A~point $b\in\R$ is a~\underline{two-sided limit point\index{limit point!two-sided|emph}} of $M\sus\R$ if for every $\ep$,
$$
\text{$P^-(b,\,\ep)\cap M\ne\emptyset\,\wedge\,P^+(b,\,\ep)\cap M\ne\emptyset$}\,. 
$$
The set of these points is denoted by \underline{$L^{\mathrm{TS}}(M)$\index{limit 
point!lmminplu@$L^{\mathrm{TS}}(M)$|emph}} ($\sus\R$).\label{LMts} Two-sided limit points 
play a~key role in the criterion of local extremes. We do not 
define one-sided neighborhoods for infinities, nor $\pm\infty$ can be 
a~one-sided limit point of a~set.

\begin{exer}\label{ex_uloZpredn}
$b\in L^-(M)$, resp. $b\in L^+(M)$ $\iff$ $\exists$ $(a_n)\sus
(-\infty,b)\cap M$, resp.  $\exists$ $(a_n)\sus(b,+\infty)\cap M$, such that $\lim a_n=b$.    
\end{exer}

\begin{exer}\label{ex_naJednostrLB}
Let $M\sus\R$ and $b\in\R$. Prove the following.
\begin{enumerate}
    \item If $b\in L^-(M)$, then $b\in L(M)$.
    \item If $b\in L^+(M)$, then $b\in L(M)$.
    \item If $b\in L(M)$, then $b\in L^-(M)$ or $b\in L^+(M)$.
    \item In general, it is not true that if $b\in L(M)$ then $b\in L^-(M)$ and $b\in L^+(M)$.
\end{enumerate}
\end{exer}
Let $M\sus\R$. By Exercise~\ref{ex_AjesteJedenPr} no finite set $M$ has a~limit point, the 
less one-sided limit point, and  
every infinite $M$ has a~limit point. This is not true for one-sided 
limit points.

\begin{exer}\label{ex_NekMnJednostr1}
Find infinite subsets of $\R$ without one-sided limit points.
\end{exer}

\begin{exer}\label{ex_NekMnJednostr}
Every infinite and bounded real set has a~one-sided limit point.
\end{exer}

\noindent
{\em $\bullet$ One-sided limits of functions. }We consider finer variants of limits of
functions, the one-sided limits.

\begin{defi}[one-sided limits]
Let\underline{\index{limit of a function!one-sided|emph}} 
$f$ be in $\mathcal{F}(M)$, $b$ be in 
$L^-(M)$ and $L$ be in $\R^*$. If for every 
$\ep$ there is a~$\de$ such that 
$$
f[P^-(b,\,\de)]\sus U(L,\,\ep)\,,
$$
we write $\lim_{x\to b^-}f(x)=L$ and say that $f$ has 
at $b$ the \underline{left-sided limit} $L$. Replacing the three signs ${}^-$ with three signs ${}^+$
we get the \underline{right-sided limit} at $b$ that is denoted by $\lim_{x\to b^+}f(x)=L$.\label{limFunOnes}
\end{defi}
Like the ordinary limit, the one-sided limit of $f$ at $b$ is not defined if $b$ is not the respective one-sided limit point of $M(f)$. The existence of
$\lim_{x\to b^{\pm}}f(x)$ again means that $b$ is the respective one-sided limit point of $M(f)$.

\begin{exer}\label{ex_zpetneVyleps}
Prove the following proposition.
\end{exer}

\begin{prop}[one-sided limits]\label{prop_ojednLimi}
The following hold.
\begin{enumerate}
\item If $\lim_{x\to a}f(x)=L$, then $\lim_{x\to a^-}f(x)=L$ 
or $\lim_{x\to a^+}f(x)=L$.
\item If $\lim_{x\to a^-}f(x)=\lim_{x\to a^+}f(x)=L$, then $\lim_{x\to a}f(x)=L$.
\item If $\lim_{x\to a^-}f(x)=K$, $\lim_{x\to a^+}f(x)=L$ and $K\ne L$, then the limit $\lim_{x\to a}f(x)$ does not exist.
\end{enumerate}
\end{prop}
For example, $\lim_{x\to 0}\sgn\,x$ does not exist because $\lim_{x\to 
0^-}\sgn\,x=-1$ and $\lim_{x\to 0^+}\sgn\,x=1$.

\begin{exer}\label{ex_jednozn_Jednostr}
Prove the following proposition.    
\end{exer}

\begin{prop}[uniqueness of one-sided limits]\label{prop_jednJednostr}
If $\lim_{x\to b^{\pm}}f(x)=K$ and $\lim_{x\to b^{\pm}}f(x)=L$ \underline{then} $K=L$, with equal signs.
\end{prop}

\begin{exer}\label{ex_Heine_Jednostr}
Prove the following proposition.    
\end{exer}
For $b\in\R$ let $I^-(b)\equiv
(-\infty,b)$ and $I^+(b)\equiv(b,+\infty)$.\label{ibPlMi}

\begin{prop}[Heine's limits]\label{prop_HeineJednostr} Suppose that 
$f$ is in $\mathcal{F}(M)$ and $b$ is in $L^{\pm}(M)$. \underline{Then} 
$\lim_{x\to b^{\pm}}f(x)=L$ $\iff$ for every sequence
$(a_n)\sus M\cap I^{\pm}(b)$ with $\lim a_n=b$ one has
$\lim f(a_n)=L$, with equal signs.
\end{prop}

Sometimes ordinary limits are replaced unnecessarily with one-sided ones. For example, we can see notation
$\lim_{x\to 0^+}\log x$. In our notation it suffices and is correct to write just
$\lim_{x\to 0}\log x$. Indeed,  
$$
\lim_{x\to 0}\log x=\lim_{x\to 0^+}\log x=-\infty\,,
$$
but $\lim_{x\to 0^-}\log x$ is not defined because 
$$
0\not\in L^-(M(\log 
x))=L^-((0,\,+\infty))=(0,\,+\infty)\,. 
$$
We conclude with one more relation between 
ordinary and one-sided limits. We use it in the proof of Corollary~\ref{cor_ekvivZde}. 

\begin{exer}\label{ex_vztObyJed}
Prove the following proposition.    
\end{exer}

\begin{prop}[using restriction]\label{prop_ojednLimi2}
Let $f\in\mathcal{F}(M)$ and $b\in L^{\pm}(M)$. \underline{Then} 
$$
\lim_{x\to b^{\pm}}f(x)=L\iff\lim_{x\to b}(f\,|\,
I^{\pm}(b))(x)=L\,, 
$$
with equal signs.
\underline{\index{limit of a function!one-sided|)}}
\end{prop}

\section[Continuity at a~point]{Continuity at a~point}\label{sec_spojFcevBodu}

\medskip\noindent
We arrive at an important definition.

\begin{defi}[pointwise continuity]\label{def_spojVbode}
Let $f\in\mathcal{F}(M)$ and $b\in M$. Then $f$ is \underline{continuous at\index{continuity of functions!at a~point|emph}}  $b$ if for 
every $\ep$ there is a~$\de$ such that 
$$
f[U(b,\,\de)]\sus U(f(b),\,\ep)\,.
$$
Else we say that $f$ is \underline{discontinuous at\index{discontinuity of $f$ at a~point|emph}} $b$.
\end{defi}

\begin{exer}\label{ex_kyJeNesp}
Let $b\in M(f)$. Explain in detail when $f$ is
discontinuous at $b$. 
\end{exer}
For example, $\sgn\,x$ is discontinuous at $0$, but it is continuous at any other point. If 
$b\not\in M(f)$ then $f$ is {\em neither continuous nor 
discontinuous} at $b$. Comparing the above definition with
$$
\lim_{x\to b}f(x)=L\iff
\forall\ep\exists\de\,\big(f[P(b,\,\de)]\sus U(L,\,\ep)\big)\,,
$$ 
we see that the last $L$ is replaced with $f(b)$, and the
deleted neighborhood $P(b,\de)$ with the full neighborhood $U(b,\de)$.

\begin{prop}[locality of continuity]\label{prop_locCont}
If\underline{\index{continuity of functions!at a~point!locality of|emph}} 
$f,g\in\mathcal{R}$, $b\in M(f)\cap M(g)$ and there is a~$\theta$ such 
that $f=g$ on $U(b,\theta)$ \underline{then} $f$ is continuous
at $b$ $\iff$ $g$ is continuous at $b$. 
\end{prop}
\duk
This is immediate from Definition~\ref{def_spojVbode} 
because we can take the $\de$ in it such that $\de\le\theta$. Then
$U(b,\de)\sus U(b,\theta)$ and $f[U(b,\de)]=g[U(b,\de)]$.
\kduk
\vspace{-3mm}
\begin{exer}\label{ex_klasDefSpovBod}
A~function $f$ is continuous at $b$ ($\in M(f)$) iff for every $\ep$ there exists a~$\de$ such that $x\in M(f)\wedge|x-b|\le\de\Rightarrow|f(x)-f(b)|\le\ep$.
\end{exer}
The continuity of $f$ at $b$ is {\em not}
equivalent to $\lim_{x\to b}f(x)=f(b)$. This only holds for
limit points of $M$. 

\begin{prop}[pointwise continuity]\label{prop_spojVbLimitou}
For any function $f\in\mathcal{F}(M)$ and any point $b\in M\cap L(M)$ the following are mutually equivalent.
\begin{enumerate}
\item The function $f$ is continuous at $b$.
\item The limit $\lim_{x\to b}f(x)=f(b)$.
\item For every sequence $(a_n)\sus M$ with $\lim a_n=b$ we have $\lim f(a_n)=f(b)$.
\end{enumerate}
\end{prop}
\duk
Let $f$ and $b$ be as stated. 
The implication $1\Rightarrow 2$. Let $f$ be continuous at $b$ by
Definition~\ref{def_spojVbode} and let an $\ep$ be given. Thus there is a~$\de$ 
such that $f[U(b,\de)]$ $\sus$ $U(f(b),\ep)$. We have $b\in L(M)$ 
and $f[P(b,\de)]$ $\sus$ $U(f(b),\ep)$. Hence $\lim_{x\to b}f(x)=f(b)$.

The implication $2\Rightarrow 3$. Suppose that $\lim_{x\to 
b}f(x)=f(b)$, that $(a_n)\sus M$ has $\lim a_n=b$ and that an $\ep$ is
given. Thus there is a~$\de$ such that 
\begin{equation}
    f[P(b,\,\de)]\sus U(f(b),\,\ep)\;.\tag{$*$} 
\end{equation}
We take an $n_0$ such that $n\ge n_0$ $\Rightarrow$ $a_n\in U(b,\de)$. Then 
also $n\ge n_0$ $\Rightarrow$ $f(a_n)\in U(f(b),\ep)$: for $a_n\ne 
b$ we use the inclusion $(*)$, and for $a_n=b$ it is automatic that 
$f(a_n)=f(b)\in U(f(b),\ep)$. Hence $\lim f(a_n)=f(b)$.  

The implication $3\Rightarrow 1$. We prove the reversal $\neg1\Rightarrow\neg3$. Suppose that $f$ is not 
continuous at $b$. Then there is an $\ep$ such that $\forall\,\de\,\exists\,a=a(\de)\in 
U(b,\de)\cap M$ with $f(a)\not\in U(f(b),\ep)$. For every $n\in\N$ we 
{\em chose} such  $a_n\equiv a(\frac{1}{n})$ and get the sequence $(a_n)\sus M$ such that $\lim a_n=b$, but for every $n$ it 
holds that $f(a_n)\not\in U(f(b),\ep)$. Hence $(f(a_n))$ does 
not converge to $f(b)$ and part~3 does not hold. 
\kduk

\noindent
The last implication again used the axiom of choice.\index{axiom!of 
choice, AC} Part~3 describes {\em Heine's definition of pointwise
continuity\index{continuity of 
functions!at a~point!Heine's definition of|emph}}. 

\begin{exer}\label{ex_HeiDefSpo}
In Proposition~\ref{prop_spojVbLimitou} one can omit in the equivalence {\em 1} $\iff$ 
{\em 3} the assumption that
$b\in L(M)$. Thus, $f$ is continuous at a~point $b\in M(f)$ $\iff$ for
every sequence $(a_n)\sus M$ with $\lim a_n=b$ one has
$\lim f(a_n)=f(b)$.
\end{exer}
The right side of this equivalence is sometimes taken as the definition
of point-wise continuity.

\medskip\noindent
{\em $\bullet$ Isolated points. }Let $M\sus\R$. The set $M\setminus L(M)$ consists of the so-called \underline{isolated points\index{isolated point|emph}} of
$M$.

\begin{exer}\label{ex_oIzBodech}
A~point $b\in M$ ($\sus\R$) is an isolated point 
of $M$ iff for some $\ep$ we have $U(b,\ep)\cap M=\{b\}$.    
\end{exer}

\begin{exer}\label{ex_limAizBody}
 Let $b\in M\sus\R$. Then $b$ is either a~limit point of $M$ or an 
 isolated point of $M$.
\end{exer}

\begin{prop}[continuity at isolated points]\label{prop_spojVizol}
Every function $f\in\mathcal{R}$ is continuous at every isolated point 
of $M(f)$.
\end{prop}
\duk
Let $f\in\mathcal{F}(M)$ and let $b\in M$ be an isolated point of $M$. By 
Exercise~\ref{ex_oIzBodech} there is a~$\de$ such that $U(b,\de)\cap M=
\{b\}$. For this $\de$ the inclusion
$f[U(b,\de)]=f[\{b\}]=\{f(b)\}\sus U(f(b),\ep)$ 
holds for every $\ep$. Hence $f$ is continuous at~$b$ by Definition~\ref{def_spojVbode}.
\kduk

\noindent
We see that every sequence $(a_n)\sus\R$, understood as 
a~function $a\in\mathcal{F}(\N)$, is  continuous at every point $n$ of the definition domain $\N$.

\begin{exer}\label{ex_ekvivNespoVbode}
Let $f\in\mathcal{F}(M)$ and $b\in M$. Then $f$ is not continuous at 
$b$ if and only if there exists a~sequence $(a_n)\sus 
M$ such that $\lim a_n=b$, $\lim f(a_n)=A$ and $A\ne f(b)$.
\end{exer}

\noindent
{\em $\bullet$ One-sided continuity at a~point. }A~function $f$ is \underline{left-continuous\index{continuity of 
functions!at a~point!left-|emph}} at $b\in M(f)$ if for every $\ep$ there is a~$\de$ such that
$$
f[U^-(b,\,\de)]\sus U(f(b),\,\ep)\,.
$$
By replacing the sign ${}^-$ with the sign ${}^+$ we get the 
\underline{right-continuity\index{continuity of 
functions!at a~point!right-|emph}}.

\begin{exer}\label{ex_ekvivSpoji}
A~function is continuous at a~point iff it is both left- and 
right-continuous at the point.
\end{exer}

\noindent
{\em $\bullet$ \underline{Riemann's function\index{Riemann, Bernhard! 
function|emph}}. }It is the function $r\in\mathcal{F}(\R)$ with values
$r(x)=0$ for $x\in\R\setminus\Q$ and $r(\frac{m}{n})=\frac{1}{n}$
for fractions $\frac{m}{n}$ in lowest terms.\label{riemannFun}

\begin{prop}[Riemann's function]\label{prop_RiemannFunkce}
Riemann's function is continuous exactly at irrational numbers.
\end{prop}
\duk
Let $x=\frac{m}{n}$ be a~fraction in lowest terms and $\ep
\le\frac{1}{n}$. For every $\de$ there is an  irrational $\al\in 
U(x,\de)$. But $r(\al)=0\not\in 
U(r(x),\ep)=U(\frac{1}{n},\ep)$, so that $r$ is discontinuous at $x$. 
Let $x\in\R\setminus\Q$ and $\ep\in(0,1)$ be given. We set 
$${\textstyle
M\equiv\{|x-\frac{m}{n}|\cc\;\frac{m}{n}\in\Q\cap U(x,\,1)\wedge \frac{1}{n}\ge\ep\}\,\text{ and }\,\de\equiv\min(M)\;.
}
$$
By Exercise~\ref{ex_ustne}, $\de$ exists and is positive. For this $\de$ we have $y\in U(x,\de)$ $\Rightarrow$ 
$r(y)\in U(r(x),\ep)=U(0,\ep)$\,---\,for every $y\in U(x,\de)$ we have $r(y)=0$ or $r(y)=\frac{1}{n}<\ep$. Hence $r$ is continuous at $x$.
\kduk
\vspace{-3mm}
\begin{exer}\label{ex_ustne}
Show that $M$ is a~nonempty finite set of positive real numbers.
\end{exer}

\section[Arithmetic of limits. Limits and order]{Arithmetic of limits. Limits and order}\label{podkap_vlLimFci}

We extend some results on limits of sequences to limits of functions. 

\medskip\noindent
{\em $\bullet$ Limits of monotone functions. }Let $f\in\mathcal{F}(M)$ 
and $X$ be any set. We say that $f$ \underline{weakly increases\index{function!weakly 
increases|emph}}, resp. \underline{weakly 
decreases\index{function!weakly decreases|emph}}, on $X$ if 
$$
\text{for 
any $x\le y$ in $X\cap M$ we have $f(x)\le f(y)$, resp. $f(x)\ge f(y)$}\,. 
$$
If $f$ weakly 
increases or weakly 
decreases on $X$, then it 
is 
\underline{monotone\index{function!monotone|emph}} on $X$. Note that $X$ 
need not be a~subset of $M$.

\begin{thm}[limits of monotone functions]\label{thm_limMonFce}
Let\index{theorem!limits of monotone functions|emph}\underline{\index{limit of a function!monotone f.|emph}} 
$f$ be in $\mathcal{F}(M)$. The following holds.
\begin{enumerate}
\item If $b\in L^-(M)$ and there is a~$\theta$ such that $f$ weakly 
increases on $P^-(b,\theta)$ \underline{then}
$$
\lim_{x\to b^-}f(x)=\sup(f[P^-(b,\,\theta)])\,.
$$
\item If $+\infty\in L(M)$ and there is a~$\theta$ such that $f$ weakly 
increases on $U(+\infty,\theta)$ \underline{then} 
$$
\lim_{x\to+\infty}f(x)=\sup(f[U(+\infty,\,\theta)])\,.
$$
\end{enumerate}
Both suprema are taken in the linear order $(\R^*,<)$.
\end{thm}
\duk
1. Let $f$, $M$, $b$ and $\theta$ be as stated, and let $\ep$ be given. We set 
$$
A\equiv\sup(f[P^-(b,\,\theta)])
$$ 
and take any $a\in U(A,\ep)$ with $a<A$. By the definition of supremum there is 
a~$c\in P^-(b,\theta)\cap M$ such that $a<f(c)\le A$. We set $\de\equiv 
b-c$. For every $d\in M$ with $c<d<b$ it holds that $a<f(c)\le f(d)\le A$. 
Hence, by Exercise~\ref{ex_nekdVlOkoli}, 
$f(d)\in U(A,\ep)$. Thus 
$$
f[P^-(b,\,\de)]\sus U(A,\,\ep)
$$ 
and $\lim_{x\to b^-}f(x)=A$. 

2. Let $f$, $M$ and $\theta$ be as stated, and let $\ep$ be given. We set
$$
A\equiv\sup(f[U(+\infty,\,\theta)])
$$ 
and take any 
$a\in U(A,\ep)$ with $a<A$. By the definition of supremum there is 
a~$c\in U(+\infty,\theta)\cap M$ such that $a<f(c)\le A$. We set 
$\de\equiv\frac{1}{c}$. For every $d\in M$ with $c<d$ it holds that 
$a<f(c)\le f(d)\le A$. Using Exercise~\ref{ex_nekdVlOkoli} we get 
that $f(d)\in U(A,\ep)$. Hence 
$$
f[U(+\infty,\,\de)]\sus U(A,\,\ep)
$$ 
and $\lim_{x\to+\infty}f(x)=A$.
\kduk

\noindent
The theorem is not valid for two-sided limits: the function $\sgn\cc\R\to
\{-1,0,1\}$ weakly increases on $\R$ but $\lim_{x\to0}\sgn\,x$ 
does not exist. We find two-sided limits of monotone functions by reducing them via 
Proposition~\ref{prop_ojednLimi} to one-sided limits. These we compute by means of
the previous theorem and the next exercise.

\begin{exer}\label{ex_variantyLimMonFce}
Describe further variants of the theorem: for locally weakly 
decreasing functions and/or the right-sided limit at $b$, resp. 
at $-\infty$. 
\end{exer}

\noindent
{\em $\bullet$ Theorem~I in Section~2.3 of Cauchy's {\em Cours d'analyse \cite{brad_sand}}. }We present this theorem in a~corrected form. In Proposition~\ref{prop_couExa} we gave counterexamples to the
original version. Our correction consists in adding the assumption of monotonicity. 

\begin{thm}[Cauchy's Theorem~I corrected]\label{thm_CauchyFixed}
Let\index{theorem!Cauchy's Theorem~I corrected|emph} $k\in\R$ and let
$f$ in $\mathcal{F}((0,+\infty))$ be a~monotone function such that 
$$
\lim_{x\to+\infty}(f(x+1)-f(x))=k\,.
$$
\underline{Then}
$$
\lim_{x\to+\infty}{\textstyle
\frac{f(x)}{x}=k\,.
}
$$
\end{thm}
\duk
We begin with the case $k=0$ and weakly increasing function $f$. The
case with weakly decreasing $f$ is similar. Let an $\ep$ be given. We take an $h>0$ such that for every $x\ge h$,
$$
-\ep\le f(x+1)-f(x)\le\ep\,.
$$
Let $n\in\N$ and $x=h,h+1,\ds,h+n-1$. By summing these inequalities and rearranging the result we get
$${\textstyle
\frac{f(h)}{n}-\ep\le\frac{f(h+n)}{n}\le\frac{f(h)}{n}+\ep\,.
}
$$
Let $c\equiv\max(\frac{|f(h)|}{\ep},\frac{1}{\ep})+1$ and $x\in\R$ be such that 
$x\ge h+c$. We take the unique $n\in\N$ such that $h+n\le 
x<h+n+1$. Then $n\ge c-1$ and
$\frac{|f(h)|}{n}\le\ep$. By the above inequalities and since $f$ weakly increases, 
$$
{\textstyle
-2\ep\le\frac{f(h)}{n}-\ep\le
\frac{f(h+n)}{n}\le
\frac{f(x)}{n}\,\text{ and }\,
\frac{f(x)}{n+1}\le
\frac{f(h+n+1)}{n+1}\le
\frac{f(h)}{n+1}+\ep\le2\ep\,.
}
$$
From $0<\frac{n}{x}\le 1$ and $0<\frac{n+1}{x}\le1+\frac{1}{c}\le1+\ep$ we get
$${\textstyle
-2\ep\le-2\ep\cdot\frac{n}{x}\le\frac{f(x)}{x}\le2\ep\cdot\frac{n+1}{x}
\le2\ep(1+\ep)\,.
}
$$
Thus $\lim_{x\to+\infty}\frac{f(x)}{x}=0=k$. 

We suppose that $k<0$ and (hence) that the function $f$ weakly 
decreases. The case of $k>0$ and weakly increasing $f$ is similar. 
Let an $\ep<\min(-\frac{k}{2},1)$ be given. We take an $h>0$ such that for every $x\ge h$, 
$$
k-\ep\le f(x+1)-f(x)\le k+\ep\,.
$$
Let $n\in\N$ and $x=h,h+1,\ds,h+n-1$. By summing these inequalities and rearranging the result we get
$$
{\textstyle
k-\ep+\frac{f(h)}{n}\le\frac{f(h+n)}{n}\le k+\ep+\frac{f(h)}{n}\,.
}
$$
Let 
$${\textstyle
c\equiv\max\big(\big\{\frac{|f(h)|}{\ep},\frac{1}{\ep},\frac{h+1}{\ep}\big\}\big)+1
}
$$ 
and $x\in\R$ be such that 
$x\ge h+c$. We take the unique $n\in\N$ such that $h+n\le x<h+n+1$. Then $n\ge c-1$ and
$\frac{|f(h)|}{n}\le\ep$. By the above inequalities and since $f$ weakly decreases,
\begin{eqnarray*}
&&{\textstyle
k-2\ep\le k-\ep+\frac{f(h)}{n+1}\le
\frac{f(h+n+1)}{n+1}\le
\frac{f(x)}{n+1}
}\\
&&{\textstyle
\text{and }\frac{f(x)}{n}\le
\frac{f(h+n)}{n}\le
k+\ep+\frac{f(h)}{n}\le k+2\ep\,.
}
\end{eqnarray*}
From $0<\frac{n+1}{x}\le1+
\frac{1}{c}\le1+\ep$, $k+2\ep<0$ and $\frac{n}{x}\ge1-\frac{h+1}{c}\ge1-\ep>0$ we get
$${\textstyle
(1+\ep)(k-2\ep)\le\frac{n+1}{x}(k-2\ep)\le
\frac{f(x)}{x}\le\frac{n}{x}(k+2\ep)
\le(1-\ep)(k+2\ep)\,.
}
$$
Thus 
$\lim_{x\to+\infty}\frac{f(x)}{x}=k$.
\kduk
\vspace{-3mm}
\begin{exer}\label{ex_uloNaCauch}
Extend this theorem to infinite limits $k=\pm\infty$.    
\end{exer}

\noindent
{\em $\bullet$ Arithmetic of limits of functions. }We extend some results on arithmetic 
of limits of sequences to functions. In proofs we use Heine's
definition of limits of functions.

\begin{thm}[arithmetic of functional limits]\label{thm_AritLimFce}
Suppose that\index{theorem!arithmetic of limits of functions|emph}\underline{\index{limit of a function!arithmetic of|emph}}
$f,g\in\mathcal{R}$, $A\in L(M(f)\cap M(g))$, $\lim_{x\to A}f(x)=K$ and $\lim_{x\to 
A}g(x)=L$. \underline{Then}
$$
\text{$\lim_{x\to A}(f+g)(x)=K+L$, $\lim_{x\to A}(fg)(x)=KL$ 
and $\lim_{x\to A}(f/g)(x)=
{\textstyle\frac{K}{L}}$}\,,
$$
if the right-hand side is not indefinite.
\end{thm}
\duk
Let $f$, $g$, $A$, $K$ and $L$ be as stated. We only consider the limit of the ratio of two functions,
proofs for sum and product are similar and easier. We assume that $\frac{K}{L}$ is not indefinite. Then 
$L\ne0$ and  $A\in L(M(f/g))$ (Exercise~\ref{ex_AjeLB}). Let 
$(a_n)\sus M(f/g)\setminus\{A\}$ be any sequence with $\lim a_n=A$. The 
implication $\Rightarrow$ in Heine's definition of limits of functions gives that
$\lim f(a_n)=K$ and $\lim\,g(a_n)=L$. Using 
Theorem~\ref{thm_ari_lim} we get that
$${\textstyle
\lim\,\frac{f(a_n)}{g(a_n)}=
\frac{\lim\,f(a_n)}{\lim\,g(a_n)}=\frac{K}{L}\,.
}
$$
Since for every sequence $(a_n)$ as above the sequence
$${\textstyle
\big(\frac{f(a_n)}{g(a_n)}\big)=((f/g)(a_n))
}
$$ 
has this limit, the 
implication $\Leftarrow$ in Heine's definition of limits of functions 
gives that also 
$\lim_{x\to A}(f/g)(x)=\frac{K}{L}$.
\kduk
\vspace{-3mm}
\begin{exer}\label{ex_AjeLB}
In the previous proof, why for $L\ne0$ is $A\in L(M(f/g))$?  
\end{exer}
Using Proposition~\ref{prop_ojednLimi2}
we easily adapt the previous theorem for one-sided limits.   

\begin{exer}\label{ex_LimFunReci0}
Deduce from the theorem the next corollary. 
\end{exer}

\begin{cor}[limits of reciprocals~1]\label{cor_LimFunReci0}
If $g\in\mathcal{R}$, $\lim_{x\to A}g(x)=B$ and $B\ne0$ \underline{then} 
$${\textstyle
\lim_{x\to A}(k_1/g)(x)=\lim_{x\to A}\frac{1}{g(x)}=\frac{1}{B}\,.
}
$$
\end{cor}

\noindent
{\em $\bullet$ Limits of functions and the linear order $(\R^*,<)$. }Recall 
that for sets $M,N\sus\R$ the notation $M<N$ means that for every $a\in M$ 
and $b\in N$ we have $a<b$. Recall that for any function $f$ and any set $X$, 
$$
f[X]=f[X\cap M(f)]=\{f(x)\cc\;x\in X\cap M(f)\}\,.
$$
In the next theorem and proposition we have $f,g\in\mathcal{R}$.

\begin{thm}[limits versus order~2]\label{thm_limFceUsp}
We assume that\index{theorem!limits versus order~2|emph}\underline{\index{limit of a function!vs. order|emph}} 
$\lim_{x\to A}f(x)=K$ and $\lim_{x\to B}g(x)=L$, where possibly $A\ne B$. The following holds.
\begin{enumerate}
\item If $K<L$ \underline{then} there is a~$\de$ such that $f[P(A,\de)]<g[P(B,\de)]$.
\item If for every $\de>0$ there exist an $x\in P(A,\de)\cap M(f)$ and a~$y\in P(B,\de)\cap M(g)$ such that $f(x)\ge g(y)$, \underline{then} $K\ge L$.
\end{enumerate}
\end{thm}
\duk
1. Since $K<L$, by Exercise~\ref{ex_ulohaNaokoli1} there is an $\ep$ such that $U(K,\ep)<U(L,\ep)$. Then by the assumption there is a~$\de$ such that 
$f[P(A,\de)]\sus U(K,\ep)$ and $g[P(B,\de)]\sus U(L,\ep)$. Hence
$f[P(A,\de)]<g[P(B,\de)]$.

2. Part~2 is the contrapositive of the implication in part~1.
\kduk

\noindent
We strengthen the theorem
in the same way as Proposition~\ref{prop_limAuspo} 
strengthens Theorem~\ref{thm_limAuspo}.

\begin{exer}\label{ex_zesLimFceUsp}
Prove the following proposition.    
\end{exer}

\begin{prop}[strengthening Theorem~\ref{thm_limFceUsp}]\label{prop_zesLimFceUsp}
Let $\lim_{x\to A}f(x)=K$ and $\lim_{x\to B}g(x)=L$, where possibly $A\ne 
B$. \underline{Then} the following hold.
\begin{enumerate}
\item If $K<L$ then there exist $\de$ and $a,b\in\R$, such that 
$$
f[P(A,\,\de)]<\{a\}<\{b\}
<g[P(B,\,\de)]\,.
$$
\item If for every $\de$ and every $a<b$ in $\R$ there is an 
$x\in P(A,\de)\cap M(f)$ and a~$y\in P(B,\de)\cap M(g)$ such that $f(x)\ge 
a$ or $g(y)\le b$, then $K\ge L$.
\end{enumerate}
\end{prop}

\begin{exer}\label{ex_verzeProJednostr}
State variants of Theorem~\ref{thm_limFceUsp} and
Proposition~\ref{prop_zesLimFceUsp} for one-sided limits and prove them.  
\end{exer}

\noindent
{\em $\bullet$ The squeeze theorem. }We generalize the squeeze theorem and the corollary of it from sequences to functions. 

\begin{thm}[unexpected squeeze~2]\label{thm_unexSque}
Suppose that\index{theorem!unexpected squeeze 2|emph} 
$f,g\in\mathcal{F}(M)$, that  
$\lim_{x\to A}f(x)=b$ is in $\R$ and that
$\lim_{x\to A}|f(x)-g(x)|=0$. 
\underline{Then} 
$$
\lim_{x\to A}g(x)=b\,.
$$
\end{thm}
\duk
Let $f$, $g$, $M$, $A$ ($\in L(M)$) and $b$ be as stated, and let an $\ep$ be given. We take a~$\de$ such that $f[P(A,\de)]\sus U(b,\frac{\ep}{2})$ and 
$(f-g)[P(A,\de)]\sus U(0,\frac{\ep}{2})$. Then for every $x\in P(A,\de)\cap M$ we have
$${\textstyle
|g(x)-b|\le |g(x)-f(x)|+|f(x)-b|\le
\frac{\ep}{2}+\frac{\ep}{2}=\ep\,.
}
$$
Hence $g[P(A,\de)]\sus U(b,\ep)$ and $\lim_{x\to A}g(x)=b$.
\kduk

Recall that
$I(a,b)=\{x\in\R\cc\;\min(a,b)\le x\le\max(a,b)\}$. 

\begin{cor}[standard squeeze~2]\label{cor_limDvaFstraz}
Suppose that
$f, g,h\in\mathcal{F}(M)$, that $\lim_{x\to K}f(x)=\lim_{x\to 
K}g(x)=b$ is in $\R$ and that for every $x\in M$ we have
$h(x)\in I(f(x),g(x))$. \underline{Then} 
$$
\lim_{x\to K}h(x)=b\,.
$$
\end{cor}
\duk
Let $f$, $g$, $h$, $M$, $K$ and $b$ be as stated. Then
$\lim_{x\to K}|f(x)-g(x)|=0$. Since for every $x\in M$ we have
$$
|f(x)-h(x)|,\,|g(x)-h(x)|\le
|f(x)-g(x)|\,,
$$
we get $\lim_{x\to K}|f(x)-h(x)|=0$ and 
$\lim_{x\to K}|g(x)-h(x)|=0$, and are done by
Theorem~\ref{thm_unexSque}.
\kduk

\noindent
The advantage of Theorem~\ref{thm_unexSque} over Corollary~\ref{cor_limDvaFstraz} is that
the theorem easily generalizes to maps between metric spaces.

\begin{exer}\label{ex_genToMS}
State and prove this
generalization. 
\end{exer}

\section[Limits  of composite functions]{Limits  of composite 
functions}\label{sec_limSloFce}

Composition\index{limit of a 
function!composite f.|(} 
of functions has no analogue for sequences, and so its interplay with limits of functions yields a~genuinely new theorem. It is usually stated as an 
implication, but we give it a~form of an equivalence.

\medskip\noindent 
{\em $\bullet$ Limits of composite functions. }Let  $f,g\in\mathcal{R}$. 
Recall that the composite function $f(g)\cc M(f(g))\to\R$ has the definition domain 
$$
M(f(g))=\{x\in M(g)\cc\;g(x)\in M(f)\}\,. 
$$
It may be a~proper subset of $M(g)$.

\begin{thm}[limits of composites]\label{thm_LimSlozFunkce}
We assume that\underline{\index{limit of a function!composite f.|emph}}
\index{theorem!limits of composites|emph}
$$
\lim_{x\to A}g(x)=K,\ \lim_{x\to K}f(x)=L\,\text{ and that }\,A\in L(M(f(g)))\,. 
$$
\underline{Then} 
$$
\lim_{x\to A}f(g)(x)=L
$$ 
if and only if one of two conditions holds.
\begin{enumerate}
\item Either $K\not\in M(f)$, or $K\in M(f)$ and $f(K)=L$.
\item For some $\theta$ we have $K\not\in g[P(A,\theta)]$.
\end{enumerate}
If neither condition holds, then $K\in M(f)$, $f(K)\ne L$ and $\lim_{x\to A}
f(g)(x)$ does not exist or equals $f(K)$.
\end{thm}
\duk
Let $A$, $g$, $K$, $f$ and $L$ be as stated and an $\ep$ be given. By the 
assumption there is a~$\de'$ such that (a) $f[P(K,\de')]\sus 
U(L,\ep)]$, and  a~$\de$ 
such that (b) $g[P(A,\de)]\sus U(K,\de')$. Suppose that condition~1 
holds. Then inclusion (a) strengthens to $f[U(K,\de')]\sus U(L,\ep)$ and  
$$
f(g)[P(A,\,\de)]=f[\,g[P(A,\,\de)]\,]\sus 
f[U(K,\,\de')]\sus U(L,\,\ep)\,.
$$
Hence $\lim_{x\to A}f(g)(x)=L$. Suppose that condition~2 holds. We 
can take the previous $\de$ such that $\de\le\theta$, where $\theta$ 
is as in condition~2. Then inclusion (b) strengthens to $g[P(A,\de)]\sus P(K,\de')$
and
$$
f(g)[P(A,\,\de)]=f[\,g[P(A,\,\de)]\,]\sus f[P(K,\,\de')]\sus U(L,\,\ep)\,.
$$
Hence again $\lim_{x\to A}f(g(x))=L$.

Suppose that neither condition holds. Hence $K$ is in $M(f)$ but 
$f(K)\ne L$, and for every $n$ there is an $a_n\in P(A,\frac{1}{n})\cap 
M(g)$ such that $g(a_n)=K$. Then $(a_n)\sus M(f(g))\setminus\{A\}$, 
$\lim\,a_n=A$ and
$$
\lim\,f(g)(a_n)=\lim\,f(g(a_n))=\lim\,f(K)=
f(K)\ \ (\ne L)\;.
$$
By Theorem~\ref{thm_HeinehoDef} the limit $\lim_{x\to 
A}f(g(x))$ either does not exist or equals to $f(K)$ (which is not $L$).
\kduk

\noindent
Condition~1 is satisfied if $K=\pm\infty$. Condition~2
is satisfied if the function $g$ is injective. We get the following 
corollary of Theorem~\ref{thm_LimSlozFunkce}.

\begin{cor}[when composite limit works]\label{cor_vzdyFung}
Let 
$$
\lim_{x\to A}g(x)=K,\ \lim_{x\to K}f(x)=L\,\text{ and }\,A\in L(M(f(g)))\,.
$$
If $K=\pm\infty$ or $g$ is injective, \underline{then}
$$
\lim_{x\to A}f(g)(x)=L\,.
$$
\end{cor} 

\begin{exer}\label{ex_slozFceHeinem}
Prove Theorem~\ref{thm_LimSlozFunkce} by means of Heine's definition of limits of 
functions.
\end{exer}

\noindent
{\em $\bullet$ Applications of Theorem~\ref{thm_LimSlozFunkce}. }The next two equivalences are
useful tools to determine limits of functions. 

\begin{cor}[shift to $0$]\label{cor_substPosun}
Let $f\in\mathcal{R}$ and $b\in L(M(f))$. \underline{Then} 
$$
\lim_{x\to b}f(x)=L\iff \lim_{x\to0}f(x+b)(x)=L\,.
$$
\end{cor}
\duk
Implication $\Rightarrow$. We assume that $\lim_{x\to b}f(x)=L$ and apply
Corollary~\ref{cor_vzdyFung} to the composite limit $\lim_{x\to0}
f(g)(x)$ where $g(x)\equiv x+b$ is injective. 

Implication $\Leftarrow$. We assume that $\lim_{x\to0}f(x+b)(x)=L$ and
apply Corollary~\ref{cor_vzdyFung} to the composite limit $\lim_{x\to b}g(h(x))$ where $g(x)\equiv f(x+b)$ and $h(x)\equiv x-b$ is injective.
\kduk

\begin{cor}[shift to $0^{\pm}$]\label{cor_ekvivZde}
Let $f\in\mathcal{R}$ and $\pm\infty\in L(M(f))$. 
\underline{Then} 
$${\textstyle
\lim_{x\to\pm\infty}f(x)=L\iff 
\lim_{x\to 0^{\pm}}
f(\frac{1}{x})(x)=L\,, 
}
$$
with equal signs. Note that the last limit is one-sided.
\end{cor}
\duk
We assume that all signs are $+$, 
the case with all signs $-$ is similar. Implication $\Rightarrow$. 
We assume that $\lim_{x\to+\infty}f(x)=L$ and apply
Corollary~\ref{cor_vzdyFung} to the composite limit $\lim_{x\to0}f(g(x)$ 
where $g(x)\equiv\frac{1}{x}\,|\,(0,+\infty)$ is injective.

Implication $\Leftarrow$. We assume that $\lim_{x\to0^+}f(\frac{1}{x})(x)=L$ and apply
Corollary~\ref{cor_vzdyFung} to the composite limit
$\lim_{x\to+\infty}g(h)(x)$ where
$g(x)\equiv f(\frac{1}{x})\,|\,(0,+\infty)$ and
$h(x)\equiv\frac{1}{x}\,|\,(0,+\infty)$ is injective.
\kduk
\vspace{-3mm}
\begin{exer}\label{ex_druhyDkPrevHod}
If $p(x)$ is 
a~non-constant polynomial, $p(a)=b$ and $f\in\mathcal{R}$ has the limit 
$\lim_{x\to b}f(x)=L$, then 
$$
\lim_{x\to a}f(p)(x)=L\,.
$$
\end{exer}
\index{limit of a function!composite f.|)}

\section[Limits of inverse functions]{Limits of inverse functions}\label{sec_limInvFci}

In this original extending section we investigate the interplay of 
functional limits and inverses of functions.

\medskip\noindent
{\em $\bullet$ Why is there nothing in
the literature about limits of inverse functions? }Maybe because it is clear that if $f\in\mathcal{R}$
is an injective function and 
$$
\lim_{x\to A}f(x)=B\,,
$$
then in general it does \underline{not} follow that 
$$
\lim_{y\to B}f^{\langle-1\rangle}(y)=A\,. 
$$
Indeed, the former limit means that all points close to $A$ are sent by $f$ 
to points close to $B$. However, 
also some points not close to $A$ may be sent by $f$ to points close to $B$. If this happens, the latter
inverse limit does not exist.

\begin{exer}\label{ex_simpForThis}
Describe a~simple example of such situation.    
\end{exer}
But even in the general situation we can say something nontrivial.

\begin{prop}[limits of inverses]\label{prop_onLimInvFun}
Suppose that $f\in\mathcal{F}(M)$ is an injective function, $A\in L(M)$ and that 
$$
\lim_{x\to A}f(x)=B\,.
$$
Then the following is true.
\begin{enumerate}
\item We have $B\in L(f[M])$ and 
$$
\lim_{y\to B}f^{\langle-1\rangle}(y)=A\,, 
$$
if this limit exists.
\item There exists a~set $N\sus M$ such that $B\in L(f[N])$ and 
$$
\lim_{y\to B}(f\,|\,N)^{\langle-1\rangle}(y)=A\,, 
$$
\end{enumerate}
\end{prop}
\duk
Let $f$, $M$ and $A$ be as stated. 1. We take any sequence $(a_n)\sus 
M\setminus\{A\}$ with $a_n\to A$ (which exists because $A\in L(M)$) 
and consider the sequence $(b_n)\equiv(f(a_n))$ ($\sus 
f[M]$). By Theorem~\ref{thm_HeinehoDef} we 
have $b_n\to B$. The injectivity of $f$ implies that $b_n=B$ for at 
most one $n$ and we get that $B\in L(f[M])$. We may assume that
$(b_n)\sus f[M]\setminus\{B\}$. 
Since $b_n\to B$ and 
$$
f^{\langle-1\rangle}(b_n)=a_n\to A\,,
$$
Theorem~\ref{thm_HeinehoDef} implies that the limit of $f^{\langle-1\rangle}(y)$ at $B$, if it exists, has to be $A$.

2. Just set 
$$
N\equiv\{a_n\cc\;n\in\N\}
$$
where $(a_n)\sus M\setminus\{A\}$ is any sequence with $a_n\to A$.
\kduk

\noindent
We present three families of theorems on limits of inverse 
functions. In the first two we present conditions under which the limit $\lim_{x\to A}f(x)=B$ does imply the limit $\lim_{y\to B}f^{\langle-1\rangle}(y)=A$.

\medskip\noindent
{\em $\bullet$ First theorem on limits of inverses. }One approach is to strengthen the 
definition of functional limit.

\begin{defi}[strong functional limit]\label{def_strongLim}
Let $f\in\mathcal{F}(M)$, $A\in L(M)$ and $B\in\R^*$. If for every $\ep$ there exist $\de$ and $\theta$ such that
$$
P(B,\,\theta)\cap f[M]\sus
f[P(A,\,\de)]\sus U(B,\,\ep)\,,
$$
we say that the \underline{strong limit\index{strong limit of 
a~function|emph}} of $f$ at $A$ is $B$ and we write $(\lim_{x\to A})f(x)=B$.\label{slimFun}
\end{defi}
The standard limit in 
Definition~\ref{def_limFce} is strengthened by adding the first inclusion.

\begin{thm}[limits of inverses~1]\label{thm_LimInvFun1}
Suppose\index{theorem!limits of inverses~1|emph} 
that $f\in\mathcal{F}(M)$ is an injective function, $A\in L(M)$, $B\in\R^*$ and that 
$$
(\lim_{x\to A})f(x)=B\,.
$$
\underline{Then} 
$$
\lim_{y\to B}f^{\langle-1\rangle}(y)=A\,. 
$$
\end{thm}
\duk
Let $f$, $M$, $A$ and $B$ be as stated. We have $B\in L(f[M])$ by 
part~1 of Proposition~\ref{prop_onLimInvFun}. Suppose for the contrary that the 
limit of $f^{\langle-1\rangle}(y)$ at $B$ is not $A$. Then there is 
a~sequence $(b_n)\sus f[M]\setminus\{B\}$ such that 
$$
b_n\to B\wedge
a_n\equiv f^{\langle-1\rangle}(b_n)\to A'\ne A\,.
$$
We take a~$\de_0$ such that 
$$
U(A,\,\de_0)\cap U(A',\,\de_0)=\emptyset\,.
$$
Since $A\in L(M)$ and the limit of $f(x)$ at $A$ is $B$, we take 
a~sequence $a_n'\sus M\setminus\{A\}$ such that $a_n'\to A$ and 
(by Theorem~\ref{thm_HeinehoDef}) $f(a_n')\to B$. We take an index $m\in\N$ such that $a_m'\in 
P(A,\de_0)$ and $f(a_m')\ne B$. Then we take an $\ep$ such that
$$
f(a_m')\not\in U(B,\,\ep)\,.
$$
Now for this $\ep$ if $\de$ is such that 
$$
f[P(A,\,\de)]\sus U(B,\,\ep)\,,
$$
then $\de<\de_0$. But for every $\de<\de_0$ and every $\theta$ the inclusion
$$
P(B,\,\theta)\cap f[M]\sus f[P(A,\,\de)]
$$
does not hold because $f$ is injective and
$$
\text{$a_n\in U(A',\,\de_0)$, hence $a_n\not\in P(A,\,\de)$, and 
$f(a_n)=b_n\in P(B,\,\theta)\cap f[M]$}
$$
for every large $n$. Thus it is not true that the strong limit of 
$f(x)$ at $A$ is $B$, contrary to the assumption. 
\kduk
\vspace{-3mm}
\begin{exer}\label{ex_fromProof}
Why is $a_n\not\in P(A,\de)$ for $\de<\de_0$?    
\end{exer}

\noindent
{\em $\bullet$ Second theorems on limits of inverses. }Another 
approach is to restrict the function $f(x)$. By 
a~\underline{strongly 
positive\index{function!strongly positive|emph}} function we mean 
any function
$$
g\cc[0,\,+\infty)\to[0,\,+\infty)
$$
such that $g(0)=0$ and that for every $c>0$ we have
$$
\inf\big(g[(c,\,+\infty)]\big)>0\,.
$$

\begin{defi}[lower regulated functions]\label{def_regFun}
We say that a~function $f$ in $\mathcal{F}(M)$ is 
\underline{lower regulated\index{function!lower regulated|emph}\index{lower regulated function|emph}} if there is a~strongly positive function 
$g(x)$ such that for every two points $x,y\in M$,
$$
|f(x)-f(y)|\ge g(|x-y|)\,.
$$
\end{defi}
For example, a~function $f\in\mathcal{R}$ is lower regulated if there is a~constant $c>0$ such that
$$
|f(x)-f(y)|\ge c|x-y|
$$ 
for every $x,y\in M(f)$. 

\begin{exer}\label{ex_isInje}
Every lower regulated function is injective.    
\end{exer}

\begin{exer}\label{ex_lowFRegFunc}
If $f\cc(a,b)\to\R$ is differentiable on $(a,b)$, 
$|f'|\ge c>0$ on $(a,b)$ and $M\sus(a,b)$, then the restriction $f\,|\,M$ is lower regulated.
\end{exer}

\begin{thm}[limits of inverses 2a]\label{thm_LimInvFun2}
Let\index{theorem!limits of inverses~2a|emph} 
$f\in\mathcal{F}(M)$, $A\in L(M)$, $b\in\R$ and let
$$
\lim_{x\to A}f(x)=b\,.
$$
If the function $f$ is lower regulated, \underline{then}
$$
\lim_{y\to b}f^{\langle-1\rangle}(y)=A\,.
$$
\end{thm}
\duk
Let $f$, $M$, $A$ and $b$ be as stated. In view of part~1 of 
Proposition~\ref{prop_onLimInvFun} and of 
Theorem~\ref{thm_HeinehoDef} it suffices to show that
$$
\text{if $(b_n)\sus f[M]\setminus\{b\}$ satisfies $b_n\to b$, then $f^{\langle-1\rangle}(b_n)\to A$}\,.
$$
We assume for the contrary that a~sequence $(b_n)\sus f[M]
\setminus\{b\}$ has $\lim b_n=b$ but $f^{\langle-1\rangle}(b_n)\not\to A$. Passing to 
a~subsequence and using part~2 of Theorem~\ref{thm_oPodposl} we 
may assume that $\lim f^{\langle-1\rangle}(b_n)=A'$ with $A'\ne A$. We denote 
$a_n\equiv f^{\langle-1\rangle}(b_n)$. Using the assumption, we 
take any sequence $(a_n')\sus M\setminus\{A\}$ with $a_n'\to A$. 
By the assumption and Theorem~\ref{thm_HeinehoDef}, $b_n'\equiv f(a_n')\to b$. There is a~constant $c>0$ such that 
$$
|a_n-a_n'|\ge c\ \ (>0)
$$
for every large $n$. But, since both 
$f(a_n)\to b$ and $f(a_n')\to b$, we have
$$
\lim_{n\to\infty}|f(a_n)-f(a_n')|=0\,.
$$
Setting 
$x\equiv a_n$ and $y\equiv a_n'$ for sufficiently large $n$, 
we get a~contradiction with the assumption that $f$ is lower regulated.
\kduk

We obtain a~similar theorem in the case when $\lim_{x\to A}f(x)=\pm\infty$. By 
an \underline{almost 
bounded\index{function!almost 
bounded|emph}} function we mean 
any function
$$
g\cc[0,\,+\infty)\to(0,\,+\infty)
$$
such that for every $c>0$ we have
$$
\sup\big(g[(c,\,+\infty)]\big)<+\infty\,.
$$

\begin{defi}[upper regulated functions]\label{def_UpregFun}
We say that a~function $f$ in $\mathcal{F}(M)$ is 
\underline{upper regulated\index{function!upper regulated|emph}\index{upper regulated function|emph}} if there is an almost bounded function 
$g(x)$ such that for every two points $x,y\in M$,
$$
|f(x)-f(y)|\le g(|x-y|)\,.
$$
\end{defi}

\begin{exer}\label{ex_UppFRegFunc}
If $f\cc(a,b)\to\R$ is differentiable on $(a,b)$, 
$|f'|\le c$ on $(a,b)$ and $M\sus(a,b)$, then the restriction $f\,|\,M$ is upper regulated.
\end{exer}

\begin{thm}[limits of inverses 2b]\label{thm_LimInvFun2b}
We\index{theorem!limits of inverses~2b|emph} 
assume that $f\in\mathcal{F}(M)$, $A$ is in $L(M)$, $B\in\{-\infty,+\infty\}$ and that
$$
\lim_{x\to A}f(x)=B\,.
$$
If the function $f$ is injective and upper regulated, \underline{then}
$$
\lim_{y\to B}f^{\langle-1\rangle}(y)=A\,.
$$
\end{thm}
\duk
Let $f$, $M$, $A$ and $B$ be as stated. We proceed as in the 
previous proof. We bring to contradiction the assumption that there is a~sequence $(b_n)$ such 
that $b_n\to B$ but 
$$
\lim f^{\langle-1\rangle}(b_n)=A'\ne A\,.
$$
We denote 
$a_n\equiv f^{\langle-1\rangle}(b_n)$. Using the assumption, we 
take any sequence $(a_n')\sus M\setminus\{A\}$ with $a_n'\to A$. 
By the assumption and Theorem~\ref{thm_HeinehoDef}, $b_n'\equiv f(a_n')\to B$. There is a~constant $c>0$ such that 
$$
|a_n-a_n'|\ge c\ \ (>0)
$$
for every large $n$. We use Exercise~\ref{ex_diifToInf}, 
pass to a~subsequence and may assume that 
$$
\lim|f(a_n)-f(a_n')|=+\infty\,.
$$
Setting 
$x\equiv a_n$ and $y\equiv a_n'$ for sufficiently large $n$, 
we get a~contradiction with the assumption that $f$ is upper regulated.
\kduk
\vspace{-3mm}
\begin{exer}\label{ex_diifToInf}
If $\lim a_n=\lim b_n=\pm\infty$ then $(b_n)$ has a~subsequence 
$(c_n)$ such that $\lim|a_n-c_n|=+\infty$.    
\end{exer}

\medskip\noindent
{\em $\bullet$ Third theorems on limits of inverses. }In these results we
associate with a~given function $f\in\mathcal{R}$ 
a~relation between the sets $L(M(f))$ and $L(f[M(f)])$.

\begin{defi}[relation $\mathcal{L}(f)$]\label{def_relLf}
Let $f\in\mathcal{F}(M)$. We define a~binary relation 
$$
\mathcal{L}(f)\sus L(M)\times L(f[M])\label{LfM}
$$ 
by putting $(A,B)$ in $\mathcal{L}(f)$ $\iff$ there exists a~sequence $(a_n)\sus M\setminus\{A\}$ such that $a_n\to A$ and $f(a_n)\to B$. 
\end{defi}

For a~relation $R\sus X\times Y$ we call an element $x\in X$, resp. 
$y\in Y$, $R$-\underline{isolated\index{relation@(
binary) relation!isolated element of|emph}} if there is no pair 
$(x,y')\in R$, resp. $(x',y)\in R$. 
We call the element $R$-\underline{unique\index{relation@(binary) relation!unique element 
of|emph}} if there is 
exactly one such pair.

\begin{prop}[relation $\mathcal{L}(f)$]\label{prop_onLf}
Let $f\in\mathcal{F}(M)$ be an injective function. The relation  $\mathcal{L}(f)$ has two properties.
\begin{enumerate}
\item No element $A\in L(M)$, resp. $B\in L(f[M])$, is $\mathcal{L}(f)$-isolated.
\item An element $A\in L(M)$, resp. $B\in L(f[M])$, is $\mathcal{L}(f)$-unique $\iff$ the limit
$$
\text{$\lim_{x\to A}f(x)$, resp. $\lim_{y\to B}f^{\langle-1\rangle}(y)$, exists}\,.
$$
\end{enumerate}
\end{prop}
\duk
Let $f$ and $M$ be as stated. 1. Let $A\in L(M)$. We take any 
sequence $(a_n)\sus M\setminus\{A\}$ with $a_n\to A$. By Theorem~\ref{thm_BolzWeier} there is a~subsequence $(a_{m_n})$
such that $\lim f(a_{m_n})=B'$. It follows that $B'\in L(f[M])$ and 
$(A,B')\in\mathcal{L}(f)$. The argument for $B$ is the same, only 
$f$ is replaced with $f^{\langle-1\rangle}$.

2. Let $A\in L(M)$. The implication $\Rightarrow$. We assume that $A$ 
is $\mathcal{L}(f)$-unique. Thus there is a~unique $B\in L(f[M])$ such that for 
a~sequence $(a_n')\sus M\setminus\{A\}$ we have $a_n'\to A$ and 
$f(a_n')\to B$. Now let $(a_n)\sus M\setminus\{A\}$ be any sequence 
with $a_n\to A$. If $f(a_n)\not\to B$ then by using part~2 of 
Theorem~\ref{thm_oPodposl} we would contradict $\mathcal{L}(f)$-uniqueness of $A$. Hence $f(a_n)\to 
B$ and $\lim_{x\to A}f(x)=B$ by 
Theorem~\ref{thm_HeinehoDef}.

The implication $\neg\Rightarrow\neg$. Now we 
assume that there exist two different elements $B,B'\in L(f[M])$
and two sequences $(a_n),(a_n')\sus M\setminus\{A\}$ such that 
$a_n,a_n'\to A$, $f(a_n)\to B$ and
$f(a_n')\to B'$. By Theorem~\ref{thm_HeinehoDef} the 
limit $\lim_{x\to A}f(x)$ does not exist.
\kduk 
\vspace{-3mm}
\begin{exer}\label{ex_proNexThm}
Prove the next theorem.    
\end{exer}

\begin{thm}[limits of inverses 3a]\label{thm_LimInvFun3a}
Let\index{theorem!limits of inverses~3a|emph} 
$f\in\mathcal{F}(M)$ be an injective function. We assume that the limit $\lim_{x\to A}f(x)$ exists for every $A\in L(M)$. \underline{Then} the relation  $\mathcal{L}(f)$ has two properties.
\begin{enumerate}
\item It is a~surjective function from 
$L(M)$ on $L(f[M])$.
\item If this function
$\mathcal{L}(f)$ is injective, then it is a~bijection and 
$$
\lim_{y\to B}f^{\langle-1\rangle}(y)=\mathcal{L}(f)^{\langle-1\rangle}(B)
$$ 
for every $B\in L(f[M])$. 
\end{enumerate}
\end{thm}
We obtain version 3b of this theorem dealing with pointwise continuity in Section~\ref{sec_aritSpoj}.

\section[Asymptotic notation]{Asymptotic notation}\label{sec_asympZnac}

Books\index{asymptotic notation|(} 
on computational complexity, or on analysis of algorithms, intersect 
with books on mathematical analysis in passages devoted to asymptotic
notation. In this section we present our version of it. In fact, what does ``asymptotic'' mean? See Definition~\ref{def_asympRel}.

\medskip\noindent
{\em $\bullet$ Asymptotic relations. }The \underline{symmetric difference\index{set operation!symmetric difference, 
$\Delta$|emph}}\label{Delta} of sets $X$ and $Y$ is
$$
X\,\Delta\,Y\equiv(X\setminus Y)\cup (Y\setminus X)\,.
$$

\begin{defi}[almost equality]\label{def_almEq}
Functions $f=\langle M(f),\R,G_f\rangle$ and $g=\langle M(g),\R,G_g\rangle$ in $\mathcal{R}$ 
are \underline{almost equal\index{almost 
equal functions, $\doteq$|emph}}, in 
symbols $f\doteq g$, if the set
$$
G_f\,\Delta\,G_g
$$ 
is finite.   
\end{defi}

\begin{exer}\label{ex_almEquaEquiv}
We have the equivalence that $f\doteq g$ iff
$$
\text{$M(f)\,\Delta\, M(g)$ and $\{x\in 
M(f)\cap M(g)\cc\;f(x)\ne g(x)\}$ are finite sets}\,.
$$
\end{exer}

\begin{exer}\label{ex_almEqu1}
The relation $\doteq$ on the set $\mathcal{R}$ is an equivalence relation.    
\end{exer} 

\begin{defi}[asymptotic relations]\label{def_asympRel}
We say that a~relation $\mathcal{A}\sus\mathcal{R}\times
\mathcal{R}$ is \underline{asymptotic\index{asymptotic relation|emph}} if for any functions $f$, $g$, $f_0$ 
and $g_0$ in $\mathcal{R}$ such that $f\doteq f_0$ and $g\doteq 
g_0$ we have the equivalence
$$
f\,\mathcal{A}\,g\iff f_0\,\mathcal{A}\,g_0\,.
$$
\end{defi}
This resembles the robustness of properties of sequences in 
Definition~\ref{def_robust}.

\medskip\noindent 
{\em $\bullet$ Asymptotic notation $O$, $\ll$ and other. }A~function 
$f\in\mathcal{R}$ is \underline{bounded\index{function!bounded|emph}} on a~set
$N$ if there is a~constant $c\ge0$ such that $|f(x)|\le c$ for every $x\in M(f)\cap N$.

\begin{defi}[$O$ and $\ll$]\label{def_velkeO}
Let $f,g\in\mathcal{R}$ and $N\sus\R$. If the function $f/g$ is bounded on the set $N$, we say that 
$$
\text{$f$ is \underline{big $O$\index{asymptotic notation!big O, $f=O(g)$ (on 
$N$)|emph}} of $g$ on $N$ and write $f=O(g)$ (on $N$)}\,.\label{O}
$$
Notation $f\ll g$ (on $N$)\index{asymptotic notation!lless or equal@$\ll$|emph} means the same. 
\end{defi}
In few exceptional cases we consider a~more general form of this notation for 
functions $f\cc M\to\C$ and $g\cc N\to\C$ with $M,N\sus\C$.

\begin{exer}\label{ex_equiBigO}
Show that $f=O(g)$ (on $N$) iff there is a~constant $c\ge0$ such 
that $|f(x)|\le c|g(x)|$  for every $x\in N\cap M(f)\cap M(g)\setminus 
Z(g)$.      
\end{exer}

\begin{prop}[$O$ is asymptotic]\label{prop_OnoChange}
Let $N\sus\R$. \underline{Then} the relation $\mathcal{A}\sus\mathcal{R}^2$,  defined by
$$
f\,\mathcal{A}\,g\iff \text{$f=O(g)$ (on $N$)}\,,
$$
is asymptotic.
\end{prop}
\duk
Let $N\sus\R$, and let $f$, $g$, $f_0$ and $g_0$ be in $\mathcal{R}$ and such that $f\doteq f_0$, $g\doteq g_0$ and $f=O(g)$ (on 
$N$). We show that also $f_0=O(g_0)$ (on $N$). By the assumption 
there is a~constant $c\ge0$ such that $\big|\frac{f(x)}{g(x)}\big|\le c$  for every $x\in M(f/g)\cap N$. It follows from the definition of the relation
$\doteq$ that the set 
$${\textstyle
X\equiv\frac{f_0}{g_0}[N]\setminus\frac{f}{g}[N]
}
$$
is 
finite (Exercise~\ref{ex_Xfini}). If $X\ne\emptyset$ then we set 
$d\equiv\max(\{|x|\cc\;x\in X\})$. For $X=\emptyset$ we set
$d\equiv1$. Then for every $x\in M(f_0/g_0)\cap N$ we 
have $\big|\frac{f_0(x)}{g_0(x)}\big|\le\max(c,d)$, as
needed.
\kduk
\vspace{-3mm}
\begin{exer}\label{ex_Xfini}
Why is the set $X$ finite?    
\end{exer}

The next exercise shows that a~simplified definition of $O$, whose versions appear in the literature, is not
asymptotic. 

\begin{exer}\label{ex_notAsym}
Suppose that we define, for $f,g\in\mathcal{R}$ and $N\sus\R$, 
that $f=O'(g)$ (on $N$) iff for some $c\ge0$ we have $|f(x)|\le c|g(x)|$ 
for every $x\in N\cap M(f)\cap M(g)$. Show that the relation $O'$ is in general not asymptotic.   
\end{exer}

\noindent
The notation 
$$
\text{$f=g+O(h)$ (on $N$)} 
$$
has the \underline{error 
form\index{asymptotic notation!error 
form|emph}} and means that 
$f-g=O(h)$ (on 
$N$). Notation like
$$
\text{$\log x=O_{\ep}(x^{\ep})$ (on $[1,+\infty)$)\index{asymptotic 
notation!big O with par@$O_{\ep}$|emph}}
$$
means that the constant $c$ in Exercise~\ref{ex_equiBigO} is a~function of $\ep$. Notation $f\gg g$ (on $N$) 
and $f=\Omega(g)$ 
(on $N$)\index{asymptotic notation!gg or equal@$\gg$, $\Omega(\cdot)$|emph} 
mean that
$g\ll f$ (on $N$). Notation $f=\Theta(g)$  (on $N$) and  $f\asymp g$ (on $N$)\index{asymptotic notation!theta@$\Theta(\cdot)$, $\asymp$|emph} 
both mean that simultaneously $f\ll g$ (on $N$) and $g\ll f$ (on $N$).\label{Omega} 

\begin{exer}\label{ex_naO}
Answer the following questions.
\begin{enumerate}
\item Is $x^2=O(x^3)$ (on  $\R\setminus(-1,1)$)?
\item Is $x^2=O(x^3)$ (on $\R$)?
\item Is $x^3=O(x^2)$ (on $\R$)?
\item Is $x^3=O(x^2)$ (on  $(-20,20)$)?
\item Is $\log x=O(x^{1/3})$ (on  $(0,+\infty)$)?
\item Is $\log x=O(x^{1/3})$ (on $(1,+\infty)$)?
\end{enumerate}
\end{exer}

\noindent 
{\em $\bullet$ Asymptotic notation $o$, $\omega$, and $\sim$. }Now the definitions employ limits of functions. 

\begin{defi}[$o$ and $\omega$]\label{def_maleo}
Let $f,g\in\mathcal{R}$ and $A\in L(M(f/g))$. If 
$$
\lim_{x\to 
A}{\textstyle\frac{f(x)}{g(x)}=0\,, 
}
$$
we say that $f$ is 
\underline{little $o$\index{asymptotic notation!little o, $f=o(g)$ ($x\to A$)|emph}} of $g$ for $x\to A$ and write 
$f(x)=o(g(x))$ ($x\to A)$. Notation $f(x)=\omega(g(x))\ (x\to 
A)$\index{asymptotic notation!little omega@$\omega(\cdot)$|emph} means 
the same.\label{o} 
\end{defi}
Like before, $f=g+o(h)$ ($x\to A$) means that $f-g=o(h)$ ($x\to A$).

\begin{defi}[$\sim$]\label{def_simi}
Let $f,g\in\mathcal{R}$ and let $A\in L(M(f/g))$. If 
$$
\lim_{x\to A}{\textstyle
\frac{f(x)}{g(x)}=1\,,
}
$$ 
we say that $f$ is \underline{asymptotically equal\index{asymptotic 
notation!asymptotic equality, $\sim$|emph}} to $g$ for $x\to A$
and write $f(x)\sim g(x)$ ($x\to A)$.
\end{defi}
For example, $x^2\sim(x-3)^2$ ($x\to+\infty$).\label{asymEqu}

\begin{prop}[$o$ and $\sim$ are asymptotic]\label{prop_onoChange}
Let $A\in\R^*$. \underline{Then} the
two relations $\mathcal{A},\mathcal{B}\sus
\mathcal{R}^2$, defined by
$$
f\,\mathcal{A}\,g\iff f=o(g)\ \ (x\to A)
\,\text{ and }\,
f\,\mathcal{B}\,g\iff f\sim g\ \ (x\to A)\,,
$$
are asymptotic.
\end{prop}
\duk
Let $A\in\R^*$ and $f$, $g$, $f_0$ and $g_0$ be 
functions in $\mathcal{R}$ such that $f\doteq f_0$, $g\doteq g_0$. Let $f=o(g)$ ($x\to A$). We show that then $f_0=o(g_0)$ $(x\to A$) as well. For
the relation $\mathcal{B}$ we argue similarly. By the 
assumption, $A$ is in $L(f/g)$ and $\lim_{x\to A}\frac{f(x)}{g(x)}=0$.
By Exercise~\ref{ex_Xfini} the set $\frac{f_0}{g_0}[\R]\Delta
\frac{f}{g}[\R]$ is finite. Hence $A\in L(f_0/g_0)$ and $\lim_{x\to A}\frac{f_0(x)}{g_0(x)}=
\lim_{x\to A}\frac{f(x)}{g(x)}=0$.
\kduk
\vspace{-3mm}
\begin{exer}\label{ex_parOtazek}
Answer the following questions.
\begin{enumerate}
\item Is $x^2=o(x^3)$ ($x\to+\infty$)?
\item Is $x^3=o(x^2)$ ($x\to 0$)?
\item Is $x^2=o(x^3)$ ($x\to 0$)?
\item Is $(x+1)^3\sim x^3$ ($x\to 1$)?
\item Is $(x+1)^3\sim x^3$ ($x\to+\infty$)?
\item Is $\mathrm{e}^{-1/x^2}=o(x^{20})$ ($x\to0$)?
\end{enumerate}
\end{exer}

\noindent
{\em $\bullet$ Properties and mutual relations of asymptotic symbols. }We cannot devote to this important
topic much space. We restrict to one proposition and ten exercises.  

\begin{prop}[$o$ implies $O$]\label{prop_oJeO}
Let $f,g\in\mathcal{R}$. If $f(x)=o(g(x))$ ($x\to A$), \underline{then} there is a~$\theta$ such that 
$$
\text{$f=O(g)$ (on $P(A,\theta)$)}\,.
$$
\end{prop}
\duk
Let $f$, $g$ and $A$ be as stated. Since $\lim_{x\to A}
\frac{f(x)}{g(x)}=0$, for $\ep=1$ there is a~$\theta$ such 
that for every $x\in M(f/g)\cap P(A,\theta)$,
$$
{\textstyle
\big|\frac{f(x)}{g(x)}\big|=\big|\frac{f(x)}{g(x)}-0\big|<1\,.
}
$$
Hence $f=O(g)$ (on $P(A,\theta)$).
\kduk

In the ten exercises, $f$, $g$ and $h$ are in $\mathcal{R}$,
$N\sus\R$ and $A\in\R^*$.

\begin{exer}\label{ex_Oslozo}
If $g=o(h)$ ($x\to A$), $f(x)=O(x)$ (on $N$) and 
$$
A\in L\big(M((f\,|\,N)(g))\big)\,, 
$$
then $f(g)=o(h)$ ($x\to A$).    
\end{exer}

\begin{exer}\label{ex_ExOne}
If $f=O(h)$ (on $N$) and $g=O(h)$ (on $N$), then $f+g=O(h)$ (on $N$).    
\end{exer}

\begin{exer}\label{ex_ExTwo}
If $f=O(h)$ (on $N$) and $g$ is bounded on $N$, then $fg=O(h)$ (on $N$).
\end{exer}

\begin{exer}\label{ex_ExThree}
If $f=O(h)$ (on $N$) and $\frac{1}{g}$ is bounded on $N$, then $f/g=O(h)$ (on $N$).
\end{exer}

\begin{exer}\label{ex_ExFour}
If $f=o(h)$ ($x\to A$), $g=o(h)$ ($x\to A$)  and $A\in L(\frac{f+g}{h})$, then $f+g=o(h)$ ($x\to A$).
\end{exer}

\begin{exer}\label{ex_ExFive}
If $f=o(h)$ ($x\to A$), $g$ is bounded on
$P(A,\theta)$ for some $\theta$ and $A\in L(\frac{fg}{h})$, then $fg=o(h)$ ($x\to A$).
\end{exer}

\begin{exer}\label{ex_ExSix}
If $f=o(h)$ ($x\to A$), $\frac{1}{g}$ is bounded on $P(A,\theta)$ for some $\theta$ and $A\in L(\frac{f}{gh})$, then $f/g=o(h)$ ($x\to A$).
\end{exer}

\begin{exer}\label{ex_ExSeven}
If $f(x)\sim h(x)$ ($x\to A$), $g(x)=o(h(x))$ ($x\to A$) and $A$ is in $L(\frac{f+g}{h})$, then $f(x)+g(x)\sim h(x)$ ($x\to A$).
\end{exer}

\begin{exer}\label{ex_ExEight}
If $f(x)\sim h(x)$ ($x\to A$), $\lim_{x\to A}g(x)=1$ and $A\in L(\frac{fg}{h})$ then
$f(x)g(x)\sim h(x)$ ($x\to A$).
\end{exer}

\begin{exer}\label{ex_ExNine}
If $f(x)\sim h(x)$ ($x\to A$), $\lim_{x\to A}g(x)=1$ and $A\in L(\frac{f}{gh})$ then
$f(x)/g(x)\sim h(x)$ ($x\to A$).
\end{exer}

Notation $o$, $O$ and $\sim$ was introduced by the German mathematicians {\em 
Paul Bachmann\index{Bachmann, Paul} (1837--1920)} and {\em Edmund 
Landau\index{Landau, Edmund} (1877--1938)}. Notation $\ll$,
$\gg$ and $\asymp$ is due to the Russian mathematician {\em 
Ivan M. Vinogradov\index{Vinogradov, Ivan M.} (1891--1983)}. 

\medskip\noindent
{\em $\bullet$ Famous asymptotics. }The \underline{prime number counting function\index{pi of x@$\pi(x)$|emph}} $\pi\cc\R\to\N_0$ is defined by
$$
\pi(x)\equiv|(-\infty,\,x]\cap\mathbb{P}|\label{piOfx}
$$
where $\mathbb{P}=\{2,3,5,7,11,\ds\}$ is the set of prime numbers. 
In 1896 the French mathematician {\em Jacques Hadamard\index{Hadamard, 
Jacques} (1865--1963)} and, in parallel with him, the Belgian 
mathematician {\em Charles Jean de la Vall\'ee Poussin\index{vallee@de la 
Vall\'ee Poussin, Charles J.} (1866--1962)} proved the famous 
\underline{Prime Number Theorem\index{theorem!PNT}} (PNT):  
$$
\pi(x)\sim\frac{x}{\log x}\ \ (x\to+\infty)\,.
$$

In Section~\ref{podkap_ctyriRobustni}
we introduced for
$k,n\in\N$ the number $r_k(n)\in\N_0$
as the size of the largest set $X\sus[n]$ containing no arithmetic
progression with length $k$. In 1975 E.~Szemer\'edi\index{Szemeredi@Szemerédi, Endre} proved the now famous
theorem, which we mentioned in Section~\ref{podkap_ctyriRobustni},
that for every $k$, 
$$
r_k(n)=o(n)\ \  (n\to+\infty)\,. 
$$

For $x\in\R$ we define
$$
D(x)\equiv|\{(m,\,n)\in\N^2\cc\;mn\le x\}|\,.\label{dOfx}
$$
For $n\in\N$ we denote by $\tau(n)$\index{tau n@$\tau(n)$|emph} the number of divisors of 
$n$. For example, 
$$
\tau(28)=|\{1,\,2,\,4,\,7,\,14,\,28\}|=6\,.
$$

\begin{exer}\label{ex_delitele}
Show that $D(x)=\sum_{n=1}^{\lfloor x\rfloor}\tau(n)$.    
\end{exer} 

\noindent
The (Dirichlet) \underline{divisor problem\index{divisor problem|emph}} 
is the problem to estimate the error in asymptotics of $D(x)$. In 1849 the 
German mathematician {\em Peter L. Dirichlet\index{Dirichlet, Peter L.} 
(1805--1859)} proved that
$$
D(x)=x\log x+(2\gamma-1)x+O(\sqrt{x})\ \ \text{(on $[1,+\infty)$)}\,,
$$
where $\ga$ is Euler's constant.
In 1903 the Russian-Ukrainian  mathematician {\em Georgij F. Voronoj\index{Voronoj, Georgij F.} (1868--1908)} improved it to
$$
D(x)=x\log x+(2\gamma-1)x+O\big(x^{1/3}\log x\big)\ \ \text{(on $[2,\,+\infty)$)}\,.
$$

\begin{exer}\label{ex_whyNotOn}
Why not on $[1,+\infty)$ as before?  
\end{exer}
The 20th century saw a~series of further advances in the divisor
problem. The current record holder is the British mathematician {\em Martin 
N. Huxley\index{Huxley, Martin N.} (1944)} who proved in 2003 that for 
every $\ep$,  
$$
D(x)=x\log x+(2\gamma-1)x+O_{\ep}\big(x^{131/416+\ep}\big)\ \ \text{(on $[1,\,+\infty)$)}\,.
$$

For $n\in\N$ and an algorithm (Turing machine)\index{Turing machine} $T$ for multiplying integers we 
define $T(n)$ as the smallest $k\in\N$ such that $T$ multiplies any two $n$-digit numbers in at most $k$ steps. The elementary school 
algorithm $T_{\mathrm{es}}$ works in 
$$
\text{$T_{\mathrm{es}}(n)=O(n^2)$ (on $\N$)}\label{Tes}
$$
steps. In 1960 the Russian mathematician {\em Anatolij A. 
Karacuba\index{Karacuba, Anatolij A.} (1937--2008)} invented an algorithm 
$T_{\mathrm{K}}$ working in
$$
\text{$T_{\mathrm{K}}(n)=O\big(n^{\log_2 
3}\big)=O\big(n^{1.585\ds}\big)$  (on $\N$)}\label{Tk}
$$
steps.
In 2021 computer scientists {\em David Harvey\index{Harvey, David}} from Australia and {\em Joris van der Hoeven\index{Hoeven@van der Hoeven, Joris} (1971)} from the Netherlands discovered an algorithm $T_{\mathrm{HH}}$ for multiplying integers with complexity
$$
T_{\mathrm{HH}}(n)=O(n\log n)\ \  \text{(on $\N$)}\,.\label{Thh}
$$

\begin{exer}\label{ex_procNe1}
But for $n=1$ the ratio $\frac{T_{\mathrm{HH}}(n)}{n\log n}$ 
is not defined?    
\end{exer}

\section{Asymptotic 
expansions}\label{sec_asyExp}

Asymptotic expansions provide infinite sequences of better and better approximations of the given function.

\medskip\noindent
{\em $\bullet$ Asymptotic expansions. }We define asymptotic scales.

\begin{defi}[asymptotic scales]\label{def_asySca}
A~sequence of functions $(f_n)\sus\mathcal{R}$ is an \underline{asymptotic scale\index{asymptotic scale|emph}} for $x\to A$ if
$A$ is in $L(\bigcap_{n=1}^{\infty}M(f_n))$, there 
is a~$\theta$ such that $f_n\ne0$ on $P(A,\theta)$ for every $n$, and for 
every $n$ we have 
$$
\text{$f_{n+1}(x)=o(f_n(x))$ ($x\to A$)}\,.
$$
\end{defi} 
For example, $(x^{-n})$ is an asymptotic scale for $x\to+\infty$, 
and $(x^n)$ for $x\to0$.

\begin{defi}[asymptotic expansions]\label{def_asyExp}
Suppose that $(a_n)\sus\R$, $f\in\mathcal{R}$ and that $(f_n)\sus\mathcal{R}$ is an asymptotic scale for $x\to A$. If for every $n$ we have
$$
{\textstyle
f(x)=\sum_{i=1}^n a_if_i(x)+o(f_{n+1}(x))\ \ (x\to A)\,,
}
$$
we say that $(a_n f_n(x))$ is an \underline{asymptotic expansion\index{asymptotic 
expansion (AE), $\approx$|emph}} of $f(x)$ for $x\to A$ and write it symbolically as
$$
{\textstyle
f(x)\approx\sum_{n=1}^{\infty} a_nf_n(x)\ \ (x\to A)}\,.\label{approx}
$$
\end{defi}
From this definition it follows that
for every $n$ there is a~$\theta_n>0$ such that $f=\sum_{i=1}^n a_if_i+O(f_{n+1})$ (on $P(A,\theta_n)$). 
Thus we have a~sequence of asymptotic approximations
$\sum_{i=1}^n a_if_i$ to $f$ with errors of order $f_{n+1}$.

\begin{exer}\label{ex_onAsyExp}
Prove it.    
\end{exer}
The assumption that $(f_n)$ is an asymptotic scale ensures that as $n$ grows, 
magnitudes of these errors get smaller and smaller. For fixed $x\in\R$, nothing is 
assumed about the convergence of the series $\sum a_nf_n(x)$; it typically diverges. Usually it is {\em not} true that $f(x)=\sum_{n=1}^{\infty}a_nf_n(x)$.

\medskip\noindent
{\em $\bullet$ Three asymptotic expansion. }We 
conclude this chapter  with three examples of asymptotic expansions. Their 
proofs will be given in {\em MA~1${}^+$}. In 1730, the Scottish mathematician {\em 
James Stirling\index{Stirling, 
James} (1692--1770)} obtained an asymptotic expansion of $\log(n!)$ 
for $n\to+\infty$. We state his expansion in a~moment. Modern theory of asymptotic expansions starts with the French
mathematician {\em Henri Poincar\'e\index{poincare@Poincar\'e, 
Henri} (1854--1912)} and his memoir of 1886. For expositions of asymptotic expansions, see \cite{cops,erde}. 

\begin{thm}[asymptotic expansion of 
$\log(n!)$]\label{thm_StirExp}
For\index{theorem!asymptotic expansion of log fac@asymptotic expansion of $\log(n"!)$} $n\in\N$,
$$
{\textstyle
\log(n!)\approx(n+\frac{1}{2})\log n-n+
\frac{1}{2}
\log(2\pi)+\sum_{k=1}^{\infty}
\frac{B_{2k}}{2k(2k-1)}\cdot 
n^{1-2k}\ \ (n\to+\infty)\,.
}
$$
\end{thm}

\begin{thm}[harmonic numbers]\label{thm_HarmExp}
For\index{theorem!asymptotic expansion of harmonic numbers}\index{harmonic number, $h_n$} 
$n\in\N$,
$$
{\textstyle
h_n=\sum_{i=1}^n\frac{1}{i}\approx 
\log n+\gamma+
\frac{1}{2n}-
\sum_{k=1}^{\infty}\frac{B_{2k}}{2k}\cdot n^{-2k}\ \ (n\to+\infty)\,.
}
$$
\end{thm}
In these two expansions, $B_k$ ($\in\Q$) for $k\in\N_0$ denote the \underline{Bernoulli numbers\index{Bernoulli number, $B_k$|emph}}.\label{bernou}
They are defined by the power series expansion
$$
{\textstyle
\frac{x}{\exp x-1}=\sum_{k=0}^{\infty}\frac{B_k}{k!}\cdot x^k
}
$$
and are named after their discoverer, the Swiss mathematician 
{\em Jacob Bernoulli\index{Bernoulli, Jacob} (1655/54--1705)}. They satisfy $B_{2k+1}=0$ for every $k\in\N$ and have initial values $B_0=1$, $B_1=-\frac{1}{2}$, $B_2=\frac{1}{6}$, $B_4=-\frac{1}{30}$, $B_6=\frac{1}{42}$, $B_8=-\frac{1}{30}$, 
$B_{10}=\frac{5}{66}$, $B_{12}=-\frac{691}{2730}$, 
$B_{14}=\frac{7}{6}$ and 
$B_{16}=-\frac{3617}{510}$.

\begin{exer}\label{ex_BernNum}
Deduce from the displayed definitoric expansion a~recurrence for $B_k$ and show that $B_{2k+1}=0$ for every $k\in\N$.    
\end{exer}
In \cite{cops} we read that the formula in Theorem~\ref{thm_StirExp} is not exactly the original expansion of 
Stirling, but a~very similar formula due to the French mathematician {\em Abraham de 
Moivre\index{moivre@de Moivre, Abraham} (1667--1754)}. The 
asymptotic expansion of harmonic numbers is due to L.~Euler.\index{Euler, Leonhard}

The third asymptotic expansion is much more recent, due to 
\cite{mont_nurl}. We say that a~graph $G=(V,E)$ is 
\underline{connected\index{graphs!connected|emph}} if for every partition 
$\{A,B\}$ of the vertex set $V$ there is an edge $e\in E$ that 
intersects both blocks $A$ and $B$.

\begin{exer}\label{ex_equiConnec}
Equivalently, a~graph $G$ is connected iff for every two vertices $u,v\in V$ 
there exists 
a~$(k+1)$-tuple of vertices 
$$
w=\langle u_0,\,u_1,\,\ds,\,u_k\rangle,\ \  k\in\N_0\,,
$$
such that $u_0=u$, $u_k=v$ and $\{u_{i-1},u_i\}\in E$ for every
$i\in[k]$. We say that $w$ is a~walk in $G$ joining $u$ and $v$.
\end{exer}

\begin{thm}[probability of connectedness]\label{thm_CommExp}
For\index{theorem!asymptotic expansion of the probability of connectedness} $n\in\N$,
\begin{eqnarray*}
&&{\textstyle
\frac{1}{2^{n(n-1)/2}}}\cdot
|\{G=([n],\,E)\cc\;\text{$G$ is a~connected graph}\}|\\
&&\approx{\textstyle 1-\sum_{k=1}^{\infty}
t_k\cdot 2^{k(k+1)/2}\cdot\binom{n}{k}2^{-kn}\ \ (n\to+\infty)\,,
}
\end{eqnarray*}
where $t_k$ is the number of irreducible tournaments of size $k$ (see below).
\end{thm}
As mentioned, this result is due to \cite{mont_nurl}. Since there are $2^{n(n-1)/2}$ graphs $G=([n],E)$, the product on the left-hand
side equals to the probability that a~random graph with the vertex set
$[n]$ is connected. By the first term of the expansion, for
$n\to+\infty$ this probability goes to~$1$. What is $t_k$? 
A~\underline{tournament\index{tournament|emph}}\label{tourn} $T=(V,E)$ is a~pair of 
a~nonempty finite set $V$ of vertices and an irreflexive relation $E\sus 
V\times V$ such that for every two distinct vertices $u,v\in V$ 
exactly one pair of $(u,v)$ and $(v,u)$ is in $E$. We say that $T$ is 
\underline{irreducible\index{tournament!irreducible|emph}} if for every partition $\{A,B\}$ of the vertices $V$ there exist pairs $(a,b)\in A\times B$ and $(c,d)\in B\times A$
such that $(a,b),(c,d)\in E$. Then $t_k$\index{tk@$t_k$|emph}\label{tk} is the number of 
irreducible tournaments $T=([k],E)$.
By \cite[A054946]{oeis} the initial values of $t_k$ are
$$
(t_1,\,t_2,\,\ds)=(1,\,0,\,2,\,24,\,544,\,22320,\,1677488,\,\ds)\,.
$$

See \cite{bori_phd,bori} for a~theory of a~calculus by which one can derive asymptotic 
expansions in certain problems of enumerative 
combinatorics\index{enumerative combinatorics} and theoretical physics.
\index{asymptotic notation|)}

\chapter[Elementary functions]{Elementary functions}\label{chap_EF}

\section[Basic elementary functions]{Basic 
elementary functions}\label{sec_elemenFce}

In this and next two sections
we\underline{\index{basic 
elementary functions, BEF|(}} 
introduce five families of functions 
$f\cc M\to\R$ with $M\sus\R$. Namely, basic elementary functions (BEF), 
elementary functions (EF), really
basic elementary functions (RBEF), polynomials (POL) and rational functions (RAC). 

\begin{defi}[BEF]\label{def_ZEF}
\underline{Basic 
elementary functions\index{basic 
elementary functions, BEF|emph}}, or {\em BEF},\label{BEF} are the following functions. 
\begin{enumerate}
\item The constant functions $k_c(x)$ for $c\in\R$.
\item The functions $\exp x$ and $\log x$.
\item The functions $a^x$ for $a>0$, $x^b$ for
$b\in\R$, $0^x$ and $x^m$ for $m\in\Z$.
\item The functions $\sin x$, $\cos x$, $\tan x$ and 
$\cot x$.
\item The functions $\arcsin x$, $\arccos x$, 
$\arctan x$ and $\mathrm{arccot}\,x$.
\end{enumerate} 
\end{defi}
We review them shortly. The functions $x^b$ and $x^m$ differ in definition domains.

\medskip\noindent
{\em $\bullet$ \underline{Constant functions\index{constant functions, $k_c(x)$|emph}}}, or constants, are functions ($c\in\R$)
$$
k_c\cc\R\to\R,\ \ k_c(x)=c\,.\label{kcx} 
$$
Instead of
$k_c(x)$ we usually write $k_c$ or just $c$.

\begin{exer}\label{ex_priklNaKon}
How many constant functions $k_c(x)$ are there?    
\end{exer}

\noindent
{\em $\bullet$  The \underline{exponential\index{exponential function, $\exp x$|emph} function} }
$$
\exp x=\exp(x)=\mathrm{e}^x\cc\R\to\R
$$ 
is for $x\in\R$ defined by the sum 
$${\textstyle
\exp x\equiv\sum_{n=0}^{\infty}\frac{1}{n!}x^n=1+x+\frac{1}{2}x^2+\frac{1}{6}x^3+\ds\,,
}\label{exp}
$$ 
where $x^0=0^0\equiv1$.

\begin{exer}\label{ex_radaJeAK} For every $x\in\R$ the series $\sum_{n=0}^{\infty}\frac{1}{n!}x^n$ is abscon.
\end{exer}
We prove the exponential identity.

\begin{thm}[exponential identity]\label{thm_expFce}
For\index{theorem!exponential identity|emph} 
every\underline{\index{exponential function, 
$\exp x$!exponential identity|emph}} 
$x,y\in\R$ we have
$$
\exp(x+y)=\exp(x)\cdot\exp(y)\,.
$$
\end{thm}
\duk
Let $x,y\in\R$. By 
Theorem~\ref{thm_cauSouRad} and Exercises~\ref{ex_radaJeAK} and \ref{ex_binomVeta}, 
the product $\mathrm{e}^x\cdot\mathrm{e}^y$ equals
$$
{\textstyle
\sum_{n=0}^{\infty}\sum_{k=0}^n\frac{x^k}{k!}\cdot\frac{y^{n-k}}{(n-k)!}=
\sum_{n=0}^{\infty}\frac{1}{n!}
\sum_{k=0}^n\binom{n}{k}x^ky^{n-k}=
\sum_{n=0}^{\infty}\frac{(x+y)^n}{n!}
}
$$
which is $\mathrm{e}^{x+y}$.
\kduk

\noindent
Nothing changes in the complex domain. Here are some more properties of 
$\exp x$.

\begin{exer}\label{ex_dok1az3}
Prove parts 1--3 of the following proposition.
\end{exer}

\begin{prop}[properties of $\mathrm{e}^x$]\label{prop_expFce}
The\underline{\index{exponential function, 
$\exp x$!properties of|emph}} 
following holds.
\begin{enumerate}
    \item For every $x\in\R$, we have $\exp x>0$ and $\exp(-x)=\frac{1}{\exp x}$. Also, $\exp 0=1$. 
    \item For all real $x<y$ we have $\exp x<\exp y$.
    \item $\lim_{x\to-\infty}\exp x=0$ and $\lim_{x\to+\infty}\exp x =+\infty$.
    \item The function $\exp$ is a~bijection from $\R$ to $(0,+\infty)$.
\end{enumerate}
\end{prop}
We prove part~4 later in Corollary~\ref{cor_imageExp}.

\medskip\noindent
{\em $\bullet$ \underline{Euler's number\index{Euler's 
number, $\mathrm{e}=2.71\ds$|emph}} }is the sum 
$${\textstyle
\mathrm{e}\equiv\exp 
1=\sum_{n=0}^{\infty}\frac{1}{n!}=2+\frac{1}{2!}+\frac{1}{3!}+\cdots=2.71828\cdots\,.\label{e}
}
$$

\begin{exer}\label{exeJeIrac}
Show\underline{\index{Euler's 
number, $\mathrm{e}=2.71\ds$!irrationality of}} 
that $\mathrm{e}$ is irrational. Hint: multiply the equality 
$\sum_{j=0}^{\infty}\frac{1}{j!}=\frac{n}{m}$ by $m!$.    
\end{exer}

\noindent
{\em $\bullet$ (Natural) \underline{logarithm\index{logarithm $\log x$|emph}} }is the inverse of the exponential function,
$$
\log\cc(0,\,+\infty)\to\R,\ \ \log x\equiv\exp^{-1}(x)\,.\label{log}
$$ 
Properties of $\log x$ follow from the properties of
$\mathrm{e}^x$.

\begin{exer}\label{ex_vlLogaritmu}
Prove the following proposition.
\end{exer}

\begin{prop}[properties of $\log x$]\label{prop_logar}
The\underline{\index{logarithm, $\log x$!properties of|emph}} 
following holds.
\begin{enumerate}
\item For all real $x,y>0$ we have 
$$
\log(xy)=\log x+\log y\,. 
$$
If $x<y$ then $\log x<\log y$. Also, $\log 1=0$.
\item $\lim_{x\to0}\log x=-\infty$ and $\lim_{x\to+\infty}\log x=+\infty$.
\item Logarithm is a~bijection from $(0,+\infty)$ to $\R$.
\end{enumerate}
\end{prop}

\noindent 
{\em $\bullet$ Real exponentiation. }In the expression $a^b$, we call $a$ the \underline{base\index{real 
exponentiation $a^b$!base|emph}} 
and $b$ the \underline{exponent\index{real 
exponentiation $a^b$!exponent|emph}}. We consider two families of real exponentials. One 
expresses $a^b$ by $\exp(b\log a)$ and extends it by the limit $\lim_{x\to-\infty}\exp x=0$. The other makes use of
iterated multiplication.

\begin{defi}[$a^b$ analytically]\label{def_aNaB}
We define\underline{\index{real exponentiation 
$a^b$!analytically|emph}} the following functions.
\begin{enumerate}
\item Let $a>0$. We set 
$$
a^x\equiv\exp(x\log a)\ \ (\in\mathcal{F}(\R))\,.\label{anax}
$$
\item Let $b>0$. We set $0^b\equiv0$ and
$$
\text{$x^b\equiv\exp(b\log x)$ for $x>0$}\,. 
$$
Thus
$x^b\in\mathcal{F}([0,+\infty))$. 
\item Let $b\le0$. We set
$$
x^b\equiv\exp(b\log x)\ \ 
(\in\mathcal{F}((0,\,+\infty)))\,. \label{xnab}
$$
\item We define 
$$
0^x\equiv k_0(x)\,|\,(0,+\infty)\ \ (\in\mathcal{F}((0,\,+\infty))\,\label{nulanax}
$$
\end{enumerate}  
\end{defi}
We do not define $0^0$ and always $a^b\ge0$. Odd roots\index{odd root} like
$\sqrt[3]{x}\equiv x^{1/3}$, $\sqrt[5]{x}\equiv x^{1/5}$, etc. are sometimes
defined for every $x\in\R$, for example, $\sqrt[3]{-8}=-2$. We
define them only for $x>0$. Note that in
parts~2 and~3 definition domains differ.

\begin{defi}[square roots]\label{def_sqRoot}
We define the \underline{(square) root\index{square root|emph}} by
$$
\sqrt{x}\equiv x^{1/2}\label{sqrt}\ \ (\in\mathcal{F}([0,+\infty)))\,.
$$
\end{defi}

\begin{exer}\label{ex_takToNej}
So $\sqrt{x}=\exp(\frac{1}{2}\log x)$, right? 
\end{exer}

\begin{exer}\label{ex_onSqRo}
Let $a\ge0$. Then $\sqrt{a}$ is the unique number $b\ge0$ such that $b^2=a$. 
\end{exer}

We proceed to the other family of real exponentials.

\begin{defi}[$a^b$ algebraically]\label{def_anaM}
Let $m\in\Z$, where now $\Z$ is understood as disjoint from $\R$.
We define\underline{\index{real exponentiation $a^b$!algebraically|emph}} 
the following functions.

\begin{enumerate}
\item If $m>0$ then 
$$
x^m\equiv\underbrace{x\cdot x\cdot\ldots\cdot x}_{\text{$m$ factors}}\ \ (\in\mathcal{F}(\R))\,.
$$  
\item We set $x^0\equiv k_1(x)$ ($\in\mathcal{F}(\R)$).\label{xnanula}
\item If $m<0$ then 
$${\textstyle
x^m\equiv\frac{k_1(x)}{x^{-m}}\ \ (\in\mathcal{F}(\R\setminus\{0\}))\,.\label{xnam}
}
$$ 
\end{enumerate}
\end{defi}
So now $0^0\equiv1$ and $a^b$ may be negative. For $m=1$ we
get the identity function $x^1=x$ ($\in\mathcal{F}(\R)$). This function also arises as the composition $\log(\mathrm{e}^x)$.

\begin{exer}\label{ex_jsouKomp}
Definitions~\ref{def_aNaB} and \ref{def_anaM} coincide on the intersection.
\end{exer}

\begin{exer}\label{ex_onE}
Show that 
$$
\mathrm{e}^x=\exp x
$$ 
for every $x\in\R$. On the left, we have the real exponentiation with the base
$\mathrm{e}=2.71\ds$ and exponent
$x$. On the right, we have a~value of the exponential function. 
\end{exer}

\begin{exer}[$(+\infty)^0$]\label{ex_laterEx}
Show that for every sequence $(a_n)\sus(0,+\infty)$ with $\lim a_n=+\infty$ and for every nonnegative
$A\in\R^*$ there exists a~sequence $(b_n)$ 
such that $\lim b_n=0$ and
$$
\lim_{n\to\infty}\big(a_n\big)^{b_n}=A\,.
$$
\end{exer}

\noindent
{\em $\bullet$ Exponential identities. }We discuss some well-known and some not so well-known identities for real
exponentiation.

\begin{thm}[three exponential identities]\label{thm_realMocn}
Let\index{theorem!three exponential identities|emph} 
$a,b>0$. For every real $x$ and $y$ the following holds.
\begin{enumerate}
\item $(a\cdot b)^x=a^x\cdot b^x$.
\item $a^x\cdot a^y=a^{x+y}$.
\item $\big(a^x\big)^y=a^{x\cdot y}$.
\end{enumerate}
\end{thm}
\duk
1. Indeed, $(ab)^x$ equals
$$
\exp(x\log(ab))=\exp(x\log a+x\log b)=\exp(x\log 
a)\cdot\exp(x\log b)=a^xb^x\,.
$$

2. Indeed,  $a^xa^y$ equals
$$
\exp(x\log a)\exp(y\log a)=\exp(x\log a+y\log a)= 
\exp((x+y)\log a)=a^{x+y}\,.
$$

3. Indeed,
$$
(a^x)^y=\exp(y\log(\exp(x\log a)))=\exp(yx\log a)=a^{xy}\,.
$$
\kduk

\noindent
If we drop the assumption that $a>0$, the third 
identity ceases to hold:
$$
\left((-1)^2\right)^{\frac{1}{2}}=1^{\frac{1}{2}}=1\ne-1=(-1)^1
=(-1)^{2\cdot\frac{1}{2}}\,,
$$
where we used both Definitions~\ref{def_aNaB} and 
\ref{def_anaM} of real 
exponentiation.

\begin{exer}\label{ex_howExa}
How exactly did we use them?    
\end{exer}

A.~Tarski\index{Tarski, Alfred} conjectured that every identity for real
exponentiation, for example
$$
x^y\cdot (x^y)^y=x^{y+y^2}\,, 
$$
can be derived from the three previous identities and from other basic properties of addition, multiplication, and exponentiation. 
In 1981 the British mathematician  {\em Alex Wilkie\index{Wilkie, Alex} 
(1948)} refuted Tarski's conjecture and showed that exponential identities like 
\begin{eqnarray*}
&&\big((1+x)^y+(1+x+x^2)^y\big)^x\cdot 
\big((1+x^3)^x+(1+x^2+x^4)^x\big)^y\\ 
&&=\big((1+x)^x+(1+x+x^2)^x\big)^y\cdot 
\big((1+x^3)^y+(1+x^2+x^4)^y\big)^x
\end{eqnarray*}
cannot be derived from the three basic identities. Identities of this kind are now called 
\underline{Wilkie's identities\index{Wilkie's 
identities|emph}}. We return to them in {\em MA~1${}^+$}.

\begin{exer}\label{ex_Wilkie}
Prove that for every real $x,y>0$ the stated Wilkie's identity holds. Hint: $(1+x)\cdot(1+x^2+x^4)=(1+x^3)\cdot(1+x+x^2)$.  
\end{exer}
The next exercise shows that $0^0$ is an indefinite expression. 
\begin{exer}\label{ex_0na0jeNeur}
Let $A\in\R^*$ with $A\ge0$. Then there exist sequences $(a_n)\sus(0,+\infty)$ and $(b_n)\sus\R$ such that 
$\lim a_n=\lim b_n=0$ and 
$$
\lim(a_n)^{b_n}=
A\,. 
$$
Could $A$ be negative?   
\end{exer}
However, the definition $0^0\equiv1$ is often useful.

\medskip\noindent
{\em $\bullet$ Cosine\index{cosine, $\cos x$|emph} and sine\index{sine, 
$\sin x$|emph}. }We define
cosine and sine for $t\in\R$ by the sums 
$${\textstyle
\cos t\equiv\sum_{n=0}^{\infty}
(-1)^n\frac{1}{(2n)!}t^{2n}\ \  (0^0=1)\,\text{ and }\,
\sin t\equiv\sum_{n=0}^{\infty}(-1)^n\frac{1}{(2n+1)!}t^{2n+1}\,.
}
$$
In other words,\label{cos}\label{sin} 
$$
{\textstyle
\cos t=1-
\frac{t^2}{2}+\frac{t^4}{24}-\cdots\,\text{ and }\,\sin t=t-\frac{t^3}{6}+\frac{t^5}{120}-\cdots\ \ (\in\mathcal{F}(\R))\,.
}
$$

\begin{exer}\label{ex_kosSinAK}
For every $t\in\R$, the series defining $\cos t$ and $\sin t$ are abscon.   
\end{exer}
The planar set
$$
S\equiv\{(x,\,y)\in\R^2\cc\;
x^2+y^2=1\}
$$ 
is the \underline{unit circle\index{unit circle|emph}}. It has radius 
$1$ and center $(0,0)$. The next theorem, whose proof we postpone to 
{\em MA~1${}^+$}, decribes the main geometric property of cosine 
and sine.

\begin{thm}[runner's]\label{thm_bezkyne}
Let\index{theorem!runner's}
$t\in\R$. A~runner starts at the point $(1,0)$ of the track $S$ and 
runs on $S$ with unit
speed. For $t>0$ she runs counter-clockwise, and for $t\le0$ 
clockwise. \underline{Then} 
$$
\text{in time $|t|$ the runner is at the point $(\cos t,\,\sin t)\in S$}\,.
$$
\end{thm}

\noindent
{\em $\bullet$ The number $\pi$.\index{number pi@number 
$\pi$|emph} }One can define it in several ways. 
\begin{defi}[the number $\pi$]\label{def_pi}
The number 
$$
\pi=3.14159\ds\label{pi}
$$ 
is twice the minimum number $\al>0$ such that $\cos\al=0$.    
\end{defi}
The existence of $\al$ follows from the continuity of cosine (see Corollary~\ref{cor_spoExpKosSin}), 
from Theorem~\ref{thm_mezihodnoty} and from the values $\cos 0=1$ and 
$\cos 2<0$ (Exercise~\ref{ex_applLeibSer}).

Or one can define $\pi$ as half of the circumference of 
$S$, which equals the time when the runner passes for the second time 
through the point $(1,0)$. This is informal because the length of a~circular arc will be introduced only in~{\em MA~1${}^+$}. 

\begin{exer}\label{ex_applLeibSer}
Using Theorem~\ref{thm_altSer} show that $\cos 1>0$ and $\cos 2<0$.
\end{exer}

\noindent
{\em $\bullet$ More on cosine and sine. }Here are some properties of these functions. 

\begin{exer}\label{ex_odvodZvety}
Deduce from Theorem~\ref{thm_bezkyne} the next proposition.
\end{exer}

\begin{prop}[on $\sin x$ and $\cos x$]\label{prop_sinKos}
The following holds.
\begin{enumerate}
\item Both functions are $2\pi$-periodic, for every $t\in\R$ we have 
$$
\cos(t+2\pi)=\cos t\,\text{ and }\,\sin(t+2\pi)=\sin t\,.
$$
\item Sine increases on the interval $[0,\frac{\pi}{2}]$ from $0$ to $1$. 
\item For every $t\in[0,\pi]$ we have 
$$
\sin t=\sin(\pi-t)
$$
and for every $t\in[0,2\pi]$ we have 
$$
\sin t=-\sin(2\pi-t)\,.
$$
\item For every $t\in\R$ we have 
$${\textstyle
\cos t=\sin(t+\frac{\pi}{2})\,\text{ and }\,\cos^2 t+\sin^2 t=1\,.
}
$$
\item For every $s,t\in\R$ we have
\begin{eqnarray*}
\sin(s\pm t)&=&\sin s\cdot\cos t\pm \cos s\cdot\sin t\text{ and}\\
\cos(s\pm t)&=&\cos s\cdot\cos t\mp \sin s\cdot\sin t\,.
\end{eqnarray*}
\end{enumerate}
\end{prop}

\noindent
{\em $\bullet$ Euler's formula. }We present an important complex identity relating the three functions $\exp x$, 
$\sin x$ and $\cos x$. We extend for it series to the complex
domain. For a~complex sequence 
$$
(z_n)\sus\C
$$ 
we define the \underline{limit\index{limit of 
a~sequence!for complex s.|emph}} $\lim z_n$ to be the unique number $z\in\C$, if it
exists, such that for every
$\ep$ there is an $n_0$ such that for every $n\ge n_0$ we have $|z_n-
z|\le\ep$. 

\begin{exer}\label{ex_limReIm}
Let $(z_n)=(a_n+b_ni)\sus\C$ with $\lim z_n=a+bi$. Then 
$$
a=\lim a_n\,\text{ and }\,
b=\lim b_n\,.
$$
\end{exer}
We extend the exponential function to $\C$ by using the same formula
$$
{\textstyle
\mathrm{e}^z=\exp z=\sum_{n=0}^{\infty}
\frac{z^n}{n!}}
\equiv\lim_{n\to\infty}{\textstyle
\sum_{j=0}^n\frac{z^j}{j!}\ \ (z\in\C)\,.
}
$$
\begin{exer}\label{ex_complExp}
This limit exists for every complex number $z$.    
\end{exer}
Since the form of the series is the same as for the real exponential, the function 
$$
\mathrm{e}^z\cc\C\to\C
$$
extends the real exponential. Recall that $i$ denotes the imaginary unit, $i^2=-1$.

\begin{thm}[Euler's\index{Euler's formula} formula]\label{thm_EulVzo}
For every $t\in\R$ we have
$$
\exp(it)=\cos t+i\sin t\,.
$$ 
\end{thm}
\duk
Let $t\in\R$. Since the natural powers of $i$ form a~$4$-periodic sequence 
$$
(i^n)=(i,\,-1,\,-i,\,1,\,i,\,-1,\,-i,\,1,\,i,\,\ds)\,, 
$$
we have for every $n\in\N_0$ that
$${\textstyle
\sum_{j=0}^n\frac{(it)^j}{j!}=\sum_{k=0}^{\lfloor n/2\rfloor}
(-1)^k\frac{1}{(2k)!}t^{2k}+
i\sum_{l=1}^{\lceil n/2\rceil}
(-1)^{l-1}\frac{1}{(2l-1)!}t^{2l-1}\,.
}
$$
We send $n\to\infty$ and get by Exercise~\ref{ex_limReIm} and by the definitions of $\exp(it)$, $\cos t$ and $\sin t$ Euler's formula.
\kduk

\noindent
{\em $\bullet$ Another beautiful exponential identity. }In \cite[Lemma 17.1]{lach} we 
found a~remarkable result on the function $\mathrm{e}^z$. In \cite{lach} 
it serves as a~lemma in proofs of limit theorems\index{central limit theorems} 
in probability theory. In the proof we follow roughly \cite{lach} and begin 
with a~lemma.

\begin{lemma}[telescoping]\label{lem_anoTeles}
Let $u_1$, $\ds$, $u_n$, $v_1$, $\ds$, $v_n$ be $2n$ formal 
variables. \underline{Then}
$$
u_1u_2\ds u_n-v_1v_2\ds v_n=
\sum_{i=1}^n V_{i-1}(u_i-v_i)U_{i+1}\,,
$$
where $V_0=1$, $V_j=v_1v_2\ds v_j$ for $j\in[n]$, $U_{n+1}=1$ and $U_j=u_ju_{j+1}\ds u_n$ for $j\in[n]$.
\end{lemma}

\begin{exer}\label{ex_proLem}
Prove the lemma.
\end{exer}

\begin{thm}[sums and products]\label{thm_lemCLT}
Let $a\in\C$, $(k_n)\sus\N$ and for $n\in\N$ let\index{theorem!sums and products|emph} 
$$
\langle a_{n,\,j}\cc\;j\in[k_n]\rangle\ \ (\sus\C)
$$
be $k_n$-tuples of complex numbers satisfying three conditions.
\begin{enumerate}
\item $\lim_{n\to\infty}\sum_{j=1}^{k_n}a_{n,j}=a$.
\item There is a~constant $c>0$ such that
$\sum_{j=1}^{k_n}|a_{n,j}|\le c$ for every $n\in\N$.
\item $\lim_{n\to\infty}\sum_{j=1}^{k_n}|a_{n,j}|^2=0$.
\end{enumerate}
\underline{Then}
$$
\lim_{n\to\infty}\prod_{j=1}^{k_n}\big(1+a_{n,\,j}\big)=\exp(a)\,.
$$
\end{thm}
\duk
Using condition~3 we take an $n_0\in\N$ such that for every $n\ge n_0$
and $i\in[k_n]$ we have $|a_{n,i}|\le\frac{1}{2}$. 
Using the exponential identity  and condition~2 we have for every $n\ge n_0$ and every set $A\sus[k_n]$ that
$$
{\textstyle
\big|\prod_{j\in A}(1+a_{n,j})\big|\le
\prod_{j\in A}(1+|a_{n,j}|)
\le\exp\big(\sum_{j\in A}
|a_{n,j}|\big)\le\mathrm{e}^c
}
$$
and, using Exercise~\ref{ex_absHexp} and condition~2, that
$${\textstyle
\big|\exp\big(\sum_{j\in A}a_{n,j}\big)
\big|=\exp\big(\mathrm{re}\big(\sum_{j\in A}a_{n,j}\big)\big)\le
\exp\big(\sum_{j\in A}|a_{n,j}|\big)
\le\mathrm{e}^c\,.
}
$$
Let $n\ge n_0$ and $i\in[k_n]$. In Lemma~\ref{lem_anoTeles} we set 
$u_i\equiv1+a_{n,i}$ and $v_i\equiv\exp(a_{n,i})$. We get 
with the help of Exercise~\ref{ex_boundExp}, the exponential identity, the above bounds and condition~3 that
\begin{eqnarray*}
&&{\textstyle\big|\prod_{j=1}^{k_n}(1+a_{n,\,j})-\exp\big(\sum_{j=1}^{k_n}a_{n,\,j}\big)\big|=\big|\sum_{i=1}^{k_n} V_{i-1}(u_i-v_i)U_{i+1}\big|\le  
}\\
&&\le{\textstyle
\sum_{i=1}^{k_n}|V_{i-1}|\cdot|u_i-v_i|\cdot|U_{i+1}|\le
\sum_{i=1}^{k_n}\mathrm{e}^c|a_{n,\,i}|^2\mathrm{e}^c\to0,\ \ n\to\infty\,.
}
\end{eqnarray*}
By condition~1 we get the stated formula.
\kduk
\vspace{-3mm}
\begin{exer}\label{ex_absHexp}
For every $z\in\C$ we have $|\exp z|=\exp(\mathrm{re}(z))$.    
\end{exer}

\begin{exer}\label{ex_boundExp}
For every $z\in\C$ with $|z|\le
\frac{1}{2}$ we have $|\exp z-1-z|\le|z|^2$. 
\end{exer}

\begin{exer}\label{ex_infProExp}
Let $a\in\C$. Deduce the infinite product
$${\textstyle
\prod_{n=1}^{\infty}\big(1+\frac{a}{n}\big)^n=\exp(a)\,.
}
$$
\end{exer}

\noindent
{\em $\bullet$ Tangent\index{tangent, $\tan x$|emph} and cotangent\index{cotangent, $\cot x$|emph}. }We define these functions by
$${\textstyle
\tan t\equiv\frac{\sin t}{\cos t}\,\text{ and }\,\cot t\equiv\frac{\cos t}{\sin t}\ \ (t\in\R)\,.\label{tan}\label{cot}
}
$$

\begin{exer}\label{ex_tanCot}
Their definition domains are
$${\textstyle
M(\tan)=
\R\setminus\{\frac{1}{2}(2m-
1)\pi\cc\;m\in\Z\}\,\text{ and }\,M(\cot)=\R\setminus
\{m\pi\cc\;m\in\Z\}\,.
}
$$
\end{exer}

\noindent
{\em $\bullet$ Arcsine\index{arcsine, $\arcsin 
x$|emph} (inverse sine)  and arccosine\index{arccosine, $\arccos 
x$|emph} (inverse cosine). }These functions
$$
\arcsin x\cc[-1,1]\to\R\,\text{ and }\,\arccos x\cc[-1,1]\to\R\label{arcsin}\label{arccos}
$$
are congruent to the inverses of the restriction $\sin x\,|\,[-\frac{\pi}{2},\frac{\pi}{2}]$ 
and $\cos x\,|\,[0,\pi]$. 

\begin{exer}\label{ex_arcsin}
Recall the meaning of ``congruence''.
\end{exer}

\noindent
{\em $\bullet$ Arctangent\index{arctangent, $\arctan x$|emph} (inverse tangent)  and arccotangent\index{arccotangent, 
$\mathrm{arccot}\,x$|emph} 
(inverse cotangent).} These functions
$$
\arctan\cc\R\to\R\,\text{ and }\,
\mathrm{arccot}\,x\cc\R\to\R\label{arctan}\label{arccot}
$$
are congruent to inverses of the restriction $\tan x\,|
\,(-\frac{\pi}{2},\frac{\pi}{2})$ and $\cot x\,|\,(0,\pi)$. 

\begin{exer}\label{ex_arctan}
Thus $\arctan x$ is congruent to a~bijection between which sets? 
\end{exer}

\section[Elementary functions]{Elementary functions}\label{sec_obecElemFce}

We\underline{\index{elementary functions, EF|(}} 
introduce a~family of elementary functions. We 
begin from operations on the set of functions $\mathcal{R}$.

\medskip\noindent
{\em $\bullet$ Five operations on $\mathcal{R}$. }Recall that 
$\mathcal{R}$ consists of functions $f\cc M\to\R$ with $M\sus\R$.

\begin{defi}[operations on 
$\mathcal{R}$]\label{def_oprOnR}
Let $f,g\in\mathcal{R}$.
We define the following operations.
\begin{enumerate}
\item  The
\underline{sum\index{function!r@$\mathcal{R}$!sum on, $+$|emph}} 
$$
f+g\cc M(f)\cap M(g)\to\R\,\text{ by }\,(f+g)(x)\equiv f(x)+g(x)\,.\label{sumFun}
$$
\item The \underline{product\index{function!r@$\mathcal{R}$!product on, $\cdot$|emph}} 
$$
fg=f\cdot g\cc M(f)\cap M(g)\to\R\,\text{ by }\,
(fg)(x)\equiv f(x)g(x)\,.\label{prodFun}
$$
\item The \underline{ratio\index{function!r@$\mathcal{R}$!ratio on, $/$|emph}} 
$$
f/g\cc 
M(f)\cap M(g)\setminus Z(g)\to\R\,\text{ by }\,(f/g)(x)\equiv{\textstyle\frac{f(x)}{g(x)}\,.} \label{diviFun}
$$
Recall that $Z(g)=\{x\in M(g)\cc\;g(x)=0\}$. 
\item Recall that the \underline{composition} 
$$
f(g)=f\circ g\cc M(f(g))\to\R,\ \  
M(f(g))=\{x\in 
M(g)\cc\;g(x)\in M(f)\}\,, 
$$
has values $f(g)(x)\equiv f(g(x))$.
\item For injective $f$ its \underline{inverse\index{function!r@$\mathcal{R}$!inverse on, $f^{-1}$|emph}} 
$$
f^{-1}\cc f[M(f)]\to\R,\ \ 
f^{-1}(y)\equiv x \iff f(x)=y\,,
$$ 
is congruent to the inverse in 
\end{enumerate}
\end{defi}
Operations in 1--4 are binary. The inverse in 5 is a~partial unary 
operation, it is not defined on non-injective functions. The ratio
of two functions is always defined, there is no forbidden division by
zero in the functional arithmetic. Recall that for $f$ in $\mathcal{R}$ and any set $X$ we have restriction 
$$
f\,|\,X\cc 
M(f)\cap X\to\R,\ \ (f\,|\,X)(x)\equiv f(x)\,. 
$$
In Chapter~\ref{chap_pr7} we add a~sixth operation on 
$\mathcal{R}$, the derivative $f\mapsto f'$.

A~semiring\index{semiring} is 
a~structure 
$$
\langle X,\,0_X,\,1_X,\,+,\,\cdot\rangle
$$ 
where $X$ is a~set, $+$ and $\cdot$ are 
commutative and associative operations on $X$, the elements $0_X$ and $1_X$ in $X$ are 
neutral to $+$ and $\cdot$, respectively, and $\cdot$ is distributive to $+$.

\begin{exer}\label{ex_aritmFunkci}
Prove the following proposition.    
\end{exer}

\begin{prop}[a~semiring on $\mathcal{R}$]\label{prop_semirFci}
Let $\mathcal{R}$ be as above, $0_{\mathcal{R}}
\equiv k_0(x)$ and $1_{\mathcal{R}}\equiv k_1(x)$. \underline{Then}
$$
\mathcal{R}_{\mathrm{smr}}\equiv\langle\mathcal{R},\,
0_{\mathcal{R}},\,1_{\mathcal{R}},\,+,\,\cdot\rangle\label{Rsmr}
$$
is a~semiring.
\end{prop}

A~function $g\in\mathcal{R}$ is \underline{invariant\index{invariant 
element in $\mathcal{R}_{\mathrm{smr}}$|emph}} in 
$\mathcal{R}_{\mathrm{smr}}$ with respect to $+$, resp. $\cdot$, if
for every function $h\in\mathcal{R}$ we have
$$
\text{$h+g=g$, resp. $h\cdot g=g$}\,.
$$
The proof of the following proposition
characterizing invariant elements in 
$\mathcal{R}_{\mathrm{smr}}$ is immediate.

\begin{prop}[the invariant function]
The real empty function $\emptyset_f$ is invariant to $+$ and $\cdot$. 
If $g\in\mathcal{R}$ with $g\ne\emptyset_f$, then 
$$
\text{$\emptyset_f+g=\emptyset_f\ne g$
and $\emptyset_f\cdot g=\emptyset_f\ne g$}\,,
$$
so that $g$ is invariant neither to $+$ nor to $\cdot$.
\end{prop}

\noindent
{\em $\bullet$ Differences of functions. }The \underline{difference\index{function!r@$
\mathcal{R}$!difference on, 
$-$|emph}} of functions $f$ and $g$ in $\mathcal{R}$ is the function 
$$
f-g\cc M(f)\cap M(g)\to\R,\ \ (f-g)(x)\equiv f(x)-g(x)\,.\label{fminusg}
$$ 

\begin{exer}\label{ex_OdecFci}
Always $f-g=f+(k_{-1}\cdot g)$.
\end{exer}

\begin{exer}\label{ex_platiTo}
Let $f,g,h\in\mathcal{R}$. Which of the implications in the (in general 
invalid) equivalence $f+g=h$ $\iff$ $f=h-g$ does gold?    
\end{exer}

\noindent
{\em $\bullet$ ``How elementary, dear Watson!''\,\index{Watson, 
John}\footnote{By \cite{dictionary} the only ``elementary'' statement in the work of A.~C. 
Doyle\index{Doyle, Arthur Conan} on Sherlock Holmes\index{Holmes, 
Sherlock} is found in the story {\em The Crooked Man}. There we read: “\ ‘Excellent!’ I [Watson] cried. ‘Elementary,’ said
he.” }} We introduce elementary functions, or EF. They are sometimes confused with basic elementary functions, or BEF, of the previous section. Recall that
\begin{eqnarray*}
\mathrm{BEF}&=&
\{\exp x,\,\log x,\,\sin x,\,\cos x,\,\tan x,\,
\cot x,\,\arcsin x,\,\arccos x,\\
&&\,\arctan x,\,\mathrm{arccot}\,x,\,0^x\}\cup\{k_c(x)\cc\;c\in\R\}
\cup\{a^x\cc\;a>0\}\cup\\
&&\cup\,\{x^b\cc\;b\in\R\}
\cup\{x^m\cc\;m\in\Z\}\,.
\end{eqnarray*}

\begin{defi}[elementary functions~1]\label{def_obecEF}
A~function $f\in\mathcal{R}$ is \underline{elementary} if 
there is a~tuple of $n\in\N$ functions 
$$
\langle f_1,\,f_2,\,\ds,\,f_n\rangle
\in\mathcal{R}^n\,,
$$
called a~\underline{generating word\index{generating word!elementary function|emph}} of $f$,
such that $f_n=f$ and for every $i\in[n]$
either $f_i\in\mathrm{BEF}$ or there exist indices $j,k<i$ for which 
$$
f_i=f_j+f_k\vee f_i=f_j\cdot f_k\vee
f_i=f_j/f_k\vee f_i=f_j(f_k)\,.
$$
\end{defi}
The set of  \underline{elementary functions\index{elementary 
functions, EF|emph}} is denoted by EF.\label{EF}

\begin{exer}\label{ex_onGenWor}
Every function $f_i$ in the
generating word of $f$ is elementary.
\end{exer} 
Said less formally,  we get EF from BEF by repeated addition, 
multiplication, division, and composition. For 
example, the \underline{identity function\index{function!identity, 
$\mathrm{id}(x)$|emph}\index{identity function, $\mathrm{id}(x)$|emph}} 
$$
\mathrm{id}(x)=x\equiv\mathrm{id}_{\R}(x)
$$ 
is elementary because $\mathrm{id}(x)=\log(\exp x)$ or 
$\mathrm{id}(x)=x^1$ with $1\in\Z$. Recall that $x^1$ with $1\in\R$ is $\mathrm{id}(x)\,|\,[0,+\infty)$.

\begin{exer}\label{ex_absHodn}
The absolute value
$|x|$ is in $\mathcal{F}(\R)$ and is 
elementary.    
\end{exer}

\begin{exer}\label{ex_deleniNulou}
What is $k_1(x)/k_0(x)$?   
\end{exer}

\begin{exer}\label{ex_prazdnaElem}
Is the empty function $\emptyset$ elementary? 
\end{exer}

\begin{exer}\label{ex_pseudoInverz}
For every $f\in\mathrm{EF}$ there is a~unique $g\in\mathrm{EF}$ such that $M(g)=M(f)$ and $f+g=k_0\,|\,M(f)$. 
\end{exer}

\begin{exer}\label{ex_EFsdefObZ}
Find functions $f,g\in\mathrm{EF}$ such that
$${\textstyle
M(f)=\Z\,\text{ and }\,M(g)=\R\setminus(\{0\}\cup\{\frac{1}{n}\cc\;n\in\Z\setminus
\{0\}\})\,.
}
$$
\end{exer}

\noindent
{\em $\bullet$ Restrictions to intervals. }We show that elementary functions are preserved by restrictions to intervals. 

\begin{prop}[getting intervals]\label{prop_intRest}
Let $a\in\R$. The functions 
$$
f_a\equiv(a-x)^{1/2},\ g_a\equiv(x-a)^{1/2},\ 
F_a\equiv\log(x-a)\,\text{ and }\,G_a\equiv\log(x-a)
$$
are elementary and have respective definition domains 
$$
(-\infty,\,a],\ [a,\,+\infty),\ 
(-\infty,\,a)\,\text{ and }\, 
(a,\,+\infty)\,.
$$
\end{prop}
\duk
This follows immediately from Exercise~\ref{ex_OdecFci}, and 
Definitions~\ref{def_ZEF}, \ref{def_oprOnR} and 
\ref{def_obecEF}. 
\kduk

Recall the real intervals as given in Proposition~\ref{prop_intervaly}: 
\begin{eqnarray*}
\mathcal{I}&\equiv&
\{\emptyset,\  \{a\},\  \R,\  (a,b),\ (-\infty,a),\ (a,+\infty),\  (a,b],\  [a,b),\  [a,b],\\\  
&&(-\infty,a],\  [a,+\infty)\cc\;a,\,b\in\R,\,a<b\}\,.
\end{eqnarray*}

\begin{prop}[interval restrictions]\label{prop_restrEF}
Let $f\in\mathrm{EF}$ be an elementary function and $I\in\mathcal{I}$ be an interval. \underline{Then} 
$$
f\,|\,I\in\mathrm{EF}\,.
$$
\end{prop}
\duk
It suffices to show that for every interval $I\in\mathcal{I}$ the 
function 
$$
h_I(x)\equiv k_0(x)\,|\,I
$$ 
is elementary, because then
$f\,|\,I=f+h_I$. 
If $I$ is $\emptyset$ or $\R$ then the result holds trivially. In the remaining cases it is easy to get $h_I$ as a~difference of two functions in Proposition~\ref{prop_intRest}, or as a~sum 
of two such differences. For example, 
$$
h_{\{a\}}=f_a-f_a+g_a-g_a,\ 
h_{(-\infty,a]}=f_a-f_a\,\text{ and }\,h_{[a,b)}=F_b-F_b+g_a-g_a\,.
$$
\kduk

\noindent
{\em $\bullet$ Really basic elementary functions }form a~subset of BEF.

\begin{defi}[RBEF]\label{def_skutecneZEF}
\underline{Really basic elementary functions\index{really basic elementary functions, RBEF|emph}} are
$$
\mathrm{RBEF}\equiv
\{\exp x,\,\log x,\,\sin x,\,
\arcsin x\}\cup\{k_c(x)\cc\;c\in\R\}\cup\,\{x^b\cc\;b\in\R\}\,.\label{RBEF}
$$
\end{defi}
We show that the other functions in BEF are redundant in generation of EF.

\begin{prop}[RBEF suffice]\label{prop_RBEFsuff}
Every function 
$$
f\in\mathrm{BEF}\setminus
\mathrm{RBEF}
$$ 
has a~generating word according to the restricted form of 
Definition~\ref{def_obecEF} in which {\em BEF} is replaced with {\em RBEF}.
\end{prop}
\duk
First we observe that the identity function is elementary in the
restricted sense because 
$$
x=\mathrm{id}(x)=\log(\exp x)\,.
$$
(We cannot use the expression
$\mathrm{id}(x)=x^1$ because functions $x^m$, $m\in\Z$, are not available.) Then we express functions in 
$\mathrm{BEF}\setminus
\mathrm{RBEF}$ by functions in $\mathrm{RBEF}$ as follows.

\begin{itemize}
\item $\cos x=\sin(x+\frac{\pi}{2})=\sin(x+k_{\pi/2}(x))$.
\item $\tan x=\sin x/\cos x$.
\item $\cot x=\cos x/\sin x$.
\item $\arccos 
x=\frac{\pi}{2}-\arcsin x=k_{\pi/2}(x)+k_{-1}(x)\cdot\arcsin x$.
\item $\arctan x=\arcsin(x/(1+x\cdot x)^{1/2})$, see Exercise~\ref{ex_arctann}.
\item $\mathrm{arccot}\,x=\frac{\pi}{2}-\arctan x=k_{\pi/2}(x)+
k_{-1}(x)\cdot\arctan x$.
\item $a^x=\exp(x\cdot\log 
a)=\exp(x\cdot k_{\log a}(x))$, for real $a>0$.
\item $0^x=k_0(x)+\log x-\log x$.
\item $x^m=x\cdot x\cdot\ldots\cdot x$, for $m\in\N$.
\item $x^0=k_1(x)$, for $0\in\Z$, and
\item $x^m=k_1(x)/(x\cdot x\cdot\ldots\cdot x)$ with $-m$ factors, for $m\in\Z$ with
$m<0$.
\end{itemize} 
\kduk

\noindent
We write $x\cdot x$ instead of the ambiguous $x^2$ which can be interpreted to be in 
$\mathcal{F}(\R)$ or in $\mathcal{F}([0,+\infty))$, depending on whether $2\in\Z$ or $2\in\R$.

\begin{exer}\label{ex_arctann}
Prove the equality of functions 
$${\textstyle
\arctan x=\arcsin\big(\frac{x}{(1+x\cdot x)^{1/2}}\big)\,.
}
$$
\end{exer}

Hence we have the next simpler definition of elementary functions.

\begin{defi}[elementary functions~2]\label{def_obecEF2}
In Definition~\ref{def_obecEF}\index{elementary functions, EF|emph} of {\em EF},
the set {\em BEF} may be replaced with the smaller set {\em RBEF}.  
\end{defi}

\begin{exer}\label{ex_jesteRedu}
In Definition~\ref{def_skutecneZEF}, which functions can be further deleted from $\{x^b\cc\;b\in\R\}$ so that Proposition~\ref{prop_RBEFsuff} still holds?
\end{exer}

The elementary function
$\frac{|x|}{x}$ ($\in
\mathcal{F}(\R\setminus\{0\})$ is $-1$ for $x<0$, and $1$ for $x>0$. It resembles the function \underline{signum\index{signum, $\sgn\,x$|emph}}\label{signum} $\sgn(x)$ in 
$\mathcal{F}(\R)$, given by $\sgn(x)=-1$ for $x<0$, $\sgn(0)=0$ and $\sgn(x)=1$ for $x>0$. 
However, signum is not elementary.

\begin{prop}[signum is not elementary]\label{prop_sgnNeniEF}
The function $\sgn\,x\not\in\mathrm{EF}$.    
\end{prop}
\duk
By Definition~\ref{def_contOnaSet}
and Theorem~\ref{thm_EFjsouSpoj} every elementary function is continuous, but 
signum is not continuous. 
\kduk

\noindent
Maybe one should reconsider the definition of elementary functions so that signum is admitted to them. We return to this matter in {\em MA~1${}^+$}.

\begin{exer}\label{ex_EFjedenBod}
Give examples of functions in {\em EF} which are not differentiable at
some points of their definition domains.\index{elementary functions, 
EF|)}    
\end{exer}

\section[Polynomials and rational functions]{Polynomials and rational functions}\label{sec_polyRac}

Restricted forms of Definition~\ref{def_obecEF} produce  
polynomials and rational functions. 

\medskip\noindent
{\em $\bullet$ Polynomials. }Our definition is not conventional. 

\begin{defi}[POL]\label{def_polynomy}
A~function $f\in\mathcal{R}$ is a~\underline{polynomial\index{polynomials, POL|emph}}\label{polynom} if 
there exists an $n$-tuple of functions 
$$
\langle f_1,\,f_2,\,\ds,\,f_n\rangle
\in\mathcal{R}^n \,,
$$
so called \underline{generating word\index{generating word!polynomial|emph}} of $f$,
such that $f_n=f$ and for every $i=1,2,\ds,n$
we have $f_i(x)=k_c(x)$ for some $c\in\R$ or 
$f_i(x)=\mathrm{id}(x)$ or there exist indices $j,k<i$ for which $f_i=f_j+f_k$ or $f_i=f_j\cdot f_k$.
We denote the set of polynomials by {\em POL}.
\end{defi}
Again, not only the last function 
$f_n=f$ but every function $f_1,f_2,\ds,f_{n-1}$ in the generating word is 
a~polynomial. In our approach, polynomials arise from the identity function 
and constants by addition and multiplication. Then it is trivial that the sum and product 
of two polynomials is again a~polynomial because it is built in in the definition. In the standard definition of polynomials this fact becomes a~nontrivial, technical result. 

However, we have to prove that our polynomials coincide with the standard
ones. For $f\in\mathcal{R}$ and $n\in\N$ we define $f^n\equiv f\cdot f\cdot\ldots\cdot f$ with 
$n$ factors $f$, and set 
$f^0\equiv k_1\,|\,M(f)$.

\begin{prop}[on zero-th power]\label{prop_zeroPower}
If $f\in\mathrm{EF}$ \underline{then} $f^0\in\mathrm{EF}$.   
\end{prop}
\duk
Indeed, $f^0=f-f+k_1=f+k_{-1}\cdot f+k_1$.
\kduk

We call the constant function $k_0(x)$ the \underline{zero 
polynomial\index{polynomials, POL!zero polynomial|emph}}. We call a~function $f\in\mathcal{R}$ a~\underline{canonical 
polynomial\index{polynomials, POL!canonical|emph}} if for some $n+1$ real numbers $a_0$, $a_1$, $\ds$, $a_n$, where $n\in\N_0$ and $a_n\ne0$, we have 
$$
f=k_{a_0}\cdot\mathrm{id}^0+k_{a_1}\cdot\mathrm{id}^1+\ds+k_{a_n}\cdot\mathrm{id}^n\,.
$$
We abbreviate it by writing 
$${\textstyle
f(x)=\sum_{j=0}^n a_jx^j\,\text{ or }\,f(x)=a_0+a_1x+\ds+a_nx^n\,,
}
$$
and call the $n+1$-tuple 
$$
\langle a_0,\,a_1,\,\ds,\,a_n\rangle
$$ 
the \underline{canonical 
form\index{polynomials, POL!canonical form of|emph}} of $f$ (we show that it is unique). It is clear that the zero 
polynomial and every canonical polynomial is a~polynomial. Now we prove the opposite.

\begin{exer}\label{ex_zeroCan}
If $f(x)=\sum_{j=0}^n a_jx^j$ is a~canonical polynomial, 
then $|Z(f)|\le n$, that is, $f(x)$ has at most $n$ zeros. 
\end{exer}
It follows that $k_0(x)$ differs from every canonical polynomial.

\begin{exer}\label{ex_rootFactor}
If $f(x)\in\mathrm{POL}$ and $f(b)=0$ then $f(x)=g(x)\cdot(x-b)$ where $g(x)\in\mathrm{POL}$.
\end{exer}

\begin{thm}[polynomials]\label{thm_oPoly}
Every\index{theorem!polynomials|emph}
nonzero polynomial is a~canonical polynomial.
Canonical polynomials have unique canonical forms.
\end{thm}
\duk
Let $f\in\mathrm{POL}$ with $f\ne k_0$ and let 
$$\langle
f_1,\,f_2,\,\ds,\,f_n\rangle,\ f=f_n\,, 
$$
be the generating word of $f$ by Definition~\ref{def_polynomy}. 
We prove by induction on $n$ that $f$ is a~canonical polynomial. If 
$n=1$ then $f=k_c$ or $f=\mathrm{id}$ and it is true. 
Let $n>1$. If $f=k_c$ or $f=\mathrm{id}$ we are again done. 
We assume that there are $j,k\in[n-1]$ such that (i) $f=f_j+f_k$ or   
(ii) $f=f_jf_k$. In case (i) both $f_j$ and $f_k$ cannot be $k_0$. If exactly one of them is $k_0$ then $f$ is a~canonical polynomial by induction. If none of $f_j$ and $f_k$ is $k_0$ then we are done by
Exercise~\ref{ex_sumCan}. In case (ii) none of $f_j$ and $f_k$ is $k_0$
and we are done by Exercise~\ref{ex_prodCan}. 

We prove uniqueness of canonical forms. If $f\in\mathcal{R}$ were 
a~canonical polynomial with two different canonical forms, by 
Exercise~\ref{ex_diffCan} $k_0=f-f$ would be a~canonical polynomial, 
which by Exercise~\ref{ex_zeroCan} is impossible. 
\kduk
\vspace{-3mm}
\begin{defi}[degree]\label{def_degPol}
The \underline{degree\index{polynomials, POL!degree of, deg|emph}} of a~nonzero 
polynomial $p(x)$ is the number $d\in\N_0$ in the canonical form 
$$
{\textstyle
p(x)=\sum_{i=0}^d a_ix^i\,.
}
$$
We denote it by $\deg p(x)$ or by 
$\deg(p(x))$.\label{degree}
\end{defi}

\begin{exer}\label{ex_sumCan}
Show that the sum of two canonical polynomials is a~canonical 
polynomial or the zero polynomial.   
\end{exer}

\begin{exer}\label{ex_prodCan}
Show that the product of two canonical polynomials is a~canonical polynomial.     
\end{exer}

\begin{exer}\label{ex_diffCan}
Show that the difference of two canonical polynomials with different canonical forms is a~canonical polynomial.     
\end{exer}

\begin{exer}\label{ex_POLjeOI}
\label{ex_polObIn}
Prove the following proposition. 
\end{exer}

\begin{prop}[POL is an integral domain]\label{prop_polJeOI}
The\index{polynomials, POL!form an integral domain} 
structure
$$
\mathrm{POL_{id}}\equiv\langle\mathrm{POL},\,k_0,\,k_1,\,+,\,\cdot\rangle\label{polynDom}
$$ 
is an integral domain.    
\end{prop} 

\begin{exer}\label{ex_polyIsom}
Show that the ring $\mathrm{POL_{id}}$ is isomorphic to the ring $\R[x]$ of real polynomials
in abstract algebra. 
\end{exer}

\noindent
{\em $\bullet$ Rational functions. }Rational functions are usually defined as ratios of
polynomials. We introduce them by augmenting
Definition~\ref{def_polynomy} with  division. Now the
resulting functions need not be everywhere defined.

\begin{defi}[RAC]\label{def_RacFce}
A~function $f\in\mathcal{R}$ is \underline{rational\index{rational functions, RAC|emph}}\label{rational} if 
there is a~tuple of $n\in\N$ functions 
$$
\langle f_1,\,f_2,\,\ds,\,f_n\rangle
\in\mathcal{R}^n\,, 
$$
called a~\underline{generating word\index{generating word!rational function|emph}} of the rational function $f$, such that $f_n=f$ such that for every $i\in[n]$
we have $f_i(x)=k_c(x)$ for some $c\in\R$ or 
$f_i(x)=\mathrm{id}(x)$ or there exist indices $j,k<i$ for which 
$$
f_i=f_j+f_k\vee f_i=f_j\cdot f_k
\vee f_i=f_j/f_k\,.
$$
The set of rational functions is denoted as 
{\em RAC}.
\end{defi}
As before, every function $f_i$, $i\in[n]$, in the generating word of 
a~rational function $f=f_n$ is rational. Our rational functions 
arise from constants and the identity by repeated addition, 
multiplication and division.

\begin{exer}\label{ex_onRatFun}
Show that $\mathrm{POL}\sus\mathrm{RAC}$ and that the sum, product and ratio of two   
rational functions is a~rational function.
\end{exer}

For instance, $\frac{1}{x}\equiv
k_1/\mathrm{id}$ is a~rational function and $M(\frac{1}{x})=\R\setminus\{0\}$. Another
rational function is the empty function:  
$\emptyset=k_1/k_0$ and $M(\emptyset)=\emptyset$. We prove for rational functions an analog of 
Theorem~\ref{thm_oPoly}, but we need for it three lemmas on computing with ratios in 
$\mathcal{R}$.

\medskip\noindent
{\em $\bullet$ The arithmetic of ratios. }We begin with addition. 
Recall the operations on $\mathcal{R}$ introduced in 
Definition~\ref{def_oprOnR}.

\begin{lemma}[the sum]\label{lem_souZlom}
If $f_1$, $g_1$, $f_2$ and $g_2$ are in $\mathcal{R}$, \underline{then}
$${\textstyle
\frac{f_1}{g_1}+\frac{f_2}{g_2}=\frac{f_1g_2+f_2g_1}{g_1g_2}\,.
}
$$
\end{lemma}
\duk
We denote the left-hand side of the equality by $F$ and the right-hand side by $G$.
Then $M(F)=M(f_1/g_1)\cap M(f_2/g_2)$ equals to
$$
\big((M(f_1)\cap M(g_1))\setminus Z(g_1)\big)\cap
\big((M(f_2)\cap M(g_2))\setminus Z(g_2)\big)
$$
and $M(G)=(M(f_1g_2+f_2g_1)\cap M(g_1g_2))\setminus Z(g_1g_2)$ equals to
$$
\big((M(f_1)\cap M(g_2))\cap(M(f_2)\cap M(g_1))\cap
(M(g_1)\cap M(g_2))\big)\setminus
(Z(g_1)\cup Z(g_2))\,.
$$
The equality $Z(g_1g_2)=Z(g_1)\cup Z(g_2)$ follows from the fact that $\R$, being a~field, is an
integral domain. Hence $M(F)=M(G)$ because both displayed sets are equal to the set
$$
\big(M(f_1)\cap M(f_2)\cap M(g_1)\cap M(g_2)\big)\setminus\big(Z(g_1)\cup Z(g_2)\big)\,.
$$
The arithmetic in $\R$ shows that
$F(x)=G(x)$ for every $x$ in $M(F)=M(G)$. Hence $F=G$.
\kduk

The proof of the second lemma is similar and is left to an exercise.

\begin{exer}\label{ex_prodRatio}
Prove the next lemma.     
\end{exer}

\begin{lemma}[the product]\label{lem_prodZlom}
If $f_1$, $g_1$, $f_2$ and $g_2$ are in $\mathcal{R}$, \underline{then}
$${\textstyle
\frac{f_1}{g_1}\cdot\frac{f_2}{g_2}=\frac{f_1f_2}{g_1g_2}\,.
}
$$
\end{lemma}

The next exercise
shows that the third lemma has to be more complicated.

\begin{exer}\label{ex_naRatioRatios}
Find four functions $f_1$, $g_1$, $f_2$ and $g_2$ in $\mathcal{R}$ such that 
$${\textstyle
\frac{f_1/g_1}{f_2/g_2}\ne\frac{f_1g_2}{f_2g_1}\,.
}
$$   
\end{exer}

\begin{lemma}[the ratio]\label{lem_ratioZlom}
If $f_1$, $g_1$, $f_2$ and $g_2$ are in $\mathcal{R}$, \underline{then}
$${\textstyle
\frac{f_1/g_1}{f_2/g_2}=
\frac{f_1(g_2)^2}{f_2g_1g_2}\,.
}
$$
\end{lemma}
\duk
We again denote the left-hand side of the equality by $F$ and the right-hand side by $G$.
Then $M(F)=(M(f_1/g_1)\cap M(f_2/g_2))\setminus Z(f_2/g_2)$ equals to
$$
\big((M(f_1)\cap M(g_1)\setminus Z(g_1))\cap(M(f_2)\cap M(g_2)\setminus Z(g_2))\big)\setminus
\big(Z(f_2)\setminus Z(g_2)\big)
$$
and $M(G)=(M(f_1g_2^2)\cap M(f_2g_1g_2))\setminus Z(f_2g_1g_2))$ equals to
$$
\big(M(f_1)\cap M(f_2)\cap M(g_1)\cap M(g_2)\big)\setminus\big(Z(f_2)\cup Z(g_1)\cup Z(g_2)\big)\,.
$$
Thus $M(F)=M(G)$ because both displayed sets are equal to
$$ 
\big(M(f_1)\cap M(f_2)\cap M(g_1)\cap M(g_2)\big)\setminus\big(Z(g_1)\cup Z(f_2)\cup Z(g_2)\big)\,.
$$
The arithmetic in $\R$ shows that
$F(x)=G(x)$ for every $x$ in $M(F)=M(G)$. Hence $F=G$.
\kduk

\noindent{\em $\bullet$ Canonical forms of rational functions. }It 
follows from Definitions~\ref{def_polynomy} and 
\ref{def_RacFce} that if $f,g\in\mathrm{POL}$ then 
$f/g\in\mathrm{RAC}$. Indeed,
suppose that 
$$
\langle f_1,\,f_2,\,\ds,\,f_m=f\rangle\,\text{ and }\,\langle g_1,\,g_2,\,\ds,\,g_n=g\rangle
$$ 
are respective generating words of the polynomials $f$ and $g$ by
Definition~\ref{def_polynomy}. 
Then 
$$
\langle f_1,\,\ds,\,f_m,\,g_1,\,\ds,\,g_n,\,
f/g\rangle
$$
is a~generating word of the rational function $f/g$ by Definition~\ref{def_RacFce}. We prove
the reverse, that every rational function is a~ratio of two 
polynomials.

\begin{thm}[rational functions]\label{thm_oRacFci}
Let $r\in\mathrm{RAC}$ with $r\ne\emptyset$.
\underline{Then} there exist polynomials $p$ and $q$ 
such that 
$$
p/q=r\,. 
$$
Hence $q\ne k_0$ and $M(r)=\R\setminus Z(q)$, where 
$Z(q)$ is a~finite set by Exercise~\ref{ex_zeroCan}.
We say that the pair 
$(p,q)$, written also as $p/q$, is a~\underline{canonical 
form\index{rational functions, RAC!canonical form of|emph}} of $r$.
\end{thm}
\duk
Let $r$ be a~nonempty rational function. By 
Definition~\ref{def_RacFce}, it has 
a~generating word 
$$
\langle f_1,\,f_2,\,\ds,\,f_n=r\rangle\,. 
$$
We proceed by induction on 
$n$. If $n=1$ then $r$ is a~constant or the identity. Thus we have canonical forms $f=k_c/k_1$, $c\in\R$, or $\mathrm{id}/k_1$. Let $n>1$. If 
$r=f_n$ is a~constant 
or the identity, we are in the
previous case. Else there exist $j,k\in\N$ with $j,k<n$ such
that (i) $r=f_j+f_k$ or (ii) $r=f_jf_k$ or (iii) $r=f_j/f_k$. 
Both $f_j$ and $f_k$ are nonempty, otherwise we would have $r=\emptyset$. By induction, $f_j$ and $f_k$ 
have canonical forms $f_j=p_1/q_1$ and $f_k=p_2/q_2$. In case (i) we get from Lemma~\ref{lem_souZlom} that
$r$ has the canonical form
$${\textstyle
r=f_j+f_k=\frac{p_1q_2+p_2q_1}{q_1q_2}\,.
}
$$
In the case (ii) we get from Lemma~\ref{lem_prodZlom} that
$r$ has the canonical form
$${\textstyle
r=f_jf_k=\frac{p_1p_2}{q_1q_2}\,.
}
$$
In the case (iii) we get from Lemma~\ref{lem_ratioZlom} that
$r$ has the canonical form
$${\textstyle
r=f_j/f_k=\frac{p_1(q_2)^2}{q_1q_2p_2}\,.
}
$$
\kduk

\noindent
Unlike for polynomials, canonical forms of rational functions are not
unique. For example, for
$r\equiv\mathrm{id}/\mathrm{id}$ ($\in\mathrm{RAC}$) every pair
$x^k/x^k$, $k\in\N$, is a~canonical form. Note that $r\ne k_1/k_1$ 
because $M(r)=\R\setminus\{0\}$ but $M(k_1/k_1)=\R$.

We introduce a~\underline{congruence\index{congruence!on rac@on $\mathrm{RAC}\setminus\{\emptyset\}$|emph}} $\sim$
on the set $\mathrm{RAC}\setminus\{\emptyset
\}$ by 
$$
r\sim s\stackrel{\text{def}}{\iff} r\,|\,M(s)=s\,|\,M(r)\;.
$$

\begin{exer}\label{ex_jetoSkutRE}
Show that $\sim$ is an equivalence relation.   
\end{exer}
For example, $k_1\sim x/x\sim(x\cdot(x-1))/(x\cdot(x-1))$.

\begin{exer}\label{ex_RCjeTel}
Prove the next proposition. 
\end{exer}

\begin{prop}[$\mathrm{RAC}$ as a~field]\label{prop_RCjeTel}
The structure
$$
\mathrm{RAC_{FI}}\equiv\langle(\mathrm{RAC}\setminus\{\emptyset\})/\!\sim,\,[k_0]_{\sim},\,[k_1]_{\sim},\,
+,\,\cdot\rangle\label{ratioFie}
$$ 
is a~field.   
\end{prop}
We call this field the \underline{field of rational 
maps\index{field of rational maps, $\mathrm{RAC_{fi}}$|emph}}.

\begin{exer}\label{ex_isomRacFun}
Show that the field  $\mathrm{RAC_{fi}}$ is isomorphic 
to the field of rational functions $\R(x)$ used in abstract algebra.
\end{exer}           
\chapter[Continuous functions]{Continuous functions}\label{chap_pr6}

In Section~\ref{sec_BlumbVeta}
we introduce dense and sparse sets and state Blumberg's Theorem~\ref{thm_BlumbergThm}: every function 
$$
f\cc\R\to\R
$$
has a~continuous restriction $f\,|\,M$ to a~set $M$ dense
in $\R$. We prove it in~{\em MA~1${}^+$}. In the extending 
Section~\ref{sec_heiZerFraSie} we prove Sierpi\'nski's\index{Sierpi\'nski, Wac\l aw} Theorem~\ref{thm_sierp}: for a~function $f\cc\R\to\R$ the equivalence
$$
\text{$f$ is continuous}\iff
\text{$f$ is sequentially continuous}
$$
holds in ZF. In Section~\ref{sec_poceSpoFun}
we show in Theorem~\ref{thm_poceSpoFun} that the set of continuous
functions $f\cc\R\to\R$ is in bijection with $\R$. The main result of
Section~\ref{sec_nabyvMezih} is Theorem~\ref{thm_mezihodnoty} 
which says that continuous functions attain every intermediate value. 

Section~\ref{sec_kompaktnost} introduces real compact, open and
closed sets. We downgrade the minimax theorem to a~mere 
Corollary~\ref{cor_priMaxi} of Theorem~\ref{thm_obrKomp} by which continuous images of compacts are
compact. Basic properties of open and closed sets are discussed. 
Theorem~\ref{thm_KompvR} 
characterizes compact sets.
We include Baire's Theorem~\ref{thm_Baire} because it is 
needed in results on real analytic functions.
Section~\ref{sec_UC} is devoted to uniform continuity.  By 
Theorem~\ref{thm_steSpo} any
continuous function with compact definition domain is uniformly
continuous. Theorem~\ref{thm_steSpoRoz} says that every uniformly
continuous function has a~unique 
continuous extension to the closure of the definition domain. In 
Theorem~\ref{thm_HMCminmax} we generalize with the 
help of the extension theorem the minimax theorem to bounded
definition domains and UC functions.

Section~\ref{sec_aritSpoj} treats interactions between continuity 
and operations on $\mathcal{R}$. Theorem~\ref{thm_aritSpojitosti} 
concerns arithmetic operations, and Theorem~\ref{thm_spjMocRady1} 
is devoted to continuity of sums of power series. 
Theorem~\ref{thm_spojSloz} deals with composition. In the 
culminating Theorem~\ref{thm_spojInverzu} we 
treat continuity of inverse functions. The chapter concludes 
with Theorem~\ref{thm_EFjsouSpoj}: every elementary 
function is continuous.

\section[Globally continuous functions]{Globally continuous functions}\label{sec_BlumbVeta}

\noindent 
{\em $\bullet$ Global continuity. }Recall that by Definition~\ref{def_spojVbode} 
a~function $f\in\mathcal{R}$ is continuous at 
a~point $a\in M(f)$ if $\forall\ep\exists\de\,
f[U(a,\de)]\sus U(f(a),\ep)$. 

\begin{defi}[global continuity]\label{def_contOnaSet}
Let $f\in\mathcal{R}$ and $X\sus\R$. We say that $f$ is 
\underline{continuous on\index{function!continuous on a 
set|emph}\index{continuity of functions!on a set|emph} $X$} if the function $f$  
is continuous at every point $b\in M(f)\cap X$. If $f\in\mathcal{R}$ is
continuous on $M(f)$, it is a~\underline{continuous\index{function!continuous|emph} function}. 
\end{defi}

\noindent
We denote the subset of continuous functions in $\mathcal{F}(M)$ by $\mathcal{C}(M)$.\index{function!cm@$\mathcal{C}(M)$|emph}\label{Cm} We denote by 
$${\textstyle
\mathcal{C}\equiv\bigcup_{M\sus\R}
\mathcal{C}(M)\label{C}
}
$$
the subset of continuous functions in $\mathcal{R}$.\index{function!c@$\mathcal{C}$|emph}  $\mathcal{R}\setminus\mathcal{C}$ are  \underline{discontinuous\index{function!discontinuous|emph}} functions.

\begin{exer}\label{ex_spojKon}
Every function $f\in\mathcal{R}$ with finite domain $M(f)$ is continuous.
\end{exer}

\begin{exer}\label{ex_spojKonst}
For every number $c\in\R$ the constant function $k_c(x)$ is continuous.    
\end{exer}

\begin{exer}\label{ex_spojIden}
The identity function $\mathrm{id}(x)=x$ is continuous.    
\end{exer}

\begin{prop}[restriction~1]\label{prop_spojRest}
Let $f\in\mathcal{C}$ and $X\sus\R$. \underline{Then} 
$$
f\,|\,X\in\mathcal{C}\,.
$$
\end{prop}
\duk
Let $b\in M(f\,|\,X)$ and let an $\ep$ be given. Thus $b\in M(f)$. Since 
$f\in\mathcal{C}$, there is a~$\de$ such that 
$f[U(b,\de)]\sus U(f(b),\ep)$.
Since $U(b,\de)\cap M(f)\cap X\sus U(b,\de)\cap M(f)$, we have 
$$
(f\,|\,X)[U(b,\,\de)]\sus f[U(b,\,\de)]\sus U(f(b),\,\ep)\,.
$$ 
\kduk

\noindent
{\em $\bullet$ Sequential continuity. }The next definition is 
sometimes used as the primary definition of continuity, both pointwise and global.

\begin{defi}[sequential continuity]\label{def_seqCont}
Let $f\in\mathcal{F}(M)$. 
\begin{enumerate}
\item If $b\in M$ and for every $(a_n)\sus M$ with $a_n\to b$ we have $f(a_n)\to f(b)$, we say that $f$ is \underline{sequentially continuous at $b$}.
\item If $f$ is sequentially continuous at every $b\in M$, we say that the function $f$ is \underline{sequentially continuous}.
\end{enumerate}
\end{defi}

Exercise~\ref{ex_HeiDefSpo} established that both definitions of continuity are equivalent:

\begin{prop}[Heine's formulation]\label{prop_HeineDefCont}
Let $f\in\mathcal{F}(M)$ and $b\in M$.
\begin{enumerate}
\item The function $f$ is continuous at $b$ $\iff$ $f$ is sequentially continuous at $b$.
\item The function $f$ is continuous $\iff$ $f$ is sequentially continuous.
\end{enumerate}    
\end{prop}
Both implications $\Rightarrow$ are easily established in ZF. Both opposite implications $\Leftarrow$  in general require 
the axiom of choice. 

\begin{exer}\label{ex_jesteJedSpoj}
See article \cite{herr} on equivalents of the axiom of choice in analysis and topology.    
\end{exer}
In the next section we prove a~surprising result due to 
W.~Sierpinski: if $M\sus\R$ is an open set, then the second implication $\Leftarrow$ in the proposition can be proven 
without the axiom of choice.

\medskip\noindent
{\em $\bullet$ Dense and sparse sets. }Let $M,N\sus\R$. The set $N$ is \underline{dense\index{set!dense|emph}} 
in the set $M$ if
for every point $a\in M$ and every $\de$ we have 
$$
U(a,\,\de)\cap M\cap N\ne\emptyset\,.
$$

\begin{exer}\label{ex_hustAposl}
$N$ is dense in $M$ $\iff$ for every point $a\in M$ there is a~sequence $(b_n)\sus M\cap N$ such that $\lim b_n=a$. 
\end{exer}

\begin{exer}\label{ex_oveJsouHust}
Show that both sets $\Q$ and $\R\setminus\Q$ are dense in $\R$.    
\end{exer}

Let $M,N\sus\R$. The set $N$ is \underline{sparse\index{set!sparse|emph}}
in the set $M$ if for every point $a\in M$ and every $\de$ there exist a~$b\in M$ and $\theta$ such that 
$U(b,\theta)\sus U(a,\de)$ and 
$$
U(b,\de)\cap M\cap N=\emptyset\,.
$$

\begin{exer}\label{ex_naRidke}
Show that $N\equiv\{\frac{1}{n}\cc\;n\in\N\}$ is sparse in $M\equiv[0,1]$. 
\end{exer}

\begin{prop}[density and continuity]\label{prop_oHusteMn}
Let $f,g\in\mathcal{C}(M)$, let $N\sus\R$ be dense in $M$ and let $f\,|\,N=g\,|\,N$. \underline{Then} 
$$
f=g\,.
$$
\end{prop}
\duk
Let $b\in M$ and $(a_n)\sus M\cap N$ have $\lim a_n=b$ (Exercise~\ref{ex_hustAposl}). Using (H) we have 
$$
f(b)=f(\lim a_n)=\lim f(a_n)=\lim g(a_n)=g(\lim a_n)=g(b)\,.
$$
\kduk

A~function $g\in\mathcal{C}$ is 
a~\underline{kernel\index{function!continuous on a set!kernel|emph}} of another 
function
$f\in\mathcal{C}$ if $g$ is a~restriction of $f$ and $M(g)$ is 
dense in $M(f)$. Using the previous proposition we easily reconstruct $f$ from $g$.

\begin{prop}[kernels]\label{prop_jadro}
Every continuous function $f$ in $\mathcal{C}$ has an at most countable kernel.    
\end{prop}
\duk
It suffices to show that every set $M\sus\R$ has an at most countable  
dense subset $N$ (Exercise~\ref{ex_whyDense}). Let 
$M\sus\R$. We obtain such set $N$ by using the axiom of choice\index{axiom!of choice, AC} and taking one element from 
every nonempty intersection $(\al,\be)\cap M$, where $\al<\be$ are 
fractions.
\kduk
\vspace{-3mm}
\begin{exer}\label{ex_whyDense}
Why does such set $N$ suffice?    
\end{exer}

\noindent{\em $\bullet$ Blumberg's theorem. }In 1922 H.~Blumberg discovered the following 
theorem.

\begin{thm}[Blumberg's]\label{thm_BlumbergThm}
For every\index{theorem!Blumberg's}
function 
$$
f\cc\R\to\R
$$ 
there exists a~set $M\sus\R$ such that
$$
\text{$M$ is dense in $\R$ and 
$f\,|\,M$ is continuous}\,.
$$
\end{thm}

\noindent
{\em Henry Blumberg\index{Blumberg, Henry} (1886--1950)} was born in 
Lithuania in the town 
Žagar$\mathrm{\dot{e}}$,\index{zagare@Žagar$\mathrm{\dot{e}}$} but the 
family emigrated to America already in 1891. We prove Blumberg's theorem 
in {\em MA~1${}^+$}.

\begin{exer}\label{ex_onBluThm}
Find such set $M$ for Riemann's function $r(x)$ (it is defined before Proposition~\ref{prop_RiemannFunkce}).  
\end{exer}

\section[Sierpi\'nski's theorem]{Sierpi\'nski's theorem}\label{sec_heiZerFraSie}

In \cite{herr_stre}, eight equivalent formulations of the axiom of countable real choice are established. One of them is the following.

\begin{prop}[AC and sequential continuity]\label{prop_LindeN}
In {\em ZF} claims~1 and~2 are equivalent.
\begin{enumerate}
\item A~function $f\cc\R\to\R$ is continuous at a~point $x$ iff $f$ is sequentially continuous at $x$.
\item The axiom of countable choice holds for subsets of $\R$. 
\end{enumerate}
\end{prop}
\duk
See \cite{herr_stre}.
\kduk

\noindent
But already in 1918 the Polish mathematician 
{\em Wac\l aw Sierpi\'nski (1882--1969)\index{Sierpi\'nski, Wac\l aw}} circumvented in \cite{sier18} the axiom of 
choice for global continuity on open real sets. We present his solution in Theorem~\ref{thm_sierp}.

\medskip\noindent
{\em $\bullet$ Choice from sets of fractions. }The existence of a~selector on the set
$$
\mathcal{P}_0(\Q)\equiv\mathcal{P}(\Q)\setminus\{\emptyset\}
$$ 
is a~theorem in ZF. In the next proposition we state a~more convenient choice principle.

\begin{exer}\label{ex_wellOrdQ}
Show that in {\em ZF} the set $\Q$ of fractions has a~well ordering.    
\end{exer}

\begin{prop}[rational choice]\label{prop_choiInQ}
We have the theorem
$$
\mathrm{ZF}\vdash\forall F\cc\N\to\mathcal{P}_0(\Q)\,\exists f\cc\N\to\Q\,\big(n\in\N\Rightarrow
f(n)\in F(n)\big)\,.
$$
\end{prop}
\duk
Let a~map $F\cc\N\to\mathcal{P}_0(\Q)$ be given. Using Exercise~\ref{ex_wellOrdQ} we take, in ZF, a~well ordering
$$
(\Q,\,\prec)\,.
$$
Then we define $f\cc\N\to\Q$ by ($n\in\N$)
$$
f(n)\equiv\min_{\prec}(F(n))\,.
$$
This is a~definition within ZF.
\kduk

\noindent
{\em $\bullet$ Sierpi\'nski's theorem. }Sierpi\'nski's method for 
eliminating AC in the proof of equivalence of continuity with 
sequential continuity uses the density of $\Q$ in any open real set.

\begin{thm}[W.~Sierpi\'nski, 1918]\label{thm_sierp}
Let $M\sus\R$ be an open set.
In\index{theorem!Sierpi\'nski's on sequential continuity|emph} 
{\em ZF}, i.e. without the axiom of choice, we have the theorem
$$
\mathrm{ZF}\vdash\forall f\in\mathcal{F}(M)\,\big(\text{$f$ is 
continuous}\iff\text{$f$ is sequentially continuous}\big)\,.
$$
\end{thm}
\duk
Let $f\cc M\to\R$ where $M\sus\R$ is an open set. The implication $\Rightarrow$ 
is easy to establish in ZF. We assume that $f$ is continuous, 
$b\in M$ and that $(a_n)\sus M$ is a~sequence with $a_n\to b$. Let an 
$\ep$ be given. We take a~$\de$ such that
$$
f[U(b,\,\de)]\sus U(f(b),\,\ep)\,.
$$
Then there is an $n_0$ such that $a_n\in U(b,\de)$ for every $n\ge 
n_0$. It follows that $f(a_n)\in U(f(b),\ep)$ for the same $n$. Hence $f(a_n)\to f(b)$ and $f$ is sequentially continuous at $b$.

We prove, in ZF, the implication $\neg\Rightarrow\neg$. This is the 
main result of the theorem. We assume that $f$ is not continuous. 
Thus there exist a~point $b\in M$ and an $\ep$ such that for every $n\in\N$ there exists a~point
$${\textstyle
\text{$b_n\in U(b,\,\frac{1}{n})\cap M$ such that $|f(b_n)-f(b)|>\ep$}\,.}
$$
We consider for $n\in\N$ the sets of fractions
$${\textstyle
A_n\equiv\{\al\in U(b,\,\frac{1}{n})\cap M\cap\Q\cc\;|f(\al)-f(b)|>\frac{\ep}{2}\}\,.
}
$$
If $A_n\ne\emptyset$ for every $n$, we use 
Proposition~\ref{prop_choiInQ} and obtain in ZF a~function
$$
\varphi\cc\N\to\Q
$$
such that $\varphi(n)\in A_n$ for every $n\in\N$. Then $(a_n)\sus\Q\cap M$ 
with $a_n\equiv \varphi(n)$ is a~sequence such that $a_n\to b$ but
$f(a_n)\not\to f(b)$, and $f$ is not sequentially continuous at the 
point $b$.

The complementary case is that $A_m=\emptyset$ for some $m\in\N$. 
This by the triangle inequality  means that
$$
{\textstyle
\text{$|f(\al)-f(b_m)|>\frac{\ep}{2}$ for every $\al\in U(b,\,\frac{1}{m})\cap M\cap\Q$}}
$$
(Exercise~\ref{ex_maleCvic}). We take large $N\in\N$ such that 
$$
{\textstyle
U(b_m,\,\frac{1}{N})\sus U(b,\,\frac{1}{m})\cap M
}
$$
and define for $n\in\N$ the sets of fractions
$${\textstyle
B_n\equiv U(b,\,\frac{1}{N+n-1})\cap\Q\,.
}
$$
Then $|f(\al)-b_m|\ge\frac{\ep}{2}$ for every $n\in\N$ and every fraction $\al\in B_n$. Using again 
Proposition~\ref{prop_choiInQ} we obtain in ZF a~function
$$
\psi\cc\N\to\Q
$$
such that $\psi(n)\in B_n$ for every $n\in\N$. Then $(a_n)\sus\Q\cap M$ 
with $a_n\equiv \psi(n)$ is a~sequence such that $a_n\to b_m$ 
but $f(a_n)\not\to f(b_m)$, and $f$ is not sequentially continuous at 
the point $b_m$.
\kduk
\vspace{-3mm}
\begin{exer}\label{ex_maleCvic}
Justify the inequality in the proof.   
\end{exer}
One might think that in the 
equivalence of both continuities the axiom of choice could be circumvented altogether by well ordering the real numbers, i.e., by replacing in Exercise~\ref{ex_wellOrdQ} the set 
$\Q$ with the set $\R$. This is impossible because $\R$ much differs 
from $\Q$. In 1970 the American mathematician {\em Robert 
M.~Solovay (1938)\index{Solovay, Robert M.}} proved in \cite{solo} that 
\begin{quote}
if the theory $\mathrm{ZF}+\mathrm{AC}+
\mathrm{I}$ is consistent, then so is the theory $\mathrm{ZF}+\mathrm{LM}$\,.    
\end{quote}
A~theory is just a~set of formulas. ZF are the axioms of ZF set theory 
(Section~...), AC is the axiom of choice 
(Axioms~\ref{axio_AC} and ...), I is the axiom 
that there exists an inaccessible cardinal 
(Definition~...) and LM is the axiom that every set of real numbers is Lebesgue measurable.
Consistency of a~theory means that one cannot derive 
from the formulas of the theory a~contradiction; else the theory is 
contradictory. Thus unless the theory 
$\mathrm{ZF}+\mathrm{AC}+
\mathrm{I}$ is contradictory we cannot prove in ZF that $\R$ has 
a~well ordering; a~well ordering of 
$\R$ yields via a~well known construction a~set of real numbers 
that is not Lebesgue measurable. 
In 1984 the Israeli mathematician {\em Saharon Shelah 
(1946)\index{Shelah, Saharon}} proved in \cite{shel} that the 
axiom I cannot be removed from Solovay's result. 

On the other hand, already in 1938/1940 K.~G\"odel published in \cite{gode38,gode40}, see also \cite[Chapter VI. ``Jetzt, Mengenlehre'']{daws}, the result that
\begin{quote}
if the theory $\mathrm{ZF}$ is consistent, then so is the theory $\mathrm{ZF}+\mathrm{AC}$\,.    
\end{quote}
Thus unless ZF is contradictory, in view of Zermelo's 
Theorem~... we cannot prove in ZF that a~well ordering of 
$\R$ does not exist.

Is it mind-boggling? Yes, it is.
It is the up-to-date answer to the question when for a~function 
$f\in\mathcal{F}(M)$ continuity is equivalent with sequential 
continuity.

\section[The cardinality of continuous functions]{The cardinality of continuous functions}\label{sec_poceSpoFun}

We show that there exists a~bijection between the sets $\mathcal{C}(\R)$ and $\R$. 

\medskip\noindent
{\em $\bullet$ The Cantor--Bernstein theorem. }We use this theorem, which we 
prove in {\em MA~1${}^+$}, to construct the mentioned bijection.

\begin{thm}[Cantor--Bernstein]\label{thm_CBv2}
Let\index{theorem!Cantor--Bernstein} 
$X$ and $Y$ be sets. If
there exist two injections,
from $X$ to $Y$ and from $Y$ to $X$, \underline{then} there exists 
a~bijection from $X$ to $Y$.
\end{thm}

\begin{exer}\label{ex_and_back}
One can extend the conclusion: there exist two bijections, from $X$ to $Y$ and from $Y$ to $X$.     
\end{exer}
We mentioned G.~Cantor earlier. {\em Felix Bernstein\index{Bernstein, 
Felix} (1878--1956)} was a~German mathematician. For example, 
$$
(m,\,n)\mapsto 2^m3^n
$$ 
is an injection from $\N\times\N$ to $\N$ and 
$$
n\mapsto(1,\,n)
$$ 
is an injection
from $\N$ to $\N\times\N$. Hence by the C.--B. theorem there exists a~bijection from $\N\times\N$ to $\N$.
The next exercise describes such 
bijection. 

\begin{exer}\label{ex_bijs}
The function $b\cc\N\times\N\to \N$, where $b(m,n)=(2m-1)\cdot 2^{n-1}$, is a~bijection. 
\end{exer}

\noindent
{\em $\bullet$ How many continuous functions from $\R$ to $\R$ are there? }As many as real numbers. To prove
that there exists a~bijection between $\mathcal{C}(\R)$ and $\R$, by
Theorem~\ref{thm_CBv2} it suffices to produce two injections, from $\R$ to $\mathcal{C}(\R)$ and from
$\mathcal{C}(\R)$ to $\R$. The former injection is clear, see
Exercise~\ref{ex_jeJasna}. We describe an injection from
$\mathcal{C}(\R)$ to $\R$. We begin
with an injection from $\R^{\N}$, the set of real sequences, to $\R$.

\begin{prop}[coding $(a_n)$ as $a$]\label{prop_injRNtoR}
There exists an injection $f\cc\R^{\N}\to\R$.    
\end{prop}
\duk
Let $X$ be the twelve-element set of ten digits $0$, $1$, $\ds$, $9$, the 
decimal point\ .\ and the minus sign $-$. We view any real number $b$ as 
a~sequence $(b_m)\sus X$, $m\in\N$. For example, 
$$
-\pi=(-,\,3,\,.,\,1,\,4,\,1,\,5,\,\ds)
$$ 
or this year is 
$$
(2,\,0,\,2,\,5,\,.,\,0,\,0,\,0,\,\ds)\,. 
$$
We fix a~coding $c$ of the elements in
$X$ by pairs of digits: $c(0)=00$, $c(1)=01$, $\ds$, $c(9)=09$, $c(.)=10$
and $c(-)=11$. We use the inverse bijection 
$$
b^{-1}=(u,\,v)\cc\N\to\N\times\N
$$ 
of Exercise~\ref{ex_bijs} and define the desired injection $f$ by
$$
f((a_n))=f(((a_{n,\,m})))
\equiv0\,.\,d_1\,d_2\,d_3\,d_4\,\ds\, 
d_{2l-1}\,d_{2l}\,\ds\ \ (\in\R)\,,
$$
where $d_i\in\{0,1,\ds,9\}$ and $d_{2l-1}d_{2l}\equiv c(a_{u(l),v(l)})$. It is clear that $f$ is 1-1.
\kduk

Now it is easy to get an injection from $\mathcal{C}(\R)$ to $\R$. 

\begin{thm}[the cardinality of $\mathcal{C}(\R)$]\label{thm_poceSpoFun}
The set $\mathcal{C}(\R)$\index{theorem!number of continuous 
functions|emph} of continuous real functions with domain $\R$ is in 
bijection with $\R$.
\end{thm}
\duk
It suffices to find an injection 
$g\cc\mathcal{C}(\R)\to\R^{\N}$ because then the composition $f(g)$ with the injection $f$ in the previous proposition is an injection from 
$\mathcal{C}(\R)$ to $\R$. This is easy, we take any bijection 
$h\cc\N\to\Q$ (Theorem) and define 
the value of $g$ on $j=j(x)\in\mathcal{C}(\R)$ as the 
real sequence 
$$
g(j)\equiv(a_n)\,\text{ where }\,a_n\equiv j(h(n))\,.
$$
By Exercise~\ref{ex_oveJsouHust} and 
Proposition~\ref{prop_oHusteMn} one can uniquely recover $j$ 
from $(a_n)$. Hence $g$ is injective.
\kduk
\vspace{-3mm}
\begin{exer}\label{ex_jeJasna}
Show that the map $a\mapsto k_a$, where $k_a$ is the constant function, is an 
injection from $\R$ to $\mathcal{C}(\R)$.  
\end{exer}

\begin{exer}\label{ex_dokazZoecn}
For every nonempty set $M\sus\R$ the sets $\mathcal{C}(M)$ and $\R$ are in bijection.     
\end{exer}

\section[Attaining intermediate values]{Attaining intermediate values}\label{sec_nabyvMezih}

We show that continuous functions map intervals to intervals. 

\medskip\noindent
{\em $\bullet$ Any continuous function attains any intermediate value. }We prove it in the next theorem. 

\begin{thm}[intermediate values]\label{thm_mezihodnoty}
Let\index{theorem!intermediate values attained|emph} $a<b$, 
$f\in\mathcal{C}([a,b])$ and $f(a)<c<f(b)$ or $f(a)>c>f(b)$. \underline{Then} 
$$
\text{$c=f(d)$ for some $d\in(a,b)$}\,.
$$
\end{thm}
\duk
Let $f(a)<c<f(b)$, the other case is similar. We set 
$$
X\equiv\{x\in[a,b]\cc\;f(x)<c\}\,\text{ and }\,d\equiv\sup(X)\ \ (\in[a,\,b])\,.
$$
The continuity of $f$ implies that
$d\in(a,b)$. We show that both $f(d)<c$ and $f(d)>c$ lead to contradiction, and hence $f(d)=c$. 
Let $f(d)<c$. By the continuity of $f$ at $d$ there is a~$\de$ such that 
for every $x\in U(d,\de)\cap[a,b]$ it holds that $f(x)<c$. But then $X$ 
contains numbers larger than $d$, which is a~contradiction. Let 
$f(d)>c$. By the continuity of $f$ at $d$ there is a~$\de$ such that 
for every $x\in U(d,\de)\cap[a,b]$ it holds that $f(x)>c$. But then every 
$x<d$ and close to $d$ lies outside $X$, which is also a~contradiction.
\kduk
\vspace{-3mm}
\begin{exer}\label{ex_obrJeInt}
For every interval $I\sus\R$ and every function $f\in\mathcal{C}(I)$ 
the image $f[I]$ is an interval.    
\end{exer}

\begin{cor}[the image of $\exp x$]\label{cor_imageExp}
We have
$$
\exp[\R]=(0,\,+\infty)\,.
$$ 
Thus the function $\exp$ is a~bijection from $\R$ to $(0,+\infty)$.
\end{cor}
\duk
Since $\exp>0$ on $\R$, we have $\exp[\R]\sus(0,+\infty)$. 
From the limits $\lim_{x\to-\infty}\exp x=0$ and 
$\lim_{x\to+\infty}\exp x=+\infty$ (part~3 of Proposition~\ref{prop_expFce}), from 
the continuity of $\exp x$ (Corollary~\ref{cor_spoExpKosSin}) 
and from Theorem~\ref{thm_mezihodnoty} it follows that 
$(0,+\infty)\sus\exp[\R]$. Thus $\exp[\R]=(0,+\infty)$. The 
exponential increases and hence it is injective and a~bijection.
\kduk
\vspace{-3mm}
\begin{exer}\label{ex_alpinista}    
Prove the following corollary. 
\end{exer}

\begin{cor}[mountaineering]\label{ex_horolezec}
A~mountain climber starts her ascend at midnight, after 24 hours she
reaches the summit and then she descends for 24 hours to the base camp. Show that 
there exists a~moment 
$$
t_0\in[0,\,24]
$$ 
when she is in both days in exactly the same altitude. 
\end{cor}

A~function $f\in\mathcal{F}(M)$ \underline{increases\index{function!increases|emph}}, respectively 
\underline{decreases\index{function!decreases|emph}}, on an 
(arbitrary) set $X$ 
if for every $x<y$ in $M\cap X$ we have  
$$
\text{$f(x)<f(y)$, respectively $f(x)>f(y)$}\,.
$$

\begin{cor}[continuity and injectivity]\label{cor_spojInjInt}
If $I\sus\R$ is an interval and $f\in\mathcal{C}(I)$ is injective, 
\underline{then} $f$ either increases or decreases on $I$.
\end{cor}
\duk
If $f$ neither increases nor decreases then $I$ contains three numbers $a<b<c$ such that $f(a)
<f(b)>f(c)$ or $f(a)>f(b)<f(c)$. In the first case, we see by
Theorem~\ref{thm_mezihodnoty} that for every $d$ such that 
$$
f(a),\,f(c)<d<f(b)
$$ 
there exist numbers $x\in(a,b)$ and a~$y\in(b,c)$ such that
$$
d=f(x)=f(y)\,, 
$$
which contradicts the injectivity of $f$. The second case 
leads to a~similar contradiction.
\kduk

Now we can easily prove 
state here as a~corollary.

\begin{cor}[Bolzano--Cauchy]\label{cor_BolzCauch}
Suppose that $a<b$, $f\in\mathcal{C}([a,b])$ and that $f(a)f(b)\le0$. \underline{Then} 
$$
Z(f)\ne\emptyset\,.
$$
\end{cor}
\duk
If $f(a)f(b)=0$ then $f(a)=0$ or $f(b)=0$. If $f(a)f(b)<0$ then $f(a)<0<f(b)$ or $f(b)<0<f(a)$, and $0=f(c)$ for some $c\in(a,b)$ by 
Theorem~\ref{thm_mezihodnoty}.
\kduk

\section[Compact, open and closed sets]{Compact, open and closed sets}\label{sec_kompaktnost}

In mathematical analysis, compactness is a~fundamental concept. We consider only compact sets in $\R$.  

\medskip\noindent
{\em $\bullet$ Real compact sets. }One possible definition is as follows.

\begin{defi}[compact sets]\label{def_compSets}
A~set $M\sus\R$ is \underline{compact\index{set!compact|emph}\index{compactness|emph}} if
every sequence $(a_n)\sus M$ has a~convergent subsequence $(a_{m_n})$ 
with 
$$
\lim a_{m_n}\in M\,.
$$
\end{defi}
An equivalent definition, which generalizes more easily, uses covers by open sets; see Theorem~\ref{thm_HeineBorelVR} below
(no proof). We introduce open sets later.

\begin{exer}\label{ex_abIsComp}
Prove that every interval $[a,b]$ is 
a~compact set.
\end{exer}

\begin{exer}\label{ex_emptyIsComp}
Is the empty set compact? Is $\R$ compact?
\end{exer}

Let $(X,\prec)$ be a~linear order and $A\sus X$. Recall that $\min(A)$, the minimum
of $A$ (relative to $\prec$), is the unique
element $a\in A$, if it exists, such that 
$$
\text{$b\succeq a$ for every $b\in A$}\,.
$$
For $\preceq$
in place of $\succeq$ we get the maximum
$\max(A)$. Statement of the next proposition employs the usual linear order 
$(\R,<)$. In the proof we use the linear order $(\R^*,<)$.

\begin{prop}[extremities of compacts]\label{prop_minMaxComp}
Any nonempty and compact set  $M\sus\R$ 
$$
\text{has both $\min(M)$ and $\max(M)$}\,.
$$
\end{prop}
\duk
We only prove that the maximum of $M$ exists, the proof for minimum is similar. 
Let $M\sus\R$, $M\ne\emptyset$, be compact and let
$$
A\equiv\sup(M)\,,
$$
in the linear order $(\R^*,<)$ (see 
$A=+\infty$ then there is a~sequence $(a_n)\sus M$ with $\lim a_n=+\infty$. But such 
$(a_n)$ has no convergent subsequence, which contradicts the compactness of $M$. Thus 
$A\in\R$. We show that $A\in M$, and then $A=\max(M)$. If $A\not\in M$, then we again take a~sequence $(a_n)\sus M$ with 
$\lim a_n=A$. It again has no convergent subsequence with limit in $M$ (the limit of every subsequence of $(a_n)$ 
is $A$), which  again contradicts the compactness of $M$.
\kduk

\noindent
{\em $\bullet$ Continuous images of compacts are compact. }This is an important interplay
between compactness and continuity. 

\begin{thm}[images of compacts]\label{thm_obrKomp}
If\index{theorem!images of compacts|emph} $f\in\mathcal{C}(M)$ and $M$ is compact \underline{then} the image 
$$
f[M]\ \  (\sus\R) 
$$
is a~compact set.    
\end{thm}
\duk
Let $f$ and $M$ be as stated and $(b_n)\sus f[M]$. Using the axiom of 
choice\index{axiom!of choice, AC} we take a~sequence $(a_n)\sus M$ such that 
$f(a_n)=b_n$. It has a~subsequence $(a_{m_n})$ with the limit 
$$
a\equiv\lim a_{m_n}\ \ (\in M)\,. 
$$
By (H) we have $\lim f(a_{m_n})=f(a)\equiv b$, so that $(b_{m_n})=(f(a_{m_n}))$ has the
limit $b\in f[M]$.
Hence $f[M]$ is compact.
\kduk

We show that continuous functions defined on compact sets
attain extremal values. If
$f\in\mathcal{R}$ and $b\in M(f)$ is such that 
$$
\text{$f(b)\ge f(x)$ 
for every $x\in M(f)$}\,, 
$$
we say that $f$ \underline{attains (at $b$) 
maximum\index{function!attains maximum|emph}}. If the opposite inequality $\le$ holds
for every $x\in M(f)$, then $f$ \underline{attains (at $b$)
minimum\index{function!attains minimum|emph}}.

\begin{cor}[minima and maxima]\label{cor_priMaxi}
Let $M\sus\R$ be a~compact set and let
$f\in\mathcal{C}(M)$. \underline{Then} 
the function $f$ attains both minimum and maximum. 
\end{cor}
\duk
Let $M$ and $f$ be as stated. By Proposition~\ref{prop_minMaxComp} 
and Theorem~\ref{thm_obrKomp} both $\min(f[M])$ and $\max(f[M])$ exist. 
Hence $f$ attains both minimum and maximum.
\kduk
\vspace{-3mm}
\begin{exer}\label{ex_neextremy}
The continuous functions
$\mathrm{id}\,|\,[0,1)$ and 
$1/(1-x)\,|\,[0,1)$ do not attain maximum.     
\end{exer}

\noindent
{\em $\bullet$ Local, global and strict extremes. }Let 
$f\in\mathcal{R}$ and $b\in M(f)$. We already defined what it means that
$f$ attains at $b$ maximum, resp. minimum. We also say that $f$ has at $b$ a~\underline{global 
maximum\index{maximum of a function!global|emph}}, resp. 
a~\underline{global minimum\index{minimum of a 
function!global|emph}}. If for every $x\in M\setminus\{b\}$ the inequality holds 
as strict, we speak of a~\underline{strict global maximum\index{maximum of a 
function!strict|emph}}, resp. 
a~\underline{strict global minimum\index{minimum of a 
function!strict|emph}}. These
(strict) minima or maxima at $b$ 
are \underline{local\index{maximum of a function!local|emph}\index{minimum of a function!local|emph}} if for some $\de$ the restriction $f\,|\,U(b,\de)$
attains them as global. 

\begin{exer}\label{ex_onExtr}
Let $f\in\mathcal{R}$. 
Prove the equivalence: $f$ attains maximum (or minimum) at every point $b\in M(f)$ 
iff $f=k_c\,|\,M(f)$ for some $c\in\R$.
\end{exer}

\noindent
{\em $\bullet$ Open and closed sets. }We introduce these important families of real sets. 

\begin{defi}\label{def_opClSets}
A~set $M\sus\R$ is 
\underline{open\index{set!open|emph}} if for every point $b\in M$ there exists
a~$\de$ such that $U(b,\de)\sus M$. The set $M$ is 
\underline{closed\index{set!closed|emph}} if its complement $\R\setminus M$ 
is open.    
\end{defi}

\begin{prop}[properties of open sets]\label{prop_openSet}
Let $\{U_i\cc\;i\in I\}$, $I\ne\emptyset$, be a~set system of open sets. \underline{Then} the following holds.
\begin{enumerate}
\item Both the empty set $\emptyset$ and the whole set $\R$ are open sets.
\item The union $\bigcup_{i\in I}U_i$
is an open set.
\item For any finite index set $I$ the intersection $\bigcap_{i\in I}U_i$ is an open set.
\end{enumerate}
\end{prop}
\duk
1. This is trivial: there is no $a\in\emptyset$ and $U(a,\de)\sus\R$ for every $a$ and $\de$. 

2. If $a\in\bigcup_{i\in I}U_i$ then $a\in U_i$ for some $i\in I$ and there is a~$\de$ such that $U(a,\de)\sus U_i$, so that $U(a,\de)\sus\bigcup_{i\in I}U_i$. 

3. If $a\in\bigcap_{i\in I}U_i$ then $a\in U_i$ for every $i\in I$ and
for every $i\in I$ there is a~$\de_i$ such that $U(a,\de_i)\sus 
U_i$. We set 
$$
\de\equiv\min(\{\de_i\cc\;i\in I\})\,.
$$
This minimum exists and is positive because $I$ is finite. Thus $U(a,\de)\sus U_i$ for every $i\in I$. 
\kduk
\vspace{-3mm}
\begin{exer}\label{ex_vlUzMno}
State and prove analogous properties of closed sets.   
\end{exer}

\begin{exer}\label{ex_opAndLimP}
If $U\sus\R$ is an open set then every element of $U$ is a~limit point of $U$.   
\end{exer}

\begin{prop}[on closed sets]\label{prop_ekvivDefUzav}
A~set $M\sus\R$ is closed $\iff$ the limit of every convergent sequence $(a_n)\sus M$ lies in $M$.
\end{prop}
\duk
Implication $\Rightarrow$.
Let $M\sus\R$ be a~closed set and let $(a_n)\sus M$ be a~sequence with 
$$
\lim a_n=a\in\R\setminus  M\,. 
$$
But then for some $\de$ we have $U(a,\de)\cap 
M=\emptyset$. This is not possible because $a_n\to a$ and $a_n\in U(a,\de)$ for every $n\ge n_0$. Hence $a\in M$. 

Reverse implication $\neg\Rightarrow\neg$. Suppose that 
$M\sus\R$ is not a~closed set. Then there is a~point 
$a\in\R\setminus M$ such that for every $n$ there exists 
$${\textstyle
a_n\in U(a,\frac{1}{n})\cap M\,. 
}
$$
Thus $(a_n)\sus M$ and $\lim 
a_n=a\not\in M$. We used the axiom of choice\index{axiom!of choice, AC}.
\kduk

\noindent
{\em $\bullet$ Characterizations of compact sets. }A~set $M\sus\R$ is 
\underline{bounded\index{set!bounded|emph}} if $M\sus[a,b]$ for some real $a<b$. We characterize real compact sets.

\begin{thm}[real compact sets]\label{thm_KompvR}
Let $M\sus\R$.\index{theorem!real compact sets|emph} \underline{Then}
$$
\text{$M$ is compact $\iff$ $M$ is bounded and closed}\,.
$$
\end{thm}
\duk
Implication $\Leftarrow$.
Let $M$ be bounded and closed and $(a_n)\sus M$. By 
Theorem~\ref{thm_BolzWeier} we have a~convergent subsequence $(a_{m_n})$.
$M$ is closed and Proposition~\ref{prop_ekvivDefUzav} 
shows that $\lim a_{m_n}\in M$. Hence $M$ is compact. 

Contrapositive implication $\neg\Leftarrow\neg$. 
Suppose that $M$ is not bounded. We then define a~sequence $(a_n)\sus M$ such 
that 
$$
m\ne n\Rightarrow|a_m-a_n|\ge1 
$$
---\,such sequence has no 
convergent subsequence. First term $a_1$ is arbitrary. 
Suppose that $a_1,a_2,\ds,a_n$ are defined
and $|a_i-a_j|\ge1$ if $i\ne j$. Since $M$ is not  bounded, there exists 
$$
\text{$a_{n+1}\in M$ such that 
$|a_{n+1}|\ge1+\max(\{|a_1|,\,|a_2|,\,\ds,\,|a_n|\})$}\,. 
$$
By the $\Delta$-inequality we have $|a_{n+1}-a_i|\ge1$ for every $i\in[n]$.
Repeating it indefinitely we get the required sequence $(a_n)$.

Suppose that $M$ is not closed. By 
Proposition~\ref{prop_ekvivDefUzav} there is a~sequence $(a_n)\sus M$ 
with limit in $\R\setminus M$. Every subsequence of $(a_n)$ has 
the same limit and therefore does not converge in $M$. 
\kduk
\vspace{-3mm}
\begin{exer}\label{ex_vyhoPrsOko}
$[a,b]\setminus P(c,\de)$ is a~compact set.   
\end{exer}

\begin{exer}\label{ex_inMS}
Show that the implication $\Rightarrow$ holds more generally in metric spaces, every compact set in a~metric space is bounded and closed.  
\end{exer}

\begin{exer}\label{ex_CountEinMS}
However, the implication $\Leftarrow$, in general, does not hold in metric
spaces. The hint is to consider discrete spaces.  
\end{exer}

We do not prove the following theorem. 

\begin{thm}[Heine--Borel]\label{thm_HeineBorelVR}
Let $M\sus\R$. \underline{Then} $M$ is compact $\iff$ for every system
$\{U_i\cc\;i\in I\}$ of open
sets the implication holds that 
$${\textstyle
\bigcup_{i\in I}U_i\supset M \Rightarrow \text{there is a~finite set $J\sus I$ such that $\bigcup_{i\in J}U_i\supset 
M$}\,.}
$$
\end{thm}
In {\em MA~1${}^+$} we prove a~more general version of this theorem for metric spaces.

\medskip\noindent
{\em $\bullet$ More on open and closed sets. }Let $A\sus M\sus\R$. Then $A$ is \underline{relatively 
closed\index{set!relatively closed|emph}} in $M$ if there is 
a~closed set $U\sus\R$ such that $A=M\cap U$.

\begin{prop}[zero sets are relatively closed]\label{prop_zeroClos}
Let $f\in\mathcal{C}$. \underline{Then} the set 
$$
Z(f)=\{x\in M(f)\cc\;f(x)=0\}
$$ 
is relatively closed in $M(f)$. 
\end{prop}
\duk
Let $N\equiv M(f)\setminus Z(f)$. Using the continuity of $f$ we see that 
$$
\text{for every $b\in 
N$ there is a~$\de_b$ such that $U(b,\,\de_b)\cap Z(f)=\emptyset$}\,. 
$$
Let
$A\equiv\bigcup_{b\in N}U(b,\de_b)$. By part~2 of Proposition~\ref{prop_openSet} the set $A$ is open. Thus $\R\setminus A$ is closed and $Z(f)=M(f)\cap(\R\setminus A)$.
\kduk

We show that continuous injective images of open sets are open. 

\begin{prop}[images of open sets]\label{prop_Obropen}
Suppose that $f\in\mathcal{C}$ is injective and both sets $M(f)$ and 
$M\sus\R$ are open. \underline{Then} the image $f[M]$ is open.
\end{prop}
\duk
Let $b\in f[M]$ and $a\equiv 
f^{-1}(b)$ ($\in M(f)\cap M$). Since $M(f)$ and $M$ are open, so is their intersection and we can take
an interval $I\equiv[a-\de,a+\de]$ such that $I\sus M(f)\cap M$. Then 
$$
f(a-\de)<b<f(a+\de)\,\text{ or }\,f(a-\de)>b>f(a+\de)
$$ 
because $f(a-
\de),f(a+\de)<b=f(a)$ and $f(a-\de),f(a+\de)>b=f(a)$ contradict by
Theorem~\ref{thm_mezihodnoty} the injectivity of $f$. We take any
$$
\ep<\min(|f(a+\de)-b|,\,|f(a-\de)-b|)\,.
$$
Using Theorem~\ref{thm_mezihodnoty} we have $U(b,\ep)\sus f[I]\sus f[M]$. Hence 
$f[M]$ is open.
\kduk

Let $A\sus M\sus\R$. Then $A$ is \underline{relatively 
open\index{set!relatively open|emph}} in $M$ if there is
an open set $B\sus\R$ such that $A=M\cap B$.

\begin{prop}[preimages of open sets]\label{prop_preiOpen}
Let $f$ be in $\mathcal{C}$ and $M\sus\R$ be open. \underline{Then} the set
$$
\text{$f^{-1}[M]=\{x\in M(f)\cc\;f(x)\in M\}$ is relatively open in $M(f)$}\,.
$$
\end{prop}
\duk
For every $b\in M$ there is an $\ep_b$ such that $U(b,\ep_b)\sus M$. For every 
$a$ in $f^{-1}[M]$ ($\sus M(f)$) there is a~$\de_a$ such that for $b\equiv f(a)$ we 
have 
$$
f[U(a,\,\de_a)]\sus U(b,\,\ep_b)\sus M\,.
$$
Thus, $U(a,\de_a)\cap M(f)\sus f^{-1}[M]$. Let 
$${\textstyle
B\equiv\bigcup_{a\in f^{-1}[M]}U(a,\,\de_a)\,. 
}
$$
By part~2 of Proposition~\ref{prop_openSet} this is an open set.
Since $f^{-1}[M]=M(f)\cap B$, the set $f^{-1}[M]$ is relatively open in $M(f)$.
\kduk

\noindent
We need Propositions~\ref{prop_zeroClos} and \ref{prop_preiOpen} in the proof of Theorem~\ref{thm_deriSEF}.

We describe the structure of any open set. An \underline{open interval\index{interval!open|emph}} 
is any nonempty interval that is an open set. These are exactly the intervals ($a<b$) 
$$
(-\infty,\,a),\  
(a,\,+\infty)\,\text{ and }\,(a,\,b)\,.
$$

\begin{thm}[structure of open sets]\label{thm_struOtMno}
A~real set\index{theorem!structure of open sets|emph} is open $\iff$ it is 
a~union of 
at most countably many mutually disjoint open intervals.
\end{thm}
\duk
Implication $\Leftarrow$ follows from part~2 of Proposition~\ref{prop_openSet}.

Implication $\Rightarrow$.
Let $M\sus\R$ be an open set. If $M=\emptyset$, the claim holds
trivially. Let $a\in M$. We define $I_a$ as the inclusion-wise maximal open interval 
$I$ such that $a\in I\sus M$ (Exercise~\ref{ex_moreDexi}). By 
Exercise~\ref{ex_explWhy} we have for any $a,b\in M$ that 
$$
I_a=I_b\,\text{ or }\,I_a\cap I_b=\emptyset\,. 
$$
Thus 
$S\equiv\{I_a\cc\;a\in\Q\cap M\}$ is an at most countable set system of mutually disjoint
open intervals such that $\bigcup S=M$ (Exercise~\ref{ex_whyUnion}).
\kduk
\vspace{-3mm}
\begin{exer}\label{ex_moreDexi}
Show that the interval $I_a$ exists. \end{exer}

\begin{exer}\label{ex_explWhy}
Explain why $I_a=I_b$ or $I_a\cap I_b=\emptyset$.
\end{exer}

\begin{exer}\label{ex_whyUnion}
Why $\bigcup S=M$?     
\end{exer}

\noindent
{\em $\bullet$ The Cantor set. }By the previous theorem any closed set $\R\setminus M$ is a~union of the ``gaps'' separating the intervals 
$I_a$. 
For $|S|=n$ we have at most $n+1$ gaps. It is not easy to imagine that for countable $S$ the set of gaps may be uncountable. Such closed sets are hard to see by our inner eye. An example is the  
\underline{Cantor set\index{Cantor set|emph}} $C$:
$${\textstyle
C\equiv\bigcap_{n=1}^{\infty}C_n\ (\sus[0,\,1]\equiv C_0)\,\text{ where }\,C_n\equiv\frac{1}{3}C_{n-1}\cup\big(\frac{1}{3}C_{n-1}+\frac{2}{3}\big)\,.\label{Cantor}
}
$$
Here we denote by $aM+b$ the set 
$\{ax+b\cc\;x\in M\}$. Thus $C$ is the leftover of the interval $[0,1]$ obtained 
 deleting first the open middle third 
$(\frac{1}{3},\frac{2}{3})$, then by deleting from the rest 
$${\textstyle
[0,\,1]\setminus(\frac{1}{3},\,\frac{2}{3})=
[0,\,\frac{1}{3}]\cup[\frac{2}{3},\,1]
}
$$ 
the open middle thirds 
$(\frac{1}{9},\frac{2}{9})$ and $(\frac{7}{9},\frac{8}{9})$, and so on.

\begin{exer}\label{ex_CabtDis}
Prove the following properties of $C$. 
\begin{enumerate}
\item $C$ is an uncountable set.
\item $C$ is closed.
\item For every $\ep$ there exist $k$ intervals $[a_i,b_i]$, $a_i<b_i$
and $i\in[k]$, such that 
$${\textstyle
\sum_{i=1}^k(b_i-a_i)\le\ep\,\text{ and }\,\bigcup_{i=1}^k[a_i,b_i]\supset C\,.
}
$$
In this sense $C$ has zero length.
\end{enumerate}  
\end{exer}

\noindent
{\em $\bullet$ Baire's theorem. }The following theorem is due to the French 
mathematician {\em Ren\'e-Louis Baire (1874--1932)\index{Baire, Ren\'e-Louis}}. In the proof we employ \underline{closed neighborhoods\index{neighborhood!closed|emph}}
$$
\overline{U}(b,\,\ep)\equiv[b-\ep,\,b+\ep]\,.\label{closedNeig}
$$
These are closed sets.

\begin{thm}[Baire's]\label{thm_Baire}
Suppose\index{theorem!Baire's|emph} 
that
$$
M=\bigcup_{n=1}^{\infty}M_n
$$
where $M\sus\R$ is a~nonempty closed set.
\underline{Then} there exists an index $n$ such that the set $M_n$ is not sparse in $M$.
\end{thm}
\duk
Suppose that $M$ is as stated and, for contradiction, that every set $M_n$ 
is sparse. We take an arbitrary point $b_0\in M$. Since 
$M_1$ is sparse, there is a~point $b_1\in M$ and $\de_1$ such that 
$\overline{U}(b_1,\de_1)\sus U(b_0,1)$, hence $\de_1\le1$, and 
$$
M_1\cap\overline{U}(b_1,\,\de_1)=\emptyset\,.
$$
Suppose that we already defined nested closed neighborhoods
$$
\overline{U}(b_1,\,\de_1)\supset
\overline{U}(b_2,\,\de_2)\supset\ds
\supset\overline{U}(b_n,\,\de_n)
$$
such that $b_i\in M$, $\de_i\le\frac{1}{i}$ and 
$M_i\cap \overline{U}(b_i,\de_i)=\emptyset$ for every $i\in[n]$.
Since $M_{n+1}$ is sparse, we can take a~point $b_{n+1}\in M$ and 
$\de_{n+1}\le\frac{1}{n+1}$ such that 
$$
\text{$\overline{U}(b_{n+1},\,\de_{n+1})\sus U(b_n,\,\de_n)$ and $M_{n+1}\cap \overline{U}(b_{n+1},\,\de_{n+1})=\emptyset$}\,.
$$
Thus we prolonged the above sequence by $n+1$-st term. Since for every $n\in\N$ 
the centers $b_n$, $b_{n+1}$, $\ds$ lie in $\overline{U}(b_n,\de_n)$, the 
sequence $(b_n)$ ($\sus M$) is Cauchy. We invoke Theorem~\ref{thm_CauchyPodm} 
and take the limit
$$
b\equiv\lim b_n\,.
$$
Since $M$ and $\overline{U}(b_n,\de_n)$ are closed sets, we see that $b\in M$ 
and $b\in\overline{U}(b_n,\de_n)$ for every $n\in\N$. Hence $b\not\in M_n$ for 
every $n\in\N$, which is a~conrtadiction.
\kduk
\vspace{-3mm}
\begin{exer}\label{ex_applBaire}
Suppose that $M\sus\R$ is a~nonempty closed set such that for every $b\in M$
and every $\de$ we have $M\cap P(b,\de)\ne\emptyset$. Then $M$ is 
uncountable. 
\end{exer}

\begin{exer}\label{ex_BaireInMS}
State and prove general version of Baire's theorem for metric spaces.    
\end{exer}

\section[Uniform continuity]{Uniform continuity}\label{sec_UC}

\noindent
{\em $\bullet$ Uniform continuity. }We define an important strengthening of continuity. 

\begin{defi}[uniform continuity]\label{def_UC}
Let $f\in\mathcal{R}$ and $M\sus\R$. We say that $f$ is   
\underline{uniformly 
continuous\index{continuity of 
functions!uniform on a~ set|emph}\index{uniform continuity|emph}} 
on $M$ if for every $\ep$ there is a~$\de$ such that for any points $a,b\in M\cap M(f)$, 
$$
|a-b|\le\de\Rightarrow|f(a)-f(b)|\le\ep\,.
$$
If $M=M(f)$, we say that $f$ is   
\underline{uniformly 
continuous\index{continuity of 
functions!uniform|emph}}, or {\em UC}.\label{UC}
\end{defi}
We denote the subset of uniformly 
continuous functions in $\mathcal{F}(M)$ by $\mathcal{UC}(M)$\index{function!ucm@$\mathcal{UC}(M)$|emph}\label{UCM} 
and set\index{function!uc@$\mathcal{UC}$|emph}$$
{\textstyle
\mathcal{UC}\equiv\bigcup_{M\sus\R}
\mathcal{UC}(M)\,.\label{UCcal}
}
$$

\begin{exer}\label{ex_stespojJespoj}
Every uniformly continuous function is continuous. 
\end{exer}

\noindent
{\em $\bullet$ Properties of {\em UC} functions. }For compact definition domains continuous functions are UC.

\begin{thm}[continuous $\Rightarrow$ UC]\label{thm_steSpo}
If\index{theorem!continuous $\Rightarrow$ UC|emph} 
$M\sus\R$ is a~compact set \underline{then}  
$$
\mathcal{C}(M)\sus\mathcal{UC}(M),\,\text{ hence }\,\mathcal{C}(M)=\mathcal{UC}(M)\,. 
$$
In words, any continuous function with compact definition domain is uniformly continuous.
\end{thm}
\duk
Suppose that $M\sus\R$ is compact and that $f\in\mathcal{F}(M)$ is not 
uniformly continuous. We prove that there is a~point 
$c\in M$ such that $f$ is not continuous at $c$. By the assumption there is an $\ep$ such that for every $\de$ there exist points  
$a,b\in M$ such that 
$$
|a-b|\le\de\,\text{ and }\,|f(a)-f(b)|>\ep\,. 
$$
Using the axiom of choice\index{axiom!of choice, AC} we get two sequences $(a_n),(b_n)\sus M$ such that 
for every $n$, 
$${\textstyle
|a_n-b_n|\le\frac{1}{n}\,\text{ and }\,|f(a_n)-f(b_n)|>\ep\,. 
}
$$
Using compactness of $M$ we pass to 
subsequences (Exercise~\ref{ex_passSubs}) and get that $\lim a_n=\lim b_n=c$ for some $c\in M$. Since $|f(a_n)-f(b_n)|>\ep$ for every $n$, it is not true that $\lim f(a_n)=\lim f(b_n)=f(c)$. By (H) $f$ is not continuous at~$c$.
\kduk
\vspace{-3mm}
\begin{exer}\label{ex_passSubs}
Explain the step where we pass to subsequences.     
\end{exer}

\begin{exer}\label{ex_nejsouSS}
The  functions $f(x)\equiv
\frac{1}{x}\,|\,(0,1]$ and $g(x)\equiv\sin(\frac{1}{x})\,|\,(0,1]$ are continuous but not {\em UC}.
\end{exer}

\begin{exer}\label{ex_HMCfun}
Let $M\equiv[0,1]\cap\Q$. Find a~function $f\in\mathcal{C}(M)\setminus\mathcal{UC}(M)$.
\end{exer}

The \underline{closure\index{set!closure of, $\overline{M}$|emph}\index{closure of a~set|emph}} of a~set $M\sus\R$ is the set
$$
\overline{M}\equiv\{b\in\R\cc\;\text{$\exists
\,(a_n)\sus M$ with $\lim a_n=b$}\}\ \ (\sus\R)\,.\label{closure}
$$
It is the set $(L(M)\cup M)\setminus\{-\infty,+\infty\}$.

\begin{exer}\label{ex_onClosure}
$M\sus\R$ is closed iff $\overline{M}=M$.   
\end{exer}

\begin{exer}\label{ex_ImUCbou}
Prove the next proposition.    
\end{exer}

\begin{prop}[UC boundedness]\label{prop_ImUCbou} 
Let $M\sus\R$ be bounded and let $f$ be in $\mathcal{UC}(M)$. 
\underline{Then} 
$$
\text{the image $f[M]$ is bounded}\,, 
$$
so that $f$ is a~bounded function.
\end{prop}

The next theorem is important.

\begin{thm}[UC extensions]\label{thm_steSpoRoz}
Let\index{theorem!extending UC functions|emph} 
$f\in\mathcal{UC}(M)$. The following holds.
\begin{enumerate}
\item Let $b\in\overline{M}$. For every sequence $(a_n)\sus M$ with $\lim a_n=b$ the limit
$$
g(b)\equiv\lim f(a_n)\ \ (\in\R)
$$
exists and does not depend on the sequence $(a_n)$.
\item The function
$$
g\cc\overline{M}\to\R
$$
so defined is the unique {\em UC} extension of $f$ to the closure of $M$.
\end{enumerate} 
\end{thm}
\duk
Let $f$ and $M$ be as stated. 

1. Let $b\in\overline{M}$ and $(a_n)\sus M$ be any sequence with $\lim 
a_n=b$. Let an $\ep$ be given. We take a~$\de$ such that 
$${\textstyle
x,\,y\in M,\,|x-y|\le\de\Rightarrow
|f(x)-f(y)|\le\ep\,.
}
$$
Since $(a_n)$ is Cauchy, there is $n_0$
such that for every $m,n\ge n_0$ we have
$|a_m-a_n|\le\de$. Thus for the same $m$
and $n$ we have
$$
|f(a_m)-f(a_n)|\le\ep
$$
and the sequence $(f(a_n))$ is Cauchy. By Theorem~\ref{thm_CauchyPodm} it has a~finite limit
$$
c\equiv\lim f(a_n)\,.
$$
Let $(a_n')\sus M$ be another sequence with $\lim a_n'=b$ and $\lim f(a_n')=c'$. If $c\ne c'$ then
$$
(b_n)\equiv(a_1,\,a_1',\,a_2,\,a_2',\,a_3,\,\ds)
$$
is a~sequence converging to $b$ such that  $\lim f(b_n)$ does 
not exist, which is a~contradiction. Hence $c$ does not depend on $(a_n)$ and
we can define a~function $g\cc\overline{M}\to\R$ by setting 
$$
g(b)\equiv\lim f(a_n)
$$
for any sequence $(a_n)\sus M$ with $a_n\to b$.

2. We prove three results: (i) $g$ extends $f$, (ii) every continuous extension of 
$f$ to $\overline{M}$ equals $g$ 
and (iii) $g$ is UC. 

(i) Let $b\in M$. We
take the constant sequence $(a_n)=(b,b,\ds)\sus M$. It converges to $b$ and we get that 
$$
g(b)=\lim f(a_n)=\lim f(b)=f(b)\,.
$$

(ii) Let $h\in\mathcal{C}(\overline{M})$ extend $f$ and $b\in\overline{M}$. We take any
sequence $(a_n)\sus M$ with $\lim a_n=b$. By (H) and the definition of $g$ we have 
$$
h(b)=\lim h(a_n)=\lim f(a_n)=g(b)\,.
$$

(iii) Let an $\ep$ be given. We take a~$\de$ such that 
$${\textstyle
a,b\in M,\ |a-b|\le\de\Rightarrow
|f(a)-f(b)|\le\frac{\ep}{3}\,. 
}
$$
Let $a,b\in\overline{M}$ with $|a-b|\le
\frac{\de}{3}$. We take two sequences $(a_n),(b_n)\sus M$ such that $\lim 
a_n=a$ and $\lim b_n=b$, and take their entries $a_m$ and $b_n$ such 
that $|a_m-a|,|b_n-b|\le\frac{\de}{3}$ and $|f(a_m)-g(a)|,|f(b_n)-g(b)|\le\frac{\ep}{3}$. Then 
\begin{eqnarray*}
&&{\textstyle|a_m-b_n|\le|a_m-a|+|a-b|+|b-b_n|\le\frac{\de}{3}+\frac{\de}{3}+\frac{\de}{3}=\de}\,\text{ and 
hence }\\
&&
|g(a)-g(b)|\le|g(a)-f(a_m)|+|f(a_m)-f(b_n)|+|f(b_n)-g(b)|\le\\
&&{\textstyle\le\frac{\ep}{3}+\frac{\ep}{3}+\frac{\ep}{3}
=\ep\,.
}
\end{eqnarray*}
\kduk

\noindent
In \cite{klaz_RAWUS} we rebuild, with the help of this extension
theorem, a~large part of univariate real analysis so that
only functions in the families 
$$
\mathcal{UC}(M),\ M\sus\Q\,, 
$$
are used. 
We conclude this section with a~generalization of 
Corollary~\ref{cor_priMaxi}.

\begin{thm}[extremes of UC functions]\label{thm_HMCminmax}
Let\index{theorem!extremes of UC functions|emph} 
$M\sus\R$ be bounded and $f$ be in $\mathcal{UC}(M)$. 
\underline{Then} there exist points $b,c\in\overline{M}$ such that for every $a\in M$,
$$
f(b)\le f(a)\le f(c)\,. 
$$
Here $f(b)$ and $f(c)$ are understood to be the values of the
extension of $f$ to the closure of $M$ provided by Theorem~\ref{thm_steSpoRoz}.
\end{thm}
\duk
Let $M$ and $f$ be as stated, and let $g\in\mathcal{C}(\overline{M})$ 
be the (uniformly) continuous extension of $f$ obtained in Theorem~\ref{thm_steSpoRoz}. We 
prove the existence of the minimum $b$, for the maximum $c$ the argument is similar. We set
$$
B\equiv\inf(f[M])\ \ (\in\R^*)
$$
(Proposition~\ref{prop_ImUCbou} shows that $B\in\R$).
Let $(a_n)\sus M$ be a~sequence such that $\lim f(a_n)=B$. Using 
Theorem~\ref{thm_BolzWeier}, we 
select a~subsequence $(b_n)$ of
$(a_n)$ such that $\lim b_n=b$ ($\in\R$). It follows that
$$
b\in\overline{M}\,\text{ and }\,g(b)=\lim f(b_n)=\lim f(a_n)=B\,. 
$$
In particular,
$B\in\R$. Since $B$ is 
a~lower bound of $f[M]$, we get that 
$$
\text{$g(b)=B\le f(a)$ for every $a\in M$}\,.
$$
\kduk
\vspace{-3mm}
\begin{exer}\label{ex_explGenee}
The inequalities in the theorem hold for every $a\in\overline{M}$.     
\end{exer}

\section[Operations on functions and continuity]{Operations on functions and continuity}\label{sec_aritSpoj}

We investigate the interplay of the five operations in
Definition~\ref{def_oprOnR} with continuity. In this respect most
interesting is the inverse $f\mapsto f^{-1}$. We prove that sums of 
power series are continuous and use it to show that the functions $\exp x$, $\cos 
x$, and $\sin x$ are continuous. The chapter concludes with the
proof of continuity of elementary functions.

\medskip\noindent
{\em $\bullet$ Arithmetic of continuity. }Recall the operations of sum, product and ratio (division) in
Definition~\ref{def_oprOnR}.

\begin{thm}[arithmetic of continuity]\label{thm_aritSpojitosti}
For function
$f,g\in\mathcal{R}$ the  following holds.
\begin{enumerate}
\item If\index{theorem!arithmetic of continuity|emph} 
$f$ and $g$ are continuous at a~point $b\in M(f)\cap 
M(g)$, \underline{then} 
$$
\text{$f+g$ and $fg$ are continuous at $b$}\,. 
$$
If $f$ and $g$ are continuous at 
a~point $b\in M(f/g)$, \underline{then} 
$$
\text{$f/g$ is continuous at $b$}\,.
$$
\item If $f,g\in\mathcal{C}$ \underline{then} $f+g,fg,f/g\in\mathcal{C}$.
\end{enumerate} 
\end{thm}
\duk
1. We consider only $f/g$, arguments for sum and product are similar
and easier. Let $f$, $g$ and $b$ be as stated and let $(a_n)\sus M(f/g)$ 
have $\lim a_n=b$. By (H) one has 
$\lim f(a_n)=f(b)$ and $\lim g(a_n)=g(b)$. 
By Theorem~\ref{thm_ari_lim} we have 
$${\textstyle
\lim (f/g)(a_n)=\lim\frac{f(a_n)}{g(a_n)}=\frac{\lim f(a_n)}{\lim g(a_n)}=
\frac{f(b)}{g(b)}=(f/g)(b)\,.
}
$$
Thus by (H) the function $f/g$ is continuous at $b$.

2. This follows from the first part.
\kduk
\vspace{-3mm}
\begin{exer}\label{ex_spojPOLaRAC}
Show that $\mathrm{POL}\sus\mathcal{C}$ and 
$\mathrm{RAC}\sus\mathcal{C}$, that is, every polynomial and every rational function is continuous. 
\end{exer}

\noindent
{\em $\bullet$ Continuity of power series. }We want to prove
continuity of all functions introduced
in Sections~\ref{sec_elemenFce} and \ref{sec_obecElemFce}. Those that were defined by composition and inverses will be discussed later. The continuity
of exponential, cosine, and sine follows from the next theorem.

\begin{thm}[continuity of $\sum_{n=0}^{\infty} a_nx^n$]\label{thm_spjMocRady1}
Suppose that\index{theorem!continuity of power series|emph} the real numbers
$a_0$, $a_1$, $\ds$, are such that
$\lim|a_n|^{1/n}=0$.
\underline{Then} for every $x\in\R$,  
$${\textstyle
S(x)\equiv\sum_{n=0}^{\infty}a_nx^n
}
$$ 
is an abscon series and the sum defines a~function
$S(x)$ that is in $\mathcal{C}(\R)$.
\end{thm}
\duk
Let the coefficients $(a_n)$ be as stated and $x\in\R$. Then
$0\le|a_n|^{1/n}|x|\le\frac{1}{2}$ for every $n\ge 
n_0$, so that $|a_nx^n|\le(\frac{1}{2})^n$ for $n\ge n_0$. The series $\sum_{n=0}^{\infty}a_nx^n$ is 
therefore abscon and converges. With $d\equiv\max(|x|,1)$ and $|c|\le1$,
\begin{eqnarray*}
|S(x+c)-S(x)|&=&{\textstyle|c|\cdot\big|\sum_{n=1}^{\infty}a_n\sum_{i=1}^n\binom{n}{i}c^{i-1}x^{n-i}\big|}\\
&\le&{\textstyle|c|\cdot\sum_{n=1}^{\infty}
|a_n|\cdot(2d)^n=O(|c|)\,\text{ (on $[-1,1]$)}\,.
}
\end{eqnarray*}
Thus the function $S(x)$ is continuous at $x$.
\kduk
\vspace{-3mm}
\begin{exer}\label{ex_rov2nerov}
Explain the displayed computation in the previous proof. 
\end{exer}

\begin{exer}\label{ex_odmFakt}
Prove that $\lim (n!)^{1/n}=+\infty$. 
\end{exer}

\begin{cor}[$\mathrm{e}^x$, cosine and sine]\label{cor_spoExpKosSin}
The functions $\mathrm{e}^x$, $\cos x$ and $\sin x$ are continuous on $\R$.     
\end{cor}
\duk
This follows from  the definitions
$${\textstyle
\mathrm{e}^x=\sum_{n=0}^{\infty}\frac{x^n}{n!},\ 
\cos x=\sum_{n=0}^{\infty}(-1)^n\frac{x^{2n}}{(2n)!}\,\text{ and }\,\sin x=\sum_{n=0}^{\infty}(-1)^n\frac{x^{2n+1}}{(2n+1)!}\,,
}
$$ 
from Theorem~\ref{thm_spjMocRady1} and from Exercise~\ref{ex_odmFakt}. 
\kduk
\vspace{-3mm}
\begin{cor}[tangent and cotangent]\label{cor_spoTamCot}
Both functions $\tan x$ and $\cot x$ are continuous on their definition domains.   
\end{cor}
\duk
As we know, $\tan x=\frac{\sin x}{\cos x}$ and $\cot x=\frac{\cos x}{\sin x}$. 
Thus continuity of these functions follows from Corollary~\ref{cor_spoExpKosSin} 
and part~2 of Theorem~\ref{thm_aritSpojitosti}.
\kduk

\noindent
{\em $\bullet$ Restriction and composition. }We treated restrictions 
in Proposition~\ref{prop_spojRest}. Now we supplement it with the
point-wise form. 

\begin{prop}[restriction 2]\label{prop_spojRest2}
Let $f\in\mathcal{R}$ and $X\sus\R$. Then the following holds.
\begin{enumerate}
\item If $f$ is continuous at $b\in M(f\,|\,X)$ \underline{then} $f\,|\,X$ is continuous at $b$.
\item If $f\in\mathcal{C}$ \underline{then} $f\,|\,X\in\mathcal{C}$.
\end{enumerate}
\end{prop}
\duk
1. Let $f$, $b$ and $X$ be as stated. We take any sequence 
$$
(b_n)\sus M(f\,|\,X)=M(f)\cap X
$$ 
with $\lim b_n=b$. By (H) we have that $\lim 
f(b_n)=f(b)$. Thus $\lim(f\,|\,X)(b_n)=\lim f(b_n)=f(b)$ and (H) gives that
$f\,|\,X$ is continuous at $b$.

2. This follows from part~1. Another proof, using neighborhoods, is given in Proposition~\ref{prop_spojRest}.
\kduk

We consider the operation of composition.

\begin{thm}[continuity of composition]\label{thm_spojSloz}
Let $f,g\in\mathcal{R}$. Then the following holds.
\begin{enumerate}
\item If\index{theorem!continuity of composition|emph} 
$g$ is continuous at $b\in M(f(g))$ and $f$ is continuous at $g(b)$ 
\underline{then} $f(g)$ is continuous at $b$.
\item If 
$f,g\in\mathcal{C}$ \underline{then} $f(g)\in\mathcal{C}$.
\end{enumerate}
\end{thm}
\duk
1. Let $f$, $g$ and $b$ be as stated,  and 
$$
(b_n)\sus M(f(g))
$$ 
be a~sequence with $\lim b_n=b$. Then $\lim g(b_n)=g(b)$ by (H). Hence, 
again by (H), $\lim f(g)(b_n)=\lim f(g(b_n))=f(g(b))=f(g)(b)$.
(H) shows that $f(g)$ is continuous at $b$.

2. This follows from the first part.
\kduk
\vspace{-3mm}
\begin{exer}\label{ex_naSpoSkla}
Prove part~1 of Theorem~\ref{thm_spojSloz} by means of neighborhoods.    
\end{exer}

\noindent
{\em $\bullet$ Inverses. }The interplay of inverses and continuity
is much more interesting than what we have seen so far. The next exercise shows that, in general, inverses do not preserve continuity. 

\begin{exer}\label{ex_discInve}
Let $f\in\mathcal{F}(\N_0)$ be defined 
by 
$${\textstyle
f(0)\equiv0,\ f(n)\equiv\frac{1}{n}\,\ds\,n\in\N\,. 
}
$$
Show that $f\in\mathcal{C}$ but that the inverse $f^{-1}$ is discontinuous. 
\end{exer}

On the other hand, in many situations inverses do preserve continuity.

\begin{thm}[continuity of inverses]\label{thm_spojInverzu}
Let\index{theorem!continuity of inverses|emph} $M\sus\R$, 
$f\in\mathcal{C}(M)$ and $f$ be injective. 
\underline{Then} in each of five situations the inverse $f^{-1}\in\mathcal{F}(f[M])$ is continuous. 
\begin{enumerate}
\item When $M$ is a~compact set.
\item When $M$ is an interval.
\item When $M$ is an open set.
\item When $M$ is a~closed set and $f$ is a~monotone function.
\item When $M\sus(a,b)$, $M$ is dense in  $(a,b)$ and $f$ is monotone and {\em UC}.
\end{enumerate}
\end{thm}
\duk
Let $M$ and $f$ be as stated.

1. Let $M$ be compact, $b\in f[M]$ and let 
$$
\text{$(b_n)\sus f[M]$ be a~sequence with $\lim b_n=b$}\,.
$$
Let $a\equiv f^{-1}(b)$ and $a_n\equiv f^{-1}(b_n)$ ($\in M$). We show that $\lim a_n=a$, 
which by (H) proves continuity of $f^{-1}$ at $b$. We show that every
subsequence of $(a_n)$ has a~subsequence with the limit $a$. Part~3 of 
Theorem~\ref{thm_oPodposl} then implies that 
$\lim a_n=a$. Let $(a_n')$ be a~subsequence of $(a_n)$. We use compactness of 
$M$ and take a~subsequence $(a_{m_n})$ of $(a_n')$ with $\lim\,a_{m_n}=c\in M$. By (H) it holds that 
$\lim\,f(a_{m_n})=f(c)=b$ because
$(f(a_{m_n}))$ is a~subsequence of $(b_n)$. Since $f$ is injective, $c=a$.

2. Let $M$ be an interval. Corollary~\ref{cor_spojInjInt} shows that $f$ increases 
or decreases. Suppose that $f$ decreases, the case of  increasing $f$ is similar. 
Theorem~\ref{thm_mezihodnoty} says that $f[M]$ is an interval. Let $b\in f[M]$ and an 
$\ep$ be given. We show that 
$f^{-1}$ is right-continuous at $b$. This is true when $b$ is the right endpoint of the interval
$f[M]$ because then 
$$
U^+(b,\,\de)\cap f[M]=\{b\}\,. 
$$
We assume that $b$ is not the right endpoint of 
$f[M]$. Since $f^{-1}$ decreases, $a\equiv f^{-1}(b)$ ($\in M$) is not 
the left endpoint of the interval $M$. We take a~small $\ep$ such that $[a-\ep,a]\sus M$. We set 
$$
\de\equiv f(a-\ep)-f(a)=f(a-\ep)-b\ \ (>0)\;. 
$$
Theorem~\ref{thm_mezihodnoty} implies that $f$ is a~decreasing 
bijection from 
$[a-\ep,a] $ to $[b,b+\de]$, and
from $(a-\ep,a]$ to $[b,b+\de)$. So $[b,b+\de)\sus f[M]$ and $U^+(b,\,\de)\cap f[M]=U^+(b,\,\de)=[b,b+\de)$. Hence
$$
f^{-1}[U^+(b,\,\de)]=U^-(a,\,\ep)\sus U(a,\,\ep) =U(f^{-1}(b),\,\ep)
$$
and $f^{-1}$ is right-continuous at $b$. The left continuity of $f^{-1}$ at $b$ is proven similarly. By Exercise~\ref{ex_ekvivSpoji}, 
$f^{-1}$ is continuous at $b$.

3. Let $M$ be open, $b\in f(M)$, $a\equiv f^{-1}(b)$ ($\in M$) and let $\ep$ be 
given. We take $\ep$ so small that $U(a,\ep)\sus M$.
Proposition~\ref{prop_Obropen} says that the set $f[U(a,\ep)]$ ($\ni b$) is open. So for some $\de$ we have that $U(b,\de)\sus f[U(a,\ep)]$. Thus
$$
f^{-1}[U(b,\,\de)]\sus U(a,\,\ep)=U(f^{-1}(b),\,\ep)\,.
$$
So $f^{-1}$ is continuous at $b$ by Definition~\ref{def_spojVbode}.

4. Let $M$ be closed and $f$ be increasing. For decreasing $f$ we 
argue similarly. We assume for contradiction that 
for some $b\in f[M]$ there is 
a~sequence $(b_n)\sus f[M]$ such that $$
\text{$\lim b_n=b$ 
but $\lim f^{-1}(b_n)$ does not exist or $\ne a\equiv f^{-1}(b)$ ($\in M$)}\,. 
$$
By part~2 
of Theorem~\ref{thm_oPodposl} 
and by Proposition 
the sequence $(b_n)$ has a~decreasing or an increasing 
subsequence 
$(c_n)$ such that $\lim f^{-1}(c_n)=B$ ($\in\R^*$) and $B\ne a$. 
We assume that $(c_n)$ decreases, 
the case of increasing $(c_n)$ is similar. Then 
$$
b<\ds<c_2<c_1\,\text{ and }\,a<\ds<f^{-1}(c_2)<f^{-1}(c_1)
$$ 
because both $f$ and $f^{-1}$ increase. By part~2 of Theorem~\ref{thm_limAuspo} we have that $B\in[a,f^{-1}(c_1))$. Thus, crucially, $B\in\R$ (here the 
argument fails for non-monotone $f$). Even $B\in M$ because $M$ is closed. Due to the continuity of $f$ in 
$B$ we have that $f(B)=\lim f(f^{-1}(c_n))=\lim c_n=b=f(a)$. 
But this contradicts the injectivity of $f$ because $B\ne a$.

5. Let $M$, $a$, $b$ and $f$ be as stated.  
We use Theorem~\ref{thm_steSpoRoz} and continuously extend $f$ to 
a~function 
$\overline{f}\cc[a,b]\to\R$. By Proposition~\ref{prop_limAuspo} and because $M$ is dense
in $(a,b)$, the function $\overline{f}$ is strictly monotone and therefore injective. By part~1 or part~2 
of this theorem  the function $(\overline{f})^{-1}$ is continuous. By Proposition~\ref{prop_spojRest}, the restriction
$(\overline{f})^{-1}\,|\,f[M]=f^{-1}$ is continuous.
\kduk

\noindent
In {\em MA~1${}^+$} we generalize part~5. 

\begin{exer}\label{ex_aniJedeNevyn}
Give examples showing that in part~4 of the theorem neither the closedness of $M$ nor the monotonicity of $f$ can be omitted.
\end{exer}

\noindent
{\em $\bullet$ Continuity of elementary functions. }We prove that $\mathrm{EF}\sus\mathcal{C}$.

\begin{exer}\label{ex_spojLogCyklFci}
Prove the following corollary.   
\end{exer}

\begin{cor}[continuity of some BEF]\label{cor_spojLogadalsich}
The functions $\log x$, $\arccos x$, $\arcsin x$, $\arctan x$ and $\mathrm{arccot}\,x$ are continuous. 
\end{cor}

\begin{prop}[continuity of $x^b$]\label{prop_spojAnaB}
For every $b\in(0,+\infty)$ the function 
$x^b\cc[0,\,+\infty)\to[0,\,+\infty)$ is continuous.
\end{prop}
\duk
Let $b>0$ and $x>0$. We find that $x^b$
is continuous at $x$ using the expression
$x^b=\exp(b\log x)$, the continuity of $\mathrm{e}^x$ (Corollary~\ref{cor_spoExpKosSin}),
continuity of logarithm (Corollary~\ref{cor_spojLogadalsich}), continuity of the constant function $k_b$ (Exercise~\ref{ex_spojKonst}), continuity of the
product (Theorem~\ref{thm_aritSpojitosti}) and continuity of the composition (Theorem~\ref{thm_spojSloz}). Continuity at $x=0$ follows with the help of Proposition~\ref{prop_spojVbLimitou}
from the limit
$$
\lim_{x\to 0}x^b=\lim_{x\to 0}\exp(b\log x)=\lim_{y\to-\infty}\exp y=0=0^b\,.
$$
Here the second equality follows from Theorem~\ref{thm_LimSlozFunkce} and from part~2 of Proposition~\ref{prop_logar}. The third equality follows from part~3 of Proposition~\ref{prop_expFce}.
\kduk
\vspace{-3mm}
\begin{thm}[$\mathrm{EF}\sus\mathcal{C}$]\label{thm_EFjsouSpoj}
Every elementary function is continuous.\index{theorem!continuity of elementary functions|emph}
\end{thm}
\duk
We proceed by induction on the length of the generating word of the given elementary function $f$ 
(Definition~\ref{def_obecEF2}). If $f$ is a~constant function, exponential, logarithm, $x^b$ 
with non-integral exponent $b>0$, sine or arcsine, it is continuous by, respectively, 
Exercise~\ref{ex_spojKonst}, Corollaries~\ref{cor_spoExpKosSin} and \ref{cor_spojLogadalsich}, 
Proposition~\ref{prop_spojAnaB} and Corollaries~\ref{cor_spoExpKosSin} and \ref{cor_spojLogadalsich}. 
If $f$ is a~sum, a~product, a~ratio or a~composition of two simpler elementary functions, it is continuous by induction and Theorems~\ref{thm_aritSpojitosti}
and \ref{thm_spojSloz}.
\kduk

\chapter[Derivatives]{Derivatives}\label{chap_pr7}

We define the derivative at any limit point in the definition domain and consider
derivatives both locally at 
points, and globally as a~unary operation on $\mathcal{R}$. This 
is more general than the standard approach, 
which defines derivatives only at inner points of domains. Our 
approach yields natural and general formulas for the derivatives of inverse functions.

Section~\ref{sec_defDerExtr} contains basic definitions. Theorem~\ref{thm_priznakExtr} is our version 
of the popular criterion of local extremes for functions with arbitrary definition domains. By 
Proposition~\ref{prop_vlDerimplSpoj}, finite derivative implies continuity. If it is additionally nonzero, the point is mapped to a~limit point of 
the image (Proposition~\ref{prop_limBodyObr} and Exercise~\ref{ex_pro0deroPla}).
Theorem~\ref{thm_nespDer} describes a~function with discontinuous derivative. Another such function is 
described in Exercise~\ref{ex_znovuTenPr}.

Section~\ref{sec_tecny} starts with Definition~\ref{def_tangentLine} of 
standard tangent lines. 
Theorem~\ref{thm_suppLine} shows
that lines touching the graph in a single point are tangents. Definition~\ref{def_limiTecn} formalizes the intuition of a~tangent as a~limit of secants. By Theorem~\ref{thm_onTangents},  
standard and limit tangents are equivalent. 
Theorem~\ref{thm_NonmainSec} shows that the tangent is a~limit of 
secants whose two intersection points enclose the point of 
contact.

Section~\ref{sec_aritmDeri} is devoted to the arithmetic 
of derivatives. In Theorem~\ref{thm_nasob} we investigate local and global derivatives of sums. We consider the relation 
between $f'$, $g'$ and $(f+g)'$ for any pair of functions 
$f,g\in\mathcal{R}$; in general, $(f+g)'\ne f'+g'$. 
Theorem~\ref{thm_LeibnizFmle} presents local and global Leibniz formulas for derivatives of products. Theorem~\ref{thm_deriPodilu} does
the same for ratios. Corollaries~\ref{cor_deriSumStan}, 
\ref{cor_Leibform}, \ref{cor_deriRat}, 
\ref{cor_derSlozF} and \ref{cor_deriInveSim} provide 
assumptions, under which the respective formulas
\begin{eqnarray*}
&&{\textstyle
(f+g)'=f'+g',\ (fg)'=f'g+fg',\  
(\frac{f}{g})'=
\frac{f'g-fg'}{g^2},\ (f(g))'=
f'(g)\cdot g'}\\
&&{\textstyle\text{and }\,(f^{-1})'=\frac{1}{f'(f^{-1})}
}
\end{eqnarray*} 
do hold. 

In Section~\ref{sec_deriSlozInve} we turn to
derivatives of composite functions (Theorem~\ref{thm_DerSlozFce}) and inverses
(Theorem~\ref{thm_DeriInvfce} and Corollary~\ref{cor_deriInve}). As before, we
give local and global formulas and allow
arbitrary definition domains. Proofs use Heine's 
definition of derivatives. 

In Section~\ref{sec_deriElemFun} 
in Theorem~\ref{thm_derMocRady} we differentiate power series. By this we get 
derivatives of the functions $\exp x$, $\sin x$, and $\cos x$. We also find the derivative of $\log x$, but the 
derivatives of other basic elementary functions (Definition~\ref{def_ZEF}) are left
to exercises. On page~\pageref{tableDer} we summarize these derivatives in a~table. 

In Section~\ref{sec_deriEF} we 
pose Problem~\ref{prob_derEF}: decide if derivatives of elementary 
functions are always elementary. In 
Theorem~\ref{thm_deriSEF} we prove that a~subfamily of so called 
simple elementary
functions\index{simple elementary 
functions, SEF} is closed to 
derivatives.

\section[${}^c$Local and global derivatives]{Local and global derivatives}\label{sec_defDerExtr}

Why derivatives? Using them we obtain local linear approximations of functions and tangent lines to their graphs. In Chapter~\ref{chap_pr9} we 
investigate more precise local polynomial approximations. The
derivative $f'$ of a~function $f$ records properties of $f$ in simpler form. For example, the increase of $f$ is recorded as the positivity of $f'$. 
If $f$ is elementary, it is straightforward 
to compute $f'$ locally. However, see Problem~\ref{prob_derEF}.

\medskip\noindent 
{\em $\bullet$ Local derivatives. }Recall that 
$\mathcal{F}(M)$ is the set of functions $f\cc M\to\R$ where $M\sus\R$, that $\mathcal{R}=\bigcup_{M\sus\R}\mathcal{F}(M)$, and that $L(M)$ ($\sus\R^*$) is the set of limit points of $M$. 

\begin{defi}[local derivatives]\label{def_Derivace}
Let $f\in\mathcal{F}(M)$ and $b\in M\cap L(M)$. The 
\underline{derivative\index{derivative of a function!at a point|emph}} of the function $f$ at the point $b$ is the limit
$$
f'(b):=\lim_{x\to b}\frac{f(x)-f(b)}{x-b}\ 
\ (\in\R^*)\,.\label{deriv}
$$
\end{defi}
This limit need not exist. If $f'(b)$ exists, 
then $b\in M(f)\cap L(M(f))$ and we often do not mention it explicitly.
The uniqueness of limits of functions implies the uniqueness of derivatives. Our definition is more general
than the standard definition, which allows for $b$ only inner points of definition domains. 

\begin{exer}\label{ex_defDeri}
Prove that also $f'(b)=\lim_{h\to0}\frac{f(b+h)-f(b)}{h}$.
\end{exer}

\begin{prop}[derivatives are local]\label{prop_ofTheDef}
Let\index{derivative of a function!at a point!locality of|emph} 
$f,g\in\mathcal{R}$ and $b$ be in $M(f)\cap M(g)$. If 
$f\,|\,U(b,\theta)=
g\,|\,U(b,\theta)$ for some $\theta>0$, \underline{then}  
$$
f'(b)=g'(b)\,, 
$$
if either side is defined.
\end{prop}
\duk
This is immediate from Proposition~\ref{prop_locLimFce}
and Definition~\ref{def_Derivace}.
\kduk

\noindent
Unlike in Proposition~\ref{prop_locLimFce}, the neighborhood $P(b,\theta)$ does not 
suffice because the definitions of $f'(b)$ and $g'(b)$ involve $f(b)$ and $g(b)$. 

Let $f\in\mathcal{R}$ and $b\in M(f)$.
If $f'(b)\in\R$, we say that $f$ is 
\underline{differentiable\index{derivative of a  function!at a point!differentiability|emph}} at $b$. 
Then $f$ has near $b$ a~local approximation by the linear function
$$
t(x)=f(b)+f'(x)\cdot(x-b)\cc\R\to\R\,.
$$
We say that $t(x)$ is the \underline{tangent\index{tangent@tangent (line)!standard|emph} (line)} to $f(x)$ at $b$. The approximation takes the form
$$
f(x)=\underbrace{f(b)+f'(b)\cdot(x-b)}_{\text{the tangent at $b$}}\,\,+\,\,\underbrace{o(x-b)}_{\text{the error}}\ \  
(x\to b)\,.
$$
Heine's definition of derivatives, abbreviated 
HDD\index{derivative of a function!at a point!Heine's 
definition of, HDD|emph}\index{HDD|emph}\label{HDD}, is a~useful tool in proofs. 

\begin{prop}[HDD]\label{prop_HeineDeri}
Let $f\in\mathcal{F}(M)$, $b\in M\cap L(M)$ and $B\in\R^*$. \underline{Then} 
$f'(b)=B$ $\iff$ for every sequence $(a_n)\sus M\setminus\{b\}$ with $\lim a_n=b$,
$$
\lim_{n\to\infty}\frac{f(a_n)-f(b)}{a_n-b}=B\,.
$$
\end{prop}
\duk
This follows from Definition~\ref{def_Derivace} and 
Theorem~\ref{thm_HeinehoDef}.
\kduk

\noindent
{\em $\bullet$ One-sided derivatives. }Let $M\sus\R$. Recall that $L^-(M)$ 
($\sus\R$) are the left-sided limit points of the set $M$. Similarly, $L^+(M)$ 
($\sus\R$) are the right-sided limit points of $M$.

\begin{defi}[one-sided derivatives]\label{def_oneSiDer}
Let $f\in\mathcal{F}(M)$ and $b$ be in $M\cap L^-(M)$. The 
\underline{left-sided 
derivative\index{derivative of a function!at a point!left-sided|emph}} of the function $f$ at the point $b$ is the left-sided limit
$$
f_-'(b):=\lim_{x\to b^-}\frac{f(x)-f(b)}{x-b}\ \  (\in\R^*)\,.\label{derivPm} 
$$
Changing the signs $-$ in the indices to
$+$, we get the \underline{right-sided
derivative\index{derivative of a function!at a point!right-sided|emph}}
$f_+'(b)$ of $f$ at $b$.    
\end{defi} 

\begin{exer}\label{ex_ojednDeri}
The following holds.
\begin{enumerate}
\item If $f'(a)=A$ then $f_-'(a)=A$ or $f_+'(a)=A$.
\item If $f_-'(a)=f_+'(a)=A$ then $f'(a)=A$.
\item If $f_-'(a)=A\ne B=f_+'(a)$ then $f'(a)$ does not exist.
\end{enumerate}

\end{exer}

\begin{exer}\label{ex_jakoVtvrz}
As in Proposition~\ref{prop_ojednLimi2}, 
we can reduce one-sided derivatives to two-sided by 
restricting the function. State 
it in detail and prove it.  
\end{exer}

\noindent
{\em $\bullet$ Zeros of derivatives and extremes. }Let 
$M\sus\R$. Recall that $a\in\R$ is a~two-sided limit point 
of $M$ if for every $\de$, both intersections $(a-\de,a)\cap M$ and $(a,a+\de)\cap M$ are nonempty.
The set of two-sided limit points
of $M$ is denoted by
$L^{\mathrm{TS}}(M)$ ($\sus\R$).

\begin{exer}\label{ex_oOLB}
$L^{\mathrm{TS}}(M)\sus L(M)$ and this inclusion is in general strict.
\end{exer}

Here is the well known criterion of local extremes in terms of vanishing derivatives for general domains. 

\begin{thm}[zeros and extremes]\label{thm_priznakExtr}
Let\index{theorem!zeros and extremes|emph} $f\in\mathcal{F}(M)$ and $b\in M\cap L^{\mathrm{TS}}(M)$. If $f'(b)\in\R^*\setminus\{0\}$ \underline{then} for every $\de$ there exist $c,d\in U(b,\de)\cap M$ such that 
$$
f(c)<f(b)<f(d)
$$
---$f$ does not have at $b$ a~local extreme.
\end{thm}
\duk
Let $f'(b)<0$, the case when $f'(b)>0$ is similar. We take an $\ep$ such 
that $U(f'(b),\ep)<0$. Let a~$\de$ 
be given. By Definition~\ref{def_Derivace}, there is a~$\theta\le\de$ such that for every $x\in P(b,\theta)\cap M$ we have 
$$
{\textstyle
\frac{f(x)-f(b)}{x-b}\in U(f'(b),\,\ep)\,\text{ and hence }\,
\frac{f(x)-f(b)}{x-b}<0\,.
}
$$
Since $b\in L^{\mathrm{TS}}(M)$, we can take points 
$c\in(b,b+\theta)\cap M$ and $d\in (b-\theta,b)\cap M$. Then, by the above inequality, $f(c)
<f(b)$ and $f(d)>f(b)$. 
\kduk
\vspace{-3mm}
\begin{exer}\label{ex_jinyExer}
In the theorem, $U(b,\theta)$ can be replaced with $P(b,\theta)$.    
\end{exer}

\begin{exer}\label{ex_Contra}
Function $f(x)=x\cc[0,1]\to\R$ has derivatives $f'(0)=f'(1)=1\ne0$
and global extremes at $0$ and $1$. Does it 
contradict the theorem?
\end{exer}

We restate the previous theorem in the contrapositive.

\begin{cor}[on local extremes]\label{cor_NPELE}
Let $f\in\mathcal{F}(M)$ and let $b\in M$. If the function $f$ has a local extreme at the point $b$, \underline{then} (exactly) one of the following three claims holds.
\begin{enumerate}
\item $b\not\in 
L^{\mathrm{TS}}(M)$.
\item $b\in
L^{\mathrm{TS}}(M)$ but the derivative $f'(b)$ does not exist.
\item $b\in 
L^{\mathrm{TS}}(M)$ and $f'(b)=0$.
\end{enumerate}  
\end{cor}

\noindent
{\em $\bullet$  An example. }We want to find the extremes of the function 
$$
f(x)=x^2\cc\Q\to\Q\,. 
$$
It is easy to see that $f'(x)=2x\cc\Q\to\Q$.
Since every $b\in M(f)$ is a~two-sided limit point of $M(f)$ and 
$f'(b)$ always exists, Corollary~\ref{cor_NPELE} says 
that $f$ may have a~local extreme only at a~zero of $f'$. We have $f'(b)=0$ iff $b=0$. So there is a~single ``suspicious'' point $b=0$. Indeed, $f$ has at $0$ a~strict global minimum.   

\begin{exer}\label{ex_whatCh}
What about the function $f(x)=x^2\cc\Z\to\N_0$?    
\end{exer}

\noindent
{\em $\bullet$ Derivatives and continuity. }In this passage we show that 
differentiability  
at a~point is a~stronger property than local continuity at that 
point.

\begin{prop}[derivatives and continuity]\label{prop_vlDerimplSpoj}
Let $f\in\mathcal{R}$. If $f'(b)\in\R$ \underline{then} $f$ is continuous at $b$.
\end{prop}
\duk
We compute 
\begin{eqnarray*}
\lim_{x\to b}f(x)&=&
\lim_{x\to b}{\textstyle
\big(f(b)+(x-b)\cdot\frac{f(x)-f(b)}{x-b}\big)}\\
&=&\lim_{x\to b}f(b)+\lim_{x\to b}(x-b)\cdot\lim_{x\to b}{\textstyle
\frac{f(x)-f(b)}{x-b}}\\
&=&f(b)+0\cdot f'(b)=f(b)\,.
\end{eqnarray*}
Note that the function $f(x)$ differs from the function $f(b)+(x-b)\cdot\frac{f(x)-f(b)}{x-b}$  (Exercise~\ref{ex_fceVleVpr}). The 
first equality is therefore nontrivial and follows from Proposition~\ref{prop_locLimFce}. 
The second equality follows from Theorem~\ref{thm_AritLimFce}. In 
the third equality, we use that $f'(b)\ne\pm\infty$. By 
Proposition~\ref{prop_spojVbLimitou}, the function $f$ is continuous at the point $b$.
\kduk
\vspace{-3mm}
\begin{exer}\label{ex_fceVleVpr}
We have $f(x)\ne f(b)+(x-b)\cdot\frac{f(x)-f(b)}{x-b}$.   
\end{exer}

\begin{exer}\label{ex_deriSignu}
Show that $\sgn'(0)=+\infty$.  
\end{exer}
So an infinite derivative at a~point does not 
imply continuity 
at that point.

\begin{exer}\label{ex_deriAbsHod}
Show that $(|x|)'_-(0)=-1$ and $(|x|)'_+(0)=+1$. 
\end{exer}
By item~3 in Exercise~\ref{ex_ojednDeri}, $(|x|)'(0)$ does not exist. Of course, continuity at 
a~point does not imply the existence of a~derivative at that 
point. 

\begin{prop}[limit points of images]\label{prop_limBodyObr}
Let $f$ be in $\mathcal{F}(M)$. If $f'(b)$ is in $\R\setminus\{0\}$, 
\underline{then} 
$$
f(b)\in L(f[M])\,.
$$
\end{prop}
\duk
Let an $\ep$ be given. Since 
$f'(b)\ne0$, by Definition~\ref{def_Derivace} there 
is a~$\de$ such that $f(x)\ne f(b)$ for every 
$x\in P(b,\de)\cap M$. By 
Proposition~\ref{prop_vlDerimplSpoj}, we can take this $\de$ so small that for the same $x$, we have $f(x)\in 
U(f(b),\ep)$. Since $b\in L(M)$, there is a~point $a\in P(b,\de)\cap M$. Then
$$
f(a)\in P(f(b),\,\ep)\cap f[M]\,. 
$$
Hence $f(b)\in L(f[M])$.
\kduk

\noindent
We use this proposition  
to differentiate inverse
functions.

\begin{exer}\label{ex_proNekoNepl}
Proposition~\ref{prop_limBodyObr} does not hold if $f'(b)=\pm\infty$.
\end{exer}

\begin{exer}\label{ex_pro0deroPla}
Proposition~\ref{prop_limBodyObr} holds if $f'(b)=0$ and $f$ is non-constant on any neighborhood $U(b,\de)$.
\end{exer}

\begin{exer}\label{ex_oneSiDiCo}
Adapt Proposition~\ref{prop_vlDerimplSpoj}
for one-sided derivatives and
one-sided continuity.
\end{exer}

\noindent
{\em $\bullet$ Derivatives of the root. }Recall that
$\sqrt{x}=x^{1/2}$ and that $M(\sqrt{x})=[0,+\infty)$. Let $a\ge0$. We compute
$$
(\sqrt{x})'(a)=\lim_{x\to a}{\textstyle\frac{\sqrt{x}-\sqrt{a}}{x-a}}
=\lim_{x\to a}{\textstyle\frac{x-a}{(x-a)(\sqrt{x}+\sqrt{a})}}
=\lim_{x\to a}{\textstyle\frac{1}{\sqrt{x}+\sqrt{a}}\,.
}
$$
The last equality again follows from Proposition~\ref{prop_locLimFce}.
Thus $(\sqrt{x})'(0)=+\infty$, and for $a>0$ we have $(\sqrt{x})'(a)=\frac{1}{2\sqrt{a}}$. 
We see that an infinite derivative at a~point does not exclude continuity at that point.

\begin{exer}\label{ex_urciJednos}
Compute $(\sqrt{x})'_-(a)$ and $(\sqrt{x})'_+(a)$ for real  $a\ge0$.   
\end{exer}

\noindent
{\em $\bullet$ Global derivative. }We extend our repertoire of operations on $\mathcal{R}$ in 
Definition~\ref{def_oprOnR} by a~unary operation of (global) derivative. For 
$f\in\mathcal{R}$ we set\index{derivative of a function!df@$D(f)$|emph}
$$
D(f):=\{b\in M(f)\cc\;f'(b)\in\R\}\ \ (\sus M(f))\,.\label{defDeri}
$$
Thus $D(f)$ is the subset of the domain where $f$ is differentiable. 
We formally define a~function $f'\cc D(f)\to\R$\index{derivative of 
a function!global, $f'$|emph} by 
$$
D(f)\ni b\mapsto f'(b)\in\R\,.\label{derivGl}
$$

\begin{defi}[global derivative]\label{def_deriFunk}
We call the unary operation on $\mathcal{R}$, 
$$
\mathcal{R}\ni f\mapsto 
f'\in\mathcal{R}\,,
$$ 
the global \underline{derivative\index{global derivative|emph}}.
\end{defi} 
$D(f)=M(f')$ may be a~proper subset of $M(f)$, which causes 
troubles. For example, $M(\sqrt{x})=[0,+\infty)$ but
$$
D(\sqrt{x})=M((\sqrt{x})')=(0,\,+\infty)\,. 
$$
Let $b\in\R$.
The notation $f'(b)$ now becomes a~little ambiguous: as the derivative of 
$f$ at $b$ it may be $\pm\infty$, but as the value of the function 
$f'$ ($\in\mathcal{R}$) at $b$ it has to be finite. We will use such 
formulations that the interpretation of the symbol 
$f'(b)$ will be clear. We will investigate interactions of derivatives 
with various operations.

\begin{prop}[derivatives and restrictions]\label{ex_deriRest}
Let $f\in\mathcal{R}$, let $X\sus\R$, and let $M=M(f)\cap X$. The following holds.
\begin{enumerate}
\item If $f'(b)$ exists and $b\in M\cap L(M)$, \underline{then} $(f\,|\,X)'(b)=f'(b)$.
\item The restriction $f'\,|\,M\cap L(M)$ is a~subfunction of $(f\,|\,X)'$.
\end{enumerate}
\end{prop}
\duk
1. By Proposition~\ref{prop_restrAlimFce},  $(f\,|\,X)'(b)$ is  
$$
\lim_{x\to b}{\textstyle\frac{(f\,|\,X)(x)-
(f\,|\,X)(b)}{x-b}}=\lim_{x\to b}{\textstyle
\big(\frac{f(x)-f(b)}{x-b}\,\big|\,X\big)(x)}=
\lim_{x\to b}{\textstyle\frac{f(x)-f(b)}{x-b}=f'(b)\,.
}
$$

2. Let $N=M\cap L(M)$, $g= f'\,|\,N$ and $c\in M(g)$. Then $c\in D(f)\cap N$ and by part~1 we have $(f\,|\,X)'(c)=f'(c)=g(c)$. Hence $g=f'\,|\,N$ is a~subfunction of $(f\,|\,X)'$.
\kduk

\begin{exer}\label{ex_derKon}
For every $c\in\R$ we have $k_c'(x)=k_0(x)$.  
\end{exer}

\begin{exer}\label{ex_derIde}
For every $c\in\R$  we have $(\mathrm{id}_{\R}(x)+k_c(x))'=k_1(x)$.   
\end{exer}

\begin{exer}\label{ex_spojResDer}
For every $f\in\mathcal{R}$ we have $f\,|\,D(f)\in\mathcal{C}$.   
\end{exer}

\noindent
{\em $\bullet$ Global one-sided derivatives. }For any 
function 
$f\in\mathcal{R}$ we introduce sets
\begin{eqnarray*}
D_-(f)&:=&\{b\in M(f)\cc\;
f_-'(b)\in\R\}
\label{defDeriPm}\text{ and}\\
D_+(f)&:=&\{b\in M(f)\cc\;f_+'(b)\in\R\}\,.
\end{eqnarray*}
We define the global \underline{left-sided 
derivative\index{derivative of 
a function!left-sided global, 
$f_-'$|emph}} of $f$
by
$$
f_-'(x)\cc D_-(f)\to\R,\ 
D_-(f)\ni b\mapsto f_-'(b)\in\R\label{derivGlpm}\,,
$$
and the (global) \underline{right-sided derivative\index{derivative of 
a function!left-sided global, 
$f_-'$|emph}} of $f$
by
$$
f_+'(x)\cc D_+(f)\to\R,\ 
D_+(f)\ni b\mapsto f_+'(b)\in\R\,.
$$

\begin{exer}\label{ex_naGlobObou}
How do $D(f)$, $D_-(f)$ and $D_+(f)$ relate?    
\end{exer}

\noindent
{\em $\bullet$ Discontinuous derivatives. }We define
a~function $f\in\mathcal{R}$ such that $M(f)=D(f)\ne\emptyset$ (thus $f\in\mathcal{C}$) but $f'\not\in\mathcal{C}$.

\begin{thm}[discontinuous derivative]\label{thm_nespDer}
There\index{theorem!discontinuous derivatives|emph} exists a~function
$f\in\mathcal{R}$ such that
$$
\text{$M(f)=D(f)\ne\emptyset$ and $f'\not\in\mathcal{C}$}\,.   
$$
\end{thm}
\duk
Let $(a_n)$ and $(b_n)$ be real sequences going to $0$ such 
that
$$
a_1>b_1>a_2>b_2>\ds>0\,\text{ and }\,
a_n-b_n=o(b_n)\ \  (n\to\infty)\,.
$$
Let
$N=
\{0\}\cup\bigcup_{n=1}^{\infty}(b_n,\,a_n)$. 
We define $f\in\mathcal{F}(N)$ by 
$$
f(0)=0\,\text{ and $f(x)=x-b_n$ for $x\in(b_n,\,a_n)$}\,. 
$$
Let $x\in(b_n,a_n)$. Then, by Proposition~\ref{prop_ofTheDef} and 
Exercise~\ref{ex_derIde}, 
$f'(x)=1$. We have $0\in L(N)$ and $${\textstyle
f'(0)=\lim_{x\to0}\frac{f(x)}{x}=0
}
$$ 
because for every $x\in (b_n,a_n)$ we have (Exercise~\ref{ex_procSe}) 
$${\textstyle
\big|\frac{f(x)}{x}\big|\le\frac{a_n - b_n}{b_n}
\to0\ \  (n\to\infty)\,. 
}
$$
Hence $D(f)=N$. Since $0\in L(N)$,  $f'(0)=0$ and $f'(x)=1$ 
for every $x\in N\setminus\{0\}$, the derivative $f'$ is a~discontinuous function.
\kduk
\vspace{-3mm}
\begin{exer}\label{ex_procSe}
Why $\lim\frac{a_n - b_n}{b_n}=0$ as $n\to\infty$? 
\end{exer}

\section[${}^c$Standard and limit tangents]{Standard and limit tangents}\label{sec_tecny}

We define two kinds of tangent lines.

\medskip\noindent
{\em $\bullet$ Standard tangents. }Tangent at a~point means a~finite derivative at that point.

\begin{defi}[tangents]\label{def_tangentLine}
Let\index{tangent@tangent (line)!standard|emph} 
$f\in\mathcal{F}(M)$, $b\in M\cap L(M)$ and $f'(b)\in\R$. 
The \underline{tangent} line to the graph $G_f$ at the point $\langle b,f(b)\rangle$ is the line 
$$
\ell=\{\langle x,\,y\rangle\in\R^2\cc\;
x\in\R,\,y=f'(b)(x-b)+f(b)\}\ \ (\sus\R^2)\,.
$$
\end{defi}
So $\ell$ goes through 
the plane point $\langle b,f(b)\rangle$. 
We rewrite the equation as 
$$
y=f'(b)x+f(b)-f'(b)b
$$ 
and see that $\ell$ 
has the slope $f'(b)$ ($\in\R$). 
We view the tangent also as the function
$\ell\in\mathcal{F}(\R)$ given by
$$
\ell(x)=f'(b)(x-b)+f(b)=f'(b)x+f(b)-f'(b)b\,.
$$

\begin{exer}\label{ex_jednStTec}
The function $\ell(x)$ is the only linear polynomial such that
$$
f(x)=\ell(x)+o(x-b)\ \ (x\to b)\,.
$$
\end{exer}

\begin{exer}\label{ex_tecOdmoc}
Let $f(x)=\sqrt{x}$.
Determine tangents at the points $\langle a,\sqrt{a}\rangle$.   
\end{exer}

\noindent
{\em $\bullet$ Touching lines are tangents. }We show that every line touching the graph 
at a~single point where the derivative exists is the tangent
at that point.

\begin{thm}[touching lines]\label{thm_suppLine}
Let\index{theorem!touching lines|emph} 
$f$ be in $\mathcal{F}(M)$, $b\in M\cap L^{\mathrm{TS}}(M)$ and 
$f'(b)\in\R^*$. If the linear function
$$
l(x)=sx+t 
$$ 
has value $l(b)=f(b)$ and if
$f(x)\ge l(x)$ for every $x\in M$, \underline{then} 
$l(x)$ is the tangent to $G_f$ at $\langle b,\,f(b)\rangle$.
\end{thm}
\duk
Let $g(x)= f(x)-
l(x)$ ($\in\mathcal{F}(M)$).
By 
Exercise~\ref{ex_whyThis},  $g'(b)=f'(b)-s$. If $g'(b)\ne0$, by
Theorem~\ref{thm_priznakExtr} there is a~point $c\in M$ near $b$ such that 
$$
\text{$g(c)<g(b)=0$ and hence $f(c)<l(c)$}\,, 
$$
in contradiction with the assumption. Thus $g'(b)=0$ and
$f'(b)=s$. It follows that $l(x)$
is the tangent line to the graph $G_f$ at the point $\langle b,f(b)\rangle$. 
\kduk
\vspace{-3mm}
\begin{exer}\label{ex_whyThis}
Show that $g'(b)=f'(b)-s$.
\end{exer}

\begin{exer}\label{ex_modifikace}
Modify the theorem for $l(x)$ touching $G_f$ from above.
\end{exer}

\noindent
{\em $\bullet$ Non-vertical lines. }\underline{Non-vertical lines\index{non-vertical lines, $\mathcal{N}$|emph}} in the plane are the lines 
$$
\ell=\{\langle x,\,y\rangle\in\R^2\cc\;x,\,y\in\R,\,y=sx+t\}\ \ (\sus\R^2)\,,
$$
where $s,t\in\R$. Note that $s$ and $t$ are unique: to different pairs $s,t$ correspond different lines $\ell$. We call $s$ the \underline{slope\index{slope|emph}\index{non-vertical lines, 
$\mathcal{N}$!slope of|emph}} of $\ell$. We denote the set of 
non-vertical lines by 
$\mathcal{N}$ ($\sus\mathcal{P}(\R^2)$).\label{Ncal}

\begin{exer}\label{ex_bijNesPri}The function\underline{\index{non-vertical lines, 
$\mathcal{N}$!parametrization by $\R^2$|emph}}
$p\cc\mathcal{N}\to\R^2$, given by 
$$
p(\ell)=\langle p_1(\ell),\,p_2(\ell)\rangle
=\langle s,\,t\rangle\,,
$$
with $s$ and $t$ as in $\ell=\ell(x)=sx+t$, is a~bijection.
\end{exer}

\begin{defi}[limits in $\mathcal{N}$]\label{def_limNescPrim}
\label{def_limvN}
Let $(\ell_n)\sus\mathcal{N}$, $\ell\in\mathcal{N}$ and let the function $p(\cdot)=\langle p_1(\cdot),p_2(\cdot)\rangle$ be as in Exercise~\ref{ex_bijNesPri}. If 
$$
\lim_{n\to\infty} p_1(\ell_n)=p_1(\ell)\,\text{ and }\,\lim_{n\to\infty} 
p_2(\ell_n)=p_2(\ell)\,,
$$ 
we say that the line $\ell$ is the \underline{limit of lines\index{limit of lines@limit of (non-vertical) lines|emph}} $\ell_n$, and write $\mathrm{Lim}\,\ell_n=\ell$\underline{\index{non-vertical lines, $\mathcal{N}$!limits of|emph}}.\label{limLin} 
\end{defi}

\begin{exer}\label{ex_defKappa}
Let $A=\langle a,b\rangle$ and $A'=\langle a',b'\rangle$ be in $\R^2$ and $a\ne a'$. There is a~unique line $\ell\in\mathcal{N}$ such that
$A,A'\in\ell$. It has slope $\frac{b'-b}{a'-a}$.    
\end{exer}
We denote this unique non-vertical line going through $A$ and $A'$ by
$$
\text{$\kappa(A,\,A')$\underline{\index{non-vertical lines, 
$\mathcal{N}$!kappa@$\kappa(A,A')$|emph}} or by $\kappa(a,\,b,\,a',\,b')$}\,.\label{kappa}
$$
If $f\in\mathcal{R}$ and $A,A'\in G_f$ with $A\ne A'$, we call the line $\kappa(A,A')$ 
a~\underline{secant\index{secant line|emph}\index{graph of 
$f$!secant|emph}} line of 
the graph $G_f$.

\medskip\noindent
{\em $\bullet$ Limit tangents. }We often read in textbooks and lecture notes that tangents
are limits of secants, but the details of this limiting process are never revealed. We explain it here. 

\begin{defi}[limit tangents]\label{def_limiTecn}
Let $f\in\mathcal{F}(M)$, $b\in M\cap L(M)$ and let $\ell\in\mathcal{N}$. The line $\ell$ is 
a~\underline{limit tangent\index{tangent@tangent
(line)!limit|emph}} to the graph $G_f$ at the point 
$\langle b,f(b)\rangle$, if for every sequence $(a_n)\sus M\setminus\{b\}$ with $\lim a_n=b$ we have, by Definition~\ref{def_limvN}, the limit
$$
\mathrm{Lim}\,\kappa(b,\,f(b),
\,a_n,\,f(a_n))=\ell\,.
$$
\end{defi}
This definition of tangents does not need derivatives.

\begin{exer}\label{ex_bodNaLT}
If $\ell$ is a~limit tangent to $G_f$ at $\langle b,f(b)\rangle$ then $\langle b,f(b)\rangle\in\ell$.
\end{exer}
We prove that the two definitions of tangents are logically 
equivalent.

\begin{thm}[standard and limit tangents]\label{thm_onTangents}
Let\index{theorem!standard and limit tangents|emph} $f\in\mathcal{F}(M)$, 
$b\in M\cap L(M)$
and let $\ell\in\mathcal{N}$. \underline{Then} $\ell$ is a~tangent 
to $G_f$ at $\langle b,f(b)\rangle$ by Definition~\ref{def_tangentLine} $\iff$ $\ell$ is a~limit tangent to $G_f$ at $\langle b,f(b)\rangle$ by Definition~\ref{def_limiTecn}.
\end{thm}
\duk
We prove the implication $\Rightarrow$. Let $f'(b)\in\R$,
$$
\ell(x)=f'(b)x+f(b)-f'(b)f(b) 
$$
and $(a_n)\sus M\setminus\{b\}$ be a~sequence with $\lim a_n=b$. We denote 
$${\textstyle
c_n=\frac{f(a_n)-f(b)}{a_n-b}\,.
}
$$
The secant $\kappa_n=\kappa_n(x)$ through the points $\langle 
b,f(b)\rangle$ and $\langle a_n,f(a_n)\rangle$ is 
$$
\kappa_n(x)=\kappa(b,\,f(b),\,a_n,\,f(a_n))(x)=
c_n(x-b)+f(b)=c_nx+f(b)-c_nb\,.
$$
By HDD we have $\lim c_n=f'(b)$. Thus also
$\lim(f(b)-c_nb)=f(b)-f'(b)b$. By Definition~\ref{def_limvN}, $\mathrm{Lim}\,\kappa_n=\ell$ .

We prove the implication $\Leftarrow$. Let $\ell(x)=sx+t$ and let $(a_n)$, $c_n$ and $\kappa_n$ be as above.
We assume that for every such sequence $(a_n)$ we have 
$\mathrm{Lim}\,\kappa_n=\ell$ 
by Definition~\ref{def_limvN}. So always $\lim c_n=s$, $\langle b,f(b)\rangle\in\ell$ by Exercise~\ref{ex_bodNaLT}, and 
$$
\lim (f(b)-c_nb)=f(b)-sb=t\,. 
$$
By HDD we have $s=f'(b)$. Hence $t=f(b)-f'(b)b$. So $\ell$ is a~tangent to $G_f$ at $\langle b,f(b)\rangle$ by Definition~\ref{def_tangentLine}.
\kduk

\noindent
{\em $\bullet$ Tangents without points of contact. }We show that 
the tangent at the point $B\in G_f$ is the limit of secants going
through pairs of points in $G_f$ that are separated by $B$ and 
converge to $B$.

\begin{exer}\label{ex_prvNex}
Prove the next lemma.    
\end{exer}

\begin{lemma}[convex combination]\label{lem_konvKomb}
Let  $s,v\in\R$ be positive, $\al=\frac{s}{s+v}$ 
and $\be=\frac{v}{s+v}$, so that $\al,\be>0$ and $\al+\be=1$. 
\underline{Then} for every $r,t\in\R$ we have 
$$
\frac{r+t}{s+v}=\al\cdot
\frac{r}{s}+\be\cdot\frac{t}{v}\,\,.
$$  
\end{lemma}

\begin{thm}[tangents and secants]\label{thm_NonmainSec}
Let\index{theorem!tangents and secants|emph} $f\in\mathcal{F}(M)$, $b\in L^{\mathrm{TS}}(M)\setminus M$, 
and $\ell\in\mathcal{N}$.  
\underline{Then} two claims are equivalent.
\begin{enumerate}
\item One can extend $f$ to $M\cup
\{b\}$ so that $\ell$ is a~tangent to $G_f$ at $\langle b,f(b)\rangle$.
\item 
For every two sequences $(x_n),(y_n)\sus M$ such that $\lim x_n=\lim y_n=b$ and $x_n<b<y_n$ for every $n$ we have the limit
$$
\mathrm{Lim}\,\kappa(x_n,\,f(x_n),\,y_n,\,f(y_n))=\ell\,.
$$ 
\end{enumerate}
\end{thm}
\duk
We prove the implication $1\Rightarrow2$. We assume that $f$ has been extended to $b$ by a~value $f(b)$, that $f'(b)\in\R$, and that 
$$
\ell=\ell(x)=f'(b)x+f(b)-f'(b)b\,.  
$$
Let $(x_n)\sus M$ with $x_n<b$ and $(y_n)\sus 
M$ with $b<y_n$ be sequences converging to $b$. We denote
$$
\text{$r_n= f(b)-
f(x_n)$, $s_n= b-x_n$, $t_n= f(y_n)-
f(b)$ and $v_n= y_n-b$}\,. 
$$
We use Lemma~\ref{lem_konvKomb} and write the slope $u_n$ of the secant 
$$
\ell_n=\kappa(x_n,\,f(x_n),\,y_n,\,f(y_n))
$$
of $G_f$ as the convex combination
$${\textstyle
u_n=\frac{f(y_n)-f(x_n)}{y_n-x_n}=
\frac{r_n+t_n}{s_n+v_n}=\al_n\cdot\frac{r_n}{s_n}+\be_n\cdot\frac{t_n}{v_n}
}
$$
of the slopes $\frac{r_n}{s_n}$ and $\frac{t_n}{v_n}$ of the respective secants 
$$
\text{$\kappa(x_n,\,f(x_n),\,b,\,f(b))$ and $\kappa(b,\,f(b),\,y_n,\,f(y_n))$}\,.
$$
Since $\lim \frac{r_n}{s_n}=\lim\frac{t_n}{v_n}=f'(b)$,  
by Corollary~\ref{cor_dvaStraz} also $\lim u_n=f'(b)$. Since 
$$
\ell_n(x)=u_nx+f(x_n)-u_nx_n\,,
$$
$\lim u_n=f'(b)$, $\lim x_n=b$ and $\lim f(x_n)=f(b)$ 
($f$ is 
continuous at $b$ due to $f'(b)\in\R$), we have by Definition~\ref{def_limvN} that $\mathrm{Lim}\,\ell_n=\ell$.

We prove the implication 
$\neg1\Rightarrow\neg2$. We assume that
$f$ cannot be extended to $b$ by any value $f(b)$ so that 
$\ell$ is tangent 
to $G_f$ at $\langle b,f(b)\rangle$. This means that if we take the unique number $f(b)\in\R$ such that 
$\langle b,\,f(b)\rangle\in\ell$ 
and if $s$ is the slope of $\ell$, then it is not true that 
$$
\lim_{x\to b}{\textstyle\frac{f(x)-f(b)}{x-b}=s\,. 
}
$$
We then obtain pairs of points in $M$ such that their components converge to $b$ from opposite sides, and it is does not hold that the 
limit of the corresponding secants is $\ell$.

The first case is that the extended function $f$ is not 
continuous at $b$. By Exercise~\ref{ex_oKL}, there exist sequences $(x_n)\sus M$ with $x_n<b$ and 
$(y_n)\sus M$  with $b<y_n$ such that $\lim x_n=\lim y_n=b$, $\lim 
f(x_n)=K$, $\lim f(y_n)=L$, and it is not true that 
$$
K=L=f(b)\,. 
$$
If $K\ne 
L$ then the slopes of the secants 
$$
\ell_n=\kappa(x_n,\,f(x_n),\,y_n,\,f(y_n))
$$
go to $\pm\infty$ and it is not true that $\mathrm{Lim}\,\ell_n=l$. 
If $K=L\ne f(b)$ then the intersections of the secants $\ell_n$ with the 
vertical line $x=b$ converge to a~(possibly infinite) point different 
from $\langle b,f(b)\rangle$. By Exercise~\ref{ex_limPrnejde}, the 
limit $\mathrm{Lim}\,\ell_n$, if it exists, cannot be a~line going 
through $\langle b,f(b)\rangle$. Again, it is not true that $\mathrm{Lim}\,\ell_n=l$.

The second case is that the extended function $f$ is continuous at $b$, but it does not hold that
$$
\lim_{x\to b}{\textstyle\frac{f(x)-f(b)}{x-b}=s\,.
}
$$
Then there is 
an $A\in\R^*\setminus\{s\}$ 
and a~sequence $(x_n)\sus M\setminus\{b\}$ lying on one side of $b$ such that 
$${\textstyle
\text{$\lim x_n=b$ 
and $\lim\frac{f(x_n)-f(b)}{x_n-b}=A$}\,. 
}
$$
We may assume that 
$x_n<b$ for every $n$, the case when always $x_n>b$ is similar. We take any sequence 
$(y_n)\sus M$ with $b<y_n$ and $\lim y_n=b$. Then $\lim f(y_n)=f(b)$ and by Exercise~\ref{ex_vybPodp} we can choose 
from 
$(y_n)$ a~subsequence $(y_{m_n})$ such that 
$${\textstyle
\lim_{n\to\infty}\frac{f(x_n)-f(y_{m_n})}{x_n-y_{m_n}}=A\,.
}
$$
Since $A\ne s$, it is not true that $\mathrm{Lim}\,\ell_n=l$.
\kduk

\noindent
Exercises~\ref{ex_oKL}--\ref{ex_vybPodp} are lemmas for the proof. Exercise~\ref{ex_uloNatecnu6} 
shows that secants through pairs of points that converge to $B\in G_f$ 
from the same side need not converge to the tangent at $B$.

\begin{exer}\label{ex_oKL}
Let $f\in\mathcal{F}(M)$ and $b\in M\cap L^{\mathrm{TS}}(M)$. If $f$ is not continuous at $b$, then for some sequences 
$(x_n),(y_n)\sus M$ with $x_n<b<y_n$ and $\lim x_n=\lim y_n=b$ we have $\lim f(x_n)=K$, 
$\lim f(y_n)=L$, but $K\ne f(b)$ or $L\ne f(b)$.    
\end{exer}

\begin{exer}\label{ex_limPrnejde}
Let $(\ell_n)\sus\mathcal{N}$, $\ell\in\mathcal{N}$, $\langle b,c\rangle\in\ell$, $\mathrm{Lim}\,\ell_n=\ell$ and 
$$
(x=b)\cap\ell_n=
\{\langle b,\,c_n\rangle\}\,. 
$$
Then $\lim c_n=c$. 
\end{exer}

\begin{exer}\label{ex_vybPodp}
Let sequences $(x_n)$, $(y_n)$, $(z_n)$ and $(u_n)$ be such that 
$\lim x_n=\lim 
z_n=b$, $x_n\ne b$, $\lim y_n=\lim u_n=c$ and $\lim
\frac{y_n-c}{x_n-
b}=A$. Then for some sequence $(m_n)\sus\N$ we have $\lim
\frac{y_n-u_{m_n}}{x_n-z_{m_n}}=A$.   
\end{exer}

\begin{exer}\label{ex_uloNatecnu6}
Find a~function $f\in\mathcal{F}(M)$ with tangent $\ell$ at $\langle b,f(b)\rangle$ and two sequences 
$(x_n),(y_n)\sus M$ such that $b<x_n<y_n$, $\lim x_n=\lim 
y_n=b$, but it does not hold that
$$
\mathrm{Lim}\,\kappa(x_n,\,f(x_n),\,y_n,\,f(y_n))=\ell\,.
$$ 
\end{exer}

\section[${}^c$Arithmetic of derivatives]{Arithmetic of derivatives}\label{sec_aritmDeri}

Global derivative is a~unary operation on $\mathcal{R}$. We 
describe its interactions with the operations of addition, 
multiplication, and division. In the next section, we treat composition and
inverse. We consider both point-wise 
and global formulas. Point-wise formulas involve finite and infinite 
values of derivatives at points. Global formulas involve only finite values of global 
derivatives.

\medskip\noindent
{\em $\bullet$ Sums. }We differentiate sums of functions.

\begin{thm}[$(f+g)'$]\label{thm_nasob}
Let\index{theorem!derivatives of sums|emph} 
$f,g\in\mathcal{R}$ and $M= M(f)\cap M(g)$.
\begin{enumerate}
\item If $b\in M\cap L(M)$, $f'(b)=K$, $g'(b)=L$, and $K+L$ is not an indefinite expression, \underline{then} 
$$
(f+g)'(b)=K+L\,.
$$
\item The function $(f'+g')\,|\,L(M)$ 
is a~subfunction of the function $(f+g)'$.
\end{enumerate}
\end{thm}
\duk
1. Let $h= f+g$, then 
$b\in M(h)\cap L(M(h))$. By Theorem~\ref{thm_AritLimFce}, 
\begin{eqnarray*}
h'(b)&=&\lim_{x\to b}{\textstyle
\frac{h(x)-h(b)}{x-b}
=\lim_{x\to b}
\big(\frac{f(x)-f(b)}{x-b}+\frac{g(x)-g(b)}{x-b}\big)}\\
&=&\lim_{x\to b}{\textstyle\frac{f(x)-f(b)}{x-b}}+\lim_{x\to b}{\textstyle\frac{g(x)-g(b)}{x-b}=f'(b)+g'(b)=K+L}\,.
\end{eqnarray*}

2. Let $h=(f'+g')\,|\,L(M)$ and 
$$
c\in M(h)=D(f)\cap D(g)\cap L(M)\,. 
$$
By part~1,  $(f+g)'(c)=f'(c)+g'(c)=h(c)$. It follows that $h$ is a~subfunction of $(f+g)'$.
\kduk
\vspace{-3mm}
\begin{exer}\label{ex_sameForMinus}
Adapt the theorem for the difference $f-g$.    
\end{exer}

\noindent
Item~1 of the theorem is standard. Item~2 is new. In our approach to derivatives, the formula 
$$
(f+g)'=f'+g'
$$ 
in general does not hold. Let $f= k_0\,|\,(-\infty,0]$ and
$g= k_0\,|\,[0,+\infty)$. Then $M=M(f)\cap M(g)=\{0\}$ and 
$$
f'+g'=k_0\,|\,(-\infty,0]+k_0\,|\,[0,+\infty)=k_0\,|\,\{0\}
\ne\emptyset_f=(k_0\,|\,\{0\})'=(f+g)'\,.
$$
The restriction to $L(M)$ therefore cannot be
omitted. 

Let $f(x)=|x|$ and $g(x)=-|x|$ ($\in\mathcal{F}(\R)$). Then $M=M(f)\cap M(g)=\R$, $L(M)=\R^*$, 
$$
(f'+g')\,|\,L(M)=k_0\,|\,(\R\setminus\{0\})\,
\text{ and }\,(f+g)'=(k_0)'=k_0\,.
$$
Thus $(f'+g')\,|\,L(M)$ may be a~proper subfunction of $(f+g)'$.

In an application of item~1 of Theorem~\ref{thm_nasob} we compute 
a~derivative. Let $f(x)=\sgn(x)$, 
$g(x)=\sqrt{x}$ and $b=0$. Then 
$M=M(f)\cap M(g)=[0,+\infty)$ and
$$
\big(\sgn(x)+\sqrt{x}\big)'(0)=\sgn'(0)+(\sqrt{x})'(0)=+\infty+(+\infty)=+\infty\,.
$$

\begin{exer}\label{ex_sgnMinusOdm}
What is $\big(\sgn(x)-\sqrt{x}\big)'(0)$?    
\end{exer}

We find assumptions for the validity of $(f+g)'=f'+g'$. We need
it for the proof of Theorem~\ref{thm_deriSEF}.

\begin{cor}[$(f+g)'=f'+g'$]\label{cor_deriSumStan}
Let $f,g\in\mathcal{R}$ and $M= M(f)\cap M(g)$. If $D(f)=M(f)$, 
$D(g)=M(g)$ and $M\sus L(M)$, \underline{then} 
$$
(f+g)'=f'+g'\,.
$$
\end{cor}
\duk
Let $h=(f'+g')\,|\,L(M)$. From the 
assumptions it follows that $M(h)=M$ and that $h=f'+g'$. Since 
$M((f+g)')\sus M$, item~2 of 
Theorem~\ref{thm_nasob} gives that $f'+g'=h=(f+g)'$. 
\kduk

\noindent
{\em $\bullet$ Products. }We differentiate products of functions.

\begin{thm}[Two Leibniz formulas]\label{thm_LeibnizFmle}
Let\index{theorem!Leibniz product formulas|emph}\index{Leibniz 
formulas|emph}\index{Leibniz, Gottfried W.} 
$f$ and $g$ be in $\mathcal{R}$ and $M= M(f)\cap M(g)$. 
\begin{enumerate}
\item \underline{Local Leibniz formula\index{Leibniz 
formulas!local one|emph}}. If $b\in M\cap L(M)$, $f'(b)=K$, $g'(b)=L$, one of $f$ and $g$ is continuous at $b$ and if $K\cdot g(b)+f(b)\cdot L$ is 
defined, \underline{then} 
$$
(fg)'(b)=Kg(b)+f(b)L\,.
$$
\item \underline{Global Leibniz formula\index{Leibniz 
formulas!global one|emph}}. The function $(f'g+fg')\,|\,L(M)$ is a~subfunction of the function $(fg)'$.
\end{enumerate}
\end{thm}
\duk
1. Let $h= fg$ and $g$ be continuous at $b$. The case when $f$ is continuous 
at $b$ is resolved in Exercise~\ref{ex_spojEf}. Then $b\in M(h)\cap L(M(h))$ and, 
by the assumptions, by Theorem~\ref{thm_AritLimFce} and by 
Proposition~\ref{prop_spojVbLimitou}, 
we have
\begin{eqnarray*}
h'(b)&=&\lim_{x\to b}{\textstyle
\frac{f(x)g(x)-f(b)g(b)}{x-b}}
=\lim_{x\to b}{\textstyle\frac{(f(x)-
f(b))g(x)+f(b)(g(x)-g(b))}{x-b}}\\
&=&\lim_{x\to b}{\textstyle\frac{f(x)-f(b)}{x-b}}\cdot\lim_{x\to b}
g(x)+f(b)\lim_{x\to b}{\textstyle\frac{g(x)-g(b)}{x-b}=Kg(b)+f(b)L}\,.
\end{eqnarray*}

2. Let $h=(f'g+fg')\,|\,L(M)$ and $c\in M(h)$. Then 
$$
c\in D(f)\cap D(g)\cap L(M)\,, 
$$
$g$ is continuous at $c$ because $g'(c)\in\R$, and by the 
first part we have that $(fg)'(c)=f'(c)g(c)+f(c)g'(c)=h(c)$.
Hence $h$  is a~subfunction of $(fg)'$.
\kduk

\noindent
Like in the case of sum, it is not hard to produce examples showing 
that the restriction to $L(M)$
cannot be omitted and that $(f'g+fg')\,|\,L(M)$ may be 
a~proper subfunction of $(fg)'$.

\begin{exer}\label{ex_spojEf}
Solve quickly the case when $f$ is  continuous at $b$.
\end{exer}

We show that the assumption of continuity of $f$ or $g$ at $b$ is 
substantial.

\begin{exer}\label{ex_uloNaLeiFor}
We define functions $f,g\in\mathcal{F}(\R)$ by $f(0)=-\frac{1}{2}$, 
$g(0)=\frac{1}{2}$ and for $a\ne0$ by
$$
f(a)=\sgn\,a\,\text{ and }\,
g(a)=-\sgn\,a\,.
$$
Show that $(fg)'(0)$ does not exist but that at $0$ the 
right-hand side of the local Leibniz formula is $+\infty$. 
\end{exer}

For future use we obtain a~simple form of the global Leibniz formula.

\begin{cor}[$(fg)'=f'g+fg'$]\label{cor_Leibform}
Let $f,g\in\mathcal{R}$ and $M= M(f)\cap M(g)$. If $D(f)=M(f)$, 
$D(g)=M(g)$ and $M\sus L(M)$, \underline{then} 
$$
(fg)'=f'g+fg'\,.
$$
\end{cor}
\duk
Let $h=(f'g+fg')\,|\,L(M)$. 
From the assumptions, it follows that $M(h)=M$ and
$h=f'g+fg'$. Since $M((fg)')\sus M$, part~2 of 
Theorem~\ref{thm_LeibnizFmle} gives that $f'g+fg'=h=(fg)'$.
\kduk

\noindent
{\em $\bullet$ Division. }We differentiate ratios of functions.

\begin{thm}[$(\frac{f}{g})')$]\label{thm_deriPodilu} Suppose\index{theorem!derivatives of ratios|emph} that 
$f,g\in\mathcal{R}$ and $M=  M(f)\cap M(g)\setminus Z(g)$. 
\begin{enumerate}
\item If $b\in M\cap L(M)$, $f'(b)=K$, $g'(b)=L$, $g$ is continuous at $b$ and if 
$$
\frac{K\cdot g(b)-f(b)\cdot L}{g(b)^2}
$$ 
is defined, \underline{then}
$$
\bigg(\frac{f}{g}\bigg)'(b)=\frac{Kg(b)-f(b)L}{g(b)^2}\,.
$$
\item The function $\frac{f'g-fg'}{g^2}\,|\,L(M)$ is a~subfunction of the function $(\frac{f}{g})'$.
\end{enumerate}
\end{thm}
\duk
1. Let $h=\frac{f}{g}$ and $g$ be continuous at $b$.
Then $b\in M(h)\cap L(M(h))$ and
$$
{\textstyle
h'(b)=\lim_{x\to b}\frac{f(x)/g(x)\,-\,f(b)/g(b)}{x\,-\,b}=\lim_{x\to b}\frac{f(x)g(b)-f(b)g(b)+f(b)g(b)-f(b)g(x)}{g(x)g(b)(x-b)}\,.
}
$$
Due to the assumptions, Proposition~\ref{prop_spojVbLimitou} and Theorem~\ref{thm_AritLimFce},  this equals
\begin{eqnarray*}
&&{\textstyle
\lim_{x\to b}\frac{f(x)-f(b)}{x-b}\lim_{x\to b}\frac{g(b)}{g(x)g(b)}-
\lim_{x\to b}\frac{f(b)}{g(x)g(b)}\lim_{x\to b}\frac{g(x)-
g(b)}{x-b}}\\
&&{\textstyle=\frac{f'(b)g(b)-f(b)g'(b)}{g(b)^2}\,.
}
\end{eqnarray*}

2. Let $h=\frac{f'g-fg'}{g^2}\,|\,L(M)$ and $c\in M(h)$. Then 
$$
c\in D(f)\cap D(g)\cap 
L(M)\setminus Z(g)\,, 
$$
$g$ is continuous at $c$, because $g'(c)\in\R$, and by the first part we have that $(\frac{f}{g})'(c)=\frac{f'(c)g(c)-
f(c)g'(c)}{g(c)^2}=h(c)$. Hence $h$  is a~subfunction of $(\frac{f}{g})'$.
\kduk

\noindent
Again, the restriction to $L(M)$ cannot in general be omitted and  
$\frac{f'g-fg'}{g^2}\,|\,L(M)$ can be a~proper restriction of $(\frac{f}{g})'$.

\begin{exer}\label{ex_uloNaDerPodilu}
Show, as in Exercise~\ref{ex_uloNaLeiFor},  that the assumption of 
continuity of $g$ at $b$ cannot be omitted.
\end{exer}

For later use we again obtain a~corollary with a~simple form of 
the formula for derivatives of ratios.

\begin{cor}[$(\frac{f}{g})'=
\frac{f'g-f'g}{g^2}$]\label{cor_deriRat}
Let $f,g\in\mathcal{R}$ and $M= M(f)\cap M(g)$. If $D(f)=M(f)$, 
$D(g)=M(g)$ and $M\sus L(M)$ \underline{then} 
$$
\bigg(\frac{f}{g}\bigg)'=\frac{f'g-f'g}{g^2}\,.
$$
\end{cor}
\duk
Let $h=
\frac{f'g-fg'}{g^2}\,|\,L(M)$. From the assumptions it 
follows that $M(h)=M\setminus Z(g)$ 
and $h=\frac{f'g-fg'}{g^2}$. Since $M((\frac{f}{g})')\sus 
M\setminus Z(g)$, part~2 of 
Theorem~\ref{thm_deriPodilu} gives that $\frac{f'g-fg'}{g^2}=h=
(\frac{f}{g})'$.
\kduk

\section[${}^c$Composite and inverse functions]{Composite and inverse functions}\label{sec_deriSlozInve}

Let $f,g\in\mathcal{R}$. Recall that the composite function 
$$
f(g)\cc M(f(g))\to\R,\,\text{ where }\,f(g)(x)= f(g(x))\,,
$$ 
has the domain
$$
M(f(g))=\{x\in M(g)\cc\;g(x)\in M(f)\}=g^{-1}[M(f)]\,.
$$ 
Thus $M(f(g))\sus M(g)$ and this inclusion may be proper. Any 
injection $f\in\mathcal{R}$ has the inverse
$$
f^{-1}\cc f[M(f)]\to\R\,\text{ where }\,f^{-1}(y)=x\iff f(x)=y\,.
$$ 
For non-injective $f$ the inverse is not defined.

\medskip\noindent
{\em $\bullet$ Composites. }We differentiate composite 
functions.

\begin{thm}[$(f(g))'$]\label{thm_DerSlozFce}
Let\index{theorem!derivatives of compositions|emph} 
$f,g\in\mathcal{R}$ and $M= M(f(g))$. 
\begin{enumerate}
\item If $b\in M\cap L(M)$, $f'(g(b))=K$, $g'(b)=L$, $g$ is 
continuous at $b$ 
and  if $K\cdot L\ne0\cdot(+\infty),(+\infty)\cdot0$, \underline{then} 
$$
f(g)'(b)=KL\,.
$$
\item The function 
$(f'(g)\cdot g')\,|\,L(M)$ is a~subfunction of the function $(f(g))'$.
\end{enumerate} 
\end{thm}
\duk
1. Let $(a_n)\sus M\setminus
\{b\}$ be any sequence with
$\lim a_n=b$. By the assumption, $\lim g(a_n)=g(b)$. We partition $(a_n)$ 
in two (possibly finite or empty) subsequences $(b_n)$ and $(c_n)$:  for every $n$ we have $g(b_n)=g(b)$ and $g(c_n)\ne g(b)$. We show that if $(x_n)=(b_n)$ or
$(x_n)=(c_n)$ is an infinite sequence, then 
$${\textstyle
\lim_{n\to\infty}\frac{f(g)(x_n)-f(g)(b)}{x_n-b}=f'(g(b))\cdot g'(b)\,.
}
$$
Then by Theorem~\ref{thm_finManyBl} it follows that $\lim\frac{f(g)(a_n)-f(g)(b)}{a_n-b}=f'(g(b))\cdot g'(b)$. By HDD, 
$$
f(g)'(b)=f'(g(b))\cdot g'(b)\,.
$$

Let the sequence $(b_n)$ be infinite. Then $\lim b_n=b$. By HDD we have
$${\textstyle
g'(b)=\lim\frac{g(b_n)-g(b)}{b_n-b}=\lim\frac{g(b)-g(b)}{b_n-b}=\lim 0=0\,.
}
$$
Thus
$${\textstyle
\lim\frac{f(g)(b_n)-f(g)(b)}{b_n-b}=
\lim\frac{f(g(b))-f(g(b))}{b_n-b}=0=f'(g(b))\cdot 0
=f'(g(b))\cdot g'(b)\,.
}
$$

Let the sequence $(c_n)$ be infinite. Then $\lim c_n=b$, $\lim g(c_n)=g(b)$
and by Theorem~\ref{thm_AritLimFce} and HDD we have 
\begin{eqnarray*}
&&{\textstyle
\lim\frac{f(g)(c_n)-f(g)(b)}{c_n-b}=
\lim\frac{f(g(c_n))-f(g(b))}{c_n-b}=
}\\
&&{\textstyle
=\lim\big(\frac{f(g(c_n))-f(g(b))}{g(c_n)-g(b)}\cdot
\frac{g(c_n)-g(b)}{c_n-b}\big)=
\lim\frac{f(g(c_n))-f(g(b))}{g(c_n)-g(b)}\cdot\lim\frac{g(c_n)-g(b)}{c_n-
b}=}\\
&&=f'(g(b))\cdot g'(b)\,.
\end{eqnarray*}

2. Let $h=(f'(g)\cdot g')\,|\,L(M)$ and $c\in M(h)$. Then 
$$
c\in M(f'(g))\cap D(g)\cap  
L(M)\,. 
$$
The functuon $g$ is continuous at $c$, because $g'(c)\in\R$, and by the 
first part $(f(g))'(c)=f'(g(c))\cdot g'(c)=h(c)$. Hence $h$ is a~subfunction of $(f(g))'$.
\kduk

\noindent
Again, the restriction to $L(M)$ cannot be omitted and $(f'(g)\cdot 
g')\,|\,L(M)$ can be a~proper restriction of $f(g)'$.

\begin{exer}\label{ex_bezSpoNrpl}
Part~1 does not hold without the continuity of $g$ at $b$. 
\end{exer}

In a~corollary, we simplify the formula for the derivative of composites.

\begin{cor}[$(f(g))'=f'(g)\cdot g'$]\label{cor_derSlozF}
Let $f,g\in\mathcal{R}$ and $M= M(f(g))$. If $D(f)=M(f)$, $D(g)=M(g)$ and 
$M\sus L(M)$ \underline{then} 
$$
(f(g))'=f'(g)\cdot g'\,.
$$
\end{cor}
\duk
Let $h=(f'(g)\cdot g')\,|\,L(M)$. From the assumptions on $f$ and $g$ it 
follows that $M(h)=M$ and $h=f'(g)\cdot g'$. Since $M((f(g))')\sus M$, part~2 of 
Theorem~\ref{thm_LeibnizFmle} gives that $$
f'(g)\cdot g'=h=(f(g))'\,.
$$
\kduk

\noindent
{\em $\bullet$ Inverses. }We take derivatives of inverse functions. A~function 
$f\in\mathcal{F}(M)$ \underline{increases\index{function!increases at a point|emph}}, 
respectively \underline{decreases,  at a point} $b\in 
M$\index{function!decreases at a point|emph} if there exists a~$\de$ such that for every $x$ and $x'$ 
with $b-\de<x<b<x'<b+\de$ we have 
$$
\text{$f(x)<f(b)<f(x')$, respectively $f(x)>f(b)>f(x')$}\,.
$$

\begin{thm}[$(f^{-1})'$]\label{thm_DeriInvfce}
Let\index{theorem!derivatives of 
inverses|emph} 
$f\in\mathcal{F}(M)$ be injective, $b$ be in $M\cap L(M)$, $f'(b)$ in $\R^*$, and $f^{-1}$ be continuous at $c= f(b)$. 
\begin{enumerate}
\item If $f'(b)\in\R\setminus\{0\}$, \underline{then} 
$$
\big(f^{-1}\big)'(c)=\frac{1}{f'(b)}=\frac{1}{f'(f^{-1}(c))}\,.
$$
\item If $f'(b)=0$ and $f$ increases, respectively decreases,
at $b$, \underline{then} 
$$
\text{$(f^{-1})'(c)=+\infty$, respectively $-\infty$}\,.
$$
\item If $f'(b)=\pm\infty$ and $c\in L(f[M])$, \underline{then} 
$$
(f^{-1})'(c)=0\,.
$$
\end{enumerate}
\end{thm}
\duk
We take a~sequence $(b_n)\sus f[M]\setminus
\{c\}$ with $\lim b_n=c$ and set $a_n= f^{-1}(b_n)$. Then 
$(a_n)\sus M\setminus\{b\}$ and by the continuity of $f^{-1}$ at $c$ we have $\lim a_n=b$.

1. Let $f'(b)$ be finite and nonzero. Then $c\in L(f[M])$ by Proposition~\ref{prop_limBodyObr}. By Theorem~\ref{thm_AritLimFce} and HDD we have
$${\textstyle
\lim\frac{f^{-1}(b_n)-f^{-1}(c)}{b_n-c}=
\lim\frac{1}{\frac{f(a_n)-f(b)}{a_n-b}}
=\frac{1}{\lim\frac{f(a_n)-f(b)}{a_n-b}}
=\frac{1}{f'(b)}\,.
}
$$
By HDD we have $(f^{-1})'(c)=\frac{1}{f'(b)}$.

2. Let $f'(b)=0$. Then $c\in L(f[M])$ by Exercise~\ref{ex_pro0deroPla}. Suppose that 
$f$ decreases (respectively increases) at $b$. Then there is an $n_0$ such that if $n\ge n_0$ then
$${\textstyle
\text{$\frac{f(a_n)-f(b)}{a_n-b}<0$ (respectively $\ds>0$)}\,. 
}
$$
The previous computation and  part~5 of Proposition~\ref{prop_dod1} show that 
$${\textstyle
\text{$(f^{-1})'(b)=\frac{1}{0^-}=-\infty$ (respectively $\ds=+\infty$)}\,.
}
$$

3. Let $f'(b)=\pm\infty$ and $c\in L(f[M])$. Then 
$${\textstyle
(f^{-1})'(c)=\frac{1}{\pm\infty}=0
}
$$ by part~1.
\kduk

\noindent
Here is a~counterexample to part~1 when the continuity of $f^{-1}$ at $c$ is 
dropped.

\begin{exer}\label{ex_CoVynechSpo}
Let
$${\textstyle
M=\big([-1,\,1]\setminus
\big\{\frac{1}{n}\cc\;n\in\N\big\}\big)\cup
\big\{1+\frac{1}{n}\cc\;n\in\N\big\}\,.
}
$$ 
We define $f\in\mathcal{F}(M)$  
by $f(x)= x$ if $x\ne
1+\frac{1}{n}$, and by $f(1+\frac{1}{n})=\frac{1}{n}$.
Then $f$ is injective and $f'(0)=1$, but for $c= 
f(0)=0$ the derivative $(f^{-1})'(c)=(f^{-1})'(0)$ does not exist. 
\end{exer}

We state the formula for the global derivative of inverses separately 
as a~corollary.

\begin{cor}[$(f^{-1})'$]\label{cor_deriInve}
Let $f\in\mathcal{R}$ be injective and 
$$
M=\{x\in 
f[M(f)]\cc\;\text{$f^{-1}$ is continuous at $x$}\}\,.
$$
\underline{Then} the function 
$\frac{1}{f'(f^{-1})}\,|\,M$ is a~subfunction of the function $(f^{-1})'$. 
\end{cor}
\duk
Let 
$h=\frac{1}{f'(f^{-1})}\,|\,M$ and $c\in M(h)$. Then 
$$
c\in M(f'(f^{-1}))\cap M\setminus Z(f'(f^{-1}))\,, 
$$
$f^{-1}$ is continuous at $c$ because $c\in M$, and by part~1 of 
Theorem~\ref{thm_DeriInvfce} it holds that $(f^{-1})'(c)=\frac{1}{f'(f^{-1}(c)}=h(c)$. Hence $h$ is a~subfunction of $(f^{-1})'$.
\kduk

\noindent
We do not need derivatives of inverses for the proof 
of Theorem~\ref{thm_deriSEF}, but for completeness, we still provide the
corollary with the simple formula for $(f^{-1})'$.

\begin{cor}[$(f^{-1})'=\frac{1}{f'(f^{-1})}$]\label{cor_deriInveSim}
Let $f\in\mathcal{R}$ be injective and $M= f[M(f)]$. If $D(f)=M(f)$, $f^{-1}\in\mathcal{C}$ and $M\sus L(M)$, \underline{then}    
$$
(f^{-1})'=\frac{1}{f'(f^{-1})}\,.
$$
\end{cor}
\duk
Let $h=
\frac{1}{f'(f^{-1})}$. It follows from the assumptions and from 
Corollary~\ref{cor_deriInve} that $M(h)=M\setminus Z(f'(f^{-1}))$ and
that $h$ is a~subfunction of $(f^{-1})'$. Suppose that $c\in M$ is such that $f'(f^{-1}(c))=0$. If $c\in M((f^{-1})')$ then part~1 of Theorem~\ref{thm_DerSlozFce} gives that
$$
1=(\mathrm{id}\,|\,M)'(c)
=(f(f^{-1}))'(c)
=f'(f^{-1}(c))\cdot(f^{-1})'(c)=
0\cdot(f^{-1})'(c)\,,
$$
which is impossible. Thus $c\not\in M((f^{-1})')$ and we deduce that 
$M((f^{-1})')\sus M\setminus Z(f'(f^{-1}))$. Hence 
$${\textstyle
\frac{1}{f'(f^{-1})}=h=(f^{-1})\,.
}
$$
\kduk

\section[${}^c$Basic elementary functions]{Basic elementary functions}\label{sec_deriElemFun}

We determine derivatives of functions in BEF (Definition~\ref{def_ZEF}). 
Constants are easy, by Exercise~\ref{ex_derKon} we have for every $c\in\R$ the derivative
$$
k_c'(x)=k_0(x)\,.
$$  

\medskip\noindent
{\em $\bullet$ The exponential function, sine and cosine. }We
strengthen Theorem~\ref{thm_spjMocRady1}
and differentiate power series. By Theorem~\ref{thm_spjMocRady1}, any real numbers $a_0$, $a_1$, $\ds$ satisfying
$$
\lim|a_n|^{1/n}=0
$$ 
determine for $x\in\R$ an abscon series 
$\sum_{n\ge0}a_nx^n$ with sum 
$${\textstyle
S(x)=\sum_{n\ge0}a_nx^n\in\mathcal{C}(\R)\,.
}
$$

\begin{thm}[derivatives of power series]\label{thm_derMocRady}
We\index{theorem!derivatives of power series|emph} have 
$${\textstyle
D(S(x))=M(S(x))=\R\,\text{ and }\,
S'(x)=\sum_{n\ge0}(n+1)a_{n+1}x^n\,.
}
$$
\end{thm}
\duk
Let $a_n\in\R$ for $n\in\N_0$ 
be such that $\lim|a_n|^{1/n}=0$. We define the function 
$${\textstyle
T(x)=\sum_{n\ge0}
(n+1)a_{n+1}x^n\,.
}
$$
By Proposition~\ref{prop_Lim_nta_odm_zn} and the assumption on $a_n$ we have
$$
\lim_{n\to\infty}\big((n+1)\cdot|a_{n+1}|\big)^{1/n}=0\,.
$$
Thus $M(T(x))=\R$. Let $x,c\in\R$ with $c\ne0$ and $y=\max(\{1,|c|,|x|\})$. Using the identity
$${\textstyle
\frac{a^{n+1}-b^{n+1}}{a-b}=\sum_{j=0}^na^j b^{n-j}\,, 
}
$$
we get the bound
\begin{eqnarray*}
&&{\textstyle
\big|\frac{1}{c}(S(x+c)-S(x))-T(x)\big|}\\
&&{\textstyle
\le\sum_{n\ge1}|a_{n+1}|\cdot
|\sum_{j=0}^n(x+c)^jx^{n-j}-(n+1)x^n| }\\
&&{\textstyle
\le|c|\sum_{n\ge1}|a_{n+1}|\cdot\sum_{j=1}^n\sum_{i=1}^j\binom{j}{i}y^{i-1}y^{n-i}}\\
&&{\textstyle
\le|c|\sum_{n\ge1}|a_{n+1}|\cdot y^{n-1}\cdot\sum_{j=1}^n2^j\le|c|\cdot
\sum_{n\ge1}|a_{n+1}|\cdot
(2y)^{n+1}\,.
}
\end{eqnarray*}
For $c\to0$ this goes in limit to $0$. Hence 
$${\textstyle
S'(x)=\lim_{c\to0}\frac{S(x+c)-S(x)}{c}=T(x)=\sum_{n\ge0}(n+1)a_{n+1}x^n\,.}
$$
\kduk
\vspace{-3mm}
\begin{exer}\label{ex_zduvVseNero}
Explain estimates in the previous proof. 
\end{exer}

\begin{cor}[$\exp x$, $\cos x$ and $\sin x$]\label{cor_derExpKosSin}
We have the global formulas
$$
\text{$(\exp x)'=\exp x$, $(\cos x)'=-\sin x$ and $(\sin x)'=\cos x$}\ \ 
(\in\mathcal{F}(\R))\,.
$$
\end{cor}
\duk
By the previous theorem, we have ($x\in\R$)
\begin{eqnarray*}
(\exp x)'&=&{\textstyle\big(\sum_{n\ge0}
\frac{1}{n!}x^n\big)'=\sum_{n\ge0}\frac{(n+1)}{(n+1)!}x^n=\exp x\,,}\\
(\cos x)'&=&{\textstyle\big(\sum_{n\ge0}(-1)^n \frac{1}{(2n)!}x^{2n}\big)'=\sum_{n\ge0}
(-1)^{n+1}\frac{(2n+2)}{(2n+2)!}}x^{2n+1}\\
&=&-\sin x\,\text{ and}\\
(\sin x)'&=&{\textstyle\big(\sum_{n\ge0}
(-1)^n\frac{1}{(2n+1)!}x^{2n+1}\big)'=\sum_{n\ge0}
(-1)^n\frac{(2n+1)}{(2n+1)!}}x^{2n}=\cos x\,.
\end{eqnarray*}
\kduk

\noindent 
{\em $\bullet$ Logarithm. }Since $\log x$ is the inverse of $\exp x$ and $(\exp x)'=\exp x$, 
Corollary~\ref{cor_deriInveSim} yields the global formula
$$
(\log x)'=\frac{1}{(\exp x)'(\log x)}=
\frac{1}{\exp(\log x)}=
\frac{k_1(x)}{x\,|\,(0,\,+\infty)}
={\textstyle\frac{1}{x}\,|\,(0,\,+\infty)\,.}
$$
We stress that $(\log x)'$ is not the function $\frac{1}{x}$ ($\in\mathcal{F}(\R\setminus\{0\})$).

\begin{exer}\label{ex_deriLogAbs}
What is $(\log|x|)'$?    
\end{exer}

\noindent
{\em $\bullet$ Real exponentiation. }In the next two exercises we differentiate specializations of the function $a^b$.

\begin{exer}\label{ex_derirealMoc}
Let $a,b\in\R$ and $m\in\Z$. Prove the global formulas below. 
\begin{enumerate}
\item For $a>0$ the derivative 
$(a^x)'=a^x\cdot\log a$ is in $\mathcal{F}(\R)$.
\item For $b>1$  the derivative $(x^b)'=bx^{b-1}$ is in 
$\mathcal{F}([0,+\infty))$.
\item For $b=1$ the derivative $(x^b)'=k_1(x)\,|\,[0,+\infty)$.
\item For $b<1$  the derivative $(x^b)'=bx^{b-1}$ is in $\mathcal{F}((0,+\infty))$.
\item We have $(0^x)'=k_0(x)\,|\,(0,+\infty)$.
\item For $m>0$  the derivative $(x^m)'=mx^{m-1}$ is in $\mathcal{F}(\R)$.
\item For $m=0$ the derivative $(x^m)'=k_0(x)$ is in 
$\mathcal{F}(\R)$.
\item For $m<0$  the derivative $(x^m)'=mx^{m-1}$ is in $\mathcal{F}(\R\setminus\{0\})$.
\end{enumerate} 
\end{exer}

\begin{exer}\label{ex_derxnabLoc}
Let $b\in(0,1)$. Prove that 
$(x^b)'(0)=+\infty$.
\end{exer}

\noindent
{\em $\bullet$ Tangent and cotangent. }

\begin{exer}\label{ex_tanCot1}
Prove the global formulas
$$
(\tan x)'=\frac{1}{\cos^2x}\,\text{ and }\,(\cot x)'=
-\frac{1}{\sin^2 x}\,.
$$
\end{exer}

\noindent
{\em $\bullet$ Inverse trigonometric functions. }

\begin{exer}\label{ex_InvTri}
Prove the global formulas below. 
\begin{enumerate}
\item The derivative $(\arcsin x)'=
\frac{1}{\sqrt{1-x^2}}$ is in $\mathcal{F}((-1,1))$.
\item The derivative $(\arccos x)'=
-\frac{1}{\sqrt{1-x^2}}$ is in $\mathcal{F}((-1,1))$.
\item The derivative $(\arctan x)'=
\frac{1}{1+x^2}$ is in $\mathcal{F}(\R)$.
\item The derivative  $(\mathrm{arccot}\,x)'=
-\frac{1}{1+x^2}$ is in $\mathcal{F}(\R)$.
\end{enumerate}
\end{exer}

\begin{exer}\label{ex_derArcLoc}
Show that 
$$
(\arcsin x)'(-1)=(\arcsin x)'(1)=+\infty
$$ 
and that 
$$
(\arccos x)'(-1)=(\arccos x)'(1)=-\infty\,.
$$
\end{exer}

\noindent
{\em $\bullet$ Overview of derivatives of $f\in\mathrm{BEF}$. }We first list 
$f$ with $D(f)=M(f)$, and then give the three 
exceptions for which $D(f)\ne M(f)$.\label{tableDer}

\begin{center}    
\begin{tabular}{l|l|l|l}
$f$ & $M(f)$ & $f'$ & $M(f')$\\
\hline\hline
$\exp x$ & $\R$ & $\exp x$ & $=M(f)$\\
\hline
$\sin x$ & $\R$ & $\cos x$ & $=M(f)$\\
\hline
$\cos x$ & $\R$ & $-\sin x$ & $=M(f)$\\
\hline
$\arctan x$ & $\R$ & $\frac{1}{x^2+1}$ & $=M(f)$\\
\hline
$\mathrm{arccot}\,x$ & $\R$ & $-\frac{1}{x^2+1}$ & $=M(f)$\\
\hline
\text{$a^x$ for $a\in(0,+\infty)$} & $\R$ & $a^x\cdot\log a$ & $=M(f)$\\
\hline
\text{$x^m$ for $m\in\N$} & $\R$ & $mx^{m-1}$ & $=M(f)$\\
\hline
\text{$x^0$ for $0\in\Z$} & $\R$ & $k_0(x)$ & $=M(f)$\\
\hline
\text{$k_c(x)$ for $c\in\R$} & $\R$ & $k_0(x)$ & $=M(f)$\\
\hline
\text{$x^m$ for negative $m\in\Z$} & $\R\setminus\{0\}$ & $mx^{m-1}$ & $=M(f)$\\
\hline
\text{$\log(|x|)$, not in BEF} & $\R\setminus\{0\}$ & $1/x$ & $=M(f)$\\
\hline
\text{$x^b$ for $b\in(1,+\infty)$} & $[0,+\infty)$ & $bx^{b-1}$ & $=M(f)$\\
\hline
\text{$x^1$ for $1\in\R$} & $[0,+\infty)$ & $k_1(x)\,|\,[0,+\infty)$ & $=M(f)$\\
\hline
$0^x$ & $(0,+\infty)$ & $k_0(x)\,|\,(0,+\infty)$ & $=M(f)$\\
\hline
\text{$x^b$ for $b\in(-\infty,0]$} & $(0,+\infty)$ & $bx^{b-1}$ & $=M(f)$\\
\hline
$\log x$ & $(0,+\infty)$ & $\frac{1}{x}\,|\,(0,+\infty)$ & $=M(f)$\\
\hline
$\tan x$ & 
\begin{tabular}{l}
\!\!\!$\R\setminus\{k\pi\,+$\\
\!\!\!$+\,\frac{\pi}{2}\cc\;k\in\Z\}$ 
\end{tabular}
& $\frac{1}{\cos^2 x}$ & $=M(f)$\\
\hline
$\cot x$ & 
$\R\setminus\{k\pi\cc\;k\in\Z\}$
 & $-\frac{1}{\sin^2 x}$ & $=M(f)$\\
 \hline\hline
\text{$x^b$ for $b\in(0,1)$} & $[0,+\infty)$ & $bx^{b-1}$ & $(0,+\infty)$\\
\hline
$\arcsin x$ & $[-1,1]$ & $\frac{1}{\sqrt{1-x^2}}$ & $(-1,1)$\\
\hline
$\arccos x$ & $[-1,1]$ & $-\frac{1}{\sqrt{1-x^2}}$ & $(-1,1)$
\end{tabular}
\end{center}

\noindent
The infinite values of pointwise derivatives of the function in BEF are not included in the table. They are as follows.
\begin{itemize}
\item $(x^b)'(0)=+\infty$ for 
$b\in(0,1)$, 
\item $(\arcsin)'(-1)=(\arcsin)'(1)=+\infty$ and
\item $(\arccos)'(-1)=(\arccos)'(1)=-\infty$.
\end{itemize}

\begin{exer}\label{ex_Mer}
Explain why we cannot merge the twelfth and thirteenth row of the 
table in one row with the functions $x^b$ for $b\in[1,+\infty)$.    
\end{exer}

Theorem~\ref{thm_nespDer} describes a~function 
$f\in\mathcal{R}$ such that $M(f)=D(f)\ne\emptyset$ and $f'\not\in\mathcal{C}$. Here is 
a~different, well-known example with $M(f)=D(f)=\R$. 

\begin{exer}\label{ex_znovuTenPr}
Let
$f\in\mathcal{F}(\R)$ be given by $f(0)=0$
and by 
$${\textstyle
f(x)= x^2\sin\big(\frac{1}{x}\big)
}
$$
for $x\ne0$. Show that $M(f)=D(f)=\R$ and $f'\not\in\mathcal{C}$.
\end{exer}

\noindent
{\em $\bullet$ Derivatives of polynomials and rational 
functions. }The zero polynomial and polynomials with degree $0$ 
are constant functions. As we know, $k_c'(x)=k_0(x)$. It is easy 
to see that for degree $d\ge1$ we have
$${\textstyle
\big(\sum_{i=0}^d a_ix^i\big)'=
\sum_{i=0}^{d-1} (i+1)a_{i+1}x^i\,.
}
$$
Let $p$ and $q$ be polynomials. By Corollary~\ref{cor_deriRat}, the derivative of the rational function $p(x)/q(x)$ is again a~rational function, 
$$
\Big(\frac{p(x)}{q(x)}\Big)'=\frac{p'(x)q(x)-p(x)q'(x)}{q(x)^2}\,.
$$
See Proposition~\ref{prop_poleRatDer} for an important property of derivatives of rational functions.

\begin{exer}\label{ex_onDerRat1}
What happens if $q(x)=k_0(x)$?    
\end{exer}

\begin{exer}\label{ex_onDerRat2}
Let $k\in\Z$. Find a~rational function $r_k(x)$ such that 
$$
r_k'(x)=x^k\,,
$$
or show that no such rational function exists.
\end{exer}

\section[${}^c$Simple elementary functions]{Simple elementary functions}\label{sec_deriEF}

\noindent
{\em $\bullet$ A~problem on elementary functions. }We introduced elementary functions in Definitions~\ref{def_obecEF} and
\ref{def_obecEF2}. In Theorem~\ref{thm_EFjsouSpoj} we established by 
induction their continuity. It was easy because continuity
of functions is preserved under sum, product, ratio, and composition. One might hope to prove similarly with the help of the formulas 
$${\textstyle
(f+g)'=f'+g',\ (fg)'=f'g+fg',\ 
(\frac{f}{g})'=\frac{f'g-fg'}{g^2}\,\text{ and }\,(f(g))'=f'(g)\cdot g'}
$$
that the derivative of an elementary function is elementary. Except that we know from Sections~\ref{sec_aritmDeri} and 
\ref{sec_deriSlozInve} that these formulas in general do not hold, not even
for elementary functions. For example, 
let 
$$
f(x)=\arcsin x\,\text{ and }\,g(x)=-\arcsin x\,. 
$$
Then 
$$
(f(x)+g(x))'=k_0(x)\,|\,[-1,1]
$$
which differs from
$${\textstyle
f'(x)+g'(x)=\frac{1}{\sqrt{1-x^2}}-
\frac{1}{\sqrt{1-x^2}}=k_0(x)\,|\,(-1,1)\,.
}
$$
Although
$f'(x),g'(x)\in\mathrm{EF}$, we cannot deduce from it that $(f(x)+g(x))'\in\mathrm{EF}$ because $(f(x)+g(x))'\ne f'(x)+g'(x)$. In this 
case, $(f(x)+g(x))'$ is elementary because
$$
(f(x)+g(x))'=k_0(x)\,|\,[-1,1]=
\sqrt{1+x}-\sqrt{1+x}+\sqrt{1-x}-\sqrt{1-x}\,,
$$
but this justification has nothing to do with the derivatives
$f'(x)=\frac{1}{\sqrt{1-x^2}}$ and $g'(x)=-\frac{1}{\sqrt{1-x^2}}$. We pose the following problem.

\begin{prob}[elementary derivatives]\label{prob_derEF}
Prove or disprove that 
$$
f\in\mathrm{EF}\Rightarrow
f'\in\mathrm{EF}\,.
$$    
\end{prob}

\noindent
{\em $\bullet$ Simple elementary functions. }We partially solve this problem in affirmative by 
introducing a~relatively large subset 
$\mathrm{SEF}$ of $\mathrm{EF}$ and showing that it is closed to 
derivatives. 

\begin{defi}[VBEF $\&$ SEF]\label{def_VBEF_SEF}
\underline{Very basic elementary functions\index{very basic elementary functions, VBEF|emph}} are 
$$
\mathrm{VBEF}=
\{k_c(x)\cc\;c\in\R\}\cup\{\exp x,\,\log x,\,\sin x,\,\arcsin_0 x\}\label{VBEF}
$$
where 
$$
\arcsin_0 x=\arcsin x\,|\,(-1,\,1)\,.\label{arcsinZero}
$$
\underline{Simple elementary functions\index{simple elementary 
functions, 
SEF|emph}},\label{SEF} 
abbreviated {\em SEF}, is the subset of {\em EF} obtained by 
replacing in Definition~\ref{def_obecEF2} the 
starting set of functions {\em RBEF} with {\em VBEF}. 
\end{defi}
We obtain simple elementary functions by starting from constants, the 
exponential function, logarithm, sine, and the restricted arkus sine, and repeatedly applying addition, multiplication, division, and composition. The next exercise explains why we do not add to VBEF the restrictions $x^b\,|\,(0,+\infty)$ for real $b>0$.

\begin{exer}\label{ex_noNeedAdd}
Express them in terms of $k_c(x)$, $\exp x$ and 
$\log x$.
\end{exer}

\begin{exer}\label{ex_ercsin0}
Show that $\arcsin_0 x\in\mathrm{EF}$.  
\end{exer}

We show that the family of functions SEF is closed to derivatives.  

\begin{thm}[derivatives in SEF]\label{thm_deriSEF}
For\index{theorem!derivatives of simple elementary functions|emph} 
every function $f\in\mathrm{SEF}$, the domain
$M(f)$ is an open set, $D(f)=M(f)$ and $f'\in\mathrm{SEF}$ as well.    
\end{thm}
\duk
First, recall 
that the identity function 
$\mathrm{id}(x)$ is in $\mathrm{SEF}$; 
it equals $\log(\exp x)$. For the constant function $k_c(x)$, we usually write just $c$.
To begin with, we check that every function in $\mathrm{VBEF}$ has the three stated properties. Clearly,  $M(k_c(x))=M(\exp x)=M(\sin x)=\R$, $M(\log x)=(0,+\infty)$ and $M(\arcsin_0 x)=(-1,1)$
are open sets. The derivatives $k_c'(x)=k_0(x)$, $(\exp x)'=\exp x$, 
$${\textstyle
(\log x)'=\frac{1}{x}\,|\,(0,\,+\infty)=\frac{1}{x}+\log x+
(-1)\cdot\log x\,, 
}
$$
$(\sin x)'=\cos x=\sin(x+\frac{\pi}{2})$ 
and
$${\textstyle
(\arcsin_0(x))'=\frac{1}{\sqrt{1-x^2}}=
1/\exp(\frac{1}{2}\log(1-x\cdot x))
}
$$
are in SEF. Each has the same domain as the original 
function.

Let $f$ be in $\mathrm{SEF}$.
We proceed by induction on the length of a~generating word of 
$f$. 
In the induction step, we show for every case (i) $f=g+h$ or 
(ii) $f=g\cdot h$ or (iii) $f=g/h$ or (iv) $f=g(h)$, where
$g,h\in\mathrm{SEF}$ are functions with the three stated properties, that $f$ has these properties as well. 
We use that $M\sus L(M)$ for every open set $M\sus\R$ (Exercise~\ref{ex_opAndLimP}) and that open sets in $\R$ are closed to finite 
intersections (part~3 of Proposition~\ref{prop_openSet}). It follows that in case (i) the set
$$
M(f)=M(g+h)=M(g)\cap M(h)
$$ 
is open. By Corollary~\ref{cor_deriSumStan}, 
$f'=g'+h'\in\mathrm{SEF}$ and 
$$
M(f')=M(g')\cap M(h')=M(g)\cap M(h)=M(f)\,. 
$$
Case (ii) is similar except that we use 
Corollary~\ref{cor_Leibform}.

In cases (iii) and (iv), we use
$\mathrm{EF}\sus\mathcal{C}$ 
(Theorem~\ref{thm_EFjsouSpoj}) and
two results on open sets and continuous functions: the zero set of a~continuous function is relatively closed
(Proposition~\ref{prop_zeroClos}) and the preimage of an open set by a~continuous function is relatively open
(Proposition~\ref{prop_preiOpen}). Hence in case (iii) there is a~closed set $U\sus\R$ such that
\begin{eqnarray*}
M(f)&=&(M(g)\cap M(h))\setminus Z(h)=
M(g)\cap M(h)\cap(\R\setminus Z(h))\\
&=&M(g)\cap M(h)\cap\big(\R\setminus(M(h)\cap U)\big)\\
&=&M(g)\cap M(h)\cap\big((\R\setminus M(h))\cup(\R\setminus U)\big)\\
&=&M(g)\cap M(h)\cap(\R\setminus U)
\end{eqnarray*}
and this is an open set. By Corollary~\ref{cor_deriRat},
$f'=\frac{g'h-gh'}{h^2}\in\mathrm{SEF}$ and
\begin{eqnarray*}
M(f')&=&(D(g)\cap M(h)\cap M(g)\cap D(h))\setminus Z(h^2)\\
&=&M(g)\cap M(h)\setminus Z(h)\\
&=&M(f)\,.
\end{eqnarray*}
Finally, in case (iv) the set $M(f)=h^{-1}[M(g)]$ is open. By Corollary~\ref{cor_derSlozF}, 
$f'=g'(h)\cdot h'\in\mathrm{SEF}$ and 
\begin{eqnarray*}
M(f')&=&h^{-1}[M(g')]\cap M(h')
=h^{-1}[M(g)]\cap M(h)\\
&=&h^{-1}[M(g)]=M(g(h))\\
&=&M(f)\,.    
\end{eqnarray*}
\kduk

\chapter[Mean value theorems]{Mean value theorems}\label{chap_pr8}

\bigskip\noindent
In Section~\ref{sec_3vetyOstrhod} we
meet three mean value theorems: Rolle's Theorem~\ref{thm_Rolle},
Lagrange's Theorem~\ref{thm_Lagrange} and
Cauchy's Theorem~\ref{thm_Cauchy}.
The first two theorems are generalized in 
Theorems~\ref{thm_rolle2} and \ref{thm_stroLagr}. The extending
Sections~\ref{sec_P_rec}--\ref{sec_tranNumbLiouv} are devoted to 
three applications of mean value theorems. In Section~\ref{sec_P_rec} we show by means 
of Rolle's theorem that the sequence 
$$
(\log n)=(0,\,\log 2,\,\log 3,\,\ds)
$$ 
is not P-recurrent. In Theorem~\ref{thm_konMnoKor} we actually 
prove a~more general result. In Sections~\ref{sec_tranNumb} and
\ref{sec_tranNumbLiouv} we show with the help of Lagrange's theorem in two 
effective ways that real
transcendental numbers exist. 

In Section~\ref{sec_monoLHR} in
Theorem~\ref{thm_derMono1} and 
Proposition~\ref{prop_deriMono2} we obtain by means of the 
first derivative results on monotonicity of functions. 
Theorems~\ref{thm_LHP1} and \ref{thm_LHP2} are l'Hospital rules
for computing functional limits of the type $\frac{0}{0}$ and 
$\frac{\infty}{\infty}$. In Section~\ref{sec_secoDeri} in Proposition~\ref{prop_DruhaDerExtr} we 
determine by the sign 
of the derivative 
$$
(f')'(b)
$$ 
the type of the local extreme of $f$ at $b$. In Theorem~\ref{thm_exOnesDer} we prove that
any convex or concave function $f$ defined on a~set $M\sus\R$ with no minimum and no maximum has 
on $L^{\pm}(M)$ finite one-sided derivatives. Thus such $f$ is 
continuous. In Theorem~\ref{thm_druDerKonvKonk} 
we determine by the sign of $(f')'(b)$ 
convexity/concavity of a~function $f$ defined on an interval. Definition~\ref{def_inflexe} introduces inflection points.
Theorems~\ref{thm_infl} and
\ref{thm_postacInflexe} provide necessary and sufficient conditions for their existence. 

In Section~\ref{sec_drawGr} we give thirteen steps for 
determining the main geometric features of the graph of a~function. Step~0 
places the function in the hierarchy 
$$
\mathrm{SEF}\sus\mathrm{EF}\sus\mathcal{R}\,.
$$
We demonstrate the procedure of thirteen steps on three examples for functions  
$\sgn\,x$, $\tan x$ and $\arcsin\big(\frac{2x}{x^2+1}\big)$. 

\section[${}^c$Mean value theorems of Rolle, Lagrange and Cauchy]{Rolle, Lagrange and Cauchy}
\label{sec_3vetyOstrhod}

The three theorems named after these mathematicians establish various relations between values of functions and 
their derivatives. In these theorems, $a<b$ are real numbers.

\medskip\noindent
{\em $\bullet$ Rolle's theorem. }We begin with this theorem and deduce the other two theorems from it.

\begin{thm}[Rolle~1]\label{thm_Rolle}
Let\index{theorem!Rolle 1|emph} 
$f\in\mathcal{C}([a,b])$ and $f(a)=f(b)$. If $f'(c)\in\R^*$ for every $c\in(a,b)$, \underline{then} $f'(c)=0$ for 
some $c\in(a,b)$.
\end{thm}
\duk
If $f(x)$ is constant, we are done as $f'(c)=0$ for every $c\in(a,b)$. Else, by Corollary~\ref{cor_priMaxi}, there exists $c\in(a,b)$
such that $f$ has at $c$ a~global extreme. Since
$c\in L^{\mathrm{TS}}([a,b])$ and $f'(c)$ exists, by
Corollary~\ref{cor_NPELE} we have $f'(c)=0$.
\kduk

\noindent
{\em Michel Rolle\index{Rolle, Michel} (1652--1719)} was a~French mathematician. 

\begin{exer}\label{ex_RolleAbsVal}
Let $f(x)$ be the absolute value $|x|$ restricted to $[-1,1]$. Then $f(1)=f(-1)=1$
and $f$ is continuous, but $f'(c)\ne0$ for every $c\in(-1,1)$. Which assumption of Rolle's theorem is not 
satisfied?
\end{exer}

The well known solution of Exercise~\ref{ex_zeroCan} employs the algebraic division of 
polynomials with remainder. Now we show an analytic solution.

\begin{cor}[the number of roots]\label{cor_bouNumRoo}
Let $f$ be a~nonzero polynomial. The set of zeros $Z(f)$ is finite and
$$
|Z(f)|\le\deg f\,.
$$
\end{cor}
\duk
We argue by contradiction. Let $g$ be a~nonzero polynomial
with the minimum degree $d\in\N_0$ such that $|Z(g)|>\deg g=d$. So there exist $d+1$ real numbers $a_1<a_2<\ds<a_{d+1}$ such that
$$
g(a_1)=g(a_2)=\ds=g(a_{d+1})=0\,.
$$
We have $d\ge1$ because nonzero constant polynomials have no roots. By Theorem~\ref{thm_Rolle}  there exist $d$ real numbers $b_i$ such that
$$
a_1<b_1<a_2<b_2<a_3<\ds<a_{d+1}\,\text{ and }\,
g'(b_1)=g'(b_2)=\ds=g'(b_d)=0\,.
$$
As we know from the end of 
Section~\ref{sec_deriElemFun}, 
$g'$ is a~nonzero polynomial with degree $d-1$. We have obtained a~contradiction with the minimality of $d$.
\kduk

The next version of Rolle's theorem does not need derivatives; 
it applies 
to the function in Exercise~\ref{ex_RolleAbsVal}. 

\begin{prop}[Rolle~2]\label{prop_Rolle}
Let 
$f\in\mathcal{C}([a,b])$ and $f(a)=f(b)$. \underline{Then} there exists $c\in(a,b)$ such that for every $\ep$ there exist $d_1,d_2\in (a,b)$ satisfying
$$
d_1<c<d_2,\ d_2-d_1\le\ep\,\text{ and }\,f(d_1)=f(d_2)\,.
$$
Geometrically, some point $\langle c,f(c)\rangle\in G_f$ can be  
arbitrarily tightly enclosed by pairs of intersection points
of horizontal secants.
\end{prop}
\duk
If $f(x)$ is constant, the result trivially holds. Else, let $f$ 
attain at $c\in(a,b)$ the maximum value (for the minimum value, we 
would argue similarly). 
We distinguish three cases. (i) $f(x)$ is constantly $f(c)$ on $[c-\de,c]$ ($\sus[a,b]$) or 
(ii) $f(x)$ is constantly $f(c)$ on $[c,c+\de]$ ($\sus[a,b]$) or 
(iii) $f(x)$ attains values smaller than $f(c)$ arbitrarily 
close to $c$ on both sides of $c$. In case (i), respectively (ii), we replace $c$
with $c-\frac{\de}{2}$, respectively $c+\frac{\de}{2}$, 
and then $c$ clearly has the stated property. In case (iii), we 
keep $c$ 
and find the required points $d_1$ and $d_2$  by means of 
Theorem~\ref{thm_mezihodnoty}. 
\kduk

The third version of Theorem~\ref{thm_Rolle} 
generalizes domains $[a,b]$ to compact sets. 
Recall that $L^{\mathrm{TS}}(M)$ denotes the set of two-sided limit points of $M$.

\begin{thm}[Rolle~3]\label{thm_rolle2}
Let\index{theorem!Rolle 3|emph}
$M\sus\R$ be a~compact set with $M\cap L^{\mathrm{TS}}(M)\ne\emptyset$
and let $f\in\mathcal{C}(M)$
satisfy two conditions.
\begin{enumerate}
\item The restriction 
$f\,|\,(M\setminus L^{\mathrm{TS}}(M))$ is constant.
\item We have
$f'(c)\in\R^*$ for every $c\in M\cap L^{\mathrm{TS}}(M)$.
\end{enumerate}
\underline{Then} $f'(c)=0$ for some $c\in M\cap L^{\mathrm{TS}}(M)$.
\end{thm}
\duk
If $f(x)$ is constant,
then $f'(c)=0$ for
every $c\in M\cap L^{\mathrm{TS}}(M)$ by the assumptions and Theorem~\ref{thm_priznakExtr}. 
Else, if $f(x)$
is not constant, by the assumptions and 
Corollary~\ref{cor_priMaxi}, $f$
attains at some $c\in 
M\cap L^{\mathrm{TS}}(M)$ global extreme. Using condition~2 and
Theorem~\ref{thm_priznakExtr}, we get $f'(c)=0$.
\kduk

\noindent
In Theorem~\ref{thm_Rolle}
we have $M=[a,b]$, $M\cap L^{\mathrm{TS}}(M)=(a,b)$ and 
$M\setminus L^{\mathrm{TS}}(M)=
\{a,b\}$. To show a~different application of Theorem~\ref{thm_rolle2}, we set 
$${\textstyle
N\equiv
\{(-1)^n\frac{1}{n}\cc\;n\in\N\}\,.
}
$$
The theorem then applies with
$M=\{0\}\cup N$ and any constant function $f$ in $\mathcal{F}(M)$. Now 
$M\cap L^{\mathrm{TS}}(M)=
\{0\}$ and 
$M\setminus L^{\mathrm{TS}}(M)=N$.

\medskip\noindent
{\em $\bullet$ Lagrange's theorem. }This mean value theorem has an interesting geometric
interpretation and generalization.

\begin{thm}[Lagrange~1]\label{thm_Lagrange}
Let\index{theorem!Lagrange 1|emph} 
$f\in\mathcal{C}([a,b])$. If $f'(c)\in\R^*$ for every $c$ in $(a,b)$, \underline{then} 
$$
f'(c)=\frac{f(b)-f(a)}{b-a}
$$
for some $c\in(a,b)$.
\end{thm}
\duk
Let $z\equiv\frac{f(b)-f(a)}{b-a}$. Then the function 
$$
g(x)\equiv
f(x)-z(x-a)\ \  (\in\mathcal{C}([a,b]))
$$
satisfies the assumptions of
Theorem~\ref{thm_Rolle}, especially $g(a)=g(b)$ ($=f(a)$). Thus $g'(c)=f'(c)-z=0$ for some $c\in(a,b)$, and $f'(c)=z$. 
\kduk

\noindent
Geometrically, there exists a~point 
$C\equiv\langle c,f(c)\rangle\in G_f$ with $a<c<b$
such that
the tangent to $G_f$ at $C$ is parallel to the secant 
$\kappa(a,f(a),b,f(b))$. 
{\em Joseph-Louis Lagrange\index{Lagrange, 
Joseph-Louis} (1736--1813)} was a~French mathematician,
physicist, and astronomer of Italian origin. 

\begin{exer}\label{ex_logicEquiv}
Let $a<b$ be real numbers. For $f\in\mathcal{R}$ let $R(f)$ be Rolle's theorem for $f$\,---\,if $f\in\mathcal{C}([a,b])$, $f(a)=f(b)$ and 
$f'(c)\in\R^*$ for every $c\in(a,b)$, then $f'(c)=0$ for some $c\in(a,b)$. Similarly, for $f\in\mathcal{R}$ let 
$L(f)$ be Lagrange's theorem for $f$\,---\,if $f\in\mathcal{C}([a,b])$ and 
$f'(c)\in\R^*$ for every $c\in(a,b)$, then $f'(c)=\frac{f(b)-f(a)}{b-a}$ for some $c\in(a,b)$.
Explain  why we know that the logical equivalence
``Rolle's theorem $\iff$ Lagrange's theorem'', formally 
$$
\big(\forall f\in\mathcal{R}\cc\;R(f)\big)\iff 
\big(\forall f\in\mathcal{R}\cc\;L(f)
\big)\,,
$$
holds before we even start proving both theorems.
\end{exer}

Tangent lines to graphs are of two kinds. Let
$f\in\mathcal{F}(M)$, $b\in M$ and let $f'(b)\in\R$. We consider the tangent
$$
l(x)\equiv f'(b)(x-b)+f(b)\ \ (\in\mathcal{F}(\R))
$$
to $G_f$ at $B\equiv\langle b,f(b)\rangle$. If for some
$\de$ we have
$$
\text{$\forall x\in U(b,\,\de)\cap M\cc\;l(x)\ge f(x)$ or 
$\forall x\in U(b,\,\de)\cap M\cc\;l(x)\le f(x)$}\,, 
$$
we say that $\ell$ is a~\underline{touching\index{tangent@tangent 
(line)!touching|emph} tangent} (see Theorem~\ref{thm_suppLine}). 
If no such $\de$ exists, we call $\ell$ a~\underline{cutting\index{tangent@tangent (line)!cutting|emph} tangent}. In the latter case, $G_f$
contains points arbitrarily close to $B$ that lie below and above the tangent. In the former case, all points in $G_f$ sufficiently close to
$B$ lie on the same side of the tangent.
The second version of Lagrange's theorem yields for any slope $r$ close to the slope of the secant a~touching tangent with slope $r$.

\begin{thm}[Lagrange~2]\label{thm_stroLagr}
Let\index{theorem!Lagrange 2|emph} $f\in\mathcal{C}([a,b])$, $f'(c)\in\R^*$
for every $c$ in $(a,b)$ and let $s\equiv\frac{f(b)-f(a)}{b-
a}$. We assume that $G_f$ is not a~straight 
plane segment. \underline{Then} there exists $\de$ such that for every $r\in(s-\de,\,s+\de)$ the graph $G_f$ 
$$
\text{has a~touching tangent with slope $r$ at some 
point $\langle c,f(c)\rangle$ with $c\in(a,b)$}\,.
$$ 
\end{thm}
\duk
Let $T$ ($\sus\R^2$) be the straight segment joining the points $\langle a,f(a)\rangle$ and $\langle b,f(b)\rangle$. We know that $G_f\setminus T\ne\emptyset$. So there is 
a~point $\langle d,f(d)\rangle\in G_f$ with $a<d<b$ lying above or below $T$. We assume
the latter case. In the former case, we argue
similarly. Let $s_1$, respectively $s_2$, be the slope of the straight plane segment joining the points $\langle a,f(a)\rangle$ and
$\langle d,f(d)\rangle$, respectively $\langle d,f(d)\rangle$ and $\langle b,f(b)\rangle$. Then $s_1<s<s_2$ and we take any $\de$ such that 
$$
s_1<s-\de<s+\de<s_2\,.
$$

We show that for every
$r\in(s-\de,s+\de)$ the graph $G_f$ has a~touching tangent with slope $r$ at 
some interior point. Let $r\in(s-\de,s+\de)$. For $t\in\R$, we define the linear function $l_t(x)\equiv rx+t$ ($\in\mathcal{F}
(\R)$) and set (in the linear order $\langle \R^*,<\rangle$)
$$
u\equiv\sup(\{t\in\R\cc\;\text{$f(x)\ge l_t(x)$ 
for every $x\in[a,\,b]$}\})\,.
$$
By Corollary~\ref{cor_priMaxi}, the function $f$ has 
a~global minimum. So
$u\in\R$. Clearly, $f(x)\ge l_u(x)$ for every $x\in[a,b]$ and for some $x=c\in[a,b]$, the equality occurs (Exercise~\ref{ex_easySee}). If $c=a$ or $c=b$, then by the choice of $\de$ and $r$, we have 
$f(d)<l_u(d)$. Hence $c\in(a,b)$. By
Theorem~\ref{thm_suppLine}, the line $l_u(x)$ is 
a~touching tangent to $G_f$ at 
$\langle c,f(c)\rangle$. By the definition of $l_u(x)$, the slope is $r$.
\kduk
\vspace{-3mm}
\begin{exer}\label{ex_easySee}
Prove the claim about $f(x)$ and $l_u(x)$.    
\end{exer}

\begin{exer}\label{ex_ifStrai}
What happens if the graph $G_f$ is a~straight plane segment?    
\end{exer}

\noindent
{\em $\bullet$ Cauchy's theorem. }A~new feature is that this mean value theorem involves two functions. 

\begin{thm}[Cauchy]\label{thm_Cauchy}
Let\index{theorem!Cauchy's mean value theorem|emph} 
$f,g\in\mathcal{C}([a,b])$, 
$g(b)\ne g(a)$ and let $f'(c)\in\R^*$ and $g'(c)\in\R$ for every $c\in(a,b)$; so now $g'(c)\ne\pm\infty$. \underline{Then} 
$$
f'(c)=\frac{f(b)-f(a)}{g(b)-
g(a)}\cdot g'(c)
$$
for some $c\in(a,b)$. 
\end{thm}
\duk
Let $z\equiv\frac{f(b)-f(a)}{g(b)-g(a)}$.
Then the function
$$
h(x)\equiv f(x)-z(g(x)-g(a))\ \  (\in\mathcal{C}([a,\,b]))
$$
satisfies the assumptions of Theorem~\ref{thm_Rolle}, especially $h(a)=h(b)=f(a)$. Hence there exists $c\in(a,b)$
with $h'(c)=f'(c)-zg'(c)=0$. So $f'(c)=zg'(c)$.
\kduk
\vspace{-3mm}
\begin{exer}\label{ex_kdeByvadd}
Where in the proof would be $g'(c)=\pm\infty$ problematic? Modify the statement of the theorem to accommodate the values $g'(c)=\pm\infty$.
\end{exer}

\section[${}^c$The sequence $(\log n)$ is not P-recurrent]{The sequence 
$(\log n)$ is not P-recurrent}\label{sec_P_rec}

This section is based on the preprint \cite{klaz_arXivLog} of the author. We show with the help of
Theorem~\ref{thm_Rolle} that the sequence
$$
(\log n)=(0,\,\log 2,\,\log 3,\,\ds)
$$ 
does not satisfy any recurrence relation from the class 
of so called P-recurrences. We begin with the definition of these recurrences.

\medskip\noindent
{\em $\bullet$ {\em P}-recurrent sequences. }These sequences generalize 
\underline{C-recurrent\index{c reccurent@C-recurrent sequence|emph}} (constantly recurrent) sequences. A~C-recurrent sequence
$(a_n)\sus\R$ satisfies for some $k\in\N$ real coefficients 
$c_1$, $c_2$, $\ds$, $c_k$,
not all of them zero, and every integer $n\ge k$ the 
relation
$${\textstyle
\sum_{i=1}^k c_i a_{n-i+1}
=c_1a_n+c_2a_{n-1}+\ds+c_ka_{n-k+1}
=0\,.
}
$$ 
A~well known C-recurrent sequence is the Fibonacci sequence
$$
(F_n)=(1,\,1,\,2,\,3,\,5,\,8,\,13,\,21,\,34,\,\ds)\,.
$$

P-recurrent sequences allow polynomials $c_i(n)$ instead of the constants $c_i$. So every C-recurrent sequence is 
P-recurrent.

\begin{defi}[P-recurrence]\label{def_PrekPosl}
Let $(a_n)\sus\R$. We say that the sequence $(a_n)$ is \underline{{\em P}-recurrent\index{P-recurrent sequence|emph}} if there
exist $k\ge1$ polynomials $p_1$, $p_2$, $\ds$, $p_k$ in $\mathrm{POL}$, not all of them zero, such that for every integer $n\ge k$ we have
$${\textstyle
\sum_{i=1}^k p_i(n)\cdot a_{n-i+1}=p_1(n)a_n+p_2(n)a_{n-1}+\ds+p_k(n)a_{n-k+1}=0\,.
}
$$
\end{defi}
For instance,  
$$
(a_n)=(n!)=(1,\,2,\,6,\,24,\,120,\,720,\,5040,\,\ds)
$$ 
is P-recurrent: for every $n\ge2$ we have
$$
1\cdot a_n+(-n)\cdot a_{n-1}=0\,.
$$

The next two exercises describe two ways to show that Definition~\ref{def_PrekPosl}
with $n\ge k$ relaxed to $n\ge n_0\ge k$ again yields 
P-recurrent sequences.

\begin{exer}\label{ex_staciNnula0}
Any relation
$${\textstyle
\sum_{i=1}^k p_i(n)\cdot a_{n-i+1}=0
}
$$
holding for every $n\ge n_0\ge k$ can be viewed as a~{\em P}-recurrence of 
an order $l\ge k$ holding for every $n\ge l$.     
\end{exer}

\begin{exer}\label{ex_staciNnula}
In any relation
$${\textstyle
\sum_{i=1}^k p_i(n)\cdot a_{n-i+1}=0
}
$$
holding for every $n\ge n_0\ge k$, 
we can replace the coefficients $p_i(x)$
with other polynomials $q_i(x)$ such that the new relation holds for every $n\ge k$.
\end{exer}
See the book \cite{stan_ECII} for uses of P-recurrent sequences in enumerative combinatorics.\index{enumerative 
combinatorics} In this section, we use Theorem~\ref{thm_Rolle} to prove the following.

\begin{thm}[logs are not P-recurrent, \cite{klaz_arXivLog}]\label{thm_logNotPrek}
The\index{theorem!logs are not P-recurrent|emph} 
sequence
$$
(\log n)=(0,\,\log 2,\,\log 3,\,\ds)
$$ 
is not {\em P}-recurrent.     
\end{thm}

\noindent
{\em $\bullet$ Poles. }We prove Theorem~\ref{thm_logNotPrek} with the help of poles. Let $f\in\mathcal{F}(M)$ and $b\in L(M)$. If  
$\lim_{x\to b}f(x)\in\R$, then $b$ is a~\underline{regular point\index{regular point|emph}} of $f$. If 
$$
f(x)\sim c(x-b)^{-k}\ \ (x\to b)
$$
for some $c\in\R\setminus\{0\}$ and $k\in\N$, we say that $b$ is a~\underline{pole\index{pole of 
a~function|emph}} of $f$ of \underline{order\index{pole of 
a~function!order of|emph}} $k$.

\begin{exer}\label{ex_onPoles}
Let $f,g\in\mathcal{R}$, $b\in L(M(f+g))$, $b$ be a~pole of $f$ of order 
$k$ and let $b$ be a~regular point of $g$ or a~pole of $g$ of order $l<k$. 
Then $b$ is a~pole of $f+g$ of order $k$.
\end{exer}

\begin{prop}[poles of derivatives]\label{prop_poleRatDer}
Let $r(x)\in\mathrm{RAC}$ and let $b\in\R$. \underline{Then} $b$ is 
a~regular point of $r'(x)$ or a~pole of $r'(x)$ with order 
at least $2$.    
\end{prop}
\duk
If $r(x)=k_0(x)$, then $r'(x)=k_0(x)$ and $b$ is 
a~regular point of $r(x)'$. We assume that $r(x)\ne k_0(x)$, use 
Theorem~\ref{thm_oRacFci} and write $r(x)$ as a~ratio of two polynomials, $r(x)=p(x)/q(x)$.
Using Exercise~\ref{ex_rootFactor} we get factorizations
$$
p(x)=p_0(x)\cdot(x-b)^k\,\text{ and }\,
q(x)=q_0(x)\cdot(x-b)^l\,,
$$
where $p_0(x),q_0(x)\in\mathrm{POL}$ are such that $p_0(b),q_0(b)\ne0$ and  $k,l\in\N_0$. We have
$$
r(x)=r_0(x)\cdot(x-b)^m\,,
$$
where $r_0(x)=p_0(x)/q_0(x)\in\mathrm{RAC}$ is such that $b\in M(r_0)$ and
$r_0(b)\ne0$, and $m=k-l\in\Z$. If $m\ge0$ then $b$ is a~regular point of 
$r'(x)$. Let $m<0$. Using results on derivatives in the previous chapter we get
$$
r'(x)=r_0'(x)\cdot(x-b)^m+
r_0(x)\cdot m(x-b)^{m-1}\,.
$$
By Exercise~\ref{ex_onPoles},  $b$ is a~pole of $r'(x)$ of order 
$1-m\ge2$.
\kduk

\noindent
{\em $\bullet$ Nonzero functions. }A~function $f\in\mathcal{R}$ is 
\underline{nonzero\index{function!nonzero|emph}} if $f(b)\ne0$ for some $b\in M(f)$. Recall that $Z(f)=\{b\in M(f)\cc\;f(b)=0\}$.

\begin{prop}[zeros of nonzero restrictions]\label{prop_nonzRestr}
Let $f\in\mathcal{R}$ be a~nonzero restriction of a~rational function.
\underline{Then} $Z(f)$ is finite.
\end{prop}
\duk
We can assume that $f(x)\in\mathrm{RAT}$ and is nonzero.
By Theorem~\ref{thm_oRacFci}, $f(x)=p(x)/q(x)$ where $p(x)$ and 
$q(x)$ are polynomials. Clearly, $Z(f)=Z(p)$. Since $f(x)$ is nonzero, $p(x)\ne k_0(x)$.
By Exercise~\ref{ex_zeroCan} or Corollary~\ref{cor_bouNumRoo}, the set $Z(p)=Z(f)$ is finite.
\kduk
\vspace{-3mm}
\begin{exer}\label{ex_zeroNonz}
The square of a~nonzero function is a~nonzero function.    
\end{exer}

\begin{cor}[nonzero derivative]\label{prop_nenulFce}
Let $r(x)\in\mathrm{RAC}$, $c_1$, $c_2$, $\ds$, $c_k$ be $k\in\N$ real numbers, not all of them zero, and let
$${\textstyle
f(x)=r(x)+\sum_{i=1}^k c_i\log(x-i+1)\,.
}
$$
\underline{Then} $f'(x)$ is a~nonzero restriction of a~rational function.
\end{cor}
\duk
We may assume that $c_k\ne0$. Clearly, 
$${\textstyle
f'(x)=\big(r'(x)+\sum_{i=1}^k\frac{c_i}{x-i+1}\big)\,|\,(k-1,\,+\infty)
}
$$
is a~restriction of a~rational function. Let $b:=k-1$. If 
$b$ is a~regular point of $r'(x)$, then by Exercise~\ref{ex_onPoles} $b$ is 
a~pole of $f'(x)$ of order $1$. Else the same exercise and 
Proposition~\ref{prop_poleRatDer} give that $b$ is a~pole of $f'(x)$ of order at 
least $2$. In either case we see that $f'(x)$ is nonzero.
\kduk

\noindent
{\em $\bullet$ Proof of Theorem~\ref{thm_logNotPrek}. }We 
prove the following stronger result.

\begin{thm}[finitely many zeros]\label{thm_konMnoKor}
Every\index{theorem!finitely many zeros|emph} function of the form
$$
{\textstyle
f(x)=r(x)\,|\,(c,\,+\infty)+
\sum_{j=1}^k p_j(x)\log(x-j+1)\,,
}
$$
where $r(x)\in\mathrm{RAC}$, $c\in\R$, $k\in\N$, $p_j(x)\in\mathrm{POL}$ and some $p_j(x)$ is nonzero, has finitely many zeros.
\end{thm}
\duk
For every such function $f=f(x)$ we define the degree  as 
$${\textstyle
\deg(f)=\min\big(\sum_{\substack{j\in[k]\\
p_j(x)\ne k_0(x)}}\deg(p_j)\label{degree2}\big)
}
$$
where the minimum is
taken over all representations of $f(x)$ in the stated form.
We argue by contradiction and consider a~function $f_0(x)$ with 
the minimum degree and infinitely many zeros. We take an infinite 
and strictly monotone sequence $(a_n)\sus Z(f_0)$
(Exercise~\ref{strMono}). We may assume that for every $n\in\N$, the
interval $(a_n,a_{n+1})$, respectively $(a_{n+1},a_n)$, is 
contained in $M(f_0)$ (Exercise~\ref{ex_procPotr}).
Using Theorem~\ref{thm_Rolle} we get a~sequence $(b_n)\sus Z(f_0')$ 
such that
$$
a_1<b_1<a_2<b_2<a_3<\ds,\,\text{ respectively }\,a_1>b_1>a_2>b_2>a_3>\ds
$$ 
(Exercise~\ref{ex_jakRolle}). 
In particular, the set $Z(f_0')$ is infinite. This is a~contradiction. If $\deg(f_0)=0$, we contradict 
Proposition~\ref{prop_nonzRestr}
and Corollary~\ref{prop_nenulFce}. If $\deg(f_0)>0$, we contradict
the minimality of $\deg(f_0)$ because $f_0'(x)$ is a function of the
considered form that has infinitely many zeros and has degree smaller than $f_0(x)$ 
(Exercise~\ref{ex_procTypuu}).
\kduk
\vspace{-3mm}
\begin{exer}\label{strMono}
How do we select in the infinite set $Z(f)$ an
increasing, or a~decreasing, sequence $(a_n)$?
\end{exer}

\begin{exer}\label{ex_procPotr}
Why can we assume that the gaps between consecutive terms in $(a_n)$
are contained in $M(f)$?
\end{exer}

\begin{exer}\label{ex_jakRolle}
How do we exactly apply Rolle's theorem to $f$ and $(a_n)$ so that we get the interleaving zeros $(b_n)$ of $f'$?    
\end{exer}

\begin{exer}\label{ex_procTypuu}
Show that for $\deg(f)>0$ the derivative $f'$ has the considered form and $\deg(f')<\deg(f)$.  
\end{exer}

\noindent
{\bf Proof of Theorem~\ref{thm_logNotPrek}. }Suppose, 
for the contrary, that the sequence $(\log n)$ is 
P-recurrent. Then there exist $k\ge1$ polynomials $p_1(x)$, $p_2(x)$, $\ds$, $p_k(x)$, not all equal to $k_0(x)$, such that for every $n\ge k$,
$${\textstyle
\sum_{j=1}^k p_j(n)\log(n-j+1)=0\,.
}
$$
Hence the function
$$
{\textstyle
f(x)=\sum_{j=1}^k p_j(x)\log(x-j+1)
}
$$
has infinitely many zeros: $Z(f)\supset\{k,k+1,\ds\}$. This 
contradicts Theorem~\ref{thm_konMnoKor} (Exercise~\ref{ex_whyIsOf}).
\kduk
\vspace{-3mm}
\begin{exer}\label{ex_whyIsOf}
Check that $\sum_{j=1}^k p_j(x)\log(x-j+1)$ belongs to the family of
functions considered in Theorem~\ref{thm_konMnoKor}.
\end{exer}

\begin{exer}\label{ex_zobecNonPrek}
Generalize Theorem~\ref{thm_logNotPrek} to sequences $(\log(n+c))$ for any real $c>-1$. 
\end{exer}

\begin{exer}\label{ex_infManZer}
Prove by means of Rolle's theorem the following proposition.    
\end{exer}

\begin{prop}\label{prop_infManZer}
Let $I\sus\R$ be a~nontrivial interval, 
$f\in\mathcal{C}(I)$, $f'(x)\in\R^*$ for every $x\in I$ and let $Z(f)$ be infinite. \underline{Then} $Z(f')$ is infinite.    
\end{prop}

\section[${}^c$Cantor's transcendental numbers]{Cantor's transcendental numbers}\label{sec_tranNumb}

We prove by means of Theorem~\ref{thm_Lagrange} the existence of transcendental numbers. In the next section 
we give a~simpler argument.  

\medskip\noindent
{\em $\bullet$ Recursive real numbers. }For simplicity of notation, we only work with real numbers in $[0,1]$. We denote by 
$(n)_{10}$ the natural number $n\in\N$ written in the base 
$10$ as a~word over the alphabet $\{0,1,\ds,9\}$. 
For example, $(2^{10})_{10}=1024$. Let 
$$
\N_{10}:=
\{(n)_{10}\cc\;n\in\N\}\,. 
$$

\begin{exer}\label{ex_onN10}
The map $\N\ni n\mapsto(n)_{10}\in\N_{10}$ is 
injective.
\end{exer}

Algorithms cannot work directly with natural numbers; they work with their codes, which are words over a~finite 
alphabet. An example of these codes are the words $(n)_{10}$.

\begin{defi}[recursive reals]\label{def_recReal}
If a~map 
$$
\mathcal{A}\cc\N_{10}\to
\{0,\,1,\ds,9\}
$$ 
is given by an algorithm, that is by a~Turing machine\index{Turing machine}, we call the sum
$$
{\textstyle
\kappa(\mathcal{A}):=\sum_{n\ge1}\mathcal{A}((n)_{10})\cdot 10^{-n}\ \ (\in[0,\,1])\label{kappaA}
}
$$
a~\underline{recursive real number\index{real numbers, $\R$!recursive ones|emph}} in the interval $[0,1]$.
\end{defi}

\begin{exer}\label{ex_someNonrec}
Show that some real numbers in $[0,1]$ are not recursive.    
\end{exer}

\noindent
{\em $\bullet$ Algebraic and transcendental numbers. }We remind these two sets of numbers. We denote by $\Q[x]$ the
set of \underline{rational polynomials\index{polynomials, POL!rational|emph}}; $\Q[x]$\label{racPol}
is a~subset of $\mathrm{POL}$ obtained in Definition~\ref{def_polynomy} by allowing only rational 
constant functions, that is, functions $k_c(x)$ for $c\in\Q$. The set $\Z[x]$\label{intePol} of \underline{integral 
polynomials\index{polynomials, POL!integral|emph}} is defined similarly.

\begin{defi}[algebraic numbers]\label{def_algNum}
A~complex number number (see 
$\al\in\C$ is \underline{algebraic\index{algebraic
number|emph}} if $p(\al)=0$ for a~nonzero polynomial $p(x)\in\Q[x]$.
\end{defi}

\begin{exer}\label{ex_ExaAlgNum}
All fractions and roots $\sqrt{n}$, $n\in\N_0$, are algebraic numbers.
\end{exer}

\begin{exer}\label{ex_ExaAlgNum1}
The polynomial $p(x)$ in Definition~\ref{def_algNum}
can always be modified so that its degree is preserved and {\em (i)} it remains rational
but becomes monic (with leading coefficient $1$) or {\em (ii)} becomes integral.
\end{exer}
Algebraic numbers for which both forms (i) and (ii) are
simultaneously achievable, that is, roots of monic integral polynomials, 
are called \underline{algebraic 
integers\index{algebraic integer|emph}}.

\begin{exer}\label{ex_FraAlgInt}
Which fractions are algebraic
integers?     
\end{exer}

\begin{exer}\label{ex_GRAlgInt}
Is the golden ratio 
$\phi\equiv\frac{1+\sqrt{5}}{2}$ an algebraic integer?     
\end{exer}

\begin{exer}\label{ex_goldFibo}
How is the golden ratio related to the Fibonacci numbers?   
\end{exer}

\begin{defi}[transcendental numbers]\label{def_traNum}
We say that a~complex number is 
\underline{transcendental\index{transcendental numbers|emph}} if it is not algebraic.    
\end{defi}

\noindent
{\em $\bullet$ Cantor's proof.
}In 1870s, G.~Cantor\index{Cantor, Georg} found a~simple proof of 
the existence of (real) transcendental numbers. This was 
a~great achievement of the nascent set theory, which we recall here.

\begin{thm}[Cantor]\label{thm_Cantor}
Real\index{theorem!Cantor's on existence of transcendental 
numbers|emph} transcendental numbers form
an uncountable set.    
\end{thm}
\duk
Let $T$ ($\sus\R$) be the set of real
transcendental numbers. We assume for the contrary that $T$ is at most countable. 
Exercise~\ref{ex_CantTrans1} shows that the set $A\cap\R$ of real algebraic numbers is countable. It follows that the set
$$
\R=T\cup(A\cap\R)
$$
is countable. This contradicts 
\kduk
\vspace{-3mm}
\begin{exer}\label{ex_CantTrans1}
Show\index{transcendental numbers!Cantor's proof of their existence} 
that the set $A$ of algebraic numbers is countable.    
\end{exer}

\noindent
{\em $\bullet$ An effective version of Cantor's proof. }We extract from the 
previous proof an algorithm 
$$
\mathcal{A}\cc\N_{10}\to
\{0,\,1,\,\ds,\,9\}
$$ 
producing a~recursive transcendental number $\kappa(\mathcal{A})$. We define 
$\mathcal{A}$ with the help of the following proposition, which in turn we obtain by means of
Theorem~\ref{thm_Lagrange}. The proposition is an effective version, for integral polynomials, of the fact that any nonzero value 
of a~continuous 
function has a~neighborhood on which the function does not vanish. 

\begin{prop}[$p(x)\ne0$ on $I$]\label{prop_nonvanPoly}
Let 
$$
p(x)=a_nx^n+\ds+a_1x+a_0,\ 
n\in\N_0,\ a_i\in\Z\text{ and }\, 
a_n\ne0\,, 
$$
be a~nonzero integral polynomial with degree $n$. Let $\al=a/10^k$, for 
$a\in\Z$ and $k\in\N_0$, be a~decimal fraction such that 
$p(\al)\ne0$. Let
$$
b=(n+1)^2\cdot
\max(\{|a_0|,\,|a_1|,\,
\ds,\,|a_n|\})\cdot(|a|+1)^n\,\text{ and }\,l=kn+b\ \ (\in\N)\,.
$$
\underline{Then} $p(x)\ne0$ for every real number $x\in[\al,
\al+10^{-l}]$.
\end{prop}
\duk
Since $\al$ has denominator $10^k$ 
and $p(x)$ has integral coefficients and degree $n$, the assumption that 
$p(\al)\ne0$ implies that $|p(\al)|\ge10^{-kn}$. For every 
$x\in(\al,\al+1]$, 
Theorem~\ref{thm_Lagrange} gives a~$y\in(\al,x)$ such that
$$
p(x)=p(\al)+p'(y)\cdot(x-\al)\,.
$$
By Exercise~\ref{ex_deriBound},  $|p'(y)|\le b$ where $b$ is defined above. If $x-\al\le10^{-l}=10^{-kn-b}<
10^{-kn}/b$ then
$${\textstyle
|p(x)|\ge|p(\al)|-|p'(y)|\cdot(x-\al)>10^{-kn}-b\cdot\frac{10^{-kn}}{b}=0\,,
}
$$
as stated.
\kduk
\vspace{-3mm}
\begin{exer}\label{ex_deriBound}
Prove the estimate $|p'(y)|\le b$.    
\end{exer}

\begin{thm}[effective Cantor's proof]\label{thm_CantTrans2}
Cantor's\index{theorem!effective Cantor's proof|emph}
proof yields an algorithm
$\mathcal{A}\cc\N_{10}\to
\{0,1,\ds,9\}$ such that the real number 
$${\textstyle
\kappa:=\kappa(\mathcal{A})=\sum_{n\ge1}\mathcal{A}((n)_{10})\cdot 10^{-n}\label{kappaNum}
}
$$
is not a~root of any nonzero integral polynomial. Hence $\kappa$ is a~recursive real transcendental number. 
\end{thm}
\duk
It is clear that there is an algorithm $\mathcal{B}\cc\N_{10}\to\bigcup_{n=1}^{\infty}\Z^n=:Z$ such that
$$
\mathcal{B}[\N_{10}]=
\{(a_0,\,a_1,\,\ds,\,a_n)\in Z\cc\;n\in\N_0\wedge a_n\ne0\}
$$
---\,$\mathcal{B}$ lists all nonzero integral polynomials by the tuples of their coefficients. For simplicity of notation we take the elements of 
$Z$ directly; in practice they are encoded as words over 
a~finite alphabet. For $m\in\N$ we set
$${\textstyle
p_m(x):=\sum_{j=0}^{n_m}a_{j,m}x^j
\,\text{ where }\,(a_{0,\,m},\,
a_{1,\,m},\,\ds,\,a_{n_m,\,m})
=\mathcal{B}((m)_{10})\,.
}
$$
The algorithm 
$\mathcal{A}$ calls the algorithm
$\mathcal{B}$ and generates a~sequence of decimal fractions
$${\textstyle
\al_1=\frac{z_1}{10^{k_1}}=\frac{0}{10^0}=0,\,
\al_2=\al_1+\frac{z_2}{10^{k_2}},\,
\al_3=\al_2+\frac{z_3}{10^{k_3}},\,\ds
}
$$
such that $k_j\in\N_0$, $0=k_1<k_2<\ds$ (so that the denominator of $\al_j$ is $10^{k_j}$), $z_j\in
\{0,1,\ds,10^{k_j-k_{j-1}}-1\}$ for every $j\in\N$ (with $k_0:=0$), and that for every $m\in\N$ we have
$$
{\textstyle
\forall x\in\big[\al_m,\,\al_m+10^{-k_m}\big]\cc p_1(x)p_2(x)\ds p_{m-1}(x)\ne0\,,
}
$$
where for $m=1$ the product is defined as $1$.

Suppose that $m\in\N$ and that $\mathcal{A}$ already generated 
$\al_1$, $\al_2$, $\ds$, $\al_m$. To generate the next decimal fraction $\al_{m+1}$, the algorithm 
$\mathcal{A}$ gets the polynomial $p_m(x)$ by calling $\mathcal{B}$ and takes a~large number $k\in\N$ such that 
$$
\text{$k>k_m$ and $10^{k-k_m}>\deg (p_m)$}\,.
$$
By Exercise~\ref{ex_zeroCan} there is 
a~number $j\in\{0,1,\ds,10^{k-k_m}-1\}$ such that 
$${\textstyle
p_m\big(\al_m+\frac{j}{10^k}\big)\ne0\,. 
}
$$
By Exercise~\ref{ex_whyAcan}, 
$\mathcal{A}$ can get this $j$. Then $\mathcal{A}$ invokes 
Proposition~\ref{prop_nonvanPoly} with $p(x)=p_m(x)$ and 
$\al=\al_m+\frac{j}{10^k}$, gets (from Proposition~\ref{prop_nonvanPoly}) the number $l\in\N$ and computes the numbers
$${\textstyle
k_{m+1}:=\max(k,\,l),\,
z_{m+1}:= j\cdot10^{k_{m+1}-k}\,\text{ and }\,\al_{m+1}:=\al_m+\frac{z_{m+1}}{10^{k_{m+1}}}\,.
}
$$ 
Since the interval $I:=\big[\al_{m+1},\al_{m+1}+10^{-k_{m+1}}\big]$ is contained in both intervals
$$
{\textstyle
\big[\al_m,\,\al_m+10^{-k_m}\big]\,
\text{ and }\,\big[\al,\,
\al+10^{-l}\big]\,,
}
$$
the value $p_j(x)\ne0$ for every $j=1,2,\ds,m$ and every $x\in I$, as required.

For $j\in\N$, we denote
by $\overline{z_j}$ the word obtained by padding in front
of $(z_j)_{10}$ zeros to the length $k_j-k_{j-1}$. For example, 
if $k_j-k_{j-1}=5$ and $z_j=(z_j)_{10}=27$, then
$\overline{z_j}=00027$. The final outputs $$
\mathcal{A}((n)_{10})\in
\{0,\,1,\,\ds,\,9\},\ n\in\N\,, 
$$
of
$\mathcal{A}$ are defined by the equality of two infinite words over 
$\{0,1,\ds,9\}$: 
$$
\mathcal{A}(1)\,
\mathcal{A}(2)\,
\ds=\overline{z_1}\,\overline{z_2}\,
\cdots\,.
$$
Since
$$
{\textstyle
\kappa=\kappa(\mathcal{A})=\sum_{n=1}^{\infty}\mathcal{A}((n)_{10})\cdot 10^{-n}\in
\bigcap_{m=1}^{\infty}\big[\al_m,\,\al_m+10^{-k_m}\big]\,,
}
$$
the value $p_m(\kappa)\ne0$ for every $m\in\N$. The number $\kappa$ is transcendental.
\kduk
\vspace{-3mm}
\begin{exer}\label{ex_whyAcan}
How does the algorithm $\mathcal{A}$ select the number $j$?    
\end{exer}

\noindent
A~natural idea is to attempt to extend the construction and obtain a~recursive real number that is not a~zero 
of any nonzero function in
$\mathrm{EF}_{\Q}$. This subfamily of elementary functions is generated by
rational constants $k_c(x)$, $c\in\Q$, and rational powers
$x^b$, $b\in\Q\cap(0,+\infty)$.

\section[${}^c$Liouville's transcendental numbers]{Liouville's transcendental numbers}\label{sec_tranNumbLiouv}

 The French mathematician and physicist {\em Joseph
Liouville\index{Liouville, Joseph} (1809--1882)} was the first who proved, in 1844, that 
transcendental numbers exist. In this section we
explain his method. It produces very simple recursive real
transcendental numbers. 

\medskip\noindent
{\em $\bullet$ Liouville's inequality. }This is another application of Theorem~\ref{thm_Lagrange}. 

\begin{thm}[Liouville's inequality]\label{thm_LiouIneq}
Let\index{theorem!Liouville's inequality|emph} 
$\al\in\R\setminus\Q$ be algebraic.
\underline{Then} there exist $n\in\N$ and $c>0$
such that for every $p/q\in\Q$ with $q>0$ we have
$$
\bigg|\al-\frac{p}{q}\bigg|\ge\frac{c}{q^n}\,.
$$
\end{thm}
\duk
Using Exercise~\ref{ex_ExaAlgNum1} we take a~nonzero
$f(x)\in\Z[x]$ with the minimum degree 
$n=\deg(f)\ge2$ such that $f(\al)=0$. Let 
$I=[\al-1,\al+1]$. If 
$\frac{p}{q}\in\Q\setminus I$ 
then 
$${\textstyle
\big|\al-\frac{p}{q}\big|\ge1\ge
1/q^n\,. 
}
$$
If $\frac{p}{q}\in I\cap\Q$ then $\frac{p}{q}\ne\al$ because $\al$ is irrational. By 
Theorem~\ref{thm_Lagrange} there is 
a~real number $x$ lying between $\al$ and $\frac{p}{q}$ such that
(recall that $f(\al)=0$)
$$
{\textstyle
f(\al)-f(\frac{p}{q})=f'(x)(\al-\frac{p}{q})\,\text{ and therefore }\,
\big|\al-\frac{p}{q}
\big|=|f(p/q)|\,/\,|f'(x)|\,.
}
$$
Crucially, $f(\frac{p}{q})\ne0$: if $f(\frac{p}{q})=0$ then
$g(x)=f(x)/(x-p/q)$ 
would be a~rational polynomial 
with $g(\al)=0$ (Exercise~\ref{ex_whyg0}) and 
$\deg g=\deg f-1$, in contradiction with the minimality 
of $\deg f$. Thus $f(\frac{p}{q})\ne0$ and, as we already argued in the previous section, we have $|f(\frac{p}{q})|\ge q^{-n}$.
We take a~real number $d>0$ such that $|f'(y)|\le d^{-1}$ for every $y\in I$ (Exercise~\ref{ex_whyd}). Then 
$$
{\textstyle
\big|\al-\frac{p}{q}
\big|\ge d/q^n\,.
}
$$
We get Liouville's inequality with the constant $c=\min(1,d)$.
\kduk
\vspace{-3mm}
\begin{exer}\label{ex_whyg0}
Why is $\al$ a~root of $g(x)$?    
\end{exer}

\begin{exer}\label{ex_whyd}
How do we get the constant $d$?    
\end{exer}

\noindent
{\em $\bullet$ Liouville's \index{transcendental numbers!Liouville's proof of their existence} 
transcendental
numbers. }The following recursive real transcendental 
number is constructed by Liouville's method.

\begin{cor}[$\lambda$ is transcendental]\label{cor_lambTran}
The\index{lambda@$\lambda=0.11000100\ds$|emph} 
real number
$${\textstyle
\lambda\equiv\sum_{n\ge1}10^{-n!}=0.110001000000000000000001000\ds
}
$$
is recursive and transcendental. 
\end{cor}
\duk
It is clear that $\lambda$ is recursive: the $j$-th decimal 
digit of $\lambda$ after the decimal point is $0$ or $1$, and 
it is $1$ iff $j=n!$
for some $n\in\N$.
The number $\lambda$ is irrational because it does not have eventually
periodic decimal expansion (Exercise~\ref{ex_perDecExp}). For 
$m\in\N$ let $q_m=10^{m!}$ and let $z_m\in\N$ be defined by 
$${\textstyle
\sum_{n=1}^m
\frac{1}{10^{n!}}=
\frac{z_m}{10^{m!}}\,.
}
$$ 
Clearly, $q_m\ge2$. Then, using a~geometric series, we get the bound
$${\textstyle
|\lambda-\frac{z_m}{q_m}|
=\sum_{n=m+1}^{\infty}\frac{1}{10^{n!}}\le
\frac{1}{10^{(m+1)!}}\frac{1}{1-10^{-
(m+1)!}}<2q_m^{-m-1}\,. 
}
$$
Fractions 
$\frac{z_m}{q_m}$  violate Liouville's inequality 
for $\lambda$ (Exercise~\ref{ex_violLiIn}) and $\lambda$ is transcendental.
\kduk
\vspace{-3mm}
\begin{exer}\label{ex_perDecExp}
Show that every rational number has an eventually periodic decimal expansion.     
\end{exer}

\begin{exer}\label{ex_violLiIn}
Explain why fractions $\frac{z_m}{q_m}$ violate for large $m$Liouville's inequality.
\end{exer}

\begin{exer}\label{ex_LamEff}
What is the complexity of the natural algorithm $\mathcal{L}$ that computes the
decimal expansion of $\lambda$?   
\end{exer}

\begin{exer}\label{ex_moreGenLam}
For every $k\in\N$, $k\ge2$, the sum $\sum k^{-n!}$ is transcendental.   
\end{exer}

\section[${}^c$Monotonicity and l'Hospital rules]{Monotonicity and l'Hospital 
rules}\label{sec_monoLHR}

We employ Theorem~\ref{thm_Lagrange} to determine intervals of 
monotonicity of functions, and we prove l'Hospital's rules by which 
one can compute limits. First, we look at the interaction of finite and 
infinite derivatives. 

\medskip\noindent
{\em $\bullet$ Dependence of finite and infinite derivatives. }Equal signs of finite derivatives force the same sign for infinite derivatives. This application
of Theorem~\ref{thm_Lagrange} 
is not so well known. We assume that $a<b$  are real numbers. 

\begin{thm}[same signs of derivatives]\label{proop_finInfder}
Let\index{theorem!same signs of derivatives|emph} 
$f\in\mathcal{C}([a,b])$. If 
$f'(c)\in\R^*$ for every
$c\in(a,b)$ and
every finite $f'(c)$ is nonnegative, \underline{then} 
$f'(c)\ge0$ for every $c\in(a,b)$.
\end{thm}
\duk
Suppose for the contrary that 
$f'(c)=-\infty$ for some $c\in(a,b)$. It follows that there 
exist numbers $c_0$ and $c_1$ such that $a<c_0<c<c_1<b$ and $f(c_0)>f(c)>f(c_1)$.  
By Theorem~\ref{thm_Lagrange} there is a~point $c_2\in(c_0,c_1)$ such that 
$$
{\textstyle
\frac{f(c_1)-f(c_0)}{c_1-c_0}=f'(c_2)\ge0\,.
}
$$
Thus $f(c_1)\ge f(c_0)$, which is a~contradiction.
\kduk
\vspace{-3mm}
\begin{exer}\label{ex_naInterac1}
State and prove the symmetric version of the theorem that forbids the derivative $f'(c)=+\infty$.    
\end{exer}

\begin{exer}\label{ex_naInterac2}
Find $f\in\mathcal{F}([-1,1])$ such that $f'(c)=1$ for every 
$c\in[-1,1]\setminus\{0\}$ and $f'(0)=-\infty$. Can $f$ be 
continuous?
\end{exer}

\noindent
{\em $\bullet$ Derivatives and intervals of monotonicity. }For any $M\sus\R$, the set
$$
M^0:=\{a\in M\cc\;\exists\,\de\cc U(a,\de)\sus M\}\label{interior}
$$ 
is the \underline{{\em interior\index{interior of a 
set|emph}}} of $M$. The interior of an interval arises by deleting 
the endpoints. We introduce the following abbreviations.

\begin{defi}[holding on]\label{def_ineqOn}
Let $A\in\R^*$, 
$f\in\mathcal{F}(M)$ and $X$ be any set. We say that 
$$
\text{\underline{$f\ge A$ (holds) 
on\index{function!inequality@(in)equality holds on|emph} $X$}} 
$$
if
$f(b)\ge A$ for every $b\in M\cap X$.
Similarly for the notation $f\le A$, $f>A$, $f<A$, $f=a$ and $f\ne a$ (holding) on $X$.    
\end{defi}

\begin{thm}[global monotonicity]\label{thm_derMono1}
Suppose that $I\sus\R$ is a~nontrivial interval,  $f\in\mathcal{C}(I)$
 and $f'(c)\in\R^*$ for every $c\in I^0$. The following holds.
\begin{enumerate}
\item If $f'\ge0$ on $I^0$ \underline{then} $f$ weakly increases on $I$. 
\item If $f'\le0$ on $I^0$ \underline{then} $f$ weakly decreases on $I$.
\item If $f'>0$ on $I^0$ \underline{then} $f$ increases on $I$. 
\item If $f'<0$ on $I^0$ \underline{then} $f$ decreases on $I$. 
\end{enumerate}
\end{thm}
\duk
We only prove the last claim and leave the rest to  Exercise~\ref{ex_otherThree}.
Let $f'<0$ on $I^0$ and let $x<y$ be in $I$. By Theorem~\ref{thm_Lagrange} there is
a~point $z\in(x,y)$ (hence $z\in I^0$) such that
$${\textstyle
\frac{f(y)-f(x)}{y-x}=f'(z)<0\,.
}
$$
Thus $f(x)>f(y)$ and $f$ decreases. 
\kduk
\vspace{-3mm}
\begin{exer}\label{ex_otherThree}
Prove parts~1--3 of the theorem.     
\end{exer} 

\begin{prop}[local monotonicity]\label{prop_deriMono2}
Let $f\in\mathcal{F}(M)$ and $b\in M$.
\begin{enumerate}
\item If $f_-'(b)<0$, there is a~$\de$ such that $f>f(b)$ on $P^-(b,\de)$.
\item If $f_-'(b)>0$, there is a~$\de$ such that $f<f(b)$ on $P^-(b,\de)$.
\item If $f_+'(b)<0$, there is a~$\de$ such that $f<f(b)$ on $P^+(b,\de)$.
\item If $f_+'(b)>0$, there is a~$\de$ such that $f>f(b)$ on $P^+(b,\de)$.
\end{enumerate}
In each case, the set $P^{\pm}(b,\de)\cap M$ is infinite.
\end{prop}
\duk
We only prove the last claim and leave the rest to Exercise~\ref{ex_otherThree2}. Let $f_+'(b)>0$. Then there is a~$\de$ such that for every $x\in P^+(b,\de)\cap M$,
$${\textstyle
\frac{f(x)-f(b)}{x-b}>0\,\text{ and therefore }\,f(x)>f(b)\,.
}
$$
The set $P^+(b,\de)\cap M$ is infinite because $b$ is a~right limit 
point of $M$.
\kduk
\vspace{-3mm}
\begin{exer}\label{ex_otherThree2}
Prove parts~1--3 of the proposition. 
\end{exer}

\begin{exer}\label{ex_jenJednop}
Can we say that in each of the four cases in the proposition the set 
$f[P^{\pm}(b,\de)]$ is infinite?    
\end{exer}

\noindent
{\em $\bullet$ Extensions of derivatives. }Let $a<b$ be real
numbers. The next proposition has the same assumptions as Theorem~\ref{thm_Lagrange} (Lagrange's mean value theorem).

\begin{prop}[extending derivatives]\label{prop_rozsDeri}
Let $f\in\mathcal{C}([a,b])$ and suppose that $f'(c)\in\R^*$  for every $c\in(a,b)$. 
\underline{Then} 
$$
a,\,b\in L(D(f)),\ f'(a)=\lim_{x\to a}f'(x)
\,\text{ and }\,
f'(b)=\lim_{x\to b}f'(x)
\,,
$$ 
if the limits exist.
\end{prop}
\duk
We show 
that $b\in L(D(f))$ and that if the limit $B=\lim_{x\to b}f'(x)$ 
exists, then $f'(b)$ exists and $f'(b)=B$. For the point $a$ and 
the derivative $f'(a)$ we argue similarly.

By Theorem~\ref{thm_Lagrange}, $(c,b)\cap 
D(f)\ne\emptyset$ for every $c\in(a,b)$. Thus $b\in 
L(D(f))$. Let $\lim_{x\to b}f'(x)=B$ and an $\ep$ be given. Then $f'[(b-\de,b)]\sus U(B,\ep)$ for some $\de$. 
By Theorem~\ref{thm_Lagrange}, for every 
$x\in(b-\de,b)$ there is a~$y\in (x,b)$ such that 
$${\textstyle
\frac{f(x)-f(b)}{x-b}=f'(y)\in U(B,\ep)\,.
}
$$
Hence $f'(b)=B$.
\kduk
\vspace{-3mm}
\begin{exer}\label{ex_jesteUlExtDe}
Find a~function $f\in\mathcal{C}([a,b])$
such that for every $c\in(a,b]$ the derivative 
$f'(c)$ exists but the limit
$\lim_{x\to b}f'(x)$ does not exist.
\end{exer}

\noindent
{\em $\bullet$ Two l'Hospital rules {\em (HR)}. }We show a~method to compute
$\lim_{x\to A}
\frac{f(x)}{g(x)}$ of the forms
$\frac{0}{0}$ and $\frac{\pm\infty}{\pm\infty}$. It is based on the transformation 
$${\textstyle
\frac{f(x)}{g(x)}=
\frac{f(x)/x}{g(x)/x}\,.
}
$$
We begin with a~local result.

\begin{thm}[local HR]\label{thm_LHP1}
Let\index{theorem!local l'Hospital rule|emph} 
$f,g\in\mathcal{R}$,
$b\in L(M(f/g))\cap M(f)\cap M(g)$ and let $f(b)=g(b)=0$. \underline{Then}
$$
\lim_{x\to b}\frac{f(x)}{g(x)}=\frac{f'(b)}{g'(b)}\ \ (\in\R^*)
$$
if the last ratio is defined.
\end{thm}
\duk
We assume that $f'(b)$ and $g'(b)$ exist and that 
$\frac{f'(b)}{g'(b)}$ is not an indefinite expression. Then by Theorem~\ref{thm_AritLimFce} we have 
$${\textstyle
\lim_{x\to b}\frac{f(x)}{g(x)}=\lim_{x\to b}
\frac{\frac{f(x)-f(b)}{x-b}}{\frac{g(x)-g(b)}{x-b}}=\frac{\lim_{x\to b}
\frac{f(x)-f(b)}{x-b}}{\lim_{x\to b}\frac{g(x)-g(b)}{x-b}}=
\frac{f'(b)}{g'(b)}\,.\hspace{3cm}\Box
}
$$

\noindent
For example, 
$${\textstyle
\lim_{x\to0}\frac{\sin x}{x}=
\frac{\cos 0}{k_1(0)}=
\frac{1}{1}=1\,\text{ and
}\,\lim_{x\to0}\frac{\exp x-1}{x}=\frac{\exp0}{k_1(0)}=
\frac{1}{1}=1\,.
}
$$

We continue with a~global result.
We assume that $a<b$ are real numbers. 

\begin{thm}[global HR]\label{thm_LHP2}
Let\index{theorem!l'Hospital rule 2|emph} 
$f,f',g,g'\in\mathcal{F}((a,b))$ and let $g'\ne0$ on 
$(a,b)$. Suppose that 
$$
\text{$\lim_{x\to a}f(x)=\lim_{x\to a}g(x)=0$ or $\lim_{x\to a}g(x)=\pm\infty$}\,.
$$
\underline{Then} 
$$
\lim_{x\to a}
\frac{f(x)}{g(x)}=\lim_{x\to a}\frac{f'(x)}{g'(x)}\ \ (\in\R^*)
$$
if the last limit exists.
\end{thm}
\duk
Note that $f,g\in\mathcal{C}$. 1. Let
$$
\lim_{x\to a}f(x)=\lim_{x\to a}g(x)=0\,\text{ and }\,
\lim_{x\to a}{\textstyle
\frac{f'(x)}{g'(x)}=K\,.} 
$$
Let an 
$\ep$ be given. We define $f(a)=g(a)=0$, then
$f,g\in\mathcal{C}([a,b))$. Theorem~\ref{thm_Rolle} implies that $g\ne0$ on $(a,b)$ (Exercise~\ref{ex_warum}). There is a~$\de$ such that
$$
(f'/g')[(a,\,a+\de)]\sus U(K,\,\ep)\,.
$$
By Theorem~\ref{thm_Cauchy}, 
for every $x\in(a,a+\de)$ there is
a~$y\in(a,x)$ such that
$${\textstyle
\frac{f(x)}{g(x)}=\frac{f(x)-f(a)}{g(x)-
g(a)}=\frac{f'(y)}{g'(y)}\in U(K,\ep)\,.
}
$$ 
Hence $\lim_{x\to a}\frac{f(x)}{g(x)}=K$.

2. We prove this case later by means of integrals.
\kduk
\vspace{-3mm}
\begin{exer}\label{ex_warum}
Why can we assume that $g\ne0$ on $(a,b)$?    
\end{exer}

\begin{exer}\label{ex_naHR2}
Let $\ep>0$. Compute $\lim_{x\to0}x^{\ep}
\log x$.    
\end{exer}

\noindent
Marquis {\em Guillaume de l'Hospital (1661--1704)\index{hospital@de 
l'Hospital, Guillaume}} published 
in 1696 the historically first textbook of differential calculus 
{\em  Analyse des Infiniment Petits pour l'Intelligence des Lignes 
Courbes}.

\begin{exer}\label{ex_otherDefDom}
Adapt the global {\em HR} for the limit at $b$ and for the definition domains $P(a,\de)$ and $U(A,\de)=U(\pm\infty,\de)$.  
\end{exer}

\noindent
{\em $\bullet$ A~non-example. }We show that the assumption that 
$M(f)=M(g)$ is an interval is
substantial. For any interval $I\sus\R$ we set $I_{\Q}:=I\cap\Q$.\label{IQ}

\begin{thm}[a~non-example]\label{prop_LHPcouExa}
There \index{theorem!counter-example to LHR~2|emph}
exist functions 
$$
f,\,g\in\mathcal{F}((0,\,1)_{\Q})
$$ 
such that 
$\lim_{x\to0}f(x)=\lim_{x\to0}g(x)=0$, $D(f)=D(g)=(0,1)_{\Q}$, 
$f'=g'=1$ on $(0,1)_{\Q}$, hence
$\lim_{x\to0}\frac{f'(x)}{g'(x)}=1$, but that
$$
\lim_{x\to0}\frac{f(x)}{g(x)}=0\,.
$$
\end{thm}
\duk
Let 
$(c_n)\sus(0,1)$ 
be a~sequence of irrational numbers such that
$${\textstyle
\text{$c_0=1>c_1>c_2>\ds>0$, $\lim c_n=0$ and $\lim\frac{c_{n-
1}}{c_n}=1$}\,. 
}
$$
For $\al\in(c_n,c_{n-1})_{\Q}$, $n\in\N$, we set
$$
f(\al)=c_n^2+\al-c_n\,\text{ and }\,g(\al)=\al\,. 
$$
The graph $G_f$
consists of short, straight ``segments'' starting on
the parabola $y=x^2$ and each has slope $1$. The other function is 
the restriction $g(x)=\mathrm{id}(x)\,|\,
(0,1)_{\Q}$. It is clear that $f'=g'=k_1(x)\,|\,
(0,1)_{\Q}$ (Exercise~\ref{ex_whyclear}) and therefore 
$$
\lim_{x\to0}{\textstyle
\frac{f'(x)}{g'(x)}}=\lim_{x\to0}
{\textstyle\frac{1}{1}=1}\,.
$$
For $\al\in(c_n,c_{n-1})_{\Q}$ we have
$$
{\textstyle
0<\frac{f(\al)}{g(\al)}\le
\frac{c_n^2+c_{n-1}-c_n}{c_n}=
c_n+\frac{c_{n-1}}{c_n}-1\to0\ \ (n\to\infty)\,.
}
$$
Thus $\lim_{x\to0}\frac{f(x)}{g(x)}=0$.
\kduk
\vspace{-3mm}
\begin{exer}\label{ex_whyclear}
Why are the derivatives of $f$ and $g$ constantly $1$?  
\end{exer}

\section[Second and higher order derivatives]{Second and higher order derivatives}\label{sec_secoDeri}

We begin with the  definition of higher-order derivatives.

\medskip\noindent
{\em $\bullet$ Derivatives of order $k\in\N_0$. }We 
simply iterate the operation of derivative in 
Definition~\ref{def_deriFunk}. Recall that $\mathcal{R}$ is the set 
of functions $f\cc M\to\R$ with $M\sus\R$.

\begin{defi}[$f^{(k)}(x)$]\label{def_kthDer}
Let $k\in\N$. We define a~unary operation $f^{(k)}$ on $\mathcal{R}$ by the iteration, 
$$
\mathcal{R}\ni f\mapsto 
f^{(k)}\equiv(\ds((f')')'\ds)'\in\mathcal{R}\,,\label{derivK}
$$ 
with $k$ applications of derivative. It is the \underline{derivative of order 
$k$\index{derivative of a function!of order $k$|emph}}. 
We set $f^{(0)}:=f$ and write $f'$ for $f^{(1)}$, $f''$ for $f^{(2)}$ and $f'''$ for $f^{(3)}$.\label{deriv2}\label{deriv3}
\end{defi}
For example, 
$$
(x\sin x)''=(\sin x+x\cos x)'=2\cos x-x\sin x\,. 
$$

\begin{exer}\label{ex_secoDeri}
What is the difference between $f''(b)$ and $(f')'(b)$?    
\end{exer}

\begin{exer}\label{ex_naDerKrad}
Determine the sequences of functions
$${\textstyle
\big((\sin x)^{(n)}\big)\,\text{ and }\,\big((\frac{1}{x})^{(n)}\big),\  
n\in\N_0\,.
}
$$
\end{exer}

\noindent
{\em $\bullet$ Derivatives of derivatives and extremes. }We can read the types of local extremes from the signs of the elements 
$(f')'(b)$ ($\in\R^*$), $b\in\R$.

\begin{prop}[$(f')'(b)$ and local extremes]\label{prop_DruhaDerExtr}
Let $f\in\mathcal{R}$. Suppose that  
$$
U(b,\,\de)\sus D(f),\ f'(b)=0\,
\text{ and that $(f')'(b)\in\R^*$}\,.
$$
If $(f')'(b)<0$ \underline{then} $f$ has at $b$ a~strict 
local maximum. If $(f')'(b)>0$ \underline{then} $f$ has at $b$ a~strict 
local minimum.
\end{prop}
\duk
Let $f$, $b$ and $\de$ be as stated and let $(f')'(b)>0$ (the case with $(f')'(b)<0$ is similar).
By parts~2 and 4 of Proposition~\ref{prop_deriMono2} there is a~$\theta<\de$ such that
for every $x\in P^-(b,\theta)$ (respectively 
$x\in P^+(b,\theta)$) we have $f'(x)<0=f'(b)$ (respectively $f'(x)>0$). By Theorem~\ref{thm_derMono1} the function $f$ decreases on 
$[b-\theta,b]$ and increases on $[b,b+\theta]$. Hence $f$ has at $b$ a~strict local minimum.
\kduk
\vspace{-3mm}
\begin{exer}\label{ex_locExtrSecD}
Show that under the assumptions of the previous proposition and with 
$f''(b)=0$, it is possible that  $f$ does not have at $b$ a~local 
extreme. 
\end{exer}

\begin{exer}\label{ex_givExa}
Find examples for Proposition~\ref{prop_DruhaDerExtr} where $(f')'(b)=\pm\infty$.    
\end{exer}

\noindent
{\em $\bullet$ Convexity and concavity. 
}Convex graphs of functions bulge downward, and concave graphs upward.

\begin{defi}[convex and concave]\label{def_konvKonk}
A~function $f\in\mathcal{F}(M)$ is \underline{convex}, respectively 
\underline{concave}\index{convexity and concavity|emph}, 
if for any three points $a<b<a'$ in $M$ 
we have
$$
\text{$f(b)\le sb+c$, respectively $f(b)\ge sb+c$}\,, 
$$
where 
$$
\text{$y=sx+c$ is the secant
$\kappa(a,\,f(a),\,a',\,f(a'))$ of $G_f$}\,. 
$$
If these inequalities hold as
strict, we speak of \underline{strict convexity}, 
respectively \underline{strict concavity}, of $f$.\index{convexity and concavity!strict|emph}
\end{defi}

\noindent
So, the strict convexity of $f$ means that the middle point $(b,f(b))\in G_f$
always lies below the secant line of $G_f$ passing through the extreme points
$(a,f(a))$ and $(a',f(a'))$. For convex $f$, the point
$(b,f(b))$ may also lie on the secant. Similarly for concavity. 

\begin{exer}\label{ex_convConc1}
Function $f(x)=x^2$ is strictly convex. Function
$f(x)=|x|$ is convex, but not strictly convex. 
Function $f(x)=\log x$ is strictly concave.    
\end{exer}

\begin{exer}\label{ex_convConc2}
(Strict) convexity and (strict) concavity are preserved under restrictions of functions.    
\end{exer}

\begin{exer}\label{ex_convConc3}
A~function $f$ is (strictly) convex $\iff$ $-f$ is (strictly) concave.   
\end{exer} 

\noindent
{\em $\bullet$ Convexity, concavity and continuity. }Convexity and concavity force the existence of one-sided derivatives. We say that a~set $M\sus\R$ is 
\underline{end-free\index{set!end-free|emph}} if it has neither
minimum nor maximum.

\begin{thm}[existence of $f_{\pm}'$]\label{thm_exOnesDer}
Let\index{theorem!existence of $f_-'$ and 
$f_+'$|emph} 
$f\in\mathcal{F}(M)$ be a~convex, respectively concave, function defined on a~nonempty 
and end-free set $M\sus\R$. \underline{Then}, with equal signs, 
for every 
$b\in M\cap L^{\pm}(M)$ there exists finite one-sided derivative 
$$
f_{\pm}'(b)\in\R\,, 
$$
and the function 
$$
f_{\pm}'\cc M\cap L^{\pm}(M)\to\R
$$ 
weakly increases, respectively weakly decreases.
\end{thm}
\duk
We prove that any convex function $f$ has at every point 
$$
b\in M\cap L^-(M)
$$ 
finite 
left-sided derivative $f_-'(b)\in\R$, and then that the function $f_-'$ weakly 
increases. The remaining three cases are treated similarly. By part~1 of Theorem~\ref{thm_limMonFce}, the limit
$$
{\textstyle
\lim_{x\to b^-}\frac{f(x)-f(b)}{x-b}=f_-'(b)
}
$$
exists and is finite because for any $a\in M$ with $a<b$ the function 
$${\textstyle
g(x)\equiv\frac{f(x)-f(b)}{x-b}\,\big|\,[a,\,b) 
}
$$
weakly increases, and 
$${\textstyle
g\le\frac{f(b')-f(b)}{b'-b}
}
$$ 
on $[a,b)$
for any fixed 
$b'\in M$ with $b'>b$\,---\,such number $b'$ exists because $\max(M)$ does not exist. These two properties of the function $g(x)$ easily follow from the convexity of $f$ because 
$${\textstyle
\frac{f(x)-f(b)}{x-b}
}
$$
is the slope of the secant $\kappa(x,f(x),b,f(b))$ of $G_f$, 
and similarly for 
$${\textstyle
\frac{f(b')-f(b)}{b'-b}\,.
}
$$

We prove monotonicity of
$f_-'$. For every two points $b<b'$ in $M\cap L^-(M)$, the inequality $f_-'(b)\le f_-'(b')$ follows again from the convexity of $f$. For every $x,y\in M$ 
with $x<b<y<b'$ we have two inequalities between slopes
$$
{\textstyle
\frac{f(x)-f(b)}{x-b}\le
\frac{f(y)-f(b)}{y-b}\le
\frac{f(y)-f(b')}{y-b'}\,,
}
$$
hence 
$${\textstyle
\frac{f(x)-f(b)}{x-b}\le
\frac{f(y)-f(b')}{y-b'}\,. 
}
$$
This inequality is preserved in the limit transitions $x\to b^-$ and $y\to(b')^-$, and we get 
$$
f_-'(b)\le f_-'(b')\,.
$$
\kduk

\noindent
But the two-sided derivative $f'(b)$ may not always exist for concave or convex function $f$ because it may be that $f'_-(b)\ne f'_+(b)$. Consider, for example, the function $|x|$ at the point $0$.

\begin{cor}[implied continuity]\label{cor_implConti}
Let $M\sus\R$ be a~nonempty end-free set and $f\in\mathcal{F}(M)$ 
be convex or concave. \underline{Then} $f\in\mathcal{C}(M)$.
\end{cor}
\duk
We show that $f$ is 
right-continuous at every point $b\in M$. The left-continuity is proven similarly. It 
follows by Exercise~\ref{ex_ekvivSpoji} that $f$ is continuous at $b$. If $b\in M$ but is 
not the right limit point of 
$M$, the function $f$ is 
right-continuous at $b$ trivially. If 
$$
b\in M\cap L^+(M)\,,
$$
then by the previous theorem the 
one-sided derivative $f_+'(b)\in\R$ exists. Thus by 
Exercise~\ref{ex_oneSiDiCo}, $f$ is 
right-continuous at $b$. 
\kduk
\vspace{-3mm}
\begin{exer}\label{ex_ConvEndpo}
Prove the next proposition.     
\end{exer}

\begin{prop}[derivatives at endpoints]\label{prop_ConvEndpo}
Let $M\sus\R$, $b=\max(M)$,
$b\in L^-(M)$ and let $f\in\mathcal{F}(M)$ be convex or concave. 
\underline{Then} there exist the one-sided and possibly infinite derivative 
$$
f_-'(b)\in\R^*\,. 
$$
The same holds if we replace max with min and the sign $-$ with $+$.
\end{prop}
\vspace{-3mm}
\begin{exer}\label{ex_ConvCont}
Is it true that if $I\sus\R$ is a~nontrivial interval and 
$f\in\mathcal{F}(I)$ is a~convex or concave function, then $f$ is continuous?
\end{exer}

\noindent
{\em $\bullet$ Convexity, concavity and $f''$. }Convex and concave 
parts of the graph $G_f$ can 
be determined by means of $f''$. Recall Definition~\ref{def_ineqOn}.

\begin{thm}[convexity, concavity and $(f')'(b)$]\label{thm_druDerKonvKonk}
Let\index{theorem!convexity, concavity and $f''$|emph} 
$I$ be a~nontrivial real interval, $f\in\mathcal{C}(I)$ and $D(f)\supset I^0$. We suppose that for every $c$ in $I^0$ the derivative  $(f')'(c)\in\R^*$ exists. The following holds.
\begin{enumerate}
\item If $f''\ge0$ on $I^0$
\underline{then} $f$ is convex.
\item If $f''>0$ on $I^0$
\underline{then} $f$ is strictly convex.
\item If $f''\le0$ on $I^0$
\underline{then} $f$ is concave.
\item If $f''<0$ on $I^0$
\underline{then} $f$ is strictly concave.
\end{enumerate}
\end{thm}
In the proof we use the next lemma whose proof we leave as an exercise.

\begin{exer}\label{ex_proofLem}
Prove the next lemma.    
\end{exer}

\begin{lemma}[on slopes]\label{lem_oSklo}
Let $(a,a')$, $(b,b')$ and $(c,c')$ be in $\R^2$, $a<b<c$ and 
$${\textstyle
\frac{b'-a'}{b-a}\le\frac{c'-b'}{c-b}\,.
}
$$
\underline{Then} the point $(b,b')$ 
lies below or on the line $\kappa(a,a',c,c')$. We get analogous results if the inequality $\le$ is 
replaced with  any of the three inequalities $\{<,\,\ge,\,>\}$. 
\end{lemma}

\noindent
{\bf Proof of Theorem~\ref{thm_druDerKonvKonk}. 
}Let $f$ and $I$ be as stated and let $f''\ge0$ on $I^0$, the remaining
three cases are treated similarly. Let
the three points $a<b<c$ be in $I$. By Theorem~\ref{thm_Lagrange} there exist numbers $y\in(a,b)$ and $z\in(b,c)$ such that 
$${\textstyle
s=\frac{f(b)-f(a)}{b-a}=f'(y)\,\text{ and }\,t=\frac{f(c)-f(b)}{c-b}=f'(z)\,.
}
$$
By Theorem~\ref{thm_derMono1}, $f'$ weakly increases on $I^0$ because
$f''\ge0$ on $I^0$. From $y<z$ it follows that $s=f'(y)\le f'(z)=t$. By
Lemma~\ref{lem_oSklo} the point $(b,f(b))$ lies below or on the line
$\kappa(a,f(a),c,f(c))$.
So $f$ is convex by Definition~\ref{def_konvKonk}.
\kduk

\noindent
{\em $\bullet$ Inflection points. }At an inflection point the graph of a~function crosses the tangent line.

\begin{defi}[inflection points]\label{def_inflexe}
Let $f\in\mathcal{F}(M)$, $b\in M\cap L^{\mathrm{TS}}(M)$, $\ell(x)$ in $\mathcal{F}(\R)$ be the tangent to $G_f$ at 
$(b,f(b))$ and $z\in\{0,1\}$. If there is a~$\de$ such that
$$
\text{$(-1)^z(f-\ell)\le0$ on $(b-\de,\,b)$ and $(-1)^z(f-\ell)\ge0$ on $(b,\,b+\de)$}\,,
$$
we call $(b,f(b))$ an \underline{inflection point\index{inflection point|emph}} 
of $G_f$. If these inequalities hold as strict,
we call $(b,f(b))$ a~\underline{strict inflection point\index{inflection 
point!strict|emph}} of $G_f$.
\end{defi}
We see that the tangent line at an inflection point is a~particular 
case of the cutting tangent of Section~\ref{sec_3vetyOstrhod}.

\begin{exer}\label{ex_onx3}
The origin $(0,0)$ is a~strict inflection point of the graph of the function
$f(x)=x^3$.     
\end{exer}

\begin{exer}\label{ex_triviInfl}
Which points of the graph of the constant function $k_1(x)$ are inflection points?
\end{exer}

\begin{thm}[no inflection]\label{thm_infl}
Suppose\index{theorem!no inflection|emph} 
that $f\in\mathcal{R}$, $D(f)\supset U(b,\de)$, that the derivative $(f')'(b)\in\R^*$ exists and that $(f')'(b)\ne0$. \underline{Then} 
$$
\text{$(b,\,f(b))$ is not an inflection point of $G_f$}\,.
$$
\end{thm}
\duk
Let $(f')'(b)>0$, the case with $(f')'(b)<0$ is similar. Let 
$\ell$ ($\sus\R^2$) be the tangent to $G_f$ at 
$(b,f(b)$. By 
Proposition~\ref{prop_deriMono2} there is a~$\theta\le\de$ such that 
for every $x\in(b-\theta,b)$ and every $x'\in(b,b+\theta)$ we have inequalities
$$
f'(x)<f'(b)<f'(x')\,.
$$
Let $x\in(b-\theta,b)$, $x'\in 
(b,b+\theta)$ and let $s$ and $t$ be the respective slopes of the 
secants
$$
\kappa(x,\,f(x),\,b,\,f(b))\;\text{ and }\;\kappa(b,\,f(b),\,x',\,f(x'))
$$
of $G_f$. The inequalities and Theorem~\ref{thm_Lagrange} give that 
$$
s<f'(b)<t\,. 
$$
Thus both points
$(x,f(x))$ and $(x',f(x'))$ lie above $\ell$. The condition in 
Definition~\ref{def_inflexe} is not satisfied.
\kduk

We obtain a~sufficient condition for existence of inflection points.

\begin{thm}[inflection exists]\label{thm_postacInflexe}
Let\index{theorem!inflection exists|emph} 
$f\in\mathcal{R}$,  $M(f'')\supset U(b,\de)$ and $z$ be in $\{0,1\}$. 
If 
$$
\text{$(-1)^zf''\ge0$ on $(b-\de,\,b)$ and $(-1)^zf''\le0$ on $(b,\,b+\de)$
} 
$$
\underline{then} 
$$
\text{$(b,f(b))$ is an inflection point of $G_f$\,.}
$$
If these inequalities hold strictly \underline{then} $(b,f(b))$ is a~strict inflection point 
of $G_f$.
\end{thm}
\duk
Let $f$, $b$, $\de$ and $z$ be as stated. We assume that $f''\le0$ on
$(b-\de,b)$ and $f''\ge0$ on
$(b,b+\de)$, the other three cases are similar. We have
$M(f')\supset U(b,\de)$ and denote the tangent line to $G_f$ at 
$(b,f(b))$ by $\ell$ ($\sus\R^2$). By Theorem~\ref{thm_derMono1} the 
derivative $f'$ weakly decreases, respectively increases, on $[b-\de,b]$, respectively
$[b,b+\de]$. Thus for every $x\in[b-\de,b)$ and $x'\in(b,b+\de]$, 
$$
f'(x)\ge f'(b)\le f'(x')\,.
$$
By Theorem~\ref{thm_Lagrange}, 
$$
{\textstyle
\frac{f(b)-f(x)}{b-x}\ge f'(b)\le
\frac{f(x')-f(b)}{x'-b}\,.
}
$$
So the slopes of the lines 
$$
\kappa(x,\,b,\,f(x),\,f(b))\,\text{ and }\,\kappa(b,\,x',\,f(b),\,f(x'))
$$ 
are at least the slope $f'(b)$ of $\ell$. So the point $(x,f(x))$ lies below or on $\ell$, and $(x',f(x'))$ above or on $\ell$. Hence
$(b,f(b))$ is an inflection point.
\kduk

\section[How to draw graphs of functions]{How to draw graphs of functions}\label{sec_drawGr}

We describe thirteen steps for determining main geometric features 
of the graph of a~function. First we define asymptotes.

\medskip\noindent
{\em $\bullet$ Asymptotes. }The graph gets arbitrarily close to these lines.

\begin{defi}[vertical asymptotes]\label{def_vertAsymp}
If for $f\in\mathcal{F}(M)$ and $b\in
L^-(M)$ we have 
$$
\lim_{x\to b^-}f(x)=\pm\infty\,, 
$$
we call the line $x=b$ a~\underline{left vertical 
asymptote\index{asymptote!left vertical|emph}} of $f$. By replacing 
the two signs ${}^-$ by two signs ${}^+$ we obtain the 
\underline{right 
vertical asymptotes\index{asymptote!right vertical}} of $f$.
\end{defi}
\vspace{-3mm}
\begin{exer}\label{ex_asymptote1}
The axis $y$ is both a~left and right vertical asymptote of $\frac{1}{x}$. It is a~right vertical asymptote of 
$\log x$.    
\end{exer}

\begin{defi}[asymptotes at infinity]\label{def_AsympInIn}
Let $s,b$ be in $\R$,  $f$ in $\mathcal{F}(M)$ and let $\pm\infty\in L(M)$. If, with equal signs,
$$
\lim_{x\to\pm\infty}(f(x)-sx-b)=0\,,
$$
we call the line $\ell(x)=sx+b$ an \underline{asymptote (of $f$) at
$\pm\infty$\index{asymptote!at infinity|emph}}. 
\end{defi}

\begin{exer}\label{ex_asymptote2}
The line 
$$
y=sx+b
$$ 
is an asymptote of a~function $f$ at $\pm\infty$ $\iff$ 
$$
\lim_{x\to\pm\infty}\frac{f(x)}{x}=s\,\text{ and }\,\lim_{x\to\pm\infty}(f(x)-
sx)=b
$$ 
(equal signs). 
\end{exer}

\begin{exer}\label{ex_asymptote3}
Find the asymptote of $\frac{1}{x}$ at $+\infty$ and at $-\infty$.
\end{exer}
Definition~\ref{def_vertAsymp} and Exercise~\ref{ex_asymptote2} imply that asymptotes are unique. 

\medskip\noindent
{\em $\bullet$ Drawing graphs of functions. 
}We recall from Definition~\ref{def_obecEF2} that an elementary function is obtained from constant functions $\{k_c(x)\cc\;c\in\R\}$ and functions 
$$
\text{$\exp x$, $\log x$, 
$\sin x$, $\arcsin x$ and $x^b$ for $b\in(0,+\infty)\setminus\N$}\,, 
$$
by
repeated addition, multiplication, division, and composition. In simple 
elementary functions we omit from the initial generators the functions 
$\arcsin x$ and $x^b$. Let $F\in\mathcal{R}$. We determine the main geometric features of $G_F$. 

\medskip\noindent
{\em {\bf Step 0. }Is it (simple) elementary? }Is $F\in\mathrm{EF}$? Is $F\in\mathrm{SEF}$?
Memberships of $F$ in EF and SEF have bearing on $M(F)$, the continuity of $F$ and $D(F)$. 

\medskip\noindent
{\em {\bf Step 1. }The definition domain. }We find $M(F)$ ($\sus\R$). If 
$F\in\mathrm{EF}$, we start from the domains
$$
M(\mathrm{e}^x)=M(\sin 
x)=M(k_c)=\R\,,
$$
$M(x^b)=[0,+\infty)$, $M(\log x)=(0,+\infty)$ and $M(\arcsin x)=
[-1,1]$, and use the domain relations 
$$
M(f+g)=M(fg)=M(f)\cap M(g),\ 
M(f/g)=M(f)\cap M(g)\setminus Z(g)
$$
and 
$$
M(f(g))=\{x\in M(g)\cc\;g(x)\in M(f)\}\,.
$$

\medskip\noindent
{\em {\bf Step 2. }Is it special? }Is $F$ \underline{even\index{function!even|emph}} 
($F(-x)=F(x)$), \underline{odd\index{function!odd|emph}} 
($F(-x)=-F(x)$) or 
$c$-\underline{periodic\index{function!c 
periodic@$c$-periodic|emph}} ($F(c+x)=F(x)$)?

\begin{exer}\label{ex_drawGr1}
Define these families of functions in more detail.     
\end{exer}

\noindent
{\em {\bf Step 3. }Derivatives and continuity. }We determine
for each point $a$ in $M(F)$ if
$F$ is continuous at $a$, and if 
$F'(a)\in\R^*$ exists and is finite. By this we determine $F'$. If $F\in\mathrm{EF}$ then $F\in\mathcal{C}$ by
Theorem~\ref{thm_EFjsouSpoj}. If $F\in\mathrm{SEF}$ then $D(F)=M(F)$ by Theorem~\ref{thm_deriSEF}.

\medskip\noindent
{\em {\bf Step 4. }Limits. }At any point 
$a\in M(F)$ where $F$ is discontinuous we investigate one-sided limits. We also investigate limits of $F$ at the 
elements of $L(M(F))\setminus 
M(F)$. For example, 
$\lim_{x\to0^{\pm}}\sgn\,x=\pm1$.
Or $\lim_{x\to-\infty}\exp x=0$ and 
$\lim_{x\to+\infty}\exp x=+\infty$.

\medskip\noindent
{\em {\bf Step 5. }Intersections of the graph with coordinate axes. }These are the points
$(x,0)$ ($\in\R^2$) where $x\in Z(F)$, 
plus the point $(0,F(0))$ if $0\in M(f)$.

\medskip\noindent
{\em {\bf Step 6. }One-sided derivatives. }At any point $a\in M(F)$ where $F'(a)$ does not exist, we investigate $F_-'(a)$ and
$F_+'(a)$. Proposition~\ref{prop_rozsDeri} is 
relevant. For example, it implies that $(|x|)'_-(0)=-1$ and
$(|x|)'_+(0)=1$. However, these one-sided derivatives of the absolute value are easily computed
directly.

\medskip\noindent
{\em {\bf Step 7. }Maximal intervals of monotonicity and extremes. }We find all inclusion-wise maximal intervals $I$ 
($\sus\R$) on which (that is, on $I\cap M(F)$) the function $F$ is monotone. 
For elementary functions,
Theorem~\ref{thm_derMono1} can be often used.
More ambitiously, one can try
to determine all inclusion-wise maximal subsets of $M(F)$ where $F$ 
is monotone, we call them
$$
\text{\underline{maximal domains of 
monotonicity\index{maximal domains of monotonicity, MDM|emph}}}\,, 
$$
or MDM.\label{MDM} 
For functions more general than elementary, this can
be difficult. We consider MDM only in Exercise~\ref{ex_GraDraw1apul}. We find local and global extremes of $F$. For this Corollary~\ref{cor_priMaxi}, 
Theorem~\ref{thm_priznakExtr}, and Proposition~\ref{prop_DruhaDerExtr} are 
relevant. 

\medskip\noindent
{\em {\bf Step 8. }The image. }The image of $F$ is the set $F[M(F)]$ ($\sus\R$). 

\medskip\noindent
{\em {\bf Step 9. }Maximal intervals of convexity and concavity. }We find 
all inclusion-wise maximal intervals $I$ ($\sus\R$) on which $F$ is convex or concave. For elementary
functions, Theorem~\ref{thm_druDerKonvKonk} can
usually be used. We will not consider the harder problem of 
finding all inclusion-wise maximal subsets of $M(f)$ where $F$ is
convex or concave. 

\medskip\noindent
{\em {\bf Step 10. }Inflection points. }We find these points of $G_F$. 
Theorems~\ref{thm_infl} and \ref{thm_postacInflexe} are 
relevant.

\medskip\noindent
{\em {\bf Step 11. }Asymptotes. }We find asymptotes of $F$. 
Definitions~\ref{def_vertAsymp} and \ref{def_AsympInIn}, and
Exercise~\ref{ex_asymptote2} are relevant.

\medskip\noindent
{\em {\bf Step 12. }Sketching the graph. }Using hand, the
computer or the Internet we sketch the usually uncountable set 
$$
G_F=\{(a,\,F(a))\cc\;a\in M(F)\}\ \ (\sus\R^2)\,.
$$

\medskip\noindent
{\em $\bullet$ First example. }Let 
$$
F=F(x)\equiv\sgn\,x\,. 
$$
Recall that $F(x)=-1$ for $x<0$, $F(x)=1$ for $x>0$ and $F(0)=0$. {\bf Step~0. }$F\not\in\mathrm{EF}$ by Proposition~\ref{prop_sgnNeniEF}.
{\bf Step~1. }$M(F)=\R$. {\bf Step~2. }The function $F$ is 
odd. {\bf Step~3. }$F$ is continuous at every point $x\ne0$ 
and $F'(x)=0$ for every $x\ne0$. At the point $0$ the function $F$ is discontinuous and 
$F'(0)=+\infty$. {\bf Step~4. }We have the one-sided limits 
$\lim_{x\to0^-}F(x)=-1$ and $\lim_{x\to0^+}F(x)=1$. Also, 
$\lim_{x\to-\infty}F(x)=-1$ and 
$\lim_{x\to+\infty}F(x)=1$. {\bf Step~5. }$G_F$ intersects 
both coordinate axes only at the origin $(0,0)$. {\bf Step~6. }Since $F'(x)$ exists for every $x\in\R$, there is nothing to compute: always 
$$
F'_-(x)=F'_+(x)=F'(x)\,. 
$$
{\bf Step~7. }We see from 
the definition of $F$ that $F$ weakly increases on $\R$, that $x$ 
is  a~global minimum, respectively maximum, of $F$ iff $x<0$, respectively $x>0$, and that $F$ has no strict local or global extreme.
{\bf Step~8. }The image of $F$ is 
$\{-1,0,1\}$. {\bf Step~9. }The maximal interval of convexity, respectively concavity, of $F$ is the interval 
$(-\infty,0]$, respectively $[0,+\infty)$. {\bf Step~10. }The function $F$ 
has no strict inflection point but it has inflection at every point $(x,F(x))$ with 
$x\ne0$. At the origin $(0,0)$ the graph $G_F$ does not have tangent.  
{\bf Step~11. }$F$ has no vertical asymptotes. The 
axis $x$ is the asymptote of $F$ at $-\infty$ and at $+\infty$. {\bf Step~12. }A~sketch of the graph $G_F$ of the function $F(x)=\sgn\,x$ is 

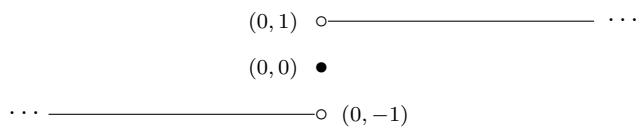
\medskip
\begin{picture}(50,50)(-40,0)
\put(0,5){$\ds$}
\put(15,5){\line(1,0){100}}\put(115,2.5){$\circ$}
\put(115,20){$\bullet$}
\put(120,40){\line(1,0){100}}\put(115,37.5){$\circ$}
\put(225,40){$\ds$}
\put(90,37.5){{\footnotesize $(0,1)$}}
\put(90,20){{\footnotesize $(0,0)$}}
\put(125,2.5){{\footnotesize $(0,-1)$}}
\end{picture}

\medskip\noindent
{\em $\bullet$ Second example. }Let 
$${\textstyle
F=F(x)\equiv\tan x=\frac{\sin x}{\cos x}=\frac{\sin x}{\sin(x+\pi/2)}\,.
}
$$
{\bf Step 0. }The last ratio shows that $F\in\mathrm{SEF}$. 
{\bf Step 1. }
\begin{eqnarray*}
M(F)&=&M(\sin x)\cap M(\cos x)\setminus Z(\cos x)=\R\setminus Z(\cos x)\\
&=&{\textstyle
\R\setminus\big\{n\pi+
\frac{\pi}{2}\cc\;n\in\Z\big\}=
\bigcup_{n\in\Z}\big(\pi n-
\frac{\pi}{2},\,\pi n+\frac{\pi}{2}\big)\,. 
}
\end{eqnarray*}
{\bf Step 2. }The function $F$ is $\pi$-periodic because
$$
\sin(\pi+x)=-\sin x\,\text{ and }\,\cos(\pi+x)=-\cos x\,. 
$$
It is an odd function because sine is odd and cosine is even. {\bf Step 3. }By 
Theorems~\ref{thm_EFjsouSpoj} and \ref{thm_deriSEF}, the 
function $F$ is continuous and $D(F)=M(F)$. By 
Exercise~\ref{ex_tanCot1}, 
$${\textstyle
F'(x)=\frac{1}{\cos^2 x}\,.
}
$$
{\bf Step 4. }For $n\in\Z$ let $b_n\equiv\pi n+\frac{\pi}{2}$. Then 
$$
\lim_{x\to 
b_n^-}F(x)=+\infty\,\text{ and }\,\lim_{x\to b_n^+}F(x)=-\infty\,. 
$$
The limits of $F(x)$ at $\pm\infty$ do not
exist. {\bf Step 5. }The graph $G_F$ intersects the axis
$y$ at the point $(0,0)$, and the axis $x$ at the points 
$${\textstyle
(b_n-\frac{\pi}{2},\,0)=
(\pi n,\,0),\ n\in\Z\,, 
}
$$
{\bf Step 6. }$D(F)=M(F)$ and 
there is nothing to compute. {\bf Step 7. }Since 
$${\textstyle
\text{$F'(x)=\frac{1}{\cos^2 x}>0$ on $M(F)$}\,, 
}
$$
the function $F$ increases on every interval 
$${\textstyle
(\pi n-\frac{\pi}{2},\,\pi n+\frac{\pi}{2}),\ n\in\Z\,. 
}
$$
These are the maximal intervals of monotonicity. $F$ has no extremes. {\bf Step~8. }Theorem~\ref{thm_mezihodnoty} and the infinite
limits in step~4 show that the image of $F$ is $\R$. {\bf Step 9. }By Corollaries~\ref{cor_deriRat} and \ref{cor_derExpKosSin}, 
and Theorem~\ref{thm_deriSEF} we have 
$${\textstyle
F''(x)=\frac{2\sin x}{\cos^3 x}\,
\text{ and }\,M(F'')=M(F')=M(F)\,. 
}
$$
Since 
$${\textstyle
\text{$F''<0$ on $(\pi n-\frac{\pi}{2},\pi n)$ and $F''>0$ on 
$(\pi n,\pi n+\frac{\pi}{2})$}\,, 
}
$$
the function $F$ is strictly concave, respectively convex, 
$$
{\textstyle
\text{on 
$(\pi n-\frac{\pi}{2},\pi n]$, respectively $[\pi n,\pi n+
\frac{\pi}{2})$, for $n\in\Z$\,. 
}}
$$
These are the maximal intervals of convexity and concavity. {\bf Step~10. }Due to the sign of $F''$ in the previous step, the 
inflection points of $F$ are exactly 
$(b_n-\frac{\pi}{2},0)=(\pi n,0)$, $n\in\Z$, and are strict. 
{\bf Step 11. }The limits in step~4 show that 
every line 
$${\textstyle
x=b_n=\pi n+\frac{\pi}{2},\ n\in\Z\,, 
}
$$
is both right and left vertical asymptote of $F$.  
There is no asymptote at 
$\pm\infty$. {\bf Step 12. } Use, for example, the drawing calculator at\\
\url{https://www.desmos.com/calculator}.

\medskip\noindent
{\em $\bullet$ Third example. }We roughly
follow the lecture notes \cite[pp. 193--194]{cern_poko}. Let $${\textstyle
F=F(x)\equiv\arcsin
\big(\frac{2x}{1+x^2}\big)\,.
}
$$
 {\bf Step 0. 
}The last expression shows that
$F\in\mathrm{EF}$. We will see that $D(F)\ne M(F)$. Therefore, by Theorem~\ref{thm_deriSEF},  $F\not\in\mathrm{SEF}$. Hence 
$F\in\mathrm{EF}\setminus
\mathrm{SEF}$. {\bf Step~1. }$M(F)=\R$ because $M(\arcsin 
x)=[-1,1]$ and $2|x|\le 1+x^2$ for every $x\in\R$ because 
$$
x^2\pm2x+1=(x\pm1)^2\ge0\,. 
$$
{\bf Step 2. }The function $F$ is odd because the functions $\sin x$, $\arcsin x$, and $\frac{2x}{1+x^2}$ are odd. It is not
periodic. {\bf Step 3. }By Theorem~\ref{thm_EFjsouSpoj} the function $F$ is 
continuous. Using part~1 of Exercise~\ref{ex_InvTri}, 
part~2 of Theorem~\ref{thm_DerSlozFce} (not
Corollary~\ref{cor_derSlozF}),
Corollary~\ref{cor_deriRat}, 
Exercise~\ref{ex_derKon}, 
Corollaries~\ref{cor_deriSumStan}
and \ref{cor_Leibform}, and part~6
of Exercise~\ref{ex_derirealMoc} we get that
$${\textstyle
D(F)=\big\{x\in\R\cc\;\frac{2x}{1+x^2}\ne\pm1\big\}=\R\setminus
\{-1,1\}=M(F)\setminus\{-1,\,1\}
}
$$ 
(in step~6 we show that the derivatives $F'(-1)$ and $F'(1)$ do not exist) and that
\begin{eqnarray*}
F'(x)&=&{\textstyle
\frac{1}{\sqrt{1-(2x/(1+x^2))^2}}\cdot\frac{2\cdot(1+x^2)-2x\cdot2x}{(1+x^2)^2}
=2\cdot\frac{(1-x^2)/(1+x^2)^2}{|(1-x^2)/(1+x^2)|}}\\
&=&{\textstyle
2\cdot\frac{1-x^2}{|1-x^2|}\cdot\frac{1}{1+x^2}=
\frac{2\cdot\sgn(1-x^2)}{1+x^2}\,|\,D(F)\,.
}
\end{eqnarray*}
{\bf Step 4. }Clearly, 
$$
\lim_{x\to-\infty}F(x)=\lim_{x\to+\infty}F(x)=
\arcsin 0=0
$$ 
because $\frac{2x}{1+x^2}\to0$ for
$x\to\pm\infty$. {\bf Step 5. }The graph $G_F$ intersects both axes 
only at the origin $(0,0)$. {\bf Step 6. }It is clear that
$\lim_{x\to1^{\pm}}F'(x)=\mp1$. So Proposition~\ref{prop_rozsDeri} gives that
$F'_{\pm}(1)=\mp1$. Since $F(x)$ is odd, $F'_{\pm}(-1)=\pm1$. By part~3 
of Exercise~\ref{ex_ojednDeri}, the derivatives $F'(-1)$ and $F'(1)$ do not exist. {\bf Step~7. }Since 
$$
\text{$F'<0$ on
$(-\infty,-1)$, $F'>0$ on $(-1,1)$ and $F'<0$ on
$(1,+\infty)$}\,,
$$
Theorem~\ref{thm_derMono1} implies that
$$
\text{$F$ decreases on $(-\infty,-1]$, increases on $[-1,1]$ 
and decreases on $[1,+\infty)$}\,. 
$$
These are the maximal 
intervals of monotonicity. Also $F(x)<0$ for $x<0$, $F(x)>0$ for 
$x>0$, and $F(0)=0$. Considering these domains of monotonicity, 
the zero limits in step~4 and the fact that $F$ is odd, we see that
$F(-1)=-\frac{\pi}{2}$ is the strict global minimum, that 
$F(1)=\frac{\pi}{2}$ is the strict global maximum, and that $F$ has no 
other (local or global) extreme. 
{\bf Step~8. }Using Theorem~\ref{thm_mezihodnoty} 
we get the 
image $F[M(F)]=[-\frac{\pi}{2},\frac{\pi}{2}]$. {\bf 
Step~9. }Using Exercise~\ref{ex_derKon}, part~1 (or 
part~2) of Theorem~\ref{thm_deriPodilu}, Corollary~\ref{cor_deriSumStan} and
part~6
of Exercise~\ref{ex_derirealMoc} we
get that
$${\textstyle
F''(x)=\frac{-4x\cdot\sgn(1-x^2)}{(1+x^2)^2}\,|\,D(F)\,.
}
$$ 
Since $F''<0$ on $(-\infty,-1)$, $F''>0$ on $(-1,0)$, $F''<0$ on $(0, 
1)$ and $F''>0$ on $(1,+\infty)$, 
Theorem~\ref{thm_druDerKonvKonk} gives that $F$ is strictly concave 
on $(-\infty,-1]$, strictly convex on $[-1,0]$, strictly concave on 
$[0,1]$ and strictly convex on $[1,+\infty)$. These are the maximal intervals of convexity 
and concavity. {\bf Step 10. }Because of the sign of $F''$ and since 
at the points $\pm1$ the graph $G_F$ does not have tangents, by 
Theorems~\ref{thm_infl} and
\ref{thm_postacInflexe} the point $(0,0)$ is the only inflection point 
of $G_F$. {\bf Step 11. }By the limits in step~4, the $x$-axis is an asymptote of $F$ at 
$-\infty$ and at $+\infty$. The function $F$ has no vertical asymptote. {\bf Step 12. }Use, for example, the drawing calculator at \url{https://www.desmos.com/calculator}. 

\begin{exer}\label{ex_GraDraw1}
Draw in steps 0--12 the graph of the  Riemann function 
$$
F=F(x)\equiv r(x)\,. 
$$
Recall that 
$r(\al)=0$ for $\al\in\R\setminus\Q$ and $r(\frac{p}{q})=\frac{1}{q}$ if the 
fraction $\frac{p}{q}$ is in lowest terms.    
\end{exer}

\begin{exer}\label{ex_GraDraw1apul}
Find for $r(x)$ maximal domains of 
monotonicity.   
\end{exer}
 
\begin{exer}\label{ex_GraDraw2}
Draw in steps 0--12 the graph of the function
$$
F=F(x)\equiv x^x\ \  (=\mathrm{e}^{x\log x})\,.
$$
\end{exer}

\chapter[Taylor polynomials]{Taylor polynomials}\label{chap_pr9}

This chapter is based on lecture~9 
\begin{quote}
\url{https://kam.mff.cuni.cz/~klazar/MAI24_pred9.pdf}    
\end{quote}
given on April 18, 2024. In Section~\ref{sec_TaylPoly}
we define for a~function $f(x)$ its Taylor polynomial $T^{f,b}_n(x)$ with 
order $n$ and center $b$. In our approach, this polynomial has degree at 
most $n$ and is determined by the condition that  
$$
f(x)=T^{f,b}_n(x)+o((x-b)^n)\ \ (x\to b)\,.
$$
By Proposition~\ref{prop_uniqTaylor} it is unique and by Theorem~\ref{thm_classTayl} if the definition domain of $f(x)$ and its derivatives is an interval $[b,c)$, then the 
coefficients of $T^{f,b}_n(x)$ are given by the standard formula. In 
Propositions~\ref{prop_withDeri1} and \ref{prop_withDeri2} and Theorem~\ref{thm_nonClasTay} 
we demonstrate that in general the coefficient of the quadratic term in $T^{f,b}_n(x)$ is independent of 
the second derivative $f''(x)$.

Section~\ref{sec_TaylorEF} is devoted to standard Taylor polynomials of functions defined on intervals.
In Propositions~\ref{prop_expSinCos}, \ref{prop_onePlusx}, \ref{prop_logs} and 
\ref{prop_arcs} we determine Taylor polynomials of functions $\mathrm{e}^x$, 
$\cos x$, $\sin x$, $(1+x)^a$, $\log(1+x)$, $\log(\frac{1}{1-x})$, 
$\arctan x$, $\arcsin x$ and $\arccos x$. We show how to compute by Taylor 
polynomials limits of the type
$$
\lim_{x\to0}\frac{f(x)}{g(x)}
$$
and in Theorem~\ref{thm_compTaylPol} and Proposition~\ref{prop_compMod} we 
describe arithmetic of Taylor polynomials. Theorem~\ref{thm_flatGraph} 
gives an example of a~nonzero and $\mathcal{C}^{\infty}(\R)$ function 
whose all derivatives at $0$ are $0$.

In Section~\ref{sec_TaylorSeries} in Theorems~\ref{thm_simpleTay} and \ref{thm_advanTay}, the simple and the advanced Taylor theorem, we give 
formulas for the Taylor remainder, which is the difference of the function and 
its Taylor polynomial. We define, in our approach, the Taylor series of 
a~function. In Theorem~\ref{thm_TaylSerExp} we derive
Maclaurin series of the functions $\mathrm{e}^x$, $\cos x$ and $\sin x$. 
Theorem~\ref{thm_NewtonBin} completely determines the domains where 
Newton's\index{Newton, Isaac} binomial series sums to the real power 
$(1+x)^a$.  

The extending Section~\ref{sec_analFunc} is devoted to functions representable as 
sums of real power series. ...

\section[Taylor polynomials]{Taylor polynomials}\label{sec_TaylPoly}

We consider local polynomial approximations of functions.

\medskip\noindent
{\em $\bullet$ Linear and constant approximations. }Let $M\sus\R$. Recall 
Definition~\ref{def_maleo} of the asymptotic notation $o(\cdot)$ and recall that if a~function $f(x)\in\mathcal{F}(M)$ is differentiable 
at a~point $b$ ($\in M\cap L(M)$), which means that the derivative $f'(b)\in\R$ exists, then $f(x)$ is approximated near $b$ as
$$
f(x)=f(b)+f'(b)(x-b)+o(x-b)\ \ (x\to b)\,.
$$
The linear function 
$$
\ell(x)\equiv f(b)+f'(b)(x-b)=f'(b)x+f(b)-f'(b)b
$$
is the tangent (line) to the graph $G_f$
at the point $(b,f(b))$.
This approximation is in fact equivalent to the differentiability.

\begin{exer}\label{ex_proveEqui}
$f(x)=f(b)+c(x-b)+o(x-b)$ ($x\to b$) $\iff$ $f'(b)=c$.  
\end{exer}

In an even simpler situation a~function $f(x)\in\mathcal{F}(M)$ continuous at $b$ is approximated 
near the point $b$ ($\in L(M)\cap M$) as
$$
f(x)=f(b)+o(1)\ \ (x\to b)
$$
by the constant function $k_{f(b)}(x)$.

\begin{exer}\label{ex_proveEqui22}
$f(x)=f(b)+o(1)$ ($x\to b$) $\iff$ $f$ is continuous at $b$.    
\end{exer}

\noindent
{\em $\bullet$ Approximation definition of Taylor polynomials. }We strengthen the above linear and constant
approximations to polynomial approximations. By 
Definition~\ref{def_polynomy}, 
polynomials are the 
$\mathcal{F}(\R)$ functions that
arise from constant functions $k_c(x)$, $c\in\R$, and the identity function 
$\mathrm{id}(x)$ by repeated addition and multiplication. 
It is clear that for every $n\in\N_0$ and real numbers $a_0$, $a_1$, 
$\ds$, $a_n$ the function
$$
{\textstyle
\sum_{j=0}^n a_j(x-b)^j\ \ \big(=\sum_{j=0}^n
k_{a_j}(x)\cdot\prod_{i=1}^j(\mathrm{id}(x)+k_{-1}(x)\cdot k_b(x))\big)
}
$$
is a~polynomial. The next definition is inspired and motivated by so called 
Peano\index{Peano, Giuseppe} derivatives,\index{Peano derivative} see survey articles \cite{evan_weil,weil} on them.

\begin{defi}[Taylor polynomials]\label{def_TaylPoly}
Let $n\in\N_0$, $M\sus\R$, $b\in M\cap L(M)$ and let $f(x)\in\mathcal{F}(M)$. If $a_0$, $a_1$, $\ds$, $a_n$ 
are $n+1$ real numbers such that
$$
{\textstyle
f(x)=\sum_{j=0}^n a_j(x-b)^j+o((x-b)^n)\ \ (x\to b)\,,
}
$$
we say that 
$${\textstyle
\text{$\sum_{j=0}^n a_j(x-b)^j$
is a~\underline{Taylor polynomial\index{Taylor 
polynomials, $T_n^{f,b}(x)$|emph}}}} 
$$
of the function $f(x)$ with order $n$ and center $b$.
\end{defi}
In the definition we assume that $b\in M(f)$ but we do not need the value $f(b)$.

\begin{prop}[uniqueness]\label{prop_uniqTaylor}
Taylor polynomials are unique.
\end{prop}
\duk
Let, for the contrary, $n\in\N_0$, $f\in\mathcal{F}(M)$, $b\in M\cap L(M)$, 
and $p(x)=\sum_{j=0}^n a_j(x-b)^j$ and $q(x)=\sum_{j=0}^n b_j(x-b)^j$ be 
two distinct polynomials such that
$$
f(x)=p(x)+o((x-b)^n)\,\text{ and }\,
f(x)=q(x)+o((x-b)^n)\ \ (x\to b)\,.
$$
Equivalently,
$$
p(x)=f(x)+o((x-b)^n)\,\text{ and }\,
q(x)=f(x)+o((x-b)^n)\ \ (x\to b)\,.
$$
Subtracting we get
$$
{\textstyle
\sum_{j=m}^n c_j(x-b)^j=o((x-b)^n)\ \ (x\to b)\,,
}
$$
for some $m\in\N_0$ with $m\le n$ and $c_m\ne0$. This is an impossible asymptotic equality.
\kduk

\noindent
The unique Taylor polynomial of a~function $f$ with order $n$ and center $b$, if it exists, is denoted by 
$$
T_n^{f,\,b}(x)\ \ (\in\mathrm{POL})\,.\label{TaylorPol}
$$

\begin{exer}\label{ex_iniSumTP}
Let $T_n^{f,b}(x)=\sum_{j=0}^n a_j(x-b)^j$. If $m\in\N_0$ with
$m\le n$, then $T_m^{f,b}(x)=\sum_{j=0}^m a_j(x-b)^j$. 
\end{exer}

Changing the variable we can move the center to $0$.

\begin{prop}[moving to $0$]\label{prop_confTo0}
Let $n\in\N_0$, $M\sus\R$, $b\in M\cap L(M)$ and $f\in\mathcal{F}(M)$. Suppose that the Taylor polynomial 
$${\textstyle
T_n^{f,b}(x)=\sum_{j=0}^na_j(x-b)^j
}
$$ 
exists. \underline{Then} the composite function 
$$
g(x)\equiv f(x+b)=f\circ (\mathrm{id}+k_b)
$$ 
has the Taylor polynomial
$${\textstyle
T_n^{g,\,0}(x)=\sum_{j=0}^na_jx^j\,.
}
$$
\end{prop}
\duk
Then $0\in M(g)\cap L(M(g))$, because $M(g)=M-b$, and using Theorem~\ref{thm_LimSlozFunkce}
we get
$$
\lim_{x\to0}{\textstyle
\frac{g(x)-\sum_{j=0}^na_jx^j}{x^n}
}=
\lim_{x\to0}{\textstyle
\frac{f(x+b)-\sum_{j=0}^na_j(x+b-b)^j}{(x+b-b)^n}
}=
\lim_{y\to b}{\textstyle
\frac{f(y)-\sum_{j=0}^na_j(y-b)^j}{(y-b)^n}
}=0\,.
$$
\kduk

\noindent
{\em $\bullet$ Classical Taylor polynomials. }In the next theorem we
obtain the classical formula for coefficients of Taylor polynomials of 
$f(x)$ in terms of values of 
derivatives $f^{(j)}(x)$. In the proof we use an auxiliary proposition. 

\begin{prop}[Taylor polynomials of $f$ and $f'$]\label{prop_TaylFandDerF}
Let $n\in\N$, $b<c$ be in $\R$ and let $f,f'\in\mathcal{F}([b,c))$.
$$\text{If }\,
T_{n-1}^{f',\,b}(x)=
\sum_{j=0}^{n-1}a_j(x-b)^j\,
\text{ \underline{then} }\,
T_n^{f,\,b}(x)=f(b)+
\sum_{j=1}^n\frac{a_{j-1}}{j}(x-b)^{j}\,.
$$
\end{prop}
\duk
Let $n$, $b$, $c$ and $f$ be as stated. Suppose that the assumption of the implication 
holds. We denote the last displayed polynomial by $p(x)$. Then $p'(x)=T_{n-1}^{f',b}(x)$. 
By the first case of Theorem~\ref{thm_LHP2} and the assumption,
$$
{\textstyle
\lim_{x\to b}\frac{f(x)-p(x)}{(x-b)^n}=
\lim_{x\to b}\frac{(f(x)-p(x))'}{((x-b)^n)'}=
n^{-1}\lim_{x\to b}\frac{f'(x)-T_{n-1}^{f',\,b}(x)}{(x-b)^{n-1}}=0\,.
}
$$
Hence $p(x)=T_n^{f,b}(x)$ by Proposition~\ref{prop_uniqTaylor}.
\kduk
\vspace{-3mm}
\begin{exer}\label{ex_whyfb}
Why is in $p(x)$ the term $f(b)$?
Does not the computation via {\em HR~2} work also 
for the simpler polynomial $q(x)\equiv\sum_{j=1}^n\frac{1}{j}a_{j-1}(x-b)^j$? Then $q'(x)=T_{n-1}^{f',b}(x)$ too.
\end{exer}

For $n=1$ we have 
$$
T^{f,b}_n(x)=f(b)+f'(b)(x-b) 
$$
whenever the Taylor polynomial exists. However, in the next passage we show that for $n\ge2$ the 
coefficient of $(x-b)^2$ in $T^{f,b}_n(x)$ has in general nothing to do with $f''(b)$. In the following theorem we describe situation when $T^{f,b}_n(x)$ 
has the classical coefficients $\frac{1}{j!}f^{(j)}(b)$, $j=0,1,\ds,n$. We assume that $n\ge2$ because the cases $n=0,1$ are resolved, for a~general definition domain $M(f)$, by Exercises~\ref{ex_proveEqui22} and 
\ref{ex_proveEqui}, respectively.

\begin{thm}[classical Taylor polynomials]\label{thm_classTayl}
Let\index{theorem!classical Taylor polynomials|emph} 
$n\in\N$ with $n\ge2$, $b<c$ be in $\R$
and let $f$, $f'$, $\ds$, $f^{(n)}$ be in
$\mathcal{F}([b,c))$. 
\underline{Then} $T_n^{f,b}(x)$ exists and
$$
T_n^{f,\,b}(x)=\sum_{j=0}^n\frac{1}{j!}f^{(j)}(b)\cdot(x-b)^j\,.
$$
\end{thm}
\duk
Let $n$, $b$, $c$ and $f$ be as stated. We proceed by induction on $n\ge2$. Let $n=2$. By Exercise~\ref{ex_proveEqui},
$$
T^{f',\,b}_1(x)=f'(b)+f''(b)(x-b)\,.
$$
Hence by Proposition~\ref{prop_TaylFandDerF},
\begin{eqnarray*}
T^{f,\,b}_2(x)&=&{\textstyle
f(b)+\frac{f'(b)}{1}(x-b)+\frac{f''(b)}{2}(x-b)^2}\\
&=&{\textstyle
\frac{f^{(0)}(b)}{0!}(x-b)^0+\frac{f'(b)}{1!}(x-b)+\frac{f''(b)}{2!}(x-b)^2\,.
}
\end{eqnarray*}
Let $n>2$. By the inductive assumption,
$$
{\textstyle
T^{f',\,b}_{n-1}(x)=\sum_{j=0}^{n-1}\frac{1}{j!}f^{(j+1)}(b)(x-b)^j\,.
}
$$
Hence by Proposition~\ref{prop_TaylFandDerF},
$$
{\textstyle
T^{f,\,b}_n(x)=f(b)+\sum_{j=1}^n\frac{1}{j}\cdot\frac{1}{(j-1)!}f^{(j)}(b)(x-b)^j=
\sum_{j=0}^n\frac{1}{j!}f^{(j)}(b)(x-b)^j\,.
}
$$
\kduk
\vspace{-3mm}
\begin{exer}\label{ex_halfNeighb}
Show that the previous theorem and proposition hold for definition domains $(c,b]$ with $c<b$ and $U(b,\de)$.   
\end{exer}

\begin{exer}\label{ex_TaylExer}
Suppose that $p(x)=\sum_{j=0}^n a_jx^j$ is a~polynomial and $m\in\N_0$. What is $T_m^{p,0}(x)$?
\end{exer}

\medskip\noindent
{\em $\bullet$ Non-classical Taylor polynomials. }Definition~\ref{def_TaylPoly} is non-standard, Taylor polynomials are usually defined by the formula 
in Theorem~\ref{thm_classTayl}. We illustrate by several examples differences between our and standard definition.

\begin{prop}[no derivatives]\label{prop_withDeri1}
Let $n\ge2$ and $b$, $c$ and $f$ be as in 
Theorem~\ref{thm_classTayl}. Let $M\sus [b,c)$ be any set such that $b\in 
M\cap L(M)$, and let
$$
g\equiv f\,|\,M\,. 
$$
\underline{Then}, no matter if for $j\ge2$ the derivatives $g^{(j)}(b)$ exist or not, we have
$$
T_n^{g,\,b}(x)=T_n^{f,\,b}(x)
=\sum_{j=0}^n\frac{1}{j!}f^{(j)}(b)(x-b)^j\,.
$$
\end{prop}
\duk
By Theorem~\ref{thm_classTayl}, 
Proposition~\ref{prop_uniqTaylor} and Definition~\ref{def_TaylPoly}.
\kduk

\noindent
Thus the Taylor polynomial $T_n^{g,b}(x)$ may exist even if none of the 
derivatives $g^{(j)}(b)$, $j\ge2$, exists. However, always 
$g^{(0)}(b)=g(b)=f(b)$ and
$g'(b)=f'(b)$.

\begin{exer}\label{ex_nonclaTp}
We define 
$${\textstyle
\text{$f\in\mathcal{F}(\{\frac{1}{n}\cc\;n\in\N\}\cup\{0\})$ by $f(0)\equiv0$ and
$f(\frac{1}{n})\equiv n^{-4}$}\,.
}
$$
Find $T_n^{f,0}(x)$ for every $n\in\N$.
What are the derivatives $f^{(j)}(0)$
for $j\in\N_0$?
\end{exer}

Definition domains of functions $g$ without derivatives $g^{(j)}(b)$, 
$j\ge2$, in Proposition~\ref{prop_withDeri1} are 
typically sparse. However, it is not hard to find a~function $f$ with the 
following properties. 
\begin{enumerate}
\item Functions $f$ and $f'$ are defined on $\R$.
\item $T_2^{f,0}(x)=0+0x+0x^2$.
\item The derivative $(f')'(0)$ does not exist.
\end{enumerate}

\begin{prop}[an example of such $f$]\label{prop_withDeri2}
The function $f$, defined by
$$
\text{$f(0)\equiv0$ and $f(x)\equiv x^3\sin(x^{-3})$ for $x\ne0$}\,,
$$
has properties {\em 1--3}. 
\end{prop}
\duk
Since $f'(x)=3x^2\sin(x^{-3})-3x^{-1}\cos(x^{-3})$ for $x\ne 0$, and $f'(0)=0$ by the definition of derivative, $f$ has property 1. The limit
$$
\lim_{x\to0}x^{-2}(f(x)-0)=0
$$
shows that $f$ has property 2. Since $f'$ is unbounded on any neighborhood of zero, $f$ has property 3.
\kduk

Recall that for every real interval $I$ we set 
$$
I_{\Q}\equiv I\cap\Q\,. 
$$
We modify the construction in Theorem~\ref{prop_LHPcouExa} and get
the next theorem. It shows that
the value $f''(b)$ is unrelated to the coefficient of $(x-b)^2$ in our Taylor 
polynomials.

\begin{thm}[independence of $a_2$ and $f''(0)$]\label{thm_nonClasTay}
For\index{theorem!independence of $a_2$ and $f''(0)$|emph} 
any $c\in\R$ there is a~function $f$  such that
$${\textstyle
\text{$f,\,f',\,f''\in\mathcal{F}([0,\,1]_{\Q})$, $T_2^{f,\,0}(x)=\sum_{j=0}^2 0x^j$ and $f''(0)=c$}\,.
}
$$
\end{thm}
\duk
Let 
$$
(c_n)\sus(0,\,1)
$$ 
be a~sequence of irrational numbers such that
$${\textstyle
\text{$c_0\equiv1>c_1>c_2>\ds>0$, $\lim c_n=0$ and $\lim\frac{c_{n-
1}}{c_n}=1$}\,. 
}
$$
We set $f(0)\equiv0$ and
for $\al\in(c_n,c_{n-1}]_{\Q}$, $n\in\N$, we define
$${\textstyle
f(\al)\equiv c_n^3+\frac{c}{2}\cdot\al^2-\frac{c}{2}\cdot c_n^2\,. 
}
$$
Thus the graph $G_f$
consists of short pierced parabolic segments starting on
the cubic $y=x^3$. Clearly, $f'(\al)=c\al$ on $(0,1]_{\Q}$ and it is easy to see that 
also $f'(0)=0$. Hence $f''(\al)=c$ on $[0,1]_{\Q}$.
For $\al\in(c_n,c_{n-1})_{\Q}$ we have
$$
{\textstyle
0<\frac{f(\al)}{\al^2}\le
\frac{c_n^3+(c/2)c_{n-1}^2-(c/2)c_n^2}{c_n^2}=
c_n+\frac{c}{2}\big(\frac{c_{n-1}}{c_n}\big)^2-\frac{c}{2}\to0\ \ (n\to\infty)\,.
}
$$
Therefore $f(x)=o(x^2)$ ($x\to0$) and 
$T_2^{f,0}(x)=0+0x+0x^2$.
\kduk

\begin{exer}\label{ex_exteCons}
Extend the previous construction to the definition domain $[-1,1]_{\Q}$ and the 
center $b=0$.    
\end{exer}

\section[Examples of Taylor polynomials]{Examples of Taylor polynomials}\label{sec_TaylorEF}

We derive 
several classical formulas for Taylor polynomials. For simplicity of
notation, we set the center to $b=0$.

\medskip\noindent
{\em $\bullet$ Taylor polynomials of the exponential, cosine, and sine. }Taylor 
polynomials are partial sums of power series defining these functions. 

\begin{prop}[$\mathrm{e}^x$, $\cos x$ and $\sin x$]\label{prop_expSinCos}
We have the following Taylor polynomials with order $n\in\N_0$ and center $0$.
\begin{enumerate}
\item If $f(x)\equiv\mathrm{e}^x$ \underline{then} 
$$
T_n^{f,\,0}(x)=\sum_{j=0}^n
\frac{1}{j!}\cdot x^j\,.
$$
\item If $f(x)\equiv\cos x$ \underline{then} 
$$
T_{2n}^{f,\,0}(x)=\sum_{j=0}^n
(-1)^j\frac{1}{(2j)!}\cdot x^{2j}\,.
$$
\item If $f(x)\equiv\sin x$ \underline{then} 
$$
T_{2n+1}^{f,\,0}(x)=\sum_{j=0}^n(-1)^j
\frac{1}{(2j+1)!}\cdot x^{2j+1}\,.
$$
\end{enumerate}
\end{prop}
\duk
These formulas follow from Theorem~\ref{thm_classTayl}, Corollary~\ref{cor_derExpKosSin} and Exercise~\ref{ex_naDerKrad}. 
\kduk

\noindent
Returning to the definitions in Section~\ref{sec_elemenFce}
we see that for each of the three functions $f(x)\in\{\mathrm{e}^x,\cos x,\sin x\}$ and every $a\in\R$ we have, remarkably,
$$
f(a)=\lim_{n\to\infty}T_n^{f,\,0}(a)\,.
$$ 

\begin{exer}\label{ex_TedToNepl}
Show that the function $f(x)\equiv\frac{1}{1-x}$ does not have this property.     
\end{exer}

\noindent
{\em $\bullet$ More examples of Taylor polynomials. }For any $a\in\R$ and $j\in\N_0$ we define 
the \underline{generalized binomial coefficient\index{generalized binomial 
coefficients|emph}} as $\binom{a}{0}\equiv1$, and for $j>0$ as
$$
\binom{a}{j}\equiv
\frac{a(a-1)\ds(a-j+1)}{j!}\,.
\label{anadje}
$$

\begin{prop}[$(1+x)^a$]\label{prop_onePlusx}
If $n\in\N_0$, $a\in\R$ and $f(x)\equiv(1+x)^a$, \underline{then}
$$
T_n^{f,\,0}(x)=
\sum_{j=0}^n\binom{a}{j}x^j\,.
$$
\end{prop}
\duk
This is immediate from Theorem~\ref{thm_classTayl} and derivatives ($j\in\N$)
$$
{\textstyle
\big((1+x)^a\big)^{(j)}=
a(a-1)\ds(a-j+1)\cdot(1+x)^{a-j}}\,,
$$
which follow from parts 2--4 of Exercise~\ref{ex_derirealMoc}.
\kduk

The following logarithmic expansions are often used. 

\begin{prop}[$\log(1+x)$ and $\log(\frac{1}{1-x})$]\label{prop_logs}
Suppose that $n\in\N_0$. If $f(x)\equiv\log(1+x)$ and  $g(x)\equiv\log(\frac{1}{1-x})$
\underline{then}
$$
T_n^{f,\,0}(x)=\sum_{j=1}^n
(-1)^{j-1}j^{-1}\cdot x^j\,\text{ and }\,T_n^{g,\,0}(x)=\sum_{j=1}^n
j^{-1}\cdot x^j\,.
$$
\end{prop}
\duk
The first formula follows from Theorem~\ref{thm_classTayl} and derivatives ($j\in\N$)
$$
\big(\log(1+x)\big)^{(j)}=
(-1)^{j-1}\cdot(j-1)!\cdot(1+x)^{-j}\,|\,(-1,\,+\infty)\,.
$$
The second formula follows from the first:
$\log(\frac{1}{1-x})=-\log(1+(-x))$.
\kduk

Taylor polynomials of inverse trigonometric functions are derived with the help of Proposition~\ref{prop_TaylFandDerF}.

\begin{prop}[$\arctan x$, $\arcsin x$ and $\arccos x$]\label{prop_arcs}
Let $n\in\N_0$.
\begin{enumerate}
\item If $f(x)\equiv\arctan x$ \underline{then} 
$$
T_{2n+1}^{f,\,0}(x)=\sum_{j=0}^n(-1)^j\frac{1}{2j+1} 
\cdot x^{2j+1}\,.
$$
\item If $f(x)\equiv\arcsin x$ \underline{then} 
$$
T_{2n+1}^{f,\,0}(x)=\sum_{j=0}^n\frac{1}{2j+1}(-1)^j\binom{-\frac{1}{2}}{j}\cdot x^{2j+1}
=\sum_{j=0}^n\frac{1}{2j+1}\binom{j-\frac{1}{2}}{j}\cdot x^{2j+1}\,.
$$
\item If $f(x)\equiv\arccos x$ \underline{then} 
$$
T_{2n+1}^{f,\,0}(x)=\frac{\pi}{2}-
\sum_{j=0}^n\frac{1}{2j+1}\binom{j-\frac{1}{2}}{j} \cdot x^{2j+1}\,.
$$
\end{enumerate}
\end{prop}
\duk
1. For $|x|<1$ we have by part~3 of Exercise~\ref{ex_InvTri} and 
Theorem~\ref{thm_geomRada} that
$${\textstyle
(\arctan x)'=\frac{1}{1+x^2}=\sum_{j=0}^{\infty}(-1)^jx^{2j}\,. 
}
$$
Exercise~\ref{ex_geomTail} gives that for every $n\in\N_0$,
$$
{\textstyle
\sum_{j=0}^{\infty}(-1)^jx^{2j}=
\sum_{j=0}^n(-1)^jx^{2j}+o(x^{2n})\ \ (x\to0)\,.
}
$$
Thus by Proposition~\ref{prop_uniqTaylor} we have for every $n\in\N_0$ that 
$${\textstyle
T_{2n}^{f',\,0}(x)
=\sum_{j=0}^n(-1)^jx^{2j}
}
$$ 
and the result follows by Proposition~\ref{prop_TaylFandDerF}. 

2. By part~1 of Exercise~\ref{ex_InvTri} we have
$${\textstyle
(\arcsin x)'=\frac{1}{\sqrt{1-x^2}}=(1-x^2)^{-1/2}\,. 
}
$$
Thus Proposition~\ref{prop_onePlusx} gives for every $n\in\N_0$ that
$$
{\textstyle
T_{2n}^{f',\,0}(x)=
\sum_{j=0}^{2n}(-1)^j\binom{-1/2}{j}x^{2j}\,.
}
$$
We are done by Proposition~\ref{prop_TaylFandDerF} and Exercise~\ref{ex_idenBinKoe}.

3. The computation is similar to the previous one. 
\kduk
\vspace{-3mm}
\begin{exer}\label{ex_geomTail}
Prove that for every $n\in\N_0$, 
$${\textstyle
\sum_{j=0}^{\infty}x^j=\sum_{j=0}^n x^j+o(x_n)\ \  (x\to0)\,.
}
$$
\end{exer}

\begin{exer}\label{ex_idenBinKoe}
Prove that for every $j\in\N_0$ and $a\in\R$,
$$
{\textstyle
(-1)^j\binom{-a}{j}=\binom{a+j-1}{j}\,.
}
$$
\end{exer}

\noindent
{\em $\bullet$ Limits
$\lim_{x\to b}\frac{f(x)}{g(x)}$. }We show how to compute them with the help of Taylor
polynomials. By Proposition~\ref{prop_confTo0} we can
restrict ourselves to the center $b=0$. An important auxiliary result concerns the limit
$$
\lim_{x\to0}\frac{p(x)\,|\,M}{q(x)\,|\,M}
$$
of the ratio of restricted polynomials. Restrictions to sets $M$ 
results from our approach, we consider Taylor polynomials of functions
with arbitrary definition domains.

\begin{prop}[ratios of restricted polynomials]
Let $n\in\N_0$, $M\sus\R$ with $0\in L(M)$ and let
$$
L\equiv\lim_{x\to0}\frac{\sum_{j=0}^n a_jx_j\,|\,M+o(x^n)}{\sum_{j=0}^n b_jx_j\,|\,M+o(x^n)}\ \ (x\to0)\,,
$$
where we assume that not all $b_j$ are zero. Let $m\in\N_0$ with $m\le n$ be the minimum 
index such that $a_m\ne0$ or $b_m\ne0$. If $b_m=0$, let $l\in\N_0$ with $m<l\le 
n$ be the minimum index such that 
$b_l\ne0$. The following holds.
\begin{enumerate}
\item If $b_m\ne0$ \underline{then} $L$ exists and 
$$
L=\frac{a_m}{b_m}\,.
$$
\item If $b_m=0$, so that $a_m\ne0$,  and $0\not\in L^{\mp}(M)$,  \underline{then} $L$ exists and, with equal signs,  
$$
L=(\pm1)^{l-m}\cdot
(\sgn\,a_m)\cdot(\sgn\,b_l)\cdot(+\infty)\,. 
$$
\item If $b_m=0$ and $0\in 
L^{\mathrm{TS}}(M)$, \underline{then} $L$ exists if $l-m$ is even and
$$
L=(\sgn\,a_m)\cdot(\sgn\,b_l)\cdot(+\infty)\,.
$$
If $l-m$ is odd \underline{then} $L$ does not exist.
\end{enumerate}
\end{prop} 
\duk
1. We define empty sums as $0$ and have 
$$
L=
\lim_{x\to0}{\textstyle
\frac{a_m+\sum_{j=m+1}^n a_jx^{j-m}+o(x^{n-m})}{b_m+\sum_{j=m+1}^n b_jx^{j-m}+o(x^{n-m})}=\frac{a_m}{b_m}\,.
}
$$

2. Now 
$$
L=\lim_{x\to0}{\textstyle
\frac{a_m+o(1)}{b_lx^{l-m}+o(x^{l-m})}
}\,|\,M=
\lim_{x\to0^{\pm}}{\textstyle
\frac{a_m+o(1)}{b_lx^{l-m}+o(x^{l-m})}
}
$$
and this equals to the stated product of signs and $+\infty$. 

3. If $l-m$ is 
even, for $x\to0$ the factor $x^{l-m}$ goes to $0$ through positive 
values and $L$ is the stated product of signs and $+\infty$. If
$l-m$ is odd then the two limits
$\lim_{x\to0^{\pm}}\frac{\cdots}{\cdots}$ are two different infinities and $L$ does not exist.
\kduk

\noindent
We mention a~simple case of the proposition for deleted neighborhoods of $0$.

\begin{cor}[$M=P(0,\de)$]\label{cor_onU0de}
For $M=P(0,\de)$ the following holds in the previous proposition.
\begin{enumerate}
\item If $b_m\ne0$ \underline{then} 
$$
L=\frac{a_m}{b_m}\,.
$$
\item If $b_m=0$, so that $a_m\ne0$, and $l-m$ is even \underline{then}
$$
L=(\sgn\,a_m)\cdot(\sgn\,b_l)\cdot(+\infty)\,.
$$
\item If $b_m=0$ and $l-m$ is odd \underline{then} $L$ does not exist.
\end{enumerate}
\end{cor}
\duk
Cases~1 and~3 of the proposition occur and case~2 does not occur.
\kduk

\noindent
We illustrate each case of the corollary by an example. Using Taylor polynomials with order $3$, we get by case~1 that
$$\lim_{x\to0}{\textstyle
\frac{\sin(2x)-2\sin x}{\cos(2x)-\cos x}}=
\lim_{x\to0}{\textstyle\frac{(2x-8x^3/6)-2(x-x^3/6)+o(x^3)}{(1-4x^2/2)-(1-x^2/2)+o(x^3)}}=\lim_{x\to0}{\textstyle\frac{-x^3+o(x^3)}{-3x^2/2+o(x^3)}
}
$$
equals $\frac{0}{-3/2}=0$. Taylor polynomials with order $4$ show by 
case~3 that
\begin{eqnarray*}
&&\lim_{x\to0}{\textstyle
\frac{\sin(2x)-2\sin x}{\cos(2x)-\cos x+3x^2/2}}=  
\lim_{x\to0}{\textstyle
\frac{-x^3+o(x^4)}{(1-4x^2/2+16x^4/24)-(1-x^2/2+x^4/24)+3x^2/2+o(x^4)}}\\
&&=\lim_{x\to0}{\textstyle
\frac{-x^3+o(x^4)}{5x^4/8+o(x^4)}
}
\end{eqnarray*}
does not exist. Finally, Taylor polynomials with order $5$ show by case~2 that
$$
\lim_{x\to0}{\textstyle
\frac{\sin(2x)-2\sin x}{\arctan x-x+x^3/3}}=
\lim_{x\to0}{\textstyle
\frac{-x^3+x^5/4+o(x^5)}{x^5/5+o(x^5)}=-\infty\,.
}
$$

\begin{exer}\label{ex_LimByTayl1}
Compute
$$
\lim_{x\to0}{\textstyle
\frac{\sqrt{1+x}-1-x/2}{\log(1+x)-x}\,.
}
$$
\end{exer}

\begin{exer}\label{ex_LimByTayl2}
Compute
$$
\lim_{x\to0}{\textstyle
\frac{\sqrt{1+x}-\sqrt{1+3x}}{\sqrt[3]{1+x}-\sqrt[3]{1+2x}}\,.
}
$$    
\end{exer}

\begin{exer}\label{ex_LimByTayl3}
Compute
$$
\lim_{x\to0}{\textstyle
\frac{\arcsin x-x}{\sin x-x}\,.
}
$$    
\end{exer}

\section[Arithmetic of Taylor polynomials]{${}^*$Arithmetic of Taylor polynomials}\label{sec_ariTayl}

In view of the usefulness of Taylor polynomials we investigate their interactions with the six operations on the set of 
functions $\mathcal{R}$, namely with addition, multiplication, division, 
composition, inverse and (global) derivative. For details of these 
operations see Definitions~\ref{def_oprOnR} and 
\ref{def_deriFunk}. Like for continuity, the most complicated and interesting 
operation turns out to be inverse. We devote to it the last
Section~\ref{sec_LIF}. We begin with the arithmetic of the error terms $o((x-b)^n)$.

\medskip\noindent
{\em $\bullet$ Arithmetic of errors. }We leave these properties of $o((x-b)^n)$ ($x\to b$) as 
exercises. We assume that $n\in\N_0$, $b\in\R$ and $x\to b$.

\begin{exer}\label{ex_arErr1}
Always $o((x-b)^n)+o((x-b)^n)=o((x-b)^n)$.    
\end{exer}

\begin{exer}\label{ex_arErr2}
Let $m\in\N_0$ and $a\in\R$. If $m\ge n+1$ then 
$$
a(x-b)^m=o((x-b)^n)\,.
$$
\end{exer}

\begin{exer}\label{ex_arErr3}
Let $m\in\N_0$ and $a\in\R$. Then 
$$
a(x-b)^m\cdot o((x-b)^n)=o((x-b)^n)\,.
$$
\end{exer}

\begin{exer}\label{ex_arErr4}
For $m,n\in\Z$, $a\in\R$ and $x\to b$ we actually have 
$$
a(x-b)^m\cdot o((x-b)^n)=o((x-b)^{m+n})\,.
$$
\end{exer}

\noindent
{\em $\bullet$ Addition and multiplication }of Taylor polynomials is easy. 

\begin{prop}[addition, multiplication]\label{prop_addiTP}
Let $n\in\N_0$, $f,g\in\mathcal{F}(M)$ and $b\in M\cap L(M)$. Suppose that the Taylor polynomials 
$T^{f,b}_n(x)$  and $T^{g,\,b}_n(x)$
exist. We write
$$
{\textstyle
p(x)=T^{f,\,b}_n(x)=\sum_{j=0}^n a_j(x-b)^j\,\text{ and }\,q(x)=T^{g,\,b}_n(x)=\sum_{j=0}^n b_j(x-b)^j\,.
}
$$
\underline{Then} the following holds.
\begin{enumerate}
\item The Taylor polynomial $T^{f+g,b}_n(x)$ exists and
$$
T^{f+g,\,b}_n(x)=T^{f,\,b}_n(x)+T^{g,\,b}_n(x)\ \ (\in\mathrm{POL})\,.
$$
\item The Taylor polynomial $T^{fg,b}_n(x)$ exists and
$$
T^{fg,\,b}_n(x)=
\sum_{\substack{i,j\in\N_0\\i+j\le n}}a_ib_j(x-b)^{i+j}\equiv r(x)\ \ (\in\mathrm{POL})\,.
$$
\end{enumerate}
\end{prop}
\duk
Let $n$, $f$, $g$, $b$, $p(x)$, $a_j$, $q(x)$ and $b_j$ be as stated, and let $x\to b$. 

1. By the assumption and Exercise~\ref{ex_arErr1} we have 
$$
f(x)+g(x)=p(x)+o((x-b)^n)+q(x)+o((x-b)^n)=p(x)+q(x)+o((x-b)^n)\,.
$$
Hence, by Proposition~\ref{prop_uniqTaylor}, 
$$
T^{f+g,\,b}_n(x)=T^{f,\,b}_n(x)+T^{g,\,b}_n(x)\,.
$$

2. By Exercises~\ref{ex_arErr1} and \ref{ex_arErr2},  $p(x)q(x)=r(x)+o((x-b)^n)$. Hence by the assumption and Exercises~\ref{ex_arErr1} and \ref{ex_arErr3},
\begin{eqnarray*}
f(x)g(x)&=&(p(x)+o((x-b)^n))\cdot(q(x)+o((x-b)^n))\\
&=&
p(x)q(x)+o((x-b)^n)=r(x)+o((x-b)^n)\,.
\end{eqnarray*}
Proposition~\ref{prop_uniqTaylor} shows that $T^{fg,b}_n(x)=r(x)$.
\kduk

\noindent
For example, if $f(x)\equiv\cos x$, $g(x)\equiv\sqrt{1+x}$, $b=0$ and $n=4$, 
then
\begin{eqnarray*}
T^{f+g,\,b}_n(x)&=&{\textstyle(1+1)+(0+\frac{1}{2})x+
(-\frac{1}{2}-\frac{1}{8})x^2+(0+\frac{1}{16})x^3+(\frac{1}{24}-\frac{5}{128})x^4
}\\
&=&{\textstyle
2+\frac{1}{2}x-\frac{5}{8}x^2+\frac{1}{16}x^3+\frac{1}{384}x^4
}
\end{eqnarray*}

\begin{exer}\label{ex_souTay}
Determine the Taylor polynomial of $\arctan x+\sin x$ with order $3$ and 
center $0$.    
\end{exer}

\noindent
For example, if $f(x)\equiv\exp x$, $g(x)\equiv\log(\frac{1}{1-x})$, $b=0$ and $n=2$, 
then
$${\textstyle
T^{fg,\,b}_n(x)=
1\cdot0+(1\cdot1+1\cdot0)x+(1\cdot\frac{1}{2}+1\cdot1+\frac{1}{2}\cdot0)x^2=x+\frac{3}{2}x^2\,.
}
$$

\begin{exer}\label{ex_prodTay}
Determine the Taylor polynomial of $\arctan x\cdot\sin x$ with order $3$ 
and center $0$.    
\end{exer}

\noindent
{\em $\bullet$ The division }of Taylor polynomials is more complicated, but manageable.

\begin{thm}[division~1]\label{thm_diviTaylPol}
Let\index{theorem!division of Taylor polynomials|emph} 
$n\in\N_0$, $f\in\mathcal{F}(M)$ and $b\in M\cap L(M)$. 
Suppose that the Taylor polynomial 
$T^{f,b}_n(x)$ exists,
$${\textstyle
p(x)\equiv T^{f,\,b}_n(x)=\sum_{j=0}^n a_j(x-b)^j
}
$$ 
and that $a_0=f(b)\ne0$. \underline{Then} for $n=0$ we have
$$
T^{1/f,\,b}_n(x)=T^{1/f,\,b}_0(x)=a_0^{-1}\,,
$$
and for $n\ge1$ we get
\begin{eqnarray*}
T^{1/f,\,b}_n(x)&=&\frac{1}{a_0}\bigg(1+
\sum_{\substack{k\in[n],\,e_1,\,e_2,\,\ds,\,e_k\in
\N\\e\equiv e_1+e_2+\ds+e_k\le n}}
(-1)^k a_{e_1}'a_{e_2}'\ds a_{e_k}'(x-b)^e\bigg)\\
&\equiv& r(x),\,\text{ where $a_j'\equiv a_j/a_0$}\,.
\end{eqnarray*}
\end{thm}
\duk
Let $n$, $f$, $b$, $p(x)$ and $a_j$ be as stated, and let $x\to b$. We 
suppose that $a_0\ne0$ and $n=0$. Using the identity $\frac{1}{1+x}=1-\frac{x}{1+x}$ we get that
$$
{\textstyle
\frac{1}{f(x)}=\frac{1}{a_0+o(1)}=a_0^{-1}\frac{1}{1+o(1)}=a_0^{-1}(1+o(1))=a_0^{-1}+o(1)\,.
}
$$
Thus $T^{1/f,b}_0(x)=a_0^{-1}$. For $n\ge1$ we use the identity 
$\frac{1}{1+x}=\sum_{k=0}^n(-1)^k 
x^k+\frac{(-x)^{n+1}}{1+x}$. By Exercises~\ref{ex_arErr1}--\ref{ex_arErr3} we have  
\begin{eqnarray*}
{\textstyle
\frac{1}{f(x)}}&=&
{\textstyle
\frac{1}{p(x)+o((x-b)^n)}=\frac{1}{a_0}\cdot
\frac{1}{1+\sum_{j=1}^n a_j'(x-b)^j+o((x-b)^n)}}\\
&&{\textstyle
\frac{1}{a_0}\cdot\Big(1+\sum_{k=1}^n(-1)^k\big(\sum_{j=1}^n a_j'(x-b)^j+o((x-b)^n)\big)^k\,+}\\
&&{\textstyle +\,o((x-b)^n)\Big)}\\
&=&{\textstyle
\frac{1}{a_0}\cdot\Big(1+\sum_{k=1}^n(-1)^k\big(\sum_{j=1}^n a_j'(x-b)^j\big)^k\Big)+
o((x-b)^n)}\\
&=&r(x)+o((x-b)^n)\,.
\end{eqnarray*}
Proposition~\ref{prop_uniqTaylor}  shows that $T^{1/f,b}_n(x)=r(x)$.
\kduk

\noindent
For example, if $n=3$, $f(x)\equiv2+\log(1+x)$ and $b=0$, then 
$p(x)=2+x-\frac{1}{2}x^2+\frac{1}{3}x^3$, $a_1'=\frac{1}{2}$, $a_2'=-\frac{1}{4}$, $a_3'=\frac{1}{6}$ and 
\begin{eqnarray*}
 T^{1/f,\,0}_3(x)&=&{\textstyle
\frac{1}{2}\big(1-\frac{1}{2}x+\frac{1}{4}x^2-\frac{1}{6}x^3+\big(\frac{1}{2}\big)^2x^2-2\big(\frac{1}{2}\big)\big(\frac{1}{4}\big)x^3-
\big(\frac{1}{2}\big)^3x^3\big)}\\
&=&{\textstyle\frac{1}{2}-\frac{1}{4}x+\frac{1}{4}x^2-\frac{13}{48}x^3\,.}
\end{eqnarray*}

\begin{exer}\label{ex_ratioTayl}
Suppose that $T^{f,b}_3(x)=5+x+x^2+x^3$. Find
$T^{1/f,b}_3(x)$.    
\end{exer}

\noindent
{\em $\bullet$ Laurent Taylor polynomials. }In this passage we
extend Theorem~\ref{thm_diviTaylPol} to the situation when $a_0=f(b)=0$.  

\begin{defi}[Laurent polynomials]\label{def_LaurPol}
A~(real) \underline{Laurent polynomial\index{Laurent polynomials|emph}} with center $b$ ($\in\R$) is a~rational function of the form
$$
l(x)=\sum_{j=m}^na_j(x-b)^j
$$
where $m\le n$ are integers and $a_j\in\R$.
\end{defi}
The adjective ``Laurent'' refers to the French mathematician {\em Pierre A. Laurent (1813--1854)\index{Laurent, Pierre A.}} who introduced power series expansions with negative exponents.

\begin{exer}\label{ex_defDomLauP}
What is $M(l(x))$?    
\end{exer}

\begin{defi}[Laurent Taylor polynomials]\label{def_LayTayPol}
Let $m,n\in\Z$ with $m\le n$, $M\sus\R$, $b\in L(M)$ and let $f(x)\in\mathcal{F}(M)$. If 
$a_j$ for $j=m,m+1,\ds,n$ are $n-m+1$ real numbers such that
$$
{\textstyle
f(x)=\sum_{j=m}^n a_j(x-b)^j+o((x-b)^n)\ \ (x\to b)\,,
}
$$
we say that 
$${\textstyle
\text{$T_{m,\,n}^{f,\,b}(x)\equiv\sum_{j=m}^n a_j(x-b)^j$
is a~\underline{Laurent Taylor polynomial\index{Laurent Taylor 
polynomials, $T_{m,n}^{f,b}(x)$|emph}}}}\label{LauTaylorPol}
$$
of the function $f(x)$ with orders $m,n$ and center $b$.    
\end{defi}
Clearly, for $n\in\N_0$ always
$$
T_{0,\,n}^{f,\,b}(x)=T_n^{f,\,b}(x)\,,
$$
if either side is defined. 

\begin{exer}\label{ex_uniqLauTayP}
Show that the Laurent polynomial $T_{m,n}^{f,b}(x)$ is uniquely 
determined by the parameters $m$, $n$, $f$ and $b$.   
\end{exer}
We generalize  Theorem~\ref{thm_diviTaylPol} as 
follows.

\begin{thm}[division~2]\label{thm_diviLauTaylPol}
Let\index{theorem!division of Laurent Taylor polynomials|emph} 
$m,n\in\Z$ with $m\le n$, $f\in\mathcal{F}(M)$ and let $b\in L(M)$. 
Suppose that the Laurent Taylor polynomial 
$T^{f,b}_{m,n}(x)$
exists,
$${\textstyle
p(x)\equiv T^{f,\,b}_{m,\,n}(x)=\sum_{j=m}^n a_j(x-b)^j
}
$$ 
and that not all coefficients $a_j$ are zero. Let 
$l\in\{m,m+1,\ds,n\}$ be the minimum index such that $a_l\ne0$. 
\underline{Then} for $l=n$ we have
$$
T^{1/f,\,b}_{-n,\,-n}(x)=a_l^{-1}(x-b)^{-l}=a_n^{-1}(x-b)^{-n}\,,
$$
and for $l<n$ we get, with $N\equiv n-l$,
\begin{eqnarray*}
&&T^{1/f,\,b}_{-l,\,n-2l}(x)=\frac{(x-b)^{-l}}{a_l}\bigg(1+
\sum_{\substack{k\in[N],\,e_1,\,e_2,\,\ds,\,e_k\in\N\\e\equiv e_1+e_2+\ds+e_k\le N}}(-1)^k a_{e_1}'a_{e_2}'\ds a_{e_k}'(x-b)^e\bigg)\\
&&\text{ where $a_j'\equiv a_{l+j}/a_l$}\,.
\end{eqnarray*}
\end{thm}
\duk
These formulas follow from
Exercises~\ref{ex_arErr1}--\ref{ex_arErr4}, and from
Theorem~\ref{thm_diviTaylPol} and its proof by means of the factorization
$${\textstyle
p(x)=a_l(x-b)^l\cdot\sum_{j=0}^N a_j'(x-b)^j=a_l(x-b)^l\cdot\big(1+\sum_{j=1}^N a_j'(x-b)^j\big)\,.
}
$$
\kduk

\noindent
For example, suppose that $m=-2$, $n=3$, $b=0$  and $f(x)\equiv\frac{1}{x}+\sin x$. We want to find $T^{1/f,0}_{?,?}(x)$. We have
$${\textstyle
T^{f,\,0}_{m,\,n}(x)=x^{-1}+x+\frac{1}{6}x^3=x^{-1}\cdot\big(1+x^2+\frac{1}{6}x^4\big)\,.
}
$$
Therefore $l=-1$, $N=4$, $a_1'=0$, $a_2'=1$,
$a_3'=0$, $a_4'=\frac{1}{6}$ and

$$
T^{1/f,\,0}_{1,\,5}(x)={\textstyle
x\big(1-1x^2-\frac{1}{6}x^4+1^2x^4\big)
=x-x^3+\frac{5}{6}x^5
}
$$

\begin{exer}\label{ex_naDiv2}
Find $T^{1/f,0}_{?,?}(x)$ if $m=0$, $n=3$ and $f(x)=\arctan x$.    
\end{exer}

\noindent
{\em $\bullet$ Composition }of Taylor polynomials is similar to division.

\begin{thm}[composition]\label{thm_compTaylPol}
Let\index{theorem!composition of Taylor polynomials|emph} 
$n\in\N_0$, $f,g\in\mathcal{R}$, $b\in M(g)\cap L(M(g))$ and $g(b)\in M(f)\cap L(M(f))$. 
Suppose that the Taylor polynomials 
$T^{f,g(b)}_n(x)$  and $T^{g,\,b}_n(x)$
exist. We write
$$
{\textstyle
p(x)\equiv T^{f,\,g(b)}_n(x)=\sum_{j=0}^n 
a_j(x-g(b))^j,\ q(x)\equiv T^{g,\,b}_n(x)=\sum_{j=0}^n b_j(x-b)^j\,,}
$$ 
so that $b_0=g(b)$.
\underline{Then} for $n=0$ we have
$$
T^{f\circ g,\,b}_n(x)=
T^{f\circ g,\,b}_0(x)=a_0\,,
$$
and for $n\ge1$ we get 
\begin{eqnarray*}
&&T^{f\circ g,\,b}_n(x)=a_0+ \sum_{\substack{k\in[n],\,e_1,\,e_2,\,
\ds,\,e_k\in\N\\
e\equiv e_1+e_2+\ds+e_k\le n}}
a_kb_{e_1}b_{e_2}\ds b_{e_k}(x-b)^e\equiv r(x)\,.
\end{eqnarray*}
\end{thm}
\duk
Let $n$, $f$, $g$, $b$, $p(x)$, $a_j$, $q(x)$ and $b_j$ be as stated, and let $x\to b$. For $n=0$ the formula is trivial. Let $n\ge1$. 
Using Exercises~\ref{ex_arErr1}--\ref{ex_arErr3} we get from 
$$
f(y)=p(y)+o((y-g(b))^n)\ \  (y\to g(b))
$$
by substituting for $y$ on the left-hand side the function $g(x)$ and on the right-hand side the expression $q(x)+o((x-b)^n)$ ($x\to b$) that for $x\to b$ we indeed have 
\begin{eqnarray*}
f(g(x))&=&p\big(q(x)+o((x-b)^n)\big)+o\big(\big(q(x)+o((x-b)^n)-g(b)\big)^n\big)\\
&=&{\textstyle
a_0+\sum_{k=1}^n a_k\big(\sum_{j=0}^n b_j(x-b)^j+o((x-b)^n)-g(b)\big)^k\,+}\\
&&{\textstyle
+\,o\big(\big(\sum_{j=0}^nb_j(x-b)^j+o((x-b)^n)-g(b)\big)^n\big)
}\\
&=&{\textstyle
a_0+\sum_{k=1}^n a_k\big(\sum_{j=1}^n b_j(x-b)^j+o((x-b)^n)\big)^k\,+}\\
&&{\textstyle
+\,o\big(\big(\sum_{j=1}^nb_j(x-b)^j+o((x-b)^n)\big)^n\big)
}\\
&=&{\textstyle a_0+\sum_{k=1}^n a_k\big(\sum_{j=1}^n b_j(x-b)^j\big)^k+o((x-b)^n)}\\
&=&r(x)+o((x-b)^n)\,.
\end{eqnarray*}
Proposition~\ref{prop_uniqTaylor}  shows that $T^{f\circ g,b}_n(x)=r(x)$.
\kduk
\vspace{-3mm}
\begin{exer}\label{ex_whyTrue1}
Why does the substitution transform 
$$
o((y-g(b))^n)\ \  (y\to g(b)) 
$$
into $o((q(x)+o((x-b)^n)-g(b))^n)$ ($x\to b$)?
\end{exer}

\begin{exer}\label{ex_whyTrue2}
Why
$$
{\textstyle
o\big(\big(\sum_{j=1}^nb_j(x-b)^j+o((x-b)^n)\big)^n\big)=o((x-b)^n)\ \ (x\to b)\,?
}
$$
\end{exer}

For example, if $n=3$, $b=0$ and $f(x)=g(x)\equiv\sin x$, then $g(b)=0$, $a_0=b_0=0$, $a_1=b_1=1$, $a_2=b_2=0$, 
$a_3=b_3=-\frac{1}{6}$ and
$$
T^{f\circ f,\,0}_3(x)={\textstyle
a_0+a_1b_1x+a_1b_3x^3+a_3b_1^3x^3=
x-\frac{1}{3}x^3\,.
}
$$

In another example we determine in two way the Taylor polynomial 
$$
T^{\sin(\cos x),\,0}_2(x)\,.
$$
Using Theorems~\ref{thm_compTaylPol} and \ref{thm_classTayl} we have 
$n=2$, $f(x)=\sin x$, $g(x)=\cos x$, $b=0$, $g(b)=1$, $a_0=\sin 1$, 
$a_1=\cos 1$, $a_2=-\frac{1}{2}\sin 1$, $b_0=1$, $b_1=0$, 
$b_2=-\frac{1}{2}$ and
$$
{\textstyle
T^{\sin(\cos x),\,0}_2(x)=a_0+a_1b_2x^2+
a_2b_1^2x^2=\sin 1-\frac{1}{2}(\cos 1)x^2\,.
}
$$
Or, directly without composing,
we set 
$$
F(x)\equiv\sin(\cos x)
$$ 
and using only Theorem~\ref{thm_classTayl} and 
formulas for differentiation
we get, of course, the same: $F'(x)=-\cos(\cos x)\cdot\sin x$ and $F''(x)=-\sin(\cos x))\cdot\sin x-\cos(\cos x)\cdot\cos x$, so that
$${\textstyle
T^{\sin(\cos x),\,0}_2(x)=F(0)+F'(0)x+\frac{1}{2}F''(0)x^2=\sin 1-\frac{1}{2}(\cos 1)x^2\,.
}
$$

\begin{exer}\label{ex_CompoTayl}
Let $f(x)\equiv T^{\sin x,0}_5(x)$ and $g(x)\equiv T^{\arcsin x,0}_5(x)$. Check by means of the formula in Theorem~\ref{thm_compTaylPol} that
$$
T^{f\circ g,\,0}_5(x)=x+o(x^5)\,.
$$
\end{exer}

\noindent{\em $\bullet$ Global derivative. }Let $n\in\N$. One might 
think that if the Taylor polynomials $T^{f',\,b}_{n-1}(x)$ and 
$T^{f,\,b}_n(x)$ exist, then they satisfy relation
$$
T^{f',\,b}_{n-1}(x)=\big(T^{f,\,b}_n\big)'\,.
$$
However, we already gave examples showing that for
general definition domains this does not hold.  

\begin{prop}[general domains are bad]\label{prop_TaylDeri}
For every $c\in\R$ there exists a~function $f$ such that $f,f',f''
\in\mathcal{F}([0,1]_{\Q})$,
$$
T^{f,\,0}_2(x)=0+0x+0x^2\
$$
and $f'(x)=cx\,|\,[0,1]_{\Q}$. Thus
$$
T^{f',\,0}_1(x)=0+cx\ne0+0x=\big(T^{f,\,0}_2(x)\big)'\,.
$$
\end{prop}
\duk
We take the function $f$ defined in the proof of 
Theorem~\ref{thm_nonClasTay}. The first equality follows from 
Exercise~\ref{ex_proveEqui}:
$$
T^{f',\,0}_1(x)=f'(0)+(f')'(0)x=0+cx\,.
$$
\kduk

\noindent
In contrast with the four previous operations on $\mathcal{R}$, global derivative 
does not yield any relation between Taylor polynomials of $f$ and $f'$. 
Recall, however, Proposition~\ref{prop_TaylFandDerF}. 
It implies that on intervals we do have 
$$
\big(T^{f,b}_n(x)\big)'=T^{f',b}_{n-1}(x)\,.
$$

\begin{cor}[intervals are good]\label{cor_intGood}
Let $n\in\N$, $b<c$ be in $\R$ and let $f,f'\in\mathcal{F}([b,c))$. 
Suppose that $T^{f',b}_{n-1}(x)$ exists. \underline{Then} also $T^{f,b}_n(x)$ exists and 
$$
\big(T^{f,\,b}_n\big)'=T^{f',\,b}_{n-1}(x)\,.
$$
\end{cor}
\duk
This is immediate from Proposition~\ref{prop_TaylFandDerF}.
\kduk
\vspace{-3mm}
\begin{exer}\label{ex_naTaylInterv}
 Extend this result to definition domains $(c,b]$ with $c<b$ and $U(b,\de)$.    
\end{exer}

\noindent{\em $\bullet$ Computing Taylor polynomials by algebra. }We give more convenient algebraic
versions of formulas in Proposition~\ref{prop_addiTP} and 
Theorems~\ref{thm_diviTaylPol} and \ref{thm_compTaylPol}.  
We work in the polynomial ring 
$$
\mathrm{POL_{id}}=\langle
\mathrm{POL},\,k_0(x),\,k_1(x),\,+,\,
\cdot\rangle
$$
which is an integral domain. 

\begin{defi}[mod $x^k$]\label{def_moduloXk}
For $p(x),q(x)\in\mathrm{POL}$ and $k\in\N_0$ we write
$$
p(x)=q(x)\ (\mathrm{mod}\;x^k)\label{equMod}
$$
and say that the polynomials $p(x)$ and $q(x)$ are \underline{equal modulo} 
$x^k$, if in the canonical form of the polynomial $p(x)-q(x)$ the coefficients 
of $x^0$, $x^1$, $\ds$, $x^{k-1}$ are zero. 
\end{defi}
In other words, (canonical forms of) $p(x)$ and $q(x)$ have identical coefficients of $x^j$ for 
$j<k$. We also write
$$
p(x)\ \mathrm{mod}\;x^k\ \ (\in\mathrm{POL})\label{redMod}
$$
to denote the unique polynomial $q(x)$ such that $\deg q(x)<k$ and $p(x)=q(x)$
($\mathrm{mod}\;x^k$).

\begin{exer}\label{ex_compMod}
Prove the next proposition.    
\end{exer}

\begin{prop}[computing mod $x^k$]\label{prop_compMod}
Let $k\in\N_0$ and let $p(x)$, $p_0(x)$,
$q(x)$ and $q_0(x)$ be polynomials such that
$$
p(x)=q(x)\ (\mathrm{mod}\;x^k)\,\text{ and }\,p_0(x)=q_0(x)\ (\mathrm{mod}\;x^k)\,.
$$
\underline{Then}
$$
p(x)+p_0(x)=q(x)+q_0(x)\ (\mathrm{mod}\;x^k)\,\text{ and }\,p(x)p_0(x)=q(x)q_0(x)\ (\mathrm{mod}\;x^k)\,.
$$
\end{prop}

The algebraic versions of Proposition~\ref{prop_addiTP} and
Theorems~\ref{thm_diviTaylPol} and \ref{thm_compTaylPol} are as follows; 
the trivial addition is omitted. For simplicity of notation, we consider
only the center $b=0$. For brevity we do not consider extensions to
Laurent Taylor polynomials.

\begin{prop}[multiplication mod $x^{n+1}$]\label{prop_multModu}
Let $p(x),q(x)\in\mathrm{POL}$ and $n\in\N_0$ be such that $n\ge\deg p(x),\deg q(x)$. \underline{Then} 
$$
(p(x)+o(x^n))\cdot(q(x)+o(x^n))=
\big(p(x)\cdot q(x)\ \mathrm{mod}\;x^{n+1}\big)+o(x^n)\ \ (x\to0)\,.
$$
\end{prop}
\duk
This follows from Exercises~\ref{ex_arErr1}--\ref{ex_arErr3}.
\kduk
\vspace{-3mm}
\begin{exer}\label{ex_explWhat}
Explain in detail what the proposition and the next two theorems 
exactly say. What is the precise meaning of the $o(\cdot)$ terms?
\end{exer}

\begin{thm}[division mod $x^{n+1}$]\label{thm_divModulo}
Let\index{theorem!division mod $x^{n+1}$|emph} 
$p(x)\in\mathrm{POL}$ with $a_0\equiv p(0)\ne0$ and let $n\in\N_0$ be such that $n\ge\deg p(x)$. \underline{Then} 
$$
\frac{1}{p(x)+o(x^n)}=\Big(\frac{1}{a_0}\sum_{k=0}^n\big(1-a_0^{-1}p(x)\big)^k\ \mathrm{mod}\;x^{n+1}\Big)+o(x^n)\ \ (x\to0)\,.
$$

\end{thm}
\duk
Using the identity 
$\frac{1}{1+x}=\sum_{k=0}^n(-1)^k 
x^k+\frac{(-x)^{n+1}}{1+x}$ and Exercises~\ref{ex_arErr1}--\ref{ex_arErr3} we get that
\begin{eqnarray*}
{\textstyle
\frac{1}{p(x)+o(x^n)}}&=&{\textstyle\frac{1}{a_0}\Big(\sum_{k=0}^n\big(1-a_0^{-1}p(x)+o(x^n)\big)^k+
\frac{(1-a_0^{-1}p(x))^{n+1}}{1+x}\Big)}\\
&=&{\textstyle
\Big(\frac{1}{a_0}\sum_{k=0}^n\big(1-a_0^{-1}p(x)\big)^k\ \mathrm{mod}\;x^{n+1}\Big)+o(x^n)\ \ (x\to0)\,.}
\end{eqnarray*}
\kduk

\begin{thm}[composition mod $x^{n+1}$]\label{thm_compModulo}
Let\index{theorem!composition  mod $x^{n+1}$|emph} 
$p(x)$ and $q(x)$ be in $\mathrm{POL}$ and let $n\in\N_0$ be such that $n\ge\deg p(x),\deg q(x)$. We write 
$p(x)=\sum_{j=0}^n a_jx^j$.
\underline{Then} 
\begin{eqnarray*}
&&(p(y)+o((y-q(0))^n))\circ(q(x)+o(x^n))\ \ (y\to q(0),\,x\to0)\\
&&=
\Big(\sum_{k=0}^n a_k\cdot q(x)^k\ \mathrm{mod}\;x^{n+1}\Big)+o(x^n)\ \ (x\to0)\,.
\end{eqnarray*}
\end{thm}
\duk
Using Exercises~\ref{ex_arErr1}--\ref{ex_arErr3} we get that ($y\to g(0)$, $x\to0$)
\begin{eqnarray*}
&&(p(y)+o((y-g(0))^n))\circ(q(x)+o(x^n))\\
&&{\textstyle
=\sum_{k=0}^na_k\big(q(x)+o(x^n)\big)^k+
o\big(\big(q(x)-q(0)+o(x^n)\big)^n\big)}\\
&&{\textstyle
=\big(\sum_{k=0}^n a_k\cdot q(x)^k\ \mathrm{mod}\;x^{n+1}\big)+o(x^n)\ \ (x\to0)\,.
}
\end{eqnarray*}
\kduk

\noindent
{\em $\bullet$ Changing the center of a~polynomial. }This is another 
technique for obtaining new Taylor polynomials from old ones.

\begin{prop}[changing the center]\label{prop_changeCen}
Let $b\in\R$, $n\in\N$ and let $a_j$ for $j=0,1,\ds,n$ be real numbers. \underline{Then} we have equality of polynomials
$$
\sum_{j=0}^n a_jx^j=\sum_{i=0}^n b_i(x-b)^i\,\text{ where }\,b_i=
\sum_{k=0}^{n-i}a_{k+i}\binom{k+i}{i}b^k
$$
\end{prop}
\duk
We obtain this identity by expanding $x^j=((x-b)+b)^j$ via binomial theorem 
(Exercise~\ref{ex_binomVeta}) and then renaming the difference $j-i$ by $k$.
\kduk

\noindent
For example,
$$
1+x+x^2=1+((x-1)+1)+((x-1)+1)^2=
3+3(x-1)+(x-1)^2\,.
$$

\begin{exer}\label{ex_chngeCen}
Change the center $0$ of $1+x+x^2+x^3$ to $-1$.    
\end{exer}

\noindent
{\em $\bullet$ Examples. }We first determine the Taylor polynomial of $\frac{1}{\cos x}$ 
with order $6$ and center $0$. We set $p(x)\equiv -\frac{1}{2}x^2+
\frac{1}{24}x^4-\frac{1}{720}x^6$ and using Theorem~\ref{thm_divModulo} we get 
\begin{eqnarray*}
&&{\textstyle
\frac{1}{\cos x}=\frac{1}{1+p(x)+o(x^6)}=
\big(\sum_{k=0}^3 (-p(x))^k\ \mathrm{mod}\;x^7\big)+o(x^6)=}\\
&&{\textstyle
=1-(-\frac{1}{2}x^2+\frac{1}{24}x^4-\frac{1}{720}x^6)+(-\frac{1}{2}x^2)^2+
2(-\frac{1}{2}x^2)\frac{1}{24}x^4\,-}\\
&&{\textstyle
-\,(-\frac{1}{2}x^2)^3+o(x^6)=1+\frac{1}{2}x^2+\frac{5}{24}x^4+\frac{61}{720}x^6
+o(x^6)\,.
}
\end{eqnarray*}
We restricted the exponent to $k\le3$ because $p(x)^k=0$ (mod $x^7$) if 
$k\ge4$.

Next we determine the Taylor polynomial of $\sqrt{1+\sin x}$ 
with order $5$ and center $0$. We set $q(x)\equiv x-\frac{1}{6}x^3+
\frac{1}{120}x^5$ and using Theorem~\ref{thm_compModulo} we get
\begin{eqnarray*}
&&{\textstyle
\sqrt{1+\sin x}=\big(1+\sum_{j=1}^5\binom{1/2}{j}\cdot q(x)^j\ \mathrm{mod}\;x^6\big)+o(x^5)=}\\
&&{\textstyle=1+\frac{1}{2}q(x)-\frac{1}{8}(x^2-2x\frac{1}{6}x^3)+\frac{1}{16}(x^3-3x^2\frac{1}{6}x^3)\,-}\\
&&{\textstyle
-\,\frac{15}{16\cdot24}x^4+\frac{105}{32\cdot120}x^5+o(x^5)=1+\frac{1}{2}x
-\frac{1}{8}x^2-\frac{1}{48}x^3+\frac{16-15}{16\cdot24}x^4\,+}\\
&&{\textstyle +\,
\frac{16-120+105}{32\cdot120}x^5+o(x^5)=}\\
&&{\textstyle=1+\frac{1}{2}x
-\frac{1}{8}x^2-\frac{1}{48}x^3+\frac{1}{384}x^4+\frac{1}{3840}x^5+o(x^5)\,.
}
\end{eqnarray*}

\begin{exer}\label{ex_najdiNaInt}
Why are these coefficients so simple? Please, work in the spirit of 
Exercise~\ref{ex_fekete2}. (I~admit that I took help from the Internet.)    
\end{exer}

Finally, we determine the Taylor polynomial of the function $\tan x$ 
with order~$6$ and center $0$. Using Proposition~\ref{prop_multModu} and the above Taylor polynomial of $\frac{1}{\cos x}$ we get that
\begin{eqnarray*}
&&{\textstyle\tan x=\frac{\sin x}{\cos x}
=\sin x\cdot\frac{1}{\cos x}=}\\
&&{\textstyle\big((x-\frac{1}{6}x^3+\frac{1}{120}x^5)\cdot(1+\frac{1}{2}x^2+\frac{5}{24}x^4+\frac{61}{720}x^6)
\ \mathrm{mod}\;x^7\big)+o(x^6)=}\\
&&{\textstyle
=x+\frac{1}{3}x^3+\frac{1-10+25}{120}x^5+o(x^6)=x+\frac{1}{3}x^3+
\frac{2}{15}x^5+0x^6+o(x^6)\,.
}
\end{eqnarray*}
There is a~formula for these coefficients in terms of the Bernoulli numbers.\index{Bernoulli number, $B_k$}

\begin{exer}\label{ex_explEven}
Show that every coefficient of $x^{2n}$ in Taylor polynomials of $\tan x$ is zero.    
\end{exer}

\begin{exer}\label{ex_sTretiOdm}
Find the Taylor polynomial with order $5$ and center $0$ of the function
$$
\sqrt[3]{1+\sin x}\,.
$$
\end{exer}

\noindent
{\em $\bullet$ Algebraic complexity of computations of Taylor polynomials. }Since this book 
is based on lectures for students in the School of Computer Science of MFF 
UK (Faculty of Mathematics and Physics of Charles University in 
Prague)\label{mffuk}, we briefly discuss algorithmic aspects of 
computing Taylor polynomials. We do not consider 
the problem to compute the Taylor polynomial of a~given function, but the problem
to compute the Taylor polynomial  of the sum, product, $\ds$ of two
functions if their Taylor polynomials are known. We consider just
the {\em arithmetic complexity\index{arithmetic complexity 
of an algorithm}} of a~computation, which is the number of 
arithmetic operations it does. The formulas in 
Theorems~\ref{thm_diviTaylPol} and \ref{thm_compTaylPol} are explicit but
do not provide efficient algorithms. Therefore we turn to formulas in 
Theorems~\ref{thm_divModulo} and \ref{thm_compModulo}. For simplicity of 
notation we treat only the case $b=0$, but everything could be easily extended to general center $b\in\R$.

\begin{thm}[computing Taylor polynomials]\label{thm_aritCompl}
Let\index{theorem!computing Taylor polynomials|emph}  $n\in\N$, $f,g\in\mathcal{R}$ and let the Taylor 
polynomials 
$$
T^{f,\,0}_n(x)=\sum_{j=0}^n a_jx^j\,\text{ and }\,T^{g,\,0}_n(x)=\sum_{j=0}^n b_jx^j
$$
exist and be given. The following holds.
\begin{enumerate}
\item The coefficients of $T^{f+g,0}_n(x)$ can be computed from $a_0$, $\ds$, $a_n$, $b_0$, $\ds$, $b_n$ in $n+1$ arithmetic operations.
\item The coefficients of $T^{fg,0}_n(x)$ can be computed from 
$a_0$, $\ds$, $a_n$, $b_0$, $\ds$, $b_n$ in $O(n^3)$ ($n\in\N$) arithmetic operations.
\item For $a_0\ne0$ the coefficients of $T^{1/f,0}_n(x)$ can be computed from 
$a_0$, $\ds$, $a_n$ in $O(n^5)$ ($n\in\N$) arithmetic operations.
\item Let in addition $g(0)=0$.
The coefficients of $T^{f\circ g,b}_n(x)$ can be computed from 
$a_0$, $\ds$, $a_n$, $b_0$, $\ds$, $b_n$ in $O(n^5)$ ($n\in\N$) arithmetic operations.
\end{enumerate}
\end{thm}
\duk
Parts~1 and~2 are left  for Exercise~\ref{ex_provItYours}. 3. Let $a_0\ne0$ and 
$$
{\textstyle
T^{1/f,\,0}_n(x)=\sum_{j=0}^n c_jx^j\,.
}
$$
We compute $c_0$, $\ds$, $c_n$ from $a_0$, $\ds$, $a_n$
by means of the formula in 
Theorem~\ref{thm_divModulo}. Let  $k\in\N$, $k\le n$, be fixed. Using the 
algorithm of part~2 we compute the 
$k$-fold product of polynomials in the summand
$${\textstyle
\big(1-a_0^{-1}\sum_{j=0}^n a_jx^j\big)^k\ \mathrm{mod}\;x^{n+1}
}
$$
in $O(kn^3)$ ($n\in\N$) arithmetic operations. Then using the 
algorithm of part~1 we compute the product of $a_0^{-1}$ with the sum of $n+1$ polynomials (with degrees at most $n$) modulo $x^{n+1}$  in Theorem~\ref{thm_divModulo} in 
$$
{\textstyle
n+2+n(n+1)+\sum_{k=1}^n O(kn^3)=O(n^5)
}
$$
arithmetic operations.

4. Let $g(0)=0$ and 
$$
{\textstyle
T^{f\circ g,0}_n(x)=\sum_{j=0}^n c_jx^j\,.
}
$$
We compute $c_0$, $\ds$, $c_n$ from $a_0$, $\ds$, $a_n$, $b_0$, $\ds$, $b_n$
by means of the formula in 
Theorem~\ref{thm_compModulo}. Using the 
algorithm of part~2 we compute the product of $a_k$ with the 
$k$-fold product of polynomials in the summand
$${\textstyle
a_k\cdot \big(\sum_{j=0}^n b_jx^j\big)^k\ \mathrm{mod}\;x^{n+1}
}
$$
in $n+1+O(kn^3)=O(kn^3)$ ($n\in\N$) arithmetic operations. Then using the 
algorithm of part~1 we compute the sum of $n+1$ polynomials (with degrees at 
most $n$) modulo $x^{n+1}$  in Theorem~\ref{thm_compModulo} in
$$
{\textstyle
n(n+1)+\sum_{k=1}^n O(kn^3)=O(n^5)
}
$$
arithmetic operations.
\kduk
\vspace{-3mm}
\begin{exer}\label{ex_provItYours}
Prove parts~1 and~2 of the theorem.    
\end{exer}
See Chapter~2 of the monograph \cite{burg_al} for more efficient polynomial arithmetic.

\chapter[Real analytic functions]{Real analytic functions}\label{chap_pr9apul}

\section[Taylor series]{Taylor series}\label{sec_TaylorSeries}

In this section we derive explicit formulas for Taylor remainder, which
is the difference of a~function and its Taylor polynomial. We define 
Taylor series and determine sums of Taylor series for several elementary functions. 

\medskip\noindent
{\em $\bullet$ Taylor remainders. }We already computed with them a~lot, but only in the $o((x-b)^n)$ form.

\begin{defi}[Taylor remainder]\label{def_TaylRema}
Suppose that a~function $f\in\mathcal{R}$ has the Taylor 
polynomial $T^{f,b}_n(x)$, so that $b\in M(f)\cap L(M(f))$. Then the function
$$
R^{f,\,b}_n(x)\equiv f(x)-T^{f,\,b}_n(x)\ \ (\in\mathcal{F}(M(f)))\label{remainder}
$$
is the \underline{Taylor remainder\index{Taylor remainders, $R^{f,b}_n(x)$|emph}} of $f$ with order $n$ and center $b$.
\end{defi}
\vspace{-3mm}
\begin{exer}\label{ex_pomocSko}
We have 
$$
R^{f,\,b}_n(x)=o((x-b)^n)\ \ (x\to b)\,.
$$
\end{exer}
This asymptotics of the Taylor remainder 
is ineffective, it provides no explicit bound on it. We improve upon it by obtaining explicit formulas for $R^{f,b}_n(x)$. 

\medskip\noindent
{\em $\bullet$ Simple Taylor theorem. }First we give 
a~simple form of the Taylor theorem. Then we explain why a~more 
complicated but stronger form is needed, and state and prove this form.

\begin{thm}[simple Taylor theorem]\label{thm_simpleTay}
Suppose\index{theorem!simple Taylor theorem|emph} 
that $n\in\N_0$, $I\sus\R$ is a~nonempty open interval, $p>0$ is 
a~real number and that $f$, $f'$, $\ds$, $f^{(n+1)}$ are in
$\mathcal{F}(I)$. \underline{Then} for every two distinct points $b,x\in I$ there is a~point $c$ between them such that
$$
\big|R^{f,\,b}_n(x)\big|=
\frac{|f^{(n+1)}(c)|}{n!\cdot p}\cdot|x-c|^{n+1-p}\cdot|x-b|^p\,.
$$
\end{thm}
\duk
Let $n$, $I$, $p$ and $f$ be as stated. We take distinct points $b,x\in I$ and 
first assume that $b<x$. We define auxiliary functions
$$
\varphi(t)\equiv(x-t)^p\,|\,[b,\,x]
$$
and (with $0^0\equiv1$)
$${\textstyle
F(t)\equiv f(x)-\sum_{i=0}^n\frac{1}{i!}f^{(i)}(t)\cdot(x-t)^i\,|\,[b,\,x]\,.
}
$$
We have $\varphi,F\in\mathcal{C}([b,x])$, $\varphi(x)=0$, $\varphi(b)=(x-b)^p$, $F(x)=0$ and $F(b)=R^{f,b}_n(x)$. Also,
\begin{eqnarray*}
F'(t)&=&{\textstyle-f'(t)-\sum_{i=1}^n\big(\frac{1}{i!}f^{(i+1)}(t)\cdot(x-t)^i-
\frac{1}{(i-1)!}f^{(i)}(t)\cdot(x-t)^{i-1}\big)}\\
&=&{\textstyle
-\frac{1}{n!}f^{(n+1)}(t)\cdot(x-t)^n\,.}
\end{eqnarray*}
By Cauchy's Theorem~\ref{thm_Cauchy} there exists a~point $c\in(b,x)$ such that
$$
{\textstyle
\frac{R^{f,\,b}_n(x)}{(x-b)^p}=\frac{F(x)-F(b)}{\varphi(x)-\varphi(b)}=\frac{F'(c)}{\varphi'(c)}=
\frac{f^{(n+1)}(c)\cdot(x-c)^n}{n!\cdot p(x-c)^{p-1}}\,.
}
$$
Solving this equation for $|R^{f,b}_n(x)|$ we get the stated formula for the (absolute value of) Taylor remainder. 

Let $x<b$. We argue as before, with the modified auxiliary function
$$
\varphi(t)\equiv(t-x)^p\,|\,[x,\,b]
$$
and the same $F(t)$ (restricted to $[x,b]$). Now we get 
$$
{\textstyle
\frac{R^{f,\,b}_n(x)}{(b-x)^p}=\frac{F(x)-F(b)}{\varphi(x)-\varphi(b)}=\frac{F'(c)}{\varphi'(c)}=
-\frac{f^{(n+1)}(c)\cdot(x-c)^n}{n!\cdot p(c-x)^{p-1}}\,.
}
$$
Solving this equation for $|R^{f,b}_n(x)|$ we get again the stated formula.
\kduk

\noindent
The formula for $|R^{f,b}_n(x)|$ in the previous and the next theorem is called 
the \underline{Schl\"omilch remainder\index{remainder of the Taylor series!schlomilch form@Schl\"omilch form|emph}}. It bears the name of the German mathematician {\em Oskar X. 
Schl\"omilch (1823--1901)\index{schlomilch@Schl\"omilch, 
Oskar X.}}. To get a~simple formula, we sacrificed the sign of $R^{f,b}_n(x)$ (we do not need it 
for determining if $R^{f,b}_n(x)\to0$
as $n\to\infty$). With some effort it can be recovered.

\begin{exer}\label{ex+morePr}
State the formula in the previous theorem in  a~more precise form as $R^{f,b}_n(x)=\cdots$.   
\end{exer}
Below, for complete treatment of the binomial series we need the Schl\"omilch 
remainder with non-integral values of the exponent $p$. For $p\in\N$ we get simple well known formulas including the sign of $R^{f,b}_n(x)$.

\begin{cor}[Lagrange and Cauchy remainders~1]\label{cor_LagCauRem}
Let $n\in\N_0$, $I\sus\R$ be a~nonempty open interval and let $f$, $f'$, $\ds$, $f^{(n+1)}$ be in
$\mathcal{F}(I)$. \underline{Then}
the following holds.
\begin{enumerate}
\item For every two distinct points $b,x\in I$ there is a~point $c$ between them such that we have the \underline{Lagrange remainder\index{remainder of the Taylor series!Lagrange form|emph}}
$$
R^{f,\,b}_n(x)=\frac{f^{(n+1)}(c)}{(n+1)!}\cdot(x-b)^{n+1}\,.
$$
\item For every two distinct points $b,x\in I$ there is a~point $c$ between them such that we have the \underline{Cauchy remainder\index{remainder of the Taylor series!Cauchy form|emph}}
$$
R^{f,\,b}_n(x)=\frac{f^{(n+1)}(c)}{n!}\cdot(x-c)^n(x-b)\,.
$$
\end{enumerate}    
\end{cor}
\duk
These remainders follow from the previous proof if we set $p=n+1$ and $p=1$, respectively. 
Then we can take (the restriction of) $\varphi(t)\equiv(x-t)^p$ in both cases
$x<b$ and $b<x$ and we get the stated formulas.
\kduk

Theorem~\ref{thm_simpleTay} has the drawback that in some cases, when we need it, it does not apply. For example, we would like to 
set in the Maclaurin series
$$
{\textstyle
(1+x)^a=\sum_{n=0}^{\infty}\binom{a}{n}x^n
}
$$
the variable $x$ to $-1$ and deduce for every real $a>0$ the summation
$$
{\textstyle
\sum_{n=0}^{\infty}\binom{a}{n}(-1)^n=0\,,
}
$$
but Theorem~\ref{thm_simpleTay} does not give it because finite derivatives 
$$
\big((1+x)^a\big)^{(n)}(-1)
$$ 
do not exist for large $n$.

\begin{exer}\label{ex_applBinom}
Prove the above summation in the case when $a\in\N$.    
\end{exer}

\noindent
{\em $\bullet$ Advanced Taylor theorem. }We fix this problem with the following theorem.

\begin{thm}[advanced Taylor theorem]\label{thm_advanTay}
Let\index{theorem!advanced Taylor theorem|emph} $n\in\N_0$, let 
$b\ne x$ and $p>0$ be real numbers, and let $I$ be the closed interval with 
endpoints $b$ and $x$. Suppose that $f\in\mathcal{C}(I)$, that $M(f^{(i)})\supset I\setminus\{x\}$ for $i=1,2,\ds,n+1$ and that
for every $i=1,2,\ds,n$ we have
$$
\lim_{t\to x}f^{(i)}(t)\cdot(t-x)^i=0\,.
$$
\underline{Then} there is a~point $c\in I^0$ such that
$$
\big|R^{f,\,b}_n(x)\big|=
\frac{|f^{(n+1)}(c)|}{n!\cdot p}\cdot|x-c|^{n+1-p}\cdot|x-b|^p\,.
$$
\end{thm}
\duk
Let $n$, $b$, $x$, $p$ and $f$ be as stated. We argue as in the proof of 
Theorem~\ref{thm_simpleTay} and use the same functions $\varphi(t)$ and $F(t)$ in $\mathcal{F}(I)$, with a~small modification. The function
$\varphi(t)$ is as in that proof. The function $F(t)$ is defined on 
$I\setminus\{x\}$ by the formula in that proof and we set $F(x)\equiv 
0$. It follows from the assumptions that then $F\in\mathcal{C}(I)$. 
Else, we argue as in that proof.
\kduk

\noindent
The advanced Taylor theorem, or rather the proof of it, again provides us
Lagrange and Cauchy remainders.

\begin{exer}\label{ex_pulkaOkoli}
Prove the next corollary.
\end{exer}

\begin{cor}[Lagrange and Cauchy remainders~2]\label{cor_LagCauRem2}
Let\index{theorem!advanced Taylor theorem|emph} $n\in\N_0$, let 
$b\ne x$ be real numbers, and let $I$ be the closed interval with 
endpoints $b$ and $x$. Suppose that $f\in\mathcal{C}(I)$, that $M(f^{(i)})\supset I\setminus\{x\}$ for $i=1,2,\ds,n+1$ and that
for every $i=1,2,\ds,n$ we have
$$
\lim_{t\to x}f^{(i)}(t)\cdot(t-x)^i=0\,.
$$ \underline{Then}
the following holds.
\begin{enumerate}
\item There is a~point $c\in I^0$ such that we have the \underline{Lagrange remainder\index{remainder of the Taylor series!Lagrange form|emph}}
$$
R^{f,\,b}_n(x)=\frac{f^{(n+1)}(c)}{(n+1)!}\cdot(x-b)^{n+1}\,.
$$
\item There is a~point $c\in I^0$ such that we have the \underline{Cauchy remainder\index{remainder of the Taylor series!Cauchy form|emph}}
$$
R^{f,\,b}_n(x)=\frac{f^{(n+1)}(c)}{n!}\cdot(x-c)^n(x-b)\,.
$$
\end{enumerate}    
\end{cor}

\noindent
{\em $\bullet$ Taylor series. }We define the Taylor series of a~function as the union of its Taylor polynomials of all orders. By Exercise~\ref{ex_iniSumTP} such 
definition is correct\,---\,the coefficient $a_j$ does not depend on 
$n$.

\begin{defi}[Taylor series]\label{def_TaylorSer}
Let $f\in\mathcal{F}(M)$ and $b\in M\cap L(M)$. Suppose that for every $n\in\N_0$ 
there exists the Taylor polynomial
$$
T_n^{f,\,b}(x)=
\sum_{j=0}^n a_j(x-b)^j\,.
$$
Then we call the series ($x\in\R$)
$$
T^{f,\,b}(x)\equiv\sum_{n=0}^{\infty}a_n(x-b)^n
\label{Tfbx}
$$
the \underline{Taylor series\index{Taylor series, $T^{f,b}(x)$|emph}} of $f(x)$
with center $b$. In the case when $b=0$ we often speak of \underline{Maclaurin series\index{Maclaurin series, $T^{f,0}(x)$|emph}} of $f(x)$.
\end{defi}

\noindent
Taylor polynomials and Taylor series are named after the English 
mathematician {\em Brook Taylor (1685--1731)\index{Taylor, Brook}}. Maclaurin 
series refer to the Scottish  mathematician {\em Colin Maclaurin (1698--1746)\index{Maclaurin, Colin}}. 
We shall investigate for which $x\in\R$ we have $T^{f,b}(x)=f(x)$. 
Thus we introduce the set
\begin{eqnarray*}
E(T^{f,b}(x))&\equiv&
\{a\in\R\cc\;\text{the series $T^{f,b}(a)$ converges and $T^{f,b}(a)=f(a)$}\}\label{ETayl}\\
&=&\{a\in\R\cc\;\lim_{n\to\infty}R^{f,\,b}_n(a)=0\}\,.
\end{eqnarray*}
For example,
always $b\in E(T^{f,b}(x))$.

\begin{exer}\label{ex_whyTaylAtb}
Why? Why do we have in the 
formula for the Taylor series that $(b-b)^0=0^0=1$?    
\end{exer}

We remind the standard notation that ($M\sus\R$)
$$
\mathcal{C}^{\infty}(M)\equiv\{f\in\mathcal{R}\cc\;\text{$f^{(j)}\in\mathcal{F}(M)$ for every $j\in\N_0$}\}\,.\label{CmInf}
$$
If $M\not\sus L(M)$ then 
$\mathcal{C}^{\infty}(M)=\emptyset$. 
Theorem~\ref{thm_classTayl} has an immediate corollary.

\begin{cor}[classical Taylor series]\label{cor_classTaySer}
Let $f\in\mathcal{C}^{\infty}(M)$ where $M\sus\R$ is a~nonempty open set. \underline{Then} for every $b\in M$ the function $f$ has the Taylor series with center $b$ and
$$
T^{f,\,b}(x)=\sum_{n=0}^{\infty}\frac{f^{(n)}(b)}{n!}\cdot(x-b)^n\,.
$$
\end{cor}

\noindent
{\em $\bullet$ Taylor series of $\exp x$, $\cos x$ and $\sin x$. }We determine 
sums of Taylor series of several elementary functions and begin with 
$\exp x$, $\cos x$ and $\sin x$. These functions are in 
$\mathcal{C}^{\infty}(\R)$ and their Taylor series are therefore given for 
every center $b\in\R$ by the formula in Corollary~\ref{cor_classTaySer}.

\begin{thm}[$\mathrm{e}^x$, $\cos x$ and $\sin x$]\label{thm_TaylSerExp}   
Let\index{theorem!Taylor series of $\mathrm{e}^x$, $\cos x$ and $\sin 
x$|emph} 
$f(x)\in\{\exp x,\cos x,\sin x\}$.
\underline{Then}
$$
f(x)=\sum_{n=0}^{\infty}\frac{f^{(n)}(b)}{n!}\cdot(x-b)^n
$$
for every $b,x\in\R$.
\end{thm}
\duk
Let $n\in\N$, $b,x\in\R$ with $b<x$ and $f(x)=\exp x$. Using Lagrange's 
remainder in Theorem~\ref{thm_simpleTay} we get
$${\textstyle
R_n^{f,\,b}(x)=(n+1)!^{-1}\cdot\mathrm{e}^c\cdot(x-b)^{n+1}\,\text{ for some }\,c=c(n)\in(b,\,x)\,.
}
$$
Thus for fixed $b$ and $x$ we have $\lim_{n\to\infty}R_n^{f,b}(x)=0$ and the sum of $T^{f,b}(x)$ equals $f(x)=\exp x$. For $x<b$ and/or $f(x)$ equal to $\cos x$ or $\sin x$ we proceed in a~similar way.
\kduk

\noindent
Note that the definitoric series of $\exp x$, $\cos x$ and $\sin x$ in Section~\ref{sec_elemenFce} are their Taylor series with center $0$:
$$
{\textstyle
\exp(a)=T^{\exp x,\,0}(a)=
\sum_{n=0}^{\infty}\frac{a^n}{n!}\,\text{ for every }\,a\in\R\,,
}
$$
and similarly
$$
{\textstyle
\cos(a)=T^{\cos x,\,0}(a)=
\sum_{n=0}^{\infty}\frac{(-1)^na^{2n}}{(2n)!},\ \sin(a)=T^{\sin x,\,0}(a)=
\sum_{n=0}^{\infty}\frac{(-1)^na^{2n+1}}{(2n+1)!}
}
$$
for every $a\in\R$.

\begin{exer}\label{ex_jeToAbscon}
All series in the theorem are abscon. 
\end{exer}

\noindent
{\em $\bullet$ Binomial series. }The next Maclaurin\index{Maclaurin, Colin} series of the real power 
$(1+x)^a$ with $x,a\in\R$ (see Definitions~\ref{def_aNaB} and 
\ref{def_anaM}) is called the \underline{binomial series\index{binomial series|emph}\index{series!binomial|emph}}. It is due to {\em Isaac Newton (1642/43--1727)\index{Newton, 
Isaac}} who was one of the main creators of modern science
(\cite{woot}).

\begin{thm}[binomial series]\label{thm_NewtonBin}
Let\index{theorem!binomial series|emph} 
$a\in\R$ and $f(x)\equiv(1+x)^a$. \underline{Then}
$$
T^{f,\,0}(x)=
\sum_{n=0}^{\infty}\binom{a}{n}\cdot x^n\,.
$$
For $x\in\R$ we have
$$
T^{f,\,0}(x)=(1+x)^a\iff x\in E(a)\equiv E(T^{f,0}(x))
$$
and the set $E(a)$ is as follows.
\begin{enumerate}
\item If $a\in\N_0$, \underline{then} $E(a)=\R$ and the binomial series is a~finite sum. If $a\not\in\N_0$, then the binomial series is an infinite series.
\item If $a>0$ but $a\not\in\N$, \underline{then} $E(a)=[-1,1]$.
\item If $a\in(-1,0)$, \underline{then} $E(a)=(-1,1]$.
\item If $a\le-1$, \underline{then} $E(a)=(-1,1)$. 
\end{enumerate}
\end{thm}
\duk
The form of $T^{f,0}(x)$ follows from Proposition~\ref{prop_onePlusx}.
Part~1 is easy, parts 2--4 are harder. If $a\in\N_0$ then $\binom{a}{n}=0$ 
for $n>a$ and the series is  a~finite sum. The equality $T^{f,0}(x)=(1+x)^a$ 
for every $x\in\R$ then follows from the binomial theorem 
(Exercise~\ref{ex_binomVeta}). If $a\not\in\N_0$ then $\binom{a}{n}\ne0$
for every $n\in\N_0$.

In the rest of the proof we assume that $a\in\R\setminus\N_0$. It is not hard 
to show (Exercise~\ref{ex_bindol}) that for some constant $c<0$,
$$
{\textstyle
\binom{a}{n}\gg n^c\ \ (n\in\N)\,.
}
$$
This bound and comparison with the geometric series give
$$
\text{$E(a)\sus[-1,1]$ for every $a\in\R\setminus\N_0$}\,. 
$$
We show that 
$$
\text{$(-1,1)\sus E(a)$ for every $a\in\R\setminus\N_0$}\,.
$$
As we know, $0\in 
E(a)$. If $x\in(-1,1)$, $x\ne0$, then by using Corollary~\ref{cor_LagCauRem} with $b=0$, $f(t)=(1+t)^a$ and the Cauchy remainder we get 
$${\textstyle
\big|R^{f,\,0}_n(x)\big|=
(n+1)\big|\binom{a}{n+1}\big|\cdot
(1+c)^{a-n-1}\cdot|x-c|^n\cdot |x|\to0,\ n\to\infty\,.
}
$$
Indeed, as for the first factor, by Exercise~\ref{ex_binhor}
$$
{\textstyle
\binom{a}{n+1}\ll n^d\ \ (n\in\N)
}
$$
for some constant $d>0$.
The number $c=c(n)\in(-1,1)$ lies between $0$ and $x$, so that 
$$
(1+c)^{a-1}\le\max\big((1-|x|)^{a-1},\,1\big)
$$
and
$$
{\textstyle
0<\frac{|x-c|}{1+c}\le\frac{|x|-|c|}{1-|c|}=1-\frac{1-|x|}{1-|c|}\le
1-\frac{1-|x|}{1}=|x|<1\,.
}
$$
We see that $R^{f,\,0}_n(x)$ goes exponentially fast to $0$ with 
$n\to\infty$. Thus $T^{f,0}(x)=f(x)$ for every $x\in(-1,1)$.

It remains to determine for $a\in\R\setminus\N_0$ if $\pm1\in E(a)$. Below we obtain for every $a\in\R\setminus\N_0$ the 
asymptotics
$$
{\textstyle
\big|\binom{a}{n}\big|=\Theta(1)\cdot n^{-a-1}\ \ (n\in\N)\,.
}
$$
The factor $\Theta(1)$ is a~function $g\cc\N\to(0,+\infty)$ such that 
$$
\text{$c\le g(n)\le d$ for every $n\in\N$ and some constants $0<c<d$}\,.
$$
We postpone the proof of the asymptotics and finish the
determination of the set $E(a)$. Let $x=1$. Then we have series
$$
{\textstyle
\sum_{n=0}^\infty\binom{a}{n}\,.
}
$$
For $a\le-1$ and $n\in\N_0$ we have the bound
$$
{\textstyle
\big|\binom{a}{n}\big|\ge\big|\binom{-1}{n}\big|=1
}
$$
and the series diverges. Thus $1\not\in E(a)$ for $a\le-1$. Let 
$a>-1$, $a\not\in\N_0$ and $x=1$. Using 
Corollary~\ref{cor_LagCauRem} with $b=0$, $f(t)=(1+t)^a$ and the Lagrange remainder we get 
$${\textstyle
\big|R^{f,\,0}_n(1)\big|=
\big|\binom{a}{n+1}\big|\cdot
(1+c)^{a-n-1}\cdot 1^{n+1}\to0,\ n\to\infty\,,
}
$$
because $0<c=c(n)<1$, $|\binom{a}{n}|=\Theta(1)\cdot n^{-a-1}$ and $-a-1<0$. Thus $1\in E(a)$ for $a>-1$. 

Let $x=-1$. Then we have series
$$
{\textstyle
\sum_{n=0}^\infty\binom{a}{n}(-1)^n\,.
}
$$
Despite the appearance it is not an alternating  series because for 
$n>a+1$ the sign of the summand does not change. For $a<0$ we have $-a-
1>-1$ and by the asymptotics of  $|\binom{a}{n}|$ the series diverges 
because $\zeta(s)$ diverges for $s\le 1$ (Exercise~\ref{zetaSmen1}). Thus 
$-1\not\in E(a)$ for $a<0$. Let $a>0$, $a\not\in\N$ and $x=-1$. We 
use Theorem~\ref{thm_advanTay} (with the Schl\"omilch remainder) with $b=0$, $f(t)=(1+t)^a$ and any 
$p\in(0,a)$. We can use this theorem because $f\in\mathcal{C}([-1,0])$ and
$$
\lim_{t\to x}f^{(i)}(t)\cdot(t-x)^i=
\lim_{t\to -1}f^{(i)}(t)\cdot(t+1)^i=0
$$
for every $i\in\N$, as $a>0$ and 
$$
f^{(i)}(t)=a(a-1)\ds(a-i+1)\cdot(1+t)^{a-i}\,.
$$
We get
$${\textstyle
\big|R^{f,\,0}_n(-1)\big|=
\frac{n+1}{p}\big|\binom{a}{n+1}\big|\cdot
(1+c)^{a-n-1}\cdot(c+1)^{n+1-p}\cdot 1^p\to0,\ n\to\infty\,,
}
$$
because the first factor is $\Theta(1)\cdot n^{-a}$ with $a>0$, we have $c=c(n)\in(-1,0)$ 
and the product $(1+c)^{a-p}$ of the second and third factor lies in 
$(0,1)$ as $a-p>0$. Thus $-1\in E(a)$ for $a>0$. This completes the proof of parts~2--4.

We derive the postponed asymptotics of 
$|\binom{a}{n}|$ ($n\in\N$). We assume that 
$a\in\R\setminus\N_0$. For $m\in\N_0$ we have
$$
{\textstyle
S(a,\,m)\equiv\binom{a}{m}^{-1}\cdot\binom{a}{m+1}=
\frac{a-m}{m+1}=-1+\frac{a+1}{m+1}\,.
}
$$
Since $\log(1+x)=x+O(x^2)$ ($|x|\le\frac{1}{2}$), we have 
$${\textstyle
\log(|S(a,m)|)=-\frac{a+1}{m+1}+O_a(m^{-2})\ \ (m\ge2|a|+2)\,.
}
$$
We set $N\equiv\lceil2|a|+2\rceil$ and $K\equiv|\prod_{m=0}^{N-1}S(a,m)|$. For every $n\ge N+1$ we have
\begin{eqnarray*}
{\textstyle
\log\big(\big|\binom{a}{n}\big|\big)}&=&{\textstyle
\log\big(\big|\prod_{m=0}^{n-1}S(a,\,m)\big|\big)=\log K-(a+1)\sum_{m=N}^{n-1}\frac{1}{m+1}\,+}\\
&&{\textstyle +\,O_a\big(\sum_{m=N}^{n-1}m^{-2}\big)
}\,.
\end{eqnarray*}
We recall the asymptotics 
$${\textstyle
\sum_{m=0}^{n-1}\frac{1}{m+1}=\log n+\ga+O(n^{-1})\ \ (n\in\N) 
}
$$
in Theorem~\ref{thm_harm_cial}. 
Hence, in view of convergence of $\zeta(2)$, 
$$
{\textstyle\log\big(\big|\binom{a}{n}\big|\big)=-(a+1)\log n+O(1)\ \ (n\in\N)\,.
}
$$
Applying $\exp x$ we get by Exercise~\ref{ex_onTheta} that
$$
{\textstyle\big|\binom{a}{n}\big|=\Theta(1)\cdot n^{-a-1}\ \ (n\in\N)\,.
}
$$
\kduk
\vspace{-3mm}
\begin{exer}\label{ex_onTheta}
Show that $\mathrm{e}^{O(1)}=\Theta(1)$.    
\end{exer}

\begin{exer}\label{ex_bindol}
 For every $a\in\R\setminus\N_0$ there is a~$c<0$ such that 
 $$
 {\textstyle
 \binom{a}{n}\gg n^c\ \ (n\in\N)\,.
 }
 $$
\end{exer}

\begin{exer}\label{ex_binhor}
 For every $a\in\R$ there is a~$d>0$ such that 
 $$
 {\textstyle
 \binom{a}{n}\ll n^d\ \ (n\in\N)\,.
 }
 $$
\end{exer}

\begin{exer}\label{ex_NewtSerAbscon}
Which of the binomial series are abscon?  
\end{exer}

\begin{exer}\label{ex_NewtSerCen}
Find $E(F^{f,b}(x))$ with $f(x)=(1+x)^a$ for general center $b$.    
\end{exer}

\noindent
{\em $\bullet$ Nice summations. }Using the binomial series we obtain families of summation identities. 

\begin{exer}\label{ex_SummIden}
Prove the following corollary.    
\end{exer}

\begin{cor}[sums of binomial series]\label{cor_sumsBinSer}
Let $a\in\R$. The next identities hold. 
\begin{enumerate}
\item If $a\in\N_0$ \underline{then}
$(1+x)^a=\sum_{n=0}^{\infty}\binom{a}{n}x^n=\sum_{n=0}^a\binom{a}{n}x^n$ for every $x\in\R$.
\item If $a\le -1$ \underline{then}
$(1+x)^a=\sum_{n=0}^{\infty}\binom{a}{n}x^n$ for every $x\in(-1,1)$.
\item If $-1<a<0$ \underline{then}
$(1+x)^a=\sum_{n=0}^{\infty}\binom{a}{n}x^n$ for every $x\in(-1,1]$.
\item If $a>0$ \underline{then}
$(1+x)^a=\sum_{n=0}^{\infty}\binom{a}{n}x^n$ for every $x\in[-1,1]$.
\end{enumerate}
\end{cor}

\noindent
{\em $\bullet$ Logarithmic series. }We determine 
Maclaurin\index{Maclaurin, Colin} series of two functions based
on logarithm.

\begin{thm}[Maclaurin series of logarithms]\label{thm_TaylSerLog}
Let\index{theorem!Maclaurin series of logarithms|emph} 
$f(x)\equiv\log(1+x)$ and $g(x)\equiv\log(\frac{1}{1-x})$. \underline{Then}
$$
\text{$T^{f,\,0}(x)=\sum_{n=1}^{\infty}\frac{(-1)^{n-1}}{n}\cdot x^n$ and $T^{f,\,0}(x)=\log(1+x)$ $\iff$ $x\in(-1,1]$}\,.
$$
Also, 
$$
\text{$T^{g,\,0}(x)=\sum_{n=1}^{\infty}\frac{1}{n}\cdot x^n$ and $T^{g,\,0}(x)={\textstyle\log\big(\frac{1}{1-x}\big)}$ $\iff$ $x\in[-1,1)$}\,.
$$
\end{thm}
\duk
It suffices to consider just the function $f(x)$ because the result for $g(x)$ follows from the identity
$$
{\textstyle
\log\big(\frac{1}{1-x}\big)=
-\log(1+(-x))\ \ (x\in[-1,\,1))\,.
}
$$
The form of $T^{f,0}(x)$ follows from Proposition~\ref{prop_logs}. For $x\in(-1,1)$, $x\ne0$, we use Corollary~\ref{cor_LagCauRem} with $b=0$, $f(t)=\log(1+t)$ and the Cauchy remainder ($n\in\N$): with some $c=c(n)$ between $0$ and $x$ we have
$$
{\textstyle
\big|R^{f,\,0}_n(x)\big|=\frac{(n-1)!}{(1+c)^{n+1}}\cdot\frac{1}{n!}\cdot|x-c|^n\cdot|x|\to0,\ n\to\infty\,,
}
$$
because we know from the proof of Theorem~\ref{thm_NewtonBin} that $0<\frac{|x-c|}{1+c}\le|x|<1$. 

For $x=1$ we use the same bound: with some $c=c(n)\in(0,1)$, $n\in\N$, we have again
$$
{\textstyle
\big|R^{f,\,0}_n(1)\big|=\frac{1}{(1+c)^{n+1}}\cdot\frac{1}{n}\cdot(1-c)^n\cdot1<n^{-1}\to0,\ n\to\infty\,.}
$$
Finally, for $x>1$ the series obviously diverges.
\kduk
\vspace{-3mm}
\begin{exer}\label{ex_TaylLogAbscon}
Which of these series are abscon?    
\end{exer}

Now we easily prove the first summation in Theorem~\ref{thm_twoSums}.

\begin{cor}[alternating harmonic series]\label{cor_alterHarm}
We have
$$
\sum_{n=1}^{\infty}\frac{(-1)^{n-1}}{n}=\log 2\,.
$$
\end{cor}
\duk
It is the instance $x=1$ of the Maclaurin series for $\log(1+x)$.
\kduk

\noindent
{\em $\bullet$ Series for inverse trigonometric functions. }We determine 
Maclaurin\index{Maclaurin, Colin} series of the functions $\arctan x$ and
$\arcsin x$. Their derivatives get more and more complicated. Therefore
we use other techniques instead of Taylor remainders to determine $E(T^{\arctan x,0}(x))$ and $E(T^{\arcsin x,0}(x))$.

\begin{thm}[Two more Maclaurin series]\label{thm_TaylSerArctan}
Let\index{theorem!Maclaurin series of arkus tangent and arkus sinus|emph}  
$f(x)\equiv\arctan x$ and $g(x)\equiv\arcsin x$. \underline{Then}
$$
\text{$T^{f,\,0}(x)=\sum_{n=1}^{\infty}\frac{(-1)^{n-1}}{2n-1}\cdot x^{2n-1}$ and $T^{f,\,0}(x)=\arctan x$ $\iff$ $x\in[-1,1]$}\,.
$$
Also,
$$
\text{$T^{g,\,0}(x)=\sum_{n=1}^{\infty}(-1)^{n-1}\frac{\binom{-\frac{1}{2}}{n-1}}{2n-1}\cdot x^{2n-1}$ and $T^{g,\,0}(x)=\arcsin x$ $\iff$ $x\in[-1,1]$}\,.
$$
\end{thm}
\duk
Both Maclaurin series follow from Proposition~\ref{prop_arcs}. It is easy 
to see that both $T^{f,0}(x)$ and $T^{g,0}(x)$ absolutely converge on $(-1,1)$. Since 
$${\textstyle
f'(x)=\frac{1}{1+x^2}=\sum_{n=0}^{\infty}(-1)^nx^{2n}
}
$$
for every $x\in(-1,1)$, using results of the next section we see that the derivative 
$$\text{$\big(f(x)-T^{f,\,0}(x)\big)'=0$
on $(-1,1)$}\,.
$$
Thus $(f(x)\,|(-1,1)=T^{f,0}(x)+k_c(x)\,|(-1,1)$ for some $c$. 
Since $f(0)=0=T^{f,0}(0)$, we get that $c=0$ and $(-1,1)\sus E(T^{f,0}(x))$. A~similar 
argument using that
$${\textstyle
g'(x)=\big(1-x^2\big)^{-1/2}=\sum_{n=0}^{\infty}(-1)^n\binom{-1/2}{n}x^{2n}
}
$$
for every $x\in(-1,1)$ shows that $(-1,1)\sus E(T^{g,0}(x))$. It is easy to see that 
$$
E(T^{f,\,0}(x)),\ E(T^{g,\,0}(x))\sus[-1,\,1]
$$
because for $|x|>1$ both Maclaurin series diverge (also, $M(\arcsin x)=[-1,1]$). Both points $\pm1$ belong to both sets
$E(T^{f,0}(x))$ and $E(T^{g,0}(x))$ due to Abel's Theorem~\ref{thm_Abel}, see 
Corollary~\ref{cor_Abel}.
\kduk

\medskip\noindent
{\em $\bullet$ A~function $f\in\mathcal{C}^{\infty}(\R)$ with $E(T^{f,0}(x))=\{0\}$. }We give an example of such function. Recall that 
if $f\in\mathcal{C}^{\infty}(\R)$ then by Theorem~\ref{thm_classTayl} for every $b\in\R$ the Taylor series
$$
T^{f,\,b}(x)
$$
exists. Trivially, always $b\in E(T^{f,b}(x))$. We have seen several 
examples of functions when the set $E(T^{f,0}(x))$ is a~nontrivial
interval. Now we define a~function
for which this set is just $\{0\}$.

\begin{thm}[minimal $E(T^{f,0}(x))$]\label{thm_flatGraph}
Let\index{theorem!minimal $E(T^{f,0}(x))$|emph} 
$F\in\mathcal{F}(\R)$ be given by
$${\textstyle
F(0)\equiv0\,\text{ and }\,F(x)\equiv\exp\left(-\frac{1}{x^2}\right)\,\ds\,x\ne0\,.
}
$$
\underline{Then}
$${\textstyle
\text{$F\in\mathcal{C}^{\infty}(\R)$, $T^{F,\,0}(x)=\sum_{n=0}^{\infty}0x^n$ and $E(T^{F,\,0}(x))=\{0\}$}\,.
}
$$
\end{thm}
\duk
Let $\R_0\equiv\R\setminus\{0\}$ and
$${\textstyle
G(x)\equiv F(x)\,|\,\R_0=\exp\big(-\frac{1}{x^2}\big)\,. 
}
$$
We prove by induction on $n\in\N_0$ that $G^{(n)}\in\mathcal{F}(\R_0)$. We actually prove by induction that for every $n\in\N_0$,
$${\textstyle
G^{(n)}(x)=\exp(-\frac{1}{x^2})\sum_{j=0}^{k_n}a_{n,\,j}x^{-j}
\ \ (\in\mathcal{F}(\R_0))\,,
}
$$
where $k_n\in\N_0$ and $a_{n,j}\in\Z$. For $n=0$ this holds with $k_0=0$ and 
$a_{n,0}=1$. Then, by the results on derivatives in Chapter~\ref{chap_pr7}, 
\begin{eqnarray*}
G^{(n+1)}(x)=(G^{(n)}(x))'&=&{\textstyle
\exp(-\frac{1}{x^2})2x^{-3}\sum_{j=0}^{k_n}a_{n,\,j}x^{-j}\,+}\\
&&{\textstyle +\,
\exp(-\frac{1}{x^2})\sum_{j=0}^{k_n}(-j)a_{n,\,j}x^{-j-1}}\\
&=&{\textstyle \exp(-\frac{1}{x^2})\sum_{j=0}^{k_{n+1}}a_{n+1,\,j}
x^{-j}\,,
}
\end{eqnarray*}
where $k_{n+1}\in\N_0$ and $a_{n+1,j}\in\Z$. It is clear that for every $n\in\N_0$ and every $x\in\R_0$,
$$
F^{(n)}(x)=G^{(n)}(x)\,.
$$
We show by induction on $n\in\N_0$ that the derivative 
$F^{(n)}(0)$ exists and is $0$. For $n=0$
this follows from the definition of $F(x)$. Suppose that it holds for an $n\in\N_0$. By Exercise~\ref{ex_jednaLim}, the above form of $G^{(n)}(x)$ and Theorem~\ref{thm_AritLimFce}, 
$$
F^{(n+1)}(0)=\lim_{x\to0}
{\textstyle
\frac{F^{(n)}(x)-F^{(n)}(0)}{x}
}=\lim_{x\to0}x^{-1}G^{(n)}(x)=0\,.
$$
We see that $F\in\mathcal{C}^{\infty}(\R)$ and $F^{(n)}(0)=0$ for every 
$n\in\N_0$. By Theorem~\ref{thm_classTayl},
$$
{\textstyle
T^{F,\,0}(x)=0+0x+0x^2+\cdots\,.
}
$$
The sum of this series is for $x\in\R$ the constant zero function $k_0(x)$.
But $k_0(x)\ne F(x)$ for every $x\ne0$, and therefore $E(T^{F,0}(x))$ is just $\{0\}$.
\kduk
\vspace{-3mm}
\begin{exer}\label{ex_jednaLim}
For every $n\in\N_0$ it is true that
$${\textstyle
\lim_{x\to0}x^{-n}\exp(-\frac{1}{x^2})=0\,.
}
$$
\end{exer}

\section[Formal power series]{${}^*$Formal power series}\label{sec_fps}

In this section we develop the theory of formal power series to 
the extent sufficient for obtaining in Theorem~\ref{thm_asymOrPar} 
rough asymptotics of the numbers 
$\mathrm{op}_n$ of nonempty ordered partitions of $[n]$ (Definition~\ref{def_numOrdPar}).

\medskip\noindent
{\em $\bullet$ Ordered partitions. }Let $k\in\N$ and $n\in\N_0$. An \underline{ordered partition\index{ordered partition|emph}} of $[n]$ with $k$ parts is any $k$-tuple
$$
\overline{A}=\langle A_1,\,A_2,\,\ds,\,A_k\rangle
$$
of (possibly empty) mutually disjoint sets $A_i$ such that $\bigcup_{i=1}^k A_i=[n]$. Recall that $[n]=\{1,2,\ds,n\}$ and $[0]=\emptyset$. We define 
$$
\mathrm{OP}(k,\,n)\equiv\{\text{ordered partitions of $[n]$ with $k$ parts}\}\,\text{ and }\,\mathrm{op}_{k,\,n}\equiv|\mathrm{OP}(k,\,n)|\,.\label{OPkn}\label{opkn}
$$
For example, $\mathrm{op}_{2,2}=4$ because
$$
\mathrm{OP}(2,\,2)=\{\langle\emptyset,\,
[2]\rangle,\,\langle[2],\,\emptyset
\rangle,\,\langle\{1\},\,\{2\}\rangle,\,
\langle\{2\},\,\{1\}\rangle\}\,.
$$
A~\underline{nonempty ordered partition\index{ordered 
partition!nonempty|emph}} $\overline{A}$ has 
all parts $A_i\ne\emptyset$.

We are interested in the asymptotics of the sequence $(\mathrm{op}_n)$ of the 
following counting numbers. 

\begin{defi}[$\mathrm{op}_n$ and $\mathrm{OP}(_)n$]\label{def_numOrdPar}
Let $n\in\N$. We denote by $\mathrm{op}_n$\label{opn} the number of nonempty ordered 
partitions of $[n]$ with any number of parts. We set $\mathrm{op}_0\equiv1$. We denote 
the set of nonempty ordered partitions of $[n]$ by $\mathrm{OP}(n)$ 
($n\in\N$).\label{OPn} Thus $\mathrm{op}_n=|\mathrm{OP}(n)|$.
\end{defi}
For example, $\mathrm{op}_2=3$ because the nonempty ordered partitions of $[2]$ are  
$$
\langle[2]\rangle,\,
\langle\{1\},\,\{2\}\rangle\,\text{ and }\,\langle\{2\},\,\{1\}\rangle\,.
$$

\begin{defi}[multinomial coefficients]\label{def_multCoef}
We suppose that $k\in\N$ and that $n,n_i\in\N_0$ for $i\in[k]$ are such that 
$n_1+n_2+\ds+n_k=n$. We define the \underline{multinomial 
coefficient\index{multinomial coefficients|emph}} $\binom{n}{n_1,n_2,\ds,n_k}$ by
$$
\binom{n}{n_1,\,n_2,\,\ds,\,n_k}\equiv
\frac{n!}{n_1!\cdot n_2!\cdot\ldots\cdot n_k!}\,.
$$
\end{defi}
Thus $\binom{n}{k}=\binom{n}{k,n-k}$ for every $n,k\in\N_0$.

\begin{exer}\label{ex_prNext}
 Prove the next proposition.   
\end{exer}

\begin{prop}[counting OP]\label{prop_countOP}
With the notation of the previous definition we have
$${\textstyle
|\{\langle A_1,\,A_2,\,\ds,\,A_k
\rangle\in\mathrm{OP}(k,\,n)\cc\;|A_i|=n_i,\;i\in[k]\}|=
\binom{n}{n_1,\,n_2,\,\ds,\,n_k}\,.
}
 $$
\end{prop}

For any set $X$ and $n\in\N_0$ we introduce the notation
$$
{\textstyle
\binom{X}{n}\equiv\{A\cc\;A\sus X\wedge|A|=n\}\,.\label{Xnadn}
}
$$
Thus $\big|\binom{[n]}{k}\big|=\binom{n}{k}$ for every $n,k\in\N_0$. 

\begin{exer}\label{ex_coJe}
What is $\binom{X}{0}$?    
\end{exer}
The following exercise shows that the sequence
$$
{\textstyle
\big(\frac{1}{n!}\cdot\mathrm{op}_n\cc\;n\in\N\big)
}
$$
is supermultiplicative.
\begin{exer}\label{ex_submOrPar}
Let $m,n\in\N$. Define an injection
$$
{\textstyle
\iota\cc\mathrm{OP}(m)\times\mathrm{OP}(n)\times\binom{[m+n]}{m}
}\to\mathrm{OP}(m+n)\,.
$$
\end{exer}

\noindent
{\em $\bullet$ Formal power series. }These are just real sequences, but 
indexed by $\N_0$ instead of $\N$. They are endowed with arithmetic 
operations. 

\begin{defi}[formal power series]\label{def_fps}
We
write
$$
\R[[x]]\equiv
\{A\cc\;\text{$A$ is a~map from $\N_0$ to $\R$}\}\label{fps}\label{Rx}
$$
for the set of \underline{formal power series\index{formal power series, fps|emph}}, or {\em fps},
with one formal variable $x$ and real coefficients $a_n\equiv A(n)$, $n\in\N_0$. We write 
the elements $A=A(x)$ in $\R[[x]]$ as formal sums
$$
{\textstyle
\sum_{n\ge0}a_nx^n\,.
}
$$
Instead of the \underline{constant term\index{formal power series, fps!constant term|emph}} $a_0x^0$ we usually write just $a_0$, and instead of $a_1x^1$ we write just $a_1x$.
\end{defi}
We define the sum and product of two fps as
$$
{\textstyle
\sum_{n\ge0}a_nx^n+\sum_{n\ge0}b_nx^n
\equiv\sum_{n\ge0}(a_n+b_n)x^n
}\index{formal power series, fps!addition of|emph}
$$
and
$$
{\textstyle
\sum_{n\ge0}a_nx^n\cdot\sum_{n\ge0}b_nx^n
\equiv\sum_{n\ge0}\big(\sum_{k=0}^n a_kb_{n-k}\big)x^n\,.
}\index{formal power series, fps!product of|emph}
$$
The neutral elements in $\R[[x]]$ are
$$
0=0_{\R[[x]]}\equiv0x^0+0x^1+\ds\,\text{ and }\,1=1_{\R[[x]]}\equiv1x^0+0x^1+0x^2+\ds\,.\index{formal power series, fps!neutral elements in|emph}
$$ 

\begin{exer}\label{ex_isAbel}
In every (commutative unital) ring $R_{\mathrm{ri}}$ the structure
$$
\langle R^{\times},\,1_R,\,\cdot\rangle
$$
is an Abelian group.
\end{exer}

\begin{prop}[ring of fps]\label{prop_fpsring}
The structure
$$
\R[[x]]_{\mathrm{ri}}\equiv
\langle\R[[x]],\,0_{\R[[x]]},\,1_{\R[[x]]},\,+,\,\cdot\rangle\label{Rxri}
$$
is a~ring. Its units are exactly the {\em fps} with nonzero constant terms.
\end{prop}
\duk
Commutativity and associativity of addition and 
multiplication of fps are left for Exercise~\ref{ex_triviCheck}. 
Neutrality of $0_{\R[[x]]}$ and $1_{\R[[x]]}$ is easy to check. We treat in detail the distributive law, additive inverses and units.

In $\R[[x]]_{\mathrm{ri}}$ the distributive law holds: 
\begin{eqnarray*}
&&{\textstyle
\sum_{n\ge0}a_nx^n\cdot
\big(\sum_{n\ge0}b_nx^n+\sum_{n\ge0}c_nx^n\big)=}\\
&&{\textstyle
=\sum_{n\ge0}a_nx^n\cdot\sum_{n\ge0}(b_n+c_n)x^n=}\\
&&{\textstyle
=\sum_{n\ge0}\sum_{j=0}^n a_j(b_{n-j}+c_{n-j})x^n=}\\
&&{\textstyle
=\sum_{n\ge0}\sum_{j=0}^n(a_jb_{n-j}+a_jc_{n-j})x^n=\ds=}\\
&&{\textstyle
=\sum_{n\ge0}a_nx^n\cdot\sum_{n\ge0}b_nx^n+
\sum_{n\ge0}a_nx^n\cdot\sum_{n\ge0}c_nx^n}\,.
\end{eqnarray*}
The additive inverse of $\sum_{n\ge0}a_nx^n$ is $\sum_{n\ge0}(-a_n)x^n$. If $A(x)$ is a~fps that has 
zero constant term, then so has every multiple of it and $A(x)$ is not a~unit.
Let $A(x)=\sum_{n\ge0}a_nx^n$ be a~fps with $a_0\ne0$. We show that there is 
a~fps $\sum_{n\ge0}b_nx^n$ such that
$$
{\textstyle
\sum_{n\ge0}a_nx^n\cdot\sum_{n\ge0}b_nx^n=a_0b_0+(a_0b_1+a_1b_0)x+\ds=1\,.
}
$$
This is equivalent with the system
$$
{\textstyle
\sum_{j=0}^na_jb_{n-j}=c_n,\ n\in\N_0\,,}
$$
where $b_n$ are unknowns, $c_0=1$ and $c_n=0$ for $n>0$. The solution starts with $b_0\equiv a_0^{-1}$, determined by the $0$-th equation. If 
$b_0$, $b_1$, $\ds$, $b_n$ are already determined, we compute $b_{n+1}$ from the $(n+1)$-st equation:
$${\textstyle
b_{n+1}\equiv-a_0^{-1}\sum_{j=0}^na_{n+1-j}b_j\,.
}
$$
Hence $A(x)$ is a~unit.
\kduk
\vspace{-3mm}
\begin{exer}\label{ex_triviCheck}
 Show that the operations $+$ and $\cdot$ on $\R[[x]]$ are commutative and associative.   
\end{exer}

\noindent
{\em $\bullet$ Formal exponential. }Formal exponential series is as 
important as its functional version $\mathrm{e}^x$.

\begin{defi}[formal exponential]\label{def_formExp}
The {\em fps}
$$
\mathrm{e}^x_{\mathrm{fps}}\equiv
\sum_{n\ge0}\frac{1}{n!}\cdot x^n
$$
is the \underline{formal exponential\index{formal exponential|emph}}.
\end{defi}

\begin{exer}\label{ex_form ExpIden}
Show that in the ring $\R[[x,y]]_{\mathrm{ri}}$ of {\em fps} with two formal variables $x$ 
and $y$  the identity
$$
\mathrm{e}^{x+y}_{\mathrm{fps}}=
\mathrm{e}^x_{\mathrm{fps}}\cdot
\mathrm{e}^y_{\mathrm{fps}}
$$
holds.
\end{exer}

\noindent
{\em $\bullet$ Infinite series in $\R[[x]]_{\mathrm{ri}}$. }We review limit
transitions and sums of infinite series in $\R[[x]]_{\mathrm{ri}}$. For a~fps
$A(x)=\sum_{n\ge0}a_nx^n$ and $m\in\N_0$ we introduce the notation
$$
[x^m]\,A(x)\equiv a_m\label{extrCoef}
$$
for the coefficient of $x^m$ in $A(x)$.

\begin{defi}[formal limit]\label{def_formConv}
Let\index{formal limit|emph} 
$(A_n(x))\sus\R[[x]]$ and let $A(x)$ be a~fps. If for every 
$m\in\N_0$ there is an $n_0=n_0(m)$ in $\N_0$ such that for every $n\ge n_0$ we have
$$
[x^m]\,A_n(x)=[x^m]\,A(x)\,,
$$
we say that the sequence of fps $(A_n(x))$ has the \underline{formal limit} $A(x)$ and we write $\mathrm{flim}\,A_n(x)=A(x)$ or 
$\mathrm{flim}_{n\to\infty}\,A_n(x)=A(x)$.\label{flim}
\end{defi}

\begin{exer}\label{ex_uniFlim}
Formal limits are unique.     
\end{exer}

For a~series of fps $\sum_{n=1}^{\infty}A_n(x)$ we define 
its sum as the formal limit 
$${\textstyle
\mathrm{flim}_{n\to\infty}
\sum_{j=1}^n A_j(x)\,,
}
$$
if this formal limit of partial sums exists. We denote the sum, which is a~fps 
$A(x)$, again by $\sum_{n=1}^{\infty}A_n(x)$. A~very 
convenient property of formal convergence is that the necessary
condition of convergence (Proposition~\ref{prop_nutnaPoKon}) is also sufficient.

\begin{thm}[NCC becomes SCC]\label{thm_NCCscc}\label{SCC}
A~series\index{theorem!NCC becomes SCC|emph} $\sum_{n=1}^{\infty}A_n(x)\sus\R[[x]]$ has a~sum $\iff$
$$
\mathrm{flim}_{n\to\infty}A_n(x)=0_{\R[[x]]}\,.
$$
\end{thm}
\duk
Let ($n\in\N$)
$$
{\textstyle
A_n(x)=\sum_{m\ge0}a_{m,\,n}x^m\,.
}
$$
We prove implication $\Rightarrow$.   
Let $m\in\N_0$. By the assumption the sequence
$$
{\textstyle
\big(\sum_{j=1}^n a_{m,\,j}\cc\;n\in\N\big)
}
$$
is eventually constant. Thus the sequence $(a_{m,1},a_{m,2},\ds)$ is eventually zero and $\mathrm{flim}\,A_n(x)=0$.

We prove implication $\Leftarrow$. Let $m\in\N_0$. Since the sequence 
$(a_{m,1},a_{m,2},\ds)$ is eventually zero, the above displayed sequence of partial sums eventually constantly equals $a_m$. 
Thus
$$
\mathrm{flim}_{n\to\infty}{\textstyle
\sum_{j=1}^n A_j(x)=\sum_{m\ge0}a_mx^m\,.
}
$$
\kduk

An important application of sums of series of fps is the formula for sums of formal geometric series.

\begin{prop}[formal geometric series]\label{prop_formGeom}
Suppose that the fps $A(x)$ has zero constant term. \underline{Then} the inverse to
the unit $1-A(x)$ equals to the sum of 
the formal geometric series with the quotient $A(x)$, 
$$
\big(1_{\R[[x]]}-A(x)\big)^{-1}=\sum_{n=0}^{\infty}A(x)^n=1_{\R[[x]]}+\sum_{n=1}^{\infty}A(x)^n\ \ (\in\R[[x]])\,.
$$
\end{prop}
\duk
It is easy to see that
$$
\mathrm{flim}\,A(x)^n=0\,.
$$
Therefore by the previous theorem
the formal geometric series has a~sum
$B(x)\equiv1+\sum_{n=1}^{\infty}A(x)^n$ ($\in\R[[x]]$). We have 
\begin{eqnarray*}
(1-A(x))\cdot B(x)&=&{\textstyle
1-A(x)+\sum_{n=1}^{\infty}A(x)^n-A(x)\cdot\sum_{n=1}^{\infty}A(x)^n}\\
&=&{\textstyle\sum_{n=0}^{\infty}A(x)^n-\sum_{n=1}^{\infty}A(x)^n}\\
&=&A(x)^0=1\
\end{eqnarray*}
(Exercise~\ref{ex_Justify}). Thus $B(x)=(1-A(x))^{-1}$
\kduk
\vspace{-3mm}
\begin{exer}\label{ex_Justify}
Justify the previous computation in detail.    
\end{exer}

\noindent
{\em $\bullet$ Exponential generating functions. }Back to combinatorics. 
A~useful family of fps are \underline{exponential generating 
functions\index{exponential generating functions, EGF|emph}}, or 
EGF.\label{EGF} These are fps of the form
$$
\sum_{n\ge0}\frac{a_n}{n!}\cdot x^n
$$
where $a_n\in\R$. Often even $a_n\in\N_0$ and the coefficients $a_n$ count some combinatorial structures of size $n$.

\begin{thm}[product formula for EGF]\label{thm_prodForm}
Let\index{theorem!product formula for EGF|emph} 
$k\in\N$. Suppose that the $k+1$ real sequences $(a_n)$ and $(b_{n,i})$ with 
$n\in\N_0$ and $i\in[k]$ are for every $n\in\N_0$ bound by the 
convolution
$$
a_n=\sum_{\substack{\overline{A}\in\mathrm{OP}(k,\,n)\\\overline{A}=\langle A_1,\,\ds,\,A_k\rangle}}
\prod_{i=1}^k b_{|A_i|,\,i}\,.
$$
\underline{Then} the corresponding {\em EGF} are related by the product
$$
\sum_{n\ge0}\frac{a_n}{n!}\cdot x^n=\prod_{i=1}^k\sum_{n\ge0}\frac{b_{n,\,i}}{n!}\cdot x^n\,.
$$
\end{thm}
\duk
This follows from the formula in Proposition~\ref{prop_countOP}:
\begin{eqnarray*}
{\textstyle
\sum_{n\ge0}\frac{a_n}{n!}}\cdot x^n&=&{\textstyle
\sum_{n\ge0}\frac{1}{n!}\big(
\sum_{\substack{\overline{A}\in\mathrm{OP}(k,\,n)\\\overline{A}=\langle A_1,\,\ds,\,A_k\rangle}}
\prod_{i=1}^k b_{|A_i|,\,i}\big)\cdot x^n}\\
&=&{\textstyle
\sum_{n\ge0}\frac{1}{n!}\big(
\sum_{\substack{n_1,\,\ds,\,n_k\in\N_0\\
n_1+\ds+n_k=n}}
\sum_{\substack{\overline{A}\in\mathrm{OP}(k,\,n)\\\overline{A}=\langle A_1,\,\ds,\,A_k\rangle\\
|A_i|=n_i}}
\prod_{i=1}^k b_{|A_i|,\,i}\big)\cdot x^n}\\
&=&{\textstyle
\sum_{n\ge0}\frac{1}{n!}\big(
\sum_{\substack{n_1,\,\ds,\,n_k\in\N_0\\
n_1+\ds+n_k=n}}\binom{n}{n_1,\,\ds,\,n_k}
\prod_{i=1}^k b_{n_i,\,i}\big)\cdot x^n}\\
&=&{\textstyle
\sum_{\substack{n,\,n_1,\,\ds,\,n_k\in\N_0\\
n_1+\ds+n_k=n}}\prod_{i=1}^k
\frac{b_{n_i,\,i}}{n_i!}\cdot x^{n_i}=
\prod_{i=1}^k\sum_{n_i\ge0}\frac{b_{n_i,\,i}}{n_i!}\cdot x^{n_i}}\,.
\end{eqnarray*}
\kduk

\begin{cor}[EGF of ordered partitions]\label{cor_EGForPar}
In the ring $\R[[x]]_{\mathrm{ri}}$
of formal power series the identity
$$
\sum_{n\ge0}\frac{\mathrm{op}_n}{n!}\cdot x^n=\frac{1}{2-\mathrm{e}^x_{\mathrm{fps}}}
$$
holds. 
\end{cor}
\duk
For $k\in\N$ and $n\in\N_0$ we denote by $a_{k,n}$ the number of nonempty 
ordered partitions of $[n]$ with $k$ parts;  $a_{k,0}=0$ for every $k$. For 
$k\in\N$, $n\in\N_0$ and $i\in[k]$ we set $b_{k,n,i}\equiv1$ if $n>0$ and  
$b_{k,0,i}\equiv0$. Then 
$${\textstyle
a_{k,\,n}=\sum_{\substack{\overline{A}\in\mathrm{OP}(k,\,n)\\\overline{A}=\langle A_1,\,\ds,\,A_k\rangle}}
\prod_{i=1}^k b_{k,\,|A_i|,\,i}
}
$$
for every $k\in\N$ and every $n\in\N_0$. We therefore use the previous theorem, on the third line below, and get
\begin{eqnarray*}
{\textstyle
\sum_{n\ge0}\frac{\mathrm{op}_n}{n!}\cdot x^n}&=&{\textstyle 1+
\sum_{n\ge1}\frac{\mathrm{op}_n}{n!}\cdot x^n}\\
&=&{\textstyle
1+\sum_{n\ge1}
\sum_{k\ge1}\frac{a_{k,\,n}}{n!}\cdot x^n=1+\sum_{k\ge1}
\sum_{n\ge0}\frac{a_{k,\,n}}{n!}\cdot x^n}\\
&=&{\textstyle
1+\sum_{k\ge1}\prod_{i=1}^k\sum_{n\ge0}\frac{b_{k,\,n,\,i}}{n!}\cdot x^n=
1+\sum_{k\ge1}\big(\mathrm{e}^x_{\mathrm{fps}}-1\big)^k}\\
&=&{\textstyle
\frac{1}{1-(\mathrm{e}^x_{\mathrm{fps}}-1)}=\frac{1}{2-\mathrm{e}^x_{\mathrm{fps}}}\,.
}
\end{eqnarray*}
On the fourth line we used Proposition~\ref{prop_formGeom}.
\kduk
\vspace{-3mm}
\begin{exer}\label{ex_zduvTo}
Justify the computation in the previous proof in detail.    
\end{exer}

\noindent
{\em $\bullet$ Formal Dirichlet series. }For the next exercise we regard the zeta series 
$\zeta(s)=\sum_{n\ge1}\frac{1}{n^s}$ as a~formal Dirichlet series, with the formal
variable $s$. For $n\in\N$ with $n>1$, let $m_n$\label{ordFac} be the number of $k$-tuples ($k\in\N$)
$$
\langle n_1,\,n_2,\,\ds,\,n_k\rangle\in
(\N\setminus\{1\})^k
$$
such that $\prod_{i=1}^k n_i=n$, and let $m_1\equiv1$. For example, $m_6=3$ 
because $6$ can be factorized as $6$, $2\cdot 3$ and $3\cdot 2$. 

\begin{exer}\label{ex_ordFac}
Show that in the ring of formal Dirichlet series the identity 
$$
\sum_{n\ge1}\frac{m_n}{n^s}=\frac{1}{2-\zeta(s)}
$$
holds.
\end{exer}

\section[Real analytic functions]{${}^*$Real analytic functions}\label{sec_analFunc}

In this section, we investigate functions that can be expressed as sums of \underline{power series\index{power series|emph}} (with the \underline{center\index{power 
series!center of|emph}} $b\in\R$). These are familiar series 
$$
{\textstyle
\sum_{n=0}^{\infty}a_n(x-b)^n,\ a_n\in\R\,,
}
$$
which look like Taylor series.
The numbers $a_0$, $a_1$, $\ds$ are the \underline{coefficients\index{power 
series!coefficient of|emph}} of the power series and $x\in\R$. Coefficients 
of Taylor series come from Taylor polynomials, but now we allow $a_n$ to 
be arbitrary real numbers.

\medskip\noindent
{\em $\bullet$ Real analytic functions. }They are sums  of power 
series and are precisely defined as follows.

\begin{defi}[RAF]\label{def_RAF}
We say that a~function 
$f\in\mathcal{R}$ is \underline{real 
analytic\index{real analytic functions, RAF|emph}}\label{RAF} if the following holds.
\begin{enumerate}
\item The definition domain $M(f)$ is an open set. 
\item For every point $b\in M(f)$ there exist real numbers $\de$ and $a_0$, $a_1$, $\ds$ such that $U(b,\de)\sus M(f)$ and for every $x\in U(b,\de)$ we have
$${\textstyle
f(x)=\sum_{n=0}^{\infty}a_n(x-b)^n\,.
}
$$
\end{enumerate}
The set of real analytic functions is denoted by {\em RAF}.
\end{defi}
For example, later we show that the function $\exp(-1/x^2)$ with the definition domain
$$
(-\infty,\,0)\cup(0,\,+\infty)
$$
is real analytic.

\begin{defi}[radius of convergence]\label{def_radiConv}
By the \underline{radius of 
convergence\index{power series!radius of convergence of|emph}} of a~power 
series $S(x)=\sum_{n=0}^{\infty}a_n(x-b)^n$ we mean the element
$$
R=R(S(x))\equiv\frac{1}{\limsup_{n\to\infty}|a_n|^{1/n}}\ \ (\in[0,\,+\infty)\cup\{+\infty\})\,,\label{radiConv}
$$
where we set $\frac{1}{0}\equiv+\infty$.    
\end{defi}
For example,
$${\textstyle
R(T^{\log(1+x),\,0}(x))=R\big(\sum_{n=1}^{\infty}(-1)^{n-1}n^{-1}x^n)=1
}
$$
because $\lim n^{1/n}=1$. One can describe by the radius of convergence of a~power series the set 
of $x\in\R$ where it converges. In this description we get the intervals 
$(b,b)=\emptyset$ and $(-\infty,+\infty)\equiv\R$.

\begin{thm}[interval of convergence]\label{thm_intConv}
Let\index{theorem!interval of convergence|emph} $S(x)=\sum_{n=0}^{\infty}a_n(x-b)^n$
be a~power series with center $b$. We define the real set 
$${\textstyle
I=I(S(x))\equiv\{x\in\R\cc\;\text{the power series $S(x)$ converges}\}}\,.\label{interConv}
$$
\underline{Then} 
$$
\text{either $I=\{b\}$ or $I$ is a~nontrivial interval}\,.
$$
In the latter case either $I=\R$ or $I$ is a~bounded interval and $b$ is its midpoint. The interior of $I$ is
$$
I^0=(b-R(S(x)),\,b+R(S(x)))
$$
and for every $x\in I^0$ the power series $S(x)$ absolutely converges. We 
call 
$\text{$I$ the \underline{interval of convergence\index{power series!interval of convergence of|emph}} of $S(x)$}$. 
\end{thm}
\duk
Let $R$, where $R\ge0$ is real or $R=+\infty$, be the radius of 
convergence of the power series $S(x)$ and $I$ be its interval of convergence. If $R=0$ then for every large $c>0$ we have $|a_n|\ge c^n$ for infinitely many $n$. Thus for every $x\ne b$ we have
$$
|a_n(x-b)^n|\ge1
$$
for infinitely many $n$ and $S(x)$ diverges. So for $R=0$ we have
$$
I=\{b\}\,\text{ and }\,I^0=\emptyset=(b-R,\,b+R)=(b,\,b)\,.
$$
If $R=+\infty$ then for every small $c>0$ we have $|a_n|\le c^n$ for every large $n$.
Thus for every $x\in\R$ and every small $c>0$ we have
$$
|a_n(x-b)^n|\le c^n
$$
for every large $n$ and, by comparison with the geometric series, $S(x)$ 
absolutely converges. In particular, for $R=+\infty$ we have
$$
I=\R=(b-R,\,b+R)=(-\infty,\,+\infty)\,.
$$

We consider the remaining case when $0<R<+\infty$. Let 
$x\in\R$ be such that $r\equiv|x-b|<R$. We show that the power series $S(x)$ 
absolutely converges. We take a~small $\de$ such that $r+\de<R$. From the 
definition of $R$ we have that
$|a_n|\le(r+\de)^{-n}$ for every large $n$. Hence
$${\textstyle
|a_n(x-b)^n|\le\big(\frac{r}{r+\de}\big)^n
}
$$
for every large $n$ and $S(x)$ absolutely converges by comparison with
the geometric series. We finally show that if $r\equiv|x-b|>R$ then the power 
series $S(x)$ diverges. We take a~small $\de$ such that $r-\de>R$. From the 
definition of $R$ we have that
$|a_n|\ge(r-\de)^{-n}$ for infinitely many $n$. Hence 
$${\textstyle
|a_n(x-b)^n|\ge\big(\frac{r}{r-\de}\big)^n>1
}
$$
for infinitely many $n$ and $S(x)$ diverges. Hence for $0<R<+\infty$ we have
$$
(b-R,\,b+R)=I^0\sus I\sus[b-R,\,b+R]\,,
$$
$S(x)$ absolutely converges on the former interval and the proof of the 
theorem is complete.
\kduk

\noindent
As we saw in the case of the binomial series, it may not be easy to determine 
if the points $b\pm  R$ belong to $I$.

\begin{exer}\label{ex_shiftI}
Intervals of convergence of power series $A(x)=\sum_{n=0}^{\infty}a_nx^n$ 
and $B(x)=\sum_{n=0}^{\infty}a_n(x-b)^n$ are related by 
$$
I(B(x))=I(A(x))+b=\{a+b\cc\;a\in I(A(x))\}\,.
$$
\end{exer}

\begin{exer}\label{ex_relEandI}
Let $T^{f,b}(x)$ be a~Taylor series. 
What is the relation between the sets
$$
E(T^{f,\,b}(x))\,\text{ and }\,I(T^{f,\,b}(x))\,?
$$
\end{exer}

\begin{exer}\label{ex_fourPS}
Find four power series with respective intervals of convergence
$$
(-1,\,1),\ (-1,\,1],\ [-1,\,1),\,\text{ and }\,[-1,\,1]\,.
$$
\end{exer}

\begin{exer}\label{ex_PoloKonv}
Find the radii of convergence of the power series with center $0$ defining
the functions $\exp x$, $\cos x$ and $\sin x$.    
\end{exer}

\noindent
{\em $\bullet$ Sum functions, power series and fps. }Suppose that 
$A(x)=\sum_{n=0}^{\infty}a_n(x-b)^n$ is a~power series with center $b$. The 
\underline{sum function\index{sum function|emph}} of $A(x)$ is the 
function
$${\textstyle
F_A\cc I(A(x))\to\R,\ 
F_A(x)\equiv\sum_{n=0}^{\infty}a_n(x-b)^n\,.\label{fax}
}
$$
We also consider the function
$$
F_A^0\cc I(A(x))^0\to\R,\ 
F_A^0\equiv F_A\,|\,I(A(x))^0\,.\label{fax0}
$$
In the next passage we prove that $F_A^0$ is real analytic. Pringsheim's 
Theorem~\ref{thm_Pringsheim} below shows that in general real analyticity 
cannot be extended to endpoints of the interval of convergence.

Notation like
$$
{\textstyle
A(x)=\sum_{n\ge0}a_nx^n\,\text{ or }\,
A(x)=\sum_{n=0}^{\infty}a_nx^n\ \ 
(a_n\in\R)
}
$$
has two meanings. We can understand it as a~power series with 
center $0$, i.e. as sequences
$$
(a_0,\,a_1x,\,a_2x^2,\,\ds)\,\text{ with }\,x\in\R
$$
that define on $I(A(x))$ the sum function $F_A(x)$. 
Or we can understand it as a~fps (formal power series) in $\R[[x]]$, i.e. just as the sequence of real coefficients
$$
(a_0,\,a_1,\,a_2,\,\ds)\,.
$$
We work simultaneously with both interpretations. For example, we speak 
of the radius of convergence of a~fps and so on.

\begin{exer}\label{ex_jesteoUlo}
Explain why a~fps $\sum_{n\ge0}a_nx^n$
is more basic than the corresponding sum function $F_A(x)=
\sum_{n=0}^{\infty}a_n(x-b)^n$, $x\in I$.  
\end{exer}

\noindent
{\em $\bullet$ Power series and real analytic functions. }We expect that
the $F_A^0$ sum function of a~power series is real analytic and we prove it.

\begin{thm}[real analyticity of $F_A^0$]\label{thm_1stIntConv}
Suppose that\index{theorem!real analyticity of $F_A^0$|emph} a~power series
$${\textstyle
A(x)=\sum_{n=0}^{\infty}a_n(x-b)^n
}
$$ 
has a~nontrivial interval of convergence $I$. \underline{Then} the sum function
$$
F_A^0\cc I^0\to\R,\ F_A^0(x)\equiv\sum_{n=0}^{\infty}a_n(x-b)^n\,,
$$
defined on the interior of $I$ is real analytic. 
\end{thm}
\duk
Let $c\in I^0$, $c\ne b$, and let $\de$ be such that $U(c,\de)\sus I^0$. Then for every $x\in U(c,\de)$ we have
\begin{eqnarray*}
F_A^0(x)&=&{\textstyle
\sum_{n=0}^{\infty}a_n(x-b)^n=
\sum_{n=0}^{\infty}a_n(x-c+(c-b))^n}\\
&=&{\textstyle
\sum_{n=0}^{\infty}a_n\sum_{j=0}^n\binom{n}{j}(x-c)^j(c-b)^{n-j}}\\
&=&{\textstyle
\sum_{j=0}^{\infty}\big(\sum_{n,\,n\ge j}\binom{n}{j}a_n(c-b)^{n-j}\big)\cdot(x-c)^j}\\
&\equiv&{\textstyle
\sum_{j=0}^{\infty}b_j(x-c)^j,\,\text{ with }\,b_j\equiv\sum_{n,\,n\ge j}\binom{n}{j}a_n(c-b)^{n-j}\,.
}
\end{eqnarray*}
On the second line, we used the binomial theorem. It remains to justify 
the change of order of summation between the second and third line, and 
the convergence of series defining the new coefficients $b_j$. We leave the 
former for Exercise~\ref{ex_necoZdukazu} and turn 
to the latter. Since $|c-b|<R(S(x))$, there is a~constant 
$d\in(0,1)$ such that for every large $n$,
$$
\big|a_n(c-b)^n\big|\le d^n\,.
$$
Also, $\binom{n}{j}=O(n^j)$ ($n\ge j$). Thus for any $d_0\in(d,1)$ we have
$${\textstyle
\big|\binom{n}{j}(c-b)^{-j}\cdot 
a_n(c-b)^n\big|\le d_0^n}
$$
for every sufficiently large $n$ ($\ge j$). Comparison with the geometric 
series shows that the series defining
$b_j$ absolutely converges.
\kduk
\vspace{-3mm}
\begin{exer}\label{ex_necoZdukazu}
Justify the change of order of summation in the proof.     
\end{exer}

\noindent
{\em $\bullet$ Arithmetic operations  on real analytic functions. }We 
investigate the interplay between real analyticity and the operations of addition, 
multiplication and division on $\mathcal{R}$.

\begin{prop}[addition of power series]\label{prop_addRAF}
Suppose that 
$${\textstyle
\text{$A(x)=\sum_{n\ge0}a_nx^n$ and $B(x)=\sum_{n\ge0}b_nx^n$}\ \ (\in\R[[x]])
}
$$
are {\em fps} 
with positive radii of convergence $R_A$ and $R_B$, respectively.
Let $C(x)\equiv A(x)+B(x)=\sum_{n\ge0}(a_n+b_n)x^n$ as a~sum of {\em fps}.
\underline{Then}
$$
F_A(x)+F_B(x)=C(x)
$$
for every $x\in(-R,\,R)$ where $R=\min(R_A,R_B)$.
\end{prop}
\duk
Let $A(x)$, $B(x)$, $R_A$, $R_B$,  $C(x)$ and $R$ be as stated, and let $x\in(-R,R)$. By linear combination of series (Proposition~\ref{prop_linComSer}), 
$$
{\textstyle
F_A(x)+F_B(x)=
\sum_{n=0}^{\infty}a_nx^n+
\sum_{n=0}^{\infty}b_nx^n 
=\sum_{n=0}^{\infty}(a_n+b_n)x^n=C(x)\,.
}
$$
\kduk

\begin{thm}[addition in RAF]\label{thm_addRAF}
If $f,g\in\mathrm{RAF}$ \underline{then} 
$f+g\in\mathrm{RAF}$. 
\end{thm}
\duk
Suppose that $f,g$ are real analytic functions and $h\equiv f+g$. The domain $M(h)=M(f)\cap M(g)$ is an open set. Let 
$b\in M(h)$. Since $f,g\in\mathrm{RAF}$, there is a~$\de$ such that $U(b,\de)\sus 
M(h)$ and for every $x\in U(b,\de)$ we have
$$
{\textstyle
f(x)=\sum_{n=0}^{\infty}a_n(x-b)^n\,\text{ and }\,g(x)=\sum_{n=0}^{\infty}b_n(x-b)^n\,.
}
$$
By linear combination of series (Proposition~\ref{prop_linComSer}) we have for the same $x$ that
$${\textstyle
h(x)=f(x)+g(x)=\sum_{n=0}^{\infty}(a_n+b_n)(x-b)^n\,.
}
$$
Thus $h$ is real analytic.
\kduk

\begin{prop}[multiplying power series]\label{prop_multRAF}
Suppose that 
$${\textstyle
\text{$A(x)=\sum_{n\ge0}a_nx^n$ and $B(x)=\sum_{n\ge0}b_nx^n$}\ \ (\in\R[[x]])
}
$$
are fps 
with positive radii of convergence $R_A$ and $R_B$, respectively.
Let $C(x)\equiv A(x)\cdot B(x)=\sum_{n\ge0}\big(\sum_{j=0}^n a_jb_{n-j}\big)x^n$ as a~product of {\em fps}.
\underline{Then}
$$
F_A(x)\cdot F_B(x)=C(x)
$$
for every $x\in(-R,\,R)$ where $R=\min(R_A,R_B)$.
\end{prop}
\duk
Let $A(x)$, $B(x)$, $R_A$, $R_B$,  $C(x)$ and $R$ be as stated, and 
let $x\in(-R,R)$. By Theorem~\ref{thm_intConv} both series 
$${\textstyle
\text{$F_A(x)=\sum_{n=0}^{\infty}a_nx^n$ and $F_B(x)=\sum_{n=0}^{\infty}b_nx^n$}
}
$$
absolutely converge. By Theorem~\ref{thm_cauSouRad} we have
$$
{\textstyle
F_A(x)\cdot F_B(x)=\sum_{n=0}^{\infty}\big(
\sum_{j=0}^n a_jb_{n-j}\big)x^n=C(x)\,.
}
$$
\kduk

\begin{thm}[products in RAF]\label{thm_multRAF}
If $f,g\in\mathrm{RAF}$ \underline{then} 
$f\cdot g\in\mathrm{RAF}$. 
\end{thm}
\duk
Suppose that $f,g$ are real analytic functions and $h\equiv 
fg$. The domain $M(h)=M(f)\cap M(g)$ is an open set. Let 
$b\in M(h)$. Since $f,g\in\mathrm{RAF}$, there is a~$\de$ such that $U(b,\de)\sus 
M(h)$ and for every $x\in U(b,\de)$ we have
$$
{\textstyle
f(x)=\sum_{n=0}^{\infty}a_n(x-b)^n\,\text{ and }\,g(x)=\sum_{n=0}^{\infty}b_n(x-b)^n\,.
}
$$
Since these series absolutely converge, by 
Theorem~\ref{thm_cauSouRad} we have for the same $x$ that
$${\textstyle
h(x)=f(x)\cdot g(x)=\sum_{n=0}^{\infty}\big(\sum_{j=0}^na_jb_{n-j}\big)(x-b)^n\,.
}
$$
Thus $h$ is real analytic.
\kduk

\begin{prop}[reciprocals of power series~1]\label{prop_diviRAF}
Suppose that 
$${\textstyle
A(x)=\sum_{n\ge0}a_nx^n\ \ (\in\R[[x]])
}
$$
is a~fps with $a_0\ne0$ and positive radius of convergence $R_A$.
Let $B(x)\equiv A(x)^{-1}=\frac{1}{A(x)}=\sum_{n=0}^{\infty}b_nx^n$ be the {\em fps} multiplicative inverse of $A(x)$.
\underline{Then}
$$
\frac{1}{F_A(x)}=B(x)
$$
for every $x\in(-R,\,R)$ where $0<R\le R_A$.
\end{prop} 
\duk
Let $A(x)$, $R_A$ and $B(x)$ be as stated. More or less by 
Proposition~\ref{prop_formGeom} we 
have 
$${\textstyle
B(x)=\sum_{n=0}^{\infty}b_nx^n=a_0^{-1}\sum_{k\ge0}(-1)^k(a_1'x+a_2'x^2+\ds)^k \ \ (\in\R[[x]])
}
$$
where $a_j'=a_0^{-1}a_j$. Our task is to deduce from this 
representation of $B(x)$ that $B(x)$ has a~positive radius of 
convergence. Comparing the coefficients we have for every $n\in\N$ the formula
$$
{\textstyle
b_n=a_0^{-1}\sum_{k\ge1}(-1)^k\sum_{\substack{m_1,\,\ds,\,m_k\in\N\\m_1+\ds+m_k=n}}a_{m_1}'\ds a_{m_k}'\,.
}
$$
Since 
$$
\limsup_{n\to\infty}|a_n'|^{1/n}=
\limsup_{n\to\infty}|a_n|^{1/n}=
{\textstyle\frac{1}{R_A}<+\infty\,,
}
$$
there is a~real constant $c\ge\frac{1}{R_A}>0$ such that $|a_n'|\le c^n$ for every 
$n\in\N$. In fact, we may set
$$
c\equiv\sup\big(\{|a_n'|^{1/n}\cc\;n\in\N\}\big)\,.
$$
Hence we bound $b_n$ for every $n\in\N$ by
\begin{eqnarray*}
|b_n|&\le&{\textstyle
|a_0^{-1}|\sum_{k\ge1}\sum_{\substack{m_1,\,\ds,\,m_k\in\N\\
m_1+\ds+m_k=n}}c^{m_1}\ds c^{m_k}}\\
&=&{\textstyle
|a_0^{-1}|c^n\sum_{k=1}^n[x^n]\big(\frac{x}{1-x}\big)^k=
|a_0^{-1}|c^n\sum_{k=1}^n[x^{n-k}](1-x)^{-k}}\\
&=&{\textstyle
|a_0^{-1}|c^n\sum_{k=1}^n
(-1)^{n-k}\binom{-k}{n-k}=
|a_0^{-1}|c^n\sum_{k=1}^n
\binom{n-1}{n-k}}\\
&=&{\textstyle
|a_0^{-1}|c^n\sum_{j=0}^{n-1}
\binom{n-1}{j}\le|a_0^{-1}|(2c)^n\,.
}
\end{eqnarray*}
Thus $\limsup|b_n|^{1/n}\le2c$ and $B(x)$ has a~positive radius of convergence $R\le\frac{1}{2c}\le R_A$. 
\kduk

\noindent
We need this result for general center $b$.

\begin{cor}[reciprocals of power series~2]\label{cor_diviRAF}
Let $b\in\R$ and let
$${\textstyle
A(x)=\sum_{n\ge0}a_n(x-b)^n
}
$$
be a~power series with center $b$, $a_0\ne0$ and positive radius of convergence $R_A$.
Let $B(x)\equiv (\sum_{n\ge0}a_nx^n)^{-1}=\sum_{n=0}^{\infty}b_nx^n$ be the {\em fps} multiplicative inverse.
\underline{Then}
$$
\frac{1}{F_A(x)}=\sum_{n=0}^{\infty}b_n(x-b)^n
$$
for every $x\in(b-R,\,b+R)$ where $0<R\le R_A$.
\end{cor} 
\duk
This follows from the previous proposition by replacing $x$ appropriately with $x-b$.
\kduk
\vspace{-3mm}
\begin{exer}\label{ex_apprRepl}
Prove the corollary in detail.    
\end{exer}

\begin{thm}[division in RAF]\label{thm_diviRAF}
If $f,g\in\mathrm{RAF}$ \underline{then} 
$f/g\in\mathrm{RAF}$. 
\end{thm}
\duk
In view of Theorem~\ref{thm_multRAF}
and the identity
$${\textstyle
f/g=f\cdot\frac{1}{g}
}
$$
it suffices to show that if $g\in\mathrm{RAF}$ then also $\frac{1}{g}\in\mathrm{RAF}$.

Suppose that $g$ is a~real analytic function and $h\equiv 
\frac{1}{g}$. The domain $M(h)=M(g)\setminus Z(g)$ is an 
open set by the proof of Theorem~\ref{thm_deriSEF}. Let 
$b\in M(h)$. Since $g\in\mathrm{RAF}$, there is a~$\de$ such that $U(b,\de)\sus 
M(h)$ and for every $x\in U(b,\de)$ we have
$$
{\textstyle
g(x)=\sum_{n=0}^{\infty}b_n(x-b)^n\,.
}
$$
Since $b\not\in Z(g)$, we have $b_0\ne0$. Hence by 
Corollary~\ref{cor_diviRAF} there is a~$\theta\in(0,\de)$ and 
a~power series $C(x)$ with center $b$ such that
$$
{\textstyle
\frac{1}{g(x)}=C(x)
}
$$
for every $x\in U(b,\theta)$, as we need.
\kduk

\noindent
{\em $\bullet$ Composition in {\em RAF} and real analyticity of simple elementary functions. }

\begin{thm}[$\mathrm{SEF}\sus
\mathrm{RAF}$]\label{thm_SEFisRAF}
Every simple elementary function is real analytic.    
\end{thm}
\duk

\kduk

\noindent
{\em $\bullet$ Abel's theorem. }This classical theorem is due to the 
Norwegian mathematician {\em Niels H. Abel (1802--1829)\index{Abel, Niels 
H.}}. We state it in not completely standard way.

\begin{thm}[Abel's]\label{thm_Abel}
For every\index{theorem!Abel's|emph}
power series
$${\textstyle
A(x)=\sum_{n=0}^{\infty}a_n(x-b)^n
}
$$ 
the sum function $F_A\cc I(A(x))\to\R$ is continuous. 
\end{thm}
In the proof we use an inequality due also to N.~H.~Abel\index{Abel, Niels 
H.}.

\begin{prop}[Abel's inequality]\label{prop_abelIneq}
Let\index{Abel's inequality|emph} 
$n\in\N$ and $a_1$, $\ds$, $a_n$, $b_1$, $\ds$, $b_n$ be $2n$ real numbers such that $b_1\ge b_2\ge\ds\ge b_n\ge0$. \underline{Then}
$$
\bigg|\sum_{j=1}^n a_jb_j\bigg|\le
{\textstyle
\max\big(\big\{
\big|\sum_{j=1}^m a_j\big|\cc\;m\in[n]
\big\}\big)\cdot b_1\equiv A\cdot b_1\,.
}
$$
\end{prop}
\duk
For $m\in[n]$ we set $A_m\equiv\sum_{j=1}^m a_j$ and define $A_0=b_{n+1}\equiv0$. Using in the right moment the triangle inequality we get
\begin{eqnarray*}
{\textstyle
\big|\sum_{j=1}^n a_jb_j\big|}&=&
{\textstyle
\big|\sum_{j=1}^n (A_j-A_{j-1})b_j\big|=\big|
\sum_{j=1}^nA_j(b_j-b_{j+1})\big|}\\
&\le&{\textstyle
\sum_{j=1}^n|A_j|(b_j-b_{j+1})
\le A\cdot\sum_{j=1}^n(b_j-b_{j+1})}\\
&=&A\cdot b_1\,.
\end{eqnarray*}
\kduk

\noindent
{\bf Proof of Theorem~\ref{thm_Abel}. }We may assume that  the center $b=0$ 
and that $I\equiv I(A(x))\ne\{0\}$. Let $y\in I$, we want to show that 
$$
\lim_{x\to y}F_A(x)=F_A(y)\,. 
$$
We may assume that $y>0$.
First we show that 
$$
\lim_{x\to y^-}F_A(x)=F_A(y)\,. 
$$
Let an $\ep>0$ be given.
Since the series $\sum_{n=0}^{\infty}a_ny^n$ converges, we can take an $N\in\N$ such that
$$
{\textstyle
\big|\sum_{j=m}^n a_jy^j
\big|\le\frac{\ep}{2}
}
$$
whenever $n\ge m\ge N$. Then we take a~small $\de\in(0,y)$ such that for every $x\in(y-\de,y)$ we have
$$
{\textstyle
\big|\sum_{j=0}^N a_jy^j-\sum_{j=0}^N a_jx^j\big|\le\frac{\ep}{2}\,.
}
$$
Now for every $x\in(y-\de,y)$ and every $n>N$ we have by Abel's inequality, used on the third line below, that
\begin{eqnarray*}
&&{\textstyle
\big|\sum_{j=0}^n a_jy^j-\sum_{j=0}^n a_jx^j\big|}\\
&&\le{\textstyle
\big|\sum_{j=0}^N a_jy^j-\sum_{j=0}^N a_jx^j\big|+\big|\sum_{j=N+1}^n a_jy^j
-\sum_{j=N+1}^n a_jx^j\big|}\\
&&\le{\textstyle\frac{\ep}{2}+
\big|\sum_{j=N+1}^n a_jy^j\cdot\big(1-\big(\frac{x}{y}\big)^j\big)\big|
\le\frac{\ep}{2}+\frac{\ep}{2}\cdot1=\ep\,.
}
\end{eqnarray*}
Thus $|F_A(y)-F_A(x)|\le\ep$ for every $x\in(y-\de,y)$ and we get the 
left-sided limit. If $y=\max(I)$, we are done. If $y\ne\max(I)$ then $y\in I^0$ and we prove the above limit
by the argument used in the proof of Theorem~\ref{thm_spjMocRady1} 
(Exercise~\ref{ex_ProvInDe}).
\kduk

\noindent
Thus, unlike real analyticity, continuity of the sum function always 
extends to endpoints of the interval of convergence. Usually one understands by 
Abel's theorem only the left-sided continuity of $F_A$. There is 
a~version of Abel's theorem for complex power series.

\begin{exer}\label{ex_necoZdukazu2}
Justify the reductions at the beginning of the proof.      
\end{exer}

\begin{exer}\label{ex_ProvInDe}
Prove in detail that the sum function of a~power series is continuous at 
every interior point of the interval of convergence.    
\end{exer}

\begin{exer}\label{ex_onAbel}
Deduce by Abel's theorem the implication
$$
{\textstyle
\text{$\sum_{n=1}^{\infty}(-1)^{n-1}\frac{1}{n}x^n=\log(1+x)$ on $(1-\de,1)$ $\Rightarrow$ $\sum_{n=1}^{\infty}(-1)^{n-1}\frac{1}{n}=\log2$}\,.
}
$$
\end{exer}
Similarly we get the second summation in Theorem~\ref{thm_twoSums}.

\begin{cor}[two applications]\label{cor_Abel}
We have, with equal signs,
$$
\sum_{n=1}^{\infty}\frac{(-1)^{n-1}}{2n-1}(\pm1)^{2n-1}=\arctan(\pm1)=\pm\frac{\pi}{4}
$$
and 
$$
\sum_{n=1}^{\infty}(-1)^{n-1}\frac{\binom{-1/2}{n-1}}{2n-1}(\pm1)^{2n-1}=\arcsin(\pm1)=\pm\frac{\pi}{2}\,.
$$
\end{cor}
\duk
Since both functions are odd, it suffices to consider only the case with sign $+$. By Theorem~\ref{thm_TaylSerArctan},
$$
{\textstyle
\text{$\sum_{n=1}^{\infty}\frac{(-1)^{n-1}}{2n-1}x^{2n-1}=\arctan x$ on $(-1,1)$}\,.
}
$$
The series $\sum_{n=1}^{\infty}\frac{(-1)^{n-1}}{2n-1}$ converges by Theorem~\ref{thm_altSer}. Thus 
$$
{\textstyle
\sum_{n=1}^{\infty}\frac{(-1)^{n-1}}{2n-1}=\sum_{n=1}^{\infty}\frac{(-1)^{n-1}}{2n-1}1^n=\lim_{x\to1}\arctan x=\arctan 1=\frac{\pi}{4}
}
$$
by Abel's theorem. The proof of the second summation is very similar.
\kduk

\noindent
{\em $\bullet$ Pringsheim's theorem. }This theorem is 
due to the German mathematician {\em Alfred Pringsheim (1850--1941)\index{Pringsheim, Alfred}} who 
was the father in law of the German novelist {\em Thomas Mann 
(1875--1955)\index{Mann, Thomas}}. Pringsheim's theorem tells us that for 
power series with nonnegative coefficients real analyticity does not
extend to endpoints of the interval of convergence.

\begin{thm}[Pringsheim's]\label{thm_Pringsheim}
Suppose\index{theorem!Pringsheim's|emph} that a~power series
$${\textstyle
A(x)=\sum_{n=0}^{\infty}a_n(x-b)^n
}
$$ 
has all coefficients $a_n\ge0$ and has 
a~nontrivial and bounded interval of convergence $I$ with $I^0=(c,d)$. \underline{Then} there does not exist real analytic function  
$$
f(x)\cc (d-\de,\,d+\de)\to\R\,\text{ with }\,\de\in(0,\,d-c)
$$
such that 
$${\textstyle
f(x)=A(x)\,\text{ on }\,(d-\de,\,d)\,.
}
$$
\end{thm}
For the proof we need the next auxiliary proposition which brings us 
back to the beginning of the book.

\begin{prop}[nonnegative double sum]\label{prop_douSum}
Let 
$$
f\cc\N_0^2\to[0,\,+\infty)
$$
be an infinite table (matrix) with nonnegative real entries. \underline{Then} we always have the 
equality of double sums
$$
\sum_{m=0}^{\infty}\sum_{n=0}^{\infty}
f(m,\,n)=\sum_{n=0}^{\infty}\sum_{m=0}^{\infty}
f(m,\,n)\ \ (\in[0,\,+\infty)\cup\{+\infty\})\,.
$$
\end{prop}
\duk
You proved it already in Exercise~... .
\kduk

\noindent
The following proof of Pringsheim's theorem is taken from \cite[pp.~240--242]{flaj_sedg}.  

\medskip\noindent
{\bf Proof of Theorem~\ref{thm_Pringsheim}. }We may assume that the center $b=0$. Let $R$ 
($\in(0,+\infty)$) be the radius of convergence of $A(x)=\sum_{n=0}^{\infty}a_nx^n$. We 
assume for the contrary that there is a~$\de\in(0,R)$ and a~real analytic function
$$
f\cc(R-\de,\,R+\de)\to\R
$$
such that $A(x)=f(x)$ for every number $x\in(R-\de,R)$. Let $h\in(0,\frac{1}{3}\de)$ and $z\equiv R-h$ ($>0$). 
By the proof of Theorem~\ref{thm_1stIntConv}
we have the expansion $f(x)=\sum_{n=0}^{\infty}b_n(x-z)^n$ around $z$ with the coefficients 
$$
{\textstyle
b_n=\sum_{m=n}^{\infty}\binom{m}{n}a_mz^{m-n}\,.
}
$$
In particular, $b_n\ge0$ for every $n\in\N_0$. Below we use  Proposition~\ref{prop_douSum} on the third line and the binomial theorem on the fourth line and get the contradiction
\begin{eqnarray*}
f(R+h)&=&{\textstyle\sum_{n=0}^{\infty}
b_n(2h)^n}\\
&=&{\textstyle 
\sum_{n=0}^{\infty}\big(\sum_{m=n}^{\infty}\binom{m}{n}a_mz^{m-n}\big)(2h)^n}\\
&=&{\textstyle 
\sum_{n=0}^{\infty}\big(\sum_{m=0}^{\infty}\binom{m}{n}a_mz^{m-n}\big)(2h)^n}\\
&=&{\textstyle
\sum_{m=0}^{\infty}\sum_{n=0}^{\infty}\binom{m}{n}a_mz^{m-n}(2h)^n}\\
&=&{\textstyle
\sum_{m=0}^{\infty}\sum_{n=0}^m\binom{m}{n}a_mz^{m-n}(2h)^n}\\
&=&{\textstyle
\sum_{m=0}^{\infty}a_m(R+h)^m=A(R+h)
}
\end{eqnarray*}
because the number $R+h$ lies outside the interval $I(A(x))$.
\kduk

\noindent
For complex power series
Pringsheim's theorem plays an even more important role. 

\begin{exer}\label{ex_necoZdukazu3}
Justify the reduction at the beginning of the proof.    
\end{exer}

\noindent
{\em $\bullet$ The Pringsheim--Boas theorem, without proof. }A.~Pringsheim\index{Pringsheim, Alfred} 
had published another theorem  in \cite{prin} in 1893. Forty years later, the American
mathematician {\em Ralph P. Boas Jr. (1912--1992)\index{Boas Jr., Ralph P.}} 
discovered that Pringsheim's proof was fallacious and gave a~correct proof, see 
\cite{boas}. The theorem says that moderate growth of all Taylor coefficients suffices for a~function to be real analytic. For the sake of brevity we relegate the proof to {\em MA~1${}^+$}.

\begin{thm}[Pringsheim--Boas]\label{thm_PriBoa}
Let\index{theorem!Pringsheim--Boas|emph} 
$a<b$ be in $\R$ and $f\in\mathcal{C}^{\infty}((a,b))$. For $j\in\N_0$ and $t\in(a,b)$ we set
$a_j(t)\equiv\frac{1}{j!}f^{(j)}(t)$ and define
$$
\rho(t)\equiv\frac{1}{\limsup_{j\to\infty}|a_j(t)|^{1/j}}\,.
$$
If $\rho(t)\ge\de>0$ for every $t\in(a,b)$, \underline{then} $f$ is real analytic. 
\end{thm}

\begin{exer}\label{ex_necoZdukazu4}
    
\end{exer}

\section[Asymptotics of ordered partitions]{${}^*$Asymptotics of ordered partitions}

After the preparation in Section~\ref{sec_fps} (definition 
of $\mathrm{op}_n$ and the formula in Corollary~\ref{cor_EGForPar}) 
and in the previous section (Pringsheim's theorem) we finally
obtain the rough asymptotics of the numbers $\mathrm{op}_n$ 
of nonempty ordered partitions of $[n]$ (Definition~\ref{def_numOrdPar}).

\medskip\noindent
{\em $\bullet$ Asymptotics of $\mathrm{op}_n$. }It is as follows. 

\begin{thm}[rough asymptotics of $\mathrm{op}_n$]\label{thm_asymOrPar}
These numbers satisfy
$$
\lim_{n\to\infty}\Big(\frac{\mathrm{op}_n}{n!}\Big)^{1/n}=\frac{1}{\log 2}\,.
$$
Equivalently (Exercise~\ref{ex_jeToEkviv}),
$$
\mathrm{op}_n=(\log 2)^{-n+o(n)}\cdot n!\ \ (n\to\infty)\,.
$$
\end{thm}
We know from
Corollary~\ref{cor_EGForPar} the formula for the EGF of the numbers $\mathrm{op}_n$, but this is only a~formal relation in the ring 
$\R[[x]]_{\mathrm{ri}}$. We convert it to a~functional relation by means of 
the next proposition whose proof makes use of Pringsheim's theorem, 
Proposition~\ref{prop_diviRAF} and Theorem~\ref{thm_diviRAF}.

\begin{prop}[determining $R$]\label{prop_deteR}
Let $a>0$ and
$${\textstyle
A(x)=\sum_{n\ge0}a_nx^n\ \ (\in\R[[x]])
}
$$ 
be a~{\em fps} with $R(A(x))>a$
such that $F_A(x)\ne0$ on $(-a,a)$, hence $a_0\ne0$, but $F_A(a)=0$.
We further assume that the {\em fps} multiplicative inverse  
$${\textstyle
B(x)=A(x)^{-1}=\frac{1}{A(x)}=\sum_{n\ge0}b_nx^n
\ \ (\in\R[[x]])}
$$
has all coefficients $b_n\ge0$. 
\underline{Then} its radius of convergence $R(B(x))=a$.
\end{prop}
\duk
If $a_0=0$ then $F_A(0)=0$, hence $a_0\ne0$ and by 
Proposition~\ref{prop_fpsring} we have the fps inverse $B(x)$. By 
Proposition~\ref{prop_diviRAF} the fps $B(x)$ has positive radius of 
convergence $R_B$. Since the
finite limit
$$
\lim_{x\to a}{\textstyle
\frac{1}{F_A(x)}
}
$$
does not exist, Proposition~\ref{prop_diviRAF} 
implies that $0<R_B\le a$. Strict inequality $R_B<a$ would contradict, 
with the help of the real analytic function
$${\textstyle
f(x)=\frac{1}{F_A(x)}\,|\,U(R_B,\,a-R_B)
}
$$
obtained via Theorem~\ref{thm},
Pringsheim's Theorem~\ref{thm_Pringsheim}. Hence $R_B=a$.
\kduk

\noindent
{\bf Proof of Theorem~\ref{thm_asymOrPar}. }Using Corollary~\ref{cor_EGForPar} and the previous proposition with $a=\log 2$ and the sum function $F_A(x)=2-\mathrm{e}^x$ we get that the EGF
$$
{\textstyle
\sum_{n\ge0}\frac{\mathrm{op}_n}{n!}\cdot x^n}
$$
has the radius of convergence $R=\log 2$. Thus
$${\textstyle
\limsup_{n\to\infty}\big(\frac{\mathrm{op}_n}{n!}\big)^{1/n}=\frac{1}{\log 2}\,.
}
$$
We know by Exercise~\ref{ex_submOrPar} that the numbers 
$\frac{\mathrm{op}_n}{n!}$ are 
supermultiplicative. Using Corollary~\ref{cor_mulFeke} we deduce that
$${\textstyle
\lim_{n\to\infty}\big(\frac{\mathrm{op}_n}{n!}\big)^{1/n}=\frac{1}{\log 2}\,.
}
$$
\kduk
\vspace{-3mm}
\begin{exer}\label{ex_jeToEkviv}
The two statements of the theorem are equivalent.    
\end{exer}

\section[Arnol'd's limits]{${}^*$Arnol'd's limits}

In the booklet \cite[p.~21]{arno}, see \cite{arno_eng} for an English 
translation, the Russian mathematician {\em 
Vladimir~I.~Arnol'd (1937--2010)\index{Arnol'd, Vladimir I.}} 
challenges the reader to compute the following limit.
\begin{quote}
Here is an example of a~problem that people like Barrow, Newton or Huygens would solve in a~couple 
of  minutes, but contemporary mathematicians, in my opinion, are unable to solve it quickly (in any case I have not yet seen a~mathematician who could cope with it quickly): evaluate
$$
\lim_{x\to0}\frac{\sin(\tan x)-\tan(\sin x)}{\arcsin(\arctan x)-\arctan(\arcsin x)}\;.
$$
\end{quote}
Indeed, for a~long time I~could not do anything with the limit. In one July night in 2025 I~only 
laboriously computed that the numerator and denominator are both 
$$
{\textstyle
-\frac{1}{30}x^7+o(x^7)\ \ (x\to0)\,,
}
$$
so that the limit equals $1$. But then a~solution occurred to me. 
Limits of the form
$$
\lim_{x\to0}\frac{f(x)-g(x)}{g^{\langle-1\rangle}(x)-f^{\langle-1\rangle}(x)}=1\,,
$$
like Arnol'd's, follow from the next result on inverse fps (formal power series).

\begin{thm}[fps and its inverse]\label{thm_formArnold}
For $n\in\N$ there exist real polynomials 
$$
p_n(x_1,\,x_2,\,\ds,\,x_n)
$$ 
with $n$ variables such that the following holds. If $A(x)=x+\sum_{n\ge2}a_nx^n$ is a~{\em fps} and
$$
B(x)=A^{\langle-1\rangle}(x)=x+\sum_{n\ge2}b_nx^n
$$
is the {\em fps} inverse of $A(x)$, \underline{then} $b_2=-a_2$ and for every $n\ge3$ we have
$$
b_n=-a_n+p_{n-2}(a_2,\,a_3,\,\ds,\,a_{n-1})\,.
$$
\end{thm}
The mentioned polynomials are functions $p_n\cc\R^n\to\R$ given 
as finite sums of \underline{monomials\index{monomial|emph}}, products of the form
$$
ax_1^{m_1}x_2^{m_2}\ds x_n^{m_n}
$$
where $a\in\R$, $m_i\in\N_0$ and $x_i$ are variables ranging in $\R$.
The fps $A(x)$ and $B(x)$ satisfy
$$
A(B(x))=B(A(x))=x. 
$$
We treat inverse fps in the next passage where we prove 
the theorem. Now we easily deduce from it Arnol'd's limits in general 
form.

\begin{cor}[Arnol'd's limits]\label{cor_ArnoldLim}
Suppose\index{theorem!Arnol'd's limits|emph} 
that
$$
f,\,g\cc(-\de,\,\de)\to\R
$$
are two injective and real analytic functions such that $f\ne g$,  
$f(0)=g(0)=0$ and $f'(0)=g'(0)=1$. \underline{Then}
$$
\lim_{x\to0}\frac{f(x)-g(x)}{g^{\langle-1\rangle}(x)-f^{\langle-1\rangle}(x)}=1\,.
$$
\end{cor}
\duk
Let
$$
{\textstyle
f(x)=x+\sum_{n\ge2}a_nx^n\,\text{ and }\,g(x)=x+\sum_{n\ge2}b_nx^n
}
$$
be the expansions of the functions around $0$ and $m\in\N$ with 
$m\ge2$ be minimum such that $a_m\ne b_m$. If $m=2$ then by 
Theorem~\ref{thm_formArnold} we have
$${\textstyle
\lim_{x\to0}\frac{f(x)-g(x)}{g^{\langle-1\rangle}(x)-f^{\langle-1\rangle}(x)}=
\lim_{x\to0}\frac{a_2x^2-b_2x^2+o(x^2)}{-b_2x^2+a_2x^2+o(x^2)}=1\,.}
$$
Let $m\ge3$. Since $a_2=b_2$, $\ds$, $a_{m-1}=b_{m-1}$, by Theorem~\ref{thm_formArnold} these equalities hold also for the coefficients of $x^2$, $\ds$, $x^{m-1}$ in the expansions of $f^{\langle-1\rangle}(x)$ and $g^{\langle-1\rangle}(x)$. Using Theorem~\ref{thm_formArnold} we again compute
\begin{eqnarray*}
&&\lim_{x\to0}{\textstyle\frac{f(x)-g(x)}{g^{\langle-1\rangle}(x)-f^{\langle-1\rangle}(x)}=}\\
&&=\lim_{x\to0}{\textstyle
\frac{a_mx^m-b_mx^m+o(x^m)}{-b_mx^m
+p_{m-2}(b_2,\,\ds,\,b_{m-1})x^m+a_mx^m-p_{m-2}(a_2,\,\ds,\,a_{m-1})x^m+o(x^m)}=}\\
&&=\lim_{x\to0}{\textstyle
\frac{a_mx^m-b_mx^m+o(x^m)}{-b_mx^m+a_mx^m+o(x^m)}=1\,.
}
\end{eqnarray*}
\kduk

\medskip\noindent
{\em $\bullet$ Inverses of {\em fps}. }We did not consider this 
operation on $\R[[x]]$ in Section~\ref{sec_fps} and we treat 
it briefly now. We look at it more thoroughly in the next final 
section. 

\begin{defi}[inverse fps]\label{def_invFps}
Let $A(x)$ and $B(x)$ be two {\em fps} with zero constant terms. We 
say that $B(x)$ is the \underline{inverse\index{formal power series, fps!inverse of|emph}} of $A(x)$, and we write $B(x)=A^{\langle-1\rangle}(x)$, if
$$
A(B(x))=x\,.
$$
\end{defi}

\begin{exer}\label{ex_invFps}
If $B(x)=A^{\langle-1\rangle}(x)$ then also $B(A(x))=x$.
\end{exer}

\begin{prop}[inverse fps]\label{prop_invFps}
The inverse of a~{\em fps} $A(x)$ is unique. It exists $\iff$ 
$[x^0]A(x)=0$ and $[x^1]A(x)\ne0$. 
\end{prop}
\duk
It is easy to see that if $[x^0]A(x)\ne0$ or if $[x^1]A(x)=0$ 
then the inverse of $A(x)$ does not exist. Suppose that 
$A(x)=\sum_{n\ge0}a_nx^n$ with $a_0=0$ and $a_1\ne0$. We are looking for a~fps $B(x)=\sum_{n\ge1}b_nx^n$ such that 
$${\textstyle
A(B(x))=
\sum_{n\ge1}a_n(b_1x+b_2x^2+\ds)^n
=x\,.}
$$
This relation is equivalent with the system of equations
$$
a_1b_1=1,\, a_1b_2+a_2b_1^2=0,\,\ds,\,
a_1b_n+q_n(a_1,\,\ds,\,a_n,\,b_1,\,\ds,\,b_{n-1})=0,\,\ds
$$
where $n\ge3$ and $q_n(\ds)$ is a~real polynomial with $2n-1$ variables. It is clear that the system has a~unique solution $(b_n)$ in terms of $(a_n)$: the sequence $(b_n)$ begins with $b_1=a_1^{-1}$ and if $b_1$, $b_2$, $\ds$, $b_n$ are already determined, then 
$$
b_{n+1}=-a_1^{-1}\cdot q_{n+1}(a_1,\,\ds,\,a_{n+1},\,b_1,\,\ds,
\,b_n)\,.
$$
\kduk

\noindent
{\bf Proof of Theorem~\ref{thm_formArnold}. }In 
view of the previous proof this proof is easy. We assume $a_1=b_1=1$  and get from the 
system that $b_2=-a_1^{-1}a_2b_1^2=-a_2$. For $n\ge3$
we see that the $n$-th equation in the system is in more detail
$$
a_1b_n+a_nb_1^n+r_n(a_1,\,\ds,\,a_{n-1},\,b_1,\,\ds,\,b_{n-1})=0
$$
where $r_n(\ds)$ is a~real polynomial with $2n-2$ variables. We get by induction
\begin{eqnarray*}
b_n&=&-a_1^{-1}a_nb_1^n-a_1^{-1}r_n(a_1,\,\ds,\,a_{n-1},\,b_1,\,\ds,\,b_{n-1})\\
&=&-a_n-r_n\big(1,\,a_2,\,\ds,\,a_{n-1},\\
&&1,\,-a_2,\,-a_3+p_1(a_2),\,\ds,\,-a_{n-1}+p_{n-3}(a_2,\,\ds,\,a_{n-2})\big)\\
&=&-a_n+p_{n-2}(a_2,\,a_3,\,\ds,\,a_{n-1})
\end{eqnarray*}
where $p_j(\ds)$ are real polynomials with $j$ variables. 
\kduk

\section[Inverses of Taylor polynomials and series]{${}^*$Inverses of Taylor polynomials and series}\label{sec_LIF}

\chapter[Newton integral]{Newton integral}\label{chap_newInt}

\chapter[Riemann integral]{Riemann integral}\label{chap_rieInt}

\chapter[Henstock--Kurzweil integral]{Henstock--Kurzweil integral}\label{chap_hensInt}

\chapter[Applications of integrals]{Applications of integrals}\label{chap_applInt}

\appendix
\chapter[ZFC set theory]{ZFC set theory}\label{doda_notation}

This and the next appendix are devoted to ZFC set theory. As we mentioned at the beginning of the 
second chapter, it is not really possible to rigorously define real 
numbers without using formalized ZFC set theory: we need it to rigorously 
define natural numbers. We need them, in turn, 
to rigorously define integers, which we use in a~rigorous definition of 
fractions. Only when fractions are rigorously defined can we 
build real numbers. We introduce a~part of ZFC set theory 
sufficient  for this undertaking in this appendix.

\section{Logical and set-theoretic notation}\label{sec_logST}

We review notation used in this book. 

\medskip\noindent
{\em $\bullet$ Fraktur and Greek alphabet. }Recall the \underline{Fraktur\index{Fraktur} hand} form of Latin letters:  
$$
\mathfrak{a},\,\mathfrak{A},\,\mathfrak{b},\,\mathfrak{B},\,\mathfrak{c},\,\mathfrak{C},\,\mathfrak{d},\,\mathfrak{D},\,\mathfrak{e},\,\mathfrak{E},\,\mathfrak{f},\,\mathfrak{F},\,\mathfrak{g},\,\mathfrak{G},\,\mathfrak{h},\,\mathfrak{H},\,\mathfrak{i},\,\mathfrak{I},\,\mathfrak{j},\,\mathfrak{J},\,\mathfrak{k},\,\mathfrak{K},\,\mathfrak{l},\,\mathfrak{L},\,\mathfrak{m},\,\mathfrak{M},\,\mathfrak{n},\,\mathfrak{N},
$$
$$
\!\mathfrak{o},\,\mathfrak{O},\,\mathfrak{p},\,\mathfrak{P},\,\mathfrak{q},\,\mathfrak{Q},\,\mathfrak{r},\,\mathfrak{R},\,\mathfrak{s},\,\mathfrak{S},\,\mathfrak{t},\,\mathfrak{T},\,\mathfrak{u},\,\mathfrak{U},\,\mathfrak{v},\,\mathfrak{V},\,\mathfrak{w},\,\mathfrak{W},\,\mathfrak{x},\,\mathfrak{X},\,\mathfrak{y},\,\mathfrak{Y},\,\mathfrak{z}\,\text{ and }\,\mathfrak{Z}\,.\label{fracAlp}
\hspace{0.5cm}
$$
The \underline{Greek alphabet} is used in mathematical notation more frequently:
$$
\alpha,\,\beta,\,\Gamma,\,\gamma,\,\Delta,\,\delta,\,\ep,\,\zeta,\,\eta,\,\Theta,\,\theta,\,\vartheta,\,\iota,
\,\kappa,\,\Lambda,\,\lambda,\,\mu,\,\nu,\,\Xi,\,\xi,\,o,\,\Pi,\,\pi,\,\rho,\,\Sigma,\,\sigma,\,
$$
$$
\tau,\,\Upsilon,\,\upsilon,\,\Phi,\,\phi,\,\varphi,\,\chi,\,\Psi,\,\psi,\,\Omega\,\text{ and }\,\omega\;\hspace{5.4cm}\label{greeAlp}
$$
---\,alpha, beta, gamma, delta, epsilon, zeta, eta, theta, iota, 
kappa, lambda, mu, nu, xi, omicron, pi, rho, sigma, tau, upsilon, phi, 
chi, psi, omega; the omitted capital letters, such as $A$ for $\al$, $B$ for $\be$, $Z$ for $\zeta$ or $H$ for $\eta$, are identical with their Latin form.

\medskip\noindent
{\em $\bullet$ Equality. }We use $=$\label{equality} 
as a set-theoretic equality, governed by the axiom of 
extensionality (Axiom \ref{axi_extens}),\index{axiom!of 
extensionality} but also as a~defining equality, or probably most often as an 
equality of identities. An example of the last type of equalities is
$$
2^3=8\,.
$$
See \cite{henk} for more information on the logic of equality. The symbols $:=$ and 
$=:$\label{defEqu} mean defining equality. For example, $[n]:=
\{1,2,\ds,n\}$ or $\{0,1,2,\ds\}=:\N_0$. 
We define $\N:=\{1,2,\ds\}$. By $\ne$ we denote non-equality. 

\medskip\noindent
{\em $\bullet$ Logical connectives. }As for logical connectives, we write $\neg$ for the 
\underline{negation\index{negation|emph}} (``it is not true that 
$\ds$''),
$\vee$\label{disju} for the \underline{disjunction\index{disjunction|emph}} (``$\ds$ or $\ds$''),\label{disjun} 
$\wedge$\label{conju} or $\&$ for the \underline{conjunction\index{conjunction|emph}} (``$\ds$ and $\ds$''),\label{conjun} 
$\to$\label{impli} for the \underline{implication\index{implication|emph}} (``if $\ds$ then 
$\ds$'') and $\leftrightarrow$\label{equi} for the \underline{equivalence\index{equivalence|emph}} (``$\ds$ if and only if 
$\ds$''). Sometimes we abbreviate the last phrase by `` $\ds$ iff $\ds$''.
We write $\exists$\label{exist} (``there is $\ds$'') and $\forall$\label{gene} 
(``for every $\ds$'') for the \underline{existential\index{quantifiers!existential|emph}} and the 
\underline{general quantifier\index{quantifiers!general|emph}}, respectively. For example, 
$$
\exists\,x\, x>0\,\text{ and }\,
\forall\,a\, a\ne b
$$
mean, respectively, that there is an element in a~domain understood 
from the context that is greater than zero, and that every element 
in the domain differs from the 
element $b$. We use symbols $\Rightarrow$ and $\Leftarrow$ for 
metamathematical implications. The symbol $\iff$ expresses 
metamathematical equivalence.

\medskip\noindent
{\em $\bullet$ Operations with sets. }We use standard set-theoretic notation: $\in$\label{member} 
for the \underline{membership\index{membership|emph}} (``$\ds$ is an element of $\ds$''), $\not\in$\label{nonmember} for the non-membership, $\cup$\label{union} 
for the \underline{union\index{set!union of two|emph}} ($a\cup b=\{c\cc\;c\in a\vee c\in b\}$), $\cap$\label{inters} 
for the \underline{intersection\index{set!intersection of two|emph}} ($a\cap 
b=\{c\in a\cc\;c\in b\}$) and $\setminus$\label{setmin} 
for the \underline{set difference\index{set!difference of two|emph}} ($a\setminus 
b=\{c\in a\cc\;c\not\in b\}$). 
By $\emptyset$\label{emptys} we denote the \underline{empty 
set\index{set!empty|emph}}, the unique set with no elements. For a~set $A$ its \underline{sum\index{set!sum|emph}} is the set 
$${\textstyle
\bigcup A=\{x\cc\;\exists\,y\,
x\in y\wedge y\in A\}\,.\label{setSum} 
}
$$
Similarly, for a~nonempty set $A$ its \underline{intersection\index{set!intersection|emph}} is the set 
$${\textstyle
\bigcap A=
\{x\cc\;\forall\,y\,
y\in A\to x\in y\}\,.\label{setInt}
}
$$

\noindent
{\em $\bullet$ Ordered tuples. }The following definition of ordered pairs is standard. 

\begin{defi}\label{def_KurPai}
Let $x$ and $y$ be sets. The \underline{ordered 
pair\index{ordered pair|emph}} $x,y$ is the set
$$
[x,\,y]:=\{\{x\},\,\{x,\,y\}\}\,.
$$
\end{defi}
This definition is due to {\em Kazimierz Kuratowski (1896--1980)\index{Kuratowski, Kazimierz}} 
in 1921. With the help of the axiom of extensionality\index{axiom!of 
extensionality}, it is easy to prove the basic property of ordered pairs.

\begin{prop}
For every four sets $a,b,c,d$, we have
$$
[a,\,b]=[c,\,d]\iff a=c\wedge b=d\,.\label{ordPair}
$$
\end{prop}
Ordered pairs are standardly generalized to ordered triples, quadruples, $\ds$ by 
iteration of ordered pairs, for example, 
$$
[a,\,b,\,c]:=[a,\,[b,\,c]],\ 
[a,\,b,\,c,\,d]:=[a,\,[b,\,[c,\,d]]],\,\ds\,.
$$
However, it is often overlooked, even in textbooks on set theory, that
this and similar definitions have the defect that the arity of the tuple 
cannot be determined solely from it, and context is needed. For 
example, if we are given the bare set $[a,b,c,d]$, we do not know if it means an ordered pair of 
the sets $a$ and $[b,c,d]$, or 
an ordered triple of the sets $a$, 
$b$, and $[c,d]$, or an ordered quadruple of the sets $a$, $b$, $c$, 
and $d$. We propose and will use an alternative definition of 
ordered tuples.

\begin{defi}\label{def_kTup}
Let $x_1$, $x_2$, $\ds$, $x_k$ be $k\ge2$ sets. The corresponding ordered 
\underline{$k$-tuple\index{ktuple@$k$-tuple|emph}}, or briefly ordered \underline{tuple}, is the set
$$
\langle x_1,\,x_2,\,\ds,\,x_k
\rangle:=\{[\{\emptyset\},\,x_1],\,
[\{\{\emptyset\}\},\,x_2],\,\ds,\,[\{\{\ds\{\emptyset\}\ds\}\},\,x_k]\}\,.\label{ordTuple}
$$
\end{defi}
In such tuples, the arity is uniquely determined. We leave the proof of the 
next proposition to the reader as an exercise.

\begin{prop}
Let $k,l\in\N\setminus\{1\}$, and let $a=\langle x_1,\,x_2,\ds,x_k
\rangle$ and 
$b=\langle y_1,\,y_2,\ds,y_l
\rangle$ be two ordered tuples. \underline{Then}
$$
a=b\iff k=l\wedge x_1=y_1\wedge x_2=y_2\wedge\ds\wedge x_k=y_k\,.
$$
\end{prop}

\section{ZFC formulas}\label{sec_ZFC}

In this appendix, we start to build the hierarchy of numeric domains
$$
\N_0-\Z-\Q-\R^*-\R-\C
$$
We build the algebraic structure of natural numbers $\N_0$ in this appendix and in 
Section~\ref{sec_ZFCtoN0}. The algebraic structures $\Z$, $\Q$, $\R^*$, and $\R$ are constructed in Sections~\ref{sec_N0toZ}--\ref{sec_RstarToR}, respectively. The last structure, the field of complex numbers $\C$, will be constructed (from $\R$) in Appendix~\ref{sec_RtoC}.

\medskip\noindent
{\em $\bullet$ Alphabets, words, and concatenations. }We have already 
introduced these and related notions in 
they form the foundation of mathematics, we review them here in more detail. An 
{\em alphabet} is any nonempty set $A$. A~{\em word} over $A$ is any finite, possibly empty, 
sequence $u=u_1u_2\ds u_m$, $m\ge0$, with every $u_i$ in $A$. We denote the 
set of words over $A$ by $A^*$. The ability to form words over an alphabet 
is really the engine that propels mathematics.

If $v=v_1v_2\ds v_n$, $n\ge 0$, is another word, and if there exists 
an index $l$ such that 
$$
v_1=u_l,\,v_2=u_{l+1},\ds,v_n=u_{l+n-1}
$$
(so that $l+n-1\le m$), we say that $v$ is a~{\em subword} of $u$, and we call the 
interval of (intuitive) 
natural numbers $\{l,l+1,\ds,l+n-1\}$ the \underline{$v$-copy\index{copy of 
a~subword|emph}} in $u$. The {\em concatenation} of $u$ and $w$ is the word
$$
w=uv:=w_1w_2\ds w_{m+n}
$$
such that $w_i=u_i$ for $i=1,2,\ds,m$ and $w_i=v_{i-m}$ for 
$i=m+1,m+2,\ds,m+n$. We know that concatenation is an associative (but not 
commutative) operation, and we speak of a word obtained by concatenating several 
words as an \underline{concatenation\index{concatenation|emph}}.

\medskip\noindent
{\em $\bullet$ {\em ZFC} formulas. }In view of the previous passage, we 
consider a~countably infinite alphabet
$$
\mathcal{A}=\{x,\,y,\,z,\,u,\,
v,\,w,\,x_0,\,x_1,\,\ds\}\cup
\{\in,\,=,\,\neg,\,\wedge,\,\vee,\,\to,\,
\leftrightarrow,\,\forall,\,\exists,\,
(,\,),\,[,\,]\}
$$
and words over $\mathcal{A}$. In the disjoint union forming $\mathcal{A}$, the former countable set $V$ consists of \underline{set 
variables\index{set 
variables|emph}}, which range in naive sets. The set 
variables $x_0$, $x_1$, $\ds$ are indexed by the intuitive 
natural numbers. The latter thirteen-element set of \underline{logical symbols\index{logical symbols|emph}} 
contains the symbols $\in$ of {\em membership}, $=$ of {\em equality}, $\neg$ of {\em negation}, $\wedge$ of {\em 
conjunction}, $\vee$ of {\em 
disjunction}, $\to$ of {\em implication}, $\leftrightarrow$ of {\em equivalence}, 
$\forall$ of the {\em general quantifier}, $\exists$ of the {\em existential 
quantifier}, and four symbols of brackets $($, $)$, $[$, and $]$. We assume that 
the meaning and use of these symbols are known to the reader. We define ZFC 
formulas as particular words over $\mathcal{A}$. An \underline{atomic formula\index{atomic formulas|emph}} is any concatenation of the form $(a\circ b)$ where $a$ and $b$ are in $V$ and $\circ$ is $\in$ or $=$.

\begin{defi}\label{def_formulas}
A~nonempty word $\varphi$ in $\mathcal{A}^*$ is called 
a~\underline{{\em ZFC} formula\index{formula@(ZFC) 
formula|emph}}, briefly a~\underline{formula}, if there exists so
called \underline{generating word\index{generating word|emph}} $g$ of $\varphi$,
$$
g\,\text{ is }\,u_1u_2\ds u_l\,\text{ in }\,(\mathcal{A}^*)^*,\ l>0\,,
$$
such that $u_l$ is $\varphi$ and for every $k\le l$, one of the following cases holds.
\begin{enumerate}
\item $u_k$ is an atomic formula.
\item $u_k$ is the concatenation $\neg u_j$ for some $j<k$.
\item $u_k$ is the concatenation 
$$
(u_i\circ u_j)
$$ 
for some indices $i,j<k$ and $\circ$ in
$\{\wedge,\vee,\to,\leftrightarrow\}$.
\item $u_k$ is the concatenation 
$$
(\circ\, a\, u_j)
$$ 
for some $\circ$ in 
$\{\forall,\exists\}$, $a$ in $V$, and $j<k$.
\end{enumerate}
\end{defi}
Clearly, $u_k$ is a~formula for every $k\le l$. This definition of ZFC 
formulas has an advantage over shorter and seemingly more elegant definitions. It provides an algorithm by which one can decide for any word 
$u\in\mathcal{A}^*$ if $u$ is a~ZFC formula. We leave the formulation of this algorithm to the interested reader as an 
exercise. The hint is to consider
the shortest generating word of $\varphi$. An example of a~formula is the axiom of extensionality
$$
(\forall x(\forall y((x=y)\leftrightarrow(\forall z((z\in x)\leftrightarrow (z\in y))))))\,.
$$
We can simplify this expression by omitting the brackets around atomic 
formulas and around the whole formula:
$$
\forall x(\forall y(x=y\leftrightarrow(\forall z(z\in x\leftrightarrow z\in y))))\,.
$$
In some cases, we replace for better readability some pairs of round brackets
by square ones. For example, we can write the axiom of extensionality as
$$
\forall x(\forall y(x=y\leftrightarrow[\forall z(z\in x\leftrightarrow z\in y)]))\,.
$$

\noindent
{\em $\bullet$ Subformulas, bound and free occurrence of set variables. }A~\underline{subformula\index{subformula|emph}} of a~formula $\varphi$ is any subword of $\varphi$ 
that is a~formula. Let $a$ be in $V$ and $\varphi$ be a~formula. An 
occurrence of $a$ in $\varphi$ is \underline{bound\index{occurrence of a~set variable!bound|emph}} (by a~quantifier) if it is included in 
a~$\psi$-copy in $\varphi$ of a~subformula $\psi$ of $\varphi$ such 
that $\psi$ has form $(\circ\, a\,\psi')$, where $\circ$ is $\forall$ or $\exists$ 
and $\psi'$ is a~formula. If the occurrence of $a$ is not bound, it is called \underline{free\index{occurrence of a~set variable!free|emph}}. A~\underline{sentence\index{sentence|emph}} is a~formula with no free occurrence of any set variable. Let 
$\varphi$ be a~formula and $a$ be a~set variable. We write $\varphi(a)$ (not to 
be confused with a~concatenation) if the only free occurrences of set variables
in $\varphi$ are those of $a$ (the set variable $a$ may not occur in 
$\varphi$). We define the more general notation $\varphi(a_1,\ds,a_k)$, where 
the $a_i$ are distinct set variables, similarly. 

\begin{defi}\label{def_qform}
By a~\underline{quasi-formula\index{quasi-formula|emph}} we 
mean any description of a~formula that is sufficiently clear and precise so 
that the formula can be easily and uniquely reconstructed from it.    
\end{defi}

\section{ZFC axioms}

Claims and statements about naive sets are expressed by 
sentences. In the second part of this section, we explain when 
a~sentence, a~claim about sets, is true. For this, we have to list the ZFC 
axioms. They form a~countably infinite set of sentences. ZFC axioms are as 
follows. 

\medskip\noindent
{\em $\bullet$ The axiom 
of extensionality. }Two sets $x$ and $y$ are equal if and only if they have the same elements $z$. 

\begin{axio}\label{axi_extens}
The \underline{axiom of extensionality\index{axiom!of extensionality|emph}} is the sentence
$$
\forall x(\forall y(x=y\leftrightarrow(\forall z(z\in x\leftrightarrow z\in y))))\,.
$$
\end{axio}

\noindent
{\em $\bullet$ The axiom of pairing. }For any two sets $x$ and $y$ there exists a~set $z=\{x,y\}$.

\begin{axio}\label{axi_pair}
The \underline{axiom of pairing\index{axiom!of pairing|emph}} is the sentence 
$$
\forall x(\forall y(\exists z(\forall u(u\in z\leftrightarrow(u=x\vee u=y)))))\,.
$$
\end{axio}

\noindent
{\em $\bullet$ The axiom of union. }For any set $x$ there exists a~set 
$${\textstyle
y=\bigcup x:=\{z\cc\;\exists u\, z\in u\in x\}\,.
}
$$

\begin{axio}\label{axi_sum}
The \underline{axiom of union\index{axiom!of union|emph}} is the sentence 
$$
\forall x(\exists y(\forall z(z\in y\leftrightarrow(\exists u(z\in u\wedge u\in x)))))\,.
$$
\end{axio}

\noindent
{\em $\bullet$ The axiom of power set. }For any set $x$ there exists a~set 
$${\textstyle
y=\mathcal{P}(x):=\{z\cc\;z\sus x\}\,.
}
$$

\begin{axio}\label{axi_power}
The \underline{axiom of power set\index{axiom!of power set|emph}} is the sentence
$$
\forall x(\exists y(\forall z(z\in y\leftrightarrow(\forall u(u\in z\to u\in x)))))\,.
$$
\end{axio}

\noindent
{\em $\bullet$ The axiom of foundation. }Every nonempty set $x$ contains an 
$\in$-minimal element $y$. 

\begin{axio}\label{axi_found}
The \underline{axiom of foundation\index{axiom!of foundation|emph}} is the sentence
$$
\forall x(\exists x_0(x_0\in x)\to
[\exists y(y\in x\wedge(\neg\exists x_1(x_1\in x\wedge x_1\in y)))])\,.
$$
\end{axio}
This axiom forbids sets $x$ such that $x\in x$, cycles like $x_0\in x_1\in\ds\in x_n\in x_0$, and infinite chains like $\ds\in x_1\in x_0$.

\medskip\noindent
{\em $\bullet$ The axiom of infinity. }There exists an inductive set $x_0$. Recall that $x_0$ is inductive if $\emptyset\in x_0$ and for every set $y\in x_0$, also $y\cup\{y\}\in x_0$.

\begin{axio}\label{axi_infi}
The \underline{axiom of infinity\index{axiom!of infinity|emph}} is the sentence
\begin{eqnarray*}
&&\exists x_0([\exists x_1((x_1\in x_0)\wedge(\neg\exists x_2\, x_2\in x_1))]\wedge[\forall x_3((x_3\in x_0)\\
&&\to(\exists x_4(x_4\in x_0\wedge(\forall x_5((x_5\in x_4)\leftrightarrow(x_5=x_3\vee x_5\in x_3))))))])\,.
\end{eqnarray*}
\end{axio}
 
\noindent
{\em $\bullet$ The axiom of choice. }
For every set $x$ of nonempty and mutually disjoint sets there exists 
a~set $y$ such that for every $z\in x$ we have $y\cap z=\{x_5\}$ for some set $x_5$. 

\begin{axio}[AC]\label{axi_AC}
The \underline{axiom of choice\index{axiom!of choice, AC|emph}} is 
the sentence
\begin{eqnarray*}
&&\forall x([(\forall x_0((x_0\in x)\to\exists x_1(x_1\in x_0)))\wedge(\forall x_2(\forall x_3((((x_2\in x)\\
&&\wedge(x_3\in x))\wedge(\exists x_4((x_4\in x_2)\wedge(x_4\in x_3))))\to x_2=x_3)))]\\
&&\to[\exists y(\forall z((z\in x)\to(\exists x_5(((x_5\in z)\wedge(x_5\in y))\\
&&\wedge(\forall x_6((x_6\in z)\wedge(x_6\in y))\to(x_6=x_5))))))])\,.    
\end{eqnarray*}
\end{axio}
The AC is the ``C'' in ZFC. The last two ZFC axioms are two axiom schemata, infinite families of axioms. 

\medskip\noindent
{\em $\bullet$ The axiom schema of separation. }For every property 
$\varphi(x)$ of sets and 
every set $y$, there exists a~set $z$ of the elements $x$ in $y$ that possess 
the property.

\begin{axio}\label{axi_separ}
By the \underline{axiom schema of separation\index{axiom!schema of separation|emph}}, for 
every {\em ZFC} formula 
$$
\varphi(x,\,x_1,\,\ds,\,x_n),\ n\ge0\,, 
$$
the sentence
$$
\forall x_1(\ds\forall x_n(
\forall y(\exists z(\forall x(x\in z\leftrightarrow(x\in y\wedge\varphi))))))
$$
is a~{\em ZFC} axiom.
\end{axio}
The last displayed word is not a formula, only a~template for it: the 
three slots, $\ds$, $n$, and $\varphi$, are to be filled in the ways clear from 
the context. The variables $x_i$ play the role of parameters. 

\medskip\noindent
{\em $\bullet$ The axiom schema of replacement. }Let $\varphi(x,y)$ be a~property of 
sets such that ``$x\mapsto y=f(x):$ $\varphi(x,y)$ holds'' is a~function. 
Then for every set $z$ there exists a~set $u$ such that $u$ arises from 
$z$ by replacing every $x\in z$ with $y=f(x)$.

\begin{axio}\label{axi_replac}
By the \underline{axiom schema of replacement\index{axiom!schema of replacement|emph}}, for every {\em ZFC} formula 
$$
\varphi(x,\,y,\,x_1,\,\ds,\,x_n),\  n\ge0\,, 
$$
the sentence
\begin{eqnarray*}
&&\forall x_1(\ds\forall x_n([\forall x(\exists v(\forall y(y=v\leftrightarrow\varphi)))]\\
&&\to[\forall z(\exists u(\forall y(y\in u\leftrightarrow (\exists x(x\in z\wedge\varphi)))))]))    
\end{eqnarray*}
is a~{\em ZFC} axiom.
\end{axio}
Again, the last displayed word is a~template for a~formula, in which the 
four slots  
$\ds$, $n$, and two $\varphi$'s, are to be filled in the ways clear from the 
context. We summarize the 
axioms of ZFC set theory as follows. 

\begin{defi}\label{def_ZFCaxio}
Let $\mathrm{ZFC}_{axi}$ ($\sus\mathcal{A}^*$) be the axioms of 
{\em ZFC} set theory. Then $\varphi\in\mathrm{ZFC}_{axi}$ $\iff$ 
$\varphi$ is Axiom~\ref{axi_extens} of extensionality or Axiom~\ref{axi_pair} 
of pairing or Axiom~\ref{axi_sum} of union or Axiom~\ref{axi_power} of power 
set or  Axiom~\ref{axi_found} of foundation or Axiom~\ref{axi_infi} of 
infinity or Axiom~\ref{axi_AC} of choice or one of  
Axioms~\ref{axi_separ} of separation  or one of Axioms~\ref{axi_replac} 
of replacement.
\end{defi}
In some popular explanations and discussions of axiomatic systems and theories, axioms 
are portrayed as simple and self-evident statements. ZFC set theory is 
not such. First, it has infinitely many axioms. Second, these axioms are far 
from straightforward. We could have written ZFC axioms with the help of 
abbreviations, for example,
the power set axiom can be stated as
$$
\forall x(\exists y(\forall z(z\in y\leftrightarrow z\sus x)))\,,
$$
where $\sus$ abbreviates the subset relation. The reason we did not do 
it is that it would be less rigorous than what we did. Although various 
abbreviations are frequently used in explanations of ZFC set theory, precise 
syntactic rules for their use are {\em never} given, which makes abbreviations
a~not-completely-rigorous device. The  matter of abbreviations is subtle, as 
one has to be aware of how they are used. For example, to be completely 
precise, we must now make clear that the axioms of logic stated 
below in the second part (we took them from \cite{soch}) are formulated for a set of 
logical symbols smaller than the one introduced above in $\mathcal{A}$. 
Namely, apart 
from $\in$, $=$, and brackets, only $\neg$, $\to$, and $\forall$ are 
allowed. Other logical symbols in $\mathcal{A}$ that appear in ZFC 
formulas have to be treated, if we want to use only these axioms of logic, as 
abbreviations: $(\varphi\wedge\chi)$ is $\neg(\varphi\to\neg\chi)$, $(\varphi\vee\chi)$ is $(\neg\varphi\to\chi)$,  
$(\varphi\leftrightarrow\chi)$ is 
$((\varphi\to\chi)\wedge
(\chi\to\varphi))$ (with the conjunction further expanded according 
to the first abbreviation), and $(\exists a\varphi)$ is $\neg(\forall 
a\neg\varphi)$, where $\varphi$ and $\chi$ are ZFC formulas, and $a$ is in 
$V$. As you see, a~rigorous syntactic theory of abbreviations is badly needed.

\section{Syntactic and semantic truth}\label{sec_howDoWe}

\begin{quote}
{\footnotesize
[Bellman:] ``Just the place for a Snark! I have said it twice:\\
That alone should encourage the crew.\\
Just the place for a Snark! I have said it thrice:\\
What I tell you three times is true.''

\noindent L.~Carroll \cite[p.~3]{carr}
}
\end{quote}

\noindent
From the previous section, we know that claims about sets are ZFC 
sentences. But which of them are true? In fact, what does it mean for 
a~sentence to be true? And what is truth? 
We first present a~syntactic Definition~\ref{def_truth} of 
truth of a~sentence. It is somewhat similar to the definition Bellman revealed to his crew 
in \cite{carr}, see the above motto. In the second part of this section, we 
discuss another kind of truth of a~sentence, the semantic one, and the 
relation of the two kinds of truth.

In order to define that a~sentence is true we need more axioms: the {\em 
axioms of logic}. As before, they are words in $\mathcal{A}^*$. Here is 
the 
\underline{propositional axiom schema\index{propositional axiom schema|emph}}. 

\begin{axio}\label{axi_propoSch}
For every {\em ZFC} formulas $\varphi$, $\chi$, and $\psi$, the three formulas
$$
(\varphi\to(\chi\to\varphi)),\ 
(\varphi\to(\chi\to\psi)\to((\varphi\to\chi)\to(\varphi\to\psi)))\,,
$$
and
$$
((\neg\varphi\to\neg\chi)\to(\chi\to
\varphi))
$$
are axioms of logic.
\end{axio}

For the next \underline{quantification axiom schema\index{quantification axiom schema|emph}}, we need the operation of substituting one set 
variable for another. Let $a$ and $b$ be set variables and $\varphi$ be a~formula. 
We say that $b$ is \underline{substitutable\index{substitutability of one set variable for another|emph}} for $a$ in $\varphi$, if no free 
occurrence of $a$ in $\varphi$ is included in a~$\psi$-copy in $\varphi$ 
of a~subformula $\psi$ of $\varphi$ that has the form $(\circ\, b\,\psi')$, with $\circ$ in 
$\{\forall,\exists\}$ and $\psi'$ a~formula. In other words, after 
substituting all free occurrences of $a$ by $b$'s, the number of bound 
occurrences of variables in $\varphi$ does not change (increase). We denote the 
formula resulting from this substitution
(it is easy to prove that it is a~formula) by $\varphi(a/b)$.

\begin{axio}\label{axi_quanSche}
Let $a$ and $b$ be set variables and $\varphi$ be any 
{\em ZFC} formula in which $a$ is substitutable for $b$. Then 
$$
((\forall b\varphi)\to\varphi(a/b))
$$
is an axiom of logic. Let $a$ be a~set variable, and $\varphi$ and $\chi$ be two {\em ZFC} formulas such that $a$ has no free occurrence in $\varphi$. Then
$$
((\forall a(\varphi\to\chi))\to(\varphi\to(\forall a\chi)))
$$
is an axiom of logic. 
\end{axio}

The last third \underline{equality axiom schema\index{equality axiom schema|emph}} concerns the logical symbol $=$ of equality.

\begin{axio}\label{axi_equaSch}
For any set variable $a$, the formula $(a=a)$ is an axiom of logic. Let $a$, $b$, $c$, and $d$ be set variable. Then the two formulas
$$
((((a=b)\wedge(c=d))\wedge(a\in c))\to(b\in d))
$$
and
$$
((((a=b)\wedge(c=d))\wedge
(a=c))\to(b=d))
$$
are axioms of logic.
\end{axio}

We summarize the axioms of logic as follows.

\begin{defi}[axioms of logic]\label{def_logicAxi}
Let 
$$
\mathrm{LOG}_{axi}\ \  (\sus\mathcal{A}^*) 
$$
be the \underline{axioms of 
logic\index{axioms of 
logic|emph}}. Then $\varphi\in\mathrm{LOG}_{axi}$ $\iff$ 
$\varphi$ is one of the propositional Axioms~\ref{axi_propoSch} or one of 
the quantification Axioms~\ref{axi_quanSche} or one of 
the equality Axioms~\ref{axi_equaSch}. 
\end{defi}

The definition of ZFC theorems is similar to Definition~\ref{def_formulas} of formulas. We define ZFC theorems, or briefly, theorems.  

\begin{defi}[ZFC theorems]\label{def_ZFCthm}
A~formula $\varphi$ is a~{\em ZFC} theorem, briefly a~theorem, if there 
exists so called \underline{proof\index{proof of a~ZFC formula|emph}} $g$ of $\varphi$,
$$
g\,\text{ is }\,u_1u_2\ds u_l\,\text{ in }\,(\mathcal{A}^*)^*,\ l>0\,,
$$
such that $u_l$ is $\varphi$ and for every $k\le l$, one of the following cases holds.
\begin{enumerate}
\item $u_k$ is in $\mathrm{ZFC}_{axi}\cup \mathrm{LOG}_{axi}$.
\item (the modus ponens rule) There exist indices $i,j<k$ such that $u_j$ 
is the concatenation $(u_i\to u_k)$.
\item (the generalization rule) There exists an index $j<k$ and a~set 
variable $b$ such that $u_k$ is the concatenation $(\forall\, b\, u_j)$.
\end{enumerate}    
\end{defi}
Again, $u_k$ is a~theorem for every $k\le l$. Similarly, the interested 
reader will easily devise an algorithm that $\ds$; no, we are kidding. Despite 
the outer similarity, there is a big difference between 
Definition~\ref{def_formulas} and Definition~\ref{def_ZFCthm}. In contrast to the 
former definition, by a ground-breaking result in mathematical logic due to {\em Alan~M.~Turing (1912--1954)\index{Turing, Alan M.}}, no algorithm exists that can decide if a formula is a theorem. What 
is the feature of Definition~\ref{def_ZFCthm} that prevents devising an algorithm in 
the style of Definition~\ref{def_formulas}? It is the modus ponens deduction rule; it 
allows one to deduce {\em shorter} formulas from longer ones.

In any case, we define truth of ZFC formulas as follows.

\begin{defi}[truth of formulas]\label{def_truth}
We say that a~{\em ZFC} formula is syntactically true if it is a~{\em ZFC} 
theorem.     
\end{defi}
For example, every ZFC axiom 
and every axiom of logic is true.
In view of Turing's mentioned result, there actually {\em is} a considerable 
difference between this definition of truth and Bellman's. We assume that informal proofs are stated so carefully that 
they can be converted with some effort into the form of 
Definition~\ref{def_ZFCthm}. 

\begin{center}
$***$    
\end{center}

\noindent
Bla

\section{What are infinite sets, really?}\label{expl_infiSets}

Bla

\chapter{The equiconsistence of AC: G\"odel}\label{app_godel}

The second goal, realized in Sections~\ref{app_godel} and \ref{app_cohen}, is to explain in 
detail how the axiom of choice, which was introduced in 
Section~\ref{sec_AC}, is situated in mathematics. In Section~\ref{app_godel}, 
we present the proof due to K.~G\"odel\index{goedel@G\"odel, Kurt} of the theorem asserting that if ZF set theory is 
consistent (i.e., free of contradictions), then so is $\mathrm{ZFC}=\mathrm{ZF}+
\mathrm{AC}$, the 
ZF set theory with the axiom of choice added. The analogous theorem for   
$\mathrm{ZF}+\neg\mathrm{AC}$, when the negation of the axiom of choice is added to ZF, was proven by P.~Cohen\index{Cohen, Paul}. We present his 
proof in Section~\ref{app_cohen}.

The theorems of G\"odel and Cohen belong to the greatest achievements in 
mathematics of the 20th century; they transformed our understanding of mathematical truth. The proofs are complicated, technically and 
conceptually, and explaining them properly and in detail is not easy. 
Fortunately, I~could lean on the works of the Czech school of set theory, the monograph \cite{soch} due to {\em 
Anton\'\i n Sochor (1942--2008)} and the textbook \cite{balc_step} due to 
{\em Bohuslav Balcar (1943--2017)} and {\em Petr \v St\v ep\'anek (1943--2012)}. I~devote this appendix to 
their memory. 

{\em Kurt G\"odel (1906--1978)\index{goedel@G\"odel, Kurt}}

\chapter{The equiconsistence of $\neg$AC: Cohen}\label{app_cohen}

{\em Paul Cohen (1934--2007)\index{Cohen, Paul}}

\chapter[From real numbers to complex numbers\checkmark]{From real numbers to complex numbers}\label{sec_RtoC}
\index{complex numbers|(}

In this appendix, we extend the hierarchy of numeric domains built in 
Chapters~\ref{sec_ZFCtoN0}--\ref{sec_RstarToR} by the complex numbers. 
We introduce their algebraic structure and prove the characterization 
Theorem~\ref{thm_Cis} which asserts that $\C$ is a~field extension of $\R$ with 
degree $[\C:\R]=2$ and is unique up to isomorphism.

Definitions~\ref{def_complNum} and  \ref{def_arithC} introduce the set of 
complex numbers and their arithmetic. 
Definition~\ref{def_CasStr} describes the algebraic structure~$\C$. In 
Definition~\ref{def_ifiel} we introduce quadratic fields. 
Theorem~\ref{thm_Cis} is the characterization theorem.
Proposition~\ref{prop_CisNormFie} shows that $\C$ is a~complete normed 
field. The basic property of the field $\C$ expressed by the {\em 
Fundamental Theorem of Algebra}, which asserts that every non-constant complex 
polynomial has a~root 
in $\C$, is proven in the following appendix in Theorem~\ref{thm_fta}.

\section{Complex numbers and their arithmetic}

We define the set 
of complex numbers.

\begin{defi}\label{def_complNum}
The set of \underline{complex numbers\index{complex numbers!set|emph}} is 
$\C:=\R\times\R$. For $a,b\in\R$ we use the standard notation 
$$
a+bi:=\langle 
a,\,b\rangle\ \ (\in\C)\,. 
$$
Sometimes, we write $a\oplus bi$ instead of $a+bi$ for better 
readability.
\end{defi}

Recall the arithmetic on $\R$ introduced in Chapter~\ref{sec_QtoRstar}. In the next 
definition, the operations $+$, $-$, and $\cdot$ on the right  of $:=$ 
are in $\R$. 

\begin{defi}\label{def_arithC}
We define the \underline{addition\index{complex numbers!addition|emph}} and 
\underline{multiplication\index{complex numbers!multiplication|emph}} of complex numbers $a+bi$ and $c+di$ as follows.
\begin{enumerate}
\item $a\oplus 
bi+c\oplus di:=(a+c)\oplus(b+d)i$.
\item $a\oplus 
bi\cdot c\oplus 
di:=(a\cdot c-b\cdot d)\oplus(a\cdot d+b\cdot c)i$.
\end{enumerate}
\end{defi}
We introduce the algebraic structure of complex numbers.

\begin{defi}\label{def_CasStr}
The algebraic structure of \underline{complex numbers\index{complex 
numbers!structure|emph}}   
$$
\C=\langle\C,\,0_{\C},\,1_{\C},\,+,\,
\cdot
\rangle
$$ 
consists of the set $\C=\R^2$ in Definition~\ref{def_complNum}, the elements $0_{\C}:=0+0i$ and $1_{\C}:=1+0i$ in $\C$, and the operations $+$ and $\cdot$ on $\C$ given in Definition~\ref{def_arithC}. 
\end{defi}
Instead of $0_{\C}$ and $1_{\C}$, we usually write just $0$ and $1$.

\section{The characterization theorem}

We introduce an algebraic structure that fits the complex numbers. We suppose 
that the reader is familiar with the basics of vector spaces over 
fields, especially with bases and dimensions. Recall 
Definition~\ref{def_field} of fields 
and recall Theorem~\ref{thm_fieEmb} that every field homomorphism is 
injective and is called a~field embedding. So let $F$ and $K$ be 
fields and
$$
f\cc F\to K
$$ 
be a~field embedding. We say that $K$  \underline{extends\index{field!extension of|emph}} $F$. We can view the field $K$ as a~vector space over $F$: the addition of 
vectors is the addition in the field $K$, and for $a\in F$ and $b\in K$, the scalar multiplication of $b$ by $a$ is defined by
$$
ab:=f(a)\cdot b\,,
$$
where $\cdot$ is multiplication in $K$. It is easy to check the axioms of vector spaces.
We denote the dimension of this vector space by $[K:F]$ ($\in\N\cup
\{\infty\}$) and call it the \underline{degree of $K$ over 
$F$\index{degree of a~field extension|emph}}.  

\begin{defi}\label{def_ifiel}
A~field $F$ is called \underline{quadratic\index{quadratic 
field|emph}} if it extends the field $\R$ and has degree $[F:\R]=2$.    
\end{defi}

We state the characterization theorem. 

\begin{thm}\label{thm_Cis}
The\index{theorem!complex numbers|emph}
structure $\C$ is a~quadratic field. Every two quadratic 
fields are isomorphic as fields.    
\end{thm}

\section{$\C$ is a~quadratic field}

We prove the former claim of the theorem right away in the next 
proposition. The latter claim is proven in 
Proposition~\ref{prop_ifie2}.

\begin{prop}\label{prop_ifie1}
The algebraic structure 
$$
\C=\langle\R^2,\,0,\,1,\,+,\,\cdot\rangle
$$ 
introduced in Definition~\ref{def_CasStr} is a~quadratic field.    
\end{prop}
\duk
The neutrality of $0$ to $+$ and $1$ to $\cdot$ follows easily from 
Definition~\ref{def_arithC}. From this definition and the fact that $\R$ is a~field, it 
follows that $+$ is commutative and 
associative, and that $\cdot$ is commutative. We show that $\cdot$ is 
associative. Let $z_1=a+bi$, $z_2=c+di$, and $z_3=e+fi$ be complex numbers. Then, indeed, 
\begin{eqnarray*}
(z_1\cdot z_2)\cdot z_3&=&((ac-bd)\oplus(ad+bc)i)\cdot(e\oplus fi)\\
&=&((ac-bd)e-(ad+bc)f)\oplus((ac-bd)f+(ad+bc)e)i\\
&=&(a(ce-df)-b(de+cf))\oplus(a(cf+de)-b(df-ce))i\\
&=&(a\oplus bi)\cdot((ce-df)\oplus(de+cf)i)=z_1\cdot(z_2\cdot z_3)\,.
\end{eqnarray*}
Similarly,
\begin{eqnarray*}
z_1\cdot(z_2+z_3)&=&(a\oplus bi)\cdot((c+e)\oplus(d+f)i)\\
&=&(a(c+e)-b(d+f))\oplus(a(d+f)+b(c+e))i\\
&=&((ac-bd)\oplus(ad+bc)i)+
((ae-bf)\oplus(af+be)i)\\
&=&z_1\cdot z_2+z_1\cdot z_3
\end{eqnarray*}
and we see that $\cdot$ is distributive to $+$.

Let $z=a+bi\in\C$. Then
$$
(a\oplus bi)+(-a\oplus(-b)i)=0\oplus0i=0
$$
and, if $a\ne0$ or $b\ne0$ and $c:=(a^2+b^2)^{-1}$, 
$$
(a\oplus bi)\cdot(ac\oplus(-bc)i)=
(a^2+b^2)c\oplus(-abc+bac)i=1\oplus0i=1\,.
$$
Hence $z$ has an additive inverse, and if $z\ne0$, it has a~multiplicative 
inverse. We have proven that $\C$ is a~field. It remains to show that $\C$ 
is an extension of $\R$ with degree $2$. It is easy to see that 
$$
f\cc\R\to\C,\ f(a)=a\oplus0i\,,
$$
is a~field embedding. Also,  
$\{1_{\C},\iota\}$, where 
$\iota:=0+1i$, is a~linear basis of the vector space $\C$ over $\R$, so that 
$[\C:\R]=2$.
\kduk

\section{$\C$ is unique}

For the proof of the isomorphism claim, we need a~lemma.

\begin{lemma}
In every quadratic field $F$ there exists an element $\iota\in F$ such that 
$$
\iota^2=\iota\cdot\iota=-1_F\,.
$$
\end{lemma}
\duk
Let $F$ be a~quadratic field. We may assume that $\R$ is a~subfield of $F$. Since 
$[F:\R]=2$, we have $F\setminus\R\ne\emptyset$ and can take 
an element $\al\in F\setminus\R$. Then the set $\{1_F,\al\}$ is linearly 
independent over $\R$, but 
$\{1_F,\al,\al^2\}$ is linearly dependent. Thus
$$
a1_F+b\al+c\al^2=0\,,
$$
where $a,b,c\in\R$ and $c\ne0$. Equivalently,  
$$
(\al+b/2c)^2=b^2/4c^2-a/c=:d\ \ (\in\R)\,.
$$
If $d\ge0$, using Lemma~\ref{lem_odmvR} we deduce 
the contradiction that $\al\in\R$. Thus $d<0$ and again using 
Lemma~\ref{lem_odmvR}, we see that the element
$${\textstyle
\iota:=\frac{\al+b/2c}{\sqrt{-d}}\ \ (\in F)
}
$$
has the property that $\iota^2=-1_F$.
\kduk

We prove that every two quadratic fields are isomorphic.

\begin{prop}\label{prop_ifie2}
Suppose that
$$
F=\langle
F,\,0_F,\,1_F,\,+,\,\cdot
\rangle\,\text{ and }\,
G=\langle
G,\,0_G,\,1_G,\,\oplus,\,\odot
\rangle
$$
are two quadratic fields. Then $F$ and $G$ are isomorphic, which 
means that there exists a~bijection $f\cc F\to G$ such that for every $a,b\in F$,
$$
f(a+b)=f(a)\oplus f(b)\,\text{ and }\,
f(a\cdot b)=f(a)\odot f(b)\,.
$$
\end{prop}
\duk
Using the previous lemma, we take elements $\iota\in F$ and $\kappa\in G$
such that $\iota^2=-1_F$ and $\kappa^2=-1_G$. For the simplicity of 
notation, we assume that $\R$ is a~subfield of $F$, and that $g\cc\R\to G$ is a~field embedding. It follows that 
$$
B=\{1_F,\,\iota\}\,\text{ and }\,
B'=\{1_G,\,\kappa\}
$$ 
is a~linear basis of $F$ over $\R$ and $G$ over $\R$, respectively. For
any $a\in F$ we write
$a=u\cdot1_F+v\cdot\iota\in F$ with $u,v\in\R$. We define a~map $f\cc F\to G$ by 
$$
f(a):=g(u)\odot1_G\oplus g(v)\odot \kappa\,.
$$

We show that $f$ is the desired isomorphism. Since $B$ and $B'$ are 
bases, $f$ is a~bijection.
Let $a,b\in F$ with 
$$
a=u\cdot 1_F+v\cdot\iota\,\text{ and }\,b=u'\cdot1_F+v'\cdot\iota,\,\text{ where }\,u,\,v,\,u',\,v'\in\R\,. 
$$
Then
\begin{eqnarray*}
f(a+b)&=&f((u+u')\cdot 1_F+(v+v')\cdot \iota)=g(u+u')\odot 1_G\oplus g(v+v')\odot\kappa\\
&=&g(u)\odot1_G\oplus g(v)\odot\kappa
\oplus g(u')\odot1_G\oplus g(v')\odot\kappa=f(a)\oplus f(b) 
\end{eqnarray*}
and
\begin{eqnarray*}
f(a\cdot b)&=&f((u\cdot u'-v\cdot v')\cdot1_F+(u\cdot v'+v\cdot u')\cdot\iota)\\
&=&g(u\cdot u'-v\cdot v')\odot1_G\oplus g(u\cdot v'+v\cdot u')\odot\kappa\\
&=&(g(u)\odot g(u')\ominus g(v)\odot g(v'))\odot 1_G\oplus\\
&&\oplus\,(g(u)\odot g(v')\oplus
g(v)\odot g(u'))\odot\kappa\\
&=&(g(u)\odot 1_G\oplus g(v)\odot\kappa)\odot(g(u')\odot 1_G\oplus g(v')\odot\kappa)\\
&=&f(a)\odot f(b)\,.
\end{eqnarray*}
In these computations, we used the definitions of $f$, $\iota$, and $\kappa$, the fact that $g\cc\R\to G$ is a~field embedding,  
and properties of the operations $+$, $\cdot$, $\oplus$, and $\odot$, 
for example, the distributive law.
\kduk

\noindent
This concludes the proof of Theorem~\ref{thm_Cis} and the entire 
buildup of the hierarchy of numeric domains that we started in 
Chapter~\ref{sec_ZFCtoN0}. A~simple field embedding of $\R$ in $\C$ is given by
$$
\R\ni a\mapsto a+0i\in\C\,.
$$ 

The previous structures $\Z$, $\Q$, $\R^*$, and $\R$ are linearly 
ordered, and the absolute value is defined for each by means of this ordering. 
Complex numbers cannot form an ordered field because, in any ordered field $F$,
squares are nonnegative and $-1_F<0_F$, but $-1$ is a~square in $\C$.
Interestingly, the \underline{absolute value\index{complex numbers!absolute 
value|emph}} on $\C$, with its standard properties, can still be defined. We set 
$$
|a+bi|:=\sqrt{a^2+b^2}\ \ (\in \R_{\ge0})\,.
$$
It clearly extends the absolute value on $\R$ (we regard 
$\R$ as a~subset of $\C$).

\begin{prop}\label{prop_CisNormFie}
The absolute value on $\C$ has the following properties.
\begin{enumerate}
\item $|z|\ge0$ and $|z|=0$ iff $z=0$. 
\item $|u\cdot v|=|u|\cdot|v|$.
\item The \underline{triangle inequality\index{complex numbers!triangle inequality|emph}} $|u+v|\le|u|+|v|$ holds.
\item The metric $|u-v|$ is \underline{complete\index{complex numbers!complete space|emph}}, i.e. every Cauchy sequence $(u_n)\sus\C$ has a~limit in the metric space $\langle\C,|u-v|\rangle$.
\end{enumerate}
\end{prop}
\duk
1. This is easy to see. 2. We compute
\begin{eqnarray*}
|(a+bi)(c+di)|&=&\sqrt{(ac-bd)^2+(ad+bc)^2}\\
&=&\sqrt{(a^2+b^2)(c^2+d^2)}\\
&=&\sqrt{(a^2+b^2)}\cdot\sqrt{(c^2+d^2)}=
|a+bi|\cdot|c+di|\,.
\end{eqnarray*}

3. We may assume that $u\ne0$. Using part~2 we see that it suffices to prove the triangle inequality for $u=1$. Let $v=a+bi$. Then
\begin{eqnarray*}
&&|1+v|\le 1+|v|\iff\sqrt{(a+1)^2+b^2}\le
1+\sqrt{a^2+b^2}\\
&&\iff(a+1)^2+b^2\le1+a^2+b^2+2\sqrt{a^2+b^2}\\
&&\iff 2a\le2\sqrt{a^2+b^2}\iff0=0\,.
\end{eqnarray*}

4. Let $(a_n+b_ni)\sus\C$ be a~Cauchy sequence. It follows that $(a_n),(b_n)\sus\R$ are real Cauchy sequences in the complete metric space 
$\langle\R,|x-y|\rangle$. We set $a:=\lim a_n$, $b:=\lim b_n$, and easily check that
in $\C$,
$$
\lim_{n\to\infty}(a_n+b_ni)=a+bi\,.
$$
\kduk

\noindent
An important property of the field $\C$ is its algebraic closedness. We prove 
it in Theorem~\ref{thm_fta}.
\index{complex numbers|)}

\chapter[Metric spaces]{Metric spaces}\label{sec_MS}

We review basics of metric spaces. After necessary definitions we treat in a~novel way the theory of their completions. 

\section{Metric spaces}

\begin{defi}\label{def_MS}
A~\underline{metric space\index{metric spaces|emph}} is a~pair $M=\langle M,d\rangle$ such that $M\ne\emptyset$ is a~set, $d\cc M\times M\to\R$ is a~function called a~\underline{metric\index{metric spaces!metric, distance|emph}} or a~\underline{distance}, and three conditions hold ($x,y,z\in M$). 
\begin{enumerate}
\item Always $d(x,y)\ge0$, with equality only for $x=y$.
\item Always $d(x,y)=d(y,x)$ (symmetry).
\item Always $d(x,y)\le d(x,z)+d(z,y)$ (the \underline{triangle inequality\index{metric spaces!triangle inequality|emph}}, briefly {\em TI}).
\end{enumerate}
\end{defi}

Let $M=\langle M,d\rangle$ 
be a~metric space, $b\in M$, and $r>0$ be a~real number. A~\underline{ball\index{metric spaces!ball in|emph}} in the space, 
centered at the point $b$ and with the radius $r$,
is the set 
$$
B_M(b,\,r)=B(b,\,r)=B=\{a\in M\cc\;d(a,\,b)<r\}\,. \label{ball}
$$
Always $B\ne\emptyset$ because $b\in B$.
We call a~set $X\sus M$ \underline{dense (in $M$)\index{dense sets|emph}} if 
$B\cap X\ne\emptyset$ for every ball $B$ in $M$. A~sequence of points $(b_n)\sus M$ in the space  
\underline{converges\index{metric spaces!convergent sequences in|emph}} to 
the \underline{limit\index{metric 
spaces!limits in|emph}} $b\in M$ if for every $\ep>0$ there is an $n_0\in\N$ 
such that for every $n\ge n_0$ we have 
$d(b,\,b_n)\le\ep$. The sequence $(b_n)$ is \underline{Cauchy\index{Cauchy 
sequence in a~metric space|emph}} if for every $\ep>0$ there is an $n_0\in\N$ 
such that for every $m,n\ge n_0$ we have 
$d(b,\,b_n)\le\ep$. The space $M$ is \underline{complete\index{complete metric spaces|emph}} if every Cauchy sequence in it has a~limit. Clearly, $X\sus M$ is dense iff for every $b\in M$ 
(equivalently, every $b\in M\setminus X$) we have a~sequence $(a_n)\sus X$ 
such that $b=\lim a_n$. 

\section{Completions of metric spaces: a~new definition}

We define completions of metric spaces. 
We say that a~metric space $M=\langle M,d\rangle$ is 
a~\underline{subspace\index{subspaces of metric spaces|emph}} of a~metric space 
$N=\langle N,e\rangle$ if $M\sus N$ and $e\,|\,M\times M=d$. The spaces $M$ and 
$N$ are \underline{isometric\index{isometric metric spaces|emph}} if there is 
a~bijection $F\cc M\to N$, called an \underline{isometry\index{isometry of 
metric spaces|emph}} of $M$ and $N$, such that $d(x,y)=e(F(x),F(y))$ 
for every $x,y\in M$.

\begin{defi}[new]\label{def_completions}
Let $M=\langle M,d\rangle$ be a~metric space. A~metric space
$C=\langle C,e\rangle$ is a~\underline{completion\index{completions of metric spaces|emph}} of $M$ if the pair of $M$ and $C$ has three properties.
\begin{enumerate}
\item $M$ is a~subspace of $C$.
\item Every Cauchy sequence $(a_n)\sus M$ has a~limit in $C$.
\item The set $C\setminus M$ is minimal to inclusion with respect to
property~2.
\end{enumerate}
\end{defi}
Our definition is new --- we do not require $C$ to be complete, nor that $M$ be dense in $C$. These properties
of $C$ follow from our definition, which is logically simpler and more 
natural than the standard definition of completions one can encounter in the literature.
We prove three results. 1. Every completion of $M$ is complete and $M$ is 
dense in it. 2. Completions are unique up to isometry. 3.~Completions exist. 

\section{Completions are complete}

We prove the first result. 

\begin{prop}\label{prop_compSp1}
Let $C=\langle C,e\rangle$ be a~completion of 
$M=\langle M,d\rangle$. Then $M$ is dense in $C$ and the metric space $C$ is complete.       
\end{prop}
\duk
The first claim is immediate by 
parts~2 and~3 in Definition~\ref{def_completions}. We show that $C$ is complete. 
Let $(a_n)\sus C$ be a~Cauchy sequence in the space $C$. Since $M$ is dense in 
$C$, we can take points $b_n\in M$ such that $e(a_n,b_n)\le n^{-1}$ for every 
$n\in\N$. We show that $(b_n)$ is a~Cauchy sequence in $M$. Let $\ep>0$ be given. Then $e(a_m,a_n)\le\ep$ for every
indices $m,n\ge n_0$. We take an index $n_1\ge n_0$ such that $e(a_n,b_n)\le \ep$ for every $n\ge n_1$. Then, using TI, we obtain for every $m,n\ge n_1$ the bound
$$
d(b_m,\,b_n)=e(b_m,\,b_n)\le e(b_m,\,a_m)+e(a_m,\,a_n)+e(a_n,\,b_n)
\le\ep+\ep+\ep=3\ep
$$
and see that $(b_n)$ is Cauchy. By part~2 of Definition~\ref{def_completions}, there is a~point $a\in C$ such that
$$
\lim_{n\to\infty}b_n=a\,.
$$
We show that $\lim a_n=a$ as well. Let $\ep>0$ be given. We take an index $n_2\ge n_1$ such that $e(b_n,a)\le\ep$ for every $n\ge n_2$.
Then for every $n\ge n_2$,
$$
e(a_n,\,a)\le e(a_n,\,b_n)+e(b_n,\,a)\le\ep+\ep=2\ep\,.
$$
We see that $\lim a_n=a$.
\kduk

\section{Completions are unique}

We prove the second result. We need a~lemma.

\begin{lemma}\label{lem_pomoLem}
Let $\ep>0$, let $M=\langle M,d\rangle$ be a~metric space, and let $a,b,c\in M$ 
be such that $d(b,c)\le\ep$. Then we have the bound
$$
|d(a,\,b)-d(a,\,c)|\le\ep\,.
$$
\end{lemma}
\duk
Indeed, by TI we have inequalities $d(a,b)\le d(a,c)+d(b,c)$ and
$d(a,c)\le d(a,b)+d(b,c)$. They are equivalent to the inequalities
$$
d(a,\,b)-d(a,\,c)\le d(b,\,c)\le\ep\,\text{ and }\,
d(a,\,b)-d(a,\,c)\ge-d(b,\,c)\ge-\ep\,,
$$
respectively, and these are equivalent to the stated bound. 
\kduk

\begin{prop}\label{prop_compSp2}
Let $C=\langle C,e\rangle$ and $C'=\langle C',e'\rangle$ be two completions of $M=\langle M,d\rangle$. Then the metric spaces 
$C$ and $C'$ are isometric.       
\end{prop}
\duk
We define a~map $F\cc C\to C'$ and show that it is an isometry of $C$ and $C'$. 
Let $a\in C$ be arbitrary. Since $M$ is dense in $C$ 
(Proposition~\ref{prop_compSp1}), $a=\lim b_n$ for some $(b_n)\sus M$. 
Clearly, $(b_n)$ is Cauchy. By part~2 of 
Definition~\ref{def_completions}, we have $\lim b_n=a'$ for 
some $a'\in C'$. We define $F\cc C\to C'$ by setting
$$
F(a):=a'\,.
$$
We show that $F$ does not depend on the choice of the sequence $(b_n)$. Let $(b_n')\sus M$ be another sequence such that $\lim b_n'=a$. By TI, 
$$
e'(b_n,\,b_n')=d(b_n,\,b_n')=e(b_n,\,b_n')\le 
e(b_n,\,a)+e(a,\,b_n')\to0\,\text{ as }n\to\infty\,. 
$$
Hence $\lim b_n=\lim b_n'$.

We show that $F$ is a~bijection and that it preserves distances. To prove the injectivity of $F$, we take $a_0,a_1\in C$, 
$a_0\ne a_1$, and sequences $(b_n^i)\sus M$ such that 
$\lim b_n^i=a_i$ for $i=0,1$. Then there is $\ep>0$ such that 
$d(b_n^0,b_n^1)\ge\ep$ for every large $n$. It follows from this 
that $F(a_0)\ne F(a_1)$. To prove the surjectivity of $F$, we 
need that $F(a)$ does not depend on the sequence $(b_n)\sus M$ chosen such that $\lim b_n=a$. We take
any point $a'\in C'$. Since $M$ is dense in $C'$, we have $a'=\lim b_n'$ 
for some $(b_n')\sus M$. Since 
$(b_n')$ is Cauchy, by part~2 of 
Definition~\ref{def_completions} we have $\lim b_n'=a$ 
for some $a\in C$. Hence $F(a)=a'$. 

Finally, to prove that $F$ preserves distances, we take two points $a_0,a_1\in C$, and sequences $(b_n^i)\sus M$ such that 
$\lim b_n^i=a_i$ for $i=0,1$. Then
\begin{eqnarray*}
&&e(a_0,\,a_1)=\lim_{n\to\infty}e(b_n^0,\,b_n^1)=\lim_{n\to\infty}d(b_n^0,\,b_n^1)=\lim_{n\to\infty}e'(b_n^0,\,b_n^1)=\\
&&=e'(F(a_0),\,
F(a_1))\,, 
\end{eqnarray*}
and $e(a_0,a_1)=e'(F(a_0),
F(a_1))$. To prove the first equality in detail, 
we take $\ep>0$ and $n_0\in\N$ such that $e(a_i,\,b_n^i)\le\ep$
for $i=0,1$ and $n\ge n_0$. Then for every $n\ge n_0$ we 
indeed have, by Lemma~\ref{lem_pomoLem}, that 
\begin{eqnarray*}
|e(a_0,\,a_1)-e(b_n^0,\,b_n^1)|&\le&
|e(a_0,\,a_1)-e(b_n^0,\,a_1)|+|e(b_n^0,\,a_1)-
e(b_n^0,\,b_n^1)|\\
&\le&\ep+\ep=2\ep\
\end{eqnarray*}
(consider triangles $a_0,b_n^0,a_1$ and $b_n^0,b_n^1,a_1$). The equality
$$
\lim e'(b_n^0,\,b_n^1)=e'(F(a_0),\,
F(a_1))
$$ 
is proven in the same way, we only work in the space $C'$ instead of $C$.
\kduk

\section{Completions exist}

We prove the third result. Two sequences $(a_n)$ and $(a_n')$ in a~metric space $M=\langle M,d\rangle$ are {\em equivalent}, written $(a_n)\sim(a_n')$, if
$\lim d(a_n,\,a_n')=0$. It is clear that $\sim$ is an equivalence relation. 
We denote by $\mathcal{C}(M)$ the set of all Cauchy sequences in $M$. We say that an {\em (equivalence) block}  
$B\in \mathcal{C}(M)/\!\sim$ is {\em bad} if there is no point $b\in M$ such that $\lim a_n=b$ for some (equivalently, any)
sequence $(a_n)\in B$. 

\begin{thm}\label{thm_compSp3}
Any\index{theorem!existence of completions of metric spaces|emph} 
metric space $M=\langle M,d\rangle$ has a~completion
$C=\langle C,e\rangle$.
\end{thm}
\duk
Let $M=\langle M,d\rangle$ be a~metric space and let
$$
M':=\{B\in \mathcal{C}(M)/\!\sim\cc\;\text{$B$ is bad}\}\,.
$$
Using the device  of disjoint union ($X,Y\leadsto(\{0\}\times X)\cup
(\{1\}\times Y)$),  
we may assume that $M\cap M'=\emptyset$. We set $C:=M\cup M'$. We first define a~map $e\cc C\times C\to\R$. Then 
we show that $\langle C,e\rangle$ is a~metric space. Finally, we show that 
the pair of $\langle M,d\rangle$ and 
$\langle C,e\rangle$ has properties 1--3 in Definition~\ref{def_completions}. 

Let $a,b\in C$. If $a,b\in M$, we define $e(a,b):=d(a,b)$. If 
$a$ is in $M$ and $b=B\in M'$, we take any sequence $(b_n)\in B$ and set 
$$
e(a,B):=\lim_{n\to\infty} d(a,\,b_n)\ \ (\in\R)\,. 
$$
This limit exists by the completeness of the metric space $\langle\R,|x-
y|\rangle$, because the 
sequence $(d(a,b_n))\sus\R$ is Cauchy by Lemma~\ref{lem_pomoLem}. Moreover, it is 
easy to show that 
the limit does not depend on the choice of the sequence $(b_n)$ in $B$. We 
similarly define values $e(A,b)$ for 
$A\in M'$ and $b\in M$. For $A,B\in M'$
and any sequences $(a_n)\in A$ and $(b_n)\in B$ we set
$$
e(A,B):=\lim_{n\to\infty} d(a_n,\,b_n)\ \ (\in\R)\,.
$$
The existence of this limit and its independence on the choice of sequences
$(a_n)$ and $(b_n)$ is justified as before.

We show that $\langle C,e\rangle$ is a~metric space. Clearly, 
$e$ is a nonnegative and symmetric function because $d$ has these properties. Let $a,b\in C$. 
Let $a=b$. If $a=b\in M$ then $e(a,b)=d(a,b)=0$ because $d$ is a~metric. If $a=A=B=b\in M'$, then 
$e(a,b)=0$ because we can take twice the 
same sequence $(a_n)=(b_n)\in A=B$ to compute the limit $e(a,b)=\lim d(a_n,b_n)=\lim d(a_n,a_n)=0$.
Let $e(a,b)=0$. If $a,b\in M$, then $a=b$ because $d$ is a metric. Let $a\in M$, $b=B\in M'$, 
and $(b_n)\in B$. Then $e(a,B)=\lim d(a,b_n)=0$ is, in fact, impossible because then $\lim b_n=a$ (in $M$) and $B$ is not a~bad block; 
same if $a=A\in M'$ and $b\in M$. Let $a=A,b=B\in M'$, and let $(a_n)\in A$ and $(b_n)\in B$. Then
$e(A,B)=\lim d(a_n,b_n)=0$ means that $(a_n)\sim(b_n)$ and 
$A=B$. We show that TI holds for $e$. This easily follows  from three facts: the definition of $e$, TI holds
for $d$, and the fact that non-strict real inequalities are preserved in limits (if $x_n\le y_n$ for every $n$, then $\lim 
x_n\le\lim y_n$, if these limits exist). Thus $\langle C,e\rangle$ is a~metric space. 

We show that $C=\langle C,e\rangle$ is a~completion of $M$ by  
Definition~\ref{def_completions}.
By the definition of $e$, the space 
$M$ is a~subspace of $C$. Let  $(a_n)\in\mathcal{C}(M)$. If $\lim a_n=b$ (in $M$) for some $b\in M$,
we are done. Else $A=[(a_n)]_{\sim}\in M'$. Then
$$
\lim_{n\to\infty}e(A,\,a_n)=\lim_{n\to\infty}\lim_{m\to\infty}d(a_m,\,a_n)=0
$$
because $(a_n)$ is a~Cauchy sequence in $M$. 
Hence $\lim a_n=A$ (in $C$). Finally, to prove property~3 in Definition~\ref{def_completions}, we take a~bad block 
$A=[(a_n)]_{\sim}\in M'$ and delete it from $C$. But then the sequence $(a_n)$
($\in\mathcal{C}(M)$) does not have a~limit in $C\setminus\{A\}$, because in $C$ it 
has the limit $A$, and property~2 is lost. The proof is complete. 
\kduk

\section[Real numbers via completions of ordered Abelian groups]{Real numbers via completions of ordered Abelian groups}

The very first idea is that by using the construction in the last subsection 
for the Euclidean metric space $\langle\Q,|x-y|\rangle$, 
one obtains real numbers. Alas, the very definition of a~metric space 
presupposes an already defined set of real numbers. In this section, we 
describe another general construction in the spirit of metric completion, 
which does produce real numbers as a~particular case.   

\section{Ordered Abelian groups}

We remind ordered  Abelian groups and show that they are endowed with 
an absolute value satisfying a~triangle inequality.

\begin{defi}
An \underline{ordered Abelian group} is a~quadruple $A=\langle A,0_A,+,
<\rangle$ such that $\langle A,0_A,+\rangle$ is an Abelian group,  
$\langle A,<\rangle$ is a~linear order, and both structures are related by the 
\underline{order axiom\index{order axiom in Abelian groups|emph}} that
$$
a<b\Rightarrow a+c<b+c\,\text{ for every }\,a,\,b,\,c\in A\,.
$$
\end{defi}
It follows from the order axiom that for every element $a\in A^*:=A\setminus\{0_A\}$ we have 
$$
-a<0_A<a\,\text{ or }\,a<0_A<-a\,.
$$
Another corollary of the order axiom is that if $a,b<0_A$, then $a+b<0_A$, and 
the same holds for $<$. Also, $a<b$ iff $-b<-a$. For any $a\in A$ we define
$$
|a|:=\max(-a,\,a)\,.
$$

\begin{prop}\label{prop_TIorAb}
Let $A=\langle A,0_A,+,
<\rangle$ be an ordered Abelian group, and let $a,b\in A$. Then
$$
|a+b|\le|a|+|b|\,.
$$
\end{prop}
\duk
If $a=0_A$ or $b=0_A$, the inequality holds trivially as an equality. Let 
$a,b<0_A$. Then $a+b<0_A$ and
$$
|a+b|=-(a+b)=-a+(-b)=|a|+|b|\,.
$$
If $a,b>0_A$, the proof is similar and even simpler. Let $a<0_A<b$. It follows 
from the order axiom that 
$$
a<a+b<b\,.
$$
If $a+b\le0_A$, then 
$$
|a+b|=-(a+b)<-a=|a|\le|a|+|b|\,.
$$
If $a+b>0_A$, then 
$$
|a+b|=a+b<b=|b|\le|a|+|b|\,.
$$
In the case $b<0_A<a$, the proof is the same.
\kduk

\section{Completions of ordered Abelian groups}

In the following, let $A=\langle A,0_A,+,<\rangle$ be an ordered Abelian 
group. Let $(a_n)\sus A$ be a~sequence of points in it. If for every $e\in A$ 
with $e>0_A$ there is an $n_0\in\N$ such that $|a_m-a_n|<e$ for every 
$m,n\ge n_0$, we say that $(a_n)$ is \underline{Cauchy\index{Cauchy 
sequence!in an ordered Abelian group|emph}}. If $b\in A$ and for every 
$e\in A$ with $e>0_A$ there is an $n_0\in\N$ such that $|a_n-b|<e$ for 
every $n\ge n_0$, we say that \underline{$(a_n)$ has the limit 
$b$\index{limits in ordered Abelian groups|emph}} and write $\lim a_n=b$ or
$\lim_{n\to\infty} a_n=b$. We check that these limits are unique.

\begin{prop}
If $\lim a_n=b$ and $\lim a_n=b'$, then $b=b'$.    
\end{prop}
\duk
Suppose that $b<b'$. The first case is that there is no $c\in A$ with $b<c<b'$.
Then we set $e:=b'-b$ ($>0_A$). 
It follows from the definition of limit that  for every large $n$ we have $a_n\le b$ and, at the same time, 
$b'\le a_n$. This is a~contradiction.

Suppose that there is a~$c\in A$ such that $b<c<b'$. We set 
$$
e:=\min(c-b,\,b'-c)\ \ (>0_A)\,. 
$$
It follows from the definition of limit and from Proposition~\ref{prop_TIorAb} that  for every large $n$, 
$$
|b-b'|\le|b-a_n|+|a_n-b'|<(c-b)+(b'-c)=
b'-b=|b-b'|\,.
$$
This is a~contradiction.
\kduk

\noindent
We say that $A$ is 
a~\underline{complete\index{complete |emph}} ordered Abelian groups, if 
every Cauchy sequence in it has a~limit.
Let $B=\langle B,0_B,\oplus,\prec\rangle$ be another 
ordered Abelian group. We say that $A$ is a~\underline{substructure} of $B$ if
$A\sus B$, $0_A=0_B$, $\oplus\,|\,A\times A$ is $+$, and 
$\prec\cap\,(A\times A)$ is $<$. We define the completion of $A$. 

\begin{defi}\label{def_complOrAb}
Let $A=\langle A,0_A,+,<\rangle$ and $B=\langle B,0_B,\oplus,\prec\rangle$ 
be two ordered Abelian groups. We say that $B$ is a~\underline{completion\index{completion of an ordered Abelian group|emph}} of $A$ if the following holds.
\begin{enumerate}
\item $A$ is a~substructure of $B$.
\item Every Cauchy sequence in $A$ has a~limit in $B$.
\item The set $B\setminus A$ is minimal 
to inclusion with respect to condition~2. 
\end{enumerate}
\end{defi}

\section{Completions are complete}

We mimic the proofs in the previous section. 

\begin{prop}\label{prop_compAb1}
Let 
$$
B=\langle B,\,0_B,\,\oplus,\,\prec\rangle
$$ 
be a~completion of 
$A=\langle A,0_A,+,<\rangle$. Then 
$A$ is dense in $B$ and the ordered Abelian group $B$ is complete.       
\end{prop}

\section{The Fundamental Theorem of Algebra}\label{sec_FTAl}

\begin{thm}\label{thm_fta}
Let $f\in\C[z]$ be a~non-constant complex polynomial. Then $f(z)$ has 
a~root in $\C$, i.e. there exists a~number $\al\in\C$ such that $f(\al)=0$.      
\end{thm}
\duk
We divide the proof in two steps. The first step is a~reduction to 
polynomials $z^n-u$. In the second step, we prove that these polynomials 
have roots.

{\em First step. If for every $n\in\N$ and every $u\in\C$ the polynomial $z^n-u$ has a~root in $\C$, then every non-constant complex polynomial has a~root in $\C$. }Let

{\em Second step. For every $n\in\N$ and every $u\in\C$ the polynomial $z^n-u$ has a~root in $\C$. }Let
\kduk

{\small \chapter[Solutions to exercises]{Solutions to exercises}\label{doda_reseni}

\centerline{{\bf 1 Prologue on functions and relations}}

\medskip\noindent
{\bf Exercise~\ref{ex_rovnFci} }Implication $\Rightarrow$ is trivial. We prove $\Leftarrow$.
Suppose that $G_f=G_g$. Then $A$ is the set of elements that appear 
as the first components of the pairs in $G_f$. Similarly for $C$. We 
deduce that $A=C$. Thus $A=C$ and $G_f=G_g$, which means that $f$ and 
$g$ are congruent.

\medskip\noindent
{\bf Exercise~\ref{ex_partFun} }Let $C$ be a~partial function from $A$ to $B$. We consider the set $A'=\{a\in A\cc\;\exists\,b\in 
B\cc a\,C\,b\}$. Then $\langle A',B,C\rangle$ is a~function 
from $A'$ to $B$.

\medskip\noindent
{\bf Exercise~\ref{ex_onEmpFun} }This follows from the definition of 
isomorphism of functions.

\medskip\noindent
{\bf Exercise~\ref{ex_obrazAvzor} } In general, neither equality holds.

\medskip\noindent
{\bf Exercise~\ref{ex_restEmpFun} }If $f$ and $g$ are functions and $f$ is empty, then $G_f=\emptyset\sus G_g$. 

\medskip\noindent
{\bf Exercise~\ref{ex_uloNeuPr} }Suppose that the elements $a$ and 
$b$ are neutral to a~commutative operation~$o$ on a~set $X$. Then 
$b=a+b=b+a=a$.

\medskip\noindent
{\bf Exercise~\ref{ex_onFamFun} }There are $3^5=243$ such words. 

\medskip\noindent
{\bf Exercise~\ref{ex_assConc} }Let $u=u_1\ds u_l$, $v=v_1\ds v_m$, and 
$w=w_1\ds w_n$ be words over $X$, with $l,m,n\in\omega$. We show that 
the concatenations $z=(uv)w=z_1\ds z_{l+m+n}$ and $z'=u(vw)=z_1'\ds 
z_{l+m+n}'$ are equal. Let $i\in[l+m+n]$. If $i\in[l]$, then 
$z_i=u_i=z_i'$. If $l+1\le i\le l+m$, then $z_i=v_{i-l}=z_i'$. Finally, if 
$l+m+1\le i\le l+m+n$, then
$z_i=w_{i-l-m}=z_i'$.

\medskip\noindent
{\bf Exercise~\ref{ex_comConc} }If $X=\{0,1\}$, $u=0$, and $v=1$, then
$uv=01\ne10=vu$.

\medskip\noindent
{\bf Exercise~\ref{ex_uloNaDrFunk} }If and only if $X=Y$.

\medskip\noindent
{\bf Exercise~\ref{ex_drFciArovn} } For injective, constant, and
identity functions, it is true. It is not true for surjective and bijective functions.

\medskip\noindent
{\bf Exercise~\ref{ex_naHomom} }Let $c\in Y$ be arbitrary. We take $b\in X$ such that $f(b)=c$. Then 
$f(a)\,q\,c=f(a)\,q\,f(b)=f(a\,p\,b)=f(b)=c$. If $f$ is not onto, we cannot 
argue about the elements in $Y\setminus f[X]$.

\medskip\noindent
{\bf Exercise~\ref{ex_naHomom1} }Let $c,d\in Y$ be arbitrary, 
$a=f^{-1}(c)$,  and $b=f^{-1}(d)$. Then 
$f^{-1}(c\,q\,d)=f^{-1}(f(a)\,q\,f(b))=f^{-1}(f(a\,p\,b))=
a\,p\,b=f^{-1}(c)\,p\,f^{-1}(d)$.

\medskip\noindent
{\bf Exercise~\ref{ex_naHomom2} }Let $a,b\in X$ be arbitrary. We see that $g(f)(a\,p\,b)=g(f(a\,p\,b))=g(f(a)\,q\,f(b))=g(f(a))\,r\,g(f(b))=g(f)(a)\,r\,g(f)(b)$. 

\medskip\noindent
{\bf Exercise~\ref{ex_jednoZnac} }No problem, for injective $f$ both meanings 
coincide.

\medskip\noindent
{\bf Exercise~\ref{ex_invInv} }Since $f^{-1}\cc f[X]\to X$ is bijective, we have $(f^{-1})^{-1}\cc X\to f[X]$. 
So for $f[X]\ne Y$, the functions $(f^{-1})^{-1}$ and $f$ differ as 
triples, but they are isomorphic. We have
$((f^{-1})^{-1})^{-1}\cc f[X]\to X$. Hence $((f^{-1})^{-1})^{-1}$ and $f^{-1}$ are equal as triples.

\medskip\noindent
{\bf Exercise~\ref{ex_rovnAskla} }Yes, they are.

\medskip\noindent
{\bf Exercise~\ref{ex_slozInjSurj} }Let $f$ and $g$ be injective
and $f(g)(x)=f(g)(y)$. 
Thus $f(g(x))=f(g(y))$ and $g(x)=g(y)$. Thus $x=y$. Let
$g\cc X\to Y$ and $f\cc Y\to B$ be onto, and let $b\in B$. Thus, there is $y\in Y$ such that
$f(y)=b$. Thus, there is $x\in X$ such that $g(x)=y$. Thus, $f(g)(x)=f(g(x))=f(y)=b$ and $f(g)$ is
onto. If $g\cc X\to Y$ and $f\cc A\to B$ are nonempty and surjective, and $Y\cap A=\emptyset$, then $f(g)\cc\emptyset\to B$ is 
not onto.

\medskip\noindent
{\bf Exercise~\ref{ex_asocSkla} }The range of $f(g(h))$ and the 
range of $f(g)(h)$ is the range $R$ of $f$. We show that
$f(g(h))$ and $f(g)(h)$ have equal graphs. Indeed, each graph is the set  of pairs $\langle x,y\rangle\in M(h)\times 
R$ such that there exist $a\in M(g)$ and $b\in M(f)$ such that $h(x)=a$, $g(a)=b$, and $f(b)=y$.

\medskip\noindent
{\bf Exercise~\ref{ex_slozNaProste} }Let $Y=h[X]$ and 
$g=\langle X,Y,G_g\rangle =\langle X,Y,G_h\rangle$. The map $f\cc Y\to Z$ is given by $f(y)=y$.

\medskip\noindent
{\bf Exercise~\ref{ex_oBijekci} }If $f$ is a bijection, then $g=
f^{-1}$ has the required properties. 
Let $g\cc Y\to X$ be as stated. Since $g(f)(x)=g(f(x))=x$, the function $f$ is
injective. Since $M(g)=Y$ and $f(g)(y)=f(g(y))=y$, the function $f$ is onto.

\medskip\noindent
{\bf Exercise~\ref{ex_kdyInvCon} }Those having empty or one-element definition domains.

\medskip\noindent
{\bf Exercise~\ref{ex_uloNaSS} }These mean $\bigcup_{i\in\N}A_i$, $\bigcup_{i\in\N}A_i$, and 
$\bigcap_{i\in\omega}A_i$, respectively. 

\medskip\noindent
{\bf Exercise~\ref{ex_equiAC} }1 $\Rightarrow$ 2. Let $X$ be as stated. We consider the set 
system $\{A_Z\cc\;Z\in X\}$ with $A_Z=Z$ and define $Y=S[X]$. 
2 $\Rightarrow$ 1. Let $\{A_i\cc\;i\in I\}$ with $A_i\ne\emptyset$ be a~given set system. We consider the set $X$ given by $a\in X$
iff $a=\{i\}\times A_i$ for some $i\in I$, and define $S(i)=y$ where 
$Y\cap X=\{\langle i,y\rangle\}$. 
1 $\Rightarrow$ 3. Let $f$ be as stated. We consider
the set system $\{A_b\cc\;b\in B\}$ with $A_b=f^{-1}[\{b\}]$ and define $g\cc B\to 
A$ by $g(b)=S(b)$. Finally, 3 $\Rightarrow$ 1. Let $\{A_i\cc\;i\in I\}$ 
with $A_i\ne\emptyset$ be a~given set system. We take the sets 
$A=\bigcup_{i\in I}\{i\}\times A_i$ and $B=I$, and the surjection $f\cc A\to B$ given by 
$f(\langle i,x\rangle)=i$. Then we define a~selector $S$ by $S(i)=x$ where $x$ is given by $g(i)=\langle i,x\rangle$. 

\medskip\noindent
{\bf Exercise~\ref{ex_oneMoreAC} }Consider the set $X=\{\{1,2\},\,\{2,3\},\,\{3,1\}\}$. 

\medskip\noindent
{\bf Exercise~\ref{ex_bloky} }Let $R$ be an equivalence relation on 
a~set $A\ne\emptyset$ (for $A=\emptyset$ 
everything trivially holds) and $[a]_R$ be an equivalence block. 
Clearly, $a\in[a]_R$. Thus the elements in $A/R$ are nonempty 
and $\bigcup A/R=A$. Let $a,b\in A$ and $c\in [a]_R\cap[b]_R$.
By the transitivity and symmetry of $R$ we have $a R b$. Thus 
$[a]_R=[b]_R$ and the elements of $A/R$ are pairwise disjoint. 

If $b,c\in[a]_R$ then $bRa$ and $cRa$. Thus $bRc$. If $bRc$, $b\in[a]_R$ and 
$c\in[a']_R$ then $bRa$, $cRa'$, so that $aRa'$. Hence $[a]_R=[a']_R$ and 
$b,c$ are in a~common block.

\medskip\noindent
{\bf Exercise~\ref{ex_rozklad} }Let $X$ be a~partition of $Y\ne\emptyset$ and $R=Y/X$. For $y\in Y$ we take 
a~block $Z\in X$ with $y\in Z$. So $yRy$ and $R$ is reflexive. For $y,y'\in Y$ with $yRy'$ 
there is a~block $Z\in X$ such that $y,y'\in Z$. So also $y'Ry$ and $R$ is 
symmetric. Let $y,y',y''\in Y$ with $yRy'$ a~$y'Ry''$. Thus there exist blocks $Z,Z'\in X$ such that 
$y,y'\in Z$ and $y',y''\in Z'$. But then $Z\cap Z'\ne\emptyset$ and therefore $Z=Z'$. Hence $yRy''$ and 
$R$ is transitive. It follows from the definition that $x,y\in Z\in X$ iff $x(Y/X)y$.

\medskip\noindent
{\bf Exercise~\ref{ex_dualita} }Let $R$, $A\ne\emptyset$ and $B$ be as stated. We prove the first equality.
We know that $C= A/R$ is a~partition of $A$, thus $S=A/C$ is an equivalence 
relation on $A$. We show that $S=R$. Let $a,b\in A$. Then $aSb$ iff there is a~block $D\in C$ such that $a,b\in 
D$. As we know from Exercise~\ref{ex_bloky}, $a,b$ lie in a~common block of $C$ iff $aRb$. Hence $S=R$.

We prove the second equality. We know that $S=A/B$ is an equivalence relation on $A$. Thus $C=A/S$ is 
a~partition of $A$. We show that $C=B$. Again, $a,b\in A$ lie in a~common block of $C$ iff $aSc$. This 
is the case iff $a,b$ lie in a~common block of $B$. Hence $C=B$.

\medskip\noindent
{\bf Exercise~\ref{ex_onLO} }Suppose that $R$ is a~linear order and $a,b$ 
are two elements such $aRb$ and $bRa$. Transitivity gives $aRa$, 
which is a~contradiction.

\medskip\noindent
{\bf Exercise~\ref{ex_tranzNeos} }Always $a\le a$ because
$a=a$. The transitivity of $\le$ follows from the transitivity of $<$. 
If $a,b$ are two distinct elements, then $a<b$ or $b<a$ by trichotomy. 
Thus $\le$ is dichotomic. If $a,b$ are two elements such that $a\le b$
and $b\le a$, then $a\ne b$ would give
$a<b$ and $b<a$, in contradiction with the previous exercise. Thus 
$a=b$ and $\le$ is weakly asymmetric.

\medskip\noindent
{\bf Exercise~\ref{ex_jenJedna} }Neither $a<b\wedge a=b$ nor $b<a\wedge a=b$ holds because $<$ is 
irreflexive. Nor $a<b\wedge b<a$ holds by Exercise~\ref{ex_onLO}.

\medskip\noindent
{\bf Exercise~\ref{ex_dalsiUl} }This follows from the trichotomy of $<$: 
if it is not the case that $a<a'$, then $a=a'$ or $a'<a$.

\medskip\noindent
{\bf Exercise~\ref{ex_minimax}. }When $m$ and $n$ are maxima of $B$, then both $n\le m$ and $m\le n$. Hence 
$m=n$. Same for minima.

\medskip\noindent
{\bf Exercise~\ref{ex_minimax2} }Proceed by induction on the size of the set. 

\medskip\noindent
{\bf Exercise~\ref{ex_infSup} }This follows from uniqueness of maxima and minima.

\medskip\noindent
{\bf Exercise~\ref{ex_dokAprVl} }Let $c=\sup(B)$. Then $c$ is an upper bound of $B$. 
Let $a<c$. Since $c$ is the minimum upper bound of $B$, the element $a$ is not an upper bound of $B$ and there 
exists the stated $b$. In the other way, let $c$ have the stated properties.
They say that $c$ is the smallest upper bound of $B$, so that $c=\min(H(B))=\sup(B)$.

Similarly one proves the equivalence that $c\in A$ is an infimum of $B$ iff $c\le b$ for every 
$b\in B$ $\&$ for every 
$a\in A$ with $c<a$ there is a~$b\in B$ such that $b<a$.

\medskip\noindent
{\bf Exercise~\ref{ex_naWO} }See Theorem~\ref{thm_wellOrdNat}.

\medskip\noindent
{\bf Exercise~\ref{ex_naWO2} }$\min(\Z)$ does not exist because $n-1<n$ for every $n\in\Z$.

\medskip\noindent
{\bf Exercise~\ref{ex_naWO3} }If a~linear order $\langle A,<\rangle$ is a~well ordering then such elements $a_n$ do not exist because $\{a_n\cc\;n\in\N\}$ would not 
have minimum. Suppose that $\langle A,<\rangle$ is not a~well ordering. We take 
any $a_1\in A$. Then we take any 
$a_2\in A$ such that $a_2<a_1$. Then any $a_3\in A$ such that $a_3<a_2$, 
and continue in this manner. In this way, we define the elements $a_n$. We 
used AC.\index{axiom!of choice, AC} 

\medskip\noindent
{\bf Exercise~\ref{ex_cvicOrder1} }See  Exercise~\ref{ex_tranzNeos}.

\medskip\noindent
{\bf Exercise~\ref{ex_cvicOrder2} }The irreflexivity and transitivity of strict inclusion are easy to check.

\medskip\noindent
\centerline{{\bf 2 From ZFC formulas to natural numbers}}

\medskip\noindent
{\bf Exercise~\ref{ex_012} }It follows from the definition that the sets $0$, 
$1$, and $2$ are elements of every inductive set.   

\medskip\noindent
{\bf Exercise~\ref{ex_minAndplus} }It follows from the previous definition 
and the axiom of extensionality.

\medskip\noindent
{\bf Exercise~\ref{ex_irrefNat} }Axiom~\ref{axi_found} of foundation forbids any set 
$x$ such that $x\in x$.

\medskip\noindent
{\bf Exercise~\ref{ex_onIndDef} }$f(n)=\lfloor n/2\rfloor$.

\medskip\noindent
{\bf Exercise~\ref{ex_nonstrOA} }If $a=b$, then $a+c=b+c$.

\medskip\noindent
{\bf Exercise~\ref{ex_semiRinf} }Additions of $1_S$ to $0_S<1_S$ 
yield infinitely many elements $0_S<1_S<1_S+1_S<\cdots$. 

\medskip\noindent
{\bf Exercise~\ref{ex_nextL} }Let $S$ be ordered and $a\in S$. Let $a\cdot0<0$. By 
adding $a\cdot0$ we obtain $a\cdot0+a\cdot0<a\cdot0$. But 
$a\cdot0+a\cdot0=a\cdot(0+0)=a\cdot0$ and we have the contradiction that 
$a\cdot0<a\cdot0$. $a\cdot0>0$ is similarly impossible. It follows from the trichotomy of $<$ that $a\cdot0=0$. If $S$ is simple, the equality $a\cdot0=0$ is easily proven by induction starting from $0\cdot0=0$.

\medskip\noindent
{\bf Exercise~\ref{ex_example} }It is clear that with the new element 
$\infty$, both operations remain commutative, and that $0_S$ and $1_S$ 
remain neutral. The three identities\,---\,for the associativity of 
$+$ and $\cdot$, and for the distributive law\,---\,remain valid: if 
$\infty$ appears, then both sides are $\infty$.

\medskip\noindent
{\bf Exercise~\ref{ex_finSiSemi} }For example, addition and multiplication modulo $m\in\N$.

\medskip\noindent
{\bf Exercise~\ref{ex_lastLem} }Let $m\in n$. Then $m\cup\{m\}\in n\cup\{n\}$ by part~6 of 
Proposition~\ref{prop_elOfOm}.

\medskip\noindent
{\bf Exercise~\ref{ex_jesteCvic} }We only need to prove that $f_S\cc\omega\to S$ is a~bijection; 
the three required properties of $f_S$ were proven in the three 
preceding lemmas. Lemma~\ref{lem_lema3} shows that 
$f_S$ is injective. To show that $f_S[\omega]=S$, we use that $S$ is 
simple and observe that $f_S[\omega]$
($\sus S$) contains $0_S$ and is closed to adding $1_S$. Indeed, $0_S=f_S(0)$ and if $x=f_S(m)$ for some $m\in\omega$, then $x\oplus1_S=f_S(m)\oplus1_S=f_S(m+1)$
by the definition of $f_S$.

\medskip\noindent
{\bf Exercise~\ref{ex_jesteDk} }Let $k=n-m\in\omega$ be the unique number
such that $m+k=n$. By the distributivity of $\cdot$ we have $l\cdot m+l\cdot k=l\cdot(m+k)=l\cdot n$. 
Thus $l\cdot k=l\cdot n-l\cdot m$.

\medskip\noindent
\centerline{{\bf 3 From natural numbers to integers}}

\medskip\noindent
{\bf Exercise~\ref{ex_disjUni} }$\langle0,m\rangle=n$ with 
$m,n\in\omega$ never holds because all elements of the nonempty left-hand side set are nonempty, but 
always $n=\emptyset$ or $\emptyset\in n$. 

\medskip\noindent
{\bf Exercise~\ref{ex_linOrdZ} }Irreflexivity is trivial from the 
definition. Let $k<l<m$ be in $\Z$. If all three are in $-\N$ or in 
$\omega$, then $k<m$ follows from the transitivity of $<$ in $\N_0$. The remaining 
case is $k\in-\N$ and $m\in\omega$, and then  $k<m$ follows from the 
definition. If $m,l\in\Z$ and are distinct, then both are in $-\N$ or 
in $\omega$ and are compared by the trichotomy of $<$ in $\N_0$, or one 
is in $-\N$ and the other in $\omega$, and $m$ and $l$ are compared from the definition. 
Thus $<$ is trichotomic. That $<$ extends the linear order in $\N_0$ is obvious.
This linear order is not a~well ordering.

\medskip\noindent
{\bf Exercise~\ref{ex_oneExZ} }$-2$, $-5$, $-6$, $-3<-2$, and $-10<1$.

\medskip\noindent
{\bf Exercise~\ref{ex_AddInvUni} }If
$\al+\be=0_R$ and $\al+\ga=0_R$ then the neutrality of $0_R$, and 
associativity and commutativity of $+$ give that $\ga=
(\al+\be)+\ga=(\al+\ga)+\be=\be$.

\medskip\noindent
{\bf Exercise~\ref{ex_2ndOAx} }If $1_R<0_R$, then $0_R<-1_R$ by 
additive inverses and the first order axiom. But then $0_R=0_R\cdot(-1_R)
<(-1_R)\cdot(-1_R)=1_R$ by the ring second order axiom, which is 
a~contradiction.

\medskip\noindent
{\bf Exercise~\ref{ex_onRiHomo} }Let $f\cc R\to R'$ be a~ring 
homomorphism. Then $0_{R'}+1_{R'}=1_{R'}=f(1_R)=f(0_R+
1_R)=f(0_R)+f(1_R)=f(0_R)+1_{R'}$ and we have equality 
$0_{R'}+1_{R'}=f(0_R)+1_{R'}$. Adding $-1_{R'}$ to it, we obtain 
$0_{R'}=f(0_R)$. Thus if $a+b=0_R$, we apply $f$ and obtain 
$f(a)+f(b)=f(0_R)=0_{R'}$. By Exercise~\ref{ex_AddInvUni}, 
$f(-a)=f(b)=-f(a)$.

\medskip\noindent
{\bf Exercise~\ref{ex_notDoma} }For every $m\in\N$, addition and 
multiplication modulo $m$ on $\omega$ yield a~finite simple ring $R_m$.

\medskip\noindent
{\bf Exercise~\ref{ex_idenProp} }Let $a,b,c\in\omega$. We used that (i) $(a+b)-c=(a-c)+b$,  (ii) $a-(b+c)=(a-b)-c$, 
(iii) $-(a-b)=(-a)+b$, (iv) $a-(b-c)=(a-b)+c$, and (v) $(a-b)-c=(a-c)-b$, 
provided that the involved differences are defined. 

\medskip\noindent
{\bf Exercise~\ref{ex_asAnExe} }We only need to prove that $f_R\cc\Z\to R$ is 
a~bijection; the three required properties of $f_R$ are proven in the 
three previous lemmas. Lemma~\ref{lem_lemaR3} shows that 
$f_R$ is injective. To prove that $f_R$ is surjective, it suffices to 
show, since $R$ is simple, that $f_R[\Z]$ contains $0_R$ and is 
closed to adding and subtracting $1_R$. The former is clear, $0_R=f_R(0)$. Let $x=f_R(m)$ for some $m\in\Z$. Then, by Lemma~\ref{lem_lemaR1},
$$
x\oplus\ominus 1_R=f_R(m)\oplus\ominus 1_R=
f_R(m\pm1)\in f_R[\Z]
$$
and we are done.

\medskip\noindent
{\bf Exercise~\ref{ex_rozd1} }Reflexivity and symmetry are 
obvious. Transitivity: if $m\ominus n\sim m'\ominus n'$ and $m'\ominus 
n'\sim m''\ominus n''$, then $m+n'=n+m'$ and $m'+n''=n'+m''$; hence 
$m+n'+m'+n''=n+m'+n'+m''$, which implies by canceling $n'+m'$ that $m+n''=n+m''$ 
and $m\ominus n\sim m''\ominus n''$.

\medskip\noindent
{\bf Exercise~\ref{ex_rozd2} }We show it here in detail only for $\cdot$. 
Let $m\ominus n\sim m'\ominus n'$ and
$k\ominus l\sim k'\ominus l'$. Then $m+n'=n+m'$, $k+l'=l+k'$, and 
$m\ominus n\cdot k\ominus l=mk+nl\ominus ml+nk$, where the variables may be primed. Then 
$mk+nl+m'l'+n'k'=ml+nk+m'k'+n'l'$
$\iff$ $m(k-l)+n(l-k)+m'(l'-k')+n'(k'-l')=0$ $\iff$ $(m+n')(k-l)+(n+m')(l-k)=0$ $\iff$ $(m+n'-n-m')(k-l)=0$ $\iff$ $0\cdot(k-l)=0$ $\iff$ T.

\medskip\noindent
{\bf Exercise~\ref{ex_rozd4} }Irreflexivity and transitivity are 
easy to see. As for trichotomy, let $k\ominus l\not\sim m\ominus n$ be 
two distinct difference integers. Thus $k+n\ne l+m$. Now $k+n<l+m$ gives
$k\ominus l<m\ominus n$, and $k+n>l+m$ gives $k\ominus l>m\ominus n$.

\medskip\noindent
{\bf Exercise~\ref{ex_rozd3} }We prove the associativity of 
$\cdot$. $(m\ominus n\cdot m'\ominus n')\cdot m''\ominus n''$ is 
$(mm'+nn')m''+(mn'+nm')n''\ominus(mn'+nm')m''+
(mm'+nn')n''$. On the other hand, $m\ominus n\cdot (m'\ominus n'\cdot m''\ominus n'')$ is $(m'm''+n'n'')m+
(m'n''+n'm'')n\ominus(m'n''+n'm'')m+(m'm''+n'n'')n$. By inspection, on 
the left-hand sides of the two $\ominus$, we have the same quadruple 
of cubic monomials, and the same holds for the right-hand sides of the two $\ominus$. 

\medskip\noindent
{\bf Exercise~\ref{ex_rozd5} }Let $m\ominus n$, $m'\ominus n'$, and 
$m''\ominus n''$ be three difference integers with $m\ominus n<m'\ominus n'$, so that $m+n'<n+m'$. Thus 
$m+m''+n'+n''<n+n''+m'+m''$ and 
$m\ominus n+m''\ominus n''<m'\ominus n'+m''\ominus n''$, which proves the 
first order axiom. 

Suppose that $m''\ominus n''>0\ominus 0$, so that 
$m''-n''>0$. Then $m\ominus n\cdot m''\ominus n''$ is $A=mm''+nn''\ominus 
mn''+nm''$, and $m'\ominus n'\cdot m''\ominus n''$ is 
$B=m'm''+n'n''\ominus m'n''+n'm''$. Now $A<B$ $\iff$ 
$mm''+nn''+m'n''+n'm''<mn''+nm''+m'm''+n'n''$ $\iff$ $m(m''-n'')+n(n''-m'')+m'(n''-m'')+n'(m''-n'')<0$ $\iff$
$(m+n'-n-m')(m''-n'')<0$ $\iff$ T, because $m+n'-n-m'<0$ and $m''-n''>0$.
This proves the second order axiom.

\medskip\noindent
\centerline{{\bf 4 From integers to fractions}}

\medskip\noindent
{\bf Exercise~\ref{ex_shodnProQ} }Reflexivity and symmetry are clear. We prove transitivity. Let $a/b\sim c/d$ and $c/d\sim e/f$. Thus 
$ad=bc$ and $cf=de$. But then $adf=bcf=bde$ and $adf=bde$. The
$d\ne0$  can be canceled ($\Z$ is a~domain)
and $af=be$. Hence $a/b\sim e/f$.

\medskip\noindent
{\bf Exercise~\ref{ex_posDeno} }$\frac{m}{n}\sim\frac{-m}{-n}$.

\medskip\noindent
{\bf Exercise~\ref{ex_lowTer0} }Let $a/b$ and $c/d$ be two equivalent 
protofractions in lowest terms, so that $ad=bc$. Since $(a,b)=1$ and $(c,d)=1$, 
we see that $a\,|\,c$ and $d\,|\,b$. Thus $c=ac'$ and $b=db'$ with $c',b'\in\Z$. 
We conclude that $1=b'c'$ and $b'=c'=\pm1$. Since $b,d>0$, the sign is $+$. Hence $a=c$ and $b=d$.

\medskip\noindent
{\bf Exercise~\ref{ex_lowTer} }Let $\al$ be a~fraction and $a/b\in\al$ be any 
protofraction in it.  By Exercise~\ref{ex_posDeno}, we may assume 
that $b>0$. Let $a=da'$ and $b=bb'$, where $d=(a,b)$. Then $a/b\sim a'/b'$ and
$a'/b'$ is in lowest terms. 

\medskip\noindent
{\bf Exercise~\ref{ex_hezkeCvic} }$\frac{0}{1}\not\sim\frac{1}{1}$ because $0\cdot1\ne1\cdot1$.

\medskip\noindent
{\bf Exercise~\ref{ex_ordFie2} }Let $a,b,c\in F\setminus\{0_F\}$ be such that $a\cdot b=1_F=a\cdot c$. Multiplying the last equality by $b$ we get that $b=(b\cdot a)\cdot c=1_F\cdot c=c$.

\medskip\noindent
{\bf Exercise~\ref{ex_neutrQ} }The element $0$ is neutral to $+$ as 
$\frac{0}{1}+\frac{a}{b}=\frac{0b+1a}{1b}=\frac{a}{b}$. Similarly, $1$ is 
neutral to $\cdot$ as $\frac{1}{1}\cdot\frac{a}{b}=\frac{1a}{1b}=
\frac{a}{b}$.

\medskip\noindent
{\bf Exercise~\ref{ex_opetLast} }We only need to prove that $f_F\cc\Q\to F$ is a~bijection; 
the three required properties of $f_F$ are proven in the three 
previous lemmas. The map 
$f_F$ is injective by Lemma~\ref{lem_tretifF}. We prove that 
$f_F$ is surjective. It 
suffices to 
show, since $F$ is simple, that $f_F[\Q]$ contains $0_F$ and is 
closed under adding and subtracting $1_F$, and under division. Clearly, $0_F=f_F(0)$. Let $x=f_F(\al)$ and $y=f_F(\be)\ne0_F$ for some $\al,\be\in\Q$ with $\be\ne0$. 
Then, by Lemma~\ref{lem_prvnifF},
$$
x\oplus\ominus 1_F=f_F(\al)\oplus\ominus 1_F=
f_F(\al\pm1)\in f_F[\Q]\,.
$$
By Lemma~\ref{lem_druhefF}, 
$$
x\oslash y=f_F(\al)\oslash f_F(\be)=f_F(\al/\be)\in f_F[\Q]\,.
$$
Thus $f_F[\Q]=F$ and we are done.

\medskip\noindent
\centerline{{\bf 5 From fractions to extended reals}}

\medskip\noindent
{\bf Exercise~\ref{ex_topBase} }(i) The union of sets in the system is the base set. (ii) The intersection of any two sets
in the system is a~union of some members of the system. Open intervals in $\langle 
X,<\rangle$ clearly satisfy (i) because every element $x\in X$ lies in an open 
interval. As for (ii), the intersection of any two open intervals is even an open interval.

\medskip\noindent

\medskip\noindent
$a\cdot b=0_F$ and $a,b\in F$ are nonzero, then the previous exercise yields $1_F=b^{-1}\cdot a^{-1}\cdot a\cdot b=
b^{-1}\cdot a^{-1}\cdot0_F=0_F$, which is a~contradiction.

\medskip\noindent
follows from the previous exercise. For example, $(a\cdot b)\cdot(a^{-1}\cdot b^{-1})=\ds=1_F$, so that $(a\cdot b)^{-1}=a^{-1}\cdot b^{-1}$.

\medskip\noindent
first two properties are 
straightforward. We prove the last two. The former follows from the 
identity $x\cdot x^{-1}=1_F$ by the multiplicativity of $|\cdot|$. The 
latter follows by applying the TI on the sum $x=(x+y)+(-y)$. 

\medskip\noindent
$d(x,y)=|x-y|\ge0_F$ is trivial, as 
is $d(x,y)=d(y,x)$ because $|z|=|-z|$. TI $d(x,y)\le d(x,z)+d(z,y)$ 
follows from applying TI for $|\cdot|$ to the sum $x-y=(x-z)+(z-y)$.

\medskip\noindent
$1_F^{-1}=1_F$, the neutrality of $1_F$, and $(a\cdot b^{-1})\cdot(b\cdot a^{-1})=1_F$. 

\medskip\noindent
$(a\cdot b^{-1})\cdot(c\cdot d^{-1})^{-1}=a\cdot b^{-1}\cdot c^{-1}\cdot d=(a\cdot d)\cdot(b\cdot c)^{-1}$.

\medskip\noindent
The implication $1\Rightarrow2$ follows  
from the inequality $x\le|x|$. The opposite implication follows from the 
inequality $|x|\le\max(-x,x)$. Suppose that $1$ holds and that $x\in 
F$ is nonzero. We take $m\in\omega$ such that $1_F/|x|\le f_F(m)$, and set 
$n=m+1$; then $3$ holds. Let $3$ hold. Then $|x|< f_F(n)$ and $1$ 
holds.

\medskip\noindent
Suppose that $\emptyset\ne X\sus F$ 
and that $X$ is lower-bounded. Then $\inf(X)=-\sup(-X)$.

\medskip\noindent
Let $(a_n)\sus F$ be Cauchy. We take 
$m\in\N$ such that if $n,n'\ge m$, then $|a_n-a_{n'}|\le1_F$. Let 
$b=\max(\{|a_1|,\ds,|a_m|\})$. Then for every $n\in\N$ we have, by TI, 
that $|a_n|\le1_F+b$.

\medskip\noindent
Let $\lim a_n=a$ and $\lim a_n=b$, where 
$a,b\in F$ are two elements. If $a\ne b$, we set $e=|a-b|/3_F$ ($>0_F$) and 
take $n_0\in\N$ such that if $n\ge n_0$, then $|a_n-a|,|a_n-b|\le e$. Then 
TI yields the contradiction that $e\le e(2_F/3_F)$. 

\medskip\noindent
See Theorem~\ref{thm_CauchyPodm}.

\medskip\noindent
Let $\al\in\Q$ and $p_{\al}=
\frac{m}{n}\in\al$ have $n>0$ and minimum $|m|+|n|$. It follows that 
$p_{\al}$ is in lowest terms. The function $\al\mapsto p_{\al}$ is 
the desired bijection. It is 
bijective and unique because two different elements in $Z_0$ are $\not\sim$. We prove it. If $\frac{k}{l},\frac{m}{n}\in 
Z_0$ are in lowest terms 
and $\frac{k}{l}\sim\frac{m}{n}$ then $kn=ml$. Thus any prime power dividing $k$ divides
$m$ and vice versa. The fundamental theorem of arithmetic gives $k=m$. 
Thus $l=n$.

\medskip\noindent
It is easy to check that $F$ is injective and preserves $+$, $\cdot$, and additive inverses. Also, $F(1)=1_{\Q}$. 

\medskip\noindent
If $m\cdot n=0$, then $F(m)\cdot F(n)=0_{\Q}$ in the field $\Q$ and we use 

\medskip\noindent
Because for every $m\in\Z$ we have $(-m)^2=m^2$.

\medskip\noindent
If $m^2$ is even, then $m$ is even 
because the squares of odd numbers are odd.

\medskip\noindent
For fractions $s,r>0$ with $s^2>2$, the inequality $(s-r)^2>2$ 
holds if $s^2-2>2sr-r^2$. Thus, for example, if $0<r<\frac{s^2-2}{2s}$. For fractions $r,s$ with $s^2<2$, $s>0$ 
and $r\in(0,1)$ (then $r^2<r$), the inequality $(s+r)^2<2$ holds if $2sr+r^2<2-s^2$. Thus, for 
example, if $0<r<\min(\{1,\frac{2-s^2}{2s+1}\})$.

\medskip\noindent
If $x>0_F$ but $x^{-1}<0_F$, then using 
the second order axiom we get $1_F=x^{-1}\odot x<0_F$, which 
contradicts Exercise~\ref{ex_2ndOAx}. The second claim follows from the 
second order axiom at once.

\medskip\noindent
It is. Two different isomorphisms would 
yield a~non-identical automorphism of the field $\Q$. From it we would 
obtain a~non-identical automorphism of the domain $\Z$, which cannot exist by induction.   

\medskip\noindent
Reflexivity and symmetry of $\sim$ are trivial. The 
transitivity easily follows from the triangle inequality.

\medskip\noindent
Every integer is HMC, even HF. Every 
fraction is a~countable set of $2$-tuples of integers and is therefore 
HMC. Since the set of fractions is countable, it is HMC.

\medskip\noindent
This follows from the previous exercise.

\medskip\noindent
The relation is irreflexive by 
definition. It is transitive because set inclusion is transitive. Let 
$X,Y\sus\Q$ be two distinct cuts. If neither $X\sus Y$ nor $Y\sus X$ 
holds, we take distinct fractions $\al\in X\setminus Y$ and $\be\in 
Y\setminus X$. Then neither $\al<\be$ nor $\al>\be$ holds, by 
property (ii) of cuts. This is a~contradiction.  

\medskip\noindent
Replace the sum of cuts with their intersection, and argue as for suprema.

\medskip\noindent
Since every Cauchy sequence is bounded, 
$\Phi(a_n)$ has property (i) of cuts. Properties (ii) and (iii) follow 
easily from the definition of $\Phi$. 

\medskip\noindent
Argue as in part~1.

\medskip\noindent
For example, $(1,0,0,0,\ds)\odot(0,1,0,0,\ds)=\overline{0}$.

\medskip\noindent
If $m$ and $n$ are so large that $|a_m-a_n|\le 
e'$ and $|a_m|,|a_n|\ge e$, then 
$|b_m-b_n|=|a_m^{-1}|\cdot|a_n^{-1}|\cdot|a_m-a_n|\le e'e^{-2}$.

\medskip\noindent
$\Phi(\al,\al,\ds)=F(\al)$.

\medskip\noindent
$a^2=b^2$ is equivalent with $(a-b)(a+b)=0_F$.

\medskip\noindent
If $a=0_F$, then $a^2=0_F\ge0_F$. If 
$a>0_F$, then $a^2>0_F\cdot a=0_F$ by the second order axiom. If $a<0_F$, 
then $0_F<-a$ by the first order axiom
and $0_F=0_F\cdot(-a)<(-a)\cdot(-a)=a^2$.

\medskip\noindent
With the help of 
in the complete ordered field $F$ we can define $a<_F b$ $\iff$ $\exists c\cc\,(c\ne0_F\wedge c^2=b-a)$.   

\medskip\noindent
Suppose that $X$ is dense in $F$ and 
that $\al\in F$. For every $n\in\N$, we select an $x_n\in X$ such that 
$|\al-x_n|\le1_F/f_F(n)$. Since $F$ is Archimedean, it follows that $\lim 
x_n=\al$. Suppose that $X$ is not dense in $F$. Then there exist $\al,e\in F$ 
with $e>0_F$ such that $|\al-x|\ge e$ for every $x\in X$. Then, clearly, 
there is no sequence $(x_n)\sus X$ such that $\lim x_n=\al$.

\medskip\noindent
See Theorem~\ref{thm_ari_lim}. 

\medskip\noindent
See Theorem~\ref{thm_ari_lim}.

\medskip\noindent
See Theorem~\ref{thm_limAuspo}.

\medskip\noindent
A~number $\ga\in\R$ is in the 
intersection iff $a_n\le\ga\le b_n$ for every $n$. This is true iff $\ga$ 
is both an upper bound of all $a_n$, and a~lower bound of all $b_n$. And 
this is true iff $\ga\in[\al,\be]$.

\medskip\noindent
Under this assumption, $\al=\be$ ($=c$).

\medskip\noindent
If $\R$ is countable, then there is a~bijection $f\cc\omega\to\R$. We have the bijection 
$g$ given by $\N\ni n\mapsto g(n)=n-1\in\omega$. Then the composite map $f(g)$ is a~surjection from $\N$ to $\R$.

\medskip\noindent
We divide $I_n$ into thirds and take a~third not containing $a_{n+1}$.

\medskip\noindent
The union of all open intervals is $X$. If $I$ and $J$ are open intervals, then so is $I\cap J$.  

\medskip\noindent
For $y\in X$ we take the identical sequence 
$(y,y,\ds)$. For $y\in Y\setminus X$ we take the subsequence because of which 
$y$ was added to $X$.

\medskip\noindent
By adding the mentioned comparisons for 
infinities, the irreflexivity, transitivity, and trichotomy are 
preserved.

\medskip\noindent
For example, $\langle[0,1],<\rangle$.

\medskip\noindent
10. The mentioned inequality turns
into the equality $C=C$ after the addition of $C$. It is easy to check 
that in other cases, the inequality is preserved (if the involved expressions 
are defined). 11. The mentioned inequality turns into the equality $C=C$ after multiplying by $C$. It is easy to check 
that in other cases, the inequality is preserved (if the involved expressions 
are defined).  

\medskip\noindent
$-\infty$, undefined, $-\infty$, and 
$-\infty$,

\medskip\noindent
$+\infty$, $+\infty$, $0$, and undefined.  

\medskip\noindent
\centerline{{\bf 7 Limits of real sequences}}

\medskip\noindent

\medskip\noindent
{\bf Exercise~\ref{ex_pocSnek} }$+\infty$, $-\infty$, $-\infty$, and undefined.

\medskip\noindent
{\bf Exercise~\ref{ex_supMinusNek} }Just $\emptyset$ and $\{-\infty\}$.

\medskip\noindent
{\bf Exercise~\ref{ex_nekdVlOkoli} }Neighborhoods of points and infinities are intervals.

\medskip\noindent
{\bf Exercise~\ref{ex_ulohaNaokoli1} }For $A,B\in\R$ we may take any 
$\ep<\frac{B-A}{2}$. If $A=-\infty$ and $B=+\infty$, we may take any 
$\ep$. If $A=-\infty$ and $B\in\R$, we may take any $\ep<\frac{1}{|B|+1}$; 
similarly if $A\in\R$ and $B=+\infty$.

\medskip\noindent
{\bf Exercise~\ref{ex_ulohaNaokol2} }This is immediate from the definition of neighborhoods.

\medskip\noindent
{\bf Exercise~\ref{ex_ulohaNaokoli3} }This is again immediate from the definition of neighborhoods.

\medskip\noindent
{\bf Exercise~\ref{ex_blowUp} }Suppose that $L=\lim a_n$. Then for any given $\ep$ there is a~$k$ such 
that for every $n\ge k$ we have $a_n
\in U(L,\ep)$. Then for every $n\ge\sum_{i=1}^k 
m_i$ we have $b_n\in U(L,\ep)$. Hence $\lim b_n=L$. In the proof of the opposite implication, we use that $(a_n)$ is a~subsequence of $(b_n)$.

\medskip\noindent
{\bf Exercise~\ref{ex_limPlInf} }$\lim a_n=+\infty$ $\iff$ for every 
$c>0$ there exists $n_0$ such that if $n\ge n_0$, then $a_n\ge c$. The 
proof is similar to the $-\infty$ case.

\medskip\noindent
{\bf Exercise~\ref{ex_jenomCvic} }If $(a_n)$ and $(b_n)$ satisfy $a_n=b_n$ 
for every $n\ge m$ and $\lim a_n=L$, then also $\lim b_n=L$. This follows 
from the fact that the index $n_0$ in the definition of limit can be taken
$\ge m$.

\medskip\noindent
{\bf Exercise~\ref{ex_robRob} }We only prove 3, the proofs for 1 and 2 
are similar. Suppose that $X$ is as 
stated, that $a_n\ne b_n$ for only finitely many $n$ 
and that $(a_n)\in\bigcap X$. Thus for every $Y\in X$ we have $(a_n)\in 
Y$. Since every $Y$ is robust, 
also $(b_n)\in Y$. Hence $(b_n)\in\bigcap X$. So $\bigcap X$ 
is robust. 

\medskip\noindent
{\bf Exercise~\ref{ex_zmenaKonmnoha} }The second and fifth properties are not robust; the others are robust.

\medskip\noindent
{\bf Exercise~\ref{ex_optVal} }Let $k\in\N$, 
$\ep=\frac{1}{k}$, and $n_0<\lceil\frac{1}{\ep}\rceil+1=k+1$. Thus $n_0\le k$ and $\frac{1}{n_0}\not\in U(0,\ep)=(-\frac{1}{k},\frac{1}{k})$ because $\frac{1}{n_0}\ge\frac{1}{k}$.

\medskip\noindent
{\bf Exercise~\ref{ex_priklNalimitu} }Clearly,  
$\frac{\sqrt[3]{n}-\sqrt{n}}{\sqrt[4]{n}}=\frac{n^{-1/6}-1}{n^{-1/4}}\to\frac{0-1}{0^+}=\frac{-1}{0^+}=-\infty$, 
due to positivity of $n^{-1/4}$.

\medskip\noindent
{\bf Exercise~\ref{ex_binomVeta} }The coefficient of the monomial $a^j b^{n-j}$ 
is the number of ways to obtain it: we choose $j$ factors
$a+b$ in the product $(a+b)^n$ from which we pick the number $a$,
and pick $b$ from the remaining $n-j$ factors. There are
$\binom{n}{j}$ ways to do it because there are $\binom{n}{j}$ $j$-element
subsets of $[n]$.

\medskip\noindent
{\bf Exercise~\ref{ex_vysvPodp} }It follows by negating the limit $n^{1/n}\to1$.

\medskip\noindent
{\bf Exercise~\ref{ex_explIndet} }In any ordered field $F$, the second order 
axiom and the transitivity of $\ge$ imply that if $a\ge b\ge0_F$ and $c\ge d\ge0_F$, then $ac\ge bd$. 

\medskip\noindent
{\bf Exercise~\ref{ex_O_mono} }It is not.

\medskip\noindent
{\bf Exercise~\ref{ex_O_mono1} }One replaces the given sequence $a_1\ge a_2\ge\ds$ with the sequence $-a_1\le 
-a_2\le\cdots$.

\medskip\noindent
{\bf Exercise~\ref{ex_jesteDodskAL} }This follows from the identity $|a_n-0|=||a_n|-0|$ ($=|a_n|$).

\medskip\noindent
{\bf Exercise~\ref{ex_preqRefTranz} }Reflexivity follows from setting 
$m_n=n$. Transitivity is clear when one views subsequences as 
obtained by omitting terms in original sequences. 

\medskip\noindent
{\bf Exercise~\ref{ex_podpNotWeakAntirefl} }For example,  $(0,1,0,1,\ds)$ and $(1,0,1,0,\ds)$.

\medskip\noindent
Hint: if $\iota_n=1\,2\,\ds\,n$ is 
the identical permutation, then the desired 
$m$-tuple is $\iota_{k-1}\ominus\iota_{k-1}\ominus\ds\ominus\iota_{k-1}$
with $l-1$ copies of $\iota_{k-1}$.

\medskip\noindent
The generalization says that in any 
linear order $(X,<)$ every sequence $(a_n)$ 
has a~monotone subsequence. The same proof works.

\medskip\noindent
{\bf Exercise~\ref{ex_slabaPodposl} }Let $(b_n)$ and $(a_n)$ be the 
stated sequences and let an $\ep$ be given. Thus there is an $n_0$ 
such that $n\ge n_0$ $\Rightarrow$ $a_n\in U(L,\ep)$. Then there is 
an $n_1$ such that $n\ge n_1$ $\Rightarrow$ 
$m_n\ge n_0$. Then for every $n\ge n_1$ we have that $b_n=a_{m_n}\in U(L,\ep)$ and $\lim b_n=L$.

\medskip\noindent
{\bf Exercise~\ref{ex_oSlPodp} }It is easy to see that the sequence $(m_n)\sus\N$ witnessing 
$(b_n)\preceq^*(a_n)$ has an increasing subsequence.

\medskip\noindent
{\bf Exercise~\ref{ex_disjointSub} }Let $B_0$ and $B_1$ be the supports 
of the two subsequences in part~1 of Theorem~\ref{thm_oPodposl}. Since
their limits differ, both sets $B_0\setminus B_1$ and $B_1\setminus 
B_0$ are infinite.

\medskip\noindent
{\bf Exercise~\ref{ex_3rdDual} }If a~sequence converges then every 
subsequence of it has the same finite limit and the right-hand side 
does not hold. If a~sequence diverges then it has no limit or an 
infinite limit. In the former case it has, by part~1 of 
Theorem~\ref{thm_oPodposl}, two subsequences with different limits. 
In the latter case the sequences itself has limit $\pm\infty$.

\medskip\noindent
{\bf Exercise~\ref{ex_whyAssu} }We can take only every other element of $B_0$.

\medskip\noindent
{\bf Exercise~\ref{ex_partBound} }Now in the proof of Theorem~\ref{thm_infinManyBl} the
subsequence corresponding to $B_0$ converges. 

\medskip\noindent
{\bf Exercise~\ref{ex_naParti} }Using Exercise~\ref{ex_disjointSub} 
we take in $(a_n)$ two disjoint subsequences $(b_n)$ and $(c_n)$ with different limits, and then proceed as in the proof of Theorem~\ref{thm_infinManyBl}.

\medskip\noindent
{\bf Exercise~\ref{ex_existHB} }By 

\medskip\noindent
{\bf Exercise~\ref{ex_lininftau} }For every $n\ge2$, it is true that $\tau(n)\ge2$ because $1$
and $n$ always divide $n$. For infinitely many $n$, namely for the prime numbers, the equality holds.

\medskip\noindent
{\bf Exercise~\ref{ex_casti3a4} }One can reduce parts 3 and 4 to parts 1 and 2 by means of the identity $\liminf a_n=-\limsup(-a_n)$. 

\medskip\noindent
{\bf Exercise~\ref{ex_naliminf1} }We let $m$ run in $\N$ and compose 
$(a_n)$ of segments $S_m$: $(a_n)=S_1S_2\ds$, where $S_m$ runs through the numbers
$-m$, $-m+\frac{1}{m}$, $-m+\frac{2}{m}$, $\ds$, $m$. 

\medskip\noindent
{\bf Exercise~\ref{ex_naliminf2} }No, it is not because $L(a_n)\cap\R$ is always a~closed set.

\medskip\noindent
{\bf Exercise~\ref{ex_naliminf3} }$L(a_n)=
\{0,+\infty\}$.

\medskip\noindent
{\bf Exercise~\ref{ex_onMonot} }The implication 
$\Leftarrow$ is trivial. The implication 
$\Rightarrow$ follows from the transitivity of $\le$. Similarly for the rest of the exercise. 

\medskip\noindent
{\bf Exercise~\ref{ex_ekvivOmez} }If such $c$ exists then for every $n$ it 
holds that $-c\le a_n\le c$ and $(a_n)$ is bounded both from below 
and from above.  Suppose that $(a_n)$ is bounded by the definition, so that 
$d\le a_n\le c$ for every $n$ 
and some numbers $d$ and $c$. Then $|a_n|\le\max(|d|,|c|)$ for every $n$.

\medskip\noindent
{\bf Exercise~\ref{ex_petRobustnich} }The last three concerning boundedness. 

\medskip\noindent
{\bf Exercise~\ref{ex_uloNaMonPos}. 
}Suppose that $(a_n)$ weakly increases for every $n\ge m$ and that 
$b_n=a_n$ 
for every $n\ge n_0$. Then $(b_n)$ weakly increases for every 
$n\ge\max(m,n_0)$. Same for weakly decreasing tails. 

\medskip\noindent
{\bf Exercise~\ref{ex_monoJekvazim} }For example, if $(a_n)$ weakly decreases, then for every $n$, the implication $a_m>a_n$
$\Rightarrow$ $m<n$ holds. Hence, $(a_n)$ goes down. Similarly, for weakly increasing sequences.

\medskip\noindent
{\bf Exercise~\ref{ex_oKvazimon1} }For example, $(1,0,2,1,3,2,4,3,5,\ds)$ goes up, 
but no tail is monotone.

\medskip\noindent
{\bf Exercise~\ref{ex_oKvazimon2} }A~sequence $(a_n)\sus\R$ is quasi-monotone iff
$$
\big(\forall\,l\,\exists\,m\cc n\ge m\Rightarrow a_n\ge a_l\big)\vee\big(\forall\,l\,\exists\,m\cc n\ge m\Rightarrow a_n\le a_l\big)\;.
$$

\noindent
{\bf Exercise~\ref{ex_robuKVMo} }Suppose for example that $(a_n)$ goes up starting from $n=m$ and that $b_n=a_n$ 
for every $n\ge n_0$. Then $(b_n)$ goes up from $n=\max(m,n_0)$.

\medskip\noindent
{\bf Exercise~\ref{ex_part3} }We reduce it to part~2 by switching from the sequence $(a_n)$ to the sequence $(-a_n)$.

\medskip\noindent
{\bf Exercise~\ref{ex_verzeBWthh} }It is easy to see that the limit $c$ of this subsequence satisfies the inequalities $a\le c\le b$.

\medskip\noindent
{\bf Exercise~\ref{ex_cauchyRobust} }Suppose that $(a_n)$ is Cauchy 
and $(b_n)$ is such that $b_n=a_n$ for $n\ge n_0$. If for a~given 
$\ep$ for every $m,n\ge n_1$ it holds that $|a_m-a_n|\le\ep$, then 
for every $m,n\ge\max(n_0,n_1)$ it holds that $|b_m-b_n|\le\ep$.
Hence $(b_n)$ is Cauchy.

\medskip\noindent
{\bf Exercise~\ref{ex_CauchyJeOmez} }Let $(a_n)$ be Cauchy. Then there 
is an $n_0$ such that for every $m,n\ge n_0$ one has that $|a_m-
a_n|\le 1$. By TI it holds for every 
$n$ that $|a_n|$ $\le$ $1+\max(\{|a_1|,\ds,|a_{n_0}|\})$. Hence $(a_n)$ is bounded.

\medskip\noindent
{\bf Exercise~\ref{ex_CauchyVeQ} }Take for example the sequence $(1,1.4,1.41,1.414,\ds)$ of truncations of the decimal expansion of $\sqrt{2}$.

\medskip\noindent
{\bf Exercise~\ref{ex_kdeuplnost} }We used it via the application of the B.--W. theorem whose proof uses the 
theorem on limits of monotone sequences. This theorem requires the 
existence of suprema and infima of sets of real numbers.

\medskip\noindent
{\bf Exercise~\ref{ex_fekete} }black

\medskip\noindent
{\bf Exercise~\ref{ex_fekete2} }Back then, there was no AI, no Internet, no mobile phones, no computers that one could use 
in everyday 
life $\ds$ In the fall of that Orwell year, I started my study of 
mathematics (and ``cybernetics'') at 
Charles University in Prague $\ds$
Now my parents are both dead, but back then my father was 
a~young, energetic man of $50$ and my mother was a~nice young lady of 
$47$ $\ds$ Do I wish I could go back to 1984? On one hand, yes, on 
the other hand, I see how stupid I was back then $\ds$

\medskip\noindent
{\bf Exercise~\ref{ex_naAFeLe} }We have $n=m+m+\ds+m+l$, with $k+1$ 
summands.

\medskip\noindent
{\bf Exercise~\ref{ex_whySuperad} }To prove the first lemma, see the 
properties of the exponential function, especially continuity. We 
prove the second lemma. If $A=+\infty$, then there is a~sequence 
$(x_n)\sus X$ with $\lim x_n=+\infty$. Then 
$\lim\exp(x_n)=+\infty=\sup(Y)$, so 
that $\sup(Y)=\exp(A)$. Let $A\in\R$. Since $x\le A$ for every $x\in X$, 
also $y=\exp x\le\exp A$ for every $y\in Y$. If $a<\exp A$, then $\log a<A$ and there is 
$x\in X$ such that $\log a<x\le A$. Then $a<\exp x\le\exp A$ and $\exp x\in Y$. Thus $\sup(Y)=\exp A$.

\medskip\noindent
{\bf Exercise~\ref{ex_exIsDef} } Let $u\ne\emptyset$. If 
$v\in[n]^*$ is $r$-sparse and $|v|\ge(\binom{n}{r}(|u|-1)+1)r$ 
then by the pigeonhole principle you can find $|u|$ intervals $\ds 
I_1\ds I_2\ds\,\ds I_{|u|}\ds$ in $v$ such that $|I_i|=r$ and that 
all $I_i$ use the same $r$-element set of letters. Taking  
appropriate terms from these intervals, one per each $I_i$, you 
easily build a~copy of $u$ in $v$. 

\medskip\noindent
{\bf Exercise~\ref{ex_abab} }We prove the upper bound $\mathrm{ex}(abab,n)\le 2n-1$ by induction on 
$n$. For $n=1$ it holds. Suppose that $v\in[n]^*$ with $n\ge2$ is a~$2$-sparse (there is 
no immediate repetition in $v$) word not containing $abab$. 
Two closest occurrences of some $j\in[n]$ (if there is any 
repetition) show that some $i\in[n]$ occurs in $u$ only 
once. Deleting this occurrence and possibly one more term of $v$, we get a~$2$-sparse word 
$v'\in([n]\setminus\{i\})^*$. Clearly, $abab\not\preceq v'$. Hence 
$|v|\le|v'|+2\le2(n-1)-1+2=2n-1$. On the other hand, words $12\ds(n-1)n(n-1)\ds 21$ show that the upper bound is tight.

\medskip\noindent
{\bf Exercise~\ref{ex_trivSzem} }These are trivialities.  

\medskip\noindent
{\bf Exercise~\ref{ex_saw1} }An automorphism sending 
a~vertex $u$ to a~vertex $v$ induces a~bijection from the set 
of edges incident with $u$ to the 
the set of edges incident with $v$.

\medskip\noindent
{\bf Exercise~\ref{ex_saw2} }In a~path starting at $v$, we 
always have a~finite number of possibilities for the next edge. An automorphism sending 
a~vertex $u$ to a~vertex $v$ induces a~bijection from the set 
of paths starting at $u$ to the set of paths starting at $v$. 

\medskip\noindent
{\bf Exercise~\ref{ex_bouNumPath} }For the first edge of a~path we have $r$ possibilities. For every 
next edge we have only at most $r-1$ possibilities because we cannot 
backtrack.

\medskip\noindent
{\bf Exercise~\ref{ex_hexag} }Let $G=\langle\Z^2,E\rangle$ be the 
$4$-regular graph of the square lattice. Let $E'=\{\{\langle a,b\rangle,\langle 
a+1,b\rangle\}\cc\;a,b\in\Z, \text{ $a$ and $b$ have the same 
parity}\}$ ($\sus E$). Then $H=G\setminus E'=\langle\Z^2,E\setminus 
E'\rangle$ is the hexagonal tiling graph. It is clearly $3$-regular 
and if $\overline{u},\overline{v}\in\Z^2$ 
are two vertices, then the shift by the vector $\overline{v}-
\overline{u}$, followed if needed by a~vertical reflection, is an 
automorphism sending $\overline{u}$ to $\overline{v}$.

\medskip\noindent
{\bf Exercise~\ref{ex_bouMeaCat} }Classically, $C_n$ is the number 
of noncrossing matchings $\langle[2n],E\rangle$.

\medskip\noindent
{\bf Exercise~\ref{ex_howReco} }We identify $e_3$ and $e_4$ as the unique two edges crossing the 
bisector of the segment $\langle 2m,0\rangle\langle 2m+1,0\rangle$.

\medskip\noindent
{\bf Exercise~\ref{ex_numsPerms} }A~word $u\in S_m$ has $m$ 
possibilities for the first letter, (for every choice of the first letter it has) $m-1$ possibilities for the second letter, $\ds$

\medskip\noindent
{\bf Exercise~\ref{ex_whereFails} }If $p_1$, $p_2$, $\ds$, $p_k$ are 
forbidden permutations, then it may happen that for some $i$ and 
$j$, $p_i$ is not $\ominus$-irreducible and $p_j$ is not 
$\oplus$-irreducible. Then we do not know if the $\ominus$- and 
$\oplus$-sum of two permutations avoiding every permutation $p_i$ still has this 
property.

\medskip\noindent
{\bf Exercise~\ref{ex_obmeTroj} }Since $|-b|=|b|$, it suffices to prove the first 
inequality. We apply to $a=(a+b)+(-b)$ the standard TI and rearrange the result.

\medskip\noindent
{\bf Exercise~\ref{ex_dodatek1} }1. Let $|a_n|\le d$ for every $n$, $L=-\infty$ and a~$c<0$ be given. It is clear that for every large
$n$ one has $b_n\le c-d$. Thus for every large $n$ we have that $a_n+b_n\le d+c-d=c$. Hence $a_n+b_n\to-\infty$. The case $L=+\infty$ is similar.

2. Let $|a_n|\le d$ for every $n$, $b_n\to0$ and an $\ep$ be given. Clearly, for every large $n$  we have $|b_n|\le\frac{\ep}{d}$.
So for every large $n$, we have $|a_nb_n|\le d\cdot\frac{\ep}{d}=\ep$. Hence $a_nb_n\to0$. 

3. Let $a_n$, $c$, $L=+\infty$ and $b_n$ be as stated and let a~$d>0$ be given. One has  for every large $n$ that $b_n\ge\frac{d}{c}$.
So for every large $n$ it holds that $a_nb_n\ge c\cdot\frac{d}{c}=d$. Hence $a_nb_n\to+\infty=L$. The other case is similar.

4. Let $|a_n|\le d$ for every $n$, $b_n\to\pm\infty$ and an $\ep$ be given. For every large $n$, we have $|b_n|\ge\frac{d}{\ep}$. So for every large $n$ we have $\big|
\frac{a_n}{b_n}\big|=|a_n|\cdot\frac{1}{|b_n|}\le d\frac{1}{d/\ep}=\ep$. Hence 
$\frac{a_n}{b_n}\to0$. 

5. Let $a_n$, $c$ and $b_n$ be as stated and 
a~$d>0$ be given. For every large $n$, one has $0<b_n\le\frac{c}{d}$.
So for every large $n$ we have $\frac{a_n}{b_n}\ge\frac{c}{c/d}=d$. Hence
$\frac{a_n}{b_n}\to+\infty$. 

6. Let $a_n$, $c$, $L=-\infty$ and $b_n$ be as stated and let a~$d<0$ 
be given. For every large $n$, one has $b_n\le dc$.
So for every large $n$ we have $\frac{b_n}{a_n}\le\frac{dc}{c}=d$. 
Hence $\frac{b_n}{a_n}\to-\infty=L$. The other case is treated similarly.

\medskip\noindent
{\bf Exercise~\ref{ex_whatPo} }Then if $\lim\frac{a_n}{b_n}$ exists, it equals $0$.

\medskip\noindent
{\bf Exercise~\ref{ex_AGnerov} }This inequality is equivalent to the inequality $(\sqrt{a}-\sqrt{b})^2\ge0$.

\medskip\noindent
{\bf Exercise~\ref{ex_OnFrec} }Instead of $k$-variable function 
$f$, we take $k+1$-variable function, where the new variable is reserved for the index $n$.

\medskip\noindent
{\bf Exercise~\ref{ex_fibon} }$F_n\ge n-1$ for every $n\in\N$.

\medskip\noindent
{\bf Exercise~\ref{ex_epDeDef} }For every $\ep$ there is $\de$ 
such that $f[U(b_1,\de)\times\ds\times 
U(b_1,\de)]\sus U(f(\overline{b}),\ep)$. The proof 
is similar to that of Proposition~\ref{prop_spojVbLimitou}.

\medskip\noindent
{\bf Exercise~\ref{ex_proGen} }The proof is similar to that of Proposition~\ref{prop_limReku}.

\medskip\noindent
{\bf Exercise~\ref{ex_fibo1} }We need to show that $F_{n+1}F_{n-1}-
F_n^2=F_{n+1}^2-F_{n+2}F_n$. This is equivalent to
$F_{n+1}(F_{n-1}-F_{n+1})=F_n(F_n-F_{n+2})$. Since the first bracket 
is $-F_n$ and the second one is $-F_{n+1}$, the identity holds. 

\medskip\noindent
{\bf Exercise~\ref{ex_uloNaRek} }This infinite continued fraction leads to the $f$-recurrent sequence $(a_n)=(0,1,\frac{1}{2},\frac{2}{3},\ds)$
for the function $f=f(x)=\frac{1}{1+x}$. Similarly to Proposition~\ref{prop_fibRat}, we show that $\lim a_n=\frac{1}{2}(\sqrt{5}-1)=1/\phi$.

\medskip\noindent
{\bf Exercise~\ref{ex_procJesiln} }The set of pairs of sequences 
$\{\langle(a_n),(b_n)\rangle:\;\exists\,n_0\,\forall\,m,n\ge n_0\cc a_m<b_n\}$ is a~proper subset of the set
$\{\langle(a_n),(b_n)\rangle:\;\exists\,n_0\,\forall\,n\ge n_0\cc a_n<b_n\}$. 

\medskip\noindent
{\bf Exercise~\ref{ex_muzePrejit}. }For example, $(a_n)=(\frac{1}{n})$ and $(b_n)=(0,0,\ds)$.

\medskip\noindent
{\bf Exercise~\ref{ex_daleZesil}. }Let $(a_n)$, $(b_n)$, $K$ and $L$ be as stated. 
We take a~number $c$ such that $K<c<L$. By Exercise~\ref{ex_nekdVlOkoli} 
there is an $\ep$ such that $U(K,\ep)<U(c,\ep)<U(L,\ep)$. Then we take any two numbers 
$a,b\in U(c,\ep)$ such that $a<b$. For every large $m$ and $n$ we have that $a_m\in U(K,\ep)$ 
and $b_n\in U(L,\ep)$. Hence $a_m\le a$ and $b\le b_n$.

Reversal of this implication is: if for every $n_0$ and every real numbers $a<b$ there exist $m$ and $n$ with $m,n\ge n_0$ such that $a_m>a$ or $b_n<b$, then $K\ge L$.

\medskip\noindent
{\bf Exercise~\ref{ex_konecInterv} }Any singleton $\{a\}$ is such an interval.

\medskip\noindent
{\bf Exercise~\ref{ex_Jedenstr} }Let $\lim a_n=-\infty$, $b_n\le a_n$ for every large $n$ 
and let a~$c<0$ be given. Then for every large $n$ one has that $b_n\le a_n\le c$. Thus $b_n\le c$ 
and $\lim b_n=-\infty$. The case of the limit $+\infty$ is similar.

\medskip\noindent
\centerline{{\bf 3 Series}}

\medskip\noindent
{\bf Exercise~\ref{ex_vysvEpdel3} }The first absolute value is at most
$\frac{\ep}{3}$ because $X'\sus X_0$. The $i$-th absolute value in
the sum is at most
$2^{-i}\frac{\ep}{3}$ because $Z_i'\sus Z_i''$. The third absolute
value is at most $\frac{\ep}{3}$ because $Y'\sus \{Z_1,\ds,Z_n\}$.

\medskip\noindent
{\bf Exercise~\ref{ex_congIsER} }Reflexivity is witnessed by the 
identity bijection, symmetry by the inverse bijection and transitivity 
by the composition of two bijections.

\medskip\noindent
{\bf Exercise~\ref{ex_congAKser} }Let
$R=\sum_{x\in X}r(x)$ and $R'=\sum_{x\in Y}s(x)$ be congruent 
AK series. If $X$ and $Y$ are finite, the equality of their sums 
is trivial. Suppose that they are infinite and that $f\cc X\to Y$ is 
a~bijection proving that $R\sim R'$.
Let $g\cc\N\to X$ be any bijection.
Then $S(R')=\lim\sum_{i=1}^n s(f(g)(i))=
\lim\sum_{i=1}^n r(g(i))=S(R)$ 
because $f(g)$ is a~bijection from $\N$ to $Y$.

\medskip\noindent
{\bf Exercise~\ref{ex_jesteTriEp} }The first absolute value is at most
$\frac{\ep}{3}$ because $Z'\sus W$. The second absolute value is at most
$\frac{\ep}{3}$ because $X'\sus X''$. The third absolute
value is at most $\frac{\ep}{3}$ because $Y'\sus Y''$.

\medskip\noindent
{\bf Exercise~\ref{ex_congAndBiSu} }Let $Q=\sum_{x\in X}r(x)$, 
$Q'=\sum_{x\in X'}r'(x)$, $R=\sum_{x\in Y}s(x)$ and 
$R'=\sum_{x\in Y'}s'(x)$ be as stated and let $f\cc X\to X'$, $g\cc 
Y\to Y'$ be bijections witnessing that $Q\sim Q'$ and $R\sim R'$. Thus
for every $x\in X$ and $y\in Y$ we have that $r(x)=r'(f(x))$ and 
$s(y)=s'(g(y))$. Let $Z\equiv X\times\{0\}\cup Y\times\{1\}$ and
$W\equiv X'\times\{0\}\cup Y'\times\{1\}$. We consider $Q+R=\sum_{z\in 
Z}t(z)$ and $Q'+R'=\sum_{z\in W}t'(z)$. We define the bijection
$h\cc Z\to W$ by $h(z)\equiv(f(x),0)$ if $z=(x,0)$, 
and by $h(z)\equiv(g(y),1)$ if $z=(y,1)$. Then it follows that for
every $z\in Z$ we have $t(z)=t'(h(z))$. Hence $Q+R\sim 
Q'+R'$.

\medskip\noindent
{\bf Exercise~\ref{ex_jesteeTriEp} }The first bound by
$\frac{\ep}{3}$ comes from the inclusion $Z\sus X''\times Y''$. The 
second bound is
$|(r+\de)\theta|\le\frac{\ep}{3}$. The third bound is $|\de s|\le\frac{\ep}{3}$.

\medskip\noindent
{\bf Exercise~\ref{ex_congAndProd} }We proceed as in the previous 
exercise, with the modifications that $Z\equiv X\times Y$, $W\equiv 
X'\times Y'$ 
and that the bijection $h\cc X\times Y\to X'\times Y'$ is given by
$h((x,y))\equiv(f(x),g(y))$. Then for every $(x,y)\in X\times Y$ we 
have that $r(x)s(y)=r'(f(x))s'(g(y))$ because
$r(x)=r'(f(x))$ and $s(y)=s'(g(y))$. Hence $Q\cdot R\sim Q'\cdot R'$. 

\medskip\noindent
{\bf Exercise~\ref{ex_jakJeBij} }For $R=\sum_{x\in X}r(x)$, $R'=\sum_{y\in 
Y}s(y)$ and $R''=\sum_{z\in Z}t(z)$ the bijection sends $(x,0)$ to 
$((x,0),0)$, $((y,0),1)$ to $((y,1),0)$ and $((z,1),1)$ to 
$(z,1)$.

\medskip\noindent
{\bf Exercise~\ref{ex_jakJeBij2} }For $R$ and $R'$ as in the previous exercise the bijection sends $(x,y)$ to $(y,x)$.

\medskip\noindent
{\bf Exercise~\ref{ex_jakJeBij3} }For $R$, $R'$ and $R''$ as in Exercise~\ref{ex_jakJeBij} the bijection sends $(x,(y,z))$ to $((x,y),z)$.

\medskip\noindent
\centerline{{\bf 4 Infinite series.  Elementary functions}}

\medskip\noindent
{\bf Exercise~\ref{ex_remaiCase} }If $A=a>0$ and $B=+\infty$, we take the 
described series $\sum a_n$ and its subseries $\sum a_{2n-1}$. Similarly
if $A\le0$ and $B=\pm\infty$. 

\medskip\noindent
{\bf Exercise~\ref{ex_robKonvDiv} }As the next solution shows, if we change 
in a~series $\sum a_n$ finitely many summands and get $\sum a_n'$ then 
there is an $m$ and a~$c$ such that for every $n\ge m$ it holds that 
$s_n'=s_n+c$. Then finite $\lim s_n$ exists iff finite $\lim s_n'$ exists.

\medskip\noindent
{\bf Exercise~\ref{ex_zmeSou} }Let 
$\sum a_n$ and $\sum b_n$ be convergent series and let there be an 
$m$ such that $b_n=a_n$ for 
$n\ne m$ and $b_m=a_m+c$ with $c\ne0$. Let $(s_n)$ and $(t_n)$ be 
respective partial sums. Then $t_n=s_n$ for $n<m$ and $t_n=s_n+c$ 
for $n\ge m$, so that $\sum b_n=\lim t_n=c+\lim s_n=c+\sum a_n$. 

\medskip\noindent
{\bf Exercise~\ref{ex_onSubser} }It is, for example, the series $1-
1+\frac{1}{2}-\frac{1}{2}+\cdots$ with sum $0$ we encountered in the 
first lecture. Its subseries 
$1+\frac{1}{2}+\ds$ has sum $+\infty$. 

\medskip\noindent
{\bf Exercise~\ref{ex_onTails} }Let $(s_n)$ be partial sums of 
$a_1+a_2+\ds$, and $(t_n)$ be partial sums of $a_m+a_{m+1}+\ds$. 
Then $t_n=s_{n+m-1}-a_1-\ds-a_{m-1}$.  

\medskip\noindent
{\bf Exercise~\ref{ex_nezapScitance} }Partial summands weakly increase (resp. decrease) for $n\ge n_0$.

\medskip\noindent
{\bf Exercise~\ref{ex_souPlusneko} }We have $\lim s_n=\lim n=+\infty$.

\medskip\noindent
{\bf Exercise~\ref{ex_vzsvRovn} }The first equality follows from the
definition of partial sums, the second from Theorem~\ref{thm_ari_lim}, the third from the assumption and Proposition~\ref{prop_LimPodpo}, and
the fourth one is trivial.

\medskip\noindent
{\bf Exercise~\ref{ex_onLinComSer} }$\sum(-1)^{n+1}n^{-2}=\sum 
n^{-2}-2\sum (2n)^{-2}=(1-\frac{1}{2})\sum n^{-2}=\frac{1}{2}\sum n^{-2}=\frac{\pi^2}{12}$.

\medskip\noindent
{\bf Exercise~\ref{ex_prNaNekSou} }$\prod_{j=1}^n(1+\frac{1}{j})=\prod_{j=1}^n\frac{j+1}{j}=n+1$ and the infinite product is $+\infty$.

\medskip\noindent
{\bf Exercise~\ref{ex_statStan} }The standard NCC is that $\lim a_n=1$.

\medskip\noindent
{\bf Exercise~\ref{ex_prvPomoc} } Due to monotonicity, we have
$\lim a_n=L$. All subsequences have this limit and hence $L=+\infty$.

\medskip\noindent
{\bf Exercise~\ref{ex_druPomoc} }Let $(s_n)$ and $(t_n)$ be partial sums 
of both series. There is a~$c$ such that for every $n\ge n_0$ we have 
$s_n\ge c+t_n$. But $\lim (c+t_n)=c+\lim t_n=c+
(+\infty)=+\infty$, and the one-policeman theorem 
shows that also $\lim s_n=+\infty$.

\medskip\noindent
{\bf Exercise~\ref{ex_HnNenicele} }Let $n\ge 2$. We assume that 
$1+\frac{1}{2}+\ds+\frac{1}{n}=m$ ($\in\N$) and deduce a~contradiction. Following the hint 
we write every denominator $j=1$, $2$, $\ds$, $n$ in the form $j=a(j)\cdot 2^{b(j)}$ 
where $a(j)\in\N$ is odd and $b(j)\in\N_0$. For $j_0=2^k$, where 
$k\in\N$ is the largest number with $2^k\le n$, this expression takes the form $j_0=1\cdot 2^k$. For every 
$j\in [n]\setminus\{j_0\}$ it holds that $b(j)<k$. Hence $1+\frac{1}{2}+\ds+\frac{1}{n}=\frac{a+b}{a\cdot 
2^k}$, where $a\equiv a(1)a(2)\ds a(n)\in\N$ is an odd number
and $b\in\N$ is even, because it is the sum of $n-1$ even numbers. The 
numerator $a+b$ is therefore odd and the power $2^k\ge2$ in the denominator
cannot be canceled. Therefore we cannot have $\frac{a+b}{a\cdot 2^k}=m$. The same argument shows that for no $n\ge2$ we have that 
$h_n=\frac{k}{l}$ with odd $l$. 

\medskip\noindent
{\bf Exercise~\ref{ex_otProGamma} }???\,---\,presently no solution is known.

\medskip\noindent
{\bf Exercise~\ref{ex_druheTvrz} }One changes the series $\sum a_n$ to  $\sum(-a_n)$.

\medskip\noindent
{\bf Exercise~\ref{ex_geneLeiSer} }If a~series $a_1+a_2+\ds$ and 
$m\in\N$ are such that $|a_m|\ge|a_{m+1}|\ge\ds$, $\lim 
a_n=0$ and $(-1)^{n-m}a_n\ge0$ 
for every $n\ge m$, then $\sum a_n$
has the sum $s\in\R$ and we have $\sum_{i=1}^n a_i\ge 
s\ge\sum_{i=1}^{n+1} a_i$ for every $n\ge m$ such that $n-m$ is even. The proof is similar to that of Theorem~\ref{thm_altSer}. 

\medskip\noindent
{\bf Exercise~\ref{ex_naGroup} }Consider the series $\sum a_n\equiv1-1+1-1+1-1+\ds$ and the
sequence $S\equiv(2,2,\ds)$. The $S$-grouping $0+0+\ds$ has the sum $0$ but the sum $\sum a_n$ does not
exist.

\medskip\noindent
{\bf Exercise~\ref{ex_naGroup2} }Consider the series $\sum 
a_n\equiv1-1+\frac{1}{2}+\frac{1}{2}-\frac{1}{2}-\frac{1}{2}+
\frac{1}{3}+\frac{1}{3}+\frac{1}{3}-\frac{1}{3}-\frac{1}{3}-
\frac{1}{3}+\ds$ and the sequence $S\equiv(2,4,6,\ds)$. The 
$S$-grouping $0+0+\ds$ has the sum $0$ but the sum $\sum a_n$ does not
exist.

\medskip\noindent
{\bf Exercise~\ref{ex_prerovAKrad} }We have a~constant $c>0$ such that $\sum_{j=1}^n|a_j|\le c$ for every
$n$. Let $f\cc\N\to\N$ be any bijection and $n\in\N$. Then
$\sum_{j=1}^n|a_{f(j)}|\le \sum_{j=1}^N|a_j|\le c$, where
$N\in\N$ is such that $[N]\supset f[\,[n]\,]$. Thus the series
$\sum a_{f(n)}$ is abscon. 

\medskip\noindent
{\bf Exercise~\ref{ex_subserAKser} }Suppose that $\sum a_n$ is abscon. 
Then we have a~constant $c>0$ such that $\sum_{j=1}^n|a_j|\le c$ for 
every $n$. Let $B\sus\N$ be infinite, $f\cc\N\to B$ be the ordering of $B$ and $n\in\N$. Then
$\sum_{j=1}^n|a_{f(j)}|\le \sum_{j=1}^N|a_j|\le c$, where
$N\in\N$ is such that $[N]\supset f[\,[n]\,]$. Thus the subseries
$\sum a_{f(n)}$ is abscon.

\medskip\noindent
{\bf Exercise~\ref{ex_exOnCompCr} }We use the comparison criterion 
with $\sum b_n$ being the stated series and $\sum 
a_n\equiv\sum\frac{c}{n(n+1)}=\sum c(\frac{1}{n}-
\frac{1}{n+1})$, where $c$ is such that $|c_n|\le c$ for every $n$. 

\medskip\noindent
{\bf Exercise~\ref{ex_onCauPro} }The coefficient of $x^n$ in $F(x)G(x)$ 
is $\sum_{j=0}^na_j b_{n-j}$. 

\medskip\noindent
{\bf Exercise~\ref{ex_zobecGeom} }This follows from the equality $q^m+q^{m+1}+\ds=q^m\cdot(1+q
+\cdots)$.

\medskip\noindent
{\bf Exercise~\ref{ex_GeomRposlU} }Every converging one, so iff $q\in(-1,1)$.

\medskip\noindent
{\bf Exercise~\ref{ex_onRooTes} }$\sum a_n\equiv\sum\frac{1}{n}$ has 
$\lim a_n^{1/n}=1$ and $\sum a_n=+\infty$. $\sum a_n\equiv\sum\frac{1}{n^2}$ has 
also $\lim a_n^{1/n}=1$ but $\sum a_n$ converges.

\medskip\noindent
{\bf Exercise~\ref{ex_onRatTes} }Use the same two series from the
previous exercise.

\medskip\noindent
{\bf Exercise~\ref{ex_limsupRatio} }For example, $\frac{1}{2}+
\frac{1}{3}+\frac{1}{2^2}+
\frac{1}{3^2}+\frac{1}{2^3}+
\frac{1}{3^3}+\ds<+\infty$, but $\frac{1}{2}(\frac{3}{2})^n\to+\infty$.

\medskip\noindent
{\bf Exercise~\ref{zetaSmen1} }This is immediate from the divergence of the harmonic series. 

\medskip\noindent
{\bf Exercise~\ref{ex_aplCKT} }Iff $s>1$, again by CCC.

\medskip\noindent
{\bf Exercise~\ref{ex_zeta2Easy} }The series $\sum_{n=2}^{\infty}\frac{1}{n(n-
1)}$ converges because $\frac{1}{n(n-1)}=\frac{1}{n-1}-\frac{1}{n}$.

\medskip\noindent
{\bf Exercise~\ref{ex_eulPro} }By the Fundamental Theorem of
Arithmetic, every natural number is expressed in a~unique way as 
a~product of prime powers. Let $s>1$ and $n\in\N$. Using the formula for sum of a~geometric series, we have
$$
{\textstyle
0\le\prod_{j=1}^n\big(
1-\frac{1}{(p_j)^s}\big)^{-1}-
\sum_{m=1}^{p_n}\frac{1}{m^s}\le
\sum_{m>p_n}\frac{1}{m^s}\to0\ \ 
(n\to\infty)\,.
}
$$

\noindent
{\bf Exercise~\ref{ex_eulPro2} }No it could not, it is an infinite product 
$\prod_{n=1}^{\infty}a_n$ such that $a_n>1$ for every $n$. 

\medskip\noindent
{\bf Exercise~\ref{ex_spanTree} }It is easy to prove by induction that 
any tree with $k\in\N$ vertices has $k-1$ edges.

\medskip\noindent
{\bf Exercise~\ref{ex_naRiemPro} }This follows from definitions. 

\medskip\noindent
{\bf Exercise~\ref{ex_ButSee} }If $A=0$ then $B=-
\infty=\lim_{n\to\infty}\sum_{j=1}^n b_{\rho(j)}$ and 
$${\textstyle
\lim_{n\to\infty}\exp(\sum_{j=1}^n b_{\rho(j)})=0=|A|\,.
}
$$ 
Similarly, if $A=\pm\infty$ then $B=+\infty=\lim_{n\to\infty}\sum_{j=1}^n 
b_{\rho(j)}$ and the displayed limit equals $+\infty=|A|$.

\medskip\noindent
{\bf Exercise~\ref{ex_naRieThmPro} }One can show that $\prod_{n=1}^{\infty}a_n$ 
has an reordering with no infinite product iff $a_n\ne0$ for every $n$ and
the following two conditions hold. (i) If $a_n<0$ for infinitely many $n$, then
it is not the case that $\sum_{|a_n|<1}\log(|a_n|)=-\infty$ and 
$\sum_{|a_n|\ge 1}\log(|a_n|)<+\infty$. (ii) If $a_n<0$ for finitely many $n$, 
then $\sum\log(a_{z_n})=-\infty$ and $\sum\log(a_{k_n})=+\infty$. 

\medskip\noindent
{\bf Exercise~\ref{ex_limBodyLimitou} }Let $M$ and $A$ be as stated and let
part~1 hold, so that $A\in L(M)$ by the given definition. We chose for 
every $n$ an $a_n\in P(A,\frac{1}{n})\cap M$ and get a~sequence $(a_n)\sus 
M\setminus\{A\}$ such that $\lim a_n=A$. Hence part~2 holds. For every $m$ there is an $\ep$ such that $a_1,\ds,a_m\not\in U(A,\ep)$. So we can choose from $(a_n)$ an injective subsequence and part~3 holds. Suppose that part~3 holds and let $(b_n)\sus M$
be an injective sequence with $\lim b_n=A$. For given $n$ we have $b_m\in 
U(A,\frac{1}{n})$ for every large $m$. From these for only one $m$ it 
holds that $b_m=A$, hence $P(A,\frac{1}{n})\cap M\ne\emptyset$ 
and part~4 holds. It is clear that part~4 implies part~1.

\medskip\noindent
{\bf Exercise~\ref{ex_AjesteJedenPr} } If $M\sus\R$ is finite, then it is
bounded and neither $-\infty$ nor $+\infty$ is a~limit point of $M$. 
Also, for every $b$ there is a~$\de$ such that $P(b,\de)\cap M=\emptyset$.
Hence $L(M)=\emptyset$. If $M\sus\R$ is infinite, then we can choose an
injective sequence $(a_n)\sus M$. By part~1 of Theorem~\ref{thm_oPodposl}, $(a_n)$
has a~subsequence $(b_n)$ with 
$\lim b_n=L$. Then, by part~3 of Proposition~\ref{prop_OlimBodech}, 
$L\in L(M)$.

\medskip\noindent
{\bf Exercise~\ref{ex_vyjmutiBodu} }This is immediate from part~2 of Proposition~\ref{prop_OlimBodech}.

\medskip\noindent
{\bf Exercise~\ref{ex_kteryDalsiLimB}. }$L(\N)=\{+\infty\}$.

\medskip\noindent
{\bf Exercise~\ref{ex_uvedtePriklad} }For instance $A=0$, $M=\{\pm\frac{1}{n}\;|\;n\in\N\}$, 
$X=\{\frac{1}{n}\;|\;n\in\N\}$, $f=0$ on $M\setminus X$ and $f=1$ on $X$.

\medskip\noindent
{\bf Exercise~\ref{ex_AC} }We select an element from each set $\{x\in P(K,\frac{1}{n})\cap M\cc\;f(x)\not\in U(L,\ep)\}$, $n\in\N$.

\medskip\noindent
{\bf Exercise~\ref{ex_triLimity}. }1. Due to the transformation ($x<0$, hence $x=-|x|$) 
$\frac{x}{\sqrt{1+x^2}-1}=\frac{1}{-\sqrt{1/x^2+1}-1/|x|}$ we get for $x\to-\infty$ 
the limit $\frac{1}{-\sqrt{1/(+\infty)+1}-0}=-1$. 

2. The transformation $\frac{1}{\sqrt{1+x}-\sqrt{x}}=\sqrt{1+x}+\sqrt{x}$ gives for $x\to+\infty$ the limit 
$\sqrt{1+(+\infty)}+\sqrt{+\infty}=+\infty$.

3 a~4. These limits are trivial, the first does not exist and the second equals $0$.

\medskip\noindent
{\bf Exercise~\ref{ex_giveExam} }For example, $(a_n)=(n\sqrt{2})$ or, more generally, 
$(a_n)=(n\al)$ where $\al>0$ is any irrational number.

\medskip\noindent
{\bf Exercise~\ref{ex_CauchyError} }The equality $\frac{f(x)}{x}=\frac{f(h)}{x}+
(1-\frac{h}{x})(k+\al)$ ((2) on p.~36 of \cite{brad_sand}), where $x=h+n$ for large 
fixed $h>0$, $n\in\N$ and $-\ep<\al<\ep$, holds only on the discrete set 
$x\in\{h+n\cc\;n\in\N\}$, but we 
need it to hold on $x\in(x_0,+\infty)$. Cauchy regards ``$x$ as a~variable quantity 
which converges towards the limit $\infty$.'' It is not too surprising 
that \cite{brad_sand} has no clear and precise definition of the limit 
of a~function at a~(finite or infinite) point.

\medskip\noindent
{\bf Exercise~\ref{ex_priklNaKon} }As many as real numbers, $k_c\mapsto c$ is a~bijection from the set of constants to $\R$.

\medskip\noindent
{\bf Exercise~\ref{ex_radaJeAK} }We take an $m$ such that $m\ge 2|x|$. Then for every $n\ge m$ we have that 
$|\frac{x^n}{n!}|\le(\frac{|x|^m}{m!})\cdot(\frac{1}{2})^{n-
m}=\frac{1}{m!}(2|x|)^m\cdot(\frac{1}{2})^n$. Then 
we use geometric series.

\medskip\noindent
{\bf Exercise~\ref{ex_dok1az3} }1. $\exp 0=1$ is trivial and the rest 
follows from the exponential identity. 2. For $x<y$ we have that 
$\mathrm{e}^y-\mathrm{e}^x=\mathrm{e}^x(\mathrm{e}^{y-x}-1)>0$, due to the exponential identity. 3. For $x>n$ it holds that 
$\mathrm{e}^x>n$, so that $\lim_{x\to+\infty}\mathrm{e}^x=+\infty$. Also, $\lim_{x\to-\infty}\mathrm{e}^x=
\frac{1}{\lim_{x\to+\infty}\mathrm{e}^x}=\frac{1}{+\infty}=0$.

\medskip\noindent
{\bf Exercise~\ref{exeJeIrac} }For contradiction, let
$\sum_{j=0}^{\infty}\frac{1}{j!}=\frac{n}{m}$ with $n,m\in\N$. Following the hint
we get that $r\equiv\sum_{j>m}\frac{m!}{j!}=n\cdot(m-1)!-\ds\in\N$. This is 
impossible because $0<r\le\frac{1}{m+1}\sum_{j=0}^{\infty}\frac{1}{(m+2)^j}=\frac{m+2}{(m+1)^2}<1$.

\medskip\noindent
{\bf Exercise~\ref{ex_vlLogaritmu} }1. $\log 1=0$ follows from $\exp 
0=1$. By flipping the graph over the line $y=x$ we get from the increasing 
function $\exp x$ the increasing function $\log x$. For
$x,y>0$ we have by the exponential identity the equality
$\exp(\log x+\log y)=x\cdot y$, so that $\log x+\log y=\log(xy)$. 2. 
These limits are again obtained by flipping the graph of the exponential
over $y=x$. 3. This follows from part~4 of the previous proposition.

\medskip\noindent
{\bf Exercise~\ref{ex_takToNej} }No, $\sqrt{x}=x^{1/2}\ne\exp(\frac{1}{2}\log x)$. The former
function has at $0$ value $0$, but the latter function is not defined at $0$. 

\medskip\noindent
{\bf Exercise~\ref{ex_onSqRo} }For $a=0$ it is trivial as $\sqrt{0}=0$ and the only solution
of $x^2=0$ is $x=0$. If $a>0$ then $\sqrt{a}=a^{1/2}=\exp(\frac{1}{2}\log a)$ and 
$(\sqrt{a})^2=\exp(2\cdot\frac{1}{2}\log a)=a$ by the exponential identity. 
The factorization $x^2-a=(x+\sqrt{a})(x-\sqrt{a})$ then shows that $\pm\sqrt{a}$ are the only solutions of the equation $x^2=a$.

\medskip\noindent
{\bf Exercise~\ref{ex_jsouKomp} }We begin with $a^x$, $a>0$. For $x=0$ we have that $\exp(x\log a)=\exp 0=1$. For $x\in\N$ it holds that
$\exp(x\log a)=\exp(\log a+\ds+\log a)$, with $x$ factors $\log a$. By 
the exponential identity this equals $\exp(\log a)\cdot\ldots\cdot\exp(\log a)
=a\cdot\ldots\cdot a$, with $x$ factors $a$. For every $x\in\N$ it holds due to the exponential identity 
that $\exp((-x)\log a)=\frac{1}{\exp(x\log 
a)}$. Thus $a^x$ agrees with $x^m$. We continue with $x^b$. Let $b\in\N$, 
$x>0$. Then again 
$\exp(b\log x)=\exp(\log x+\ds+\log x)=\exp(\log x)\cdot\ldots\cdot\exp(\log x)
=x\cdot\ldots\cdot x$, with $b$ factors $x$. Also $0^b=0=0\cdot\ldots\cdot0$.  
Let $b=0$ and $x>0$. Then $x^b=1$.
Let $b\in\Z$ with $b<0$ and $x>0$.
Then we again get by the exponential identity that  $x^b=\frac{1}{x^{-b}}$.
Hence $x^b$ agrees with $x^m$. Finally, $0^x=0$ for $x\in\N$ also agrees with $x^m$.

\medskip\noindent
{\bf Exercise~\ref{ex_onE} }By Definition~\ref{def_aNaB} one has 
that $\mathrm{e}^x=
\exp(x\log\mathrm{e})$. Since $\mathrm{e}=\exp 1$ and $\log x$ is 
inverse to $\exp x$, it equals to $\exp(x\log(\exp 
1))=\exp(x\cdot1)=\exp x$.

\medskip\noindent
{\bf Exercise~\ref{ex_laterEx} }Note that $a_n=\mathrm{e}^{\log a_n}$.
For $A=0$ we set $b_n=-\frac{1}{\sqrt{\log(2+|a_n|)}}$. For $0<A<+\infty$
we set $b_n=\frac{\log A}{\log(2+|a_n|)}$.
For $A=+\infty$
we set $b_n=\frac{1}{\sqrt{\log(2+|a_n|)}}$.

\medskip\noindent
{\bf Exercise~\ref{ex_howExa} }$(-1)^2=1$ and $(-1)^1=1$ are
computed by the latter definition, $1^{\frac{1}{2}}=1$ by the former.

\medskip\noindent
{\bf Exercise~\ref{ex_Wilkie} }With $A\equiv1+x$, $B\equiv1+x+x^2$, 
$C\equiv1+x^3$ and $D\equiv1+x^2+x^4$, for which by the 
hint $AD=BC\equiv E$, we should show that $(A^y+B^y)^x\cdot(C^x+D^x)^y=
(A^x+B^x)^y\cdot(C^y+D^y)^x$. Equivalently, that $E^{xy}(1+(\frac{B}{A})^y)^x(1+
(\frac{C}{D})^x)^y=E^{xy}(1+(\frac{A}{B})^x)^y(1+(\frac{D}{C})^y)^x$. But 
this holds because $\frac{B}{A}=\frac{D}{C}$ and 
$\frac{C}{D}=\frac{A}{B}$, and multiplication is commutative.

\medskip\noindent
{\bf Exercise~\ref{ex_0na0jeNeur} }For $A=0$ we set $a_n\equiv\frac{1}{n^n}$ 
and $b_n\equiv\frac{1}{n}$. For $0<A<1$ we set $a_n\equiv A^n$ and $b_n\equiv\frac{1}{n}$. 
For $A=1$ we set $a_n=b_n\equiv\frac{1}{n}$. For $1<A<+\infty$ we set 
$a_n\equiv\frac{1}{A^n}$ and $b_n\equiv-\frac{1}{n}$. For 
$A=+\infty$ we set $a_n\equiv\frac{1}{n^n}$ 
and $b_n\equiv-\frac{1}{n}$. No, it could not, $a^b<0$ only for $b\in\Z\setminus\{0\}$.

\medskip\noindent
{\bf Exercise~\ref{ex_kosSinAK} }Proceed as in Exercise~\ref{ex_radaJeAK}.

\medskip\noindent
{\bf Exercise~\ref{ex_applLeibSer} }By the theorem, $\cos 1\ge 1-\frac{1^2}{2!}=\frac{1}{2}$ and $\cos 2\le1-\frac{2^2}{2!}+\frac{2^4}{4!}=-\frac{1}{3}$.

\medskip\noindent
{\bf Exercise~\ref{ex_odvodZvety} }1. The runner runs one lap in time $2\pi$ and gets in the same position. 
2. This is the behavior of the $y$-coordinate of the runner in the first 
quarter of the lap. 3. The track is symmetric according to the $y$-axis, 
and according to the origin $(0,0)$. 4. The counter-clockwise rotation of 
$S$ around the origin by $\frac{\pi}{2}$ is equivalent to the exchange of 
the coordinate axes. The second relation says that the points on $S$ have 
distance $1$ from the origin. 5. Search the Internet for pictures for ``geometric proof of summation formulae for sine and cosine''.

\medskip\noindent
{\bf Exercise~\ref{ex_limReIm} }This follows from the bounds $|a-a_n|,|b-b_n|\le|a+bi-z_n|$. 

\medskip\noindent
{\bf Exercise~\ref{ex_complExp} }It follows from the fact that the 
metric space $(\C,|\cdot|)$ is complete (every Cauchy sequence
in it converges) and the fact that the sequence $(\sum_{j=0}^n
\frac{z^j}{j!})$, $n\in\N$, is Cauchy, which is easy to prove by 
means of (complex) geometric series. 

\medskip\noindent
{\bf Exercise~\ref{ex_proLem} }Iterate the identity $u_1u_2\ds u_n-v_1v_2\ds 
v_n=(u_1-v_1)u_2\ds u_n+v_1(u_2\ds u_n-v_2\ds v_n)$.

\medskip\noindent
{\bf Exercise~\ref{ex_absHexp} }$\mathrm{e}^{a+bi}=\mathrm{e}^a\mathrm{e}^{bi}$. 

\medskip\noindent
{\bf Exercise~\ref{ex_boundExp} }Using geometric series we have $|\exp z-1-z|\le\frac{1}{2}\sum_{n\ge2}|z|^n=\frac{|z|^2}{2(1-|z|)}\le|z|^2$.

\medskip\noindent
{\bf Exercise~\ref{ex_infProExp} }For every $n\in\N$ we set $k_n=n$ and 
$a_{n,j}=\frac{a}{n}$, and use the theorem. 

\medskip\noindent
{\bf Exercise~\ref{ex_tanCot} }We know from the properties of cosine 
and sine what zeros they have.

\medskip\noindent
{\bf Exercise~\ref{ex_arcsin} }By the definition of inverses in 
the functional inverse of $\sin x\,|\,[-
\frac{\pi}{2},\frac{\pi}{2}]$ goes from $[-1,1]$ to $[-\frac{\pi}{2},\frac{\pi}{2}]$, but
$\arcsin\cc[-1,1]\to\R$. So arcsine is not equal to this inverse and is
only congruent to it, in the sense of congruence of functions in Definition~\ref{def_rovnostFci}.  Similarly for arccosine.

\medskip\noindent
{\bf Exercise~\ref{ex_arctan} }Between $\R$ and $(-\frac{\pi}{2},\frac{\pi}{2})$.

\medskip\noindent
{\bf Exercise~\ref{ex_aritmFunkci} }This follows from the commutativity, associativity
and distributivity of the operations $+$ and $\cdot$ on $\R$ and from the
fact that the operation of intersection of two sets enjoys 
these properties too: $M\cap N=N\cap M$, $(M\cap N)\cap P=M\cap(N\cap P)$ 
and $M\cap(N\cap P)=(M\cap N)\cap(M\cap P)$. In $\R$ the number
$0$, respectively $1$, is neutral to addition, respectively 
multiplication, and 
always $\R\cap M=M$. No function $f\in\mathcal{R}$ with $M(f)\ne\R$ 
has additive or multiplicative inverse.

\medskip\noindent
{\bf Exercise~\ref{ex_OdecFci} }The functions $f-g$ and $f+f_{-1}\cdot g$ 
have equal values and also equal definition domains: $M(f)\cap 
M(g)=M(f)\cap(\R\cap M(g))$.

\medskip\noindent
{\bf Exercise~\ref{ex_platiTo} }We set $g=h\equiv\emptyset_f$ and see that 
$\Rightarrow$ does not hold. Similarly, $f=g\equiv\emptyset_f$ show that $\Leftarrow$ does not hold. 

\medskip\noindent
{\bf Exercise~\ref{ex_onGenWor} }For every $i\in[n]$ the initial segment 
$\langle f_1,f_2,\ds,f_i\rangle$ is a~generating word of $f_i$.

\medskip\noindent
{\bf Exercise~\ref{ex_absHodn} }$|x|=(\mathrm{id}\cdot\mathrm{id})^{1/2}$.

\medskip\noindent
{\bf Exercise~\ref{ex_deleniNulou} }It is the empty function $\emptyset$.

\medskip\noindent
{\bf Exercise~\ref{ex_prazdnaElem} }Yes, it is, for example by the previous exercise. 

\medskip\noindent
{\bf Exercise~\ref{ex_pseudoInverz} }$g\equiv k_{-1}\cdot f$.

\medskip\noindent
{\bf Exercise~\ref{ex_EFsdefObZ} }For instance $f(x)\equiv\sqrt{\sin(\pi 
x)}+\sqrt{-\sin(\pi x)}$ and $g(x)\equiv\frac{1}{\sin(\pi/x)}$; here $\pi$ is $k_{\pi}$ and $x$ is the identity $\mathrm{id}$.

\medskip\noindent
{\bf Exercise~\ref{ex_arctann} }On $(-\frac{\pi}{2},\frac{\pi}{2})$ we 
have $\tan x=\frac{\sin x}{\sqrt{1-\sin^2 x}}$ and $\sin x=
\frac{\tan x}{\sqrt{1+\tan^2 x}}$. Hence $\arctan x=\arcsin(\sin(\arctan x))=
\arcsin\big(\frac{x}{\sqrt{1+x^2}}\big)$. 

\medskip\noindent
{\bf Exercise~\ref{ex_jesteRedu} }$x^b$ for $b\le 0$ and $b\in\N$ ($\sus\R$).

\medskip\noindent
{\bf Exercise~\ref{ex_EFjedenBod} }For instance $\arcsin x$, $|x|$ or $\arcsin(\sin x)$.

\medskip\noindent
{\bf Exercise~\ref{ex_zeroCan} }Suppose that $p$ has the 
canonical from $\sum_{j=0}^na_jx^j$ and that $n$ is minimum. We show by induction on $\deg p\equiv n$ 
that $|Z(p)|\le n$. For $n=0$ it holds, then $p=k_{a_0}$ with 
$a_0\ne0$ and $p$ has no zero. Let $n>0$. If $p$ has no zero, the 
inequality holds. Let $a\in Z(p)$. Then we divide with remainder the 
polynomial $p$ in the above canonical form by the polynomial 
$x-a=\mathrm{id}-k_a$ and get the expression $p=(x-a)q$, for some 
canonical polynomial $q$ with degree at most $n-1$. For every $b\ne a$ with $p(b)=0$ we have $q(b)=0$. By induction, $|Z(p)|=1+|Z(p)\setminus
\{a\}|\le1+|Z(q)|\le 1+\deg q\le 1+(n-1)=n$.

\medskip\noindent
{\bf Exercise~\ref{ex_rootFactor} }We proved it in the previous proof.

\medskip\noindent
{\bf Exercise~\ref{ex_sumCan} }Let $f=\sum_{j=0}^m a_jx^j$ and 
$g=\sum_{j=0}^n m_jx^j$. Then $f+g=\sum_{j=0}^p(a_j+b_j)x^j$, 
where $p\equiv\max(m,n)$, if $j>m$ then $a_j\equiv0$ and if $j>n$ then
$b_j\equiv0$, is a~canonical polynomial, provided that not all sums $a_j+b_j$ are zero. If they are all zero, then $f+g$ is the zero polynomial.

\medskip\noindent
{\bf Exercise~\ref{ex_prodCan} }Let $f=\sum_{j=0}^m a_jx^j$ and 
$g=\sum_{j=0}^n m_jx^j$. Then the product $fg=\sum_{j=0}^{m+n}
(\sum_{i=0}^j a_ib_{j-i})x^j$ is a~canonical polynomial. 

\medskip\noindent
{\bf Exercise~\ref{ex_diffCan} }Let $f=\sum_{j=0}^m a_jx^j$ and 
$g=\sum_{j=0}^n m_jx^j$ be such that $m\ne n$ or $a_j\ne b_j$ for some 
$j\le\min(m,n)$. It follows that $f-g=\sum_{j=0}^p(a_j-b_j)x^j$, 
where $p\equiv\max(m,n)$, if $j>m$ then $a_j\equiv0$ and if $j>n$ then
$b_j\equiv0$, is a~canonical polynomial because not all 
differences $a_j-b_j$ is zero.

\medskip\noindent
{\bf Exercise~\ref{ex_POLjeOI} 
}By Proposition~\ref{prop_semirFci} it only remains to prove that additive inverses exist and that 
the product of two nonzero polynomials is nonzero. The former
is provided by Exercise~\ref{ex_pseudoInverz} and
the latter by Exercise~\ref{ex_prodCan}.

\medskip\noindent
{\bf Exercise~\ref{ex_polyIsom} }The map that sends $k_0$ to the abstract zero polynomial, and a~canonical 
polynomial $f=\sum_{j=0}^n a_jx^j$ to the abstract polynomial 
$\sum_{j=0}^n a_jx^j$ in $\R[x]$, is the required isomorphism. 

\medskip\noindent
{\bf Exercise~\ref{ex_onRatFun} }This is a~trivial consequence of 
the definitions.

\medskip\noindent
{\bf Exercise~\ref{ex_prodRatio} }Now the domain of the left-hand 
side is $\big(M(f_1)\cap M(g_1)\setminus 
Z(g_1)\big)\cap\big(M(f_2)\cap M(g_2)\setminus Z(g_2)\big)$, and of 
the right-hand side it is $(M(f_1)\cap M(f_2))\cap(M(g_1)\cap 
M(g_2))\setminus (Z(g_1)\cup Z(g_2))$. They both are equal to
$M(f_1)\cap M(f_2)\cap M(g_1)\cap M(g_2)\setminus(Z(g_1)\cup Z(g_2))$. 

\medskip\noindent
{\bf Exercise~\ref{ex_naRatioRatios} }Let $f_1=f_2=g_1\equiv k_1$ 
and $g_2\equiv\mathrm{id}$. Then $\frac{f_1/g_1}{f_2/g_2}$ is
$\mathrm{id}\,|\,\R\setminus\{0\}$ but $\frac{f_1g_2}{f_2g_1}$ is $\mathrm{id}$.

\medskip\noindent
{\bf Exercise~\ref{ex_jetoSkutRE} } The reflexivity and symmetry of $\sim$ 
are clear. We prove transitivity. Let
$r\sim s$ and $s\sim t$. We take these rational functions in
canonical forms: $r=\frac{a}{b}$, $s=\frac{c}{d}$ and $t=\frac{e}{f}$. 
By the assumption $r=s$ on $M(r)\cap M(s)=\R\setminus Z(bd)$ and $s=t$ on 
$\R\setminus Z(df)$. Thus $r=t$ on $\R\setminus(Z(bd)\cup Z(df))=
(\R\setminus Z(bf))\setminus Z(d)=(M(r)\cap M(t))\setminus Z(d)$. By 
the continuity of $r$ and $t$ on $M(r)\cap M(t)\cap 
Z(d)$ the functions $r$ and $t$ 
are equal on the whole $M(r)\cap M(t)$. Hence $r\sim t$.

\medskip\noindent
{\bf Exercise~\ref{ex_RCjeTel} }Let $r,s,r',s'\in\mathrm{RAC}\setminus\{\emptyset\}$.
It is not hard to see, using continuity of rational functions, that $r\sim r'$ and $s\sim s'$ 
imply that also $r+s\sim r'+s'$ 
and $r\cdot s\sim r'\cdot s'$. Thus addition and multiplication of equivalence blocks of
rational functions is correctly defined. Neutrality of $[k_0]_{\sim}$ and $[k_1]_{\sim}$ to $+$ and $\cdot$, respectively, is immediate. The commutativity and associativity of addition and multiplication and the distributive
law in $\mathrm{RAC_{fi}}$ follow from the arithmetic in $\R$. Additive and multiplicative inverses require
equivalence blocks. Let
$r=p/q\in\mathrm{RAC}\setminus\{\emptyset\}$ 
have the canonical form $(p,q)$. Then
$[(-p)/q]_{\sim}$ is the additive inverse: the sum is $[k_0/q]_{\sim}=[k_0]_{\sim}$. If $p\ne k_0$ then $[q/p]_{\sim}$ is the 
multiplicative inverse: the product is
$[pq/qp]_{\sim}=[k_1]_{\sim}$.

\medskip\noindent
{\bf Exercise~\ref{ex_isomRacFun} }Recall that the abstract rational functions $\R(x)$ are
the equivalence classes 
$[\frac{p}{q}]_{\sim}$ where $p,q\in\R[x]$, $q\ne0$, are abstract (real) 
polynomials and the equivalence relation (congruence) $\sim$ on $\R[x]\times(\R[x]\setminus\{0\})$ is given by 
$\frac{p_1}{q_1}\sim
\frac{p_2}{q_2}$ iff $p_1q_2=p_2q_1$ (with the product and equality in $\R[x]$). It is not 
too hard to show that the map sending the zero rational function $[k_0]_{\sim}\in\mathrm{RAC}$ to the zero 
abstract rational function $[\frac{0}{1}]_{\sim}$ in $\R(x)$, and the nonzero rational function 
$[r]_{\sim}=[p/q]_{\sim}$, where 
$(p,q)\in(\mathrm{POL}\setminus\{k_0\})^2$ is a~canonical form of 
$r\in(\mathrm{RAC}\setminus\{\emptyset\})\setminus\{[k_0]_{\sim}\}$, to the nonzero abstract rational function
$[\frac{p}{q}]_{\sim}\in\R(x)\setminus\{0\}$, is an isomorphism of fields. The notation $p=\sum_{j=0}^n a_jx^j$, where $n\in\N_0$, $a_j\in\R$ and $a_n\ne0$, is 
interpreted in two ways: as the canonical form of a~nonzero polynomial in $\mathrm{POL}$, and
as a~nonzero abstract polynomial in $\R[x]$.

\medskip\noindent
\centerline{{\bf 5 Limits of functions. Asymptotic notation}}

\medskip\noindent
{\bf Exercise~\ref{ex_uloZpredn} }Let $b\in L^-(M)$. For every 
$n\in\N$ we chose from $P^-(b,\frac{1}{n})\cap M$ a~point $a_n$ and get
the required sequence $(a_n)$. If 
$b\not\in L^-(M)$ then for some 
$\ep$ we have that $P^-(b,\ep)\cap M=\emptyset$ and the required 
sequence does not exist. For $L^+(M)$, the proof is similar.

\medskip\noindent
{\bf Exercise~\ref{ex_naJednostrLB} }The first two implications follow 
from the fact that any one-sided deleted neighborhood is contained in the two-sided one. We prove the third implication. Let $b$ be a~limit 
point of $M$. Then there is a~sequence $(a_n)\sus M\setminus
\{b\}$ such that $\lim a_n=b$. It has a~subsequence $(a_{m_n})$ such 
that for every $n$ it holds that $a_{m_n}>b$ or for every $n$ it 
holds that $a_{m_n}<b$. Thus $b$ is a~one-sided limit point of $M$. 4. 
For instance $0\in L([0,1))$ but $0\not\in L^-([0,1))$.

\medskip\noindent
{\bf Exercise~\ref{ex_NekMnJednostr1} }It is, for example, the set $\N\sus\R$.

\medskip\noindent
{\bf Exercise~\ref{ex_NekMnJednostr} }We find in the set a~strictly monotone sequence.

\medskip\noindent
{\bf Exercise~\ref{ex_zpetneVyleps} }1. Here $a$ is a~limit point of 
$M$. Thus, as we know, it is the left or right limit point of $M$.
2. This follows from the inclusion $P(a,\de)\sus P^-(a,\de')\cup 
P^+(a,\de'')$ for $\de=\min(\de',\de'')$. 3. If 
$\lim_{x\to a}f(x)=A$, it suffices to take $\ep$ so small that 
$U(A,\ep)$ and $U(K,\ep)$, or $U(A,\ep)$ and $U(L,\ep)$, 
are disjoint and we get a~contradiction.

\medskip\noindent
{\bf Exercise~\ref{ex_jednozn_Jednostr} 
}The proof is similar to the proof of Proposition~\ref{prop_jednLimi}.

\medskip\noindent
{\bf Exercise~\ref{ex_Heine_Jednostr} 
}The proof is similar to the proof of Theorem~\ref{thm_HeinehoDef}.

\medskip\noindent
{\bf Exercise~\ref{ex_vztObyJed} }This follows from the equality $f[P^{\pm}(b,\de)]=(f\,|\,I^{\pm}(b))
[P(b,\de)]$.

\medskip\noindent
{\bf Exercise~\ref{ex_kyJeNesp} }There is a~sequence 
$(b_n)\sus M(f)$ such that $\lim b_n=b$, but $\lim f(b_n)$ does not exist or is not equal to 
$f(b)$. Or, by Exercise~\ref{ex_ekvivNespoVbode}, iff there is a~sequence $(b_n)\sus 
M(f)$ with $\lim b_n=b$ such that $\lim f(b_n)=A\ne f(b)$.

\medskip\noindent
{\bf Exercise~\ref{ex_klasDefSpovBod} }The solution of the inequalities 
$|x-b|<\de$, respectively $|f(x)-f(b)|<\ep$, is exactly $U(b,\de)$, 
respectively $U(f(b),\ep)$. Non-strict inequalities are equivalent: 
decrease $\de$ or $\ep$ a~little.

\medskip\noindent
{\bf Exercise~\ref{ex_HeiDefSpo} }If $b\in M\setminus L(M)$, there is a~$\de$ such that $U(b,\de)\cap 
M(f)=\{b\}$. Then the only sequences in $M(f)$ with the limit $b$ are the eventually constant 
sequences $(a_n)$ with $a_n=b$ for $n\ge n_0$. Then $\lim f(a_n)=f(b)$ which agrees with continuity of $f$ 
in such point $b$ (trivially or by Proposition~\ref{prop_spojVizol}).

\medskip\noindent
{\bf Exercise~\ref{ex_oIzBodech} }If $b\in M$ is isolated, it is not 
limit and there is a~$\ep$ such that  $M\cap P(b,\ep)=\emptyset$. Then 
$M\cap U(b,\ep)=\{b\}$. If $b\in M$ is not isolated, it is limit and for 
every $\ep$ we have that $M\cap P(b,\ep)\ne\emptyset$. Hence for 
every $\ep$ we have that $M\cap U(b,\ep)\ne\{b\}$.

\medskip\noindent
{\bf Exercise~\ref{ex_limAizBody} }This follows from the definition of isolated points.

\medskip\noindent
{\bf Exercise~\ref{ex_ekvivNespoVbode} }This is immediate from Heine's 
definition of continuity at a~point
and from part~3 of Theorem~\ref{thm_oPodposl}.

\medskip\noindent
{\bf Exercise~\ref{ex_ekvivSpoji} }Let $f$ be continuous at $b\in 
M(f)$ and let an $\ep$ be given. Then for some $\de$ it holds that $f[U(b,\de)]\sus U(f(b),\ep)$. Thus $f[U^-(b,\de)]\sus 
U(f(b),\ep)$ and $f[U^+(b,\de)]\sus U(f(b),\ep)$ because $U^-(b,\de)$ 
and $U^+(b,\de)$ are contained in $U(b,\de)$. So $f$ is both left- and 
right-continuous at $b$.

Let $f$ be both left- and 
right-continuous at $b\in M(f)$ and let an $\ep$ be given. Then there 
exist $\de'$ and $\de''$ such that $f[U^-(b,\de')]\sus U(f(b),\ep)$
and $f[U^+(b,\de')]\sus U(f(b),\ep)$. We set $\de
\equiv\min(\de',\de'')$. We get 
that $f[U(b,\de)]\sus U(f(b),\ep)$, because 
$U(b,\de)\sus U^-(b,\de')\cup U^+(b,\de'')$. Hence $f$ is continuous
at $b$.

\medskip\noindent
{\bf Exercise~\ref{ex_ustne} }Every number in $M$ is positive because 
$x$ is irrational. $U(x,1)$, which is an interval of length $2$, 
contains only finitely many fractions
with bounded denominators. Hence $M$ is a~finite set. But $U(x,1)$
contains at least one integer and so $M\ne\emptyset$.

\medskip\noindent
{\bf Exercise~\ref{ex_variantyLimMonFce} }Let $f\in\mathcal{F}(M)$. For $b\in 
L^-(M)$ and $f$ that weakly decreases on $P^-(b,\de)$ we replace
supremum with infimum. For $b\in L^+(M)$ and $f$ that weakly increases, 
resp. weakly decreases, on $P^+(b,\de)$ we take infimum, 
resp. supremum. For $+\infty\in L(M)$ 
and $f$ that weakly decreases on $U(+\infty,\de)$ we replace supremum
with infimum. For $-\infty\in L(M)$ and $f$ that weakly increases, 
resp. weakly decreases, on $U(-\infty,\de)$ we take infimum, 
resp. supremum.

\medskip\noindent
{\bf Exercise~\ref{ex_uloNaCauch} }Let $k=-\infty$, hence $f$ weakly decreases, and let a~$c<-1$ be 
given. We take an $h>0$ such that for every $x\ge h$ we have $f(x+1)-
f(x)\le c$. Let $n\in\N$ and $x=h,h+1,\ds,h+n-1$. By summing
these inequalities and rearranging the result we get
$\frac{f(h+n)}{n}\le c+
\frac{f(h)}{n}$. Let $d\equiv\max(|f(h)|,2(h+1))+1$ and $x\in\R$ be such that $x\ge h+d$. 
We take the unique $n\in\N$ such that $h+n\le x<h+n+1$. Then $n\ge|f(h)|$. By these inequalities and since $f$ weakly decreases, we have $\frac{f(x)}{n}\le
\frac{f(h+n)}{n}\le c+
\frac{f(h)}{n}\le c+1$. Since $\frac{n}{x}\ge1-\frac{h+1}{d}\ge\frac{1}{2}$ and $c+1<0$, we have 
$\frac{f(x)}{x}\le
(c+1)\frac{n}{x}\le\frac{c+1}{2}$. Hence $\lim_{x\to+\infty}
\frac{f(x)}{x}=-\infty$. For $k=+\infty$ the argument is similar. 

\medskip\noindent
{\bf Exercise~\ref{ex_AjeLB} }$M(f/g)=M(f)\cap M(g)\setminus 
Z(g)$ and there is a~$\de$ such that $Z(g)\cap U(A,\de)=\emptyset$. 

\medskip\noindent
{\bf Exercise~\ref{ex_LimFunReci0} }This is a~particular case of part~3 
of the theorem with $f=k_1$.

\medskip\noindent
{\bf Exercise~\ref{ex_zesLimFceUsp} }1. The previous proof of the theorem is easily modified, 
for $K<L$ there is an $\ep$ and real numbers $a,b$ such that $U(K,\ep)<
\{a\}<\{b\}<U(L,\ep)$. 2. 
Again the reversal of an implication. 

\medskip\noindent
{\bf Exercise~\ref{ex_verzeProJednostr} }The ordinary limits are only 
changed to one-sided. The proofs are basically reductions to ordinary 
limits by means of Proposition~\ref{prop_ojednLimi2}.

\medskip\noindent
{\bf Exercise~\ref{ex_genToMS} }First we generalize the 
notion of limits from $\mathcal{R}$ to maps between MS (metric space(s)). If 
$f\cc M\to N$, where $(M,d)$ and $(N,e)$ are MS, and $a\in M$, $b\in 
N$, then $\lim_{x\to a}f(x)=b$ means that for every $\ep$ there is 
a~$\de$ such that if $x\in M$ satisfies $0<d(x,a)\le\de$ then
$e(f(x),b)\le\de$. Theorem~\ref{thm_unexSque} 
generalizes as follows. If $a\in M$, $b\in N$ and $f,g\cc M\to N$ are 
maps between MS such that $\lim_{x\to a}f(x)=b$ and 
$\lim_{x\to a}e(f(x),g(x))=0$ (this is a~map from the MS $M$ to the MS $(\R,|x-y|)$), then also $\lim_{x\to 
a}g(x)=b$. 

\medskip\noindent
{\bf Exercise~\ref{ex_slozFceHeinem} }Let a~sequence $(a_n)\sus M(f(g))\setminus\{A\}$ 
have $\lim a_n=A$. By Heine's definition of limits of functions 
(HDLF), $\lim g(a_n)=K$. Suppose that condition~$1$ holds. Then for 
the $n$ with $g(a_n)=K$ we have that $f(g)(a_n)=f(g(a_n))=f(K)=L$. If 
there are infinitely many $n$ such that $g(a_n)\ne K$ then for the 
corresponding subsequence $(a_{m_n})$ we have by 
HDLF that $\lim f(g(a_{m_n}))=L$. Thus $\lim f(g)(a_n)=L$ and HDLF 
says that $\lim_{x\to A}f(g)(x)=L$. 
Suppose that condition~$2$ holds.
Then, deleting from $(g(a_n))$ finitely many terms, we may assume that $(g(a_n))\sus M(f)\setminus\{K\}$.
By HDLF, $\lim f(g(a_n))=L$. Again by HDLF, $\lim_{x\to A}f(g)(x)=L$. 
The case that none of the conditions holds is resolved via 
HDLF already in the original 
proof.

\medskip\noindent
{\bf Exercise~\ref{ex_druhyDkPrevHod} }The inner function $p(x)$ now in 
general is not injective, but we have the finiteness bound $|p^{-1}
(y)|\le\deg p$ for every $y\in\R$, and Theorem~\ref{thm_LimSlozFunkce} 
can be used because condition~2 
holds. 

\medskip\noindent
{\bf Exercise~\ref{ex_simpForThis} }For example, take $M\equiv\{\frac{1}{n}\cc\;n\in\N\}\cup\{1+\frac{1}{n}\cc\;n\in\N\}$ and define $f\in\mathcal{F}(M)$ by $f(\frac{1}{n})=f(1+\frac{1}{n})\equiv\frac{1}{n}$.

\medskip\noindent
{\bf Exercise~\ref{ex_fromProof} }Because $a_n\in U(A',\de_0)$, for 
$\de<\de_0$ we have $P(A,\de)\sus U(A,\de_0)$ and $U(A,\de_0)\cap 
U(A',\de_0)=\emptyset$.

\medskip\noindent
{\bf Exercise~\ref{ex_isInje} }If $x\ne y$ then $|f(x)-f(y)|\ge g(|x-y|)>0$ and $f(x)\ne f(y)$.

\medskip\noindent
{\bf Exercise~\ref{ex_lowFRegFunc} }Let $x<y$ be in $M$. By Lagrange's theorem we have for some $d\in(x,y)$ that
$$
|f(y)-f(x)|=|f'(d)|\cdot|y-x|\ge
c|y-x|\,.
$$

\medskip\noindent
{\bf Exercise~\ref{ex_UppFRegFunc} }This is similar to the previous 
proof, only the inequality is reversed.

\medskip\noindent
{\bf Exercise~\ref{ex_diifToInf} }It is easy to see that there is an
increasing sequence $(m_n)\sus\N$ such that $|a_n-b_{m_n}|\ge n$ for every $n\in\N$ (the subsequence of $(b_n)$ is running fast ahead of $(a_n)$).

\medskip\noindent
{\bf Exercise~\ref{ex_proNexThm} }Both properties follow easily from Proposition~\ref{prop_onLf}.

\medskip\noindent
{\bf Exercise~\ref{ex_almEquaEquiv} }A~pair $(a,b)$ is in $G_f\Delta 
G_g$ iff $a\in M(f)\Delta M(g)$ or $f(a)\ne g(a)$. 

\medskip\noindent
{\bf Exercise~\ref{ex_almEqu1} }Reflexivity and symmetry of 
$\doteq$ are clear. The transitivity follows from the fact that for any
three sets $A$, $B$ and $C$ we have $A\Delta C\sus A\Delta B\cup B\Delta 
C$.  

\medskip\noindent
{\bf Exercise~\ref{ex_equiBigO} }This is just an explicit unfolding of the boundedness of $\frac{f}{g}$ 
on $N$.

\medskip\noindent
{\bf Exercise~\ref{ex_Xfini} }Because if $f\doteq f_0$ and $g\doteq g_0$, then $\frac{f}{g}\doteq\frac{f_0}{g_0}$. 

\medskip\noindent
{\bf Exercise~\ref{ex_notAsym} }If $f=O'(g)$ (on $N$) and $a\in 
M(f)\cap M(g)\cap N\cap Z(g)$, then necessarily $a\in Z(f)$. Thus if we 
change $f$ to $f_0$ by replacing $f(a)=0$ with any nonzero value 
$f_0(a)$, then $f_0\ne O'(g)$ (on $N$). Hence $O'$ in general does not 
stand change in a~single value of the function and is not an asymptotic relation.

\medskip\noindent
{\bf Exercise~\ref{ex_naO} }1. Yes. 2. No (problem near $0$). 3. No  
(problem near $\pm\infty$). 4. Yes. 5. No (problem near $0$). 6. Yes.

\medskip\noindent
{\bf Exercise~\ref{ex_parOtazek} }1. Yes. 2. Yes. 3. No. 4. No. 5. Yes. 6. Yes. 

\medskip\noindent
{\bf Exercise~\ref{ex_Oslozo} }This follows from the fact that the limit 
$\lim_{x\to A}\frac{u(x)}{v(x)}=0$ is preserved when we multiply the numerator 
$u(x)$ by a~bounded factor.

\medskip\noindent
{\bf Exercise~\ref{ex_ExOne} }We assume that, for a~constant $c\ge0$, for every $x\in 
M(f)\cap M(h)\cap N\setminus Z(h)$ it holds that $|\frac{f}{h}(x)|\le 
c$, and the same holds with $g$ in place of $f$. Hence it holds for every $x\in 
M(f)\cap M(g)\cap M(h)\cap N\setminus Z(h)$ that $|\frac{f+g}{h}(x)|\le |\frac{f}{h}(x)|+
|\frac{g}{h}(x)|\le2c$.

\medskip\noindent
{\bf Exercise~\ref{ex_ExTwo} }For $f$ we have a~bound as previously and 
for every $x\in M(g)\cap N$ it holds that $|g(x)|\le c$. Hence it holds for every $x\in 
M(f)\cap M(g)\cap M(h)\cap N\setminus Z(h)$ that $|\frac{fg}{h}(x)|=|\frac{f}{h}(x)|\cdot
|g(x)|\le c^2$.

\medskip\noindent
{\bf Exercise~\ref{ex_ExThree} }This is similar to the previous exercise.

\medskip\noindent
{\bf Exercise~\ref{ex_ExFour} }We have that $\lim_{x\to A}
\frac{f}{h}(x)=0$ and $\lim_{x\to A}\frac{g}{h}(x)=0$. Hence the limit 
$\lim_{x\to A}\frac{f+g}{h}(x)$, which is defined because $A$ is 
a~limit point of $M((f+g)/h)$, equals 
$\lim_{x\to A}\frac{f}{h}+\lim_{x\to A}\frac{g}{h}=0+0=0$.

\medskip\noindent
{\bf Exercise~\ref{ex_ExFive} }Similarly, 
from $\lim_{x\to A}\frac{f}{h}(x)=0$ and a~bound that $|g(x)|\le c$ for 
every $x\in M(g)\cap P(A,\theta)$ we easily deduce that also $\lim_{x\to 
A}\frac{fg}{h}(x)=0$.

\medskip\noindent
{\bf Exercise~\ref{ex_ExSix} }This is similar to the previous exercise.

\medskip\noindent
{\bf Exercise~\ref{ex_ExSeven} }Like before we get from the assumptions on $f$, $g$, $h$ and 
$A$ that $\lim_{x\to A}\frac{f+g}{h}(x)=\lim_{x\to A}\frac{f}{h}(x)+\lim_{x\to A}
\frac{g}{h}(x)=1+0=1$.

\medskip\noindent
{\bf Exercise~\ref{ex_ExEight} }We get from the assumptions on $f$, $g$, $h$ and $A$ that $\lim_{x\to A}\frac{fg}{h}(x)=\lim_{x\to 
A}\frac{f}{h}(x)\cdot\lim_{x\to A}g(x)=1\cdot1=1$.

\medskip\noindent
{\bf Exercise~\ref{ex_ExNine} }This is similar to the previous exercise.

\medskip\noindent
{\bf Exercise~\ref{ex_delitele} }For $k\in\N$ with $k\le x$ the number of 
pairs $(m,n)\in\N^2$
with $mn=k$ equals $\tau(k)$.

\medskip\noindent
{\bf Exercise~\ref{ex_whyNotOn} }Because $\lim_{x\to1^+}x^{1/3}\log x=0$.

\medskip\noindent
{\bf Exercise~\ref{ex_procNe1} }No problem in our definition of big O, only we get no upper bound on $T_{\mathrm{HH}}(1)$. 

\medskip\noindent
{\bf Exercise~\ref{ex_onAsyExp} }This is an application of Proposition~\ref{prop_oJeO}.

\medskip\noindent
{\bf Exercise~\ref{ex_BernNum} }Moving the denominator $\exp x-1$ to the right, we get for $n\in\N_0$ that the coefficient $\sum_{k=0}^{n-1}\frac{B_k}{k!}\cdot
\frac{1}{(n-k)!}$ of $x^n$ equals $0$ for $n\ne1$, and $1$ for $n=1$. So for $n=1$ we get that $B_0=1$, and for $n\ge2$ that 
$B_{n-1}=-\frac{1}{n}\sum_{k=0}^{n-2}\binom{n}{k}B_k$. The second claim follows from the identity $f(-x)=f(x)$ for the formal power series
$f(x)\equiv\frac{x}{\mathrm{e}^x-1}+\frac{x}{2}$.

\medskip\noindent
{\bf Exercise~\ref{ex_equiConnec} }Suppose that $G$ is connected by
the new equivalent definition and
that $\{A,B\}$ is a~partition of $V$. We take any two vertices
$u\in A$ and $v\in B$. The walk joining $u$ and $v$ provides the 
edge that intersects both $A$ and 
$B$. If $G$ is not connected by
the new equivalent definition, we take any two (distinct) vertices 
$u,v\in V$ that cannot be joined by any walk, define $A$ as the set 
of vertices that can be reached from $u$ by a~walk, and set 
$B\equiv V\setminus A$ ($\ni v$). Then $\{A,B\}$ is a~partition of $V$
such that no edge in $E$ intersects both $A$ and $B$. Hence $G$ is not 
connected by the original definition.

\medskip\noindent
\centerline{{\bf 6 Continuous functions}}

\medskip\noindent
{\bf Exercise~\ref{ex_spojKon} }Use Proposition~\ref{prop_spojVizol}, every point of the definition domain is isolated.

\medskip\noindent
{\bf Exercise~\ref{ex_spojKonst} }For every 
$x\in\R$, $\de$ and $\ep$ we have $k_a[U(x,\de)]=\{a\}\sus U(k_a(x),\ep)=U(a,\ep)$.

\medskip\noindent
{\bf Exercise~\ref{ex_spojIden} }For every $a\in\R$ and given 
$\ep$ it suffices to set $\de=\ep$ 
because $x[U(a,\de)]=U(a,\de)$.

\medskip\noindent
{\bf Exercise~\ref{ex_hustAposl} }Let $N$ be dense in $M$ and $a\in 
M$. Using the axiom of choice\index{axiom!of choice, AC} 
we take for every 
$n$ a~point $b_n$ in $N\cap U(a,\frac{1}{n})$ and get a~sequence 
$(b_n)\sus N$ 
with $\lim b_n=a$. If $(b_n)\sus N$ has $\lim b_n=a\in M$ then for every 
$\de$ for every large $n$ it holds that 
$b_n\in U(a,\de)$.

\medskip\noindent
{\bf Exercise~\ref{ex_oveJsouHust} }Let $a<b$ be in $\R$. We 
take $n\in\N$ so large that $\frac{2\sqrt{2}}{n}<b-a$. It 
follows that for some $m\in\Z$ we have $\frac{m}{n}\in (a,b)$, and 
that for some $m\in\Z\setminus\{0\}$ we have $\frac{m\sqrt{2}}{n}\in 
(a,b)$. Thus every nontrivial interval contains both a~fraction 
and an irrational number. Hence $\Q$ and $\R\setminus\Q$ are dense in 
$\R$.

\medskip\noindent
{\bf Exercise~\ref{ex_naRidke} }Let $0\le a<b\le 1$ and $n\in\N$ be 
maximum such that $\frac{1}{n}\ge\frac{a+b}{2}$. Then the 
nontrivial interval $(\max(
\{a,\frac{1}{n+1}\}),\frac{a+b}{2})$ is contained in $(a,b)$ and is 
disjoint to $N$.

\medskip\noindent
{\bf Exercise~\ref{ex_whyDense} }Because then $f\,|\,N$ is an at most countable kernel of $f$.

\medskip\noindent
{\bf Exercise~\ref{ex_onBluThm} }$M=\R\setminus\Q$.

\medskip\noindent
{\bf Exercise~\ref{ex_wellOrdQ} }In ZF we easily define a~bijection $f\cc\N_0\to\Q$. Now $\N_0=\omega$ is well ordered by the relation $\in$.

\medskip\noindent
{\bf Exercise~\ref{ex_maleCvic} }If $A_m=\emptyset$ then $|\al-b|\le\frac{\ep}{2}$ for every $\al\in U(b_m,
\frac{1}{m})\cap M\cap\Q$. If we had for any of these fractions $\alpha$ that $|\al-b_m|\le
\frac{\ep}{2}$, we would have $|b_m-b|\le|b_m-\al|+|\al-b|\le\frac{\ep}{2}+\frac{\ep}{2}=\ep$, which is not the case.

\medskip\noindent
{\bf Exercise~\ref{ex_and_back} }The inverse of any bijection is a~bijection.  

\medskip\noindent
{\bf Exercise~\ref{ex_bijs} }The function $s$ is onto $\N$ because 
every natural number 
is a~product of an odd number and a~power of two. These expressions are unique: if 
$(2k-1)2^{l-1}=(2m-1)2^{n-1}$ then $l=n$, else $2$ would divide an odd 
number, so also $k=m$. Hence $s$ is injective.

\medskip\noindent
{\bf Exercise~\ref{ex_jeJasna} }
By Exercise~\ref{ex_spojKon} constants are continuous. Clearly, 
$k_a=k_b$ implies that $a=b$.

\medskip\noindent
{\bf Exercise~\ref{ex_dokazZoecn} }An injection from $\R$ to
$\mathcal{C}(M)$ is given by restrictions of constant functions. To get an
injection from $\mathcal{C}(M)$ to $\R^{\N}$ we use Proposition~\ref{prop_jadro}.

\medskip\noindent
{\bf Exercise~\ref{ex_obrJeInt} }Let $a<c<b$ with $a,c\in f[I]$. Then by 
Theorem also $c\in f[I]$. Hence $f[I]$ is an interval. 

\medskip\noindent
{\bf Exercise~\ref{ex_alpinista} }We prove that for every two
continuous functions $f,g\cc[0,1]\to\R$ satisfying
$f(0)=g(1)=0$ and $f(1)=g(0)=1$ there is a~$t\in(0,1)$ such that
$f(t)=g(t)$. We set $h\equiv 
f-g\cc[0,1]\to\R$ and use the theorem on intermediate values.

\medskip\noindent
{\bf Exercise~\ref{ex_abIsComp} }This follows from Theorems~\ref{thm_BolzWeier} and \ref{thm_limAuspo}.

\medskip\noindent
{\bf Exercise~\ref{ex_emptyIsComp} }The empty set is trivially compact. The set $\R$ is not compact.

\medskip\noindent
{\bf Exercise~\ref{ex_neextremy} }Continuity of both functions is easy to show. Neither function has maximum because for every
$x\in[0,1)$ and every $y\in(x,1)$ the value at $y$ is greater than that at $x$.

\medskip\noindent
{\bf Exercise~\ref{ex_onExtr} }If
$f$ attains maximum or minimum
everywhere, then for every $a,b$ in $M(f)$ we have both $f(a)\le f(b)$
and $f(a)\ge f(b)$, and $f$ is constant. Clearly, constant functions
attain maximum and minimum
everywhere. 

\medskip\noindent
{\bf Exercise~\ref{ex_vlUzMno} }$\emptyset$ and $\R$ are closed sets, 
and closed sets are closed to finite unions and arbitrary intersections. The
proof follows from passing to complements and using de Morgan formulas.

\medskip\noindent
{\bf Exercise~\ref{ex_opAndLimP} }This is clear from the fact that if
$b\in M$ then $U(b,\de)\sus M$ for some $\de$, and therefore
$P(b,\theta)\sus M$ for every $\theta\le\de$.

\medskip\noindent
{\bf Exercise~\ref{ex_vyhoPrsOko} }$[a,b]\setminus P(c,\de)=[a,b]\cap(\R\setminus((c-\de,c)\cup(c,c+
\de)))$ which is a~closed and bounded set.

\medskip\noindent
{\bf Exercise~\ref{ex_inMS} }The arguments in the proof of Theorem~\ref{thm_KompvR} easily extend to metric spaces.

\medskip\noindent
{\bf Exercise~\ref{ex_CountEinMS} }Any discrete metric space $(X,d)$ is closed 
and bounded, but if $X$ is infinite then $X$ is not compact because sequences 
$(x_n)\sus X$ with $x_m\ne x_n$ for $m\ne n$ do not have convergent subsequences.

\medskip\noindent
{\bf Exercise~\ref{ex_moreDexi} }$I_a$ is the  union of all open intervals that contain $a$ and are contained in $M$.

\medskip\noindent
{\bf Exercise~\ref{ex_explWhy} }The union of two intersecting open intervals is an open interval. 

\medskip\noindent
{\bf Exercise~\ref{ex_whyUnion} }Because for every $a\in M$ there is 
a~$b\in M\cap\Q$ near $a$ such that $a\in I_b$.

\medskip\noindent
{\bf Exercise~\ref{ex_CabtDis} }$C$ is closed because it is an intersection of closed sets. It is uncountable because it is exactly the set of those points in $[0,1]$ whose $3$-adic expansions use only digits $0$ and $2$. $C$ has ``length'' $0$ because $1-\frac{1}{3}-\frac{2}{9}-
\frac{4}{27}-\frac{8}{81}-\ds=0$.

\medskip\noindent
{\bf Exercise~\ref{ex_applBaire} }The stated property of points of $M$ means 
that every singleton set $\{b\}$, $b\in M$, is sparse in $M$. If $M$ were at 
most countable, the union $M=\bigcup_{b\in M}\{b\}$ would 
contradict Baire's theorem.

\medskip\noindent
{\bf Exercise~\ref{ex_BaireInMS} }Baire's theorem for metric spaces 
concerns complete spaces $(X,d)$. In such space every Cauchy sequence of 
points converges. A~set $Y\sus X$ is sparse (in $X$) if every ball $B\sus X$ 
has a~subball $B'\sus B$ such that $Y\cap B'=\emptyset$. The theorem then 
says that if $(X,d)$, $X\ne\emptyset$, is complete and $X=\bigcup_{n=1}^{\infty}X_n$, then 
there is an index $n$ such that $X_n$ is not sparse. The proof is similar to the 
one we gave in the case $(X,d)=(M,|x-y|)$ for nonempty closed $M\sus\R$.

\medskip\noindent
{\bf Exercise~\ref{ex_stespojJespoj} }Let $f\cc M\to\R$ be uniformly continuous, $c\in M$ and an 
$\ep$ be given. We take for this $\ep$ the $\de$ guaranteed by the uniform continuity. Then certainly 
$f[U(c,\de)]\sus U(f(c),\ep)$, so that $f$ is continuous at $c$. 

\medskip\noindent
{\bf Exercise~\ref{ex_passSubs} }If $(a_n),(b_n)\sus M$ for a~compact set $M$, we have convergent subsequences $(a_{k_n})$ and $(b_{m_n})$. For simpler notation we denote them again by $(a_n)$ and $(b_n)$.

\medskip\noindent
{\bf Exercise~\ref{ex_nejsouSS} }Let $a_n\equiv\frac{1}{n}$ and $b_n\equiv\frac{2}{\pi(2n-1)}$. Then 
$\lim(a_n-a_{n+1})=\lim(b_n-b_{n+1})=0$ but for every $n$ one has that $f(a_{n+1})-f(a_n)=1$ and $f(b_{n+1})-
f(b_n)=2$.

\medskip\noindent
{\bf Exercise~\ref{ex_HMCfun} }For example $f\equiv0$
on $[0,\frac{1}{\sqrt{2}}]\cap\Q$ and $f\equiv1$ 
on $[\frac{1}{\sqrt{2}},1]\cap\Q$.

\medskip\noindent
{\bf Exercise~\ref{ex_onClosure} }If $M$ is closed then the set of finite 
limits of sequences in $M$ coincides with $M$. If $M$ is not closed then
a~point outside of $M$ is in $\overline{M}$. 

\medskip\noindent
{\bf Exercise~\ref{ex_ImUCbou} }We take a~$\de$ for $\ep=1$ in the UC 
property of $f$, and a~finite set $N\sus M$ with the property that for 
every $a\in M$ there is a~$b_a\in N$
with $|a-b_a|\le 1$. Then for every $a\in M$ we have
$|f(a)|\le|f(a)-f(b_a)|+|f(b_a)|
\le1+\max(\{|f(x)|\cc\;x\in N\})$.

\medskip\noindent
{\bf Exercise~\ref{ex_explGenee} }This follows easily from the theorem on the 
relation of limits of functions and ordering. 

\medskip\noindent
{\bf Exercise~\ref{ex_spojPOLaRAC} }It follows from the arithmetic of continuity, from the definition of 
polynomials and rational functions, and from continuity of constants and identity.

\medskip\noindent
{\bf Exercise~\ref{ex_rov2nerov} }In the difference $a_n(x+c)^n-a_nx^n$ we use the binomial theorem, 
cancel $a_nx^n$ and take out $c$. The triangle inequality gives that 
$|a_n\sum_{i=1}^n\binom{n}{i}c^{i-1}x^{n-i}|$ 
is at most $|a_n|\sum_{i=1}^n\binom{n}{i}|c|^{i-
1}|x|^{n-i}$. We replace the numbers $|c|$ and $|x|$ by $d$ which is not 
smaller and is at least $1$, so that we can increase $i-1$ to $i$. Then we
bound the inner sum by $(d+d)^n$.

\medskip\noindent
{\bf Exercise~\ref{ex_odmFakt} }It suffices to show that for every $c>0$ we have that $\lim\frac{c^n}{n!}=0$. The sequence $(\frac{c^n}{n!})$ is 
nonnegative and, for large $n$, decreasing. Therefore it has the limit $d\ge0$ and there is an $m$ such that
$d=\inf(\{\frac{c^n}{n!}\cc\;n\ge m\})$. Suppose for the contrary that $d>0$. Then for large enough $n\ge m$ 
one has that $\frac{d(n+1)}{c}>\frac{c^n}{n!}$. Hence $d>\frac{c^{n+1}}{(n+1)!}$, which is a~contradiction, and $d=0$.

\medskip\noindent
{\bf Exercise~\ref{ex_naSpoSkla} }Let $b\in M(f(g))$ and let an $\ep$ be given. There is a~$\de$ such that
$f[U(g(b),\de)]\sus U(f(g)(b),\ep)$. There is a~$\theta$ such that $g[U(b,\theta)]\sus 
U(g(b),\de)$. Hence $f(g)[U(b,\theta)]\sus\ds\sus U(f(g)(b),\ep)$ and $f(g)$ is continuous at $b$ according to the definition.

\medskip\noindent
{\bf Exercise~\ref{ex_aniJedeNevyn} }The function $f$ given as $f(x)=x$ on $(0,1)$ and with the value 
$f(2)=1$ is continuous and increasing but the inverse $f^{-1}\cc(0,1]\to(0,1)\cup\{2\}$ is not continuous 
at $1$. The function $f$ 
with the values $f(0)=1$ and $f(n)=1-\frac{1}{n}$ for $n\in\N$ has the closed definition domain $\N_0\sus\R$ and 
is injective and continuous, but the inverse 
$f^{-1}\cc\{1-\frac{1}{n}\cc\;\;n\in\N\}\cup\{1\}\to\N_0$ is not continuous at $1$. 

\medskip\noindent
{\bf Exercise~\ref{ex_spojLogCyklFci} }$\log x$ is continuous by any of parts 2--4 of the theorem, $\arccos x$ and $\arcsin x$ by any of parts 1, 2 and 
4, and $\arctan x$ and $\mathrm{arccot}\,x$ by any of parts 2--4.

\medskip\noindent
\centerline{{\bf 7 Derivatives}}

\medskip\noindent
{\bf Exercise~\ref{ex_defDeri} }Use Corollary~\ref{cor_substPosun}.

\medskip\noindent
{\bf Exercise~\ref{ex_ojednDeri} }It is an instance of  Proposition~\ref{prop_ojednLimi}.

\medskip\noindent
{\bf Exercise~\ref{ex_jakoVtvrz} }Let $f\in\mathcal{F}(M)$ 
and $b\in M\cap L^{\pm}(M)$. Then we have the equivalence that 
$f'_{\pm}(b)=A$ 
$\iff$ $(f\,|\,I^{\pm}(b))'(b)=A$, with equal signs. In fact, it is an 
instance of  Proposition~\ref{prop_ojednLimi2}.

\medskip\noindent
{\bf Exercise~\ref{ex_oOLB} }The first claim follows easily from 
definitions. $0$ is a~limit point of the interval $(0,1)$ but it is not 
its two-sided limit point.

\medskip\noindent
{\bf Exercise~\ref{ex_jinyExer} }The displayed inequalities show that $c,d\ne b$. 

\medskip\noindent
{\bf Exercise~\ref{ex_Contra} }No, it does not. Neither $0$ nor $1$ is a~two-sided limit point of the domain $[0,1]$.

\medskip\noindent
{\bf Exercise~\ref{ex_whatCh} }Now $M(f)=\Z$ and 
$L^{\mathrm{TS}}(M(f))=\emptyset$. By Corollary~\ref{cor_NPELE}, every 
point $b$ in $M(f)$ is ``suspicious''.
In fact, $f$ has both local minimum and local maximum at every point of its domain. 

\medskip\noindent
{\bf Exercise~\ref{ex_fceVleVpr} }The function $f(x)$ is defined at 
$b$, the other function is not. This is the only point where they 
differ.

\medskip\noindent
{\bf Exercise~\ref{ex_deriSignu} }$\sgn_-'(0)=\lim_{x\to 
0^-}\frac{-1-0}{x-0}=+\infty$. Similarly $\sgn_+'(0)=+\infty$.
By item~2 in Exercise~\ref{ex_ojednDeri}, 
$\sgn'(0)=+\infty$.

\medskip\noindent
{\bf Exercise~\ref{ex_deriAbsHod} }$(|x|)'_-(0)=\lim_{x\to 
0^-}\frac{-x-0}{x-0}=-1$ and similarly for $(|x|)'_+(0)$.

\medskip\noindent
{\bf Exercise~\ref{ex_proNekoNepl} }For example, 
$\sgn'(0)=+\infty$ and 
$\sgn\,0=0\not\in L(\sgn[\R])=\emptyset$.

\medskip\noindent
{\bf Exercise~\ref{ex_pro0deroPla} }Use continuity of the function 
$f$ at $b$. 

\medskip\noindent
{\bf Exercise~\ref{ex_oneSiDiCo} }If $f\in\mathcal{R}$ and $f_-'(b)\in\R$ then $f$ is left-continuous at $b$, and similarly for $f_+'(b)\in\R$. Proofs are the same as for Proposition~\ref{prop_vlDerimplSpoj}.

\medskip\noindent
{\bf Exercise~\ref{ex_urciJednos} }The derivative 
$(\sqrt{x})'_-(0)$ is not defined because $0$ is not a~left limit point 
of the definition domain $[0,+\infty)$. Remaining 
one-sided derivatives are equal to the ordinary ones.

\medskip\noindent
{\bf Exercise~\ref{ex_derKon} }For every $x\ne b$ we have $\frac{k_c(x)-k_c(b)}{x-b}=\frac{0}{x-b}=0$. Hence $k_c'=k_0$.

\medskip\noindent
{\bf Exercise~\ref{ex_derIde} }For every $a$ we have $\lim_{x\to a}\frac{x+c-(a+c)}{x-a}=1$.

\medskip\noindent
{\bf Exercise~\ref{ex_spojResDer} }This is what Proposition~\ref{prop_vlDerimplSpoj} says.

\medskip\noindent
{\bf Exercise~\ref{ex_naGlobObou} }$D(f)\sus D_-(f)\cup D_+(f)$ and this inclusion is in general strict.

\medskip\noindent
{\bf Exercise~\ref{ex_procSe} }This is exactly the assumption that $a_n-b_n = o(b_n)$ ($n\to\infty$).

\medskip\noindent
{\bf Exercise~\ref{ex_jednStTec} }This follows from the limit 
$\lim_{x\to b}\frac{f(x)-f(b)}{x-b}=f'(b)$. 

\medskip\noindent
{\bf Exercise~\ref{ex_tecOdmoc} }For 
$a>0$ the equation is $y=\frac{x}{2\sqrt{a}}+\frac{\sqrt{a}}{2}$. At $\langle 0,0\rangle$ the tangent is not defined.

\medskip\noindent
{\bf Exercise~\ref{ex_whyThis} }Use the definition
of derivatives or results on derivatives of sums of functions.

\medskip\noindent
{\bf Exercise~\ref{ex_modifikace} }The required modifications are 
straightforward.

\medskip\noindent
{\bf Exercise~\ref{ex_bijNesPri} }We only need to show that for every 
non-vertical line
$\ell\in\mathcal{N}$ there exists a~unique pair $\langle s,t\rangle\in\R^2$ such
that $\ell=\{\langle x,sx+t\rangle\cc\;x\in\R\}$. Suppose that $\langle s',t'\rangle\in\R^2$ is 
another pair determining $\ell$. Then $t=s0+t=s'0+t'=t'$, hence $t=t'$. From $s+t=s1+t=s'1+t'=s'+t'$
we get $s=s'$.

\medskip\noindent
{\bf Exercise~\ref{ex_defKappa} }The system $sa+t=b$ $\&$ 
$sa'+t=b'$, with given $a\ne a'$, $b$, $b'$,
and unknowns $s$ and $t$, has a~unique solution whose component $s$ is given by the 
stated formula.

\medskip\noindent
{\bf Exercise~\ref{ex_bodNaLT} }We assume that $\ell(x)=sx+t$ is a~limit tangent to 
$G_f$ at $\langle b,f(b)\rangle$. We take any sequence $(b_n)\sus M(f)\setminus
\{b\}$ such that $b_n\to 
b$. Let the line $\kappa(b_n,f(b_n),b,f(b))$ be given 
as $y=s_nx+t_n$. Then $s_n\to s$, $t_n\to t$ and always 
$s_nb+t_n=f(b)$. Hence $f(b)=s_nb+t_n\to sb+t$. We see 
that $sb+t=f(b)$ and $\langle b,f(b)\rangle\in\ell$.

\medskip\noindent
{\bf Exercise~\ref{ex_prvNex} }This is a~straightforward check of a~rational identity. 

\medskip\noindent
{\bf Exercise~\ref{ex_oKL} }Let $b$, $M$ and $f$ be as stated. Using Heine's 
definition of pointwise continuity and Theorem~\ref{thm_oPodposl}, we take a~sequence $(x_n)$ converging to $b$ from one side
and such that $\lim f(x_n)\ne f(b)$. Then we take any 
sequence $(y_n)$ converging to $b$ from the other side and such that 
$\lim f(y_n)$ exists.

\medskip\noindent
{\bf Exercise~\ref{ex_limPrnejde} }Let $\ell_n(x)=s_nx+t_n$ and $\ell(x)=sx+t$. By the assumption, $\lim s_n=s$ and $\lim t_n=t$. But $c=sb+t$ and 
$c_n=s_nb+t_n$, so that $\lim c_n=\lim (s_nb+t_n)=(\lim s_n)\cdot b+\lim t_n=sb+t=c$.

\medskip\noindent
{\bf Exercise~\ref{ex_vybPodp} }Let $d_n$ be the infimum of the $\ep$ such that $\frac{y_n-c}{x_n-
b}\in U(A,\ep)$. We may clearly assume that always $y_n\ne c$ or that $A=0$. Then we easily define by induction an increasing sequence $(m_n)\sus\N$ such that always $\frac{y_n-
u_{m_n}}{x_n-z_{m_n}}\in U(A,d_n+\frac{1}{n})$. We are done because $\lim d_n=0$.

\medskip\noindent
{\bf Exercise~\ref{ex_uloNatecnu6} }We set $b=0$, $M=\R$, $f(x)= 
x^2\sin(\frac{1}{x})$ for $x\ne0$, $f(0)=0$, 
$x_n=\frac{2}{(4n-1)\pi}$ and $y_n=\frac{2}{(4n-3)\pi}$, for $n$ running in $\N$. 
Then $f'(0)=0$ and $x_n,y_n\to0$, but the secants $\kappa(x_n,f(x_n),y_n,f(y_n))$ 
have slopes $\frac{y_n^2+x_n^2}{y_n-x_n}\gg \frac{n^2}{n^2}=1$ ($n\in\N$).

\medskip\noindent
{\bf Exercise~\ref{ex_sameForMinus} }The adaptation is straightforward.

\medskip\noindent
{\bf Exercise~\ref{ex_sgnMinusOdm} }$(\sgn(x)-
\sqrt{x})'(0)=\lim_{x\to0}\frac{1-\sqrt{x}}{x}=\lim_{x\to0}\frac{1/\sqrt{x}-1}{\sqrt{x}}=\frac{+\infty}{0^+}=+\infty$.

\medskip\noindent
{\bf Exercise~\ref{ex_spojEf} }Use that $fg=gf$.

\medskip\noindent
{\bf Exercise~\ref{ex_uloNaLeiFor} }Thus $f'(0)=+\infty$, $g'(0)=-\infty$ and $f'(0)g(0)+f(0)g'(0)=(+\infty)\frac{1}{2}+(-\frac{1}{2})(-\infty)=
+\infty$. On the other hand, 
$(fg)(x)=-1$ for $x\ne0$ and $(fg)(0)=-\frac{1}{4}$ give $(fg)'_-(0)=+\infty$, 
$(fg)'_+(0)=-\infty$ and $(fg)'(0)$ does not exist.

\medskip\noindent
{\bf Exercise~\ref{ex_uloNaDerPodilu} }We modify the previous exercise 
and change the value $f(0)$ 
to $f(0)=\frac{1}{2}$. Then $\frac{f'(0)g(0)-f(0)g'(0)}{g(0)^2}=
((+\infty)\cdot\frac{1}{2}-\frac{1}{2}\cdot(-\infty))/\frac{1}{4}=+\infty$, but $(f/g)(x)=-1$ for 
$x\ne0$ and $(f/g)(0)=1$. The derivative $(f/g)'(0)$ again does not exist.

\medskip\noindent
{\bf Exercise~\ref{ex_bezSpoNrpl} }With $g$ not continuous at $b$ the 
formula $(f(g))'(b)=f'(g(b))\cdot g'(b)$
may not hold in the sense that the right-hand side is defined but the 
left-hand side is undefined. We set $f(x)\equiv x^2$, take $g(x)$ as
the modified signum with the value $g(0)=\frac{1}{2}$ and $b\equiv0$ 
($M=\R\sus L(M)=\R^*$). Then $f'(g(b))=(2x)(\frac{1}{2})=1$, 
$g'(b)=+\infty$ and the right-hand side is $1\cdot(+\infty)=+\infty$. However, $f(g)(x)=1$ for every $x\ne0$ and $f(g)(0)=\frac{1}{4}$, so that $(f(g))'(b)$ does not exist.

\medskip\noindent
{\bf Exercise~\ref{ex_CoVynechSpo} } All stated properties of $f$, except the last one, are very easy to see. 
Clearly, $f^{-1}(c)=f^{-1}(0)=0$.
We show that $(f^{-1})'(0)$ does not
exist by providing two sequences $(y_n)\sus f[M]\setminus\{c\}$ 
with $y_n\to c$ such that the two limits 
$\lim\frac{f^{-1}(y_n)-f^{-1}(c)}{y_n-c}=\lim\frac{f^{-1}(y_n)}{y_n}$
are different. For the first sequence $(y_n)\equiv(\frac{1}{n})$
the limit is $\lim
\frac{1+1/n}{1/n}=+\infty$. For the second sequence 
$(y_n)\equiv(\frac{2}{2n+1})$
the limit is $\lim
\frac{2/(2n+1)}{2/(2n+1)}=1$.

\medskip\noindent
{\bf Exercise~\ref{ex_zduvVseNero} }For $n=0$ we have
$\frac{1}{c}(a_n(x+c)^n-a_nx^n)=0$. For $n\ge1$ this difference equals
$a_n\sum_{j=0}^{n-1}(x+c)^j x^{n-1-j}$. We use the transformation
$\sum_{j=0}^n(x+c)^jx^{n-j}-(n+1)x^n=c\sum_{j=1}^n\sum_{i=1}^j\binom{j}{i}c^{i-1}x^{n-i}$ (by the binomial theorem), where $n\ge1$. 
Then we use the triangle inequality and the estimates $|c|,|x|\le y$.
Finally, the sum of binomial coefficients in the $j$-th row 
of the Pascal\index{Pascal, Blaise} triangle is $2^j$ and $y^{n-1}\le 
y^{n+1}$.

\medskip\noindent
{\bf Exercise~\ref{ex_deriLogAbs} }$(\log|x|)'=\frac{1}{x}$ ($\in\mathcal{F}(\R\setminus
\{0\})$).  

\medskip\noindent
{\bf Exercise~\ref{ex_derirealMoc} }1. The derivative $(a^x)'=(\exp((x\log a))'=\ds=a^x\cdot\log
a$ follows from Corollaries~\ref{cor_derSlozF} and \ref{cor_derExpKosSin}, 
and from the derivative $(k_c(x)\cdot\mathrm{id}(x))'=k_c(x)$. 2. For 
$b>1$ and $x>0$ we get the derivative from the expression 
$x^b=\exp(b\log x)$; for $x=0$ we use the definition of the derivative
and the limit 
$$
\lim_{x\to0}x^{-1}\exp(b\log x)
=\lim_{x\to0}\exp((b-1)\log x)=
\lim_{y\to-\infty}\exp y=0\,.
$$
3. This is trivial. 4. For $b<1$ we
use the expression $x^b=\exp(b\log x)$. 5. This is trivial. 6. We use 
the Leibniz formula. 7. This is trivial. 8. We use the expression 
$x^m=\frac{1}{x^{-m}}$ and the derivative of ratios. 

\medskip\noindent
{\bf Exercise~\ref{ex_derxnabLoc} }Now we have the limit 
$$
\lim_{x\to0}x^{-1}\exp(b\log x)
=\lim_{x\to0}\exp((b-1)\log x)=
\lim_{y\to+\infty}\exp y=+\infty\,.
$$

\medskip\noindent
{\bf Exercise~\ref{ex_tanCot1} }The derivative of $\tan x=\frac{\sin x}{\cos x}$ follows from the
derivatives $(\sin x)'=\cos x$ and
$(\cos x)'=-\sin x$, from the identity $\sin^2 x+\cos^2 x=k_1(x)$ and from
Corollary~\ref{cor_deriRat}. Similarly for $\cot x$.

\medskip\noindent
{\bf Exercise~\ref{ex_InvTri} }Let $f\equiv\sin x\,|\,
[-\frac{\pi}{2},\frac{\pi}{2}]$. Then $\arcsin x$ is the inverse of 
$f$. By part~1 of Theorem~\ref{thm_DeriInvfce} 
we have for $a\in(-1,1)$  that 
$${\textstyle
(\arcsin x)'(a)=\frac{1}{f'(\arcsin a)}
=\frac{1}{(\cos x\,|\,
[-\frac{\pi}{2},\frac{\pi}{2}])(\arcsin a)}=
\frac{1}{(\sqrt{1-f(\arcsin a)^2})}=
\frac{1}{\sqrt{1-a^2}}\,.
}
$$ 
We can differentiate $\arctan x$ similarly  or we can use the 
expression in terms of acsine in Exercise~\ref{ex_arctann}.  
Similarly for the derivatives of $\arccos x$ and $\mathrm{arccot}\,x$.

\medskip\noindent
{\bf Exercise~\ref{ex_derArcLoc} }This is an application, at the points $-1$ and $1$, of part~2 of Theorem~\ref{thm_DeriInvfce} to the inverse of $\sin x\,|\,
[-\frac{\pi}{2},\frac{\pi}{2}]$.

\medskip\noindent
{\bf Exercise~\ref{ex_Mer} }Because for $0\in\R$ we have $x^0=k_1(x)\,|\,(0,+\infty)$. 

\medskip\noindent
{\bf Exercise~\ref{ex_znovuTenPr} }For $x\ne0$ we have $f'(x)=2x\sin(\frac{1}{x})-
\cos(\frac{1}{x})$;  $f'(0)=\lim_{x\to0}\frac{1}{x}\cdot x^2\sin(\frac{1}{x})=\lim_{x\to0}x\sin(\frac{1}{x})=0$. This is a~discontinuous
function.

\medskip\noindent
{\bf Exercise~\ref{ex_onDerRat1} }Both sides of the equality are $\emptyset$, the empty rational function.

\medskip\noindent
{\bf Exercise~\ref{ex_onDerRat2} }If $k\ne-1$ then $r_k(x)=
\frac{x^{k+1}}{k+1}+k_c(x)$. For $k=-1$ no rational function $r_{-1}(x)$ with
$r_{-1}'(x)=1/x$ exists by Proposition~\ref{prop_poleRatDer}.

\medskip\noindent
{\bf Exercise~\ref{ex_noNeedAdd} }$\exp(b\log x)$. 

\medskip\noindent
{\bf Exercise~\ref{ex_ercsin0} }This follows from Proposition~\ref{prop_restrEF}.

\medskip\noindent
\centerline{{\bf 8 Applications of mean value theorems}}

\medskip\noindent
{\bf Exercise~\ref{ex_RolleAbsVal} }$f'(0)$ does not exist. 

\medskip\noindent
{\bf Exercise~\ref{ex_logicEquiv} }If we assume the left-hand side of 
the equivalence, then we can deduce from it the right-hand side with the help of the reduction 
of Lagrange's theorem to Rolle's that is described in the proof of 
Theorem~\ref{thm_Lagrange}. If we assume the right-hand side of 
the equivalence, then the left-hand side follows because Rolle's theorem 
is a~particular case of Lagrange's.

\medskip\noindent
{\bf Exercise~\ref{ex_easySee} }If $f(x)<l_u(x)$ for some $x\in[a,b]$, then by the continuity of $l_t(x)$ in $t$ (for fixed $x$) we would have
$f(x)<l_{u'}(x)$ for every
$u'\in(u-\ep,u]$. This is a~contradiction with the definition
of $u$ as a~supremum. We get, by 
increasing $u$ a~bit, a~similar contradiction if $f(x)>l_u(x)$ for 
every $x\in[a,b]$. 

\medskip\noindent
{\bf Exercise~\ref{ex_ifStrai} }Then every tangent line to $G_f$ has slope $s=\frac{f(b)-f(a)}{b-a}$.

\medskip\noindent
{\bf Exercise~\ref{ex_kdeByvadd} } If $f(a)=f(b)$ and
$g'(c)=\pm\infty$ then the product $z\cdot g'(c)=0\cdot(\pm\infty)$ is 
indefinite. However, one can modify the statement of the theorem by 
distinguishing the cases $f(b)\ne f(a)$ and $f(b)=f(a)$. In the former case we 
can allow $g'(c)\in\R^*$, and in the latter case we have a~$c\in(a,b)$ such that 
$f'(c)=0$.

\medskip\noindent
{\bf Exercise~\ref{ex_staciNnula0} }Add the dummy zero coefficients $p_{k+1}(x)=0$, $\ds$, $p_{n_0}(x)=0$.

\medskip\noindent
{\bf Exercise~\ref{ex_staciNnula} }Multiply the coefficients $p_i(x)$ by the polynomial
$(x-k)(x-k-1)\ds(x-n_0+1)$.

\medskip\noindent
{\bf Exercise~\ref{ex_onPoles} }Use Exercise~\ref{ex_ExSeven}.

\medskip\noindent
{\bf Exercise~\ref{ex_zeroNonz} }If $f\in\mathcal{R}$ is nonzero, then $f(b)\ne0$ for some $b\in M(f)$. Then $f^2=f\cdot f$, $M(f^2)=M(f)\cap M(f)=M(f)$, $b\in M(f^2)$ and $f^2(b)=f(b)^2\ne0$.

\medskip\noindent
{\bf Exercise~\ref{strMono} }First we select an injective sequence. 
Then we use Proposition.

\medskip\noindent
{\bf Exercise~\ref{ex_procPotr} }The denominator of $r(x)$ has finitely 
many zeros $\{z_1<z_2<\ds<z_l\}$ in $(k-1,+\infty)$, and a~gap $(z_{i-
1},z_i)$, $i\in[l+1]$, between them, where $z_0=k-1$ and 
$z_{l+1}=+\infty$,  contains a~tail $(a_r,a_{r+1},\ds)$ of $(a_n)$.

\medskip\noindent
{\bf Exercise~\ref{ex_jakRolle} }Suppose that $(a_n)$ decreases.
Then we can apply to any restiction
$f\,|\,[a_{n+1},a_n]$ Rolle's theorem, because $f$ is
differentiable on $[a_{n+1},a_n]$ and $f(a_{n+1})=f(a_n)=0$. We get
a~point $b_n\in(a_{n+1},a_n)$ with $f'(b_n)=0$. Similarly for increasing sequence $(a_n)$.

\medskip\noindent
{\bf Exercise~\ref{ex_procTypuu} }Clearly, $r'(x)\in\mathrm{RAC}$.
Suppose that $p_j(x)\ne0$.
Then we have $(p_j(x)\log(x-j+1))'=p_j'(x)\log(x-j+1)+p_j(x)\cdot\frac{1}{x-
j+1}\,|\,(j-1,+\infty)$. The latter summand is absorbed in $(r(x)\,|\,(c,+\infty))'$. If $\deg
p_j=0$, the former summand disappears. If $\deg p_j>0$, then
in the former summand we have $\deg p_j'=\deg p_j-1$.

\medskip\noindent
{\bf Exercise~\ref{ex_zobecNonPrek} }We just replace the function $\log(x-j+1)$ with $\log(x-j+1+c)$.

\medskip\noindent
{\bf Exercise~\ref{ex_infManZer} }Proceed as in the proof of Theorem~\ref{thm_konMnoKor}.

\medskip\noindent
{\bf Exercise~\ref{ex_onN10} }This follows from uniqueness of expansions of 
natural numbers in base $10$. See Theorem~... .

\medskip\noindent
{\bf Exercise~\ref{ex_someNonrec} }The set of recursive real numbers is countable.

\medskip\noindent
{\bf Exercise~\ref{ex_ExaAlgNum} }Any fraction $\al=\frac{a}{b}$ is 
a~root of $x-\al$ and of $bx-a$. The number $\sqrt{n}$ is a~root of $x^2-n$.

\medskip\noindent
{\bf Exercise~\ref{ex_ExaAlgNum1} }We get form (i) by dividing the polynomial 
by the leading coefficient. We get form (ii) by multiplying it 
by the product of denominators of its coefficients.

\medskip\noindent
{\bf Exercise~\ref{ex_FraAlgInt} }Exactly the integers $\Z$.

\medskip\noindent
{\bf Exercise~\ref{ex_GRAlgInt} }It is, $\phi^2-\phi+1$.  

\medskip\noindent
{\bf Exercise~\ref{ex_goldFibo} }We have the (Binet\index{Binet, Jacques P. M.}\index{Binet, Jacques P. M.!formula}) formula $F_n=\frac{1}{\sqrt{5}}\big(\phi^n-\psi^n\big)$ where $\psi$ is the other root of $x^2-x+1$.

\medskip\noindent
{\bf Exercise~\ref{ex_CantTrans1} }Every nonzero complex polynomial with degree
$d\in\N_0$ has at most $d$ complex roots (basically Exercise~\ref{ex_zeroCan}). Since the set of
integer polynomials is countable, the set of algebraic numbers is
also countable because it is a~countable union of finite sets. 
But the set $\R$ of real numbers is uncountable
(Theorem) and 
therefore also
the set of real transcendental numbers is uncountable, in 
particular nonempty. 

\medskip\noindent
{\bf Exercise~\ref{ex_deriBound} }$|(\sum_{j=0}^na_jx^j)'(x)|$, where 
$x\in[\al,\al+1$ with $\al=\frac{a}{10^k}$ for $a\in\Z$ and $k\in\N_0$, is by the 
triangle inequality at most $\sum_{j=0}^nj|a_j|(|a|+1)^{j-1}$. This is at most $(n+1)^2\max(\ds)
(|a|+1)^n$. We replace $n$ by $n+1$ in order to have only positive 
factors in the bound, and $n-1$ by $n$ in order to have nonnegative 
exponents.

\medskip\noindent
{\bf Exercise~\ref{ex_whyAcan} }The algorithm  computes the values 
$p_m(\al_m+\frac{j}{10^k})$ for $j=0$, $1$, $\ds$, $10^{k-k_m}-1$, which are always fractions, 
and selects the first $j$ for which the value is nonzero. Such $j$ exists 
because $p_m(x)$ is a~nonzero integral polynomial with sufficiently small 
degree.

\medskip\noindent
{\bf Exercise~\ref{ex_whyg0} }In
the division $g(x)=\frac{f(x)}{x-
\frac{p}{q}}$ the polynomial $f(x)$ loses one multiplicity of the root
$\frac{p}{q}$, but $\al\ne\frac{p}{q}$ as it is irrational. 

\medskip\noindent
{\bf Exercise~\ref{ex_whyd} }The constant $d$ exists because $f'$ is continuous and $I$ is compact.

\medskip\noindent
{\bf Exercise~\ref{ex_perDecExp} }Let $\frac{a}{b}$ with
$a\in\Z$ and $b\in\N$ be any fraction. We may assume that
$a\in\N$. Recall the school algorithm for computing the decimal
expansion of $\frac{a}{b}$, for example $\frac{1}{7}=1:7=0.142\ds$, 
with the residues $1,3,2,6,\ds$. Once a~residue repeats, the residues and the expansion start to repeat.

\medskip\noindent
{\bf Exercise~\ref{ex_violLiIn} }Let $n\in\N$ and $c>0$ be
arbitrary. We take $m\in\N$ large enough so that
$\frac{2}{q_m^{(m+1)/2}}<c$ and $\frac{m+1}{2}\ge n$. Then $|\lambda-\frac{z_m}{q_m}|<\frac{2}{q_m^{m+1}}=\frac{2}{q_m^{(m+1)/2}}\cdot\frac{1}{q_m^{(m+1)/2}}<\frac{c}{q_m^n}$ and Liouville's inequality is violated.

\medskip\noindent
{\bf Exercise~\ref{ex_LamEff} }Algorithm $\mathcal{L}$ determines for every input $n\in\N$ if $n=m!$ 
for some $m\in\N$. If it is the case, then $\mathcal{L}$ outputs the
digit $\mathcal{L}(n)=1$, and else it outputs $0$. In more detail, 
$\mathcal{L}$ multiplies the numbers $1,2,\ds,m$ as long as $m!\le n$. 
Since $m!\ge 2^{m-1}$, $\mathcal{L}$ knows the digit 
$\mathcal{L}(n)$ at the latest for $m\le\log_2(n+2)+1\le\log_2(6n)$.
Multiplying two numbers $\le n$ takes time 
$O(\log_2(n+1)^2)=O(\log_2(6n)^2)$.
Thus $\mathcal{L}$ computes the decimal digit $\mathcal{L}(n)$ ($\in
\{0,1\}$) in time $O(\log_2(6n)^3)$, which is time
polynomial in the size of $n$  (number of its digits). This 
(probably) cannot be achieved by the algorithm $\mathcal{A}$ in Theorem~\ref{thm_CantTrans2}. In 
any case, the description of $\mathcal{A}$ is much more 
complex than that of~$\mathcal{L}$. 

\medskip\noindent
{\bf Exercise~\ref{ex_moreGenLam} }The proof is very similar to the 
proof of Corollary~\ref{cor_lambTran}.

\medskip\noindent
{\bf Exercise~\ref{ex_naInterac1} }In the assumptions we just change $\ge$ to $\le$. In the conclusion we change, of course, $-\infty$ to $+\infty$ The proof of this version is similar to the original one. 

\medskip\noindent
{\bf Exercise~\ref{ex_naInterac2} }We define $f(x)\equiv x+1$ for 
$-1\le x<0$, $f(0)\equiv0$ and 
$f(x)\equiv x-1$ for $0<x\le1$. 

\medskip\noindent
{\bf Exercise~\ref{ex_otherThree} }These proofs are very similar to the one shown. 

\medskip\noindent
{\bf Exercise~\ref{ex_otherThree2} }These proofs are very similar to the one shown.

\medskip\noindent
{\bf Exercise~\ref{ex_jenJednop} }We cannot. The signum function $\sgn\,x$ shows, for $b\equiv0$, that these sets may have just
one element.

\medskip\noindent
{\bf Exercise~\ref{ex_jesteUlExtDe} }Take, for example, 
$k_0(x)\,|\,[a,b]$ and, to get $f(x)$, for any sequence $a<a_1<a_2<\ds<b$ with
$\lim a_n=b$ deform the flat $G_{k_0}$ around every point $(a_n,0)$ 
by a~small upward  bump. Thus a~function 
$f\in\mathcal{F}([a,b])$ arises that satisfies the assumptions of 
the proposition. Each bump is so low that the value
$f'(b)=k_0'(b)=0$ remains, but at the same time it is so steep that, say, $f'(a_n)=1$ for every $n$. The bumps are also so narrow that
$f'\big(\frac{a_n+a_{n+1}}{2}\big)=0$. Then $f'(b)=0$ but 
the limit $\lim_{x\to b}f'(x)$ does not exist.

\medskip\noindent
{\bf Exercise~\ref{ex_naHR2} }We compute by part~2 of HR~2 on the interval $(0,1)$ that
$$
{\textstyle
\lim_{x\to0}x^{\ep}\log x=
\lim_{\ds}\frac{\log x}{x^{-\ep}}=
\frac{1}{-\ep}\lim_{x\to0^+}
\frac{1/x}{x^{-\ep-1}}=
\frac{1}{-\ep}\lim_{x\to0^+}x^{\ep}=0\,.
}
$$

\medskip\noindent
{\bf Exercise~\ref{ex_otherDefDom} }We replace $f(x)$ with $f(a+b-x)$. The definition domain $P(b,\de)$ is the
union $(b-\de,b)\cup(b,b+\de)$. We move from the interval
$U(+\infty,\de)=(\frac{1}{\de},+\infty)$ and $x\to+\infty$ to the interval $(0,\de)$ and $y\to0$ 
by means of the substitution $y=\frac{1}{x}$.

\medskip\noindent
{\bf Exercise~\ref{ex_whyclear} }Because of the definitions of the functions 
$f(x)$, $g(x)$ and the derivative.

\medskip\noindent
{\bf Exercise~\ref{ex_secoDeri} }No, the former vale is always in $\R$, but the 
latter element may be $\pm\infty$.

\medskip\noindent
{\bf Exercise~\ref{ex_naDerKrad} }The former sequence is $4$-periodic: $(\overline{\sin x,\cos x,-\sin x,-\cos x})$. The $n$-th term of the latter sequence is
$(\frac{1}{x})^{(n)}=(-1)^nn!x^{-n-1}$. 

\medskip\noindent
{\bf Exercise~\ref{ex_locExtrSecD} }Consider $f(x)\equiv x^3$ and $b\equiv0$.

\medskip\noindent
{\bf Exercise~\ref{ex_givExa} }We set $b\equiv0$ and define $g(x)\equiv\sqrt{x}$ for $x\ge0$ and 
$g(x)\equiv-\sqrt{-x}$ for $x\le0$. Then the function
$${\textstyle
f(x)\equiv\int_{-\de}^x g(t)\,\mathrm{d}t\cc U(0,\,\de)\to\R
}
$$
has at $0$ a~strict local minimum, $f'=g$ on $U(0,\de)$ and $(f')'(0)=g'(0)=+\infty$.

\medskip\noindent
{\bf Exercise~\ref{ex_convConc1} }The strict convexity of $x^2$ 
follows from Theorem~\ref{thm_druDerKonvKonk} 
because $(x^2)''=2>0$. The claim about $|x|$ is clear
from the graph which is a~union of two half-lines. The strict concavity of $\log x$ follows from
Theorem~\ref{thm_druDerKonvKonk} because $(\log x)''=-x^{-2}<0$.

\medskip\noindent
{\bf Exercise~\ref{ex_convConc2} } This is logically clear from the definition. 

\medskip\noindent
{\bf Exercise~\ref{ex_convConc3} } This is clear from the definition by applying the symmetry $(x,y)\mapsto(x,-y)$ of the plane. 

\medskip\noindent
{\bf Exercise~\ref{ex_ConvEndpo} }The argument, based on Theorem~\ref{thm_limMonFce}, is the same as in the 
proof of Theorem~\ref{thm_exOnesDer}. Only the upper bound is now not available.

\medskip\noindent
{\bf Exercise~\ref{ex_ConvCont} }It is not true because the one-sided 
derivative at the endpoint may be infinite, and then continuity at the point is not guaranteed. For example, the function 
$f\in\mathcal{F}([0,1])$, given as $f(x)=0$ for $x\ne 1$ and $f(1)=1$, 
is convex. 

\medskip\noindent
{\bf Exercise~\ref{ex_proofLem} }It 
follows that $(c,c')$ lies above or on the line going through the first 
two points. Thus the second point $(b,b')$ lies below or on the line
going through the first point $(a,a')$ and the third point $(c,c')$. 

\medskip\noindent
{\bf Exercise~\ref{ex_onx3} }Indeed, the tangent at $(0,0)$ is the $x$-axis. 

\medskip\noindent
{\bf Exercise~\ref{ex_triviInfl} }Every point of the graph is an inflection point. 

\medskip\noindent
{\bf Exercise~\ref{ex_asymptote1} }These functions have one-sided 
limits $\pm\infty$ at $0$.

\medskip\noindent
{\bf Exercise~\ref{ex_asymptote2} }Let 
$\lim_{x\to+\infty}
(f(x)-sx-b)=0$. By adding the limit $\lim_{x\to+\infty}b=b$ we 
get that $\lim_{x\to+\infty}(f(x)-sx)=b$. Dividing by the limit $\lim_{x\to+\infty}x=+\infty$ we get that
$\lim_{x\to+\infty}(\frac{f(x)}{x}-s)=0$, thus
$\lim_{x\to+\infty}\frac{f(x)}{x}=s$.

Suppose that $\lim_{x\to+\infty}\frac{f(x)}{x}=s$
and $\lim_{x\to+\infty}(f(x)-sx)=b$. Subtracting from the latter 
limit the limit $\lim_{x\to+\infty}b=b$, we get 
that $\lim_{x\to+\infty}
(f(x)-sx-b)=0$.

\medskip\noindent
{\bf Exercise~\ref{ex_asymptote3} }It is the axis $x$.

\medskip\noindent
{\bf Exercise~\ref{ex_drawGr1} }$F\in\mathcal{F}(M)$ is even, resp. odd, if $M=-M$ 
($=\{-x:\;x\in M\}$) and for every $x\in M$ we have 
$F(-x)=F(x)$, resp. $F(-x)=-F(x)$. The function $F$ is
$c$-periodic ($c\in\R$) if $M=M\pm c$ ($=\{x\pm c:\;x\in M\}$) and for every $x\in M$ we have 
$F(x+c)=F(x)$. 

\medskip\noindent
{\bf Exercise~\ref{ex_GraDraw1} }{\bf 0. }$F\not\in\mathrm{EF}$. 
{\bf 1. }$M(F)=\R$. {\bf 2. }$F$ is an even and $1$-periodic function. 
The only periods are integers. {\bf 3. }From 
Proposition~\ref{prop_RiemannFunkce} we know that $F$ is continuous 
exactly at irrational numbers. We show that $F'(\al)$ does not exist 
for any $\al\in\R$. For rational $\al=\frac{m}{n}$ in lowest terms, 
$F(\al)=\frac{1}{n}$ and $F(x)=0$ for $x$ arbitrarily close to $\al$ 
both from the left and 
the right. This gives differential ratios at $\al$ going 
to $+\infty$ on the left, and to 
$-\infty$ on the right, of $\al$. If $\al$ is irrational then these 
zero values of $F(x)$ give diff. ratios at $\al$ equal 
to $0$ and arbitrarily close to $\al$. But by the theorem of 
Dirichlet\index{Dirichlet, Peter L.} there exist infinitely many 
different fractions $\frac{m}{n}$ in lowest terms such that $|\al-
\frac{m}{n}|<\frac{1}{n^2}$. These fractions give at $\al$ diff. 
ratios in absolute value  $>\frac{1/n}{1/n^2}=n\to+\infty$. 
{\bf 4. }It is not hard to see that $\lim_{x\to\al}F(x)=0$ for every 
$\al\in\R$ and that the limits $\lim_{x\to\pm\infty}F(x)$ do not 
exist. {\bf 5. }The function $F$ intersects the axis~$x$ exactly at 
the points $(\al,0)$ for irrational $\al$, and the axis $y$ at the point
$(0,1)$. {\bf 6. }We show that for irrational $\al\in\R$ the derivative 
$F'_+(\al)$ does not exist, and that 
$F'_+(\al)=-\infty$ for $\al\in\Q$. 
Suppose that $\al$ is irrational and that $\be\to\al^+$ via irrational 
$\be$. Then 
$\frac{F(\be)-F(\al)}{\be-\al}=0$. For rational Dirichlet's 
approximations $\be$ of $\al$ we get $\frac{F(\be)-F(\al)}{\be-\al}
\to-\infty$. Thus $F'_+(\al)$ does not exist. Let $\al=\frac{k}{l}$ be 
rational and in lowest terms. For irrational $\be$ we get 
$\frac{F(\be)-F(\al)}{\be-\al}=-\frac{1}{l}\cdot\frac{1}{\be-\al}
\to-\frac{1}{l}\cdot\frac{1}{0^+}=-\infty$. For rational 
$\be=\frac{m}{n}$ in lowest terms we similarly get 
$\frac{F(\be)-F(\al)}{\be-\al}=\frac{n^{-1}-l^{-1}}{\be-
\al}\to\frac{0-l^{-1}}{0^+}=-\infty$. Similarly, 
$F'_-(\al)$ does not exist for irrational $\al\in\R$ and 
$F'_-(\al)=+\infty$ for $\al\in\Q$.

{\bf 7. }Since $F$ is not monotone on any non-trivial interval, the
only maximal intervals of monotonicity are the singletons
$\{a\}$, $a\in\R$. As for the extremes, the function $F$ has at every irrational $\al$ 
non-strict global minimum, and at every $\al\in\Q$ strict local
maximum. Also, at every $\al\in\Z$ the function $F$ has non-strict
global maximum. $F$ has no other 
(local or global) extreme. 
{\bf 8. }The image of $F$ is the set $\{\frac{1}{n}\cc\;n\in\N\}\cup\{0\}$. {\bf 9. }Like in step~7, maximal intervals of convexity and concavity  are the singletons
$\{a\}$, $a\in\R$. {\bf 10. }Since there are no tangents, there are no
inflection points. {\bf 11. }$F$ has no asymptotes. {\bf 12. }Here are eleven topmost  points of $G_F$ in $[0,1]\times\R$: 

\medskip
\begin{picture}(50,50)(-40,0)
\put(-5,5){\line(1,0){240}}
\put(0,5){\line(0,1){5}}
\put(230,5){\line(0,1){5}}
\put(-2,15){{\footnotesize $0$}}
\put(228,15){{\footnotesize $1$}}
\put(-2,45){$\bullet$} 
\put(228,45){$\bullet$}
\put(113,25){$\bullet$}
\put(75,18){$\bullet$}
\put(151,18){$\bullet$}
\put(55,15){$\bullet$}
\put(170,15){$\bullet$}
\put(44,13){$\bullet$}
\put(90,13){$\bullet$}
\put(136,13){$\bullet$}
\put(182,13){$\bullet$}
\end{picture}

\medskip\noindent
{\bf Exercise~\ref{ex_GraDraw1apul} }The description of MDMs of $F(x)=r(x)$ is 
not simple. It suffices to describe those on which $F$ weakly increases. 
Since $F$ is even, from them we easily get those on which $F$ weakly
decreases. It is easy to see that the MDMs on which $F$ is constant 
is exactly the set 
$\{\al\in\R\cc\;\text{$\al$ is irrational}\}$, together with the 
sets $\{\frac{m}{n}
\in\Q\cc\;m\in\Z\wedge\mathrm{GCD}(m,n)=1\}$, $n\in\N$.

We describe the MDMs where $F$ weakly increases and is not constant. We introduce sequences of fractions in lowest terms on which $F$ is constant. Let 
$n\in\N$. A~$C(n,-\infty)$-sequence is a~left-infinite sequence 
of fractions $S=(\ds<\frac{m_3}{n}<\frac{m_2}{n}<\frac{m_1}{n})$ such 
that $S$ contains all 
lowest terms fractions 
$\frac{m}{n}\le\frac{m_1}{n}$. We set $l(S)\equiv\frac{m_1}{n}$ to be the last term of $S$. Similarly, 
$S=(\frac{m_1}{n}<\frac{m_2}{n}<\frac{m_3}{n}<\ds)$ is a~$C(n,
+\infty)$-sequence if $S$ contains all lowest terms fractions $\frac{m}{n}\ge \frac{m_1}{n}$. We set 
$f(S)\equiv\frac{m_1}{n}$ to be the first term of $S$. For $k\in\N$ with 
$k\ge2$ we say that $S=(\frac{m_1}{n}<\ds<\frac{m_k}{n})$ is a~$C(n,k)$ 
sequence if $S$ contains all lowest terms fractions 
$\frac{m}{n}$ with 
$\frac{m_1}{n}\le\frac{m}{n}\le\frac{m_k}{n}$. Again, $f(S)$ and $l(S)$ are the first and last term of $S$, respectively. Note that for every $S$ the restriction  $F\,|\,S$ 
is constantly $\frac{1}{n}$ and that every $S$ is internally MDM, 
no $\al\in\R$ lying inside $S$ but not in $S$ can be added to $S$ because then $F\,|\,S\cup\{\al\}$ is not monotone. 

We introduce similar sequences of fractions in lowest terms on which $F$ increases. 
An $I(-\infty)$-sequence is a~left-infinite sequence 
of fractions $S=(\ds<\frac{m_3}{n_3}<\frac{m_2}{n_2}<\frac{m_1}{n_1})$ 
such that $\ds>n_3>n_2>n_1$ 
and there is no fraction $\frac{m}{n}$ in lowest terms with 
$\frac{m_{i+1}}{n_{i+1}}<\frac{m}{n}<\frac{m_i}{n_i}$ and 
$n_{i+1}>n>n_i$. We set
$l(S)\equiv\frac{m_1}{n_1}$ to be the last term of $S$. We similarly define finite $k$-term, $k\in\N$ with $k\ge2$, $I(k)$-sequences 
$S=(\frac{m_1}{n_1}<\ds<\frac{m_k}{n_k})$ and their 
first and last terms $f(S)$ and $l(S)$. Again every such
$S$ is internally MDM. Note that that there is no 
$I(+\infty)$-sequence. An $I(-\infty)$-sequence can 
run to the left to $-\infty$, an example is given by the sequence
$(\ds<-\frac{n^2+1}{n}<\ds<-\frac{5}{2}<-\frac{2}{1})$. 

Now we describe the structure of any MDM $X$ ($\sus\R$) of $F$ on which 
$F$ weakly increases and is not constant. It has the form 
$$
X=AB_1B_2\ds B_kC
$$ 
where $A$ is the part going from $-\infty$, $k\in\N_0$, $B_1B_2\ds B_k$ 
is the alternating part consisting of finite segments, which is empty 
for $k=0$, and $C$ is the part going to $+\infty$. For the set $A$ we 
have three possibilities: (i) 
$A=C(n,-\infty)$ for some $n\in\N$, (ii) $A=I(-\infty)$ going to the left to $-\infty$ or (iii) $A=
\{\al\in\R\cc\;\text{$\al$ is irrational and $\al<c$}\}$ for some 
$c\in\R$. In $B_1B_2\ds B_k$, constant and increasing rational sequences alternate. The set $C$ is always a~$C(n,+\infty)$-sequence.

In case (i), we have $k=2j-1$ with $j\in\N$, $B_1$, $B_3$, 
$\ds$ are $I(k_1)$-, $I(k_3)$-, $\ds$
sequences and $B_2$, $B_4$, $\ds$ are $C(n_2,k_2)$-, $C(n_4,k_4)$-, 
$\ds$ sequences such that $l(A)=f(B_1)$, $l(B_1)=f(B_2)$, 
$\ds$, $l(B_{k-1})=f(B_k)$, and $C$ is a~$C(n_{k+1},+\infty)$-sequence 
such that $l(B_k)=f(C)$. The concatenation $AB_1B_2\ds B_kC$ is done in such a~way that the listed equal first and last elements are
identified. In case (ii), we have $k=2j$ with $j\in\N_0$, the 
alternating part may be empty, and the constant and increasing rational
sequences alternate in it in the other order. In case (iii), we have 
the subcase (iv) with $c$ irrational
and the subcase (v) with $c\in\Q$. 
In subcase (iv) we have $k=2j-1$ with $j\in\N$, $B_1$ is an $I(-
\infty)$-sequence with infimum $c$, $B_3$, $B_5$, $\ds$ are finite 
increasing rational sequences and $B_2$, $B_4$, $\ds$ are finite constant 
rational sequences. In subcase (v) we may have both $k=2j-1$ with 
$j\in\N$ and $k=2j$ with $j\in\N_0$,
but $f(B_1)$, respectively $f(C)$, equals $c$, and in $B_1B_2\ds B_k$
finite constant and increasing rational sequences alternate in the 
appropriate order.

\medskip\noindent
{\bf Exercise~\ref{ex_GraDraw2} }{\bf 0. }$F\in\mathrm{SEF}$.
{\bf 1. }$M(F(x))=M(\log x)=
(0,+\infty)$. {\bf 2. }This function is not of a~special form. {\bf 3. }We have $F\in\mathcal{C}$ 
and $D(F)=M(F)$ because $F\in\mathrm{SEF}$. The 
derivative is $F'(x)=(1+\log x)F(x)=(1+\log x)x^x$. {\bf 4. 
}$\lim_{x\to0}F(x)=\ds=1$ because $\lim_{x\to0}x\log x=0$. Clearly,  
$\lim_{x\to+\infty}F(x)=+\infty$. {\bf 5. }$G_F$ is disjoint to 
both coordinate axes. {\bf 6. }Since $F\in\mathrm{SEF}$, we 
have $D(F)=M(F)$ and there is nothing to compute. {\bf 7. }We 
equate the $F'$ found in step~3 to $0$ and get that $F'<0$ on 
$(0,\frac{1}{\mathrm{e}})$ and $F'>0$ on 
$(\frac{1}{\mathrm{e}},+\infty)$. Thus the maximal intervals of 
monotonicity of $F$ are the intervals $(0,
\frac{1}{\mathrm{e}}]$ and 
$[\frac{1}{\mathrm{e}},+\infty)$. On the former $F$ decreases and on 
the latter it increases. At $x=
\frac{1}{\mathrm{e}}$ the function $F$ has a~strict global minimum with the value 
$1/\mathrm{e}^{1/\mathrm{e}}$. It follows that $F$ has no other (local or global) extreme. In particular, there is no local maximum:
if $x\in(0,\frac{1}{\mathrm{e}}]$ then $F(y)>F(x)$ for every 
$y\in(0,x)$, and if $x\in[\frac{1}{\mathrm{e}},+\infty)$ then 
$F(y)>F(x)$ for every $y\in(x,+\infty)$. {\bf 8. }By 
steps~4 and~7 and continuity of $F$ the image of $F$ is the interval 
$[1/\mathrm{e}^{1/\mathrm{e}},+\infty)$. {\bf 9. }We have 
$F''(x)=(\frac{1}{x}+(1+\log x)^2)x^x$. Since $F''>0$ on 
$(0,+\infty)$, $F$ is convex on its definition domain. {\bf 10. }It 
follows from the previous step that $F$ has no inflection. {\bf 
11. }It is clear that $F$ has no vertical asymptotes. Since 
$\lim_{x\to+\infty}\frac{F(x)}{x}=+\infty$, by 
Exercise~\ref{ex_asymptote2} the function $F$ does not have asymptote at $+\infty$. {\bf 12. }\url{https://www.desmos.com/calculator}.

\medskip\noindent
\centerline{{\bf 9 Taylor polynomials and series. Real analytic functions}}

\medskip\noindent
{\bf Exercise~\ref{ex_proveEqui} }If $f(x)=f(b)+c(x-b)+o(x-b)$ ($x\to 
b$), then, by the definition of the symbol $o$, we have $\lim_{x\to 
b}\frac{f(x)-(f(b)+c(x-b))}{x-b}=0$.
This by arithmetic of functional limits is equivalent with 
$\lim_{x\to b}\frac{f(x)-f(b)}{x-b}=c$. In the other directions, if
$\lim_{x\to b}\frac{f(x)-f(b)}{x-b}=c$, then arithmetic of functional
limits gives that 
$${\textstyle
\lim_{x\to 
b}\frac{f(x)-(f(b)+c(x-b))}{x-b}=0\,.
}
$$
Hence $f(x)=f(b)+c(x-b)+o(x-b)$ ($x\to b$). 

\medskip\noindent
{\bf Exercise~\ref{ex_proveEqui22} }If $f$ is continuous at $b$, then 
$\lim_{x\to b}f(x)=f(b)$ (we are assuming that $b$ is not an isolated
point of $M$). Arithmetic of functional limits gives that
$\lim_{x\to b}\frac{f(x)-f(b)}{1}=0$.
Thus $f(x)=f(b)+o(1)$ ($x\to b$). The opposite implication is proven similarly. 

\medskip\noindent
{\bf Exercise~\ref{ex_iniSumTP} }This follows from the fact that if $m<n$ then $\sum_{j=m+1}^n a_j(x-b)^j+o((x-b)^n)=o((x-b)^m)$.

\medskip\noindent
{\bf Exercise~\ref{ex_whyfb} }The assumption of HR~2 requires that the limit is of the type $\frac{0}{0}$.

\medskip\noindent
{\bf Exercise~\ref{ex_halfNeighb} }l'Hospital rule HR~2 works for 
these one-sided neighborhoods,

\medskip\noindent
{\bf Exercise~\ref{ex_TaylExer} }$T_m^{p,0}=\sum_{j=0}^m a_jx^j$, 
where for $j>n$ we set $a_j\equiv0$.

\medskip\noindent
{\bf Exercise~\ref{ex_nonclaTp} }The coefficients $a_j$ in 
$T_n^{f,0}(x)=\sum_{j=0}^n a_jx^j$ are $a_0=a_1=a_2=a_3=0$, $a_4=1$ and 
$a_j=0$ for $j\ge 5$. As for the derivatives, $f^{(0)}(0)=f(0)=0$,
$f'(0)=0$ and $f^{(j)}(0)$ does not exist for any $j\ge2$ because $L(D(f))=L(\{0\})=\emptyset$. 

\medskip\noindent
{\bf Exercise~\ref{ex_TedToNepl} }Now $T_n^{f,0}(x)=\sum_{j=0}^n x^j$ and 
$\lim T_n^{f,0}(2)=+\infty$, but $f(2)=-1$.

\medskip\noindent
{\bf Exercise~\ref{ex_geomTail} }For $|x|<1$ we have 
$\sum_{j=0}^{\infty}x^j=\sum_{j=0}^n x^j+\frac{x^{n+1}}{1-x}$.

\medskip\noindent
{\bf Exercise~\ref{ex_idenBinKoe} }Just multiply $(-1)^j$ and the numerator of 
the generalized binomial coefficient on the left-hand side.

\medskip\noindent
{\bf Exercise~\ref{ex_LimByTayl1} }

\medskip\noindent
{\bf Exercise~\ref{ex_arErr1} }

\medskip\noindent
{\bf Exercise~\ref{ex_arErr2} }

\medskip\noindent
{\bf Exercise~\ref{ex_arErr3} }

\medskip\noindent
{\bf Exercise~\ref{ex_souTay} }This Taylor polynomial equals $(0+0)+(1+1)x+(0+0)x^2+(-\frac{1}{3}-\frac{1}{6})x^3=2x-\frac{1}{2}x^3$.

\medskip\noindent
{\bf Exercise~\ref{ex_prodTay} }

\medskip\noindent
{\bf Exercise~\ref{ex_ratioTayl} }

\medskip\noindent
{\bf Exercise~\ref{ex_najdiNaInt} }Near $0$ we have $\sqrt{1+\sin x}=\sin(\frac{x}{2})+\cos(\frac{x}{2})$.

\medskip\noindent
{\bf Exercise~\ref{ex_explEven} }$\tan x$ is an odd function.

\medskip\noindent
{\bf Exercise~\ref{ex_sTretiOdm} }

\medskip\noindent
{\bf Exercise~\ref{ex_pulkaOkoli} }

\medskip\noindent
{\bf Exercise~\ref{ex_coJe} }$\{\emptyset\}$.

\medskip\noindent
{\bf Exercise~\ref{ex_submOrPar} }Let $X\in\binom{[m+n]}{m}$. We define the 
maps $\kappa\cc[m]\to[m+n]$ and $\lambda\cc[n]\to[m+n]$ by $\kappa(j)\equiv$ the $j$-th element of 
$X$ (in the standard ordering of $\N$), and similarly by $\lambda(j)\equiv$  the $j$-th element of 
$[m+n]\setminus X$. Now we can
define the injection $\iota$ by
$$
\langle\langle A_1,\,\ds,\,A_k\rangle,\,\langle B_1,\,\ds,\,B_l\rangle,\,X\rangle
\mapsto
\langle\kappa[A_1],\,\lambda[B_1],\,
\kappa[A_2],\,\lambda[B_2],\,\ds\rangle\,.
$$

\medskip\noindent
{\bf Exercise~\ref{ex_invFps} }

\medskip\noindent
\centerline{{\bf 10 Primitives of UC functions}}

\medskip\noindent
{\bf Exercise }For example, $F(x)\equiv2\sqrt{x}$ and $f(x)\equiv1/\sqrt{x}$. Then $F'=f$ but $M(F)\ne M(f)$.

\medskip\noindent
\centerline{{\bf 14 More applications of Riemann integrals}}

\medskip\noindent
{\bf Exercise }For example, $g(t)=(\cos t,\sin t)$ for 
$t$ in $[0,\pi]$ and $h(x)=(x,f(x))$
for $x$ in $[-1,1]$.

\medskip\noindent
{\bf Exercise }For example, $f(t)=t(\overline{b}-\overline{a})
+\overline{a}$ for $t\in[0,1]$.

\medskip\noindent
{\bf Exercise }For example, $\varphi(t)=(t,2\cos(10\pi t),2\sin(10\pi t))$
for $t\in[0,1]$.

\medskip\noindent
{\bf Exercise }These are exactly constant maps from 
intervals $[a,b]$ to $\R^d$.

\medskip\noindent
{\bf Exercise }This length is $b-a$.

\medskip\noindent
{\bf Exercise }If $c\in(a,b)$ is such that 
$\varphi(c)\not\in\overline{a}\overline{b}$, then $\|\varphi(c)-
\varphi(a)\|+\|\varphi(b)-\varphi(c)\|>\|\overline{b}-
\overline{a}\|$. 

\medskip\noindent
{\bf Exercise }

}

\newpage

\addcontentsline{toc}{chapter}{Notation index}

\centerline{{\bf Page numbers of definitions of symbols}}

\label{znaceni}
\bigskip\noindent
$\langle A,<\rangle$ \dotfill\pageref{linearOrd}\\
$\langle a_1,a_2,\ds,a_k\rangle$ \dotfill\pageref{ordTuple}\\
$\mathfrak{a},\mathfrak{A},\ds,\mathfrak{Z}$ \dotfill \pageref{fracAlp}\\
$[A,B]$ \dotfill\pageref{ordPair}\\
$A\times B$ \dotfill\pageref{cartProd}\\
$\al,\be,\ds,\omega$ \dotfill \pageref{greeAlp}\\
AC \dotfill\pageref{AC}\\
$\binom{a}{j}$ \dotfill\pageref{anadje}\\
$[a]_R$ \dotfill\pageref{block}\\
$\arccos x$ \dotfill\pageref{arccos}\\
$\mathrm{arccot}\,x$  \dotfill\pageref{arccot}\\ 
$\arcsin x$ \dotfill\pageref{arcsin}\\
$\arcsin_0 x$ \dotfill\pageref{arcsinZero}\\
$\arctan x$ \dotfill\pageref{arctan}\\
$\sim$ (asymptotic equality) \dotfill\pageref{asymEqu}\\
$\approx$ (asymptotic expansion) \dotfill\pageref{approx}\\
$a^x$ \dotfill\pageref{anax}\\
$B(b,r)$ \dotfill\pageref{ball}\\
BEF \dotfill\pageref{BEF}\\
$B_k$ \dotfill\pageref{bernou}\\
$(b_n)\preceq(a_n)$ \dotfill\pageref{preceq}\\
$\mathcal{C}$ \dotfill\pageref{C}\\
$C$ \dotfill\pageref{Cantor}\\
$\odot$ (Cauchy product) \dotfill\pageref{CauchyPr}\\ 
$\mathcal{C}(M)$ \dotfill\pageref{Cm}\\
$\mathcal{C}^{\infty}(M)$ \dotfill\pageref{CmInf}\\
$C_n$ \dotfill\pageref{Catalan}\\
$\wedge$, $\&$ (conjunction) \dotfill\pageref{conjun}\\
$\cos t$ \dotfill\pageref{cos}\\
$\cot t$ \dotfill\pageref{cot}\\
$D(B)$ \dotfill\pageref{lowerB}\\
$:=$, $=:$ (define) \dotfill\pageref{defEqu}\\
$\deg p(x)$ \dotfill\pageref{degree}\\
$\deg f(x)$ \dotfill\pageref{degree2}\\
$\Delta(A,B)$ \dotfill\pageref{Delta}\\
$\vee$ (disjunction)  \dotfill\pageref{disjun}\\
$D(f)$ \dotfill\pageref{defDeri}\\
$D_{\pm}(f)$ \dotfill\pageref{defDeriPm}\\
$D(x)$ \dotfill\pageref{dOfx}\\
$\mathrm{e}$\dotfill\pageref{e}\\
EF \dotfill\pageref{EF}\\
EGF \dotfill\pageref{EGF}\\
$=$ (equality)\dotfill\pageref{equality}\\
$E(T^{f,b}(x))$ \dotfill\pageref{ETayl}\\
$\mathbb{E}\,X$ \dotfill\pageref{expec}\\
$\exp(x)$, $\mathrm{e}^x$ \dotfill\pageref{exp}\\
$\mathrm{ex}(u,n)$ \dotfill\pageref{exuen}\\
$\emptyset_{f}$ \dotfill\pageref{empReaFun}\\
$f_{\emptyset}$ \dotfill\pageref{empFun}\\
$f^{-1}$ \dotfill\pageref{inverse}\\
$f'$ \dotfill \pageref{derivGl}\\
$f_{\pm}'$ \dotfill \pageref{derivGlpm}\\
$f''$ \dotfill \pageref{deriv2}\\
$f'''$ \dotfill \pageref{deriv3}\\
$f\cc A\to B$ \dotfill\pageref{function}\\
$F_A(x)$ \dotfill \pageref{fax}\\
$F_A^0(x)$ \dotfill \pageref{fax0}\\
$f'(b)$ \dotfill\pageref{deriv}\\
$f_{\pm}'(b)$ \dotfill \pageref{derivPm}\\
$f[C]$ \dotfill\pageref{image}\\
$f^{-1}[C]$ \dotfill\pageref{preimage}\\
$f\,|\,C$ \dotfill\pageref{restric}\\
$f(g),f\circ g$ \dotfill\pageref{compos}\\
$f+g$ \dotfill\pageref{sumFun}\\
$f\cdot g$, $fg$ \dotfill\pageref{prodFun}\\
$f/g$ \dotfill\pageref{diviFun}\\
$f-g$ \dotfill\pageref{fminusg}\\
$f^{(k)}$ \dotfill \pageref{derivK}\\
$\mathrm{flim}\,A_n(x)$, $\mathrm{flim}_{n\to\infty}\,A_n(x)$ \dotfill \pageref{flim}\\
$\mathcal{F}(M)$ \dotfill \pageref{FM}\\
fps \dotfill \pageref{fps}\\
$\ga$ \dotfill \pageref{gamma}\\
$G_f$ \dotfill \pageref{graphFun}\\
$G=\langle V,E\rangle$ (graph) \dotfill \pageref{graph}\\
$H(B)$ \dotfill \pageref{upperB}\\
HDD \dotfill \pageref{HDD}\\ 
$h_n$ \dotfill\pageref{hen}\\
$I(b)^+$, $I(b)^-$ \dotfill\pageref{ibPlMi}\\
$\mathrm{id}_X$ \dotfill\pageref{identityFun}\\
$\inf$ \dotfill\pageref{inf}\\
$I_{\Q}$ \dotfill\pageref{IQ}\\
$I(S(x))$ \dotfill\pageref{interConv}\\
$\kappa$ \dotfill\pageref{kappaNum}\\
$\kappa(\mathcal{A})$ \dotfill\pageref{kappaA}\\
$\kappa(A,A')$, $\kappa(a,b,a',b')$ \dotfill\pageref{kappa}\\
$k_c(x)$ \dotfill\pageref{kcx}\\
$L(a_n)$ \dotfill \pageref{Lan}\\
$\mathcal{L}(f)$ \dotfill \pageref{LfM}\\
$\mathrm{LFP}(f)$ \dotfill \pageref{LFP}\\
$\mathrm{LFT}$ graph \dotfill \pageref{LFT}\\
$\lim a_n$ \dotfill \pageref{lim}\\
$\liminf$ \dotfill \pageref{liminf}\\
$\mathrm{Lim}\, \ell_n$ \dotfill \pageref{limLin}\\
$\limsup$ \dotfill \pageref{limsup}\\
$\lim_{x\to A}f(x)$ \dotfill \pageref{limFun}\\
$\lim_{x\to a^{\pm}}f(x)$ \dotfill \pageref{limFunOnes}\\
$(\lim_{x\to A})f(x)$ \dotfill \pageref{slimFun}\\
$L(M)$ \dotfill \pageref{LM}\\
$L^{\pm}(M)$ \dotfill \pageref{LMpm}\\
$\log(x)$ \dotfill\pageref{log}\\
$L^{\mathrm{TS}}(M)$ \dotfill \pageref{LMts}\\
$\overline{M}$ \dotfill\pageref{closure}\\
$M^0$ \dotfill\pageref{interior}\\
$\max(\cdot)$, $\min(\cdot)$  \dotfill\pageref{minmax}\\
MDM \dotfill\pageref{MDM}\\
$M(f)$ \dotfill\pageref{Mf}\\
MFF UK \dotfill\pageref{mffuk}\\
$m(n)$ \dotfill\pageref{meaNum}\\
$m_n$ \dotfill\pageref{ordFac}\\
$\cdots=\cdots\ (\mathrm{mod}\;x^k)$ \dotfill\pageref{equMod}\\
$\cdots\ \mathrm{mod}\;x^k$ \dotfill\pageref{redMod}\\
$M(w)$ \dotfill\pageref{Mw}\\
$\N$ \dotfill\pageref{natural}\\
$\mathcal{N}$ \dotfill\pageref{Ncal}\\
$\N_0$ \dotfill\pageref{natZer}\\
NCC \dotfill\pageref{NCC}\\
$\binom{n}{j}$ \dotfill\pageref{binom}\\
$O$, $\ll$ \dotfill\pageref{O}\\
$o$, $\omega$ \dotfill\pageref{o}\\
$\Omega$, $\Theta$, $\asymp$ \dotfill\pageref{Omega}\\
$\mathrm{OP}(k,n)$ \dotfill\pageref{OPkn}\\
$\mathrm{op}_{k,n}$ \dotfill\pageref{opkn}\\
$\mathrm{OP}(n)$ \dotfill\pageref{OPn}\\
$\mathrm{op}_n$ \dotfill\pageref{opn}\\
$P(A,\ep)$ \dotfill \pageref{PAeps}\\
$P^{\pm}(b,\,\ep)$ \dotfill \pageref{delNeiPm}\\
$\phi$ \dotfill \pageref{phi}\\
$\pi$ \dotfill\pageref{pi}\\
$\prod_{n=1}^{\infty}a_n$ \dotfill\pageref{infiProd}\\
$\pi(x)$ \dotfill\pageref{piOfx}\\
$(p_n)$ \dotfill \pageref{pen}\\
POL \dotfill\pageref{polynom}\\
$\pi(p,n)$ \dotfill\pageref{piePn}\\
$\Q[x]$ \dotfill \pageref{racPol}\\
$\mathcal{R}$ \dotfill \pageref{funcR}\\
RAC \dotfill \pageref{rational}\\
$\mathrm{RAC_{fi}}$ \dotfill \pageref{ratioFie}\\
RAF \dotfill \pageref{RAF}\\
RBEF \dotfill \pageref{RBEF}\\
$R^{f,b}_n(x)$ \dotfill\pageref{remainder}\\
$r_k(n)$ \dotfill \pageref{erken}\\
$\R^{\N}$ \dotfill \pageref{RtoN}\\
$\sqrt{x}$ (root of $x$) \dotfill\pageref{sqrt}\\
$\mathcal{R}_{\mathrm{smr}}$ \dotfill \pageref{Rsmr}\\
$R(S(x))$ \dotfill \pageref{radiConv}\\
$r(x)$ \dotfill \pageref{riemannFun}\\
$\R[[x]]$ \dotfill \pageref{Rx}\\
$\R[[x]]_{\mathrm{ri}}$ \dotfill \pageref{Rxri}\\
SCC \dotfill\pageref{SCC}\\
SEF \dotfill\pageref{SEF}\\
$\bigcap$ (set intersection) \dotfill\pageref{setInt}\\
$\bigcup$ (set sum) \dotfill\pageref{setSum}\\
$\sgn(x)$\dotfill\pageref{signum}\\
$\sum a_n$, $\sum_{n=1}^{\infty}a_n$,
$a_1+a_2+\cdots$ \dotfill\pageref{seriesSum}\\
$\sin t$ \dotfill\pageref{sin}\\
$(s_n)$ \dotfill\pageref{partSum}\\
$\sup$ \dotfill\pageref{sup}\\
$\tan t$  \dotfill\pageref{tan}\\
$\tau(n)$ \dotfill\pageref{tau}\\
$T^{f,b}(x)$\dotfill\pageref{Tfbx}\\
$T_{m,n}^{f,b}(x)$\dotfill\pageref{LauTaylorPol}\\
$T_n^{f,b}(x)$\dotfill\pageref{TaylorPol}\\
$T_{\mathrm{HH}}(n)$ \dotfill\pageref{Thh}\\
$t_k$ \dotfill\pageref{tk}\\
$T_{\mathrm{K}}(n)$ \dotfill\pageref{Tk}\\
$T=(V,E)$ (tournament) \dotfill\pageref{tourn}\\
$\mathfrak{T}_{\mathrm{SR}}$ \dotfill\pageref{tsr}\\
$U(b,\ep)$ \dotfill\pageref{neighb}\\
$U^{\pm}(b,\ep)$ \dotfill\pageref{neighbPm}\\
$\overline{U}(b,\ep)$ \dotfill\pageref{closedNeig}\\
UC \dotfill\pageref{UC}\\
$\mathcal{UC}$ \dotfill\pageref{UCcal}\\
$\mathcal{UC}(M)$ \dotfill\pageref{UCM}\\
$U(\pm\infty,\ep)$ \dotfill\pageref{neighbInf}\\
VBEF \dotfill\pageref{VBEF}\\
$x^0$ \dotfill\pageref{xnanula}\\
$0^x$ \dotfill\pageref{nulanax}\\
$x^b$\dotfill\pageref{xnab}\\
$x^m$ \dotfill\pageref{xnam}\\
$[x^m]A(x)$ \dotfill\pageref{extrCoef}\\
$\binom{X}{n}$ \dotfill\pageref{Xnadn}\\
$\zeta(s)$ \dotfill\pageref{zeta}\\
$Z(f)$ \dotfill\pageref{zeroFun}\\
$\Z[x]$ \dotfill \pageref{intePol}
\newpage

\addcontentsline{toc}{chapter}{Author and subject index}

{\small\printindex}
\end{document}